\documentclass[a4paper]{book}
\pdfoutput=1
\usepackage[latin1]{inputenc}
\title{Sheaf Theory Through Examples \\ (\textbf{Abridged Version})}
\author{Daniel Rosiak}
\date{December 12, 2020}
\usepackage[papersize={7in, 9in}, margin=1in]{geometry}
\usepackage{amsmath}
\usepackage{amsfonts}
\usepackage{amssymb}
\usepackage{graphicx, graphics}
\graphicspath{ {figures/} }
\usepackage[dvipsnames]{xcolor}
\usepackage[font=small, labelfont=bf]{caption}
\usepackage{gfsporson}
\usepackage{multirow}
\usepackage{forest}
\usepackage{makeidx}
\makeindex
\usepackage[cal=boondoxo]{mathalfa}
\usepackage[integrals]{wasysym}
\usepackage{amsthm}
\usepackage{textcomp}
\usepackage{mathrsfs}
\usepackage{pgfplots}
\usepgfplotslibrary{patchplots}
\pgfplotsset{compat=1.15}
\usepackage{tikz-cd}
\usepackage{tikz}
\usepackage[backend=biber,
style=alphabetic,
]{biblatex}
\usetikzlibrary{decorations.pathmorphing}
\usetikzlibrary{decorations.markings, decorations.shapes}
\usetikzlibrary{hobby}
\usetikzlibrary{arrows, patterns, shapes, automata, petri, positioning, calc, backgrounds, decorations, snakes, trees, bending, arrows.meta, fadings, intersections}
\usepackage{blkarray}
\let\amsamp=&
\usepackage{setspace}
\usepackage{mathtools}
\DeclarePairedDelimiter{\abs}{\lvert}{\rvert}

\begingroup
\catcode`\&=13
\gdef\pampmatrix{%
	\begingroup
	
	\begin{psmallmatrix}%
	}
	\gdef\endpampmatrix{\end{psmallmatrix}\endgroup}
\endgroup

\usepackage{rotating}

\usepackage{skak}
\usepackage{chessboard}
\storechessboardstyle{4x4}{maxfield=d4}
\setchessboard{boardfontsize=6pt, labelfontsize=4pt}

\usepackage{diagram}
\def\checkmark{\tikz\fill[scale=0.4](0,.35) -- (.25,0) -- (1,.7) -- (.25,.15) -- cycle;}

\theoremstyle{definition}
\newtheorem{example}{Example}[section]
\theoremstyle{definition}
\newtheorem{definition}{Definition}[section]
\theoremstyle{definition}
\newtheorem{proposition}{Proposition}[section]
\newtheorem{theorem}{Theorem}[section]
\theoremstyle{theorem}
\theoremstyle{definition}
\newtheorem{exercise}{Exercise}[section]

\newenvironment{boxedtikzcd}
{\begin{lrbox}{\boxedtikzcdbox}\begin{tikzcd}}
		{\end{tikzcd}\end{lrbox}\fbox{\usebox{\boxedtikzcdbox}}}
\newsavebox{\boxedtikzcdbox}

\usepackage{hyperref}
\begin{document}
	\maketitle 
	\chapter*{Preface} 
	After circulating an earlier version of this work among colleagues back in 2018, with the initial aim of providing a gentle and example-heavy introduction to sheaves aimed at a less specialized audience than is typical, I was encouraged by the feedback of readers, many of whom found the manuscript (or portions thereof) helpful; this encouragement led me to continue to make various additions and modifications over the years. \par 
	The project is now under contract with the MIT Press, which would publish it as an open access book in 2021 or early 2022. In the meantime, a number of readers have encouraged me to make available at least a portion of the book through arXiv. The present version represents a little more than two-thirds of what the professionally edited and published book would contain: the fifth chapter and a concluding chapter are missing from this version. The fifth chapter is dedicated to toposes, a number of more involved applications of sheaves (including to the ``$n$-queens problem" in chess, Schreier graphs for self-similar groups, cellular automata, and more), and discussion of constructions and examples from cohesive toposes.\par 
	
	Feedback or comments on the present work can be directed to the author's personal email, and would of course be appreciated.    
	\tableofcontents
	\newpage
	
	\chapter*{Introduction} 
	\addcontentsline{toc}{chapter}{Introduction}
	\section{An Invitation} 
	In many cases, events and objects are given to observation as extended through time and space, and so the resulting data is local and distributed in some fashion. For now, we can think of this situation in terms of data being indexed by, or attached to (``sitting over"), given regions or domains of some sensors. In saying that the data is \textit{local}, we just mean that it holds only within, or is only defined over, a certain region, i.e., its validity is restricted to a prescribed region or domain or reference context, and we expect that whenever a property holds at a point of its extended domain, then it also holds at ``nearby" points. We collect temperature and pressure readings and thus form a notion of ranges of possible temperatures and pressures over certain geographical regions; we record the fluctuating stockpile of products in a factory over certain business cycles; we accumulate observations or images of certain patches of the sky or the earth; we gather testimonies or accounts about particular events understood to have unfolded over a certain region of space-time; we build up a collection of test results concerning various parts of the human body; we amass collections of memories or recordings of our distinct interpretations of a certain score of music; we develop observations about which ethical and legal principles or laws are respected throughout a given region or network of human actors; we form a concept of our kitchen table via various observations and encounters, assigning certain attributes to those regions of space-time delimiting our various encounters with the table, where we expect that the ascribed properties or attributes are present throughout the entirety of a region of their extension. Even if certain phenomena are not intrinsically local, frequently its measurement or the method employed in data collection may still be local. \par 
	But even the least scrupulous person does not merely accumulate or amass local or partial data points. From an early age, we try to understand the various modes of \textit{connections} and \textit{cooperations} between the data, to patch these partial pieces together into a larger whole whenever possible, to resolve inconsistencies among the various pieces, to go on to build coherent and more global visions out of what may have only been given to us in pieces.
	As informed citizens or as scientists, we look at the data given to us on arctic sea-ice melting rates, on temperature changes in certain regions, on concentrations of greenhouse gases at various latitudes and various ocean depths, etc., and we build a more global vision of the changes to our entire planet on the basis of the connections and feedbacks between these various data. As investigators of a crime, we must ``piece together" a complete and consistent account of the events from the partial accounts of various witnesses. As doctors, we must infer a diagnosis and a plan of action from the various individual test results concerning the parts of a patient's body. We take our many observations concerning the behavior of certain networks of human actors and try to form global ethical guidelines or principles to guide us in further encounters. \par 
	Yet sometimes information is simply not local in nature. Roughly, one might think of such non-locality in terms of how, as perceivers, certain attributes of a space may appear to us in a particular way but then cease to manifest themselves in such a way over subparts of that space, in which case one cannot really think of the perception as being built up from local pieces. For a different example: in the game of Scrabble$^{\text{TM}}$,\index{Scrabble} one considers assignments of letters, one by one, to the individual squares in a lattice of squares, with the aim of building words out of such assignments. One might thus suspect that we have something like a ``local assignment" of data (letters in the alphabet) to an underlying space (15 $\times$ 15 grid of squares). Yet this assignment of letters to squares in order to form words is not really local in nature, since, while we do assign letters one by one to the grid of squares, the smallest unit of the game is really a \textit{legal word}, but not all sub-words or parts of words are themselves words, and so a given word (data assignment) over some larger region of the board may cease to be a word (possible data assignment) when we restrict attention to a subregion. \par  
	Even when information is local, there are many instances where we cannot synthesize our partial perspectives into a more global perspective or conclusion. As investigators, we might fail to form a coherent version of events because the testimonies of the witnesses cannot be made to agree with what other data or evidence tells us regarding certain key events. As musicians, we might fail to produce a compelling performance of a score because we have yet to figure out how to take what is best in each of our ``trial" interpretations of certain sections or parts of the entire score and splice them together into a coherent single performance or recording of the entire score. A doctor who receives conflicting information from certain test results, or testimony from the patient that conflicts with the test results, will have difficulty making a diagnosis. In explaining the game of rock-paper-scissors to children, we tell them that rock beats scissors, scissors beats paper, and paper beats rock, but we cannot tell the child how to win \textit{all the time}, i.e., we cannot answer their pleas to provide them with a global recipe for winning this game.\par 
	For distinct reasons, differing in the gravity of the obstacle they represent, we cannot always ``lift" what is local or partial up to a global value assignment or solution. A problem may have a number of viable and interesting local solutions but still fail to have even a single global solution. When we do not have the ``full story," we might make faulty inferences. Ethicists might struggle with the fact that it is not always obvious how to pass from the instantiations or particular variations of a seemingly locally valid prescription, valid or binding for (say) a subset of a network of agents, to a more global principle, valid for a larger network. In the case of the doctor attempting to make a diagnosis out of conflicting data, it may simply be a matter of either collecting more data, or perhaps resolving certain inconsistencies in the given test results by ignoring certain data in deference to other data. Other times, as in the case of rock-paper-scissors, there is simply nothing to be done to overcome the failed passage from the given local ranking functions to a global ranking function, for the latter simply does not exist. The intellectually honest person will eventually want to know if their failure to lift the local to the global is due to the inherent particularity or contextuality of the phenomena being observed or whether it is simply a matter of their own inabilities to reconcile inconsistencies or repair discrepancies in data-collecting methods so as to patch together a more global vision out of these parts. \par 
	\textit{Sheaf theory} is the roughly 70-year old collection of concepts and tools designed by mathematicians to tame and precisely comprehend problems with a structure exactly like the sorts of situations introduced above. The reader will have hopefully noticed a pattern in the various situations just described. We produce or collect assignments of data indexed to certain regions, where whenever data is assigned to a particular region, we expect it to be applicable throughout the entirety of that region. In most cases, these observations or data assignments come already distributed  in some way over the given network formed by the various regions; but if not, they may become so over time, as we accumulate and compare more local or partial observations. In certain cases, together with the given value assignments and a natural way of decomposing the underlying space, revealing the relations between the regions themselves, there may emerge correspondingly natural ways of restricting assignments of data along the subregions of given regions. In such cases, in this movement of decomposition and restriction, the glue or system of translations binding the various data together, permitting some sort of transit between the partial data items, becomes explicit; in this way, an internal consistency among the parts may emerge, enabling the controlled gluing or binding together of the local data into an integrated whole that now specifies a solution or system of assignments over a larger region embracing all of those subregions. Such structures of coherence emerging among the partial patches of local data, once explicitly acknowledged and developed, may enable a unique \textit{global} observation or solution, i.e., an observation that no longer refers merely to yet another local region but now extends over and embraces all of the regions at once; as such, it may even enable predictions concerning missing data or at least enable principled comparisons between the various given groups of data. Sheaves provide us with a powerful tool for precisely modeling and working with the sort of local-global passages indicated above. Whenever such a local-global passage is possible, the resulting global observations make transparent the forces of coherence between the local data points by exhibiting to us the principled connections and translation formulas between the partial information, making explicit the glue by which such partial and distinct clumps of data can be ``fused" together, and highlighting the qualities of the distribution of data. And once in this framework, we may even go on to consider systematic passages or translations between distinct such systems of local-to-global data. \par 
	On the other hand, when faced with \textit{obstructions} to such a local-global passage, we typically revise our basic assumptions, or perhaps the entire structure of our data, or maybe just our manner of assigning the data to our regions. We are usually motivated to do this in order to allow precisely such a global passage to come into view. When we can satisfy ourselves that nothing can be done to overcome these obstructions, we examine what the failure to pass from such local observations to the global in this instance can tell us about the phenomena at hand. \textit{Sheaf cohomology} is a tool used for capturing and revealing precisely obstructions of this sort. \par     
	The purpose of this book is to provide an inviting and (hopefully) gentle introduction to sheaf theory, where the emphasis is on explicit constructions, applications, and a wealth of examples from many different contexts. Sheaf theory is typically presented as a highly specialized and advanced tool, belonging mostly to algebraic topology and algebraic geometry (the historical ``homes" of sheaves), and sheaves accordingly have acquired a somewhat intimidating reputation. And even when the presentation is uncharacteristically accessible, emphasis is typically placed on abstract results, and it is left to the reader's imagination (or ``exercises") to consider some of the things they might be used for or some of the places where they can be found. This book's primary aim is to dispel some of this fear, to demonstrate that sheaves can be found all over (and not just in highly specialized areas of advanced math), and to give a wider audience of readers a more inviting tour of sheaves. Especially over the last few years, the interest in sheaves among less and less specialized groups of people appears to be growing immensely; but, whenever I spoke to newcomers to sheaves, I invariably heard that the existing literature was either too specialized or too forbidding. This book accordingly also aims to fill a gap in the existing literature, which for the most part tends to either focus exclusively on a particular use of sheaves or assumes a formidable pre-existing background and high tolerance for abstraction. I do not share the view that applications or concrete constructions are mere corollaries of theorems, or that examples are mere illustrations with no power to inform ``deeper" conceptual advances. I am not sure if I would go as far as to endorse Vladimir Arnold's \index{Vladimir Arnold} idea that ``The content of a mathematical theory is never larger than the set of examples that are thoroughly understood," but I do believe that one barrier to the wider recognition of the immense power of sheaf theory lies in the tendency to present much of sheaf theory as if it were a forbiddingly abstruse or specialized tool, or as belonging mainly to one area of math. One thing this book aims to show is that it is no such thing. Moreover, well-chosen examples are not only useful, both pedagogically and ``psychologically," in helping newcomers get a better handle on the abstract concepts and advance forwards with more confidence, but can even jostle experts out of the rut of the `same old examples' and present interesting challenges both to our fundamental intuitions of the underlying concepts and to preconceptions we might have about the true scope of applicability of those concepts. \par 
	Before outlining the contents of the book, the next section offers a more detailed, but still ``naive," glimpse into the \textit{idea} of a sheaf via a toy construction, with the aim of better establishing intuitions about the underlying sheaf idea.  
	\section{A First Pass at the Idea of a Sheaf}
	Suppose we have some `region', which, for the moment, we can represent very naively and abstractly as 
	\begin{center}
		\includegraphics[scale=0.4]{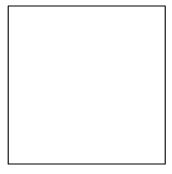}
	\end{center}
	We are less interested in the ``space itself" and more in how the space serves as a site where various things \textit{take place}. In other words, we think of this region as really just an abstract domain supporting various \textit{happenings}, where such happenings carry information for appropriate sensors\index{sensors} or ``measuring instruments" (in a very generalized sense), so that interrogating the space becomes a matter of asking the sensors about what is happening on the space.\footnote{The description of sheaves as ``measuring instruments" or the ``meter sticks" on a space that we are invoking---so that the set of all sheaves on a given space supply one with an arsenal of all the meter sticks measuring it, yielding ``a kind of `superstructure of measurement'"---ultimately comes from Grothendieck,\index{Alexander Grothendieck} who was largely responsible for many of the key ideas and results in the early development of sheaf theory. In speaking of (another early sheaf theorist) Jean Leray's\index{Jean Leray} work in the 40s, Grothendieck said this: 
		\begin{quote}
			The essential novelty in his ideas was that of the (Abelian) sheaf over a space, to which Leray associated a corresponding collection of cohomology groups (called ``sheaf coefficients"). It is as if the good old standard ``cohomological metric" which had been used up to then to ``measure" a space, had suddenly multiplied into an unimaginably large number of new ``meter sticks" of every shape, size and form imaginable, each intimately adapted to the space in question, each supplying us with very precise information which it alone can provide. This was the dominant concept involved in the profound transformation of our approach to spaces of every sort, and unquestionably one of the most important mathematical ideas of the 20th century. (\cite{grothendieck_recoltes_1986}, Promenade 12) 
	\end{quote}
Then the sheaves on a given space will incorporate ``all that is most essential about that space...in all respects a lawful procedure [replacing consideration of the space by consideration of the sheaves on the space], because it turns out that one can ``reconstitute" in all respects, the topological space by means of the associated ``category of sheaves" (or ``arsenal" of measuring instruments)...[H]enceforth one can drop the initial space...[W]hat really counts in a topological space is neither its ``points" nor its subsets of points, nor the proximity relations between them; rather, it is the \textit{sheaves on that space, and the category that they produce}" (\cite{grothendieck_recoltes_1986}, Promenade 13). The reader for whom this is overwhelming should press on and rest assured that we will have a lot more to say about all this later on in the book, and the notions and results alluded to in the above will be motivated and discussed in detail.} For instance, the region might be the site of some happenings that supply \textit{visual information}, so that as a sensor monitors the happenings over a region (or some part of it), it collects specifically visual information about whatever is going on in the area of its purview:   
	\begin{center}
		\includegraphics[scale=0.25]{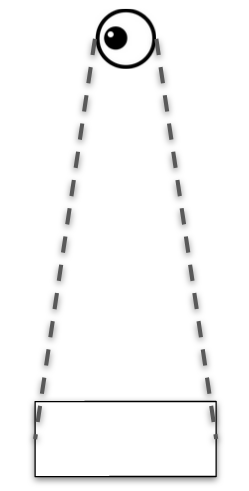}
	\end{center}
	There might then be another sensor, taking in visual information about another region or part of some overall `space', offering another ``point of view" or ``perspective" on another part of the space; and it may be that the underlying regions monitored by the two sensors overlap in part:
	\begin{center}
		\includegraphics[scale=0.25]{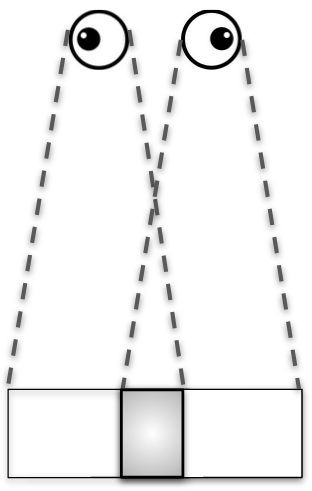}
	\end{center}
	Since we are ultimately interested in the informative happenings on the space, we want to see how the distinct ``perspectives" on what is happening throughout the space are themselves related; to this end, a very natural thing to do is ask how the data collected by such neighboring sensors are related. Specifically, a very natural thing to ask is whether and how the perspectives are \textit{compatible} on such overlapping sub-regions, whenever there are such overlaps between the underlying regions over which they, individually, collect data. \par 
	A little more explicitly: if we assume the first sensor collects visual data about its region (call it $U_1$), we may imagine, for concreteness, that the particular sort of data available to the sensor consists of sketches, say, of characters or letters (so that the underlying region acts as some sort of generalized sketchpad or drawing board)   
	\begin{center}
		\includegraphics[scale=0.25]{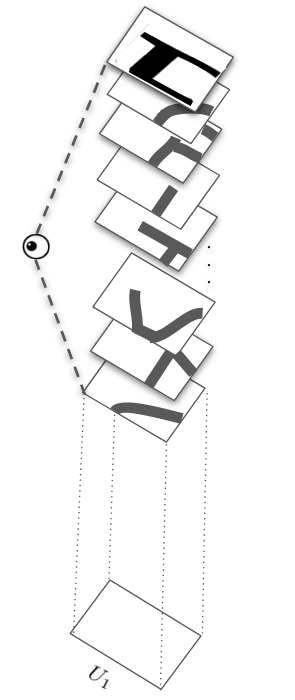}
	\end{center}
	While not really necessary, the sensor might even be supposed to be equipped to ``process" the information it collects, translating such visual inputs into reasonable guesses about which possible capital letter or character the partial sketch is supposed to represent. In any event, attempting to relate the two ``points of views" by considering their compatibility on the region where their two surveyed regions overlap, we are really thinking about first making a selection from each of the collections of data assigned to the individual sensors:   
	\begin{center}
		\includegraphics[scale=0.25]{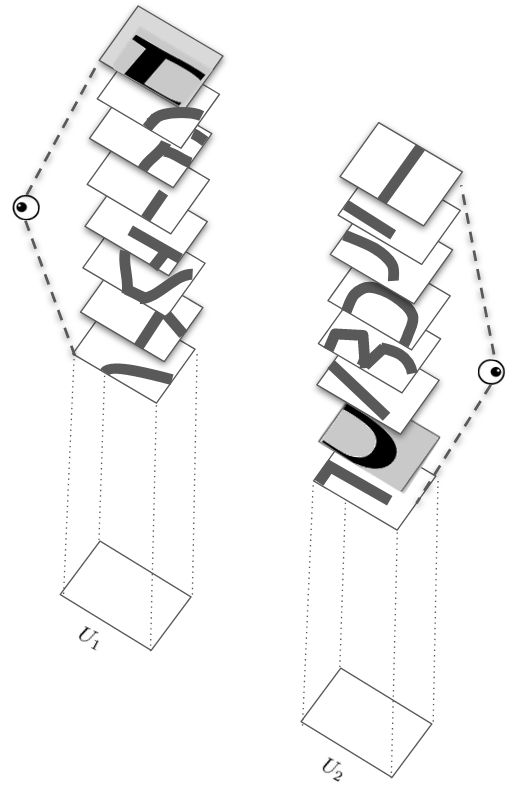}
	\end{center}
Corresponding to how the underlying regions are naturally related by an ``inclusion" relation, the compatibility question, undertaken at the level of the selections (highlighted in gray above) from the collections of all informative happenings on the respective regions, will involve looking at whether those data items ``match" (or can otherwise be made ``compatible") when we restrict attention to that region where the individual regions monitored by the separate sensors overlap:  
	\begin{center}
		\includegraphics[scale=0.25]{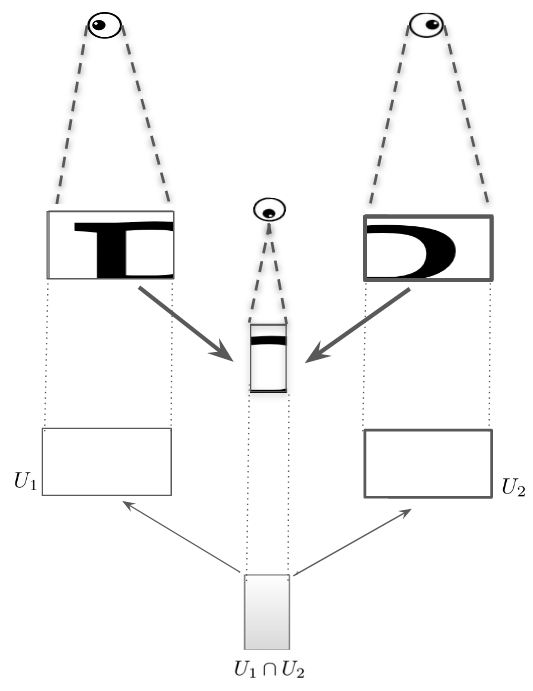}
	\end{center}
	If the given selection from what they individually ``see" does match on the overlap, then, corresponding to how the regions $U_1$ and $U_2$ may be joined together to form a larger region, 
	\begin{center}
		\includegraphics[scale=0.25]{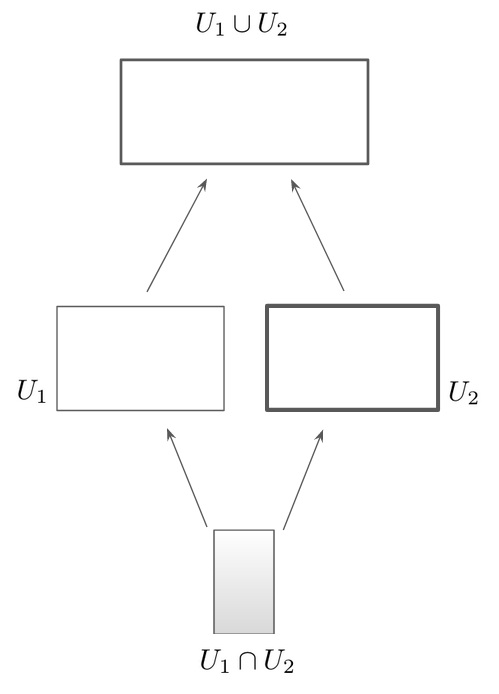}
	\end{center}
	at the level of the data on the happenings over the regions, we can pull this data back into an item of data given now over the entire space $U_1 \cup U_2$, with the condition that we expect that restricting this new, more comprehensive, perspective back down to the original individual regions $U_1$ and $U_2$ will give us back whatever the two individual sensors originally ``saw" for themselves: 
	\begin{center}
		\includegraphics[scale=0.3]{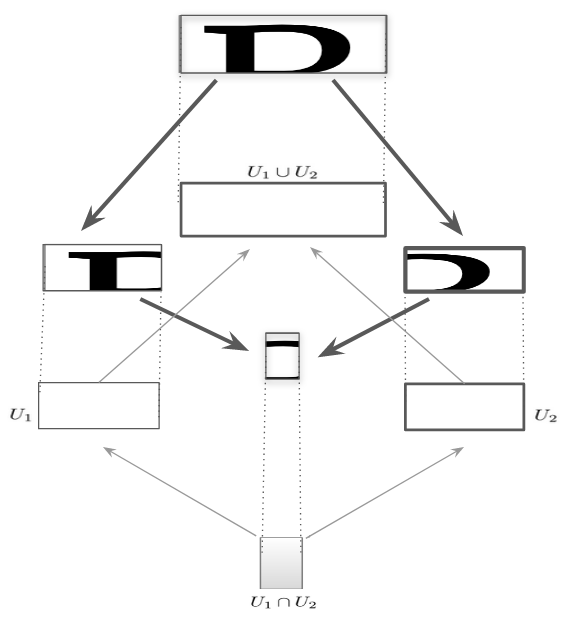}
	\end{center}
	In other words, given some selection from what sensor 1 ``sees" as happening in its region $U_1$ and from what sensor 2 ``sees" as happening in its region $U_2$, provided their ``story" agrees about what is happening on the overlapping region $U_1 \cap U_2$, then we can paste their individual visions into a single and more global vision or story about what is happening on the overall region $U_1 \cup U_2$ (and we expect that this story ultimately ``comes from" the individual stories of each sensor, in the sense that restricting the ``global story" down to region $U_1$, for instance, will recover exactly what sensor 1 already saw on its own). \par 
	Another way to look at this is as follows: while the sensor on the left, when left to its own devices, will believe that it may be seeing a part of any of the letters $\{B, E, F, P, R\}$, checking this assignment's compatibility with the sensor on the right amounts to constraining what the left sensor believes by what the sensor on the right ``knows," in particular that it cannot be seeing an $E$ or an $F$. Symmetrically, the sensor on the right will have its own ``beliefs" that might, in the matching with the left sensor, be constrained by whatever the left sensor ``knows." In matching the two sensors along their overlap, and patching their perspectives together into a single, more collective, perspective now given over a larger region (the union of their two regions), we are letting what each sensor individually ``knows" constrain and be constrained by what the other ``knows." \par 
	In this way, as we cover more and more of a `space' (or, alternatively, as we decompose a given `space' into more and more pieces), we can perform such compatibility checks at the level of the data on the happenings on the `site' (our collection of regions covering a given space), and then ``glue together," piece by piece, the partial perspectives represented by each sensor's local data collection into more and more embracing or ``global" perspectives. More concretely, continuing with our present example, suppose there are two additional regions, covering now some southwest and southeast regions, respectively, so that, altogether, the four regions cover some region (represented by the main square):  
	\begin{center}
		\includegraphics[scale=0.25]{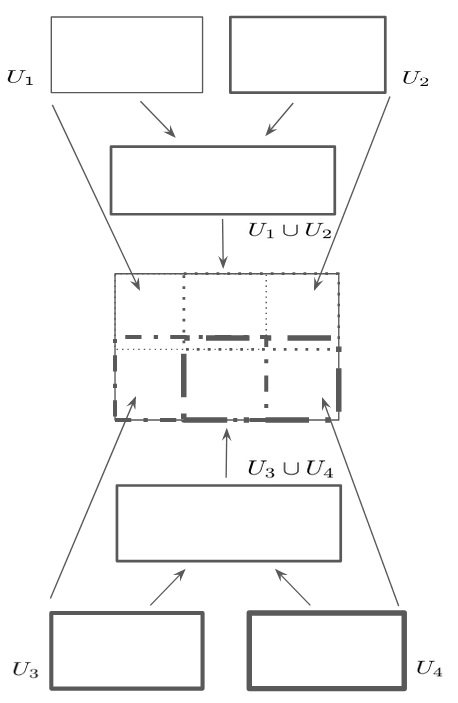}
	\end{center} 
	where we have left implicit the obvious intersections ($U_1 \cap U_2$, $U_3 \cap U_4$, $U_1 \cap U_3$, etc.). With the four regions $U_1, U_2, U_3,$ and $U_4$, to each of which there corresponds a particular sensor, we have the entire central region $U = U_1 \cup U_2 \cup U_3 \cup U_4$ `covered'. Part of what this means is that, were you to invite \textit{another} sensor to observe the happenings on some further portion of the space, in an important sense, this extra sensor would be superfluous---since, together, the four regions monitored by the four individual sensors already have the overall region `covered'. \par 
	For concreteness, suppose we have the following further selections of data from the data collected by each of these new (southwest and southeast) sensors, so that altogether, having performed the various compatibility checks (left implicit), the resulting system of ``points of view" on our site can be represented as follows:
	\begin{center}
		\includegraphics[angle=16, scale=0.3]{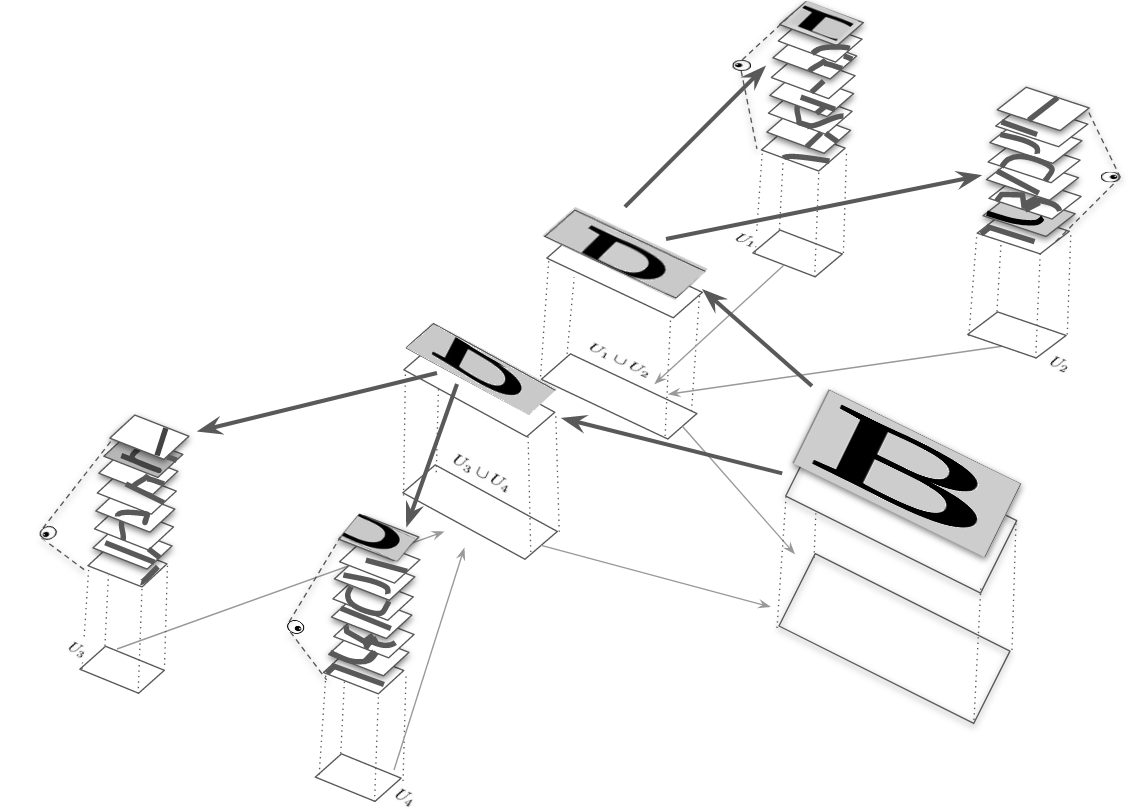}
	\end{center} 
	This system of mutually compatible local data assignments or ``measurements" of the happenings on the space---where the various data assignments are, piece by piece, constrained by one another, and thereby patched together to supply an assignment over the \textit{entire} space covered by the individual regions---is, in essence, what constitutes our \textit{sheaf}. The idea is that the data assignments are being ``tied together" in a natural way  
	\begin{center}
		\includegraphics[angle=16, scale=0.3]{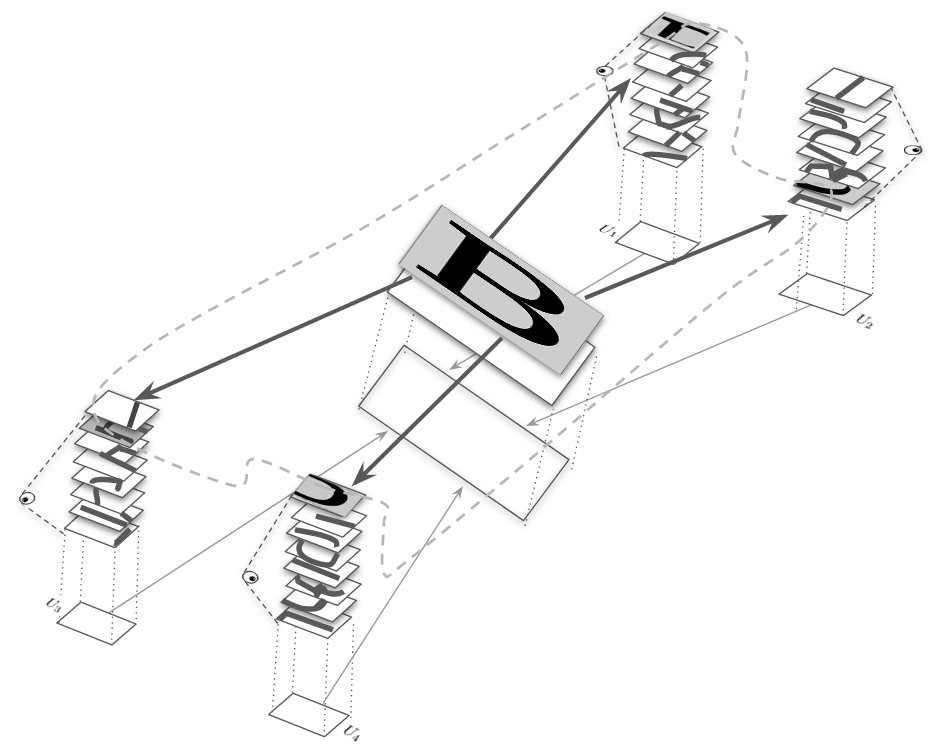}
	\end{center} 
	where this last picture is meant to serve as motivation or clarification regarding the agricultural terminology of `sheaf':  
	\begin{center}
		\includegraphics[scale=0.2]{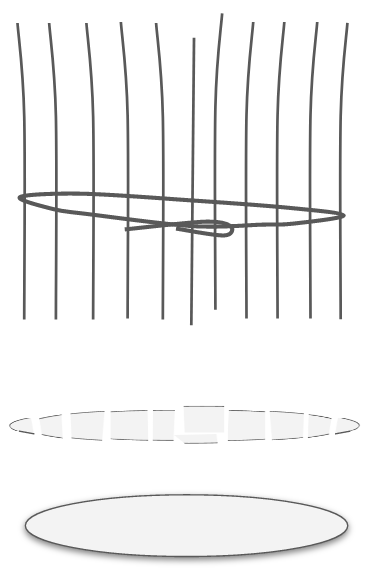}
	\end{center}
	Here one thinks of various `regions' as the parcels of an overall `space' covered by those pieces, the collection of which then serves as a `site' where certain happenings are held to take place, and the abstract sensors capturing local snapshots or measurements of all that is going on in each parcel are then regarded as being collected together into `stalks' of data, regarded as sitting over (or growing out of) the various parts of the ground space to which they are attached. A selection of a particular snapshot made from each of the individual stalks (collections of snapshots) then amounts to a cross section and the process of restriction (along intersecting regions) and collation (along unions of regions) of these sections then captures how the various stalks of data are ``bound together." \par  
	To sum up, then: the first characteristic feature of this construction is that some information is received or assigned \textit{locally}, so that the records or observations made by each of the individual sensors are understood as being ``about," or indexed to, the entirety of some limited region, so that whenever something holds or applies at a ``point" of that region, it will hold nearby as well. Next, since together the collection of regions monitored by the individual sensors may be seen as \textit{collectively covering} some overall region, we can check that the individual sensors that cover regions that have some overlap can ``communicate" their observations to one another, and a natural expectation is that, however different their records are on the non-overlapping region, there should be some sort of \textit{compatibility} or \textit{agreement} or \textit{mutual constraining} of the data recorded by the sensors over their shared, overlapping region; accordingly, we ask that each such pair of sensors covering overlapping regions ``check in" with one another. Finally, whenever such compatibility can be established, we expect that we can bind the information supplied by each sensor together, and regard them as patching together into a \textit{single sensor supplying data over the union} of the underlying (and partially overlapping) individual regions, in such a way that were we to ``restrict" that single sensor back down to one of the original regions, we would recover exactly the partial data reported by the original sensor assigned to that individual region. \par 
	While most of the more fascinating and conspicuous examples of such a construction come from pure and applied math, something very much like the sheaf construction appears to be operative in so many areas of ``everyday life." For instance, related to the toy example discussed above, even the way our binocular vision systems work appears to involve something like the collation of images into a single image along overlapping regions whenever there is agreement (from the input to each separate eye).\footnote{That visual information processing itself appears to fundamentally involve some sort of sheaf-like process appears even more acutely in other species, such as certain insects like the dragonfly, whose compound eyes contain up to 30,000 facets, each facet within the eye pointing in a slightly different direction and taking in light emanating from only one particular direction, resulting in a mosaic of partially overlapping images that are then integrated in the insect brain.} More generally, image and face recognition appears to operate, in a single brain (where clusters of neurons play the role of individual sensors), in something like the patchwork ``sum of parts" way described above. Moving beyond the individual, collective knowledge itself appears to operate in a fundamentally very similar way: a society's store of knowledge consists of a vast patchwork built up of partial records and data items referring to particular (possibly overlapping) regions, each of which data items can be (and often are!) checked for compatibility whenever they involve data that both refer to, or make claims about, the same underlying domain. \par 
	The very simple and naive presentation given to it above admittedly runs the risk of downplaying the power and scope of this construction; it would be difficult to overstate just how powerful the underlying idea of a sheaf is. An upshot of the previous illustration, though, is that while sheaves are often regarded as highly abstract and specialized constructions, whose power derives from their sophistication, the truth is that the underlying idea is so ubiquitous, so ``right before our eyes," that one might even be impressed that it was finally named explicitly so that substantial efforts could then be made to refine our ideas of it. In this context, one is reminded of the old joke about the fish, where an older fish swims up to two younger fish, and greets them ``morning, how's the water?" After swimming along for some time, one of the younger fishes turns to the other and says   
	\begin{quote}
		``What the hell is water?"
	\end{quote} 
	 In this same spirit, Grothendieck \index{Alexander Grothendieck} would highlight precisely this ``simplicity" of the fundamental idea behind sheaves (and, more generally, toposes):  
	\begin{quote}
		As even with the idea of sheaves (due to Leray), or that of schemes, as with all grand ideas that overthrow the established vision of things, the idea of the topos had everything one could hope to cause a disturbance, primarily through its ``self-evident" naturalness, through its simplicity (at the limit naive, simple-minded,  ``infantile") -- through that special quality which so often makes us cry out: ``Oh, that's all there is to it!", in a tone mixing betrayal with envy, that innuendo of the ``extravagant", the ``frivolous", that one reserves for all things that are unsettling by their unforeseen simplicity, causing us to recall, perhaps, the long buried days of our infancy.... (\cite{grothendieck_recoltes_1986}, Promenade 13) 
	\end{quote}
	\section{Outline of Contents}   
	The rest of the book is structured as follows. The first chapter is dedicated to exposition of the most important category-theoretic concepts, tools, and results needed for the subsequent development of sheaves. Category theory is indispensable to the presentation and understanding of the notions of sheaf theory. While in the last decade there have appeared a number of accessible introductions to category theory,\footnote{The general reader without much, or any, background in category theory is especially encouraged to have a look at the engaging and highly accessible \cite{spivak_category_2014}. Readers with more prior mathematical experience may find \cite{riehl_category_2016}, displaying the ubiquity of categorical constructions throughout many areas of mathematics, a compelling introduction. \cite{lawvere_sets_2003} is also highly recommended, especially for those readers content to be challenged to work many things out for themselves through thought-provoking exercises, often giving one the feeling of ``re-discovering" things for oneself.} feedback from readers of earlier drafts of this book convinced me that the best approach to an introduction to sheaves that aims to reach a much wider audience than usual would need to be as self-contained as possible. In this first chapter, all the necessary categorical fundamentals are accordingly motivated and developed. The emphasis here, as elsewhere in the book, is on explicit constructions and creative examples. For instance, the concept of an \textit{adjunction}, and key abstract properties of such things, is introduced and developed first through an extended example involving ``dilating" and ``eroding" an image, then again through the development of ``possibility" and ``necessity" modalities applied to both modeling the consideration of attributes of a person applied to them \textit{qua} the different ``hats" they wear in life, and then applied to graphs of traveling routes. While the reader already perfectly comfortable with category theory is free to skip this chapter or just skim through it, or refer back to later cited examples as needed, there are a few novel examples and (hopefully enlightening or at least mildly entertaining) philosophical discussions of important results such as the Yoneda lemma that may interest the expert as well. \par 
	Chapter 2 returns to presheaves (introduced in Chapter 1) to consider them in more depth. It discusses four main perspectives on \textit{presheaves}, develops a few notable examples of each of these, and develops some useful ways of understanding such constructions more generally. This is done both for its own sake and in order to build up to the following chapter dedicated to the initial development of the sheaf concept.\par 
	Chapter 3 introduces sheaves (specifically on topological spaces) and some key sheaf concepts and results---as always, through a diverse collection of examples. Throughout this chapter, some of the vital conceptual aspects of sheaves in the context of topological spaces are motivated, teased out, and illustrated through the various examples, and sometimes the same aspect is revisited from new perspectives as the level of complexity of the examples increases. \par 
	Chapter 4 is dedicated to a ``hands on" introduction to sheaf cohomology. The centerpiece of this chapter is an explicit construction, with worked-out computations, involving sheaves on complexes. There is also a brief look at \textit{cosheaves} and an interesting example relating sheaves and cosheaves. \par 
	Chapter 5 revisits and revises a number of earlier concepts, and develops sheaves from the more general perspective of \textit{toposes}. The important notions in topos theory (especially as this relates to sheaves) are motivated and developed through a variety of examples. We move through various layers of abstraction, from sheaves on a site (with a Grothendieck ``topology") or Grothendieck toposes to elementary toposes. The last few sections are devoted to illustrations, through concrete examples, of some slightly more advanced topos-theoretical notions and examples. The book concludes with an abridged presentation of some special topics, including a brief introduction to \textit{cohesive toposes}. There are many other directions the book could have taken at this point, and more advanced sheaf-theoretical topics that might have been considered, but in the interest of space, attention has been confined to this short final section on the special topic of cohesive toposes. \par 
	Throughout each chapter, I occasionally pause for a few pages to highlight, in a more ``philosophical" fashion (in what I call ``Philosophical Passes"), some of the important conceptual features to have emerged from the preceding technical developments. The overall aim of the ``Philosophical Pass" sections is to periodically step back from the technical details and examine the contributions of sheaf theory and category theory to the broader development of ideas. These sections may provide some needed rest for the reader, letting the brain productively ``switch modes" for some time, and giving one something to think about ``beyond the formal details." A lot of category theory, and the sheaf theory built on it, is deeply philosophical, in the sense that it speaks to, and further probes, questions and ideas that have fascinated human beings for millenia, going to the heart of some of the most lasting and knotty questions concerning, for instance, what an individual object is, the nature of the concept of `space', and the dialectics of continuity and discreteness. I hope it is not entirely due to my bias as someone who doubles as a professional philosopher that I believe that this sort of ``behind the scenes" reflection is an indispensable part not just of doing good mathematics but also of advancing our inquiry, as human beings, into some of these fundamental questions. 
	\chapter{Categorical Fundamentals for Sheaves}
	\section{Categorical Preliminaries}
	The language of category theory is indispensable to the presentation and understanding of the notions of sheaf theory. It is likely that any reader of this book has at least already \textit{heard} of categories, and may already be familiar with at least the basics of category theory. However, we will motivate and develop the necessary notions, and do so in a way that emphasizes connections with later constructions and perspectives that will emerge in our development of sheaves. The rest of this first section of the chapter supplies the definition of a category, then considers some notable examples of categories, and then presents an alternate perspective on categories. \par 
	Fundamentally, the specification of a category involves two main components: establishing some \textit{data} or givens, and then ensuring that this data conforms to two simple axioms or laws. To define, or verify that one has, a category, one should first make sure the right data is present. This first main step of establishing the data of a category really involves doing four things. First of all, it means identifying a collection of \textit{objects}. Especially when one is assembling a category out of already established mathematical materials, these objects will typically already go by another name, like vertices, sets, vector spaces, topological spaces, types, various algebras or structured sets, and so on. \par 
	Second, one must assemble or specify a collection of ``morphisms" or mappings, which is just some principled way of establishing connections between the objects of the first step. Again, when dealing with already established structures, these will usually already have a name, like arrows or edges, functions, linear transformations, continuous maps, terms, homomorphisms or structure-preserving maps, and so on. Many of the categories one meets in practice have sets with some structure attached to them for objects and (the corresponding) ``structure-preserving" mappings or connections between those sets for morphisms, so this is a good ``model" to keep in mind. \par 
	Third, and perhaps most importantly, one must specify an appropriate notion of \textit{composition} for these mappings, where for the moment this can be thought of in terms of specifying an operation that enables us to form a ``composite" map that goes directly from object $A$ to $C$ whenever there is a mapping from $A$ to $B$ juxtaposed with a mapping from $B$ to $C$. This composition operation in fact already determines the fourth requirement: that for each object, there is assigned a unique ``identity" (the ``do nothing") morphism that starts out from that object and return to itself. These four constituents---objects, morphisms, composites, and identities---supply us with the data of the category. \par 
	Next, one must show that the data given above conforms to two very ``natural" laws or axioms. First, if we have a morphism from one ``source" object to another ``target" object, then following that morphism with the identity morphism on the ``target" object should be the same thing as ``just" traveling along the original morphism; and the same should be true if we first travel along the identity morphism on the source object and then apply the morphism. In short, the identity morphisms cannot do anything to change other morphisms---this was why we referred to them above as the ``do nothing" morphisms. \par 
	Finally, a category must satisfy what is called the associative law, where this can be thought of as follows: if you have a string of morphisms from $A$ to $B$ and from $B$ to $C$ and from $C$ to $D$, then it should make no difference whether you choose to first go directly from $A$ to $C$ (using the composite map that we have by virtue of the third step in the data construction) followed by the map from $C$ to $D$, or if you go from $A$ to $B$ and then go directly from $B$ to $D$ (using the composite map). \par 
	An entity that has all the data specified above, data that in turn conforms to the two laws described in the preceding two paragraphs, is a category. The informal description given in the preceding paragraphs is given more formally in the following definition. 
\begin{definition}
	A \textit{category} $\textbf{C}$ \index{category!definition} consists of the following data:\footnote{Throughout this document, categories are generally designated with bold font. However, sometimes we may use script font instead, especially when dealing with things like pre-orders (discussed below), where each individual order is already a category. We will always make it clear what category we are working with, so this shouldn't be a problem.} 
	\begin{itemize}
		\item A collection $Ob(\textbf{C})$, whose elements are \textbf{objects}; 
		\item For every pair of objects $x, y \in Ob(\textbf{C})$, a collection $\text{Hom}_{\textbf{C}}(x,y)$ (or just $\textbf{C}(x,y)$) of \textbf{morphisms} from $x$ to $y$;\footnote{The term ``morphism" comes from \textit{homomorphism}, which is how one refers to a structure-preserving function in algebra, and which explains the notation ``Hom." Morphisms are also commonly referred to as ``arrows" or ``maps."} 
		\item To each object $x \in Ob(\textbf{C})$ is assigned a specified \textbf{identity morphism} on $x$, denoted $\text{id}_x$ $\in Hom_{\textbf{C}}(x,x)$;
		\item For every three objects $x,y,z \in Ob(\textbf{C})$, a function 
		\begin{equation*}
		\circ: \text{Hom}_{\textbf{C}}(y,z) \times \text{Hom}_{\textbf{C}}(x,y) \rightarrow \text{Hom}_{\textbf{C}}(x,z),
		\end{equation*}
		called the \textbf{composition formula}, which acts on elements to assign, to any morphism $f: x \rightarrow y$ and any $g: y \rightarrow z$, the composite morphism\footnote{One reads this right-to-left: first apply $f$, then run $g$ on the result.} $g \circ f: x \rightarrow z$: 
		\begin{equation*}
		\begin{split} 
		& \circ: \text{Hom}_{\textbf{C}}(y,z) \times \text{Hom}_{\textbf{C}}(x,y) \rightarrow \text{Hom}_{\textbf{C}}(x,z) \\
		& \circ( \hspace*{1em} g \hspace*{2em}, \hspace*{1.5em}f) \hspace*{0.5em}\mapsto (g \circ f)  
		\end{split} 
		\end{equation*}
	\end{itemize}
	This data gives us a category provided it further satisfies the following two axioms:  
	\begin{itemize}
		\item \textbf{Associativity} (of composition): if $x \xrightarrow{f} y \xrightarrow{g} z \xrightarrow{h} w$, then $h \circ (g \circ f) = (h \circ g) \circ f$. 
		\begin{center}
			\begin{tikzcd}
				x \arrow[r, "f"] \arrow[rr, bend left = 50, "{g \circ f}"] & y \arrow[r, "g"] \arrow[rr, bend right = 50, "{h \circ g}", swap] & z \arrow[r, "h"] & w
			\end{tikzcd}
		\end{center}
		\item \textbf{Identity}: if $f: x \rightarrow y$, then $f = f \circ id_x$ and $f = id_y \circ f$. 
	\end{itemize} 
\end{definition} 
\begin{example}
	The category $\textbf{Set}$ consisting of sets for objects and functions (with specified domain and codomain) for morphisms is in fact a category.\footnote{While this comment may not make sense to the reader right now, set theory can be thought of as ``zero-dimensional" category theory.}    
\end{example}
\begin{example}
	(\textit{Category of Pre-orders (Posets)}) \index{category!of pre-orders} \index{category!of posets}
	Recall that a relation between sets $X$ and $Y$ is just a subset $R \subseteq  X \times Y$, and that a \textit{binary relation} on $X$, is a subset $R \subseteq  X \times X$. It is customary to use infix notation for binary relations, so that, for instance, one writes $a \leq b$ for $(a,b) \in R$. We define a \textit{pre-order} \index{pre-order!defined} as a set with a binary relation (call it `$\leq$') that further satisfies the properties of being \textit{reflexive} and \textit{transitive}. In other words, it is a pair $(X, \leq_X)$ where we have 
	\begin{itemize}
		\item $x \leq x$ for all $x \in X$ (reflexivity); and 
		\item if $x \leq y$ and $y \leq z$, then $x \leq z$ (transitivity).
	\end{itemize} 
Then a \textit{poset} \index{poset!defined} is a pre-order that is additionally \textit{anti-symmetric}, where this means that $x \leq y$ and $y \leq x$ implies that $x = y$. \par 
It is often useful to represent a given poset (or pre-order) with a diagram. For instance, supposing we have an order-structure on $P = \{a,b,c,d\}$ given by $a \leq c, b \leq c, b \leq d$, together with the obvious identity (reflexivity) $x \leq x$ for all $x \in P$. The data of this poset may be displayed by the diagram:    
\begin{center} 	
	\begin{tikzpicture}[yscale=0.75, xscale=0.5]
	\node (c) at (2.4,0) {$c$};
	\node (d) at (6.4,0) {$d$};
	\node (a) at (2.4,-2) {$a$};
	\node (b) at (6.4,-2) {$b$};
	
	\draw[->] (a) -- (c);
	\draw[->] (b) -- (c);
	\draw[->] (b) -- (d);
	\end{tikzpicture} \end{center} 
Pre-orders (posets) can themselves be related to one another, and the right notion here is one of a monotone (or order-preserving) map. 
\begin{definition}
	A \textit{monotone} (\textit{order-preserving}) map \index{monotone map!defined} between pre-orders (or posets) $(X, \leq_X)$ and $(Y, \leq_Y)$ is a function $f: X \rightarrow Y$ satisfying that for all elements $a, b \in X$, 
	\begin{equation*}
	\text{if } a \leq_X b, \text{ then } f(a) \leq_Y f(b).  
	\end{equation*} 
\end{definition}
\noindent 
\textbf{Pre} is the category having pre-orders for objects and order-preserving functions for morphisms. \textbf{Pos} is the category having posets for objects and order-preserving functions for morphisms.\footnote{As one can see from the examples given thus far, it is common for a category to be named after its objects. However, this widespread practice is not really in accord with the ``philosophy" of category theory, which gives primacy to the morphisms (or at least demands that objects be considered together with their morphisms). We will explore this point further in section \ref{sec: No Objects}.} Each identity arrow will just be the corresponding identity function, regarded as a monotone map. It is easy to verify that for two monotone maps $X \xrightarrow{f} Y$ and $Y \xrightarrow{g} Z$ between orders, the function composition $g \circ f$ is also monotone. \par 
If we further add the property that for all $x, x' \in X$, either $x \leq x'$ or $x' \leq x$, i.e., any two objects are \textit{comparable}, then we get what are called \textit{linear orders}.\index{linear order} In particular for $n \in N$ a natural number, we can consider the linear order $[n] = (\{0,1,\dots, n\}, \leqslant)$, where every finite linear order may be represented pictorially  
	\begin{center} 
		\begin{tikzpicture}\small
		\foreach \x in {0,...,3}
		\fill (0,0) circle (1.5pt);
		\fill (1,0) circle (1.5pt); 
		\fill (2,0) circle (1.5pt); 
		\fill (3,0) circle (1.5pt); 
		\node[draw=none] (ellipsis1) at (4.15,0) {$\cdots$};
		\fill (5.3,0) circle (1.5pt); 
		
		\foreach \x in {0,...,3}
		\node[above=2pt] at (\x,0) {\x};
		\node[above=2pt] at (5.3,0) {n}; 
		
		\foreach \x in {0,...,3} 
		\draw [-to] (\x + 0.15, 0) -- (\x + 0.85, 0);  
		\draw [-to] (4.35, 0) -- (5.15, 0); 
		\end{tikzpicture} 
	\end{center} 
	Together with morphisms Hom($[m]$, $[n]$) defined as all the functions $f: \{0,1,\dots,m\} \rightarrow \{0,1,\dots,n\}$ such that, for every pair of elements $i,j \in \{0,1,\dots, m\}$, if $i \leqslant j$, then $f(i) \leqslant f(j)$, i.e., monotone functions, this produces another category: \textbf{FLin}, the category of finite linear orders. \par    
	We also have \textit{cyclic orders}, \index{cyclic order} defined not as a binary relation, but as a ternary relation $[a,b,c]$ (read ``after $a$, one arrives at $b$ before $c$"). More formally, a cyclic order on a set is a ternary relation that satisfies: 
	\begin{enumerate}
		\item cyclicity: if $[a,b,c]$, then $[b,c,a]$; 
		\item asymmetry: if $[a,b,c]$, then not $[c,b,a]$; 
		\item transitivity: if $[a,b,c]$ and $[a,c,d]$, then $[a,b,d]$; 
		\item totality: if $a,b,$ and $c$ are distinct, then we have either $[a,b,c]$ or $[c,b,a]$. 
	\end{enumerate}
	You can think of a cyclic order on a set as an arrangement of the objects of that set around a circle, so that a cyclic order on a set with $n$ elements can be pictured as an (evenly spaced) arrangement of the objects of the set on an $n$-hour clock face. 
	\begin{center}
		\begin{tikzpicture}[scale = 0.6]
		\small 
		\foreach \ang\lab\anch in {90/1/north, 45/2/{north east}, 0/3/east, 270/i/south, 180/{n-1}/west, 135/n/{north west}}{
			\draw[fill=black] ($(0,0)+(\ang:3)$) circle (.08);
			\node[anchor=\anch] at ($(0,0)+(\ang:2.8)$) {$\lab$};
		}
		
		\foreach \ang\lab in {90/1,45/2,180/{n-1},135/n}{
			\draw[->,shorten <=7pt, shorten >=7pt] ($(0,0)+(\ang:3)$) arc (\ang:\ang-45:3);
			\node at ($(0,0)+(\ang-22.5:3.5)$) {};
		}
		
		\draw[->,shorten <=7pt] ($(0,0)+(0:3)$) arc (360:325:3);
		\draw[->,shorten >=7pt] ($(0,0)+(305:3)$) arc (305:270:3);
		\draw[->,shorten <=7pt] ($(0,0)+(270:3)$) arc (270:235:3);
		\draw[->,shorten >=7pt] ($(0,0)+(215:3)$) arc (215:180:3);
		\node at ($(0,0)+(0-20:3.5)$) {};
		\node at ($(0,0)+(315-25:3.5)$) {};
		\node at ($(0,0)+(270-20:3.5)$) {};
		\node at ($(0,0)+(225-25:3.5)$) {};
		
		\foreach \ang in {310,315,320,220,225,230}{
			\draw[fill=black] ($(0,0)+(\ang:3)$) circle (.02);
		}
		
		\end{tikzpicture}
	\end{center} 
	Such finite cyclically ordered sets are sometimes designated $\Lambda_n$, for each natural number $n$. If we take as objects, for each $n \in \mathbb{N}$, the object $\Lambda_n$, and for morphisms $\text{Hom}_{\bf{\Lambda}}(\Lambda_m, \Lambda_n)$ monotone functions, i.e., functions from $\{0,1,\dots, m\}$ to $\{0,1,\dots, n\}$ such that whenever $[f(a), f(b), f(c)]$, we have $[a,b,c]$ for all $a,b,c \in \{0,1,\dots, m\}$, then we get the \textit{cyclic category} \index{category!cyclic} $\bf{\Lambda}$.\footnote{Another usual way of defining the morphisms of this category is in terms of the increasing functions $f: \mathbb{Z} \rightarrow \mathbb{Z}$ satisfying $f(i + m + 1) = f(i) + n + 1$.} \par 
	Orders, especially pre-orders and posets, are very important in category theory, and we will see a lot more of them throughout the book. 
\end{example}
\begin{example}
	A graph \index{graph} is typically represented by a bunch of dots or vertices together with certain edges or arrows linking a pair of vertices and defining what is called a relationship of \textit{incidence} between the vertices and edges. More formally, a (simple) \textit{graph} $G$ consists of a set $V$ of \textit{vertices}, together with a collection of two-element subsets $\{x,y\}$ of $V$ (or sometimes just represented by a set $E$ that consists of the ``names" of such pairings, via stipulating an additional mapping that interprets edges as pairs of vertices), called the \textit{edges}. A graph morphism $G \rightarrow H$ is then a function $f: V \rightarrow V'$ on the vertices such that $\{f(x), f(y)\}$ is an edge of $H$ whenever $\{x,y\}$ is an edge of $G$.\par
	As the pairs of vertices above are defined to be \textit{unordered}, the resulting graphs are undirected. If we are assuming that the map interpreting edges as unordered pairs of vertices does so in a one-to-one way, we are requiring that the graph be ``simple" in the sense of having at most one edge between two vertices. In this case, we will have constructed the category of undirected (simple) graphs, $\textbf{UGrph}$, or more commonly $\textbf{SmpGrph}$. This is often what the graph theorist means, by default, by `graph'. Note that if we allowed instead, for each unordered pair of distinct vertices, an entire set of edges between these, we would generalize this to \textit{multigraphs}.\par   
	We can further define \textit{directed graphs} (which often go under the name of \textit{quivers} by category theorists).\label{example: graph} A (directed) \textit{graph} $G = (V, A, s, t)$ consists of a set $V$ of vertices, a set $A$ of directed edges, or arrows (or arcs), and two functions 
		\begin{center} 
		\begin{tikzcd}
			A \arrow[r, shift left = 1ex, "s"] \arrow[r, "t", swap] & V 
		\end{tikzcd}
	\end{center} 
that act to pick out the \textit{source} and \textit{target} of an arc. \par 
Then if $G = (V, A, s,t)$ and $G' = (V', A', s', t')$ are two graphs, a \textit{graph homomorphism} $f: G \rightarrow G'$ requires that (the ordered pair) $(f(x), f(y))$ is an arc of $G'$ whenever $(x,y)$ is an arc of $G$. More explicitly, a graph morphism $g: G \rightarrow G'$ is a pair of morphisms $g_0: V \rightarrow V'$ and $g_1: A \rightarrow A'$ such that sources and targets are preserved, i.e., 
\begin{equation*}
s' \circ g_1 = g_0 \circ s \text{ and } t' \circ g_1 = g_0 \circ t.
\end{equation*}
In general, there may exist several parallel arrows, i.e., with the same source and same target, in which case we are dealing with directed \textit{multigraphs}. If we allow closed arrows, or loops, i.e., arrows whose source are target are identical, then we are dealing with \textit{looped} (or \textit{reflexive}) graphs. There is a lot more to say about distinctions between different graphs, the distinct categories for each, and their categorical features of interest; but we will postpone this until later chapters. \par 
In the case of the above directed graphs, this is the category $\textbf{dGrph}$ (or just $\textbf{Grph}$), which has directed graphs as objects, and (directed) graph homomorphisms (i.e., source and target preserving morphisms) as morphisms.\index{category!of graphs}    
\end{example}
\begin{example}
	The category \textbf{Mon} (\textbf{Group}) of monoids (groups) \index{category!of monoids} \index{category!of groups} has monoids (groups) for objects and monoid (group) homomorphisms for morphisms. (This example, together with the necessary definitions, will be discussed in much more detail in a moment.)
\end{example}
\begin{example}
	The category $\textbf{Vect}$ \index{category!of vector spaces} is the category of $k$-vector spaces (for a given field $k$, dropping the $k$ when this is understood), which has vector spaces for objects and linear transformations for morphisms. Restricting attention to just finite-dimensional vector spaces yields the category $\textbf{FinVect}$, which is where most of linear algebra takes place.   
\end{example}
The previous four examples are just a few of the many examples of \textit{categories of structures}, or sets with some structure on them. When, in 1945, Eilenberg\index{Eilenberg} and MacLane\index{MacLane} first defined categories and the related notions (introduced below) allowing categories to be compared, they stressed how it provided ``opportunities for the comparison of constructions[...]in different branches of mathematics." But with Grothendieck's\index{Alexander Grothendieck} \textit{Tohoku} paper a decade later, it became more and more evident that category theory was not just a convenient way of comparing different mathematical structures, but was itself a significant mathematical structure of its own intrinsic interest. One way of starting to appreciate this is to realize that we do not just have categories consisting \textit{of} mathematical objects/structures, but equally important are those categories that allow us to view categories themselves \textit{as} mathematical objects/structures. The following important examples supply simple examples of this perspective of \textit{categories as structures} (the first two of which reveal crucial features of categories in general and are accordingly often said to supply us with a means of doing ``category theory in the miniature").  
\begin{example}
	(\textit{Each order is already a category}) Let $(X, \leq_X)$ be a given pre-order\index{pre-order!as category} (or, less generally, a poset).\index{poset!as category} It is easy to check that we can form the category $\mathcal{X}$ by taking 
	\begin{itemize}
		\item the elements of $X$ as the objects of $\mathcal{X}$; and  
		\item for elements $a, b \in X$, there exists a morphism in $\mathcal{X}$ from $a$ to $b$ exactly when $a \leq b$ (and there is \textit{at most} one such arrow, so this morphism will necessarily be unique). 
	\end{itemize}
Notice how transitivity of the relation $\leq$ automatically gives us the required composition morphisms, while reflexivity of $\leq$ just translates to the existence of identity morphisms. Thus, we can regard any given poset (pre-order) $(X, \leq_X)$ as a category $\mathcal{X}$ in its own right.  
\end{example}
\begin{example}
 (\textit{Each monoid is already a category}) A \textit{monoid} \index{monoid!defined} $\mathcal{M} = (M, \cdot, e)$ is a set $M$ equipped with 
	\begin{itemize} 
		\item an associative binary multiplication operation $\cdot: M \times M \rightarrow M$, i.e., $\cdot$ is a function from $M \times M$ to $M$ (a \textit{binary operation on} $M$) assigning to each pair $(x,y) \in M \times M$ an element $x \cdot y$ of $M$, where this operation is moreover associative in the sense that 
		\begin{equation*}
		x \cdot (y \cdot z) = (x \cdot y) \cdot z
		\end{equation*} 
		for all $x, y, z \in M$; and  
		\item a two-sided ``identity" element $e \in M$, where this satisfies 
		\begin{equation*}
		e \cdot x = x = x \cdot e
		\end{equation*}  
		for all $x \in M$.  
	\end{itemize} 
	Comparing this definition to that of a category, it is straightforward to see how any monoid $\mathcal{M}$ can be regarded as a category of its own.\index{monoid!as category} Specifically, it is a category with just one object. Explicitly, a monoid $(M, e, \cdot)$ can be considered as a category $\mathcal{M}$ with one object and with hom-set equal to $M$ itself, where the identity morphism comes from the monoid identity $e$ and the composition formula from the monoid multiplication $\cdot:M \times M \rightarrow M$. In other words, the category $\mathcal{M}$ is defined as consisting of 
	\begin{itemize}
		\item objects: the single object $M$ itself; 
		\item morphisms: the members of $M$ (where each monoid element represents a distinct endomorphism, i.e., map from $M$ to $M$, on the single object). 
	\end{itemize}
Then, the identity $\text{Id}_{\mathcal{M}}$ is given by $e$ and composition of arrows $x, y \in M$ is just given by the monoid multiplication 
\begin{equation*}
x \circ y = x \cdot y. 
\end{equation*}
Conversely, notice that if $\textbf{C}$ is a category with only one object $a$ and $M$ is its collection of morphisms, then $(M, \circ, \text{Id}_a)$ will be a monoid. \par 
Finally, an element $m \in M$ of a monoid is said to have an \textit{inverse} provided there exists an $m' \in M$ such that $m \cdot m' = e$ and $m' \cdot m = e$. Recall that a \textit{group} is just a monoid for which every element $m \in M$ has an inverse. Similar to the above, then, any group itself gives rise to a category in which there is just one object, but where every morphism (given by the group elements) is now an isomorphism. 
\end{example}
	The previous two examples are not just examples of \textit{any old} categories, but in an important sense, categories in general may be regarded as a sort of fusion of preorders on the one hand and monoids on the other. Over and above the fact that each monoid and each preorder is itself already a category, these two examples are ``special" in that categories more generally are exceptionally ``monoid-like" and ``preorder-like." We saw that every monoid is a single-object category. Seen from the other side, categories in general may be regarded as the ``many-object" version of monoids. We saw that every preorder is a single-arrowed category, as between any two objects there is at most one arrow. Seen from the other side, categories may be regarded as the ``many-arrowed" version of preorders. Monoids furnish us with not just a study of composition ``in the miniature" (by collapsing down to a single object), but in a sense the associative binary operation and neutral or identity element that comprise the data of a monoid seem to provide a prototype for the general associativity and identity \textit{axioms} of a category. Preorders, for their part, furnish us not just with a study of comparison of objects via morphisms ``in the miniature" (by collapsing down to at most one morphism from any object to another), but in a sense the reflexivity and transitivity of the order seems to provide the model for the key \textit{data} specifying a category, i.e., the assignment of an identity arrow to each object (via reflexivity) and the composition formula (via transitivity).\par 
Before continuing with examples, there is another important (if somewhat ``philosophical") way in which monoids can shed light on categories. This has the added benefit of introducing the interesting notion of \textit{oidification} and an alternative (philosophically appealing, if somewhat less useful) definition of categories.
\subsection{Aside on ``No Objects"}
\label{sec: No Objects}
The following presents an alternative, single-sorted or ``no objects" version of the definition of a category.
\begin{definition}
	 (\textit{Category Definition Again (``No Objects" version)}) A \textit{category} (single-sorted) \index{category!definition} is a collection $C$, the elements or `individuals' of which are called \textit{morphisms}, together with two endofunctions $s, t: C \rightarrow C$ (think ``source" and ``target") on $C$ and a partial function $\circ: C \times C \rightarrow C$, where these satisfy the following axioms: 
	\begin{enumerate}
		\item $x \circ y$ is defined iff $s(x) = t(y)$
		\item $s(s(x)) = s(x) = t(s(x))$ and $t(t(x)) = t(x) = s(t(x))$ (so $s$ and $t$ are \textit{idempotent} endofunctions on $C$ with the same image)
		\item if $x \circ y$ is defined, then $s(x \circ y) = s(y)$ and $t(x \circ y) = t(x)$
		\item $(x \circ y) \circ z = x \circ (y \circ z)$ (whenever either is defined)
		\item $x \circ s(x) = x$ and $t(x) \circ x = x$. 
	\end{enumerate}
\end{definition}
Notice how the elements of the shared image of $s$ and $t$, i.e., the $x$ such that $s(x) = x$ (equivalently, $t(x) = x$), are the \textit{identities} (or \textit{objects}). \par 
Probably the ``punch" of this definition is lost on a reader seeing it for the first time. Behind this presentation is the idea that each object in the usual definition of a category can in fact be \textit{identified} with its identity morphism, allowing us to realize an ``arrows-only" (or ``object-free") definition of a category. It is in the context of such an ``arrows-only" version that we can even more easily see how monoids are just one-object categories (so that, ultimately, categories in general are just many-object monoids). From a given monoid, we obtain a category by defining $s(x) = t(x) = e$, where $e$ is the monoid constant (``identity") element. Going the other way, given a (non-empty) category satisfying any of 
\begin{itemize}
	\item $s(x) = s(y)$,
	\item $s(x) = t(y)$, or 
	\item $t(x) = t(y)$,\footnote{In other words, $s$ is a constant function (and thus, so is $t$, and they are in fact equal).}
\end{itemize} 
we can define $e$ as the (unique) identity morphism, and thus obtain a monoid. \par 
In this way, single-sorted categories appear, in an especially clear way, via what is sometimes called the ``\textit{oid}ification" \index{oidification} of mon\textit{oid}s, where this describes a twofold process whereby 
\begin{enumerate}
	\item a concept is realized as equivalent to a certain category with a \textit{single} object; and then 
	\item the concept is generalized (``oidified") by moving to categories of that type that now have more than one object.  
\end{enumerate} 
We will see more examples of this process later on. \par 
For now, let us remark briefly on the significance of this ``object-free" perspective. Consider how, in the context of graphs and graph theory, the novice will likely see arcs (arrows) as secondary to vertices (objects), for the arcs are frequently construed as just pairs of vertices. It also seems a valid observation that ``psychologically" it is somehow more natural for many of us to begin with objects (as the irreducible ``simples") and then move on to \textit{relations} between those objects. But in more general treatments of graphs, dealing with directed multigraphs or quivers for instance, one begins to appreciate that this proclivity really gets things backwards: in fact, in more general settings, arcs are more naturally seen as primary and vertices can be seen as ``degenerate" sorts of arcs, or as equivalence classes of arcs under the relations ``has the same source (target) as." \par 
In a similar fashion, one might argue that our default ``object-oriented" mindsets can get things backwards, in terms of what is really fundamental conceptually. It is often said in category theory that ``what matters are the arrows/relations, not objects," for by the above line of reasoning, it is the algebra of morphisms that really determines the category. This is a very powerful idea, one that seems to permeate many aspects of category theory, and even resurfaces in a particularly poignant way with one of the key results in category theory (Yoneda lemma and embedding). The ``object-free" definition of a category (as above) is not standard, perhaps because it seems to complicate the presentation of many classical examples of categories, whose presentation is comparatively more straightforward using the classical definition of a category. However, the ``object-free" approach is arguably even more fundamental conceptually, and well-attuned to the core ``philosophy" of much of the categorical approach (which insists, in many contexts, that what matters is how objects and structures interact or relate), so it is worthwhile to at least be familiar with the existence of such an alternative definition. 
\subsection{A Few More Examples}
\begin{example}
 	Suppose we are given $V$ a vector space. Then we can define a category $\textbf{V}$ as follows: 
 	\begin{itemize}
 		\item for objects: \textbf{V} has only one object, called $*$; 
 		\item for morphisms (arrows from $*$ to $*$): the vectors $v$ in $V$
 		\item for the identity arrow for $*$: the zero vector; and 
 		\item for composition of vectors $v$ and $v'$: their sum.  
 	\end{itemize} 
\end{example}
We turn now to a very important example, one that starts to make better sense of the idea that category theory, in being visualized as it is by arrows between dots, might lead one to want to regard category theory as some sort of graph theory, but with ``extra data," where this involves some ``extra structure" (specifically, the \textit{composition} of arrows). In our definition of a graph earlier, observe that there were no other conditions placed on arrows/edges and vertices, other than that involving the source and target functions, picking out the source vertex and the target vertex of a given arc $a$; in particular, there was no requirement regarding the composition of arrows/edges. Thus, it is not the case that a category \textit{is} a graph, for a directed graph in general has no notion of composition of edges/arrows (and not even a notion of identity arrows). However, any category (well, any ``small" category, on which more below) does have an underlying graph. While the converse does not hold, it is an important fact that every directed graph can be \textit{made} into a certain category, via a special construction, discussed in the following.     
\begin{example}
	Given a directed graph $G$, we first describe the notion of a \textit{path} in $G$, as any sequence of successive arrows where the target of one arrow is the source of the other. More explicitly, for each $n \in \mathbb{N}$, we define a \textit{path} \index{path!defined} through $G$ of length $n$ as a list of $n$ edges, 
	\begin{center}
		\begin{tikzcd}
			i(0) \arrow[r, "e(1)"] & i(1) \arrow[r, "e(2)"] & i(2) \arrow[r] & \cdots & \arrow[r, "e(n)"] & i(n)
		\end{tikzcd}
	\end{center}
	where the target of each edge is the source of the next one. A path of length 1 would then be a single edge, while a path of length 0 would be a vertex. We can create a category $\textbf{Pth}(G)$, the category of paths through $G$, \index{category!of paths} with objects the nodes of $G$ and for morphisms from objects $x$ to $y$ all the paths through $G$ from $x$ to $y$. Given two paths, 
	\begin{center}
	\begin{tikzcd}
		i(0) \arrow[r] & i(1) \arrow[r] & i(2) \arrow[r] & \cdots \arrow[r] & i(n)
	\end{tikzcd}
\end{center}
and 
	\begin{center}
	\begin{tikzcd}
		j(0) \arrow[r] & j(1) \arrow[r] & j(2) \arrow[r] & \cdots \arrow[r] & j(m), 
	\end{tikzcd}
	\end{center} 
with the end node of the first equal to the start node of the second, i.e., $i(n) = j(0)$, we form the \textit{composite path} by concatenating or sticking the two paths together along this identical node, i.e., 
	\begin{center}
	\begin{tikzcd}
		i(0) \arrow[r] & i(1) \arrow[r] & \cdots \arrow[r] & i(n) = j(0) \arrow[r] & j(1) \arrow[r] & \cdots \arrow[r] & j(m)
	\end{tikzcd}
\end{center}
resulting in a new path from $i(0)$ to $j(m)$. Then, concatenating paths end to end is associative, making composition in $\textbf{Pth}(G)$ associative. As for ensuring that each object (vertex of $G$) has an identity arrow in $\textbf{Pth}(G)$, we can observe that each vertex has an associated ``length 0" path, and sticking such a path at the end of another path does nothing to change that other path. Thus, we can just take the paths of length 0 as our identity arrows, i.e., the identity of an object $x$ is given by the path of length 0 from $x$ to $x$. \par 
	 Moreover, given a graph homomorphism $f: G \rightarrow G'$, every path in $G$ will be sent under $f$ to a path in $G'$. We will have more to say about this category, and the construction that generates it, in a subsequent section. 
\end{example}
There are many more categories that we might mention, and that are important to mathematicians. However, we will instead move forward and let the categories that will be of particular use to us emerge organically throughout the book. For now, here is a curtailed and rather arbitrary list of just a few more categories of general interest. 
\begin{itemize}
	\item $\textbf{Top}$: the category that has topological spaces for objects and continuous functions for morphisms.\index{category!of topological spaces} 
	\item $\textbf{Measure}$:\index{category! of measure spaces} the category that has measure spaces as objects and (on one definition) for morphisms appropriate equivalence classes of measurable functions. 
	\item $\textbf{Cat}$: the category of...categories! This has categories for objects and functors (defined and discussed shortly) for morphisms.\index{category! of categories}  
\end{itemize}
\subsection{Some New Categories From Old}
	Finally, there are many important things one can do \textit{to} categories, to generate new categories from old ones. Attention is confined, for the moment, to those that will be most important for our purposes.  
	\begin{definition}
		Let $\textbf{C}$ be a category. The \textit{dual} (or \textit{opposite}) category $\textbf{C}^{op}$ \index{category!opposite} is then defined as follows: 
		\begin{itemize}
			\item objects: same as the objects of $\textbf{C}$; 
			\item morphisms: given objects $A, B$ the morphisms from $A$ to $B$ in $\textbf{C}^{op}$ are exactly the morphisms from $B$ to $A$ in \textbf{C}. (In other words, just reverse the direction of all the arrows in \textbf{C}.) 
		\end{itemize} 
	Identities for $\textbf{C}^{op}$ are defined as before, and composites are formed by reversing arrows as one would expect, yielding a category. In more detail, for each $\textbf{C}$-arrow $f: A \rightarrow B$, introduce an arrow $f^{op}: B \rightarrow A$ in $\textbf{C}^{op}$, so that ultimately, these give all and only the arrows in $\textbf{C}^{op}$. Then the composite $f^{op} \circ g^{op}$ will be defined precisely when $g \circ f$ is defined in $\textbf{C}$, where for 
\begin{center} 
	\begin{tikzcd}[column sep =large]
		A \arrow[bend left = 30]{r}[name=U]{f} & B \arrow{l}[name=L]{f^{op}} \arrow[bend left = 30]{r}[name=Q]{g} & C \arrow{l}[name=R]{g^{op}},
	\end{tikzcd}
\end{center}   
we have that $f^{op} \circ g^{op} = (g \circ f)^{op}$.
	\end{definition} \noindent 
In short, and in slogan-form, 
\begin{quote}
	\textit{Given a category, just reverse all its morphisms, and you'll get another category (its} dual\textit{)!}
\end{quote}
	With this seemingly innocuous construction, every result in category theory will have a corresponding dual, essentially got ``for free" by simply formally ``reversing all arrows." In other words, when a statement is true in a category $\textbf{C}$, then its dual will be true in the dual category $\textbf{C}^{op}$. Such duality not only can clarify and simplify relationships that are often hidden in applications or particular contexts, but it also ``reduces by half" the proof of certain statements (since the other, dual statement will ``follow by duality")---or, to see things another way, it ``multiplies by two" the number of results, as each theorem will have its corresponding dual. 
Finally, for any $\textbf{C}$, note that we will have that $(\textbf{C}^{op})^{op} = \textbf{C}$. \par 
Next, we consider how, given a category $\textbf{C}$, we can form a new category by taking as our objects all the \textit{arrows} of $\textbf{C}$. 
	\begin{definition}
		For a category \textbf{C}, we define the \textit{arrow category} \index{category!arrow} of \textbf{C}, denoted $\textbf{C}^{\rightarrow}$, as having for 
		\begin{itemize}
			\item objects: morphisms $A \rightarrow B$ of $\textbf{C}$; and for 
			\item morphisms: from the object $A \rightarrow B$ to the object $A' \rightarrow B'$, a morphism is a couple $(A \rightarrow A', B \rightarrow B')$ of morphisms of $\textbf{C}$ making the diagram 
			\begin{center}
				\begin{tikzcd}
					A \arrow[d] \arrow[r] & B \arrow[d] \\ 
					A' \arrow[r] & B'
				\end{tikzcd} 
			\end{center}
			\par \noindent 
			commute. 
		\end{itemize}
	Composition of arrows is then carried out in the obvious way. 
	\end{definition} \noindent 
With the arrow category, we are seeing \textit{all} the arrows of the old category as our objects in the new category. The next construction instead looks at just \textit{some} of the old arrows, where we restrict attention to arrows that have fixed domain or codomain.
		\begin{definition}
			Given a category $\textbf{C}$, and an object $A$ of $\textbf{C}$, we can form the two categories called the \textit{slice} and \textit{co-slice} categories, \index{category!slice} \index{category!co-slice} respectively denoted 
			\begin{equation*}
			(\textbf{C} \downarrow A) \hspace*{3em} (A \downarrow \textbf{C}),
			\end{equation*}
			also called the category of
			\begin{equation*} 
			\text{objects over } A \hspace*{3em} \text{objects under } A,
			\end{equation*}  
			respectively.\footnote{These categories are also particular cases of a more general construction, known as \textit{comma categories}.\index{category!comma} It is not uncommon to see the slice category of objects over $A \in \text{Ob}(\textbf{C})$ referred to as $\textbf{C} / A$, and the co-slice category of objects under $A$ referred to as $A / \textbf{C}$.} The objects of the new category are given by 
			\begin{equation*}
				\text{arrows to } A \hspace*{3em} \text{arrows from } A.
			\end{equation*}  
			In other words, objects of the slice category are given by all pairs $(B, f)$, where $B$ is an object of $\textbf{C}$ and $f: A \rightarrow B$ an arrow of $\textbf{C}$, and of the co-slice category by all pairs $(B, f)$ such that $f: B \rightarrow A$ is an arrow of $\textbf{C}$. \par 
			Morphisms in the new category are given by $h: (B, f) \rightarrow (B', f')$ where this is an arrow $h: B \rightarrow B'$ of $\textbf{C}$ for which the respective triangles 
				\begin{center}
					\begin{tikzcd}
						B \arrow[dr, "f", swap] \arrow[rr, "h"] & & B' \arrow[dl, "g"] & &  A \arrow[dl, "f", swap] \arrow[dr, "g"] \\ 
						& A & & B \arrow[rr, "h", swap] & & B'
					\end{tikzcd}
				\end{center} 
				\par \noindent 
				\textit{commute} in the sense that, for instance, for the triangle on the left, $g \circ h = f$. \par 
			Composition in $(\textbf{C} \downarrow A)$ and $(A \downarrow \textbf{C})$ is then given by composition in $\textbf{C}$ of the base arrows $h$ of such triangles.  
		\end{definition}  
	Categories of this type play an important role in advancing some of the general theory, in addition to being of some intrinsic interest. For now, the slice category of \textit{objects over} $A$ might be thought of as giving something like a view of the category \textit{seen within the context of} $A$ (and the corresponding dual statement for the category of objects under $A$). \par  
Finally, we define the following notion of a \textit{subcategory}. 
\begin{definition}
	A \textit{subcategory} \index{category!subcategory} $\textbf{D}$ of a category $\textbf{C}$ is got by restricting to a subcollection of the collection of objects of $\textbf{C}$ (i.e., every $\textbf{D}$-object is a $\textbf{C}$-object), and a subcollection of the collection of morphisms of $\textbf{C}$ (i.e., if $A$ and $B$ are any two $\textbf{D}$-objects, then all the $\textbf{D}$-arrows $A \rightarrow B$ are present in $\textbf{C}$), where we further require that
	\begin{itemize}
		\item if the morphism $f: A \rightarrow B$ is in $\textbf{D}$, then $A$ and $B$ are in $\textbf{D}$ as well. 
		\item if $A$ is in $\textbf{D}$, then so too is the identity morphism $\text{Id}_A$. 
		\item if $f: A \rightarrow B$ and $g: B \rightarrow C$ are in $\textbf{D}$, then so too is the composite $g \circ f: A \rightarrow C$. 
	\end{itemize}
\end{definition} \noindent 
Moreover, we can also define the following: 
\begin{definition}
	Let $\textbf{D}$ be a subcategory of $\textbf{C}$. Then we say that $\textbf{D}$ is a \textit{full subcategory} \index{category!full subcategory} of $\textbf{C}$ when $\textbf{C}$ has no arrows $A \rightarrow B$ other than the ones already in $\textbf{D}$, i.e., for any $\textbf{D}$-objects $A$ and $B$, 
	\begin{equation*}
	\text{Hom}_{\textbf{D}}(A, B) = \text{Hom}_{\textbf{C}}(A, B).
	\end{equation*}
\end{definition}
\begin{example}
	The\index{category!of finite sets} category $\textbf{FinSet}$ of finite sets---the category whose objects are all finite sets and whose morphisms are all the functions between them---is a subcategory of $\textbf{Set}$. In fact, it is a \textit{full} subcategory. \par 
	The category of abelian groups is a (full) subcategory of the category of groups. \par 
	If $\textbf{C}$ is the category that has as objects those parts of $\mathbb{R}^n$ that are open, and for morphisms those mappings between objects that are continuous, then a subcategory $\textbf{D}$ of $\textbf{C}$ is formed by restricting to mappings that have a derivative, where a rule of basic calculus shows that $\textbf{D}$ has composition. A further subcategory of $\textbf{D}$ could be got by further restricting to those mappings that have all derivatives (i.e., the \textit{smooth} ones). There are many other important examples of subcategories that we will encounter throughout this book. 
\end{example}
There are a number of other useful things one can do with categories, not to mention the important things one can do and find within categories. Discussion of such matters is left to emerge organically throughout the book. \par 
	The real power of category theory, however, only really comes into its own once it is realized how, by putting everything on the same ``plane," we can consider principled relations \textit{between} categories. This is what we discuss in the next section. 
	\section{Prelude to Sheaves: Presheaves} 
	\subsection{Functors}
	It is often said that category theory privileges relations over objects. But a category itself can be considered as an object, and then a natural question is ``what do relations between categories look like?" If a category is a context for studying a specific type of mathematical object and the network of relations entertained between those objects, a \textit{functor} is a principled way of comparing categories, translating the objects and actions of one category into objects and actions in another category in such a way that certain structural relations are preserved through this translation. As a way of moving in a controlled way \textit{between} categories, one can initially think of a functor as doing any of the following things: specifying data locally; producing a picture of the source category in the target category, modeling one category or some aspect of that category within another; ``realizing" an abstract theory of some structured notion (such as a `group') in a certain background or on a specific ``stage"; taking advantage of the methods available in the target category to analyze the source category; converting a problem in one category into another where the solution might be more readily apparent; forgetting or deliberately losing some information, perhaps in order to examine or identify those features more robust to variations or to ease computation. But underneath these different interpretations or uses is a very simple requirement: a functor just transforms objects and maps in the domain category into objects and maps in the codomain category, in such a way that two equations are satisfied. Functors also come in two ``flavors," depending on their direction or variance. Formally, 
	\begin{definition}
		A \textit{(covariant) functor} \index{functor!defined} $F: \textbf{C} \rightarrow \textbf{D}$ between categories \textbf{C} and \textbf{D} is an assignment of 
		\begin{enumerate}
			\item an object $F(c) \in \textbf{D}$ for every object $c \in Ob(\textbf{C})$; and 
			\item a morphism $F(f): F(c) \rightarrow F(c')$ in $\textbf{D}$ for every morphism $c \rightarrow c'$ in $\textbf{C}$, 
		\end{enumerate}
	which assignments moreover satisfy the following two axioms: 
		\begin{enumerate}
			\item For any object $c$ in \textbf{C}, $F(id_c) = id_{Fc}$ (`$F$ of the identity on $c$ is the identity on $Fc$'); 
			\item For any composable pair $f,g$ in \textbf{C}, $F(g) \circ F(f) = F(g \circ f)$.
		\end{enumerate} 
		A \textit{(contravariant) functor}, i.e., a functor $F: \textbf{C}^{op} \rightarrow \textbf{D}$, is defined in the same way on objects, but differently on morphisms (reversing the direction of all arrows). Explicitly, to each morphism $f: c \rightarrow c' \in \textbf{C}$ it assigns a morphism $Ff: Fc' \rightarrow Fc \in \textbf{D}$. This assignment must satisfy the same identity axiom as above, but for any composable pair $f,g$ in \textbf{C}, we now have $F(f) \circ F(g) = F(g \circ f)$ (note the change in direction). \par 
		All the information of this definition is displayed below (the covariant case on the left and contravariant case on the right, and with identity maps omitted except for one object): 
		\par 
		\vspace*{0.5cm} 
	\small
	\begin{tikzcd}
		c \arrow[d, "f", swap] \arrow[dd, bend left, "{g \circ f}"] \arrow[loop left,distance=1em]{}{\mathrm{id}_c} & & F(c)\arrow[d, "F(f)", swap] \arrow[dd, bend left = 50, "{F(g \circ f)}"] \arrow[loop left,distance=1em]{}{{F(id_c)}} & & & & c \arrow[d, "f",swap] \arrow[loop left,distance=1em]{}{\mathrm{id}_c} \arrow[dd, bend left, "{g \circ f}"] & & F(c) \arrow[loop left,distance=1em]{}{{F(id_c)}}  \\
		c' \arrow[d, "g", swap] & &  F(c')\arrow[d, "F(g)",swap] & & & & c' \arrow[d, "g",swap] & &  F(c')\arrow[u, "F(f)"]  \\
		c'' & & F(c'') & & & & c''  & &  F(c'')\arrow[u, "F(g)"] \arrow[uu, bend right = 50, "{F(g \circ f)}",swap] \\
		\textbf{C} \arrow[rr,"F"] & & \textbf{D} & & & & \textbf{C}^{op} \arrow[rr,"F"] & &\textbf{D} \\
	\end{tikzcd}
	\end{definition} \noindent 
Functors will usually be denoted with upper-case letters, such as $F, G, P$, etc., though we may occasionally use a more evocative name, to indicate what the functor does. \par 
The truth of the frequently-cited claim of Eilenberg and Maclane that ``the whole concept of a category is essentially an auxiliary one; our basic concepts are essentially those of a functor and of a natural transformation"\footnote{See \cite{maclane_general_1945}, 247.} proves itself with time to anyone who works with categories. In addition to their intrinsic interest, functors are of special interest to us because of their role in the definition of \textit{presheaves}.   
\begin{definition}
	A (set-valued) \textit{presheaf} \index{presheaf!defined} on \textbf{C}, for \textbf{C} a small category, is a (contravariant) functor $\textbf{C}^{op} \rightarrow \textbf{Set}$.\footnote{By `small', one means that the category has no more than a set's worth of arrows.\index{category!small} \par 
	Incidentally, the reader who wonders why, if a presheaf is just a (contravariant) functor, we bother giving it two names, might be satisfied by the fun notion, used by the nLab authors, of a \textit{concept with an attitude}, meant to capture those situations in math when one and the same concept is given two different names, one of the names indicating a specific perspective or ``attitude" suggesting what to do with the objects, or the sorts of things one might expect to be able to do with them. In renaming a (set-valued) contravariant functor as a presheaf, then, we have a concept with an attitude, specifically looking forward to \textit{sheaves}.}   
\end{definition} \noindent 
A presheaf will often be thought of as consisting of some specification or assignment of local data, according to the ``shape" of the domain category; a sheaf will emerge as a special sort of presheaf in that its local data can be glued or patched together (locally). Before addressing in more detail the nature of presheaves, we pause to provide some examples of functors in general.
\subsection{Examples of Functors}
Functors appear all over mathematics. But perhaps the ``lowest hanging" examples of functors is supplied by those mathematical entities of a certain type each of which individually assembles into a category, thus giving us a particularly simple ``way in" to categories as objects of study in their own right. In these cases, we would expect that a functor between such objects, now each regarded as individual categories in their own right, would recover the usual important relations that are expected to obtain among such objects. This is what is illustrated by the following two examples.   
\begin{example}
	Recall that a preorder $\mathcal{X} := (X, \leq)$ is defined as a set $X$ together with a reflexive and transitive binary relation $\leq$. Recall also that we can transform a given preorder into a category $\mathcal{X}$ by defining, for every pair of objects $x,x' \in X$, the hom-set $\text{Hom}_{\mathcal{X}}(x,x')$ as either empty (in case the pair $(x,x')$ is not related by $\leq$) or as consisting of the unique morphism $x \rightarrow x'$ (just in case $x \leq x'$), making the composition formula completely determined. \par 
	For a morphism of preorders $f: (X, \leq_X) \rightarrow (Y, \leq_Y)$, we take an object $x \in Ob(\mathcal{X}) = X$ and assign it the object $f(x) \in Y = Ob(\mathcal{Y})$. Given a morphism $f: x \rightarrow x'$ in $\mathcal{X}$, we have automatically that $x \leq x'$ and by the definition of a morphism of preorders as order-preserving, we will have $f(x) \leq f(x')$, so all we have to do is assign to $f$ the unique morphism $f(x) \rightarrow f(x')$ in $\mathcal{Y}$. In other words, the usual preorder maps, i.e., monotone maps, are nothing other than functors between such categories. \par 
	In short, regarding the preorders (or posets) $\mathcal{X}$ and $\mathcal{Y}$ as categories, then a covariant functor from $\mathcal{X}$ to $\mathcal{Y}$ is nothing other than a monotone (order-preserving) function, while a contravariant functor is just an anti-tone (order-reversing) function (i.e., whenever $x \leq x'$ in $\mathcal{X}$, then $f(x') \leq f(x)$ in $\mathcal{Y}$).\index{antitone map}  
	\end{example} 
	\begin{example}  
	Recall that each monoid $\mathcal{M}$ (and each group $\mathcal{G}$) can be regarded as its own category. Explicitly, we saw earlier that a monoid $(M, e, \cdot)$ can be considered as a category $\mathcal{M}$ with one object and with hom-set equal to $M$ itself, where the identity morphism comes from the monoid identity $e$ and the composition formula from the monoid multiplication $\cdot: M \times M \rightarrow M$.\footnote{Notice that, getting a little ``meta," when we said that each monoid could be regarded as its own category, we were really just constructing a \textit{functor} $i: \textbf{Mon} \rightarrow \textbf{Cat}$, one that takes a monoid to its corresponding category!} \par 
	Then if we have monoids \index{monoid} regarded as categories, we might hope that a covariant functor from the one to the other would just be a monoid homomorphism (the ``usual" notion of a morphism between monoids), where this is a map $\phi: M \rightarrow N$ between two monoids that respects the structure in the sense that 
	\begin{equation*}
	\phi(m \cdot_M m') = \phi(m) \cdot_N \phi(m') \text{ and } \phi(e_M) = e_N.  
	\end{equation*} 
	So each monoid homomorphism is indeed just a functor between such one-object categories. It is easy to see that the above equations are the same as defining a covariant functor between $\mathcal{M}$ and $\mathcal{N}$, where these monoids are each regarded as categories. A contravariant functor from $\mathcal{M}$ to $\mathcal{N}$, for its part, is then just a monoid morphism that flips the elements, i.e., $\phi(m \cdot m') = \phi(m') \cdot \phi(m)$.  \par 
	Since a group is just a monoid in which every element is invertible, a similar story can of course be told of groups, \index{group} i.e., a group can be regarded via a functor $\textbf{Grp} \rightarrow \textbf{Cat}$ as a category with one object such that every morphism is an isomorphism and a functor between such categories is just a group homomorphism. In this context, we might also mention that there exists a functor $Core: \textbf{Mon} \rightarrow \textbf{Grp}$ that\index{functor!Core} takes a monoid $M$ and returns the subset of invertible elements of $M$, which of course forms a group, typically called the \textit{core} of $M$. There is a related functor $\textbf{Cat} \rightarrow \textbf{Grpd}$ sending a category $\textbf{C}$ to the largest groupoid\index{groupoid} inside \textbf{C}, also called its core.\footnote{A groupoid is just like a group except that it can have more than one object, as our earlier discussion of ``oidification" \index{oidification} would have suggested. More formally, a \textit{groupoid} is a category such that every morphism is an isomorphism; a morphism between groupoids is also just a functor.}
\end{example}
Moving beyond examples where the functors pass between categories that are fundamentally the same sort of structure, we need to begin to appreciate some of the other important things functors do. 
\begin{example}
	In many settings, one might want to transfer one system of objects that present themselves in a certain way in one context to another context where irrelevant or undesirable (e.g., noisy) features are suppressed, while simultaneously preserving certain basic qualitative features. In the definition of a category, and indeed in the definition of many mathematical objects, typically one specifies (i) underlying data, together with (ii) some extra structure, which in turn may satisfy (iii) some properties. One obvious thing to do when considering some category $\textbf{C}$ is to deliberately ``forget" or ignore some or all of the structure or the properties carried by the source category by passing, via a functor, to another category. This process informally describes what are usually called \textit{forgetful functors}, \index{functor!forgetful} which provide us with a large source of examples.\par
	There are many examples of this where \textbf{Set} is the target category, since many important categories are sets with some structure; however, forgetful functors need not have \textbf{Set} for the target category. For instance, since a group is just a monoid $(M, e, \cdot)$ with the extra property that every element $m \in M$ has an inverse, this means that to every group we can assign its underlying monoid, and every group homomorphism gets assigned to a monoid homomorphism between its underlying monoids, simply by ``forgetting" the extra conditions on a group. Thus, there is a forgetful functor $U: \textbf{Grp} \rightarrow \textbf{Mon}$. \par 
	While the ``forgetting" terminology might vaguely suggest some sort of (possibly pejorative) loss of information, another way of looking at the same process is as extracting and emphasizing only the ``important" features of the objects under study. An illustration of this comes from \textit{detectors}, \index{detector} which do indeed act to forget or lose information, while preserving fundamental features of the underlying signal, and this is regarded as exactly what is useful about such tools, since what is removed is ``clutter," leaving us with a compressed representation of the original information (with the effect that the result of applying the functor might be more robust to variations, more relevant to a particular application, simpler for computation purposes, etc.). We know that a \textit{signal} is a collection of (local) measurements related to one another, and the topology on these measurements tells us how a measurement responds to noise, e.g., a signal over a discrete set is typically either not changed by noise at all or it changes drastically, while a signal over a smoother space may depend less drastically on perturbations. As already anticipated, as a forgetful functor, a detector acts to \textit{remove} something---specifically, it acts to remove the topological structure from the signal (which may have the effect of quantizing those signals that previously varied smoothly). \par 
	As a specific and simple instance of this, consider a \textit{threshold detector}.\footnote{This idea for this threshold example comes from \cite{robinson_topological_2014}.} A threshold detector takes a continuous function $f \in \textbf{Cont}(\mathbb{R})$ and returns the open set on which $f(x) > T$ for some threshold $T \in \mathbb{R}$. The domain of this functor is the category $\textbf{Cont}(\mathbb{R})$ with continuous (real-valued) functions for objects and for morphisms the functions $f \rightarrow g$ whenever $f(x) > g(x)$ for all $x \in \mathbb{R}$. The threshold detector is then the functor $D$ that assigns to each $f \in \textbf{Cont}(\mathbb{R})$ the open set $D(f) = \{x \in \mathbb{R}: f(x) > T\}$, i.e., it lands in the category \textbf{Op} of open sets of $\mathbb{R}$ with morphisms given by subset inclusion.  Moreover, one can see that if $f \rightarrow g$, then we will have $D(g) \subseteq D(f)$, making $D$ a contravariant functor $\textbf{Cont}(\mathbb{R}) \rightarrow \textbf{Op}(\mathbb{R})$. \par 
	Forgetful functors in general frequently can tell us interesting things about the source category. For instance, we have a functor $U: \textbf{Cat} \rightarrow \textbf{Grph}$, which informs us of the fact that categories have underlying graphs. While we did not explicitly adopt this description above, we could have equally defined a category by saying that the data of a (small) category involves a set of objects (sometimes denoted $\textbf{C}_0$), a set of morphisms (denoted $\textbf{C}_1)$, and a diagram $\begin{tikzcd} \textbf{C}_1 \arrow[r, shift left = 1ex, "s"] \arrow[r, "t", swap] & \textbf{C}_0 \end{tikzcd}$, together with some structure (composition and identities) and properties (unit and associativity axioms). Comparing this formulation to the definition of a directed graph (see \ref{example: graph}), it should come as no surprise that there is such a functor sending categories to their underlying directed graph (or ``1-globular set," on which more below). \par  
	Forgetful functors often come paired with corresponding \textit{free functors}.\index{functor!free} For instance, corresponding to the ``underlying graph" functor $U$, there exists the \textit{free category functor} $F: \textbf{Grph} \rightarrow \textbf{Cat}$, which we have in fact already described. Recall that given a directed graph $G$, we can create a category $\textbf{Pth}(G)$,\index{category!of paths} the category of paths through $G$, with objects the nodes of $G$ and arrows the paths through $G$. The resulting category of paths of a graph $G$, $\textbf{Pth}(G)$, gives us the \textit{free category generated by} $G$, which can be thought of as being the result of freely adding all paths (all possible composite arrows) as well as all the identity arrows to a given graph. The resulting category has the same set of objects, i.e., vertices, as the original graph, but it will in general have a larger set of morphisms (the hom-set $Hom(v,v')$ in the resulting category will consist of all the paths in the graph $G$ from $v$ to $v'$), and graph homomorphisms then extend into unique (covariant) functors. \par 
	In short, this path construction gives rise to a functor $F: \textbf{Grph} \rightarrow \textbf{Cat}$, called the \textit{free category functor}. $\textbf{Pth}(G)$ is always the \textit{largest} category generated by $G$. On the other hand, $G$ also generates a \textit{smallest} category by taking the quotient of $\textbf{Pth}(G)$ by the relation that identifies two paths that share the same source and the same target. In this connection, any category $\textbf{C}$ can be obtained as a quotient of the corresponding category of paths of its underlying graph, under the equivalence relation identifying two paths if and only if they have the same composite in $\textbf{C}$.\footnote{It is worth adding, in passing, that the pair of functors \begin{tikzcd}[ampersand replacement=\&] \textbf{Cat} \arrow[r, shift left = 1ex, "U"] \& \textbf{Grph} \arrow[l, "F"] \end{tikzcd} is extremely significant, for it forms an important adjunction which gives rise to a particularly special ``monad" that is a starting point for the generalization to $n$-categories\index{category!n-category} (on which more below).} \par 
	In short, constructions and results established in the context of graphs can be applied to categories, once we `forget' about composition; and conversely, results concerning categories can be applied to graphs by simply replacing a graph by its category of paths. 
\end{example}
\begin{example} 
	We just saw that categories have underlying graphs. There is the important related notion of a \textit{diagram} \index{functor!diagram} in a category \textbf{C}, a notion that in some sense captures a generalized idea of a subgraph of a given category's underlying graph. A diagram is defined as a functor $F: \textbf{J} \rightarrow \textbf{C}$ where the domain, called the \textit{indexing category} \index{category!indexing} or \textit{template}, is a small category. Typically, one thinks of the indexing category as a directed graph, i.e., some collection of nodes and edges, serving as a template defining the shape of any realization of that template in \textbf{C}, and which may also specify some commutativity conditions on the edges which are to be respected by \textbf{C}. Then a diagram can be regarded as something like an instantiation or realization in \textbf{C} of a particular template \textbf{J}. Each node in the underlying graph of the indexing category is instantiated with the objects of \textbf{C}, while each edge is instantiated with a morphism of \textbf{C}. If we write the objects in the index category \textbf{J} as $i, j, ...$, and the values of the functor $F: \textbf{J} \rightarrow \textbf{C}$ in the form $F_i, F_j, ...$ (or $F(i), F(j),...$), then a diagram is a family of objects $F(i)$ of \textbf{C} indexed by the nodes of \textbf{J} and a family of arrows $F(e)$ of \textbf{C} indexed by the edges of \textbf{J}; accordingly, one sometimes speaks of a diagram $F$ as a $\textbf{J}$-indexed set, or $\textbf{J}^{op}$-parametrized set (depending on the variance of the functor). Functoriality demands that any of the composition relations (in particular, commutative diagrams) that obtain in \textbf{J} carry over (under the action of $F$) to the image in $\textbf{C}$.\par 
	We will have a lot more to say about this perspective in the next chapter. For now, let us look at a few concrete illustrations of this. First consider the category 
	\begin{center} 
		\begin{tikzcd}[framed]
			\bullet{0} \arrow[loop above,"{id_0}"] \arrow[r,"f"] & \bullet{1} \arrow[loop above,"{id_1}"]
		\end{tikzcd}
	\end{center} 
often called \textbf{2}.\footnote{$\textbf{2}$ is isomorphic to the linear order $\textbf{[1]}$, so one will occasionally also see it go by that name.} With such a category for our indexing category, \index{category!indexing} a (set-valued) diagram yields a category that has as objects all the functions from one set to another set, and as morphisms the commutative squares between those arrow-objects. In more detail: a morphism from object $f:a \rightarrow b$ to object $g:c \rightarrow d$ will be a pair of functions $\langle h, k \rangle$ such that 
	\begin{center} 
		\begin{tikzcd}
			a \arrow[r, "h"] \arrow[d, "f", swap]
			& c \arrow[d, "g"] \\
			b \arrow[r, "k", swap]
			& d
		\end{tikzcd} 
	\end{center} commutes. Composition is componentwise, i.e., $\langle j, l \rangle \circ \langle h, k \rangle = \langle j \circ h, l \circ k \rangle$, and the identity arrow for $f: a \rightarrow b$ will be the function pair $\langle id_a, id_b \rangle$. Look familiar? It should. This is just the \textit{arrow category} introduced earlier!\par 
	Suppose instead we take for our indexing category $\textbf{3}$, or \textbf{[2]}, the linear order category with length 2 
	\begin{center} 
		\begin{tikzcd}[framed]
			\bullet{0} \arrow[r, "f"] \arrow[dr, "{g \circ f}", swap]  & \bullet{1} \arrow[d, "g"] \\
			& \bullet{2}
		\end{tikzcd} 
	\end{center} 
Then a diagram on this category just acts to pick out as objects commutative triangles. As a final example, taking the category $\textbf{2} \times \textbf{2} \times \textbf{2}$ as our indexing category just serves to pick out as objects commutative cubes in the target category: 
	\begin{center} 
		\begin{tikzcd}[row sep=1.5em, column sep = 1.5em]
			\bullet \arrow[rr] \arrow[dr] \arrow[dd] &&
			\bullet \arrow[dd] \arrow[dr] \\
			& \bullet \arrow[rr] &&
			\bullet \arrow[dd] \\
			\bullet \arrow[rr] \arrow[dr] && \bullet \arrow[dr] \\
			& \bullet \arrow[rr] \arrow[uu,<-] && \bullet 
		\end{tikzcd} 
	\end{center}  
	As we will see, this ``diagram approach" can be significantly generalized and can even be used to provide definitions of $n$-categories, specifying the data for an $n$-category\index{category!n-category} as a diagram (presheaf) $A: \Sigma^{op} \rightarrow \textbf{Set}$, where $\Sigma$ is some category of shapes, and the functor yields, for each shape, a set of ``cells" of that shape. The globular shapes (to be defined below) are the most basic cell shape.\footnote{Moreover, this diagram approach already suggests a more general definition of presheaves: for categories \textbf{C} and \textbf{J}, a \textbf{C}-presheaf on \textbf{J} can be defined as a contravariant functor from \textbf{J} to \textbf{C}. While this more general definition is perfectly coherent (and is useful for achieving greater generality), presheaves are ordinarily regarded as valued in \textbf{Set}. While this is not entirely necessary, there is also good reason for it. In brief, it has to do with the fact that the category of sets is somewhat special: the usual categories are \textit{enriched} over sets, by which we mean that given a pair of objects $X, Y \in \textbf{C}$, we can form $Hom_{\textbf{C}}(X,Y)$, an object of \textbf{Set}. Moreover, another factor here has to do with the fact that only set-valued functors are \textit{representable} (on which, more below). However, in most categories $\textbf{C}$, the hom-sets $Hom_{\textbf{C}}(X,Y)$ are richer than just sets, e.g., for $\textbf{C}$ a category of chain complexes, $Hom_{\textbf{C}}(X,Y)$ is an abelian group.}
\end{example} 
\begin{example}
	\label{example: diagram 1}
	Expanding on the previous perspective, assume we are given as indexing category $\textbf{J} :=$
	\begin{center} 
		\begin{tikzpicture}[framed, scale=0.7]
		\tikzset{vertex/.style = {shape=circle,draw, fill=black, minimum size=3pt, inner sep =0pt}}
		\tikzset{edge/.style = {->,> = latex'}}
		\node[vertex] (a) [label=above:{$a$}] at  (0,2) {};
		\node[vertex] (a5) [label=above:{$d$}] at  (4.6,1.2) {};
		\node[vertex] (a1) [label=below:{$b$}] at (1.5,1.2) {};
		\node[vertex] (a2) [label=above:{$c$}] at (3,2) {};
		\node[vertex] (a6) [label=above:{$q$}] at (4.5,0) {};
		\node[vertex] (a7) [label=above:{$r$}] at (6.5,0) {};
		
		\path[-latex] (a6) edge node[auto] {$i$} (a7);
		
		\path[-latex] (a) edge node[pos=0.5, above] {$f$} (a1);
		
		\path[-latex] (a2) edge node[above] {$g$} (a1);
		\path[-latex] (a2) edge node[pos=0.5, above] {$h$} (a5);
		\end{tikzpicture} 
	\end{center} 
	Now let the diagram $F: \textbf{J} \rightarrow \textbf{Set}$ be given on objects by 
	\begin{align*}
	& F(a) = \{1,2\}, \hspace*{1em} F(b) = \{1,2\}, \hspace*{1em} F(c) = \{1,2,3\}, \\
	& F(d) = \{1,2,3,4\}, \hspace*{1em} F(q) = \{1,2,3\}, \hspace*{1em} F(r) = \{1,2\}. 
	\end{align*} 
	and on morphisms by 
	\begin{align*}
	& F(f) = 1 \mapsto 1, 2 \mapsto 2; \\
	& F(g) = 1 \mapsto 1, 2 \mapsto 2, 3 \mapsto 1; \\
	& F(h) = 1 \mapsto 1, 2 \mapsto 2, 3 \mapsto 4; \\
	& F(i) = 1 \mapsto 2, 2 \mapsto 1, 3 \mapsto 1. 
	\end{align*} 
	This can be pictured as follows: 
	\par 
	\begin{center} 
		\begin{tikzpicture}[framed, scale=0.7]
		\tikzset{vertex/.style = {shape=circle,draw, fill=black, minimum size=3pt, inner sep =0pt}}
		\tikzset{edge/.style = {->,> = latex'}}
		\node[vertex] (a) [label=left:{\footnotesize $(a,1)$}] at  (0,2) {};
		\node[vertex] (a2) [label=left:{\footnotesize$(a,2)$}] at  (0,1) {};
		\node[vertex] (d) [label=right:{\footnotesize$(d,1)$}] at  (6.5,1.2) {};
		\node[vertex] (d2) [label=right:{\footnotesize$(d,2)$}] at  (6.5,0) {};
		\node[vertex,gray] (d3) [label=right:{\footnotesize $(d,3)$}] at  (6.5,-1) {};
		\node[vertex] (d4) [label=right:{\footnotesize$(d,4)$}] at  (6.5,-2) {};
		
		\node[vertex] (b) [label=above:{\footnotesize $(b,1)$}] at (2,1) {};
		\node[vertex] (b2) [label=below:{\footnotesize $(b,2)$}] at (2,-0.1) {};
		\node[vertex] (c) [label=above:{\footnotesize $(c,1)$}] at (4,2) {};
		\node[vertex] (c2) [label=below:{\footnotesize $(c,2)$}] at (4,1) {};
		\node[vertex] (c3) [label=below:{\footnotesize $(c,3)$}] at (4,-0.1) {};
		\node[vertex] (q) [label=left:{\footnotesize $(q,1)$}] at (6.5,-3) {};
		\node[vertex] (q2) [label=left:{\footnotesize $(q,2)$}] at (6.5,-4) {};
		\node[vertex] (q3) [label=left:{\footnotesize $(q,3)$}] at (6.5,-5) {};
		\node[vertex] (r) [label=right:{\footnotesize $(r,1)$}] at (8.5,-3) {};
		\node[vertex] (r2) [label=right:{\footnotesize $(r,2)$}] at (8.5,-4.5) {};
		
		\path[-latex, dashed] (q) edge node[auto] {} (r2);
		\path[-latex, thick] (q2) edge node[auto] {} (r);
		\path[-latex, thick] (q3) edge node[auto] {} (r);
		
		\path[-latex,ultra thick] (a) edge node[pos=0.5, above] {} (b);
		\path[-latex] (a2) edge node[pos=0.5, above] {} (b2);
		\path[-latex,ultra thick] (c) edge node[above] {} (b);
		\path[-latex] (c2) edge node[above] {} (b2);
		\path[-latex,ultra thick] (c3) edge node[above] {} (b);
		\path[-latex,ultra thick] (c) edge node[pos=0.5, above] {} (d);
		\path[-latex] (c2) edge node[pos=0.5, above] {} (d2);
		\path[-latex,ultra thick] (c3) edge node[pos=0.5, above] {} (d4);
		\end{tikzpicture} 
	\end{center} 
	This realization affords us a concrete illustration of another important construction, the \textit{category of elements}, which will be used later on in explaining why there are different thicknesses of arrows in this picture.   
\end{example} 
There are also many functors that recover important established constructions that appear within the context of more specialized study of certain structures. The following examples show a few of these.
\begin{example} 
	\label{example: ncoloring}
	Consider the set of vertex colorings of an undirected connected graph subject to the condition that no adjacent vertices are assigned the same color. There is a functor $nColor: \textbf{UCGraph}^{op} \rightarrow \textbf{Sets}$\index{functor!nColoring} that takes an undirected connected graph to the set of $n$-colorings of its vertices, i.e., it colors an undirected connected graph with \textit{at most} $n$ colors. For any $G$, $nColor(G)$ will be the set of all $n$-colorings of $G$. Note that if graph $G$ has an $n$-coloring, then clearly each of its subgraphs will have an $n$-coloring. Moreover, for any $f: G \rightarrow G'$, $nColor(f)$ will be the function restricting the colorings of $G'$ to those of $G$.\footnote{Looking ahead, a sheaf will be defined as a particular presheaf satisfying certain properties with respect to ``covers" of the objects of the domain category. Anticipating this, the $nColor$ functor will in fact turn out to be a sheaf since, if $\{G_i \hspace*{0.3em}| \hspace*{0.3em} i \in I\}$ covers $G$, and if $\{c_i \in nColor(G_i) \hspace*{0.3em}| \hspace*{0.3em} i \in I \}$ is a family of colorings such that the colorings agree on intersections among the $G_i$, then the $c_i$'s induce a unique coloring of the entire graph $G$.}
\end{example}  
\begin{example}
	Suppose we have two parallel arrows 
	\begin{center}
		\begin{tikzcd}
		 X \arrow[r, shift right = 1.5, "g",swap] \arrow[r, "f"] & Y,
		\end{tikzcd} 
	\end{center}
which, for simplicity, we may take as living in \textbf{Set}. We can then define something called the \textit{equalizer of}\index{equalizer} $f$ and $g$ as a set
	\begin{equation*}
	Eq(f,g) := \{x \in X \hspace*{0.25em}|\hspace*{0.25em} f(x) = g(x) \}. 
	\end{equation*}
	This set of elements of $X$ for which the two functions coincide will of course be a subset of $X$, so we can consider the \textit{inclusion} $e: Eq(f,g) \rightarrow X$. Altogether, this will result in the diagram
	\begin{center} 
		\begin{tikzcd}
			Eq(f,g) \arrow[r, dashed, "e"] & X \arrow[r, shift right = 1.5, "g",swap] \arrow[r, "f"] & Y
		\end{tikzcd} 
	\end{center} 
that in fact commutes, i.e., $f \circ e = g \circ e$. The equalizer construction really involves both the set $Eq(f,g)$ and this map making the diagram commute. This construction can be defined in more general categories. In this setting, the equalizer emerges as a special sort of object, namely as the \textit{universal} object. This means that the \textit{equalizer} of the parallel arrows (in this case, $f, g$) is a morphism $e: \text{Eq} \rightarrow X$ such that $f \circ e = g \circ e$, where this is \textit{universal} with this property, meaning that given \textit{any other} morphism $u: Z \rightarrow X$ in the category such that $f \circ u = g \circ u$, there exists a unique $v: Z \rightarrow \text{Eq}$ such that $u$ factors through $v$, i.e., $e \circ v = u$.\footnote{This is a particular instance of something called a \textit{limit}, which will be discussed in more detail further on.}\par 
Since a graph can be defined, as we have seen, as a pair of functions \begin{tikzcd} A \arrow[r, shift left = 1ex, "s"] \arrow[r, "t", swap] & V, \end{tikzcd} where $A$ stands for arrows and $V$ for vertices, and where $s$ just picks out the source vertex of an arrow and $t$ the target vertex, consider that for any graph $G$ we can find its set of length one loops via the equalizer construction $Eq(G)$:
	\begin{center} 
	\begin{tikzcd}
		Eq(s,t) \arrow[r, dashed, "e"] & A \arrow[r, shift right = 1.5, "s",swap] \arrow[r, "t"] & V.
	\end{tikzcd} 
	\end{center} 
	This equalizer assignment is functorial, since given a graph homomorphism $G \rightarrow G'$, there is an induced function $Eq(G) \rightarrow Eq(G')$. \par 
	 In a similar fashion, various categories---including graphs, reflexive graphs, discrete dynamical systems, simplicial sets, the category of elements (see below)---support a construction that allows us to count the \textit{connected components} in that category. For concreteness, we stick with the case of the category of graphs, \textbf{Grph}, and consider the ``connected components" \index{functor!connected components} functor $\Pi_0: \textbf{Grph} \rightarrow \textbf{Set}$. It is obtained via the dual to the equalizer, namely the \textit{coequalizer} construction \index{coequalizer}:
	 \begin{center} 
	 \begin{tikzcd}
	 	A \arrow[r, shift right = 2, "{s}", swap] \arrow[r, "{t}"] & V \arrow[r, "q",dashed] & Coeq(s,t). 
	 \end{tikzcd}
	 \end{center}
	 We define $Coeq(s,t)$ as $V \big/$$\sim$, where $\sim$ is an equivalence relation on $V$, i.e., $s(x) \sim t(x)$ for all $x \in A$, and where $q$ is the quotient function $q: V \rightarrow V \big/$$\sim$. This construction accordingly acts to identify all arrows where the source of one arrow is equal to the target of the other. In other words, all we are doing is picking out the connected components of the graph. This assignment of the set of connected components of a graph can be shown to be functorial as well. In later examples, we will see instances of this functor in action. 
\end{example}
Here is an example of a different flavor, one that also ties together a number of constructions introduced thus far. 
\begin{example} \label{qua}
	There are natural language expressions that we use all the time to express that someone or something has a certain property \textit{qua} (or \textit{as}) one thing but not \textit{qua} some other thing. For instance, one might say  
	\begin{quote}
		John is fair \textit{qua} father, but not \textit{qua} politician,  
	\end{quote}
or 
\begin{quote}
	Maria is inspirational \textit{qua} teacher, but not \textit{qua} basketball player.
\end{quote}
We often make use of judgments involving the logic of \textit{qua}. Suppose someone asks you whether `Abe is honest'. You might intelligibly answer, ``Well, yes and no." You might elaborate by saying ``It depends: in some respects/aspects, Abe is honest; in other respects/aspects, not so much." You may debate about which of the aspects are relevant, or most relevant, and also how Abe's behavior, under a particular aspect (in a particular respect), should be interpreted (as honest or dishonest). However, in general, this sort of answer and the ensuing discussion make sense. Once agreement about these matters (which aspects are relevant, etc.) has been achieved, one can arrive at a judgment about Abe's honesty, one that accommodates the fact that the answer depends on the aspect under consideration. \par 
	We might conceptualize this situation category-theoretically, using a particular  ``category of aspects" or \textit{qua} category.\footnote{The idea for this, and the key definitions provided below (as well as the example, with mostly trivial modifications), comes from \cite{la_palme_reyes_models_1999}.} In this setting, we will be able to model things like `honesty of Abe under aspect $A$', and moreover model the assembly of global judgments of the type `Abe is honest', `Abe is not honest', `Abe is dishonest', etc., from this data of judgments about Abe's honesty \textit{qua} the various relevant aspects. \par 
	Let us first define something we could call the \textit{nominal category} $\textbf{CN}$. 
	\begin{definition}
		The \textit{nominal category},\index{category!nominal} \textbf{CN}, has for
		\begin{itemize}
			\item objects: CNs (count nouns) relevant to discussion, e.g., `a student', `a coworker', `a husband', `a parent', `a family man', `a student, a coworker, and a family man', written as 
			\begin{center} 
				\begin{tikzcd}
					\boxed{\text{a s}}, \boxed{\text{a c}}, \boxed{\text{a h}}, \boxed{\text{a p}}, \boxed{\text{a f}}, \boxed{\text{a scf}}
				\end{tikzcd}
			\end{center}
			\item morphisms: ``identification" postulates of the form (copula connecting nouns)
			\begin{center} 
				\begin{tikzcd}
					\boxed{\text{a f}} \arrow[r, "is"] & \boxed{\text{a p}}, \boxed{\text{a scf}} \arrow[r, "is"] & \boxed{\text{a f}} 
				\end{tikzcd}
			\end{center} 
			where these are meant to capture the identifications frequently used in natural languages, such as `a family man is a parent', `a student, a coworker, and a family man is a family man', `a dog is an animal', etc. \par 
			Identity morphisms are those particular axioms of the form 
			\begin{center} 
				\begin{tikzcd}
					\boxed{\text{a f}} \arrow[r, "is"] & \boxed{\text{a f}} 
				\end{tikzcd}
			\end{center} 
		\end{itemize}
		Composition is given by stringing together identifications in the obvious way, i.e., whenever we have two arrows, the co-domain of one as the domain of the other, we complete the graph by adding an arrow that is the composition of the two arrows, and where this corresponds to the rule of inference, normal in natural languages, from things like `a human is a primate' and `a primate is a mammal' to `a human is a mammal'.   
	\end{definition} 
	The axiom arrows of this category can be thought of as supplying a \textit{system of identifications}, where this replaces a notion of \textit{equality} between kinds (since equality is a relation that might obtain only between the members of a given kind). This category is assumed to be posetal, with at most one arrow between two objects.\par 
Then, we know that for a category \textbf{C}, we can define the \textit{arrow category} of \textbf{C}, denoted $\textbf{C}^{\rightarrow}$.\index{category!arrow} There is actually an embedding (on which more later) from $\textbf{C}$ to $\textbf{C}^{\rightarrow}$, which allows us to \textit{identify} any object $A$ of $\textbf{C}$ with the object $A \xrightarrow{1_A} A$ in $\textbf{C}^{\rightarrow}$. And with this arrow category construction, we can define our important category. 
	\begin{definition}
		The \textit{qua category}\index{category!qua} (called the \textit{aspectual category} in \cite{la_palme_reyes_models_1999}) of $\textbf{CN}$, \textbf{Qua(CN)} (or just \textbf{Qua}), is defined as $(\textbf{CN}^{\rightarrow})^{op}$. 
	\end{definition}
	Because of how $\textbf{CN}$ was defined, we will have at most one morphism from an object $A$ to an object $B$. When such a morphism exists, we can see $A \rightarrow B$ as $A \text{ qua } B$, e.g., 
	\begin{center} 
		\begin{tikzcd}
			\boxed{\text{a f}} \arrow[r] & \boxed{\text{a p}}
		\end{tikzcd}
	\end{center} 
	will be to look at `a family man \textit{qua} a parent'. Identifying $A \xrightarrow{1_A} A$ with $A$, we are thus identifying the count noun $A$ with its ``global" aspect $A \text{ qua } A$. \par 
	In modeling consideration of the various aspects of people, we will be interested in a particular subcategory of $(\textbf{CN}^{\rightarrow})^{op}$, namely the co-slice category of objects under the global object. For instance, the category 
	\begin{center} 
		\begin{tikzcd}
			\boxed{\text{a scf}} \downarrow \textbf{CN}
		\end{tikzcd}
	\end{center} 
	of objects under `a student, a coworker, a family man'---where this is identified with the global aspect $\boxed{\text{ a scf } } \xrightarrow{qua} \boxed{\text{ a scf }}$---forms a subcategory $\textbf{A}$ of the \textit{qua} category $(\textbf{CN}^{\rightarrow})^{op}$. \par 
	For concreteness, suppose we have 
	\begin{center}
		\begin{tikzcd}
			& \boxed{\text{a scf}}qua\boxed{\text{a scf}} \\
			\boxed{\text{a scf}}qua\boxed{\text{a s}} \arrow[ur] & \boxed{\text{a scf}}qua\boxed{\text{a c}} \arrow[u] & \boxed{\text{a scf}}qua\boxed{\text{a f}} \arrow[ul] \\
			& & \boxed{\text{a scf}}qua\boxed{\text{a p}} \arrow[u] & \boxed{\text{a scf}}qua\boxed{\text{a h}} \arrow[ul] 
		\end{tikzcd}
	\end{center}
	In this way, such an $\textbf{A}$ will thus serve as a way of representing the aspects of a person, including for instance Abe, relevant to whether or not a certain predicable holds of them, e.g., `honesty'. This is something that will be evaluated `\textit{qua} scf' (in terms of all their ``hats," via the global aspect), `\textit{qua} student', `\textit{qua} family man', etc., where this latter has two subaspects: `\textit{qua} parent' and `\textit{qua} husband'. Abbreviating these aspects, then, we could have displayed $\textbf{A}$ as 
	\begin{center}
		\begin{tikzcd}
			& G \\
			S \arrow[ur] & C \arrow[u] & F \arrow[ul] \\
			& & P \arrow[u] & H \arrow[ul] 
		\end{tikzcd}
	\end{center}
Before defining the relevant functor, let us also record the following definition of a concept we will use here and throughout this book. 
\begin{definition}
	A functor $F: \textbf{C}^{op} \rightarrow \textbf{Set}$ is a \textit{subfunctor} \index{functor!sub} of $G: \textbf{C}^{op} \rightarrow \textbf{Set}$, denoted $F \subseteq  G$, if for any $f: b \rightarrow a$ in $\textbf{C}^{op}$ there exists a commutative diagram 
	\begin{center}
		\begin{tikzcd}
			F(a) \arrow[d, hookrightarrow] \arrow[r, "F(f)"] & F(b) \arrow[d, hookrightarrow] \\
			G(a) \arrow[r, "G(f)", swap] & G(b).
		\end{tikzcd}	
	\end{center} 
\end{definition}
	\noindent 
Then, given a \textit{qua} category and $\mathscr{P}$ a set of predicables that are applicable to the count nouns of $\textbf{CN}$---where ``predicables" may be thought of for now as just involving grammatical expressions consisting of adjectives, verbs, or adjectival and verb phrases, including expressions such as ``mortal" or ``honest" but also ``to be a person," where this is derived from or sorted by a count noun---we define an \textit{interpretation} of $(\textbf{Qua}, \mathscr{P})$ to be a functor 
		\begin{equation}
		X: \textbf{Qua}^{op} \rightarrow \textbf{Set}
		\end{equation}
		together with a set 
		\begin{equation}
		\{X_{\phi} \hookrightarrow X | \phi \in \mathscr{P}\}
		\end{equation}
		of subfunctors of $X$ that satisfy the following conditions: 
		\begin{enumerate}
			\item $X(\boxed{\text{ A}}qua \boxed{\text{ B}}) = X(\boxed{\text{ A}}qua \boxed{\text{ A}})$; and 
			\item $X_{\phi}(\boxed{\text{ A }}qua \boxed{\text{ B }}) = X_{\phi}(\boxed{\text{ B}}qua \boxed{\text{ B}}) \circ X(\boxed{\text{ B }}qua\boxed{\text{ B }} \rightarrow \boxed{\text{ A }}qua\boxed{\text{ A }})$
		\end{enumerate} \noindent 
	Notice that since $\textbf{Qua} = (\textbf{CN}^{\rightarrow})^{op}$, its dual, $\textbf{Qua}^{op}$ is just $\textbf{CN}^{\rightarrow}$, making $X$ equivalently expressible as a functor 
	\begin{equation}
	X: \textbf{CN}^{\rightarrow} \rightarrow \textbf{Set}. 
	\end{equation}
	The possibility of comparing the interpretations of count nouns forms the basis of the possibility of comparing the corresponding interpretations of predicables deemed functorial. The functoriality here can be understood as saying, for instance, if Abe is honest \textit{qua} family man, then he is honest \textit{qua} parent. Moreover, $\boxed{\text{ a p }}$ is interpreted as the set of parents, $\boxed{\text{ a scf }}$ as the set of students who are also coworkers and family men. We will return to this example throughout the book. 
\end{example}
The next example is immensely important for the general theory that will be developed in coming sections. 
\begin{example} 
	\label{example: hom-functor}
	Let \textbf{C} be an arbitrary category, and fix an object $a$ of \textbf{C}. Then we can form the (covariant)\index{functor!Hom} \textit{Hom-functor} $\text{Hom}_{\textbf{C}}(a,\--): \textbf{C} \rightarrow \textbf{Set}$, which takes each object $b$ of \textbf{C} to the set $\text{Hom}_{\textbf{C}}(a,b)$ of $\textbf{C}$-morphisms from $a$ to $b$, and takes each \textbf{C}-morphism $f: b \rightarrow c$ to the following map between hom-sets: \begin{equation} \text{Hom}_{\textbf{C}}(a,f): \text{Hom}_{\textbf{C}}(a,b) \rightarrow \text{Hom}_{\textbf{C}}(a,c),
	\end{equation} 
	which outputs $f \circ g: a \rightarrow c$ for input $g: a \rightarrow b$. In other words, the action on morphisms is given by \textit{post-composition}. This hom-functor will be defined for any object whenever the hom-sets of $\textbf{C}$ are \textit{small}.\footnote{The other way of saying this is `whenever \textbf{C} is \textit{locally small}.'} Intuitively, the set $\text{Hom}_{\textbf{C}}(a,b)$ can be thought of as the set of ways to pass from $a$ to $b$ within \textbf{C}, or the set of ways $a$ ``sees" $b$ within the context or framework of \textbf{C}. Then, refraining from ``filling in" the object $b$, it should be obvious how $\text{Hom}(a,\--)$ can be thought of as representing in a rather general fashion `where and how $a$ goes elsewhere' or `how $a$ sees its world'. Given an object $a \in \textbf{C}$, we say that the covariant functor $\text{Hom}(a, \--)$ is \textit{represented by} $a$; for reasons we will see below, this functor is also denoted $Y^a$ (or sometimes $h^a$). It will turn out to be an important observation that instead of restricting ourselves to the hom-functor on a given $a$, we can assign to \textit{each} object $c \in \textbf{C}$ its hom-functor $\text{Hom}(c, \--)$, and then collect all these together.\par 
	Dually, we can also form the \textit{contravariant Hom-functor} $\text{Hom}_{\textbf{C}}(\--,a): \textbf{C}^{op} \rightarrow \textbf{Set}$, for a fixed object $a$ of \textbf{C}, which takes each object $b$ of \textbf{C} to the set $\text{Hom}_{\textbf{C}}(b,a)$ of $\textbf{C}$-arrows from $b$ to $a$, and takes each \textbf{C}-arrow $f: b \rightarrow c$ to $\text{Hom}_{\textbf{C}}(f,a): \text{Hom}_{\textbf{C}}(c,a) \rightarrow \text{Hom}_{\textbf{C}}(b,a)$, i.e., outputting $g \circ f: b \rightarrow a$ for input $g: c \rightarrow a$, acting by \textit{pre-composition}. This functor can be thought of as representing `how $a$ is seen by its world'. Given an object $a \in \textbf{C}$, we say that the contravariant functor $Y_a := \text{Hom}(\--,a)$ is \textit{represented by} $a$. As above, instead of restricting ourselves to the hom-functor on $a$, we can ultimately let this functor vary over all the objects of \textbf{C}.
\end{example}
Throughout this book, we will see many more examples of functors. For now, though, we can continue to develop the main concepts that will let us ``ascend" once more in generality, regarding functors (presheaves) as objects in a category, with morphisms given by certain \textit{transformations} between the functors.  
\subsection{Natural Transformations}
Functors are important for many reasons. In particular, as we will see below, special ``universal" properties are given in terms of functors. Moreover, it is possible to use two functors to do a variety of important things, such as produce a new category from old categories.\footnote{Via the \textit{comma category} construction.} However, perhaps most important for our purposes is the fact that functors can be composed, and there is a nice notion of comparing functors.\par  
There may exist a variety of ways of embedding or modeling or instantiating one category within another, i.e., there may exist many functors from one category to another. Sometimes these will be equivalent, but sometimes not. Moreover, the same blueprint may be realized in different ways, i.e., there can be different functors that act the same way on objects. \textit{Natural transformations} enable us to compare these realizations.\index{natural transformation} If functors allow us to systematically import or transform objects from one category into another and thus translate between different categories, natural transformations allow us to compare the different translations in a controlled manner.  
\begin{definition} \label{defnattrans}
	Given categories \textbf{C} and \textbf{D} and functors $F, G: \textbf{C} \rightarrow \textbf{D}$, a \textit{natural transformation}\index{natural transformation!defined} $\alpha: F \Rightarrow G$, depicted in terms of its boundary data by the ``globular" diagram 
	\par 
	\begin{center}
	\begin{tikzcd}[column sep=large]
			\textbf{C}
			\arrow[bend left=50]{r}[name=U,label=above:$F$]{}
			\arrow[bend right=50]{r}[name=D,label=below:$G$]{} &
			\textbf{D}
			\arrow[shorten <=8pt,shorten >=8pt,Rightarrow,to path={(U) -- node[label=right:$\alpha$] {} (D)}]{}
	\end{tikzcd} 
	\end{center}  
	\par \noindent 
	consists of the following: 
	\begin{itemize}
		\item for each object $c \in \textbf{C}$, an arrow $\alpha_c: F(c) \rightarrow G(c)$ in \textbf{D}, called the $c$-component of $\alpha$, the collection of which (for all objects in \textbf{C}) define the \textit{components} of the natural transformation; and
		\item for each morphism $f: c \rightarrow c'$ in \textbf{C}, the following square of morphisms, called the \textit{naturality square} for $f$, must commute in \textbf{D}: \par 
		\centering 
		\begin{tikzcd}
				F(c) \arrow[r, "{\alpha_c}"] \arrow[d, "F(f)", swap]
				& G(c) \arrow[d, "G(f)"] \\
				F(c') \arrow[r, "{\alpha_{c'}}", swap]
				& G(c')
		\end{tikzcd}
	\end{itemize}
	The set of natural transformations $F \rightarrow G$ is sometimes denoted $Nat(F,G)$. 
\end{definition}
Composition of natural transformations is a little more complicated than the ``usual" composition, for there are in fact two types of composition: vertical and horizontal: \par 
	\begin{tikzcd}[column sep=huge]
	\textbf{C}
	\arrow[bend left=70]{r}[name=U,label={[yshift=-0.2cm] above:\footnotesize $F$}]{}
	\arrow{r} [name = M, label={[yshift=-0.27cm]: \footnotesize $G$}]{}
	\arrow[bend right=70]{r}[name=D,label=below:\footnotesize $H$]{} &
	\textbf{D}
	\arrow[shorten >= 4.5pt, shorten <= 4.5pt, Rightarrow,to path={(U) -- node[label=right:\footnotesize $\alpha$] {} (M)}]{}
	\arrow[shorten <=4.5pt, shorten >= 4.5pt, Rightarrow,to path={(M) -- node[label=right: \footnotesize $\beta$] {} (D)}]{}
	& & 
		\textbf{C}
	\arrow[bend left=50]{r}[name=Q,label={[yshift=-0.2cm]above:\footnotesize $F_1$}]{}
	\arrow[bend right=50]{r}[name=R,label=below:\footnotesize $G_1$]{} &
	\textbf{D} 
	\arrow[bend left=50]{r}[name=S,label={[yshift=-0.2cm]above:\footnotesize $F_2$}]{}
	\arrow[bend right=50]{r}[name=T,label=below:\footnotesize $G_2$]{}
	\arrow[shorten <=9pt,shorten >=7pt,Rightarrow,to path={(Q) -- node[label=right:$\alpha$] {} (R)}]{} & 
	\textbf{E} 
	\arrow[shorten <=9pt,shorten >=7pt,Rightarrow,to path={(S) -- node[label=right:$\beta$] {} (T)}]{}
\end{tikzcd} \par \noindent 
Vertical composition uses the symbol `$\circ$', giving $\beta \circ \alpha: F \Rightarrow H$ for the diagram on the left. Componentwise, this is defined by $(\beta \circ \alpha)_c := \beta_c \circ \alpha_c$. Horizontal composition uses the symbol `$\diamond$', giving $\beta \diamond \alpha: F_2 \circ F_1 \Rightarrow G_2 \circ G_1$ on the right, whose component at $c \in \textbf{C}$ is defined as the composite of the following commutative square: 
\par 
\begin{center} 
\begin{tikzcd}
	F_2 F_1 (c) \arrow[r, "{\beta_{F_1 c}}"] \arrow[d, "{F_2 (\alpha_c )}", swap] \arrow[dr, dashed, anchor=center, "{(\beta \Diamond \alpha)_c}"] & G_2 F_1 (c) \arrow[d, "{G_2 (\alpha_c )}"] \\ 
	F_2 G_1 (c) \arrow[r, "{\beta_{G_1 c}}", swap] & G_2 G_1 (c)
\end{tikzcd}   
\end{center}
\begin{example}
An \textit{endofunctor}\index{functor!endo} is a functor whose domain is equal to its codomain, i.e., a functor from a category to itself. Using this notion, given a category $\textbf{C}$, we can form the category $\textbf{End}(\textbf{C})$ that\index{category!of endofunctors} has 
\begin{itemize}
	\item for objects: the endofunctors $F: \textbf{C} \rightarrow \textbf{C}$; and 
	\item for morphisms: the \textit{natural transformations} between such endofunctors. 
\end{itemize}
\end{example}
\begin{example}
	For \textbf{J} an arbitrary category viewed as a template or indexing category \index{category!indexing} for \textbf{C}, we can produce the category $\textbf{C}^{\textbf{J}}$ of $\textbf{J}$-diagrams in \textbf{C}, where each object is a functor $F: \textbf{J} \rightarrow \textbf{C}$, and for two such objects $F, G$, a morphism of $\textbf{C}^{\textbf{J}}$ from $F$ to $G$ is a natural transformation between the functors. Moreover, given three such functors and two such natural transformations, we can form the composite natural transformation (either vertical or horizontal), and such composition is associative.
\end{example}
As such examples suggest, via the notion of natural transformation, we can more generally form the \textit{functor category},\index{category!functor} for natural transformations can be thought of as morphisms between functors. This is defined on categories \textbf{C} and \textbf{D} as having for objects all the functors from \textbf{C} to \textbf{D} and for morphisms all the natural transformations between such functors. There are clearly identity natural transformations and a well-defined composition formula (see Definition \ref{defnattrans} for the definition of vertical composition) for the natural transformations, and the category laws hold more generally, so we indeed have defined a category: the category of functors, denoted $Fun(\textbf{C}, \textbf{D})$, or more commonly, $\textbf{D}^{\textbf{C}}$. For our purposes, the most important thing to note here is that since presheaves are just (contravariant) functors, and a morphism of presheaves from $F$ and $G$ is just a natural transformation $\alpha: F \Rightarrow G$, i.e., a family of functions (one for each object of the domain category) subject to the naturality square commuting condition, we can form the presheaf functor category.  
\begin{definition}
	The \textit{presheaf category},\index{category!presheaf} denoted $\textbf{Set}^{\textbf{C}^{op}}$, is the (contravariant) functor category having for objects all functors $F: \textbf{C}^{op} \rightarrow \textbf{Set}$, and for morphisms $F \rightarrow G$ all natural transformations $\theta: F \Rightarrow G$ between such functors. Such a $\theta$ assigns to each object $c$ of \textbf{C} a function $\theta_c: F(c) \rightarrow G(c)$, and does so in such as way as to make all diagrams 
	\begin{center} 
		\begin{tikzcd}
			F(c) \arrow[r, "{\theta}_c"] \arrow[d, "F(f)", swap]
			& G(c) \arrow[d, "G(f)"] \\
			F(d) \arrow[r, "{\theta}_d", swap]
			& G(d)
		\end{tikzcd}
	\end{center} 
	commute for $f: d \rightarrow c$ in \textbf{C}.
\end{definition}
Because of the importance that this category will have in our story, we will see many more examples of natural transformations, in a variety of contexts. For now, here is an important example of natural transformations that has the additional benefit of introducing a number of other pivotal concepts in category theory, in particular that of \textit{limit} and \textit{colimit}.    
\begin{example}
	Recall the discussion of \textbf{J}-shaped diagrams $D: \textbf{J} \rightarrow \textbf{C}$. Using this construction and natural transformations, we can introduce the concepts of cones and cocones of a diagram, thereby characterizing the limit\index{limit} and colimit\index{colimit} of a diagram as the universal such (co)cone. \par 
	In a category \textbf{C}, a \textit{terminal object}\index{terminal object} is a special object, usually denoted 1 (owing to the fact that in \textbf{Set}, it is just a 1-element set), satisfying a certain universal property: 
	\begin{quote}
		for every object $x$ of $\textbf{C}$, there exists a unique morphism $!: x \rightarrow 1$. 
	\end{quote} 
If such a terminal object exists, it will be unique (up to unique isomorphism).\footnote{Dually, an \textit{initial object}\index{initial object} in a category $\textbf{C}$ is an object $\emptyset$ such that for any object $x$ of $\textbf{C}$, there is a unique morphism $!: \emptyset \rightarrow x$. Similarly, an initial object, if it exists, will be unique up to unique isomorphism, letting us speak of \textit{the} initial object. Note that an initial object in $\textbf{C}$ is the same as a terminal object in $\textbf{C}^{op}$.} \par 
But $\textbf{Cat}$ is a category, and we thus speak of the terminal object in $\textbf{Cat}$ as the \textit{terminal category}. \index{category!terminal} This is just the unique (up to isomorphism) category with a single object and a single morphism (necessarily the identity morphism on that object). We denote this $\textbf{1}$ (or sometimes $\underline{1}$). \par   
	Let $t: \textbf{J} \rightarrow \underline{1}$ denote the unique functor to the terminal category. Suppose given an object $c \in Ob(\textbf{C})$, which is represented by the functor $c: \underline{1} \rightarrow \textbf{C}$. Then, precomposing with $t$ to get $c \circ t: \textbf{J} \rightarrow \textbf{C}$ will just give us the \textit{constant functor}\index{functor!constant} at $c$, where this sends each object in \textbf{J} to the same \textbf{C}-object $c$ and every morphism in \textbf{J} to the identity $id_c$ on that object. Thus, composing with $t$ induces a functor $\textbf{C} \cong Fun(\underline{1}, \textbf{C}) \rightarrow Fun(\textbf{J}, \textbf{C})$, which is commonly denoted $\Delta_t: \textbf{C} \rightarrow Fun(\textbf{J}, \textbf{C}) = \textbf{C}^{\textbf{J}}$. Thus, we have an embedding $\Delta: \textbf{C} \rightarrow \textbf{C}^{\textbf{J}}$ that takes an object $c$ to the constant functor at $c$ and a morphism $f: c \rightarrow c'$ to the \textit{constant natural transformation},\index{natural transformation!constant} where each component is defined to be the morphism $f$. One can observe that each arrow $f: c \rightarrow c'$ in \textbf{C} induces a natural transformation $\Delta(c) \xrightarrow{\Delta(f)} \Delta(c')$ such that 
	\begin{center} 
		\begin{tikzcd}
			(\Delta c)(i) \arrow[r, "{\Delta(f)_i}"] \arrow[d, "{(\Delta c)(e)}", swap] & (\Delta c')(i) \arrow[d, "{(\Delta c')(e)}"] & & i \arrow[d, "e"] \\ 
			(\Delta c)(j) \arrow[r, "{\Delta (f)_j}", swap] & (\Delta c')(j) & & j  
		\end{tikzcd}   
	\end{center} \par \noindent 
	commutes for each edge $e$ of the indexing category $\textbf{J}$. But recall that the constant functor just sends every object $c \in \textbf{C}$ to $c$ and assigns $id_c$ to each edge, so the previous diagram reduces to 
	\begin{center} 
		\begin{tikzcd}
			c \arrow[r, "{\Delta(f)_i}"] \arrow[d, "{id_c}", swap] & c' \arrow[d, "{id_{c'}}"] \\ 
			c \arrow[r, "{\Delta (f)_j}", swap] & c'  
		\end{tikzcd}   
	\end{center} \par \noindent 
	which obviously commutes. \par 
	If we now consider, for an arbitrary $\textbf{J}$-diagram $F: \textbf{J} \rightarrow \textbf{C}$ and for $c \in \textbf{C}$, the arrows (which are in fact natural transformations)
	\begin{equation*}
	\Delta c \longrightarrow F \hspace*{4em} F \longrightarrow \Delta c,  
	\end{equation*}
	we get that a typical arrow in $\textbf{C}^{\textbf{J}}$ corresponding to these arrows is just a natural transformation, i.e., a family of arrows of $\textbf{C}$, 
	\begin{equation*}
	(\Delta c)(i) \xrightarrow{\xi(i)} F(i) \hspace*{4em} F(i) \xrightarrow{\xi (i)} (\Delta c)(i)  
	\end{equation*}
	indexed by the various objects or nodes of $\textbf{J}$ and such that 
	\begin{center} 
		\begin{tikzcd}
			(\Delta c)(i) \arrow[r, "{\xi(i)}"] \arrow[d, "{(\Delta c)(e)}", swap] & F(i) \arrow[d, "{F(e)}"] & i \arrow[d, "e"] & & F(i) \arrow[d, "{F(e)}", swap]  \arrow[r, "{\xi(i)}"] & (\Delta c)(i) \arrow[d, "{(\Delta c)(e)}"] \\ 
			(\Delta c)(j) \arrow[r, "{\xi(j)}", swap] & F(j) & j & & F(j) \arrow[r, "{\xi(j)}", swap] & (\Delta c)(j)  
		\end{tikzcd}   
	\end{center} \par \noindent 
	commute for each such edge $e: i \rightarrow j$ in \textbf{J}. But when we apply the functor $\Delta$, these commutative squares collapse to the commutative triangles
	\begin{center} 
		\begin{tikzcd}
			& F(i) \arrow[dd, "{F(e)}"]  & & i \arrow[dd, "e"] & & F(i) \arrow[dd, "{F(e)}", swap]  \arrow[dr, "{\xi(i)}"] \\ 
			c \arrow[ur, "{\xi(i)}"] \arrow[dr, "{\xi (j)}", swap] & & & & & & c \\ 
			& F(j) & & j & & F(j) \arrow[ur, "{\xi(j)}", swap]  
		\end{tikzcd}   
	\end{center} \par \noindent 
	The definitions guarantee that whenever the indexing category has composable edges, the corresponding composite triangles commute. The natural transformations represented by the triangles on the left give a \textit{left solution} for the diagram in $\textbf{C}$, sometimes also called a \textit{cone over}\index{cone} the diagram $F$ with \textit{summit} vertex $c$. The natural transformations represented by the triangles on the right give a \textit{right solution} for the diagram, also called a \textit{cocone for} (or \textit{cone under}) the diagram $F$ with \textit{nadir} $c$.\footnote{Hopefully this is already clear, but in case not: the terminology of `over' and `under' has to do with the fact that the above triangles can be presented as rotated clockwise 90 degrees.} We can then form the category of cones, where an object in the category of cones over $F$ will be a cone over $F$, with some summit, while a morphism from a cone $\xi: c \Rightarrow F$ to a cone $\mu: d \Rightarrow F$ is a morphism $f: c \rightarrow d$ in \textbf{C} such that for each index $j \in \textbf{J}$, $\mu_j \circ f = \xi_j$, i.e., a map between the summits such that each leg of the domain cone factors through the corresponding leg of the codomain cone. \par Using these notions, we can define the \textit{limit}\index{limit} of $F$ in terms of a universal cone, where a cone $\alpha: L \rightarrow F$ with vertex $L$ is universal with respect to $F$ provided for every cone $\Delta c \rightarrow F$, there is a unique map $g: \Delta c \rightarrow F$ making 
	\begin{center} 
		\begin{tikzcd}[row sep = large]
			c \arrow[dr, "{\xi(i)}"] \arrow[ddr, "{\xi (j)}", swap, bend right = 20] \arrow[rr, dashed, "g"] & & L \arrow[dl, "\alpha(i)", swap] \arrow[ddl, "\alpha(j)", bend left =20] \\ 
			& F(i) \arrow[d, "{F(e)}", swap]  & & i \arrow[d, "e"]  \\ 
			& F(j) & & j
		\end{tikzcd}   
	\end{center} 
	commute. In such a case, one usually refers (somewhat improperly) to the universal cone by just the vertex $L = \varprojlim F$, and calls this the limit of $F$.
	 \par 
	We can see a limit for a diagram $F: \textbf{J} \rightarrow \textbf{C}$ as a representation for the corresponding functor $Cone(\--, F): \textbf{C}^{op} \rightarrow \textbf{Set}$, sending $c \in \textbf{C}$ to the set of cones over $F$ with summit $c$.\footnote{While we could consider limits and colimits in any category, by something called the Yoneda lemma (on which much more below) we can be assured that the constructions of (co)limits of diagrams valued in the category $\textbf{Set}$ suffice to provide formulae for (co)limits in any category. To ensure that we have a \textit{set} of cones, we need only assume that the diagram is indexed by a small category $\textbf{J}$ and that $\textbf{C}$ is locally small, thereby guaranteeing that the functor category $\textbf{C}^{\textbf{J}}$ is locally small. Under certain conditions, we can weaken these restrictions.} The limiting cone will be \textit{universal} in the sense that for any other cone over $F$, there will exist a unique arrow from the summit of that cone to the summit of the limiting cone, i.e., it must pass uniquely through the limiting cone if it wants to pass down to $F$. 
	\par 
	The dual construction produces a category of cocones $\textbf{CoCones}(F)$, wherein the universal cocone emerges as the \textit{colimit}\index{colimit} of the diagram $F$, denoted $\text{colim } F$ (or sometimes $\varinjlim F$), as a representation for $Cone(F, \--)$, forcing all cocones to receive maps from the colimiting cone if they want to receive maps from $F$. We will have more to say about this in a moment. \par 
	First, let us consider an alternative way to understand this construction. To understand this, we first need another concept, one that we will discuss via a few examples, before coming back to what it has to do with limits (and colimits). 
	\begin{definition}
		Let \textbf{C} be a category, and let $F: \textbf{C} \rightarrow \textbf{Set}$ be a (covariant) functor. Then the \textit{category of elements of} $F$,\index{category!of elements} denoted $\int_{\textbf{C}} F$ (or just $\int F$ if the context is clear), is defined: 
	\begin{equation}
		\begin{split}
		Ob(\int F) = \{(c, x) \hspace*{0.3em}| \hspace*{0.3em} c \in \textbf{C}, x \in F(c) \} \\
		Hom_{\int F}((c,x), (c', x')) = \{f: c \rightarrow c' \hspace*{0.3em}| \hspace*{0.3em} F(f)(x) = x'\}.
		\end{split}
		\end{equation}
		Similarly for the contravariant case: for $F: \textbf{C}^{op} \rightarrow \textbf{Set}$ the \textit{category of elements of} $F$, denoted $\int_{\textbf{C}^{op}} F$ (or just $\int F$), is defined: 
		\begin{equation}
		\begin{split}
		Ob(\int F) = \{(c, x) \hspace*{0.3em}| \hspace*{0.3em} c \in \textbf{C}, x \in F(c) \} \\
		Hom_{\int F}((c,x), (c', x')) = \{f: c \rightarrow c' \hspace*{0.3em} | \hspace*{0.3em} F(f)(x') = x\}.
		\end{split}
		\end{equation}
		Associated to these constructions are the natural functors $\pi_{F}: \int F \rightarrow \textbf{C}$, called the projection functors, sending each object $(c, x) \in Ob(\int F)$ to the object $c \in Ob(\textbf{C})$ or $Ob(\textbf{C}^{op})$, and each morphism $f: (c, x) \rightarrow (c', x')$ to the morphism $f: c \rightarrow c'$, i.e., $\pi(f, (c, x), (c', x')) = f$. 
	\end{definition} \noindent 
	As a concrete instance of this, recall the ``vertex coloring" functor $nColor$. \index{functor!nColoring} An object in the category of elements $\int nColor$ of this functor $nColor$ will be a graph together with a chosen $n$-coloring, i.e., objects are $n$-colored graphs. A morphism $\phi: G \rightarrow G'$ between a pair of $n$-colored graphs will be a graph homomorphism $\phi: G \rightarrow G'$ so that the induced function $nColor(\phi): nColor(G') \rightarrow nColor(G)$ takes the chosen coloring of $G'$ to the chosen coloring of $G$, i.e., the graph homomorphism $\phi$ will preserve the chosen colorings in the sense that each red vertex of $G$ will be carried to a red vertex of $G'$. In short, then, $\int nColor$ is the category of $n$-colored graphs and the color-preserving graph homomorphisms between them. \par 
	For another example, recall the hom-functors, first introduced in \ref{example: hom-functor}. Objects in the category of elements of $\text{Hom}_{\textbf{C}}(c, \--)$ are the morphisms $f: c \rightarrow d$ in \textbf{C}. A morphism from $f: c \rightarrow d$ to $g: c \rightarrow e$ is then a morphism $h: d \rightarrow e$ such that $g = h \circ f$. $h$ is said to be a morphism under $c$ because of the diagram attached to this condition: 
	\begin{center}  
		\begin{tikzcd}
			& c \arrow[dl, "f", swap] \arrow[dr, "g"] \\
			d \arrow[rr, "h", swap] && e
		\end{tikzcd} 
	\end{center} 
	This category is none other than the \textit{co-slice category} of objects \textit{under} the $c \in \textbf{C}$.\index{category!co-slice} Note that the forgetful functor $U: c/\textbf{C} \rightarrow \textbf{C}$ sends an object $f: c \rightarrow d$ to the codomain, and takes a morphism (a commutative triangle) to the arrow opposite the object $c$, i.e., to $h$ in the above instance. We could also construct the dual category of elements $\int \text{Hom}_{\textbf{C}}(\--, c)$ in terms of the \textit{slice category} $\textbf{C}/c$ \textit{over} the object $c \in \textbf{C}$.\footnote{In this connection, we can mention the important result that for $\textbf{C}$ small and $P$ a presheaf on \textbf{C}, one can show an equivalence of categories 
		\begin{equation*}
		\textbf{Set}^{\textbf{C}^{op}}/P \cong \textbf{Set}^{(\int_{\textbf{C}} P)^{op}} . 
		\end{equation*} } \par 
	The category of elements is extremely significant because any universal property can be seen as defining an initial or terminal object in this category. In particular, it turns out that for any small functor\footnote{A functor or diagram is small if its indexing category is small. \index{functor!small}} $F: \textbf{C} \rightarrow \textbf{Set}$, we have 
	\begin{equation}
	\text{colim } F \cong \Pi_0 (\int F),
	\end{equation}
	where $\Pi_0$ operates by picking out the connected components\index{functor!connected components} and $\int F$ is the category of elements\index{category!of elements} of $F$. So, in other words, the set of connected components of the category of elements of a functor $F$, $\Pi_0 (\int F)$, is isomorphic to the colimit\index{colimit} of $F$.\footnote{Alternatively, in general, we could just have said that a colimit is an initial object\index{initial object} in the category $\int Cone(F, \--)$, and we note that the forgetful functor $\int Cone(F, \--) \rightarrow \textbf{C}$ will take a cone to its nadir.} \par 
	To see this in action, recall the functor (diagram) 
	\begin{center} 
		\begin{tikzpicture}[framed, scale=0.7]
		\tikzset{vertex/.style = {shape=circle,draw, fill=black, minimum size=3pt, inner sep =0pt}}
		\tikzset{edge/.style = {->,> = latex'}}
		\node[vertex] (a) [label=left:{\footnotesize $(a,1)$}] at  (0,2) {};
		\node[vertex] (a2) [label=left:{\footnotesize$(a,2)$}] at  (0,1) {};
		\node[vertex] (d) [label=right:{\footnotesize$(d,1)$}] at  (6.5,1.2) {};
		\node[vertex] (d2) [label=right:{\footnotesize$(d,2)$}] at  (6.5,0) {};
		\node[vertex,gray] (d3) [label=right:{\footnotesize $(d,3)$}] at  (6.5,-1) {};
		\node[vertex] (d4) [label=right:{\footnotesize$(d,4)$}] at  (6.5,-2) {};
		
		\node[vertex] (b) [label=above:{\footnotesize $(b,1)$}] at (2,1) {};
		\node[vertex] (b2) [label=below:{\footnotesize $(b,2)$}] at (2,-0.1) {};
		\node[vertex] (c) [label=above:{\footnotesize $(c,1)$}] at (4,2) {};
		\node[vertex] (c2) [label=below:{\footnotesize $(c,2)$}] at (4,1) {};
		\node[vertex] (c3) [label=below:{\footnotesize $(c,3)$}] at (4,-0.1) {};
		\node[vertex] (q) [label=left:{\footnotesize $(q,1)$}] at (6.5,-3) {};
		\node[vertex] (q2) [label=left:{\footnotesize $(q,2)$}] at (6.5,-4) {};
		\node[vertex] (q3) [label=left:{\footnotesize $(q,3)$}] at (6.5,-5) {};
		\node[vertex] (r) [label=right:{\footnotesize $(r,1)$}] at (8.5,-3) {};
		\node[vertex] (r2) [label=right:{\footnotesize $(r,2)$}] at (8.5,-4.5) {};
		
		\path[-latex, dashed] (q) edge node[auto] {} (r2);
		\path[-latex, thick] (q2) edge node[auto] {} (r);
		\path[-latex, thick] (q3) edge node[auto] {} (r);
		
		\path[-latex,ultra thick] (a) edge node[pos=0.5, above] {} (b);
		\path[-latex] (a2) edge node[pos=0.5, above] {} (b2);
		\path[-latex,ultra thick] (c) edge node[above] {} (b);
		\path[-latex] (c2) edge node[above] {} (b2);
		\path[-latex,ultra thick] (c3) edge node[above] {} (b);
		\path[-latex,ultra thick] (c) edge node[pos=0.5, above] {} (d);
		\path[-latex] (c2) edge node[pos=0.5, above] {} (d2);
		\path[-latex,ultra thick] (c3) edge node[pos=0.5, above] {} (d4);
		\end{tikzpicture} 
	\end{center} 
from \ref{example: diagram 1}.
The thicknesses and colorings in this picture can now be explained. The picture above is in fact a representation of the category of elements of $F$. The various thicknesses depict the action of taking its connected components. By inspection, one can verify that this is just the set 
\begin{equation*} \{[(a,1)], [(a,2)], [(d,3)], [(q,1)], [(q,2)]\},
\end{equation*} 
where each element is a representative of one of the components, which is in turn isomorphic to a set of cardinality $5$, entailing that $colim_{\textbf{J}} F \cong$ a set with $5$ elements. \par 
As for limits:\index{limit} we could also show that the \textit{limit} of any small functor $F: \textbf{C} \rightarrow \textbf{Set}$ is isomorphic to the set of functors $\textbf{C} \rightarrow \int F$ that define a \textit{section} to the canonical projection $\pi: \int F \rightarrow \textbf{C}$.\footnote{We will discuss sections later on.} Alternatively, we can define the limit as a terminal object in the category of elements of cones over $F$, i.e., in $\int Cone(\--,F)$. Note that the forgetful functor $\int Cone(\--, F) \rightarrow \textbf{C}$ will send a given cone to its summit. \par 
 In short, both the limiting and the colimiting cones are universal in the sense of acting as a kind of gateway through which all other cones must pass. In any particular case, such universal objects need not exist. However, we noted above that a diagram is said to be \textit{small} if its indexing category is a small category, which further allows us to define a category $\textbf{C}$ as \textit{complete} \index{category!complete} if it admits limits of all small diagrams valued in \textbf{C}, and as \textit{cocomplete} \index{category!cocomplete} if it admits all colimits of all small diagrams valued in \textbf{C}.      
\end{example}
\section{Yoneda: The Most Important Idea in Category Theory}
We are now in a position to consider what is perhaps the most important idea in category theory, the Yoneda results.\index{Yoneda} But in the coming sections, we will motivate this idea through a simplified special case, its analogue for posets (in fancier language, its ``$\textbf{2}$-enriched" analogue). This motivation requires that one first understand \textit{enrichment}, the introduction of which also gives us a chance to refine our understanding of categories in general.  
\subsection{First, Enrichment!}
Not all categories were created equal. For instance, in certain categories, there is a natural way of combining elements of the category, i.e., making use of an operation that takes two elements and ``adds" or ``multiplies" them together. Not all categories admit such a thing. Those that do are called \textit{symmetric monoidal}. 
\begin{definition} 
	A \textit{symmetric monoidal structure} on a category $\mathcal{V}$ consists of the following data: 
	\begin{enumerate}
		\item a bifunctor $\-- \otimes \--: \mathcal{V} \times \mathcal{V} \rightarrow \mathcal{V}$, called the \textit{monoidal product}; 
		\item a unit object $I \in \text{Ob}(\mathcal{V})$, called the \textit{monoidal unit}, 
	\end{enumerate} 
	subject to the following specified natural isomorphisms: 
	\begin{equation}
	v \otimes w \cong_{\gamma} w \otimes v \hspace*{2em} u \otimes (v \otimes w) \cong_{\alpha} (u \otimes v) \otimes w \hspace*{2em} I \otimes v \cong_{\lambda} v \cong_{\rho} v \otimes I
	\end{equation}
	that witness symmetry, associativity, and unit conditions on the monoidal product. There are then standard ``coherence conditions" that these natural transformations are expected to obey. \par 
	A category equipped with such a symmetric monoidal structure is then called a \textit{symmetric monoidal category},\index{category!symmetric monoidal} denoted, e.g., $(\mathcal{V}, \otimes, I)$. A \textit{monoidal category} is similarly defined, except the symmetry natural isomorphism displayed on the left above is left out. If the natural isomorphisms involving associativity and the unit are replaced by \textit{equalities}, then the monoidal structure is said to be \textit{strict}.
	\end{definition}  
This is defined on categories in general, but an especially simple special case comes from restricting the definition to preorders (as categories). 
\begin{definition}
	A \textit{symmetric monoidal structure} on a preorder $(X, \leq)$ consists of
	\begin{itemize}
		\item an element $I \in X$ called the monoidal unit, and 
		\item a function $\otimes: X \times X \rightarrow X$, called the monoidal product.
	\end{itemize}
	These must further satisfy the following, for all $x_1,x_2, y_1, y_2, x, y, z \in X$, where we use infix notation, i.e., $\otimes(x_1,x_2)$ is written $x_1 \otimes x_2$: 
	\begin{itemize}
		\item \textbf{monotonicity}: if $x_1 \leq y_1$ and $x_2 \leq y_2$, then $x_1 \otimes x_2 \leq y_1 \otimes y_2$
		\item \textbf{unitality}: $I \otimes x = x$ and $x \otimes I = x$
		\item \textbf{associativity}: $(x \otimes y) \otimes z = x \otimes (y \otimes z)$
		\item \textbf{symmetry}: $x \otimes y = y \otimes x$. 
	\end{itemize}
	Then a preorder equipped with a symmetric monoidal structure, $(X, \leq, I, \otimes)$, is called a \textit{symmetric monoidal preorder}.\index{pre-order!symmetric monoidal}
\end{definition} \noindent  
Monoidal units may be, e.g., $0, 1,$ true, false, $\{*\}$, etc. Monoidal ``products" include $\otimes, +, *, \wedge, \vee, \times$, etc.
\begin{example}
	The simplest nontrivial preorder is $\textbf{2} = \{0 \xrightarrow{\leq} 1\}$. Alternatively, you might think of this as $\textbf{2} = \{false, true\}$ with the single non-trivial arrow $false \leq true$. There are two different symmetric monoidal structures on it. To consider one of these: let the monoidal unit be $true$ and the monoidal product be $\wedge$ (AND), giving a monoidal preorder $(\textbf{2}, \leq, true, \wedge)$.\footnote{\cite{fong_seven_2018} sensibly calls this $\textbf{Bool}$, but we may just stick with calling it $\textbf{2}$, after its carrier pre-order. The reader who desires a more in-depth treatment of enrichment, or who is intrigued by any of these matters, will surely enjoy the recent \cite{fong_seven_2018}. Readers with a higher tolerance for abstraction might also find \cite{kelly_basic_2005} useful.} 
\end{example}
\begin{example}
	For a set $S$, the powerset $\mathcal{P}(S)$ of all subsets of $S$, with the order $A \leq B$ given by subset relation $A \subseteq  B$, in fact has a symmetric monoidal structure on it: $(\mathcal{P}(S), \leq, S, \cap)$ is a symmetric monoidal preorder. \par 
	In particular, taking $S$ equal to a two-element set, this is isomorphic to $(A_4, \leq_k, B, \otimes)$, Belnap's\index{Belnap} four-valued ``knowledge lattice"\index{lattice!knowledge} (or ``approximation lattice") $A_4 = (\{\bot, t, f, \top\}, \leq_k)$, often used by relevance and paraconsistent logicians, where the values are the various subsets of $\{t,f\}$. Here, $\top$ (or sometimes $B$) is `\textit{both} true and false', $\otimes$ is a `consensus' connective corresponding to meet, and $(A_4, \leqslant_k)$ is the (complete) lattice corresponding to an ordering on epistemic states (`how much info/knowledge')
	\begin{center}
		\includegraphics[scale=0.4]{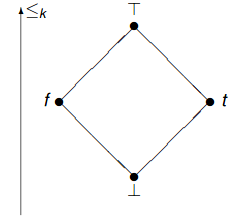}
	\end{center}
This structure has four `truth values': the classical ones ($t$ and $f$); a truth value $\bot$ that intuitively captures the ``lack of information" (``neither $t$ nor $f$"); and a truth value $\top$ that can be deployed to represent contradictions (``\textit{both} $t$ and $f$"). The underlying partial order of the lattice has $t$ and $f$ as its intermediate truth values, $\bot$ as the $\leq_k$-minimal element, and $\top$ as the $\leq_k$-maximal element. Overall, the partial order $(\{\bot, t, f, \top\}, \leq_k)$ of the lattice is often regarded as ranking the ``amount of knowledge or information," or ``approximates the information in," i.e., if $x \leq_k y$, then $y$ gives us \textit{at least as much} information as $x$ (possibly more). A move up in the lattice represents an ``increase in the amount of information," with $\otimes$ taking the uppermost element below both $x$ and $y$.\footnote{For more on this lattice, see \cite{belnap_how_1992}.} 
\end{example}
\begin{example}
	Let $[0,\infty]$ be the set of non-negative real quantities, together with $\infty$. Consider the preorder $([0, \infty], \geq)$, with the natural order $\geq$, e.g., $\pi \geq 0.8$, $14.\overline{33} \geq 11$, and of course $\infty \geq x$ for all $x \in [0, \infty]$. There is a symmetric monoidal structure here, with monoidal unit $0$ and monoidal product $+$ (where in particular $x + \infty = \infty$ for any $x \in [0, \infty]$). After \cite{fong_seven_2018}, we can call this symmetric monoidal preorder $\textbf{Cost} := ([0,\infty], \geq, 0, +)$, since we think of the elements of $[0, \infty]$ as costs. 
\end{example} 
In the standard definitions of a category that we have seen thus far, the hom-sets are \textit{sets}, i.e., objects of the category $\textbf{Set}$. On this approach, with an ordinary category $\textbf{C}$, given a \textit{set} of ``points" (objects), that the hom-sets are specifically sets effectively means that the task or question of getting from (or relating) one point to another has a \textit{set} of approaches or answers or names. But what if we generalized this story and let the hom-sets of a category come from some category other than $\textbf{Set}$?  \par  
Symmetric monoidal categories are important, in large part, because of something we can \textit{do} with them: we can \textit{enrich} an (arbitrary) category in them! What does that mean? \cite{fong_seven_2018} suggests a very nice intuitive way of thinking of this: enriching in, say, a monoidal preorder $\mathcal{V} = (V, \leq, I, \otimes)$ just means ``letting $\mathcal{V}$ structure the question of the relations or paths between the objects of the underlying category." In this general context, enriching in different monoidal categories often recovers (while generalizing) important entities in math. For instance, it emerges that categories ``enriched in \textbf{Cost}," or  $\textbf{Cost}$-categories, provide a powerful generalization of the notion of metric space.
\begin{definition}
	Let $\mathcal{V} = (V, \leq, I, \otimes)$ be a symmetric monoidal preorder. A $\mathcal{V}$-category \index{category!enriched} $\mathcal{X}$ consists of
	\begin{itemize}
		\item specification of a set $\text{Ob}(\mathcal{X})$, elements of which are objects; 
		\item for every two objects $x, y$, specification of an element $\mathcal{X}(x,y) \in V$, called the hom-object, 
	\end{itemize}
	and where these satisfy the two properties 
	\begin{itemize}
		\item for every object $x \in \text{Ob}(\mathcal{X})$, we have $I \leq \mathcal{X}(x,x)$, and 
		\item for every three objects $x, y,z \in \text{Ob}(\mathcal{X})$, we have $\mathcal{X}(x,y) \otimes \mathcal{X}(y,z) \leq \mathcal{X}(x,z)$. 
	\end{itemize}
	In this case, we call $\mathcal{V}$ the \textit{base of the enrichment} for $\mathcal{X}$, or just say $\mathcal{X}$ is \textit{enriched in} $\mathcal{V}$. 
\end{definition}
\begin{example} 
	What happens if we enrich in \textbf{Cost}$= ([0, \infty], \geq, 0, +)$? Following the definition: a \textbf{Cost}-category $\mathcal{X}$ consists of 
	\begin{itemize}
		\item (i) a collection $Ob(\mathcal{X})$, and 
		\item (ii) for every $x, y \in Ob(\mathcal{X})$ an element $\mathcal{X}(x,y) \in [0, \infty]$.
	\end{itemize}
	The idea here is that $Ob(\mathcal{X})$ provides the ``points", while $\mathcal{X}(x,y) \in [0, \infty]$ plays the role of supplying the ``distances." Still just following the definition, the properties of a category enriched in \textbf{Cost} are given by:  
	\begin{itemize}
		\item $0 \geq \mathcal{X}(x,x)$ for all $x \in Ob(\mathcal{X})$, and 
		\item $\mathcal{X}(x,y) + \mathcal{X}(y,z) \geq \mathcal{X}(x,z)$ for all $x,y,z \in Ob(\mathcal{X})$.  
	\end{itemize} 
	Note that since $\mathcal{X}(x,x) \in [0, \infty]$, the property $0 \geq \mathcal{X}(x,x)$ implies that $\mathcal{X}(x,x) = 0$. So this is in fact equivalent to the first condition $d(x,x) = 0$ describing a metric. And the second condition here is clearly the usual triangle inequality! We have thus defined, with the notion of a \textbf{Cost}-category, an extended (Lawvere) metric space.
\end{example} \noindent 
The usual definition of a metric space goes as follows:  
\begin{definition}
	A \textit{metric space} \index{metric space} ($X, \rho$) consists of $X$ a non-empty set, the elements of which are ``points," and a function $\rho: X \times X \rightarrow \mathbb{R}_{\geq 0}$ called a \textit{metric}, where this means for all $x,y,z \in X$: 
	\begin{itemize}
		\item (i) $0 \leq \rho(x,y)$ (or just $\rho(x,x) = 0$); 
		\item (ii) if $ \rho(x,y) = 0$, then $x = y$;
		\item (iii) $\rho(x,y) = \rho(y,x)$; 
		\item (iv) $ \rho(x,y) + \rho(y,z) \geq \rho(x,z)$. 
	\end{itemize}
	If we instead take a function $\rho: X \times X \rightarrow [0, \infty] = \mathbb{R}_{\geq 0} \cup \{\infty\}$, then we have an \textit{extended metric space}. 
\end{definition} \noindent 
From the categorical viewpoint, the generalized construction of a $\textbf{Cost}$-category already suggests that 
\begin{itemize}
	\item (ii) if $\rho(x,y) = 0$, then $x = y$;
	\item (iii) $\rho(x,y) = \rho(y,x)$ 
\end{itemize}
are somehow not as ``natural" as the other two conditions (triangle inequality and that points are at ``zero distance" from themselves). Indeed, there are contexts in which (ii) is not satisfied, yet we would still like to have a metric. Also, requiring (iii) or symmetry prevents us from regarding a number of constructions we would like to regard as ``metrics" as legitimate metrics, so abandoning this condition is also desirable. 
\begin{example}
	Now take the symmetric monoidal preorder $\textbf{2} = (\{false, true\}, \leq, true, \wedge)$. Enriching in \textbf{2} recovers the notion of a preorder, since for any $x,y \in \mathcal{P}$, there is either $0$ (``false") or $1$ (``true") arrow from $x$ to $y$. Accordingly, the ``homs" here will be objects of $\textbf{2}$, not \textbf{Set}. More formally, a \textbf{2}-category consists of 
	\begin{itemize}
		\item a specification of a set of objects
		\item for every object $x, y$, and element $\mathcal{X}(x,y) \in \textbf{2}$
	\end{itemize} 
	where this data satisfies 
	\begin{enumerate} 
		\item for every element $x \in Ob(\mathcal{X})$, $true \xrightarrow{\leq} \mathcal{X}(x,x)$, so $\mathcal{X}(x,x) = true$
		\item for every $x,y,z$, $\mathcal{X}(x,y) \wedge \mathcal{X}(y,z) \xrightarrow{\leq} \mathcal{X}(x,z)$ 
	\end{enumerate}
	The first condition above just amounts to reflexivity and the second to transitivity; understanding $\mathcal{X}(x,y) = true$ to just mean that $x \leq y$, clearly this just recovers the notion of a preorder. Thus, the theory of $\textbf{2}$-enriched categories just recovers precisely the theory of ordered sets and monotonic maps between them.
\end{example}
\begin{example}
	Returning to $A_4$,\index{lattice!knowledge} we can understand `t' as `told True', `f' as `told False', `N' (or $\bot$) as `told nothing (i.e., \textit{neither} told True nor told False)', `B' (or $\top$) as `\textit{both} told True and told False'. $\bot$ is at the bottom of the lattice as it gives no information at all, while $\top$ is at the top since it gives ``too much" (or inconsistent) information. \par 
	When we enrich in $A_4$, the resulting $A_4$-category $\mathcal{X}$ will describe, for any two objects $x,y$ of $\mathcal{X}$, all the (true, false, null, both true and false) information that has been received/inputed (perhaps from several independent sources) about whether or not you can ``get from" $x$ to $y$. \par 
	Enriching in $A_4$ implies that the issue of passing from $x$ to $y$ is structured by how much information/knowledge we (or some system, like a computer, prepared to receive and reason about inconsistent information) might have about the question. For instance, `I have been told that `yes' (`no') one can (cannot) pass from $x$ to $y$' or `I have been told both that you can and you cannot pass from $x$ to $y$' or `I have not been told anything about whether or not you can pass from $x$ to $y$', etc. 
\end{example}
In the next few sections, we will make use of this notion of enrichment to give a particularly simple presentation of the powerful abstract Yoneda results towards which we are leading.
\subsection{Downsets and Yoneda in the Miniature}
Given a poset $\mathcal{P} = (P, \leq)$, we saw how we can regard $\mathcal{P}$ as a category.  
\begin{definition}\label{definition: downset}
	Let $\mathcal{P}$ be a poset, and $A \subseteq  \mathcal{P}$ a subset. Then, we call the subset $A$ a \textit{downset}\index{downset} if for each $p \in A$ and $q \in \mathcal{P}$, we have that $p \in A$ and $q \leq p$ implies that $q \in A$. Dually, a subset $U \subseteq  \mathcal{P}$ is an \textit{upper set} (or \textit{up-set})\index{upper set} provided: if $p \in U$ and $p \leq q$, then $q \in U$.  \par \noindent 
	We can further define, for each element $p \in \mathcal{P}$, the down-set generated by $p$---called its \textit{principal downset},\index{downset!principal} denoted $\mathcal{D}_p$ (or sometimes just $\downarrow p$)--- as\footnote{Dually, as one would expect, we can also define, for each point $p$, its \textit{principal upper set}\index{upper set!principal} $U_p$ (or sometimes $\uparrow p$) as 
		\begin{equation*}
		\uparrow p := \{q : p \leq q \}. 
		\end{equation*}} 
	\begin{equation*}
	\downarrow p := \{q \in \mathcal{P}: q \leq p \}. 
	\end{equation*}
\end{definition} 
For instance, consider the following poset $\mathcal{P}$ on $P = \{a,b,c,d\}$ given by $a \leq c$; $b \leq c$; $b \leq d$; and the obvious identity (reflexivity) $x \leq x$ for all $x \in P$. The data of this poset is more helpfully displayed in the picture:   
\begin{center} 	
	\begin{tikzpicture}[yscale=0.75, xscale=0.5]
	\node (c) at (2.4,0) {$c$};
	\node (d) at (6.4,0) {$d$};
	\node (a) at (2.4,-2) {$a$};
	\node (b) at (6.4,-2) {$b$};
	
	\draw[->] (a) -- (c);
	\draw[->] (b) -- (c);
	\draw[->] (b) -- (d);
	\end{tikzpicture} \end{center} 
\begin{exercise}
	Is $\{a,b,c\}$ a downset? How about $\{a,b\}$? And $\{a,c,d\}$?
\end{exercise} \par \noindent 
\textit{Solution}: Yes, $\{a,b,c\}$ is a downset; same with $\{a,b\}$. But $N = \{a,c,d\}$ is not a downset, for in particular $d \in N$ and thus, considering $b \in P$, as $b \leq d$, in order for $N$ to be a downset, we should have that $b \in N$, but $b \notin N$. \par \vspace*{1em}
In general, we denote by $\mathcal{D}(\mathcal{P})$ the collection of all down-sets of the poset $\mathcal{P}$. Observe that $\mathcal{D}(\mathcal{P})$, the collection of all the downsets of $\mathcal{P}$, has a natural order on it---namely, $U \leq V$ if $U$ contained in $V$. Then, $(\mathcal{D}(\mathcal{P}), \subseteq )$ is itself an order under inclusion, one that we will denote by $\mathcal{D}(\mathcal{P})$, or $\textbf{Down}(\mathcal{P})$ when we want to emphasize that we are regarding this as a \textit{category}.\footnote{Dually, we write $\mathcal{U}(\mathcal{P})$ for the collection of all the upper sets of $\mathcal{P}$; this also has a natural order on it---namely, $U \leq V$ if $U$ contained in $V$, making $(\mathcal{U}(\mathcal{P}), \subseteq )$ an order as well.} This poset that consists of the collection of all down-sets of $\mathcal{P}$, ordered by inclusion, is sometimes called the \textit{down-set completion}. The following diagram displays the information of all the downsets of our given $\mathcal{P}$, naturally ordered by inclusion: 
\begin{center} 	
	\begin{tikzpicture}[yscale=0.85, xscale=0.73]
	\node (e) at (4.4, -3) {$\emptyset$};
	\node (c) at (2.4,0) {$\{a,b,c\}$};
	\node (d) at (6.4,0) {$\{a,b,d\}$};
	\node (t) at (4.4, 1) {$\{a,b,c,d\}$};
	\node (a) at (2.4,-2) {$\{a\}$};
	\node (b) at (6.4,-2) {$\{b\}$};
	\node (b1) at (4.4,-1) {$\{a,b\}$};
	\node (d1) at (8.4,-1) {$\{b,d\}$};
	
	\draw[->] (b1) -- (c);
	\draw[->] (a) -- (b1);
	\draw[->] (b) -- (b1);
	\draw[->] (b) -- (d1);
	\draw[->] (d1) -- (d);
	\draw[->] (b1) -- (d);
	\draw[->] (c) -- (t);
	\draw[->] (d) -- (t);
	\draw[->] (e) -- (a);
	\draw[->] (e) -- (b);
	
	\node (n) at (10.4, -1) {$=$};
	\node (e3) at (14.4, -3) {$\downarrow \emptyset$};
	\node (c3) at (12.4,0) {$\downarrow c$};
	\node (d3) at (16.4,0) {$\downarrow a \hspace*{0.2em}\cup \downarrow d$};
	\node (t3) at (14.4, 1) {$\downarrow c \hspace*{0.2em}\cup \downarrow d$};
	\node (a3) at (12.4,-2) {$\downarrow a$};
	\node (b3) at (16.4,-2) {$\downarrow b$};
	\node (b4) at (14.4,-1) {$\downarrow a \hspace*{0.2em}\cup \downarrow b$};
	\node (d4) at (18.4,-1) {$\downarrow d $};
	
	\draw[->] (b4) -- (c3);
	\draw[->] (a3) -- (b4);
	\draw[->] (b3) -- (b4);
	\draw[->] (b3) -- (d4);
	\draw[->] (d4) -- (d3);
	\draw[->] (b4) -- (d3);
	\draw[->] (c3) -- (t3);
	\draw[->] (d3) -- (t3);
	\draw[->] (e3) -- (a3);
	\draw[->] (e3) -- (b3);
	\end{tikzpicture} \end{center} 
There are a couple of valuable general observations to note at this point, which can be illustrated via this particular example. The first observation will allow us to construe downsets in terms of monotone maps (functors) from $\mathcal{P}^{op}$ to the order $\textbf{2}$.\par 
First, consider that any given element $A$ of $\mathcal{D}(\mathcal{P})$ represents something like a ``choice" of elements from the underlying set $P$, with the further requirement that, as a downset, whenever $x \in A$, then any $y \in P$ such that $y \leq x$ in $\mathcal{P}$ is also in $A$. But this requirement is the same as saying that for any $x, y$ such that $y \leq x$ in $\mathcal{P}$, if we have that `it is \textit{true} that $x \in A$', then we must also have that `it is \textit{true} that $y \in A$'. And this is just to say that 
\begin{equation*}
y \leq x \text{ implies } \phi(y) \geq \phi(x),
\end{equation*}  
where $\phi$ is an antitone map from the order $\mathcal{P}$ to order $\textbf{2}$; or, equivalently, it is a monotone map from the opposite order $\mathcal{P}^{op}$ to $\textbf{2}$. Such maps are themselves ordered under the pointwise inclusion ordering. If we designate such a poset of monotone maps, ordered by inclusion, by $\textbf{2}^{\mathcal{P}^{op}}$ or $\text{Monot}(\mathcal{P}^{op}, \textbf{2})$, then we can see that there is a map between the orders 
\begin{align*}
\mathcal{D}(\mathcal{P}) & \rightarrow \text{Monot}(\mathcal{P}^{op}, \textbf{2}) \\ 
D & \mapsto \phi_D,
\end{align*} 
where $\phi_D$ acts as the characteristic (or indicator) function, mapping to $1$ on $D$ and $0$ elsewhere. In other words, given a downset $D$ of $\mathcal{P}$, we define $\phi_D: \mathcal{P}^{op} \rightarrow \textbf{2}$ by setting $\phi_D(x) = 1$ precisely when $x \in D$ (i.e., assigns it to the characteristic function of $D$). Conversely, given a monotone map in $\text{Monot}(\mathcal{P}^{op}, \textbf{2})$, we can send this to the inverse image $\phi^{-1}(1) \in \mathcal{D}(\mathcal{P})$, recovering a unique downset (you can verify for yourself that the subset $\phi^{-1}(1)$ is a downset). In order theory, in general, a map $F: P \rightarrow Q$, where $P$ and $Q$ are posets, is said to be an \textit{order-embedding} \index{order embedding} provided 
\begin{equation*}
x \leq y \text{ in } P \text{ iff } F(x) \leq F(y) \text{ in } Q,
\end{equation*}
and then such an order-embedding yields an \textit{order-isomorphism} \index{order isomorphism} between $P$ and $Q$.
But since, in our case, $A \subseteq  B$ iff $\phi_A \leq \phi_B$, altogether we have thus described an order-embedding, giving us an order-isomorphism   
\begin{equation*}
\mathcal{D}(\mathcal{P}) \cong \text{Monot}(\mathcal{P}^{op}, \textbf{2}). 
\end{equation*}   
Notice how, included among the maps $\mathcal{P}^{op} \rightarrow \textbf{2}$ are those $\phi_p$ for any given element $p \in \mathcal{P}$. These send $q \mapsto 1$ iff $p \leq q$ in $\mathcal{P}$ (or, equivalently, but perhaps more clearly, iff $q \leq p$ in $\mathcal{P}^{op}$, which is our domain), for such a map acts as the indicator function of the set of all $a \leq b$, i.e., the principal downset $\downarrow b$ of $b$. \par 
Observe that the principal downset $\downarrow p : = \{p' \in P : p' \leq p \}$ is itself a downset (i.e., will belong to $\mathcal{D}(\mathcal{P})$), \textit{for any} $p$ in our set $\mathcal{P}$. (This follows from the transitivity of $\mathcal{P}$.) The principal downsets are actually rather special objects among the downsets. To see this, consider again the diagram of the downsets of our particular $\mathcal{P}$. You can see that every object in $\mathcal{D}(\mathcal{P})$ is a principal downset or a union of principal downsets, and that the principal downsets $\downarrow x$ run through all $x \in \mathcal{P}$. Observe also that for $A \subseteq  \mathcal{P}$, $\downarrow A$ will be the \textit{smallest} downset that contains $A$, and moreover $A = \downarrow A$ iff $A$ is a downset. We will return to these facts, and make better sense of them, shortly. \par 
For now, we need to realize that $\downarrow$ can actually be regarded as a monotone map from $\mathcal{P}$ to $\mathcal{D}(\mathcal{P})$, which in our particular case may be pictured as 
\begin{center}
	\includegraphics*[scale=0.25]{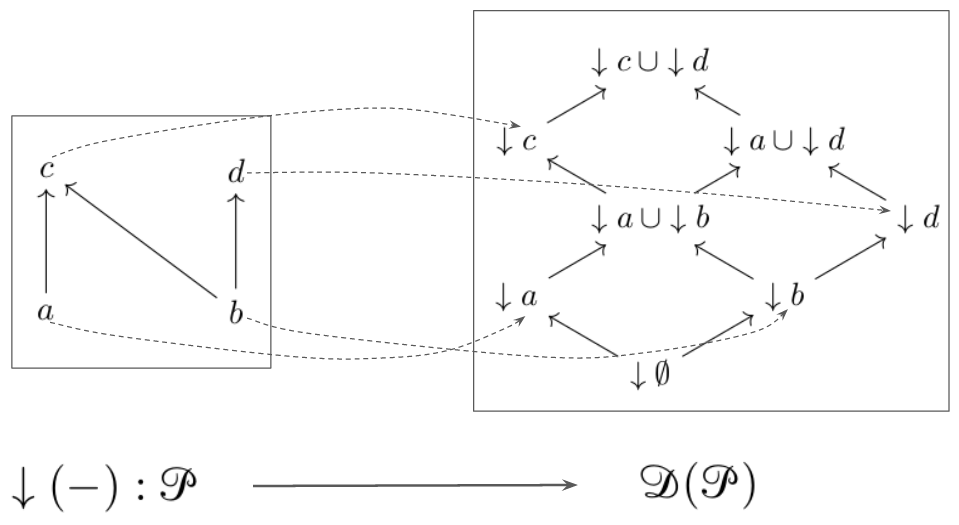}
\end{center}
 That this $\downarrow$ really just defines a monotone map, i.e., $x \leq y$ implies $\downarrow x \subseteq  \downarrow y$, is easy to see in the general case as well. For, let $x \leq y$ in $\mathcal{P}$. We want $\downarrow x \subseteq  \downarrow y$ in $\mathcal{D}(\mathcal{P})$. Take any $x' \in \downarrow x$. Then we must have that $x' \leq x$. By transitivity of the order, $x' \leq x$ and the assumed $x \leq y$ yield $x' \leq y$. Thus, $x' \in \downarrow y$. Altogether, this shows that $\downarrow x \subseteq  \downarrow y$. \par 
 We have the other direction as well, namely $\downarrow x \subseteq  \downarrow y$ in $\mathcal{D}(\mathcal{P})$ implies $x \leq y$ in $\mathcal{P}$. Altogether, then, we actually have another \textit{order-embedding}, embedding any poset into its ``down-set completion":
\begin{equation}
\downarrow(\--): \mathcal{P} \rightarrow \mathcal{D}(\mathcal{P}) \hspace*{1em}
p \mapsto \downarrow p . 
\end{equation}  
This is all part of a much bigger story, so the main results are set aside for emphasis.
\begin{proposition}
	(\textit{Yoneda lemma for posets}) \index{Yoneda lemma!for posets} Given $\mathcal{P}$ a poset, $x \in \mathcal{P}$, and $A \in \mathcal{D}(\mathcal{P})$, then 
	\begin{equation*}
	x \in A \text{ iff } \downarrow x \subseteq  A. 
	\end{equation*}
\end{proposition}
\begin{proof}
	($\Rightarrow$) Take $y \in \downarrow x$. But then $y \leq x$, and so $y \in A$ since $A$ is a downset. \par \noindent 
	($\Leftarrow$) If $x \in A$, then all $y \leq x$ is also in $A$ (as $A$ is a downset), and in particular $x \in \downarrow x$, so $\downarrow x \subseteq  A$.
\end{proof} \noindent 
Applying the Yoneda lemma to two principal downsets, we will have that $y \leq z$ iff for all $x$\footnote{Readers familiar with real analysis may recognize in this the construction of \textit{Dedekind cuts}!}
\begin{equation*}
x \leq y \Rightarrow x \leq z. 
\end{equation*} 
The most important corollary, or application, of the lemma is the following: 
\begin{proposition}
	(\textit{Yoneda embedding for posets}) \index{Yoneda embedding!for posets} This $\downarrow$ defines an order-embedding (embedding any poset into its down-set completion): 
	\begin{equation}
	\downarrow(\--): \mathcal{P} \rightarrow \mathcal{D}(\mathcal{P}) \hspace*{1em}
	p \mapsto \downarrow p.  
	\end{equation}  
\end{proposition}
\begin{proof}
	Let $x \leq y$. Then $x \in \downarrow y$ (conversely, if $x \in \downarrow y$, clearly we must have $x \leq y$). Applying the previous lemma (taking $A = \downarrow y$), $x \in \downarrow y$ holds precisely when $\downarrow x \subseteq  \downarrow y$. On the other hand, the converse holds as well, i.e., $\downarrow x \subseteq  \downarrow y$ implies $x \in \downarrow y$, which implies $x \leq y$. Altogether then, 
	\begin{equation}
	x \leq y \text{ iff } (\downarrow x) \subseteq  (\downarrow y ). 
	\end{equation}
\end{proof} \noindent 
The Yoneda results thus assure us, in a slogan, that 
\begin{quote}
	\textit{To know everything ``below" an element is just to know that element}. 
\end{quote}
Let us start to generalize this story. Given a poset $\mathcal{P}$, by considering $\mathcal{P}$ as a category, we might have looked at presheaves $\mathcal{P}^{op} \rightarrow \textbf{Set}$. But since preorders are the same thing as $\textbf{2}$-enriched categories, supposing we want not \textit{sets} for our hom-sets, but rather ``truth values," it is natural instead to consider $\textbf{2}$-enriched ``presheaves" on $\mathcal{P}$. Instead of arbitrary set-valued data, then, such a $\textbf{2}$-presheaf assigns to each $x \in \mathcal{P}$ a \textit{property} or \textit{truth-value}; and this will just recover the monotone maps (functors) $\mathcal{P}^{op} \rightarrow \textbf{2}$ (or, equivalently, an anti-tone map, or contravariant functor (presheaf),\index{antitone map} from $\mathcal{P}$ to $\textbf{2}$). \par 
We saw how such a monotone map (from the opposite order) effectively acts as the characteristic (or indicator) function $\phi_A$ of a downset $A \subseteq  P$, forming part of the important order-isomorphism 
\begin{equation*}
\mathcal{D}(\mathcal{P}) \cong \text{Monot}(\mathcal{P}^{op}, \textbf{2}). 
\end{equation*}   
We know how to convert any poset into a category, so---renaming $\mathcal{D}(\mathcal{P}) := \textbf{Down}(\mathcal{P})$ and $\text{Monot}(\mathcal{P}^{op}, \textbf{2}) := \textbf{Monot}(\mathcal{P}^{op}, \textbf{2})$, to emphasize that we are now dealing with categories---the above actually describes
\begin{equation*}
\textbf{2-PreSh}(\mathcal{P}) := \textbf{Monot}(\mathcal{P}^{op}, \textbf{2}) \cong \textbf{Down}(\mathcal{P}).
\end{equation*}  
In the order setting, via the principal downsets, we were then able to construct an embedding $\mathcal{P} \rightarrow \textbf{Down}(\mathcal{P})$. Similarly, but in much greater generality, we will see that there is an embedding $\textbf{C} \rightarrow \textbf{PreSh}(\textbf{C})$, taking a general category $\textbf{C}$ to its category of presheaves. Before defining the Yoneda lemma and embedding in the general case, we need to take a step back for a moment and discuss \textit{representability}. To motivate this, consider that in a poset, regarded as a category, a principal downset on an element $p \in \mathcal{P}$ is just all the arrows into $p$ (``all eyes on $p$"), where these amount to all elements that $p$ ``looks down on" (i.e., all elements that ``look up to" $p$). This identification of ``arrows into $p$" and ``\textit{elements} below $p$" can be made, since a poset is precisely a category for which there is \textit{at most one} arrow $q \rightarrow p$ for any $p, q$, allowing us to \textit{identify} $\text{Hom}_{\mathcal{P}}(q, p)$ with the element $q$. In other words, 
\begin{equation*} 
\text{Hom}_{\mathcal{P}}(\--, p) = \downarrow p . 
\end{equation*} 
Before cashing in on the power of this statement, we will discuss the ``representability" operative here by considering a sort of miniaturized version of this phenomenon. 
\subsection{Representability Simplified} 
Consider a general map \index{representability}  
\begin{equation*}
T \times X \rightarrow Y
\end{equation*}
that has a product for its domain (you can think of this as involving sets and functions for now). We cannot typically expect to be able to reduce this to a specification of what is happening on $T$ and $X$ separately, as the interaction of the two factors of the product is essentially involved in supplying the values of the mapping itself. Thus, in considering a map that has for domain a product (all three objects different, for the most general case), 
\begin{equation*}
T \times X \xrightarrow{f} Y, 
\end{equation*}
we can ask (non-trivial) questions about the relations between any of the separate objects involved in the product and the codomain object.  \par 
Observe that if we use the terminal object $\textbf{1}$ to pick out ``points" of $X$, via 
$\textbf{1} \xrightarrow{x} X$, any such point will give rise, via $f$, to the map $f_x$
\begin{center}
	\begin{tikzcd}
		& T \times X \arrow[dr, "f"] \\ 
		T \arrow[ur, "{\langle id_T, \overline{x}\rangle}"] \arrow[rr, "f_x", swap] & & Y 
	\end{tikzcd}
\end{center}  
where $\overline{x}$ is defined as the composite constant map $T \rightarrow \textbf{1} \xrightarrow{x} X$. Thus, 
\begin{equation*}
f_x(t) = f(t, x)
\end{equation*}
for all $t$. In this way, we are now regarding the map $f$ as an $X$-parameterized family of maps $T \rightarrow Y$, one for each of the points of $X$. In this setting, one possible question would be to consider, for a pair of sets $T, Y$, whether there is a set $X$ \textit{large enough} for its points to supply, via the maps $f(\--, x)$, \textit{all} maps $T \rightarrow Y$. As a very simple illustration of this, in the simple case of sets described by their cardinality (number of elements), if $T$ is a set with $4$ elements and $Y$ a set with $3$ elements, then $X$ would need to have 
\begin{equation*}
3^4 = 81
\end{equation*} 
elements, for that is the number of maps $T \rightarrow Y$. \par 
This is in fact part of a much more general story than the above discussion might suggest, one that gets at ``representability."\index{representability} For, we do not need to restrict attention to sets and their cardinal number properties, or even make this a matter of \textit{size}. More generally, for a function $g: T \rightarrow Y$, whenever there is at least one $\textbf{1} \xrightarrow{x_0} X$ such that 
\begin{equation*}
g(\--) = f(\--, x_0),
\end{equation*}  
i.e., for all $t \in T$, 
\begin{equation*}
g(t) = f(t, x_0),
\end{equation*} we say that $g$ is \textit{representable} by $x_0$, or $f$-\textit{represented by} $x_0$. This might be seen as a ``toy" instance of the much more general notion of representable functors. Recall the hom-functor $Hom_{\textbf{C}}(\--, c)$ (and its dual $Hom_{\textbf{C}}(c, \--)$), where this may be given by 
\begin{equation*}
Hom_{\textbf{C}}(\--, \--): \textbf{C}^{op} \times \textbf{C} \rightarrow \textbf{Set}. 
\end{equation*}
For any categories $\textbf{C}, \textbf{D}, \textbf{E}$, an isomorphism can be demonstrated between 
\begin{equation*}
\textbf{E}^{\textbf{C} \times \textbf{D}} \cong (\textbf{E}^{\textbf{D})^{\textbf{C}}},
\end{equation*}
allowing us to move freely between functors $\textbf{C} \times \textbf{D} \rightarrow \textbf{E}$ and $\textbf{C} \rightarrow \textbf{E}^{\textbf{D}}$. Thus, if we fix one of the variables of this Hom$_{\textbf{C}}$, then we get the important\index{functor!representable} \textit{representable functors}: 
\begin{align*} 
Hom_{\textbf{C}}(c, \--) : & \textbf{C} \rightarrow \textbf{Set} \\ 
& d \mapsto Hom_{\textbf{C}}(c, d) \\ 
& (f: b \rightarrow c) \mapsto (f \circ (\--): Hom_{\textbf{C}}(a,b) \rightarrow Hom_{\textbf{C}}(a,c)),
\end{align*} 
and 
\begin{align*} 
Hom_{\textbf{C}}(\--, c) : & \textbf{C}^{op} \rightarrow \textbf{Set} \\ 
& d \mapsto Hom_{\textbf{C}}(d, c) \\ 
& (f: c \rightarrow d) \mapsto ((\--) \circ f: Hom_{\textbf{C}}(b,a) \rightarrow Hom_{\textbf{C}}(c,a)).
\end{align*} 
But this is just to describe the Yoneda-embedding functors taking, for instance in the contravariant case, 
\begin{align*}
\textbf{C} & \rightarrow \textbf{Set}^{\textbf{C}^{op}} \\ 
& x \mapsto Hom_{\textbf{C}}(\--, x). 
\end{align*}
Functors (presheaves) of this form are then said to be \textit{representable}. More formally, 
\begin{definition} 
	For a locally small category $\textbf{C}$, we say that a functor $F:\textbf{C} \rightarrow \textbf{Set}$ is a \textit{representable functor} \index{functor!representable} if there exists an object $c \in \textbf{C}$ (sometimes called the \textit{representing object}) together with a natural isomorphism $Hom_{\textbf{C}}(c, \--) \cong F$; or, equivalently, one speaks of a \textit{representation} for a (covariant) functor $F$ as an object $c \in \textbf{C}$ together with a specified natural isomorphism $Hom_{\textbf{C}}(c, \--) \cong F$.\footnote{Ordinarily, one requires that the domain \textbf{C} of a representable functor be locally small so that the hom-functors $Hom_{\textbf{C}}(c, \--)$ and $Hom_{\textbf{C}}(\--, c)$ are valued in the category of sets. Usually it is said that only set-valued functors are representable, however this may be generalized through the use of Grothendieck universes (or large cardinals).} \par 
	If $F$ is a contravariant functor, then the desired natural isomorphism is given between $Hom_{\textbf{C}}(\--, c) \cong F$.
\end{definition}
 In the covariant case, the representable functor can be thought of, intuitively, as encoding how a category ``is seen" or ``is acted on" by a certain object; in the contravariant case, how the category ``sees" or ``acts on" the chosen object. For instance, in the category of topological spaces $\textbf{Top}$, if we regard all the maps from \textbf{1} (the one-point space) to a space $X$, this just produces the points of $X$, i.e., ``\textbf{1} sees points."\footnote{This example is lifted from \cite{leinster_basic_2014}.} \par 
It is worth lingering a bit with this notion of representability.\index{representability} It might be useful to mention, moreover, that most functors (valued in \textbf{Set}) are \textit{not} representable---if you were to pick a functor randomly, odds are it would not be representable. Thus, examples of \textit{non}-representable functors abound; but perhaps a few concrete non-examples are in order. 
\begin{example}
	The \textit{covariant} powerset functor $\mathbb{P}: \textbf{Set} \rightarrow \textbf{Set}$ is not representable. This functor $\mathbb{P}$ is such that $\mathbb{P}(X)$ is just the power set of $X$, and for any function $X \rightarrow Y$, the map $\mathbb{P}(X) \rightarrow \mathbb{P}(Y)$ takes $A \subseteq  \mathbb{P}(X)$ to the image under $f$, i.e., $f(A)$. \par 
	To see that it is not representable, suppose we have a representing object $X$, i.e., $X \in \textbf{Set}$ represents $\mathbb{P}$. Then, in particular, we will need that 
	\begin{equation*}
	|\text{Hom}_{\textbf{Set}}(X, \--)| = |\mathbb{P}(\--)|
	\end{equation*}
	for all sets, i.e., 
	\begin{equation*}
	|\text{Hom}_{\textbf{Set}}(X, Y)| = |\mathbb{P}(Y)|
	\end{equation*}
	for all $Y \in \textbf{Set}$. But then we can take $Y = \{*\}$, a singleton set. For any non-trivial $X$, there can be only one map to the singleton set, so $|\text{Hom}_{\textbf{Set}}(X, Y)| = 1$. Yet the powerset of a singleton set is, of course, of cardinality $2$. Thus 
	\begin{equation*}
	|\text{Hom}_{\textbf{Set}}(X, Y)| \neq |\mathbb{P}(Y)|
	\end{equation*} 
	for our given $Y = \{*\}$. This contradiction tells us that there can be no such representing object $X$ in $\textbf{Set}$ for the covariant powerset functor. \par 
	On the other hand, the contravariant powerset functor (presheaf) is representable!\footnote{In general, it seems to be easier to find representables among \textit{contravariant} functors.}  
	Specifically, $\mathbb{P}$ (contravariant now) is representable by the 2-element set $2 = \{0,1\}$, so that for each set $Y$, we have the isomorphism 
	\begin{align*}
	Hom_\textbf{Set}(Y, 2) & \cong \mathbb{P}(Y) \\
	f & \mapsto f^{-1}(\{1\}).
	\end{align*} 
	Effectively, this says that $2$ is a set that contains a universal subset $\{1\}$ that pulls back to any other subset, via the characteristic function of that subset. 
\end{example}
The great importance of representable functors is in part due to the fact that representable functors can encode a universal property of its representing object. For instance, a category $\textbf{C}$ will have an \textit{initial object}\index{initial object} precisely when the constant functor $*: \textbf{C} \rightarrow \textbf{Set}$ is representable, i.e., an object $c \in \textbf{C}$ will be initial iff the functor $Y^c$ is naturally isomorphic to the constant functor sending every object to the singleton set. Dually, an object $c \in \textbf{C}$ will be \textit{terminal}\index{terminal object} iff the functor $Y_c$ is naturally isomorphic to the constant functor $*: \textbf{C}^{op} \rightarrow \textbf{Set}$. Put otherwise: an object $c \in \textbf{C}$ is \textit{initial} if, for all objects $d \in \textbf{C}$, there exists a unique morphism $c \rightarrow d$; while an object $c \in \textbf{C}$ is \textit{terminal} if, for all objects $d \in \textbf{C}$, there exists a unique morphism $d \rightarrow c$.\par
The absence of such universal properties can be used, as we effectively did in dealing with the (covariant) powerset functor above, to show that a candidate non-representable functor is in fact not representable. The general idea here---which method you might use to convince yourself of the non-representability of the functors described in the coming examples---is to (i) assume the functor is representable, (ii) consider a possible ``universal" element for the functor, and then (iii) produce a contradiction by showing that this element cannot actually have the special universal property that it needs to have. \par 
Recall that in $\textbf{Set}$, every one-element (singleton) set is a terminal object (a special object with a special universality property). Thus, in our discussion of the covariant powerset functor, another way of saying that 
\begin{equation*}
|\text{Hom}_{\textbf{Set}}(X, \{*\})| \neq |\mathbb{P}(\{*\})|
\end{equation*} 
would accordingly have been to say that $\mathbb{P}$ does not preserve the terminal object. 
\begin{example}
	The (covariant) functor $\textbf{Grp} \rightarrow \textbf{Set}$ that takes a group to its set of subgroups is not representable. 	
\end{example} 
\begin{example}
	The (covariant) functor $\textbf{Rng} \rightarrow \textbf{Set}$ that takes a ring $R$ to its set of squares, i.e., $\{r^2: r \in R\}$, is not representable. 
\end{example}
It turns out that all universal properties themselves can be captured by the fact that certain data defines an initial or terminal object in an appropriate category, specifically the category of elements of the representable functor, a fact that can be rather useful (but that we simply record, without proof, before giving examples). 
\begin{proposition} \label{representable}
	A covariant (contravariant) set-valued functor is representable\index{functor!representable} iff its category of elements\index{category!of elements} has an initial\index{initial object} (respectively, terminal)\index{terminal object} object.
\end{proposition} 
\begin{example}
	Applied to a poset $\mathcal{P}$, consider $\mathcal{P}$ as a category. For an arbitrary element $p \in \mathcal{P}$, first check that the slice category\index{category!slice} $p / \mathcal{P}$ (also denoted $(p\downarrow \mathcal{P})$)\footnote{Hopefully this latter notation is not too confusing in this context, given that we are also talking about principal downsets!} is just the principal downset generated by $p$; dually, the co-slice category\index{category!co-slice} $(\mathcal{P} \downarrow p)$ is the principal upper set of $p$. Recall that we can construct the category of elements in terms of the slice category. Thus, a $\textbf{2}$-presheaf (i.e., downset $A \subseteq  P$) is representable iff it has a greatest element.  
\end{example}
\begin{example}
	Leaving the details to the reader, we indicate that the functor $nColor$ \index{functor!nColoring} described earlier is represented by the functor $K_n$ (the complete graph\index{graph!complete} on $n$ nodes, where a \textit{complete graph} is just a graph in which every pair of distinct vertices is connected by a unique edge), which basically says that if we want to know how many homomorphisms there are from a graph $G$ to the complete graph on $n$-vertices, we could just as well consider all the $n$-colorings of $G$.  
\end{example} \noindent 
In a moment, we will see the most important category-theoretic result, which morally shows how an object is defined completely by its functorial (relational) properties. The next proposition tells us that even if a functorial definition does not correspond to an object, i.e., if the particular functor is not representable, it is still ``built out of" the representables (in particular, it is the colimit of a diagram of representables). 
\begin{proposition}
	Every object $P$ in the presheaf category \index{category!presheaf} $\textbf{Set}^{\textbf{C}^{op}}$, i.e., every contravariant functor on $\textbf{C}$, is a colimit\index{colimit} of a diagram of representable objects, in a canonical way, i.e., 
	\begin{equation*}
	P \cong \text{colim}(\int P \xrightarrow{\pi_P} \textbf{C} \xrightarrow{\textbf{y}} \textbf{Set}^{\textbf{C}^{op}})
	\end{equation*} 
	where $\pi$ is the projection functor and $\textbf{y}$ is the Yoneda embedding.\index{functor!representable}  
\end{proposition} 
This proposition states that given a functor $P: \textbf{C}^{op} \rightarrow \textbf{Set}$, there will be a canonical way of constructing a (small) indexing category $\textbf{J}$ and a corresponding diagram $A: \textbf{J} \rightarrow \textbf{C}$ of shape $\textbf{J}$ such that $P$ is isomorphic to the colimit of $A$ composed with the Yoneda embedding. The indexing category that serves to prove the proposition is the category of elements of $P$.\footnote{A proof of this fact can be found in \cite{riehl_category_2016}. That \textit{every presheaf is a colimit of representable presheaves} is closely related to another construction, namely the \textit{Cauchy completion} (or \textit{Karoubi envelope}) of a category, in which the fact that representable presheaves are \textit{continuous} in a precise sense is exploited. The main idea here is that while we have the powerful Yoneda (full and faithful) embedding sending a category \textbf{C} to the category of presheaves $\textbf{Set}^{\textbf{C}^{op}}$, in general, a category \textbf{C} cannot be recovered from $\textbf{Set}^{\textbf{C}^{op}}$, so a natural question to ask is how or to what extent, given $\textbf{Set}^{\textbf{C}^{op}}$, it can be said to determine \textbf{C}. Basically, if a category \textbf{C} (or $\textbf{C}^{op}$) can be shown to be \textit{Cauchy complete}, then it can not only be recovered (up to equivalence) from the presheaf category (or covariant functor category of variable sets $\textbf{Set}^{\textbf{C}}$), but it can be shown to \textit{generate} the original presheaf (variable set) category.} \par 
As representability often seems to baffle the newcomer, the next (optional) section offers an elaboration on the phenomenon of (non-)representability. The reader eager to press on to the main Yoneda results can skip ahead a few pages.  
\subsection{More on Representability, Fixed points, and a Paradox}
A moment ago, when we were thinking of general morphisms $T \times X \rightarrow Y$ as a family of morphisms $T \rightarrow Y$ parameterized or indexed by the elements of $X$---in which setting we were considering a sort of ``miniature" version of \index{representability} representability---this was effectively to look at arbitrary maps
\begin{equation*}
\hat{f}: X \rightarrow Y^T,
\end{equation*}   
which led to the question of when (and which) $X$ can ``parameterize" all the maps from $T$ itself to some $Y$. This is effectively the same as asking when such $\hat{f}$ are onto (surjective). In particular, though, we can ask this for $X = T$, so that we are considering 
\begin{equation*}
\hat{f}: T \rightarrow Y^T
\end{equation*}   
or 
\begin{equation*}
f: T \times T \rightarrow Y,
\end{equation*}     
which are effectively ``$Y$-valued" \textit{relations} on $T$ (or $Y$-\textit{attributes} of type $T$).
Via the object (function) $Y^T$, $Y$-valued predicates can be thought of as `talking about' $T$. A special circumstance would be where \textit{all} the ways of `talking about itself' can be said by $T$ itself! This is captured by the surjectivity of the map $g$, i.e., when \textit{every} element $f: T \rightarrow Y$ of $Y^T$ is representable in $T$. The next result concerns when this can occur.   
\begin{theorem}
	(\textit{Lawvere's Fixed-Point Theorem}) \index{Lawvere!fixed-point theorem} If
	\begin{equation*}
	\hat{f}: X \rightarrow Y^X
	\end{equation*}
	is surjective, i.e., every $g: X \rightarrow Y$ is representable (in the above sense), then $Y$ will have the \textit{fixed-point property}, i.e., \textit{every} endomap $\tau: Y \rightarrow Y$ has at least one \textit{fixed point},\index{fixed point} where this of course means some $y \in Y$ such that 
	\begin{equation*}
	\tau(y) = y. 
	\end{equation*} 
\end{theorem}
\begin{proof}
	Consider $p: X \rightarrow Y$, an arbitrary ``predicate," i.e., element of $Y^X$. Since any endomap $\alpha: Y \rightarrow Y$ just ``shuffles around" the elements of $Y$, we can define $p$ as the composite of the diagonal map, the function $f$ (got from $\hat{f}$ via the standard exponential conversion), and an endomap, 
	\begin{center}
		\begin{tikzcd}
			X \times X \arrow[r, "{f}"] & Y \arrow[d, "{\alpha}"] \\ 
			X \arrow[u, "{\delta}"] \arrow[r, "p", swap] & Y. 
		\end{tikzcd}
	\end{center}  
	By assumption, moreover, there will be an $x \in X$ that \textit{represents} $p$. Thus, 
	\begin{equation*}
	p(x) = \alpha(f(\delta(x))) = \alpha(f(x,x)) = \alpha(p(x)),
	\end{equation*}
	making $p(x)$ a fixed point of $\alpha$. 
\end{proof} \noindent 
Notice that $Y$ has the fixed-point property provided \textit{every} endofunction on $Y$ has a fixed point. But any set with more than one element clearly has an endofunction on it that does not have a fixed point (hint: the simplest example is a two-point set, where the points are `true' and `false'; then, an endomap without fixed points is given by the familiar `negation' map); thus, no set with more than one element will have the fixed-point property. In order to appreciate the importance of the theorem in $\textbf{Set}$, we can present the theorem in another light, namely via the contrapositive. 
\begin{theorem}
	(\textit{Cantor's Theorem})\index{Cantor} If $Y$ has at least one endomap $\tau$ that has no fixed points (i.e., for all $y \in Y, \tau y \neq y$), then for every object $X$ and for every 
	\begin{equation*}
	X \xrightarrow{\hat{f}} Y^X 
	\end{equation*} 
	$\hat{f}$ is not surjective. \par 
	In other words, $\hat{f}$ not being surjective means that for every attempt $\phi: X \times X \rightarrow Y$ to parameterize maps $X \rightarrow Y$ by the points of $X$, there must be at least one map $g: X \rightarrow Y$ that gets left out, i.e., is not representable by $\phi$ (meaning, does not occur as $\phi(\--, x)$ for any point $x$ in $X$). 
\end{theorem}  
\begin{proof}
	Again, define $g$ as the composite of the diagonal map, the function $f$ (got from $\hat{f}$ via the exponential conversion), and an endomap, 
	\begin{center}
		\begin{tikzcd}
			X \times X \arrow[r, "{f}"] & Y \arrow[d, "{\alpha}"] \\ 
			X \arrow[u, "{\delta}"] \arrow[r, "g", swap] & Y. 
		\end{tikzcd}
	\end{center}  
	In other words, 
	\begin{equation*}
	g(x) = \alpha(f(x,x)).
	\end{equation*}
	Then, for all $x \in X$, 
	\begin{equation*}
	g(\--) \neq f(\--, x)
	\end{equation*}
	as functions of one variable. For, if we \textit{did} have $g(\--) = f(\--, x_0)$ for some $x_0 \in X$, then by evaluation at $x_0$, 
	\begin{equation*}
	f(x_0, x_0) = g(x_0) = \alpha(f(x_0, x_0))
	\end{equation*}
	where the leftmost equality follows from $g$ being representable, and the second equality is by definition. But then, $\alpha$ has a fixed point. Contradiction.  
\end{proof}
Cantor's famous result that there is no surjective map from a set to its powerset, i.e., 
	\begin{equation*}
	X \rightarrow 2^X,
	\end{equation*}
	is a special case of the above.\par 
It is best, though, to see how this special result is part of something more general. Given our way of thinking about maps $Y^X$ as providing a particular way (or name for how) $X$ ``speaks about" or describes itself, the generalized version of the above can be regarded as saying that, provided the ``truth-values" or properties of $X$ are non-trivial, there will be no way that the elements of object $X$ can ``talk about" themselves (in the sense of talking about their own truthfulness or their own properties). The result appeals to an observation concerning the fundamental limitations in how an object $X$ can address its own properties. Many apparent ``paradoxes" of the past seem to play off of this. For instance, the Liar paradox was an ancient way of exhibiting the trouble one can get into when natural languages attempt to construct self-referential statements that speak about their own truthfulness (if you permit this, it seems one must open the door to certain inconsistencies in natural language). Russell's famous paradox was basically a simplified version of something Cantor\index{Cantor} himself already found, one that did not involve the notion of \textit{size}: namely that if we take $T$ as the set of \textit{all} sets, then by Cantor's theorem, there is a set larger than $T$, namely the powerset of $T$, yet $T$ is assumed to contain all sets, so we are saying that $T$ contains a subset that is larger than itself. G\"{o}del's famous incompleteness results revealed limitations in formal systems and provability statements within those systems. Brandenburger-Keisler's paradox (a sort of two-person or interactive version of Russell's paradox) about the description of a belief situation in which ``Ann believes that Bob believes that Ann believes that Bob believes something false about Ann"---the paradox is: does Ann believe that Bob has a false belief about Ann?---suggests that not every description of beliefs can be ``represented." There are a variety of other results\footnote{See \cite{yanofsky_universal_2003} for more; \cite{abramsky_contextual_2014} is also of interest, in this connection.} that one could enumerate as further examples of what are arguably all variations on the same theme: 
\begin{quote}
	\textit{Letting things address their own properties, without limitations, can lead to problems.} 
\end{quote}    
The phenomenon of (non-)representability is really at the core of such problems. The following (apparently ``paradoxical") example is meant mostly ``for fun," to get the reader thinking more about some of the subtleties in issues of representability.  
\begin{example} 
The issue underlying the following example sometimes goes under the name of ``Grelling's paradox."\index{Grelling's paradox}\footnote{\cite{yanofsky_universal_2003} has a very nice discussion of this and a number of other such ``paradoxes."} Consider the set of all English words. Some of these words \textit{describe themselves}, while others (most) do not. Adjectives, perhaps more than any other type of word, are used to describe things. So let us restrict attention to the set of adjectives, which we may denote $Adj$. Among the adjectives, certain of them \textit{describe themselves}, while others (most) do not. Those that describe themselves are said to be \textit{autological} (or \textit{homological}). For instance, the following adjectives are homological: ``English" (is English!); ``polysyllabic" (is polysyllabic); ``Hellenic" (is ``of Greek origin"); ``Unhyphenated." Those adjectives, by contrast, that \textit{do not describe themselves} are said to be \textit{heterological}. For instance, the following are heterological: ``Spanish" (not a Spanish word!); ``Misspelled" (is spelled correctly!); ``Long" (is hardly long); ``Monosyllabic"; ``Hyphenated." \par
It seems plausible that all adjectives will be either homological or heterological. However, consider the adjective ``heterological." Is it heterological? Suppose not. Then, it might naturally be assumed, it will be homological. So it describes itself. Thus ``heterological" (which says that it does not describe itself) must be heterological after all. So if ``heterological" is not heterological, then it is heterological. On the other hand, then, we suppose that the answer to the question is affirmative, i.e., that ``heterological" is heterological. Then ``heterological," being heterological, does not describe itself. But this implies that it is \textit{not} heterological after all (since ``heterological" says is that it is of the sort that does not describe itself, and we just said that ``heterological" does not describe itself, so it is not described by the description ``does not describe itself")! \par  
We might formalize this seemingly ``paradoxical" situation by first considering that we are dealing with a function
\begin{equation*}
f: Adj \times Adj \rightarrow 2
\end{equation*}
defined on all adjectives $a_1, a_2$ by 
\begin{equation*}
f(a_1, a_2) = 
\begin{cases}
1 & \text{if }  a_2 \text{ describes } a_1 \\
0 & \text{if  } a_2 \text{ does not describe } a_1. 
\end{cases}
\end{equation*}
Then we know there is a predicate (a map $Adj \rightarrow 2$) that can be defined on $Adj$ that is not representable by any element of $Adj$. We get this by applying the fixed-point theorem, with $\alpha$ the negation map $\neg: 2 \rightarrow 2$, setting $\alpha(0) = 1$ and $\alpha(1) = 0$. More explicitly, using the idea from before, 
\begin{center}
	\begin{tikzcd}
	Adj \times Adj \arrow[r, "{f}"] & 2 \arrow[d, "{\alpha}"] \\ 
		Adj \arrow[u, "{\delta}"] \arrow[r, "g", swap] & 2, 
	\end{tikzcd}
\end{center} 
we know how to construct $g$ as a (non-representable) function naming a particular property of adjectives, namely as the characteristic function of a subset of adjectives that cannot be described by any adjectives. In particular, the adjective ``heterological" will be in this subset. In terms of the above, that $g$ is such a characteristic function just says that we must have that 
\begin{equation*}
g(\--) \neq f(\--, a)
\end{equation*} 
for all adjectives $a$, since if there \textit{were} an adjective $a_0$ that satisfied $g(\--) = f(\--, a_0)$, evaluating at $a_0$ would give 
\begin{equation*}
f(a_0, a_0) = g(a_0) = \alpha(f(a_0, a_0)),
\end{equation*} 
the first equality from the (assumed) representability of $g$ and the second by definition of $g$. But this is certainly false, due to the nature of the map $\alpha$. \par 
Observe that the hypothetical 
\begin{equation*}
f(a_0, a_0) = g(a_0) = \alpha(f(a_0, a_0)),
\end{equation*}
which yields a contradiction (whether we choose $f(a_0, a_0) = 1$, when ``$a_0$ describes itself," or $f(a_0, a_0) = 0$, when ``$a_0$ does not describe itself"), makes precise exactly the ``paradox" described at the beginning. \par 
Altogether, it is perhaps more telling to consider such a situation in terms of the non-representability of the $g$ given above, where this just means that the property of ``not being described by" an adjective (which applies to ``heterological" in particular) is \textit{not representable}, for there is no adjective that might represent itself via $f$.  
\end{example} 
\subsection{Yoneda in the General}
Let us now make good use of the notions of representability and the model of the Yoneda results for posets. In the special case of posets, we saw that we can identify a principal down set $\downarrow p$ with a representable functor $Hom_{\mathcal{P}}(\--, p)$. For any element $p \in \mathcal{P}$, there will be a representable $\textbf{2}$-presheaf (think: it is ``represented" by $p$) 
\begin{equation*}
\phi_p: \mathcal{P}^{op} \rightarrow \textbf{2} 
\end{equation*}
that takes $q \mapsto 1$ iff $p \leq q$. In this way, the representable presheaves act as the ``characteristic maps" of the principal downsets of $\mathcal{P}$, and the $\textbf{2}$-enriched version of the Yoneda embedding taking each $p \mapsto \phi_p$ is the same as the inclusion of the elements of the poset into its downsets (which is, in turn, the same as considering the $\textbf{2}$-enriched presheaves on $\mathcal{P}$) 
\begin{equation*}
\mathcal{P} \hookrightarrow \textbf{Down}(\mathcal{P}) \cong \textbf{2-}\textbf{PreSh}(\mathcal{P})
\end{equation*} \noindent 
The Yoneda results in the case of categories more generally, i.e., in the $\textbf{Set}$-enriched setting, are effectively a far-reaching generalization of this idea, and supply perhaps the most important and well-utilized results in category theory.   
\begin{proposition}
	(\textit{Yoneda Lemma}) \index{Yoneda lemma} For any functor $F: \textbf{C} \rightarrow \textbf{Set}$, where $\textbf{C}$ is a locally small category, and for any object $c \in \textbf{C}$, the natural transformations $Y^c \Rightarrow F$ are in bijection with elements of the set $F(c)$, i.e.,\footnote{Recall that by $Y^c$ we just mean $Hom_{\textbf{C}}(c,\--)$, while $Y_c$ is used for $Hom_{\textbf{C}}(\--, c)$.}  
	\begin{equation}
	\text{Nat}(Y^c, F) \cong F(c).
	\end{equation} 
	Moreover, this correspondence is natural in both $F$ and $c$. In the contravariant case, i.e., for $F: \textbf{C}^{op} \rightarrow \textbf{Set}$, things are as above, except we have 
	\begin{equation}
	\text{Nat}(Y_c, F) \cong F(c).
	\end{equation} 
\end{proposition} 
We are not going to prove this (it is a good exercise to actually attempt to prove this yourself!), but instead will unpack it and then discuss its significance at a more general level. We will confine attention to the contravariant version in what follows (but dual statements can be made for the covariant version). \par 
The idea is that for a fixed category \textbf{C}, given an object $c \in \textbf{C}$ and a (contravariant) functor $F: \textbf{C}^{op} \rightarrow \textbf{Set}$, we know that the object $c$ gives rise to another special (representable) functor $Y_c: \textbf{C}^{op} \rightarrow \textbf{Set}$. A very natural question to ask, then, is about the maps $Y_c \Rightarrow F$,
\begin{center} 
	\begin{tikzcd}[column sep=huge]
		\textbf{C}^{op}
		\arrow[bend left=50]{r}[name=U,label=above:$Y_c$]{}
		\arrow[bend right=50]{r}[name=D,label=below:$F$]{} &
		\textbf{Set}
		\arrow[shorten <=10pt,shorten >=10pt,Rightarrow,to path={(U) -- node[label=right:$?$] {} (D)}]{} . 
	\end{tikzcd}
\end{center}  
The functors we are comparing both ``live" in $\textbf{Set}^{\textbf{C}^{op}}$, so the collection of maps from $Y_c$ to $F$ are just the natural transformations that belong to $Hom_{\textbf{Set}^{\textbf{C}^{op}}}(Y_c, F)$. But what is this set? Notice that from the input data $F$ and $c$ we were given (``given an object $c$ and a functor $F$"), we could have also constructed the set $F(c)$, by simply applying $F$ on the given object $c$. The Yoneda lemma just assures us that these two sets are the same! Moreover, all the generality of natural transformations is encoded in the particular case of identity maps (used in the proof of the lemma).\par \noindent 
``Naturality" in $F$ in the definition just means that, given any $\upsilon: F \rightarrow G$, the following diagram commutes:\footnote{Note: all the Hom's are $Hom_{\textbf{Set}^{\textbf{C}^{op}}}$.} 
\begin{center}  
	\begin{tikzcd}
		Hom(Y_c, F) \arrow[r, "\cong"] \arrow[d, "{Hom(Y_c, \upsilon)}", swap] & Fc \arrow[d, "\upsilon_c"] \\
		Hom(Y_c, G) \arrow[r, "\cong", swap] & Gc 
	\end{tikzcd} 
\end{center} 
On the other hand, naturality in $c$ means that, given any $h: c \rightarrow c'$, the following diagram commutes: 
\begin{center}  
	\begin{tikzcd}
		Hom(Y_c, F) \arrow[r, "\cong"] & Fc  \\
		Hom(Y_c', F) \arrow[u, "{Hom(Y_h, F)}"] \arrow[r, "\cong", swap] & Fc' \arrow[u, "Fh", swap]
	\end{tikzcd} 
\end{center} 
The most significant application of the Yoneda lemma is given by the \textit{Yoneda embedding}, which tells us that any (locally small) \textbf{C} will be isomorphic to the full subcategory of $\textbf{Set}^{\textbf{C}^{op}}$ spanned by the contravariant representable functors, while $\textbf{C}^{op}$ will be isomorphic to the full subcategory of $\textbf{Set}^{\textbf{C}}$ spanned by the covariant representable functors. We have seen that for each $c \in \textbf{C}$, we have the covariant functor $Y^c$ going from \textbf{C} to \textbf{Set} and the contravariant functor $Y_c$ going from $\textbf{C}^{op}$ to \textbf{Set}. If we let this functor vary over all the objects of \textbf{C}, the resulting functors can be gathered together into the (for example, covariant) functor $Y^{\bullet}: \textbf{C}^{op} \rightarrow Hom(\textbf{C}, \textbf{Set})$; dually, we have the contravariant functor $Y_c$ going from \textbf{C}$^{op}$ to \textbf{Set}, and collecting these functors together as we let $c$ vary will give a functor $Y_{\bullet}: \textbf{C} \rightarrow Hom(\textbf{C}^{op}, \textbf{Set})$.\footnote{It is not unusual to rename these functors, as we do in the following definition, with a lowercase (bold) \textbf{y} in both cases, leaving the appropriate variance to context.}
\begin{definition} The \textit{Yoneda embedding} \index{Yoneda embedding} of \textbf{C}, a locally small category, supplies functors 
	\begin{center} 
		\begin{tikzcd}
			\textbf{C} \arrow[r, hookrightarrow, "{\textbf{y}}"] & \textbf{Set}^{\textbf{C}^{op}} & & \textbf{C}^{op} \arrow[r, hookrightarrow, "{\textbf{y}}"] & \textbf{Set}^{\textbf{C}} \\
			c \arrow[d, "f", swap] \arrow[r, maps to] & Hom(\--, c) \arrow[d] & & c \arrow[d, "f", swap] \arrow[r, maps to] & Hom(c, \--) \\
			d \arrow[r, maps to] & Hom(\--, d) & & d \arrow[r, maps to] & Hom(d, \--) \arrow[u]
		\end{tikzcd}
	\end{center} 
	defining full and faithful embeddings.\footnote{An embedding in the categorical sense is a full and faithful functor. A functor is faithful (full) provided it is injective (surjective) when restricted to each set of morphisms that have a given source and target. Note, moreover, that a full and faithful functor is necessarily injective on objects (up to isomorphism). \par 
	The proof is this result can be found in any text on category theory.} 
\end{definition}
The Yoneda embedding $\textbf{y}$ gives us a representation of \textbf{C} in a category of set-valued functors and natural transformations. An important consequence of the embedding is that any pair of isomorphic objects $a \cong b$ in \textbf{C} are representably isomorphic, i.e., $Y^a \cong Y^b$. The Yoneda lemma supplies the converse, namely if either the (co- or contravariant) functors represented by $a$ and $b$ are naturally isomorphic, then $a$ and $b$ will be isomorphic; so in particular, if $a$ and $b$ represent the same functor, then $a \cong b$. In many cases, it will be easier or more revealing to give such an arrow $Y^a \rightarrow Y^b$ or $Y_a \rightarrow Y_b$ than to supply $a \rightarrow b$, for the category \textbf{Set}$^{\textbf{C}^{op}}$ in general has more structure than does \textbf{C}, i.e., it is complete, cocomplete, and cartesian closed (basically, any morphism defined on a product of two objects can be identified with a morphism defined on one of the factors). Thus we can use the more advanced tools and universal properties (like the existence of limits) that come with the presheaf category, and be sure that an arrow of the form $Y_a \rightarrow Y_b$, for instance, comes from a unique $a \rightarrow b$ even if \textbf{C} on its own may not allow the advanced constructions. Analogously, representing a rational number in terms of downward (upward) closed sets under the standard ordering results in a Dedekind cut, and altogether this embeds the rationals into the reals, allowing for solutions to more equations. Passing from a category \textbf{C} to its presheaf category can also be regarded as adjoining colimits (think generalized sums) to \textbf{C}, and doing so in the most ``free" way.\footnote{This is a powerful and general idea, but the reader who desires a more concrete way of thinking about the previous statement, might consider the unions (\textit{colimits}) that showed up in the downset poset, after we embedded $\mathcal{P}$ into $\mathcal{D}(\mathcal{P})$, where these were not present in $\mathcal{P}$ itself.} In general, in passing to the presheaf category, many non-representable presheaves will show up as well. But, as we saw, the representables have a very special role to play.  
\subsection{Philosophical Pass: Yoneda and Relationality}
Occasionally, an idea is powerful enough that it seems, almost effortlessly, to transcend its local and native context of application and speak articulately to many other contexts. The idea underlying the Yoneda results is like this. In the most general sense, it might be regarded as saying something like  
\begin{quote} 
	\textit{To understand an object it suffices to understand all its relationships with other things}.
\end{quote} 
Or, to say it in another way,
\begin{quote}
	If you want to know whether two objects $A$ and $B$ are the same, just look at whether all the ways of probing $A$ with other things (or, dually, probing other things with $A$) is the same as all the ways of probing $B$ with other things (or, dually, probing other things with $B$).
\end{quote}
The previous two ``slogans" are an attempt to convey what is sometimes called the ``Yoneda philosophy." \index{Yoneda philosophy} It tells us, fundamentally, that if you want to understand what something \textit{is}, there's no need to chase after some ``object \textit{in itself}"; instead, just consider all the ways the (candidate) object transforms, perturbs, and constrains other things (or, dually, all the ways it is transformed, perturbed, and constrained by others)---this will tell you what it is! \par 
This idea that an object is determined by its totality of behaviors in relation to other entities capable of affecting it (or being affected by it), in addition to appealing to many of our intuitions, is an idea that one can find versions of throughout a number of different contexts, for instance in the 17th century philosopher Spinoza's idea that what a body \textit{is} (its ``essence") is inseparable from all the ways that the body can affect (causally influence) and be affected (causally influenced) by \textit{other bodies}. In this general approach, what an object \textit{is} can be entirely encapsulated by regarding ``all at once" (generically) all of its interrelations or possible interactions with the other objects of its world. In the covariant case, we do this by regarding its ways of affecting other things; dually, in the contravariant case, by its particular ways of being affected by the other objects that inhabit its world. To use another metaphor: if you want to understand if two ``destination points" are the same, just inspect whether, ranging over all the addresses with which they can communicate, the networks of routes connecting them to the addresses are the same. For a given object $c$, the representable functor just captures, all at once, the most generic and universal ``picture" of that object, supplying a ``placeholder" for each of the possible attributes of that object (and then Yoneda's lemma just says that to specify an actual object of type $c$, it suffices to fill in all the placeholders for every attribute found in the generic thing of type $c$).\par
If one regards $\text{Hom}_{\textbf{C}}(\--, A)(U) = \text{Hom}_{\textbf{C}}(U, A)$ as telling us about ``$A$ viewed from the perspective of $U$," then the fact (from the Yoneda embedding) that 
\begin{equation*}
(\text{Hom}_{\textbf{C}}(\--, A) \cong  \text{Hom}_{\textbf{C}}(\--, B)) \text{ iff } (A \cong B) \text{ iff } (\text{Hom}_{\textbf{C}}(B, \--) \cong  \text{Hom}_{\textbf{C}}(A, \--))
\end{equation*} 
can be glossed as saying that two objects $A$ and $B$ will be the same precisely when they ``look the same" from all perspectives $U$.\footnote{It is not uncommon to hear such interpretations, namely that we can retrieve the object itself via all the ``perspectives on it." While this is perfectly ``correct" and useful metaphor, when understood in the right way, it is potentially a little misleading. For if taken too literally, it might sound to the unscrupulous listener like a kind of relativism, where objects are reducible to the ``points of view" on them (where it is typically implicitly understood that ``perspectives" are something only the type `human being' is capable of having). But the ``relationism" of Yoneda should not be confused with such a ``relativism." Not all transformations are of the ``perspectival" or ``seeing" type (in fact, probably most are not). Thankfully, such misunderstandings can be blocked by Yoneda itself. For the (inadequate) ``relativist" interpretation would need to assume that there is an object (namely those of the type `human being') that can represent \textit{any} functor whatsoever, that thereby itself mediates all possible exchanges between objects. But Yoneda tells us no such thing.} This paradigm seems especially natural in many contexts, beyond mathematics. It seems especially appropriate to an adequate description of \textit{learning}, wherein an object comes to be known and recognized through probing it with other things, varying the ``perspectives" on it. For instance, suppose you are tasked with having a robot learn how to identify objects that it has never been exposed to before, without relying on much training data or manual supervision, with the aim of having it come away with an ability to correctly discriminate between (what you naively take to be) different objects and readily recognize other instances of objects ``of the same type" in other contexts. \par
You might place the robot in a room with a number of other objects, say, a ceramic cup, some red rubber balls, a steel cable, a small plant, a plastic bottle filled with water, a worm, and a thermostat. If the robot cannot interact with the objects in any way, and the observable interactions and changes that would unfold without its intervention are rather uneventful or slow to unfold, it is not clear how the robot could learn anything at all. On the other hand, suppose you have enabled the robot to inspect, manipulate, or otherwise instigate or probe the objects in a number of ways, and observe the outcomes. At first, these action attempts might be more or less randomized. The robot might simply locally perform simple action sequences or gestures such as 
\begin{quote}
	\textit{Grasp, Release, Put, Pull, Push, Rotate, Twist, Throw, Squeeze, Bend, Stack, See, Hear, Locate}.   
\end{quote}  
It may grasp a rubber ball. It twists the plastic water bottle. It attempts to bend the steel cable. It sees something move (without the robot having performed any other action that might explain this), or hears it wriggling. Such actions can provide the robot with much information about the objects populating the room, and sometimes even the mere successful implementation of a certain isolated action can provide further information, such as about the location (``within reach radius") of an object-candidate that is engaged by \textit{Grasp}. And once the robot has a decent working sense of some possible object-candidates, it can use one to probe others, and learn even more. It may \textit{Push} on a number of objects, with little to no effect, and then pushing on (what we know to be) the thermostat, alter the room's temperature, upon which change it may observe different effects this has throughout the room (e.g., the water evaporates, the plant withers, the worm moves more, other things remain unchanged in certain relevant respects). In this way, our robot goes around the room and probes ``possible object" regions in different ways, and observes the effects of these variations, giving them their own name. Via such probings of the objects of the room, and ``composite" action-sequences, such as \textit{Grasp, then Release, then Hear}, it seems like our robot will have a chance at learning ``what is what." \par  
Later, if we take our robot and place it in a new room, one that has just (what we know to be) a red balloon and a (similarly-shaped) blue rubber ball; and if, in the previous room, all that the robot came to ``know" of (what we know to be) the red rubber ball is the visual information (gathered through \textit{See}) of its color and shape; then it may assume that the red balloon ``is" the object it knew as the rubber ball (which will go under the name of a mapping from some $A_i$ to color data), while it may assume the blue rubber ball is some entirely new thing. On the other hand, if our robot \textit{had} probed the red rubber ball in more ways---say, having subjected it to $Bounce, Twist, Throw, Hear, Grasp$---you can be sure that it would take much less for it to come to recognize that the blue rubber ball was something like the object it knew in the previous room, while the red balloon was something very different. \par  
The idea of Yoneda is that we can be assured that if the robot wants to learn whether some object $A$ is the same thing as object $B$, it will suffice for it learn whether
\begin{equation*}
(\text{Hom}_{\textbf{C}}(\--, A) \cong  \text{Hom}_{\textbf{C}}(\--, B))  
\end{equation*} 
or, dually, 
\begin{equation*}
(\text{Hom}_{\textbf{C}}(B, \--) \cong  \text{Hom}_{\textbf{C}}(A, \--)).
\end{equation*}
In terms of the discussion above, this ``means" having the robot explore whether
\begin{quote}
	all the ways of probing $A$ with objects of its environment amount to the same as all the ways of probing $B$ with objects of its environment. 
\end{quote}
This is a fascinating idea, philosophically, and one that we think has much in its favor even beyond the narrower context of mathematics. Apparently, in Japanese, the word for human being, \textit{ningen}, is made up of two characters, the first of which means something like a human or person, while the second is a representation of the doors of a gate, and means something like ``betweenness," so that the literal meaning of the word ``human" is ``the relation between persons." This---rather than the narrative of ``atomistic individualism," that often ignores or glosses over the immense load of relational constraints and determinations that come with differential obligations, stresses, and opportunities---seems more attuned to the ``Yoneda" way of thinking. \par   
  The fundamental intuition behind the Yoneda philosophy, then, is that to know or access an object it suffices to know or access how it can be transformed by different objects, or how other objects transform into it. More exactly, Yoneda tell us that if there is a natural way of passing an object $c$'s vision of its world (or how it is seen by its world) on to a functor $F$ on that same category, then to recover this vision it suffices to ask $F$ how it acts on $c$. While, mathematically speaking, the usefulness of the lemma often boils down to the fact that we are able to reduce the computation of natural transformations (which can be unwieldy) to the simple evaluation of a (set-valued) functor on an object, in a sense the full philosophical significance of the lemma points in the other direction. Given a category and an object in that category, rather than regard the object ``on its own" (moreover, treating the entire category in a ``detached" manner, as delimiting the outer boundaries of our consideration), via Yoneda we can regard that object as entirely characterized by its perspective or action on its world (or its world's perspective or action on it), and moreover place the category in which it lives in the wider category of all presheaves or sets varying over that category. Via Yoneda, we can perform this sort of passage from the detached consideration of a given object to the consideration of all its interrelations with the other objects of its world \textit{for every object of a given category}. In doing so, we can think of ourselves as taking an entire category $\textbf{C}$ that previously was itself being regarded in a ``detached" manner, and placing it in the more ``continuous" (in a loose sense) context of the category of all the presheaves over $\textbf{C}$. The category of presheaves over $\textbf{C}$ into which \textbf{C} is embedded not only has certain desirable properties that the original category may lack, like possessing all categorical limits, but it can be understood (in both intuitive and in various technical ways) as providing the continuous counterpart to the  ``detached" consideration of the original category.
\section{Adjunctions}
In this final section, we turn to the last of the really fundamental notions in category theory: that of adjointness, or adjoint functors.\index{functor!adjoint} \index{adjunction} The notion of an adjunction is in a sense a generalization or weakening of the notion of an equivalence of categories,\footnote{We have basically already been using this notion of \textit{equivalence of categories},\index{category!equivalence of categories} which is the category-theoretic notion of ``sameness" of categories. More formally, an equivalence of two categories $\textbf{C}$ and $\textbf{D}$ consists of a pair of functors $F: \textbf{C} \rightarrow \textbf{D}$, $G: \textbf{D} \rightarrow \textbf{C}$ that are inverse to each other (up to natural isomorphism of functors), in the sense that we have the natural isomorphisms $F \circ G \cong Id_{\textbf{D}}$ and $G \circ F \cong Id_{\textbf{C}}$.} where we are interested not so much in a relation (or isomorphism) between two categories but in the relation between specific functors between those categories. Another perspective would be to say that adjunctions represent something like a broadening of the notion of an inverse, involving \textit{unique} ``reversal attempts," which moreover frequently exist even when inverses do not.  \par 
Like so many other important notions in category theory, these notions arise in an especially simple form in the special context of orders, so that adjoint relations between orders allow one to display the features of adjointness in a particularly accessible form. Adjoint functors between orders first appeared under the name of \textit{Galois connections}. The next few examples will illustrate and motivate, via very explicit examples in the context of orders, some of the many fundamental general features and properties of adjunctions. 
\subsection{Adjunctions through Morphology}
\begin{example}
	Suppose you receive a dark photocopy of some text, where the pen or marker appears to be ``bleeding": 
	\begin{center}
		\includegraphics*[scale=0.25]{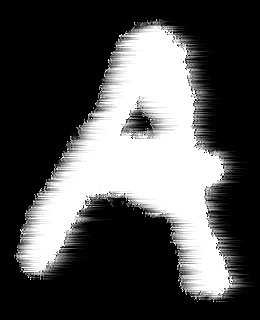}
	\end{center}
	With the help of your favorite programming language, you might perform what the image processing community would call an ``erosion" \index{erosion} of the image. After doing this (perhaps a few times), you would be left with something like
	\begin{center}
		\includegraphics*[scale=0.25]{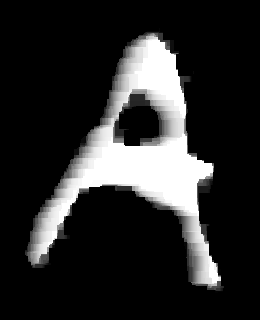}
	\end{center} 
	As one can see, erosion effectively acts to make thicker lines skinnier and detects, or enhances, the holes inside the letter `A'. \par 
	Suppose, instead, that you attempted a \textit{dual} operation, called ``dilation," \index{dilation} the effect of which is to thicken the image, so that lightly drawn figures are presented as if written with a thicker pen, and holes are (gradually) filled. In the case of dilating the original image, you would be left with something like 
	\begin{center}
		\includegraphics*[scale=0.25]{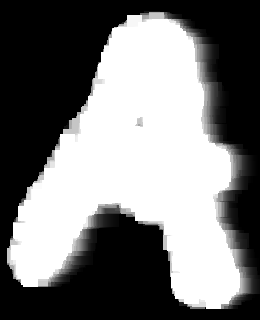}
	\end{center} 
	Now suppose that, for instance, after eroding the image you received, certain things have become harder to read. You decide that you would like to `undo' what you have done, perhaps because you have lost some important information. It seems sensible to hope you might undo it, and get back to the original image, by dilating the result of your erosion. But, in general, erosion and dilation\index{erosion}\index{dilation} do not admit inverses---in particular, they are thus not one another's inverse---and there is no way of determining precisely one image from the dilated image or eroded image. If an image is eroded and then dilated (or conversely), the resulting image will not be the original image. Such operations function to \textit{discard} information, so perhaps it is not so surprising after all that one would not get back to the original by `undoing' an erosion, for instance, by dilating the result.  \par 
	However, in the failed search for an inverse to each operation, you will very quickly alight upon something new, the basic properties of which appear to be useful and interesting in their own right. In particular, eroding after we have dilated an image yields a very different result than dilating after eroding, even though neither composite gives back the original image. However, for an arbitrary image $I$, you will notice that it is ``bounded" in both directions, in the sense of containing the result of one of the two composite operations while being contained by the result of the other, i.e.,  
	\begin{quote}
	\hspace*{2em} Dilating\index{dilation}\index{erosion} after eroding $I$ $\subseteq  I \subseteq $ Eroding after dilating $I$.	
	\end{quote}
	If we erode an image and then dilate the eroded image (making use of the same ``structuring element" through both operations, on which more below), we arrive at a subset of the original image, sometimes called the \textit{opening} \index{opening} of the image by the image processing community. For instance, if we start with the following image
	\begin{center}
		\includegraphics*[scale=0.25]{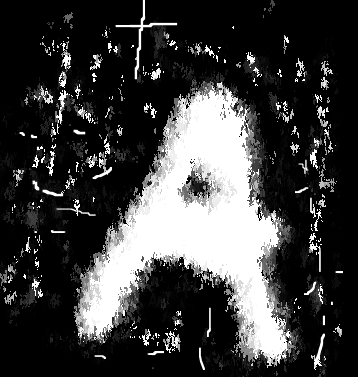}
	\end{center}
	then its opening yields 
	\begin{center}
		\includegraphics*[scale=0.25]{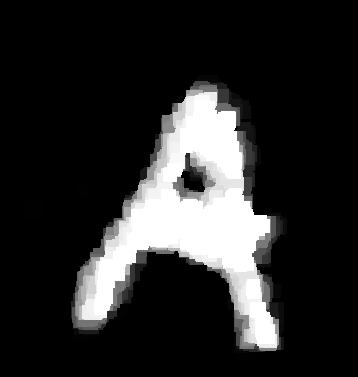}
	\end{center}
	Opening an image will leave one with an image that is generally smaller than the original, as it removes noise and protrusions and other small objects from an image, while preserving the shape and size of the more substantial objects in the image. On the other hand, dilating an image and then eroding\index{dilation}\index{erosion} the dilated image (with the same structuring element throughout), sometimes called \textit{closing} \index{closing} the image, leaves one with an image that is generally larger than the original. Starting with the following image
	\begin{center}
		\includegraphics*[scale=0.27]{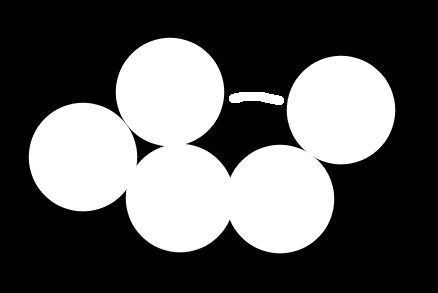}
	\end{center}
	closing it will get rid of small holes, fill gaps in contours, smooth sections of contours, and fuse thin gulfs or breaks between figures. The result of closing the above yields
	\begin{center}
		\includegraphics*[scale=0.27]{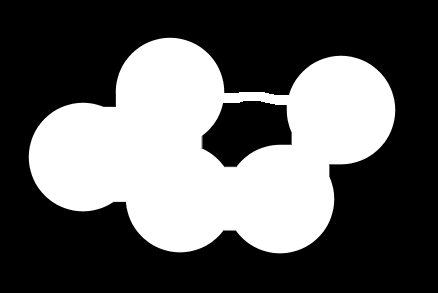}
	\end{center} 
	However, you may quickly learn that opening (or closing) an image twice leaves one with the same image as opening (or closing) it once. In other words, opening and closing are \textit{idempotent} operations. \par 
	Altogether, the two basic operations of dilation and erosion are not quite inverses of one another; yet, as may already be evident from the discussion of the ways their idempotent composites ``bound" the original above and below, they are nevertheless related in a special way---and are, in a sense, the ``closest thing" to an inverse (when this does not exist). We will explore that notion more closely and formally now.\par 
	\textit{Mathematical morphology} \index{morphology} is a field that deals with the processing of binary, gray-level, and other signals, and has proven useful in image processing. The majority of its tools are built on the two fundamental operators of dilation and erosion, and combinations thereof. It has a variety of applications involving image processing and feature extraction and recognition, including applications in X-ray angiography, feature extraction in biometrics, text restoration, etc.\index{dilation}\index{erosion}\par 
	A fundamental idea, in this setting, is to ``probe" an image with a basic, pre-defined shape, and then examine how this fixed shape relates to the shapes comprising the image. One calls the probe the \textit{structuring element}, which is itself a subset of the space, so that, e.g., in the simple case of binary images, such a structuring element is itself just a binary image; these are, moreover, taken to have a defined origin. For instance, in the digital space $E = \mathbb{Z}^2$ (imagine a grid of squares), one might take for structuring element a $3 \times 3$ square, i.e., the set 
	\begin{equation*} 
	\{(-1,-1), (-1,0), (-1,1), (0,-1), (0,0), (0,1), (1,-1), (1,0), (1,1) \},
	\end{equation*} or perhaps a $3 \times 1$ rectangle with a ``central" square, or a diamond or disk-shaped element, and so on. Going back to our earlier example of the blurry binary image letter `A', we may take $0$s to represent background and $1$ for foreground, so that you basically have something like: 
	\begin{equation*} 
	\begin{smallmatrix}
	0 & 0 & 0 & 0 & 1 & 0 & 0 & 0 & 0\\
	0 & 0 & 0 & 1 & 1 & 1 & 0 & 0 & 0\\ 
	0 & 0 & 0 & 1 & 0 & 1 & 0 & 0 & 0\\ 
	0 & 0 & 1 & 1 & 0 & 1 & 1 & 0 & 0 \\
	0 & 0 & 1 & 1 & 1 & 1 & 1 & 0 & 0\\
	0 & 1 & 1 & 1 & 0 & 0 & 1 & 1 & 0\\ 
	0 & 1 & 1 & 0 & 0 & 0 & 1 & 1 & 0\\ 
	1 & 1 & 0 & 0 & 0 & 0 & 0 & 1 & 1 
	\end{smallmatrix}
	\end{equation*} 
	Then, we might take for structuring element $B \subseteq  E$ a ``diamond" (with origin boxed off)
	\begin{equation*} 
	\begin{smallmatrix}
	0 & 1 & 0 \\ 
	1 & \fbox{1} & 1 \\ 
	0 & 1 & 0
	\end{smallmatrix} 
	\end{equation*}
	In the simple case of our running example of a binary image in a bounded region, we took the structuring element $B$, placing $B$'s origin at each pixel of our image $X$ as we ``scan" over all of $X$,\index{dilation}\index{erosion} and compute at each pixel the dilation and erosion by taking the maximum or minimum value, respectively, of all pixels within the ``window" or neighborhood covered by the structuring element (so that, e.g., in the case of dilation, a pixel is set to 1 if any of its neighboring pixels have the value 1). \par  
	More specifically, assuming we have fixed an origin in $E$, to each point $p$ of $E$ there will correspond the translation map that takes the origin to $p$; such a map will then takes $B$ in particular onto $B_p$, the translate of $B$ by $p$. In general, translation by $p$ is a map $E \rightarrow E$ that takes $x$ to $x + p$; thus, it takes any subset $X$ of $E$ to its \textit{translate} by $p$, 
	\begin{equation*}
	X_p = \{x+p \hspace*{0.25em}| \hspace*{0.25em} x \in X \}. 
	\end{equation*} 
	For a structuring element $B$, then, we can consider all its translate $B_p$. Given a subset (image) $X$ of $E$, we can examine how the translates $B_p$ of a given structuring element $B$ interact with $X$. In the simple case of Boolean images (as subsets of a Euclidean or digital space), we carry out this examination via two operations: 
	\begin{equation}
	X \oplus B = \{x + b \hspace*{0.25em}| \hspace*{0.25em} x \in X, b \in B \} = \bigcup_{x \in X} B_x = \bigcup_{b \in B} X_b,
	\end{equation}
	called Minkowski addition, and its dual, 
	\begin{equation}
	X \ominus B = \{p \in E \hspace*{0.25em}| \hspace*{0.25em} B_p \subseteq  X \} = \bigcap_{b \in b} X_{-b}.
	\end{equation}
	The former transformation (taking $X$ into $X \oplus B$) is in fact what gives us a \textit{dilation}, \index{dilation} the basic property of which is that it distributes over union, while the latter (taking $X$ into $X \ominus B$) is an \textit{erosion}, \index{erosion} the basic property of which is that it distributes over intersection. In the simple set-theoretical binary image case, dilations coincide with Minkowski addition; yet erosion of an image is the intersection of all translations by the points $-b$. In short,  
	\begin{quote}
		Dilation of $X$ by $B$ is computed as the union of translations of $X$ by the elements of $B$, 
	\end{quote}
while 
\begin{quote}
	Erosion of $Y$ by $B$ is the intersection of translations of $Y$ by the reflected elements of $B$.
\end{quote} \noindent 
While we will see that we can give more general definitions of these operations, the particular behaviors of these operations in fact already follow from a general relationship underlying these operations.
\begin{proposition} 
	For every subset $X, Y, B$ of our space $E$, where $B$ is any structuring element, we have  
	\begin{equation}
	X \oplus B \subseteq  Y \text{ iff } X \subseteq  Y \ominus B. 
	\end{equation}
\end{proposition}
	\begin{proof}
		$(\Rightarrow)$ Suppose $X \oplus B \subseteq  Y$ and let $z \in X$ and $b \in B$. Then $z + b \in X \oplus B$, and thus $z + b \subseteq  Y$. And $z + b \subseteq  Y$ for any $b \in B$ implies that $z \in Y \ominus B$.\index{dilation}\index{erosion} \par \noindent  
		$(\Leftarrow)$ Suppose $X \subseteq  Y \ominus B$ and let $z \in X \oplus B$. Then there exists $x \in X$ and $b \in B$ such that $z = x + b$. But $x \in X$ and $X \subseteq  Y \ominus B$ entails that $x \in Y \ominus B$. Thus, for every $b' \in Y$, we have $x + b' \in Y$, and in particular, $b \in B$, so $x + b \in Y$. But $z = x + b$, so $z \in Y$. 
	\end{proof}
	Before discussing the significance of this more generally, it is worth noting that there is no need to restrict attention, as we have thus far, to consideration of dilation and erosion in the simple case of Boolean images in digital space. We can also define the dilation and erosion of a \textit{function} by a structuring element (itself regarded as a \textit{structuring function}). For instance, given a function $f: E \rightarrow T$ from a space $E$ to $T$ a set of, e.g., grey-levels (i.e., a complete lattice that comes from a subset of $\overline{\mathbb{R}} = \mathbb{R} \cup \{-\infty, +\infty\}$), and given a point $p \in E$, the \textit{translate} of $f$ by $p$ is the function $f_p$ whose graph is obtained by translating the graph $\{(x, f(x)) \hspace*{0.25em}| \hspace*{0.25em} x \in E\}$ by $p$ in the first coordinate, i.e., $\{(x+p, f(x)) \hspace*{0.25em}| \hspace*{0.25em} x \in E\}$, so that for all $y \in E$, we have 
\begin{equation}
f_p(y) = f(y-p). 
\end{equation}
This defines the translation of a function by a point. Of course, if extend this to the translation by a pair $(p, t)$, we get that for all $y \in E$, 
\begin{equation}
f_{(p,t)}(y) = f(y-p) + t.
\end{equation}
This approach allows us to define Minkowski addition and the dual operation for two \textit{functions} $E \rightarrow T$. Where $f$ plays the role of an (grey-level) image, and $b$ is the functional analogue of a structuring element (i.e., a \textit{structuring function}), for all $p \in E$, these operations will take on values 
\begin{equation}
(f \oplus b)(p)  = \sup_{y \in E} (f(y) + b(p-y)) 
\end{equation}
and 
\begin{equation}
(f \ominus b)(p) = \inf_{y \in E} (f(y) - b(y-p)).
\end{equation}
Then the operator $\delta_g: T^E \rightarrow T^E$ taking $f \mapsto f \oplus g$ is \textit{dilation by} \index{dilation} $g$, and $\epsilon_g: T^E \rightarrow T^E$ taking $f \mapsto f \ominus g$ is \textit{erosion by} \index{erosion} $g$. In a low-dimensional case, using a portion of a disk for structuring function, these operators might do something like:\footnote{This image (and the one that follows) is taken, with slight modifications, from \cite{hlavac_grayscale_nodate}.} 
\begin{center}
	\includegraphics*[scale=0.25]{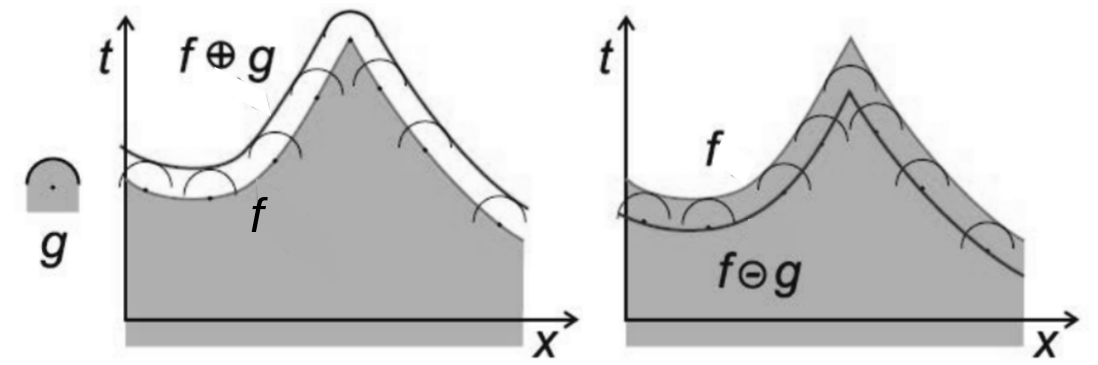}
\end{center}
When we use a ``flat structuring" element instead, such as that represented by a line, things are simplified even further, and we are effectively applying max and min filters, i.e., for dilation
\begin{equation*}
(f \oplus g)(x) = \sup_{y \in E, x-y \in B} f(y) = \sup_{y \in B_x} f(y) 
\end{equation*}
and for erosion 
\begin{equation*}
(f \ominus g)(x) = \inf_{y \in E, x-y \in B} f(y) = \inf_{y \in \breve{B}_x} f(y), 
\end{equation*}
where $\breve{B}$ is the transpose or symmetrical of $B$, i.e., $\{-b \hspace*{0.25em}|\hspace*{0.25em} b \in B\}$. Erosion by a flat structuring function acts to shrink peaks and flatten valleys, while dilation acts in a dual fashion, flattening or ``rounding" peaks and accentuating valleys. Taking, for instance, some price data, we might then have something like 
\begin{center}
	\includegraphics*[scale=0.25]{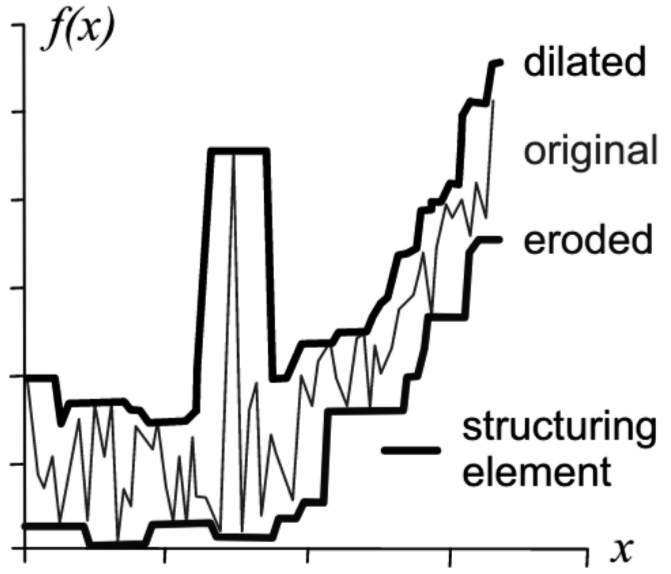}
\end{center}
Whether in its functional treatment, or in terms of the special relation $X \oplus B \subseteq  Y$ iff $X \subseteq  Y \ominus B$, the relation between dilation and erosion\index{dilation}\index{erosion} is part of a much more general and powerful story, exemplifying the notion of a Galois connection,\index{Galois connection} itself an instance of the more general notion of an adjunction. We first supply the relevant definitions, and then further explore some of the powerful abstract features of this notion through the particular case of the operations of dilation and erosion. 
\begin{definition} 
Let $\mathcal{P} = (P, \leq_P)$ be a preordered set, and $\mathcal{Q} = (Q, \leq_Q)$ another preorder. Suppose we have a pair of monotone maps $F: P \rightarrow Q$ and $G: Q \rightarrow P$, 
\begin{center}  
	\begin{tikzcd}
		P \arrow[r, shift left = .7ex, "F"] &
		Q \arrow[l, shift left = .2ex, "G"]
	\end{tikzcd} 
\end{center} 
such that for all $p \in P$ and $q \in Q$, we have the two way rule 
\begin{center} 
	\begin{equation*}
	\frac{F(p) \leq q}{p \leq G(q)}
	\end{equation*}
\end{center} 
where the bar indicates `iff'. If such a condition obtains, the pair $(F, G)$ is said to form a monotone \textit{Galois connection} between $\mathcal{P}$ and $\mathcal{Q}$. 
\end{definition} \noindent 
When such a connection obtains, we also say that $F$ is the \textit{left (or lower) adjoint} and $G$ the \textit{right (or upper) adjoint} of the pair, and write (for reasons we will see in a moment) $F \dashv G$ to indicate the relation. \par 
Such a situation can be expressed in terms of the behavior of certain special arrows associated to each object of $P$ and $Q$. In particular, let $p$ be an object of $P$, and set $q = F(p)$. Then 
\begin{center} 
	\begin{equation*}
	\frac{F(p) \leq F(p)}{p \leq G(F(p))}
	\end{equation*}
\end{center} 
where the top is the identity arrow on $F(p)$ in $Q$, indicating reflexivity (which holds in any order). If we use $\theta$ to designate the bijection that realizes the `if and only if' situation, then $\theta(\text{id}_{F(p)})$ is a special arrow of $Q$, called the \textit{unit} \index{adjunction!unit of} of $p$, where this arrow enjoys a certain universality property. There is a corresponding dual notion of a \textit{counit}.\index{adjunction!counit of} In short, for each $p \in P$, we call the \textit{unit} an element $p \leq GFp$ that is least among all $x$ with $p \leq Gx$; dually, for each $q \in Q$, the \textit{counit} is an element $FGq \leq q$ that is greatest among all $y$ with $Fy \leq q$. It can be shown that, given order-preserving maps $F: P \rightarrow Q$ and $G: Q \rightarrow P$, it is in fact \textit{equivalent} to say (i) $F \dashv G$ and (ii) $p \leq GFp$ and $FGq \leq q$.  
\par 
Changing the variance of the functors involved gives us a slightly different notion. 
\begin{definition}
	 Let $\mathcal{P} = (P, \leq_P)$ and $\mathcal{Q} = (Q, \leq_Q)$ be orders. Suppose we have a pair of anti-tone (order-reversing) \index{antitone map} maps $F: P \rightarrow Q$ and $G: Q \rightarrow P$, 
	 \begin{center}  
	 	\begin{tikzcd}
	 		P \arrow[r, shift left = .7ex, "F"] &
	 		Q \arrow[l, shift left = .2ex, "G"]
	 	\end{tikzcd} 
	 \end{center} 
	 such that for all $p \in P$ and $q \in Q$, we have the two way rule 
	 \begin{center} 
	 	\begin{equation*}
	 	\frac{q \leq F(p)}{p \leq G(q)} . 
	 	\end{equation*}
	 \end{center}  
 If such a condition obtains, the pair $(F, G)$ is said to form an \textit{anti-tone Galois connection} between $\mathcal{P}$ and $\mathcal{Q}$. 
\end{definition}
We have seen many times now how any order $\mathcal{P}$ can be regarded as a category by taking
\begin{equation*}
x \leq_P y \text{ iff there exists an arrow } x \rightarrow y.
\end{equation*}
And in this setting, we further know that covariant functors between such categories are just monotone (order-preserving) functions, and that contravariant functors are anti-tone (order-reversing) functions. Thus, it is entirely natural to attempt to regard Galois connections in a more general, categorical guise. Doing so gives us the notion of an adjunction, which can accordingly be seen as a straightforward categorical generalization of the notion of a (monotone) Galois connection. 
	\begin{definition}
		An \textit{adjunction} \index{adjunction!defined} is a pair of functors $F: \textbf{C} \rightarrow \textbf{D}$ and $G: \textbf{D} \rightarrow \textbf{C}$ such that there is an isomorphism
		\begin{equation*}
		Hom_{\textbf{D}}(F(c), d) \cong Hom_{\textbf{C}}(c, G(d)),
		\end{equation*}
		for all $c \in \textbf{C}, d \in \textbf{D}$, which is moreover natural in both variables. When this obtains, we say $F$ is left adjoint to $G$, or equivalently $G$ is right adjoint to $F$, denoted $F \dashv G$.\footnote{Sometimes the morphisms $F(c) \xrightarrow{f^{\sharp}} d$ and $c \xrightarrow{f^{\flat}} G(d)$ of the bijection given above are said to be \textit{adjunct} or \textit{transposes} of each other.} \par \noindent 
		In saying that the isomorphism is ``natural in both variables", we mean that for any morphisms with domain and codomain as below, the square on the left commutes (in \textbf{D}) iff the square on the right commutes (in \textbf{C}): 
		\begin{center}
			\begin{tikzcd}[row sep = small]
				F(c) \arrow[dd, "F(h)", swap] \arrow[r, "{f^{\sharp}}"] & d \arrow[dd, "k"] & & c \arrow[dd, "h", swap] \arrow[r, "{f^{\flat}}"] & G(d) \arrow[dd, "G(k)"] \\
				& & \iff \\
				F(c') \arrow[r, "{g^{\sharp}}", swap] & d' & & c' \arrow[r, "{g^{\flat}}", swap] & G(d') 
			\end{tikzcd}
		\end{center} 
	\end{definition} 
As one might expect, by considering functors of different variance, corresponding to \textit{antitone} Galois connections, there is another notion, namely that of \textit{mutually right adjoints} (and further, mutually left adjoints). 
\begin{definition}
	Given a pair of functors $F: \textbf{C}^{op} \rightarrow \textbf{D}$ and $G: \textbf{D}^{op} \rightarrow \textbf{C}$, if there exists a natural isomorphism 
	\begin{equation}
	Hom_{\textbf{D}}(F(c), d) \cong Hom_{\textbf{C}}(G(d), c),
	\end{equation}
	then we say that $F$ and $G$ are \textit{mutually left adjoint}. Given the same functors, if there exists a natural isomorphism 
	\begin{equation}
	Hom_{\textbf{D}}(d, F(c)) \cong Hom_{\textbf{C}}(c, G(d)),
	\end{equation}
	then we say that $F$ and $G$ are \textit{mutually right adjoint}.
\end{definition} \noindent 
An antitone Galois connection just names a mutual right adjoint situation between preorders (posets). \par 
	Following the example of the Galois connection definition, we can also recover the notions of the unit and counit of an adjunction.
	\begin{definition}
		Given an adjunction $F \dashv G$, there is a natural transformation 
		\begin{equation*}
		\eta: \text{id}_{\textbf{C}} \Rightarrow GF
		\end{equation*} 
		called the \textit{unit} of the adjunction.\index{adjunction!unit of} Its component 
		\begin{equation*}
		\eta_c: c \rightarrow GFc
		\end{equation*}
		at $c$ is the transpose of the identity morphism $\text{id}_{Fc}$. \par 
		Dually, there is a natural transformation $\mu: FG \Rightarrow \text{id}_{\textbf{D}}$, called the \textit{counit} of the adjunction, \index{adjunction!counit of} with component 
		\begin{equation*}
		\mu_d: FGd \rightarrow d
		\end{equation*}
		at $d$ defined as the transpose of the identity morphism $\text{id}_{Gd}$. 
	\end{definition}
Any adjunction comes with a unit and a counit. In fact, conversely, given opposing functors $F: \textbf{C} \rightarrow \textbf{D}$ and $G: \textbf{D} \rightarrow \textbf{C}$, supposing they are equipped with natural transformations $\eta: \text{id}_{\textbf{C}} \Rightarrow GF$ and $\mu: FG \Rightarrow \text{id}_{\textbf{D}}$ satisfying a pair of conditions, then this data can be used to exhibit $F$ and $G$ as adjoint functors. In other words, we can use the natural transformations exemplifying the counit and unit maps, together with some conditions on these, to actually define an adjunction. 
\begin{definition}
	(\textit{Adjunction, again}) \index{adjunction!defined} An \textit{adjunction} consists of a pair of functors $F: \textbf{C} \rightarrow \textbf{D}$ and $G: \textbf{D} \rightarrow \textbf{C}$, equipped with further natural transformations $\eta: \text{id}_{\textbf{C}} \Rightarrow GF$ and $\mu: FG \Rightarrow \text{id}_{\textbf{D}}$ satisfying what are sometimes called the \textit{triangle identities}:\index{triangle identities}
	\begin{center}
		\begin{tikzcd}
			F \arrow[dr, Rightarrow, "{\text{id}_F}", swap] \arrow[r, Rightarrow, "F\eta"] & FGF \arrow[d, Rightarrow, "\mu F"] & & 	G \arrow[dr, Rightarrow, "{\text{id}_G}", swap] \arrow[r, Rightarrow, "\eta G"] & GFG \arrow[d, Rightarrow, "G \mu"] \\
			& F & & & G
		\end{tikzcd}
	\end{center}  
\end{definition} \noindent 
Then, the isomorphism $Hom_{\textbf{D}}(F(c), d) \cong Hom_{\textbf{C}}(c, G(d))$ realizing $F$ and $G$ as an adjoint pair, will exist precisely where there exists a pair of natural transformations, as above, satisfying the triangle identities.\par \vspace*{1em} 
	Let us now return to the dilation and erosion of images and breath some life into these ideas. The special relation relating Minkowski addition to its dual did not really depend on the particular form given to it in the translation-invariant case of a binary image, but exemplifies a more general notion of dilation and erosion\index{dilation}\index{erosion} on an arbitrary complete lattice, using the operations of supremum and infimum, where we have the adjunction\index{dilation}\index{erosion}
	\begin{equation*}
	\text{dilate } \dashv \text{ erode}. 
	\end{equation*}
	 To see how this works, we can first observe that both operations, dilation and erosion, are order-preserving (monotone), in the sense that $X \subseteq  Y$ implies $X \oplus Z \subseteq  Y \oplus Z$ and also $X \ominus Z \subseteq  Y \ominus Z$. Moreover, while the order of an image intersection (union) and a dilation (erosion) cannot be interchanged freely, the dilation of the union of two images is indeed equal to the union of the dilations of the images, so the order can be interchanged; likewise, erosion of the intersection of two images yields the intersection of their erosions. \par 
	In dealing with the pair (dilate, erode), there is no need to be restricted to the poset of subsets of a digital space. There are clearly many choices for the underlying object space, i.e., of where the images in question are held to reside. It is common to consider $\mathbb{P}(E)$ the space of all subsets of $E$ (where $E$ is the $d$-dimensional Euclidean space $\mathbb{R}^d$ or the digital space $\mathbb{Z}^d$). But we might also consider Conv($E$) the space of all convex subsets of $E$; or $\mathcal{P}^E$ the space of ``image" functions from $E$ a discrete space to $\mathcal{P}$ a pixel lattice; or the space $(T^3)^{E}$ of RGB color images (where RGB colors are triples ($r,g,b$) of numerical values, $T^3$ the lattice of RGB colors under componentwise order, and an RGB image as a function $E \rightarrow T^3$ taking each point $p \in E$ to a triple $(r(p), g(p), b(p))$ representing the RGB coloration of $p$); and so on. The take-away, though, is that despite the differences in these underlying spaces, all such spaces form \textit{complete lattices} (where this means it has all joins and meets, not just finite or binary ones). So if we consider a complete lattice $\mathcal{L}$, with the order $\leq$, supremum $\bigvee$, infimum $\bigwedge$, least element $0$, and greatest element $I$, such a lattice can be thought of as our ``image lattice," corresponding to a particular set of images we are working with. Traditionally, one then defines dilations and erosions as follows. 
	\begin{definition}  Let $\mathcal{L}$ and $\mathcal{M}$ denote complete lattices. Then for $\delta: \mathcal{L} \rightarrow \mathcal{M}$ and $\epsilon: \mathcal{M} \rightarrow \mathcal{L}$, we say that 
		\begin{itemize}
			\item $\delta$ is a \textit{dilation} \index{dilation!defined} provided for every $S \subseteq  \mathcal{L}$, 
			\begin{equation}
			\delta(\bigvee S) = \bigvee_{X \in S} \delta(X). 
			\end{equation}
			\item $\epsilon$ is an \textit{erosion} \index{erosion!defined} provided for every $T \subseteq  \mathcal{M}$, 
			\begin{equation}
			\epsilon(\bigwedge T) = \bigwedge_{Y \in T} \epsilon(Y). 
			\end{equation}
		\end{itemize}
		Note that this also applies in the case of $S, T$ empty, in which case a dilation is held to preserve $0$, while an erosion preserves $I$. 
	\end{definition} 
Our dilate-erode pair is actually an antitone Galois connection, where $\mathcal{M} = \mathcal{L}^{op}$, which just means that for all $S, T$ in $\mathcal{L}$ 
\begin{center} 
	\begin{equation*}
	\frac{T \leq^{op} \delta(S)}{S \leq \epsilon(T)}
	\end{equation*}
\end{center}   
which is, of course, the same as 
\begin{center} 
	\begin{equation*}
	\frac{T \geq \delta(S)}{S \leq \epsilon(T)}
	\end{equation*}
\end{center}    
or, equivalently, 
\begin{center} 
	\begin{equation*}
	\frac{\delta(S) \leq T}{S \leq \epsilon(T)}
	\end{equation*}
\end{center}    
where the order here is now the same, that given on $\mathcal{L}$, above and below the line. Thus, we have recovered the usual notion of an adjunction with $\delta: \mathcal{L} \rightarrow \mathcal{L}$ order-preserving and $\epsilon: \mathcal{L} \rightarrow \mathcal{L}$ order-preserving! Dilations and erosions\index{dilation}\index{erosion} are then precisely just the order-preserving (monotone) transformations on a complete lattice that moreover commute with the union (supremum) and intersection (infimum), respectively.\footnote{Note that, if we were to regard the pair as comprising an anti-tone Galois connection, then we would be saying that the operators exchanged suprema and infima, in the sense that, e.g., $\delta(\bigvee_i x_i) = \bigwedge \delta(x_i)$.} 
\par 
Morphological operators are thereby given a unified treatment in the general framework of an adjoint pair on complete lattices. A number of well-established properties concerning the interaction of these operators then fall out immediately from the general framework of adjunctions. Conversely, we can illustrate such general facts via the present operators on images. \par 
Suppose we have an adjunction $\delta \dashv \epsilon$ on a complete lattice $\mathcal{L}$.\footnote{Looking ahead to what will have to be true of such functors, they have the names they do; however, we do not yet require anything about the maps $\delta$ and $\epsilon$, except that they form an adjoint pair moving between $\mathcal{L}$ and itself.} Then a number of morphologically significant facts come ``for free" as corollaries of general categorical truths about an adjoint pair. Even the fact that $\delta$ is a dilation and $\epsilon$ is an erosion in the first place can be derived from the existence of this special adjoint relationship. In what follows, we explore some of these general truths through the lens of some notable particular truths about dilations and erosions.  
\begin{enumerate}
	\item (\textbf{Uniqueness of Adjoints}) 
	\begin{proposition} To each dilation\index{dilation}\index{erosion} $\delta$ there corresponds a unique erosion $\epsilon$, namely 
	\begin{equation*}
	\epsilon(X) = \bigvee \{S \in \mathcal{L}| \delta(S) \leq X\},
	\end{equation*}
	and to each erosion $\epsilon$ there corresponds a unique dilation: 
	\begin{equation*}
	\delta(X) \bigwedge \{S \in \mathcal{L}| \epsilon(S) \geq X\}.
	\end{equation*}
	\end{proposition} 
	This ultimately derives from a general result that assures us that, like inverses, adjoints are unique (well, actually ``unique up to unique isomorphism," but we can ignore this in our special case): 
	\begin{proposition}
		Adjoint maps are unique.
	\end{proposition} 
In the case of orders, with order-preserving maps between them, this just means 
\begin{enumerate}
	\item if $F_1$ and $F_2$ are left adjoints of $G$, then $F_1 = F_2$. 
	\item if $G_1$ and $G_2$ are right adjoints of $F$, then $G_1 = G_2$.
\end{enumerate}
\begin{proof}
	(We focus on the simple case of orders, and prove (a); (b) follows by duality) From the adjointness assumptions, we have both  
	\begin{equation*}
	\frac{F_1(p) \leq q}{p \leq G(q)}
	\end{equation*}
	and 
	\begin{equation*}
	\frac{p \leq G(q)}{F_2(p) \leq q}, 
	\end{equation*}
	so immediately we have that $F_1(p) \leq q$ iff $F_2(p) \leq q$. Set $q = F_1(p)$, making $F_1(p) \leq q$ trivially true, forcing $F_2(p) \leq F_1(p)$ to be true as well. Similarly, set $q = F_2(p)$ and use the trivial truth $F_2(p) \leq F_2(p)$ to force $F_1(p) \leq F_2(p)$. In a poset, this entails that $F_1(p) = F_2(p)$, $p$ arbitrary. 
\end{proof}
The adjunction then gives rise to the formulas 
\begin{equation*}
G(q) = \bigvee \{p|F(p) \leq q\}
\end{equation*}
and 
\begin{equation*}
F(p) = \bigwedge \{q | p \leq G(q)\}
\end{equation*}
which displays the uniqueness of the adjoints, and so explains the unique erosion (dilation) corresponding to each dilation (erosion),\index{dilation}\index{erosion} as written above. \par 
In general, a given map may or may not have a left (or right) adjoint; the map may have one without the other, neither, or both (where these may be the same or different). But if it does have a left (or right) adjoint, we can be confident that, even though they are not quite inverses, the adjoint is unique up to isomorphism. \par 
Adjoint functors also interact in particularly interesting and useful ways with the limit and colimit constructions, a connection we now explore. \par 
	\item (\textbf{Limit and Colimit Preservation})\index{limit}\index{colimit} 
	\begin{proposition} $\delta$ is a dilation and $\epsilon$ is an erosion, and both are order-preserving. 
	\end{proposition} 
This follows immediately from a very important category-theoretic result, namely that 
	\begin{proposition}
		Right adjoints preserve limits (\textit{RAPL}); Left adjoints preserve colimits (\textit{LAPC}).\footnote{Terminologically, a general functor that is limit-preserving is often said to be a \textit{continuous} functor, while a colimit-preserving functor is a \textit{cocontinuous} functor. Another related concept we will make use of later on in the book is the following: a functor is said to be \textit{left exact} if it preserves \textit{finite} limits, and \textit{right exact} if it preserves \textit{finite} colimits.\par 
		Speaking of (co)limits, it is worthwhile noting that entities exhibiting universality, like colimits and limits, initial objects and terminal objects, can themselves be phrased entirely in terms of adjoint functors. Then, one of the advantages of this adjunction perspective is that the (co)limit of \textit{every} \textbf{J}-shaped diagram in $\textbf{C}$ can be defined all at once, rather than just taking the (co)limit of a particular $\textbf{J}$-shaped diagram $X: \textbf{J} \rightarrow \textbf{C}$. We will return to this in a later chapter.} 
	\end{proposition} 
Instead of proving this in the general case, we will show how it obtains in our special case of maps between orders, in which setting limits are infima and colimits are suprema.
\begin{proposition}
	\begin{enumerate}
		\item If $f: \mathcal{Q} \rightarrow \mathcal{P}$ has a right adjoint (i.e., is a left adjoint), then it preserves the suprema that exist in $\mathcal{Q}$. 
		\item If $g: \mathcal{P} \rightarrow \mathcal{Q}$ has a left adjoint (i.e., is a right adjoint), then it preserves the infima that exist in $\mathcal{P}$.  
	\end{enumerate}
\end{proposition}
\begin{proof}
	(Of (a), since (b) follows by duality) Assume $S = \{q_i\}_{i \in I}$ is a family of elements of $\mathcal{Q}$ with a supremum $\bigvee S$ in $\mathcal{Q}$. Claim: $f(\bigvee S)$ is the supremum in $\mathcal{P}$ of the family $\{f(q_i)\}_{i \in I}$, i.e., 
	\begin{equation*}
	f(\bigvee S) = \bigvee f(S). 
	\end{equation*} 
	But $f(\bigvee S) \leq p$ iff $\bigvee S \leq g(p)$ (since, by assumption, $f$ has a right adjoint, call it $g$). And this latter inequality holds iff for all $q_i \in S$, we have $q_i \leq g(p)$. But then we can again use the assumed adjoint relation $f \dashv g$, and see that this latter inequality will hold iff $f(q_i) \leq p$ for all $q_i \in S$, and this in turn will hold iff for all $t \in f(S)$, we have $t \leq a$. In sum, then, we have that $f(\bigvee_i q_i) \leq p$ if and only if $\bigvee_i f(q_i) \leq p$, or that $f$ preserves any suprema that exist in $\mathcal{Q}$.  
\end{proof}
But a dilation (erosion) was just \textit{defined} as an order-preserving map that commutes with colimits (limits).	So $\delta$ being a left adjoint suffices to tell us that $\delta$ must be a dilation (and dually, for an erosion).\par 
\item (\textbf{Adjoints Compose}) 
\begin{proposition} Given two dilations $\delta: \mathcal{L} \rightarrow \mathcal{M}, \delta':\mathcal{M} \rightarrow \mathcal{N}$ and two erosions $\epsilon: \mathcal{M} \rightarrow \mathcal{L},  \epsilon': \mathcal{N} \rightarrow \mathcal{M}$ such that $\delta \dashv \epsilon$ and $\delta' \dashv \epsilon'$, then their composition forms an adjunction $\delta' \circ \delta \dashv \epsilon \circ \epsilon'$. 
\end{proposition}
This exemplifies a general result in category theory, namely: 
	\begin{proposition}
		Left (right) adjoints are closed under composition,  i.e., given the adjunctions
		\begin{center} 
			\tikzset{
				,no line/.style={%
					,draw=none
					,commutative diagrams/every label/.append style={/tikz/auto=false}
				}
			}
			\begin{tikzcd}[column sep =large]
				\textbf{C} \arrow[bend left = 30]{r}[name=U]{F} & \textbf{D} \arrow{l}[name=L]{G} \arrow[bend left = 30]{r}[name=Q]{F'} \arrow[from=L, to=U, no line, pos=.5]{}{\perp} & \textbf{E} \arrow{l}[name=R]{G'} \arrow[from=R, to=Q, no line, pos=.5]{}{\perp},
			\end{tikzcd}
		\end{center}  
		the composite $F' \circ F$ is left adjoint to the composite $G \circ G'$: 
		\begin{center} 
			\tikzset{
				,no line/.style={%
					,draw=none
					,commutative diagrams/every label/.append style={/tikz/auto=false}
				}
			}
			\begin{tikzcd}[column sep =large]
				\textbf{C} \arrow[bend left = 30]{r}[name=U]{F' \circ F} & \textbf{E} \arrow{l}[name=L]{G \circ G'} \arrow[from=L, to=U, no line, pos=.5]{}{\perp}. 
			\end{tikzcd}
		\end{center}  
		In this way, arbitrarily long strings of adjoints can be produced. 
	\end{proposition}
Moreover, another fact from morphology follows from the facts that adjoints compose (and using LAPC and RAPL), namely that for dilations and erosions on the same complete lattice, if $\delta_j \dashv \epsilon_j$ forms an adjoint pair for every $j \in J$, then $(\bigvee_j \delta_j, \bigwedge_j \epsilon_j)$ is an adjunction. 
\item The ``opening" operator bounds an image on the left, while its ``closing" bounds it on the right, i.e., 
\begin{proposition} 
	\begin{equation*} 
	\delta \epsilon \leq \text{id} \leq \epsilon \delta.
	\end{equation*}
\end{proposition} 
This is immediate from the unit natural transformation, $\text{id} \leq \epsilon \delta$, and the counit $\delta \epsilon \leq \text{id}$. 
\item (\textbf{Fixed Point Formulae})\index{fixed point!formulae} 
\begin{proposition} 
	$\delta \epsilon \delta = \delta$ and $\epsilon \delta \epsilon = \epsilon$.
\end{proposition} 
The unit and counit maps, satisfying the triangle identities, give the following general ``fixed point formulae" result underlying the above: 
\begin{proposition}
	If $\mathcal{P}$ and $\mathcal{Q}$ are posets and $F: \mathcal{P} \rightarrow \mathcal{Q}$ and $G: \mathcal{Q} \rightarrow \mathcal{P}$ form a (monotone) Galois connection (adjunction), with $F \dashv G$, then the following \textit{fixed point formulae} will hold for $F$ and $G$: 
	\begin{equation*}
	FGF = F \text{ and } GFG = G.
	\end{equation*}
\end{proposition} 
\begin{proof}
	The triangle identities give $F(p) \leq FGF(p) \leq F(p)$ for all $p \in \mathcal{P}$, so $F = FGF$. The second formula follows similarly. 
\end{proof}
\end{enumerate} 
In the last two items, we saw how the unit and counit maps determine two important endomaps, namely $\delta \circ \epsilon$ (``opening")\index{opening} and $\epsilon \circ \delta$ (``closing")\index{closing}. The presence of unit and counit further give us the fixed point formulae, which translates to the morphologically-significant fact, 
\begin{equation}
\delta \epsilon \delta = \delta \hspace*{2em} \text{ and } \hspace*{2em} \epsilon \delta \epsilon = \epsilon. 
\end{equation} 
This ``stability" property of openings and closings means, in terms of the interpretation of such operations as filters, that they effectively ``complete their task" (unlike many other filters, where repeated applications can involve further modifications of the image, with no guarantee of the outcome after a finite number of iterations). In general, the above fixed point formula further entails, in particular, that $\epsilon \delta$ and $\delta \epsilon$ are each idempotent. Thus, altogether, the composite monotone map $\epsilon \delta$, for its part, has the properties that 
\begin{itemize}
	\item $p \leq \epsilon \delta (p)$, and 
	\item $\epsilon \delta \epsilon \delta (p) = \epsilon \delta (p)$. 
\end{itemize}
But this is exactly to say that $\epsilon \delta$ is a \textit{closure} operator, in the following general sense. 
\begin{definition} 
	A \textit{closure operator}\index{closure} on a poset (typically some poset of subobjects, e.g., on the powerset poset) $\mathcal{P}$ is an endomap $K: \mathcal{P} \rightarrow \mathcal{P}$ such that
	\begin{enumerate}
		\item for each $p \leq p' \in \mathcal{P}$, $K(p) \leq K(p')$ (monotonicity);  
		\item for each $p \in \mathcal{P}$, $p \leq K(p)$ (extensivity); 
		\item for each $p \in \mathcal{P}$, $K(K(p)) = K(p)$ (idempotence).
	\end{enumerate}
\end{definition} \noindent 
There is an important dual notion to closure, called the kernel operator (or dual closure), where this is an endomap that, like $K$, arises from a Galois connection, and is both monotone and idempotent, yet satisfies the dual of the extensivity property.  
\begin{definition}
	A \textit{kernel operator} (or \textit{dual closure})\index{closure!dual} is an endomap $L$ satisfying 
	\begin{enumerate}
		\item for each $p \leq p' \in \mathcal{P}$, $L(p) \leq L(p')$ (monotonicity);  
		\item for each $p \in \mathcal{P}$, $L(p) \leq p$ (contractivity);
		\item for each $p \in \mathcal{P}$, $L(L(p)) = L(p)$ (idempotence). 
	\end{enumerate}
\end{definition}
In short, these notions are all part of a much more general story, namely that for a Galois connection or adjunction on posets such as $\delta \dashv \epsilon$ as above, the composite $\epsilon \circ \delta$ will automatically be monotone, extensive, and idempotent, i.e., a closure operator on the underlying poset (or lattice) $\mathcal{P}$; dually, $\delta \circ \epsilon$ will be monotone, contracting, and idempotent, i.e., a kernel operator on $\mathcal{Q}$. \par 	
Before leaving this example, we will explore a few last notions via morphology. We will let the induced kernel operator, called ``opening" by the morphology community, be denoted $\phi = \delta \epsilon$, while $\kappa = \epsilon \delta$ will denote the induced closure (or ``closing" operator, to be consistent with the mathematical morphology literature). In the binary case, opening and closing are typically defined, respectively, as 
	\begin{equation}
	X \text{o} B = (X \ominus B) \oplus B = \bigcup \{B_p | p \in E, B_p \subseteq  X\}
	\end{equation}
	\begin{equation}
	X \text{\textbullet} B = (X \oplus B) \ominus B.
	\end{equation}
	Morphological closing\index{closing} is just dilation\index{dilation} (by some $B$) followed by erosion of the result by $B$, while morphological opening\index{opening} is the erosion\index{erosion} (by some $B$) followed by dilation of the resulting image by $B$. Closing acts to fill out narrow holes. In terms of translations with the structuring element, the opening of an image $A$ by $B$ is the complement of the union of all translations of $B$ that fall outside (do not overlap) $A$. As extensive (i.e., larger than the identity mapping), closings of an image are generally ``larger" than the original image. Opening, for its part, acts to remove noise, narrow connections between regions, and parts of objects, generally attenuating ``peaks" and other small protrusions or components. If you have a note where the writing appears to be growing tiny roots from its edges, opening effectively acts to remove these outer leaks at the boundary, rounding the edges. In terms of translations with the structuring element, the opening of $A$ by $B$ is the union of all translations of $B$ that fit completely within $A$. As anti-extensive, openings of an image are generally ``smaller" than the original. 
	\begin{exercise}
		Composing dilations and erosions, we found the composite operations of opening $(\phi = \delta \epsilon$) and closing $(\kappa = \epsilon \delta$), which were moreover idempotent. Further composing openings and closings with one another (e.g., $\kappa \circ \phi$), how many more distinct operations can we produce? Describe, in terms of their effect on images, at least one of these ``image filters." Finally, consider how the composite operators must be related to one another.   
	\end{exercise}
	\noindent
	\textit{Solution}: By (alternately) composing openings $\phi (= \delta \epsilon)$ and closings $\kappa (= \epsilon \delta)$, we can obtain four new filters in total, each four of which are idempotent. 
	\begin{enumerate}
		\item closing-after-opening: $\kappa \phi$
		\item opening-after-closing: $\phi \kappa$
		\item opening-after-closing-after-opening: $ \phi \kappa \phi$
		\item closing-after-opening-after-closing: $\kappa \phi \kappa$.
	\end{enumerate} 
	You can easily convince yourself that no other operator can be obtained by further composition with any combination of $\phi$s and $\kappa$s. Any attempt to produce further new operations by pre- or post-composing the above four with $\kappa$ or $\phi$ will just reduce back to one of those four, by the idempotence of these operators (together with the idempotence of $\kappa$ and $\phi$ themselves). \par 
	These composites are used in the course of various image processing tasks, such as ``smoothing" an image or performing image segmentation. An example of $\kappa \phi$ is given by the following:\footnote{This image is taken from \cite{bobick_binary_2014}.}
	\begin{center}
		\includegraphics*[scale=0.25]{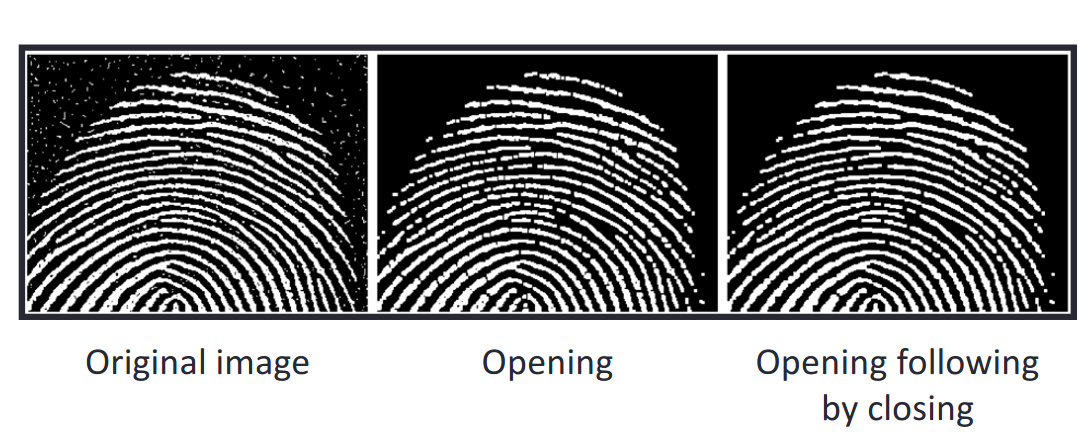}
	\end{center}
	In terms of the relations between these four (and the original opening and closing operators as well), it is easy to show that 
	\begin{equation*}
	\phi \leq \phi \kappa \phi \leq \{ \begin{smallmatrix}
	\kappa \phi \\ 
	\phi \kappa 
	\end{smallmatrix} \} \leq \kappa \phi \kappa \leq \kappa,  
	\end{equation*}
	and moreover $\phi \kappa \phi$ will be the greatest filter smaller than $\phi \kappa \wedge \kappa \phi$, while $\kappa \phi \kappa$ will be the smallest filter greater than $\phi \kappa \vee \kappa \phi$.
\end{example} 
\subsection{Adjunctions through Modalities}\label{section: modalities}
It may help to get an ever better handle on adjunctions by looking at another application of such notions. At least since Aristotle's attempt to understand certain statements containing the words ``necessary"\index{necessity} and ``possible,"\index{possibility} philosophers and logicians have been interested in the ``logic" of different operators describing different ``ways of being true." Modal logic began as the study of \textit{necessary} and \textit{possible} truths, but in at least the last 100 years it has been recognized that modalities\index{modalities} abound in both natural and formal languages; so these days modal logic is more commonly regarded as the much broader study of a variety of constructions that modify the truth conditions of statements (which includes, most notably, statements concerning knowledge, belief, ethics, temporal happenings, how computer programs behave, and the dynamical properties of state transitions). There are in fact a number of close connections between erosion and dilation and the modal operators $\Box$ of necessity and $\Diamond$ of possibility, respectively, and similarly one can define modal operators in terms of adjunctions. \par 
The following realizes these ideas more concretely. 
\subsubsection{On What is Not}
\begin{example}	
Since at least the time of one of the first Western philosophical texts, the Parmenides (around 500 B.C.E), the nature of \textit{negation}\index{negation} has been on people's minds. This includes a number of issues, such as
	\begin{itemize}
		\item  would a complete description of \textit{what is} need to include any description of \textit{what is not}? In other words, what is the ``ontological status" of the negated entities or negative states of affairs?\footnote{One position, in this context, might articulate the view that everything is what it is---as the individual thing it is---only on account of how it is \textit{not} some other things, and accordingly try to take very seriously the idea that ``all determination arises from a kind of negation." An opposing position might argue that negations always just describe \textit{privations}, and a complete and accurate description of reality would not need to involve mention of any ``negative entities."}
		\item when is the negation of a negation (negated entity) the identity (the original entity)? 
	\end{itemize}
	One might further motivate such concerns as follows. One might try to argue that, philosophically, holes, shadows, fissures, boundaries---and other such ``negative" or derivative entities---seem somehow less real or fundamental than (or at least not to be on the same footing with) the ``positive" objects that produce or surround or support them. At the very least, this sort of observation seems to have some validity in that it does seem somehow more difficult to supply identity criteria for holes, for instance, compared to ordinary material objects (for holes appear to be \textit{made of nothing}), or even to speak of what holes are (what are the \textit{parts} of a hole?). \par     
	Today, we are most accustomed to thinking of negation as a linguistic or logical operator on a language, where the operation leads from an expression to the contradictory expression. Typically, the expression in question is a sentence or a part of a sentence. But we might attempt to regard negation, more broadly, as an operation that can also take place on larger wholes, on entire structures or theories. Moreover, one might argue that, however one approached it, the ``right" understanding (and description) of negation would need to capture, above all, the relation or dependence between what (the structure) is being negated and the result of this operation of negation. \par   
	Such a perspective on negation is arguably exemplified in the facts that (i) negation is a contravariant functor on a particular category (to itself), and (ii) this functor has special relations to itself, in that it is adjoint, in fact self-adjoint. This perspective, developed formally in the following example, might even suggest, informally, that one think of the action of certain contravariant functors as a generalized sort of negation of structures. \par 
	First recall that a \textit{Heyting algebra}\index{algebra!Heyting} $H$ is a poset with all finite products and coproducts, and that is moreover cartesian closed. Another way of describing such an $H$ is as a distributive lattice with a least element $0$ and a greatest element $1$, expanded with an operation $\Rightarrow$, where this means that for any two elements $p, q$ of the lattice, there exists an exponential $q^p$, usually written 
	\begin{equation*}
	p \Rightarrow q.
	\end{equation*}
	This operation is characterized by the adjunction \index{adjunction} 
	\begin{equation*}
	r \leq (p \Rightarrow q) \text{ iff } r \wedge p \leq q. 
	\end{equation*}	  
	In other words, $\Rightarrow$ is a binary operation on a lattice with a least element, such that for any two elements $p, q$ of the lattice, $\max \{r \hspace*{0.25em}| \hspace*{0.25em} r \wedge p \leq q \}$ exists (where this latter set contains an element greater than or equal to every one of its elements, and such a least upper bound for all those elements $r$ where $r \wedge p \leq q$ is what is denoted by `$p \Rightarrow q$'). Another way to think of this $p \Rightarrow q$ is in the setting of the propositional calculus, where it is the weakest condition needed for the inference rule of modus ponens to hold, i.e., to enforce that from $p \Rightarrow q$ and $p$ we can infer $q$.\par 
	Heyting algebras\index{algebra!Heyting} serve as models for intuitionistic propositional calculus, and in this setting, with variables regarded as propositions, $\wedge$ as 'and', $\vee$ as 'or', and $\Rightarrow$ as implication, we can also define \textit{negation} of $p$ as 
	\begin{equation*}
	\neg p := (p \Rightarrow 0). 
	\end{equation*}  
	On account of how $\Rightarrow$ is defined, we can rewrite this as 
	\begin{equation*}
	q \leq \neg p \text{ iff } q \wedge p = 0,
	\end{equation*}
	revealing $\neg p$ to be the union of all those $q$ whose meet with $p$ in the lattice is $0$, the least element. \par 
	In any Heyting algebra $H$, we have not just that 
	\begin{equation*}
	p \leq \neg \neg p,
	\end{equation*}
	but also that 
	\begin{equation*}
	p \leq q \text{ implies } \neg q \leq \neg p.
	\end{equation*} 
	But this reveals how negation is just a contravariant functor from the Heyting algebra\index{algebra!Heyting} to itself! More explicitly, 
	\begin{proposition} 
		$\neg$ is a functor $\neg: H \rightarrow H^{op}$ (and also $\neg: H^{op} \rightarrow H$). This functor is, moreover, adjoint to itself, since $p \leq \neg q$ iff $q \leq \neg p$. 
	\end{proposition} \noindent 
Let us spell out the above, specifically its self-adjointness, more explicitly. The first inequality, $p \leq \neg \neg p$, is immediate from $q \leq \neg p \text{ iff } q \wedge p = 0$, using the further fact that, for any Heyting algebra,\index{algebra!Heyting} $p \wedge \neg p = 0$ (this follows from the adjunctive definition of $\Rightarrow$ together with the definition of negation as $\neg p = (p \Rightarrow 0)$). For the second, suppose $p \leq q$ in $H$. Then, $p \wedge \neg q \leq q \wedge \neg q$ and the right-hand side of this inequality is $0$. So $p \wedge \neg q = \neg q \wedge p = 0$, and so by $q \leq \neg p \text{ iff } q \wedge p = 0$, we have that $\neg q \leq \neg p$. Incidentally, this moreover shows that
	\begin{equation*}
	\neg p = \neg \neg \neg p
	\end{equation*} which we might call the `$1 = 3$' fact. This is a result of the contravariant functoriality of $\neg$ and that $p \leq \neg \neg p$. For, suppose $p \leq \neg \neg p$. Then, by the contravariant functoriality inequality, we also have that $\neg \neg \neg p \leq \neg p$. And since $p \leq \neg \neg p$ holds for all $p$ in $H$, it holds in particular for $\neg p$. Thus, $\neg p \leq \neg \neg \neg p$, giving the other side of the equality, and so, altogether, $1 = 3$. \par 
	That $\neg$, as a functor $H \rightarrow H^{op}$ and also $H^{op} \rightarrow H$, is adjoint to itself means that for all $p \in H$ and $q \in H^{op} (= H)$, we have the two way rule 
	\begin{center} 
		\begin{equation*}
		\frac{\neg p \leq^{op} q}{p \leq \neg q}
		\end{equation*}
	\end{center} 
	or 
	\begin{center} 
		\begin{equation*}
		\frac{\neg p \geq q}{p \leq \neg q}
		\end{equation*}
	\end{center} 
	where the top (holding in $H^{op}$) holds if and only if the bottom (holding in $H$) does. Another way of seeing the truth of the fact that for all $x \in H$, $x \leq \neg \neg x$, notice that as an adjoint, letting $q = \neg p$, we must have 
	\begin{center} 
		\begin{equation*}
		\frac{\neg p \geq \neg p}{p \leq \neg \neg p}
		\end{equation*}
	\end{center} 
	and since the top is always true, the bottom must be as well. \par 
	An adjoint is a kind of generalized inverse, and as such, an adjunction describes a kind of loosening or weakening of the notion of equivalence. In the present situation, asking when it is the case that $p = \neg \neg p$ (when the `do nothing' functor is equivalent to applying the negation functor twice), is like asking when the above adjunction happens to be an equivalence (isomorphism). If the relations $\leq$ are replaced by $=$, then we get isomorphisms. This captures the following well-known relation between Heyting algebras and Boolean algebras: 
	\begin{proposition}
		A Heyting algebra $H$ is Boolean\index{algebra!Boolean} (i.e., $\neg \neg x = x$ for all $x \in H$) if and only if the above adjunction is an equivalence. 
	\end{proposition} \noindent 
This is a stricter requirement, and in general we need not have $x = \neg \neg x$. This requirement is something one might not always want to impose, and this is in fact one of the merits or utility of working with the more general Heyting algebras.\index{algebra!Heyting} For any topological space $X$, the set $\mathscr{O}(X)$ of open sets of $X$ forms a Heyting algebra; for instance, the opens in the real line accordingly form a Heyting algebra, but one that is not Boolean,\index{algebra!Boolean} since the complement of an open set is not necessarily open. \par 
	We could go on to dualize things and describe a dual notion, namely that of \textit{co-Heyting algebras},\index{algebra!co-Heyting} which support a corresponding but different notion of negation. The utility of considering such things can be motivated with another problem related to natural language. In many natural languages, for instance in English, one often has resort to forms of negation that do not seem to be captured by, or behave as, the single negation operator of classical logic. Suppose someone is described to you as `not honest', after they act in a particular way in a particular situation. This is not necessarily to say that they are `dishonest.' In natural language, we can deny that a person is honest in at least two distinct ways: (i) by asserting that someone is not honest (negating the predicable `to be honest'); or (ii) by asserting that they are dishonest (negating the adjective `honest'). It is easy to appreciate, intuitively, how the second, e.g., `dishonest', is a stronger form of negation than the first (`not honest'). Moreover, suppose your friend Abe is someone you would be willing to describe as `not dishonest.' This does not seem to convey the same thing as describing Abe as `honest.' Finally, while we expect it to be the case that Abe is either honest or not honest, it seems plausible to assert, as well, that he is neither honest nor dishonest.  We would like to know how to capture these observations more formally, and describe formal relations between the different forms of negation, as applied to natural language. \par
	More generally, compare the sentence
	\begin{quote}
		It is false that not $p$,
	\end{quote}  
	with the sentence 
	\begin{quote}
		It is not false that $p$. 
	\end{quote}
	These sentences clearly do not say the same thing. The first indicates the \textit{necessity} of\index{necessity} $p$, while the second indicates its \textit{possibility}.\index{possibility} We can make sense of this in the context of a particular algebra called a bi-Heyting algebra, by interpreting the ``it is false" in such sentences as the Heyting negation and the ``not" as the corresponding co-Heyting negation. This setting will also allow us to define modal operators in terms of pairings of both negations. \par 
	A \textit{co-Heyting algebra}\index{algebra!co-Heyting} is a poset whose dual is a Heyting algebra. Unpacking this, we can equivalently observe that a co-Heyting algebra will be a (bounded) lattice expanded with a binary operation $\smallsetminus$ such that for every $r, p, q$, we have the adjunction rule 
	\begin{equation*}
	(p \smallsetminus q) \leq r \text{ iff } p \leq q \vee r. 
	\end{equation*}
	In other words, $p \smallsetminus q = \bigwedge \{r | p \leq q \vee r \}$. 
	A corresponding unary negation operation $\sim$ can then be defined by 
	\begin{equation}
	\sim p := (1 \smallsetminus p). 
	\end{equation}	  
	It follows that we have the following adjunction rule for this `negation' $\sim$: 
	\begin{center} 
		\begin{equation*}
		\frac{\sim p \leq q}{1 =p \vee q}
		\end{equation*}
	\end{center} 
	Notice how, similar to how we have $p \wedge \neg p = 0$ in a Heyting algebra (though not necessarily $p \vee \neg p = 1$), in a co-Heyting algebra we have 
	\begin{equation*}
	p \vee \sim p = 1.
	\end{equation*}
	Similar to how, in a Heyting algebra, the negation $\neg$ is order-reversing and satisfies $x \leq \neg \neg x$, in a co-Heyting algebra the negation $\sim$ is also order-reversing and it satisfies $\sim \sim x \leq x$. \par   
	Moreover, just as, in a Heyting algebra, $p \vee \neg p$ is not necessarily the ``top" (\textit{true}), so in a co-Heyting algebra, $p \wedge \sim p$ is not necessarily the ``bottom" (\textit{false}). In particular, then, in a co-Heyting algebra, we can thus define the generally non-trivial notion of the \textit{boundary}\index{boundary} of $p$, as 
	\begin{equation*}
	\partial p : = p \wedge \sim p. 
	\end{equation*}
	Incidentally, this recovers, in purely algebraic terms, the geometric notion of \textit{boundary}, thus suggesting certain deep connections between logic and geometry. 
	Bi-Heyting\index{algebra!bi-Heyting} algebras--a bounded distributive lattice that is both a Heyting and a co-Heyting algebra---can accordingly be deployed to shed further light into some of these deep connections between logic and geometry. \par 
	We can thus work in the context of a bi-Heyting algebra and combine these negations to form our modal operators.\index{modalities} In particular, we will let
	\begin{equation*} 
	\Diamond p = \sim \neg p,
	\end{equation*} read as `possibly $p$'. This begins to suggest how we might formalize the fact that, for instance, John not being dishonest does not let us conclude, in general, that John is honest, but rather only that he is \textit{possibly} honest.\index{possibility} Similarly, we will let 
	\begin{equation*} 
	\Box p = \neg \sim p,
	\end{equation*} read as `necessarily $p$'.\index{necessity} \par 
	Given such definitions, we can show that
	\begin{proposition}
		We have the adjunction\index{adjunction} $\Diamond \dashv \Box$. 
	\end{proposition} 
	\begin{proof}
		To show that $\Diamond \dashv \Box$, we need only verify some equivalences, each of which more or less follows automatically from definitions. By definition, $\Diamond = \sim \neg$ and $\Box = \neg \sim$, so the adjunction just says that $\sim \neg \dashv \neg \sim$. But 
		\begin{equation*}
		\sim \neg p \leq q
		\end{equation*}
		just says that 
		\begin{equation*}
		1 \leq \neg p \vee q
		\end{equation*}
		which is of course equivalent to 
		\begin{equation*}
		1 \leq q \vee \neg p, 
		\end{equation*}
		which itself is just to say that 
		\begin{equation*}
		\sim q \leq \neg p.
		\end{equation*}
		Using this last inequality, and the definition of $\neg$, we have 
		\begin{equation*}
		\sim q \wedge p \leq 0
		\end{equation*}
		or, equivalently, 
		\begin{equation*}
		p \wedge \sim q \leq 0. 
		\end{equation*}
		This last line can be written 
		\begin{equation*}
		p \leq (\sim q \Rightarrow 0),
		\end{equation*}
		which is the same as saying that 
		\begin{equation*}
		p \leq \neg \sim q. 
		\end{equation*}
		Altogether, this string of equivalences shows that 
		\begin{equation}
		\sim \neg p \leq q \text{ iff } p \leq \neg \sim q,
		\end{equation}
		which is precisely what is needed to show that $\sim \neg \dashv \neg \sim$ (or $\Diamond \dashv \Box$). 
	\end{proof}
	Working in a bi-Heyting algebra and using the above definitions of $\Box$ and $\Diamond$, we can moreover consider repeated applications of these operators. As a matter of simplifying computations involving repeated applications of the operators,\index{modalities} we will define $\Box_0 = \Diamond_0 = $ Id, and
	\begin{equation*}
	\Box_{n+1} := \neg \sim \Box_n, \hspace*{2em} \Diamond_{n+1} :=  \sim \neg \Diamond_n,
	\end{equation*}
	where $\Box_n$ is of course just the result of iterating ($n$ times) the composition of $\neg \sim$, and $\Diamond_n$ by iterating ($n$ times) the composition of $\sim \neg$. $\Box_n$ and $\Diamond_n$ are clearly both order-preserving, for all $n$, as is evident from the double fact that, in a Heyting algebra, $\neg$ is order-reversing and satisfies $p \leq \neg \neg p$ for all $p$, and that, in a co-Heyting algebra, $\sim$ is order-reversing and $\sim \sim p \leq p$. Moreover, we have 
	\begin{enumerate}
		\item $\Box_{n+1} \leq \Box_n \leq \text{Id} \leq \Diamond_n \leq \Diamond_{n+1}$ for all $n$; and 
		\item $\Diamond_n \dashv \Box_n$ for all $n$. 
	\end{enumerate}  
	\begin{proof}
		1. First, we know that, for any $p$ in the bi-Heyting algebra, we have that $\neg p \leq \sim p$. From this, taking $p = \sim p$, we have that 
		\begin{equation*}
		\neg \sim p \leq \sim \sim p
		\end{equation*}
		and since $\sim \sim p \leq p$, this gives that 
		\begin{equation*}
		\neg \sim p \leq p. 
		\end{equation*}
		Moreover, $p \leq \neg \neg p$, and, applying the result $\neg p \leq \sim p$ again, now to $p = \neg p$, we get that 
		\begin{equation*}
		\neg \neg p \leq \sim \neg p.
		\end{equation*}
		Altogether, this gives that 
		\begin{equation*}
		\neg \sim p \leq p \leq \sim \neg p. 
		\end{equation*} 
		By definition, then, this reads as 
		\begin{equation*}
		\Box p \leq \text{ Id} p \leq \Diamond p, 
		\end{equation*}
		and further iterating this, letting $p = \Box_n p$, and then $\Diamond_n p$, gives the main result.  
	\end{proof}
	\begin{proof}
		2. We want that $\Diamond_n \dashv \Box_n$ for all $n$. But since adjoints compose, this follows from the preceding result, by iterating.
	\end{proof}
	In a bi-Heyting algebra\index{algebra!bi-Heyting} where countable suprema and infimi exist satisfying 
	\begin{equation*}
	\frac{b \leq \bigwedge_n a_n}{b \leq a_n \text{ for all } n}
	\end{equation*}
	and 
	\begin{equation*}
	\frac{\bigvee_n a_n \leq b}{a_n \leq b \text{ for all } n}
	\end{equation*}
	we can further define 
	\begin{definition}
		\begin{equation}
		\Box p = \bigwedge_n \Box_n p, 
		\end{equation}
		\begin{equation}
		\Diamond p = \bigvee_n \Diamond_n p. 
		\end{equation}
	\end{definition} \par 
	Then, for a bi-Heyting algebra\index{modalities} that has countable suprema and infima satisfying the above rules, the modal operators $\Box$ and $\Diamond$ are such that $\Box p$ will be the largest complemented $x$ such that $x \leq p$, and $\Diamond p$ will be the smallest complemented $x$ such that $p \leq x$. This moreover realizes the fact that $\Box$ and $\Diamond$ are both order-preserving and that $\Box p \leq p \leq \Diamond p$.\par 
	More generally, for any bounded distributive lattice and operators $\Box$ and $\Diamond$ defined just as above, the same properties will hold of $\Box$ and $\Diamond$. In particular, it can be shown that $\Diamond \dashv \Box$, and thus
	\label{5 inequalities}
	\begin{enumerate}
		\item $\Box \leq \text{ id } \leq \Diamond$,
		\item $\Box \Box  = \Box , \Diamond \Diamond = \Diamond$, 
		\item $\text{id } \leq \Box \Diamond$, 
		\item $\Diamond \Box \leq \text{ id}$,
		\item $\Diamond (\phi \wedge \Box \psi) = \Diamond \phi \wedge \Box \psi$. 
	\end{enumerate}
	\subsubsection{Adjoint Modalities in the Qua Category}
	As a particular realization of some of these ideas,\index{modalities} recall the \textit{interpretation} functor on the \textit{qua}\index{category!qua} category $\textbf{Qua}$ from Example \ref{qua}.\footnote{Again, this material on the \textit{qua} category is ultimately derived from \cite{la_palme_reyes_models_1999}; the reader who desires to pursue these matters further than what is discussed here can find many more interesting details in that paper.} Applying all of this to the particular subcategory $\textbf{A}$ from earlier, an interpretation amounts to a set $X$ together with a set of predicates of $X$, where by `predicate of $X$' is just meant a family $\{\phi_A \}_{A \in Ob(\textbf{A})}$ of subsets of $X$ with the functorial property
	\begin{equation*}
	\text{ if } x \in \phi_A (\text{ i.e., } x \in_A \phi) \text{ and } A' \rightarrow A \in \textbf{A}, \text{ then } x \in \phi_{A'}. 
	\end{equation*}
	And since $\textbf{A}$ was just the comma category $\boxed{\text{ a scf }} \downarrow \textbf{CN}$, so that all aspects are of the form $\boxed{\text{a scf}}qua \boxed{\text{ B }}$, an interpretation of $\textbf{A}$ just associates to every aspect the same set $X$. Notice that we can collect together such families of predicates (as subfunctors\index{functor!sub} of $X$) into the set $\mathcal{P}(X)$ of all predicates of $X$. These predicates in fact form a bounded distributive lattice with 2 negations---in fact, a bi-Heyting algebra\index{algebra!bi-Heyting}---as
	\begin{equation*}
	(\mathcal{P}(X), \leq, \vee, \wedge, 1, 0, \neg, \sim),
	\end{equation*} 
	where there is the natural ordering 
	\begin{equation*}
	\phi \leq \psi \text{ iff } \forall A \in \textbf{A}, \forall x \in X, x \in_A \phi \Rightarrow x \in_A \psi.
	\end{equation*} 
	Also, as expected, we have 
	\begin{equation*}
	x \in_A (\phi \vee \psi) \text{ iff } x \in_A \phi \text{ or } x \in_A \psi 
	\end{equation*}
	and 
	\begin{equation*}
	x \in_A (\phi \wedge \psi) \text{ iff } x \in_A \phi \text{ and } x \in_A \psi . 
	\end{equation*}
	$0$ (or $\bot$) is the predicate `false', the bottom element of the order, while $1$ (or $\top$) is the predicate `true', the top element of the order. Given a predicate $\phi$, we define two negations, $\neg \phi$ and $\sim \phi$, which are in fact two new predicates (i.e., have the requisite functorial property):\index{modalities} 
	\begin{equation*}
	x \in_A \neg \phi \text{ iff } \forall A' \rightarrow A \in \textbf{A} \hspace*{0.5em} x \notin_{A'} \phi
	\end{equation*}
	and 
	\begin{equation*}
	x \in_A \sim \phi \text{ iff } \exists A \rightarrow A' \in \textbf{A} \hspace*{0.5em} x \notin_{A'} \phi
	\end{equation*}
	Applied to the predicate `honest', for instance, we have the natural reading: `$\neg$ honest' as `dishonest', and `$\sim$ honest' as `not honest'. 
	\par 
	As we saw in our earlier discussion of these same negations, we have the following adjunctions for our negations (which hold for arbitrary properties $\phi, \psi$): 
	\begin{enumerate}
		\item \begin{equation*}
		\frac{\psi \leq \neg \phi}{\psi \wedge \phi = 0}
		\end{equation*}
		\item \begin{equation*}
		\frac{\sim \phi \leq \psi}{1 = \psi \vee \phi}
		\end{equation*}
	\end{enumerate}
	Notice also that for every property $\phi$, it is a consequence of the above that $\phi \wedge \neg \phi = 0$ and $\phi \vee \sim \phi = 1$. However, $\phi \wedge \sim \phi$ need not be $0$, and similarly, $\phi \vee \neg \phi$ is not necessarily $1$. \par 
	Applying all this, suppose we want to model a discussion and decision regarding Abe's honesty. We assume the aspects considered relevant to Abe's honesty have been agreed upon, and likewise, agreement has been achieved on Abe's honesty under each of the relevant aspects, i.e., for every aspect $A \in \textbf{A}$, we know whether $Abe \in_A \textit{ honest}$ or $Abe \notin_A \textit{ honest}$. If Abe fails to be honest under every sub-aspect of one of the aspects, say $F = $ family man, then we would say that `Abe is dishonest' under that aspect, i.e., \textit{qua} family man. By contrast, we would say that `Abe is not honest' under the aspect $F$ precisely when he fails to be honest with respect to one of the super-aspects of $F$. The ultimate judgment regarding Abe's honesty is then obtained by restricting to the global aspect, $G$ (or `\textit{qua} scf'), where this means that Abe is `honest', `not honest', or `dishonest' precisely when $Abe \in_G honest$, $Abe \in_G \sim honest$, $Abe \in_G \neg honest$, respectively. In more detail, 
	\begin{quote}
		$Abe \in_G honest$ iff $\forall A$ Abe $\in_A honest$, i.e., `Abe is honest iff Abe is honest under any aspect'. 
	\end{quote}
	\begin{quote}
		$Abe \in_G \sim honest$ iff $\exists A$ Abe $\notin_A honest$, i.e., `Abe is not honest iff Abe fails to be honest under at least one of the aspects'. 
	\end{quote}
	\begin{quote}
		$Abe \in_G \neg honest$ iff $\forall A$ Abe $\notin_A honest$, i.e., `Abe is dishonest iff Abe fails to be honest under every one of the aspects'. 
	\end{quote}
	If Abe is honest, then he cannot be dishonest (and conversely), i.e., $\phi \wedge \neg \phi = 0$. However, that $\phi \vee \neg \phi$ is not necessarily $1$ means, of course, that it is not always the case that Abe is either honest or dishonest. One scenario in which this might occur would be where Abe is honest under the aspect $S$ (Abe is honest \textit{qua} student) but fails to be honest in all the other aspects. This reveals how the negation is not, in general, Boolean. Yet we can note that $\sim$ is Boolean globally. This means that for the `global aspect', we will have that Abe is either honest or not honest, but not both, i.e., $honest \vee \sim honest = 1$ and $honest \wedge \sim honest = 0$. However, it may occur that Abe is both honest and not honest under the very same aspect (as long as this is not the global aspect).\par 
	Thus, notice that even though Abe is `not honest' under a particular aspect precisely when he is not honest under every aspect, for all aspects other than the global one, there is a difference between `not honest' under an aspect and \textit{failing} to be honest under that aspect, for the former has the functoriality property: if Abe is not honest under aspect $A$, then he is not honest under any subaspect $A' \rightarrow A$. Failing to be honest under a given aspect, by contrast, is simply the absence of Abe's honesty under that aspect. Absence of honesty under an aspect is not functorial, so there may be an absence of Abe's honesty under an aspect and Abe's honesty under another.  \par 
	Suppose Abe is not dishonest under an aspect. What can we conclude from this? Not that Abe is honest; rather, only that he is \textit{possibly} honest. Abe is not dishonest under a given aspect precisely when he is honest under at least one aspect, as the following string of equivalences reveal: 
	\begin{equation*}
	\frac{\frac{Abe \in_A \sim \neg honest}{\exists A \rightarrow A' \text{ } Abe \notin_{A'} \neg honest}}{\frac{\exists A \rightarrow A' \exists A'' \rightarrow A' \text{ } Abe \in_{A''} honest}{\exists A' \text{ } Abe \in_{A'} honest}}
	\end{equation*}
	Restricting to the global level $A = G$, gives 
	\begin{equation*}
	\frac{Abe \in \sim \neg h}{\exists A' Abe \in_{A'} h},
	\end{equation*}
	or `Abe is not dishonest iff Abe is honest under at least one aspect'. Since this works for any predicate, this suggests we define our modal operator\index{modalities}
	\begin{equation*}
	\Diamond \phi : = \sim \neg \phi,
	\end{equation*}
	read as `possibly $\phi$'.\index{possibility} \par \noindent 
	Similarly, we can calculate 
	\begin{equation*}
	\frac{\frac{Abe \in_A \neg \sim honest}{\forall A' \rightarrow A \text{ } Abe \notin_{A'} \sim honest}}{\frac{\forall A' \rightarrow A \exists A' \rightarrow A'' \text{ } Abe \in_{A''} honest}{\forall A' \text{ } Abe \in_{A'} honest}}.
	\end{equation*}
	Again, letting $A = G$ the global aspect, this becomes 
	\begin{equation*}
	\frac{Abe \in \neg \sim honest}{Abe \in honest},
	\end{equation*}
	i.e., Abe is necessarily honest under a given aspect iff he is honest under every aspect. In other words, $x \in_A \Box \phi$ iff for all aspects $A'$, $x \in_{A'} \phi$ iff $x \in_G \phi$.  
	This suggests we define a further modal operator
	\begin{equation*}
	\Box \phi : = \neg \sim \phi,
	\end{equation*}
	read as `necessarily $\phi$'.\index{necessity} Observe that $\Diamond \phi$ and $\Box \phi$, defined thus, are clearly themselves predicates of $X$. Moreover, they are themselves adjoint, $\Diamond \dashv \Box$, and as such satisfy the (in)equalities we isolated in \ref{5 inequalities}.
\subsubsection{Adjoint Modalities in Graphs}
	The modalities\index{modalities} discussed in the last few sections are particularly well-illustrated and tangible in the context of graphs and their subgraphs.\index{graph} If we have a directed multi-graph $G$, the lattice of subgraphs of $G$ constitutes a bi-Heyting algebra,\index{algebra!bi-Heyting}\footnote{This is actually a special case of something that will be discussed in a later chapter, namely that \textit{any presheaf topos is bi-Heyting}. The category of (multi-)graphs, we will see, can be represented as $\textbf{Set}^{\textbf{C}^{op}}$ for $\textbf{C}$ the category of two objects and two non-trivial morphisms between those objects; moreover, the lattice of subobjects of any object of a bi-Heyting topos is itself a (complete) bi-Heyting algebra.} where a subgraph $X$ of $G$ is a directed multi-graph, i.e., consists of subsets $X_0$ of the vertices $G_0$ of $G$ and a subset $X_1$ of the edges $G_1$ of $G$, such that every edge in $X_1$ has both its source and target in the vertex set $X_0$. It is clear that we can take unions and intersections of subgraphs, but what it means to take ``complements" is not as evident. The set-theoretical complement $c(X)$ of a subgraph $X$ will not work, since in general it will not even be a graph, as we might wind up with edges whose source or target is missing in the set $c(X)$. There are two obvious ways to ``repair" this complement operation, though: we can either discard such problem edges; or, on the other hand, we can retain them and ``complete" them by adding their sources and targets in the underlying graph. The first option in fact gives rise to our Heyting negation $\neg X$, and the second to the co-Heyting $\sim X$. \par 
	It is instructive to see these notions ``at work" and to make explicit computations with them. So observe that given the graph $G$ together with its subgraph $Y$, as displayed below, then for $\sim Y$ we will have 
	\begin{center}
		\includegraphics*[scale=0.3]{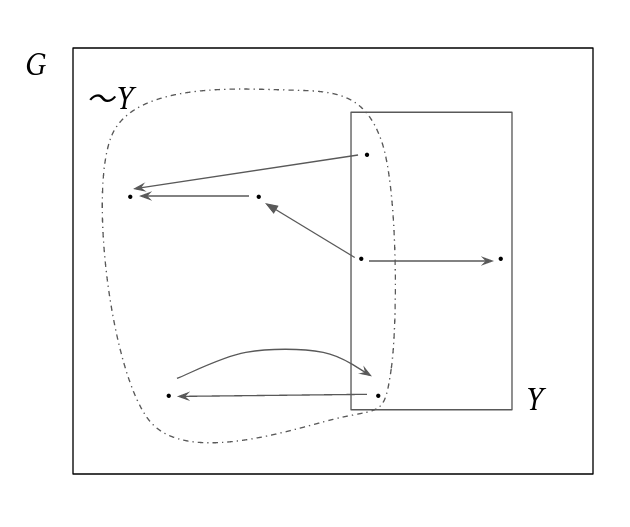}
	\end{center}
	Running this again, we compute $\sim \sim Y$ (which is, importantly, not identical to $Y$!),
	\begin{center}
		\includegraphics*[scale=0.3]{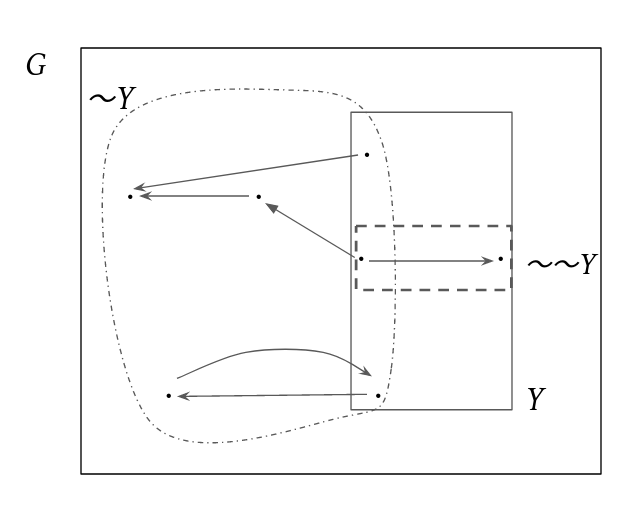}
	\end{center}
	The reader can also verify for themselves that for such a $G$ and $Y$, $\neg Y \neq \sim Y$.\footnote{In general, though, a bi-Heyting algebra for which $\neg x = \sim x$ for all $x$ is necessarily a Boolean algebra.} \par 
	Let us now look closely at an example that computes modalities in graphs.\index{modalities} Suppose we have the following graph\index{graph} $G$, and the subgraph $X$ of $G$, as follows: 
	\begin{center}
		\includegraphics*[scale=0.25]{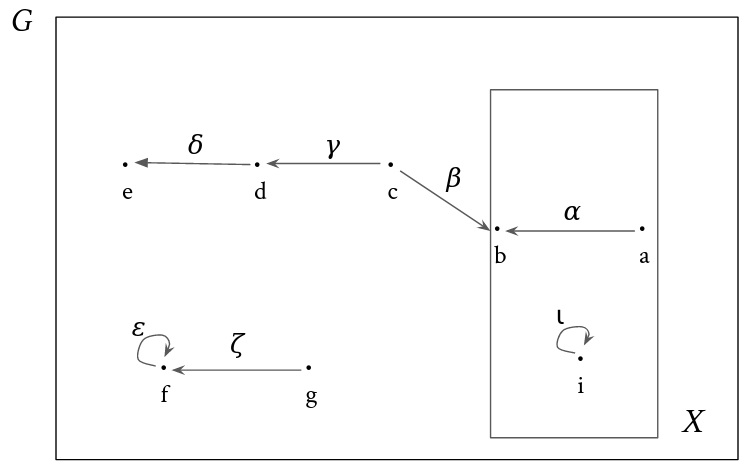}
	\end{center}
	Then, for $\neg X$ we will get the largest subgraph disjoint from $X$:
	\begin{center}
		\includegraphics*[scale=0.25]{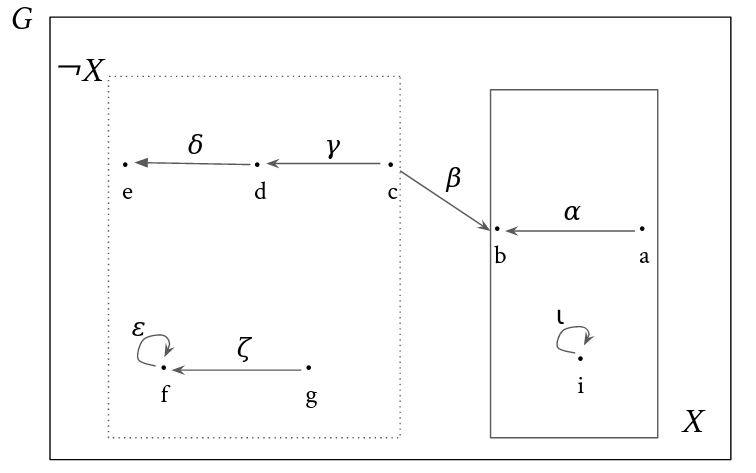}
	\end{center}
	$\sim X$, for its part, is here different from $\neg X$, and yields the smallest subgraph whose union with $X$ gives all of $G$: 
	\begin{center}
		\includegraphics*[scale=0.25]{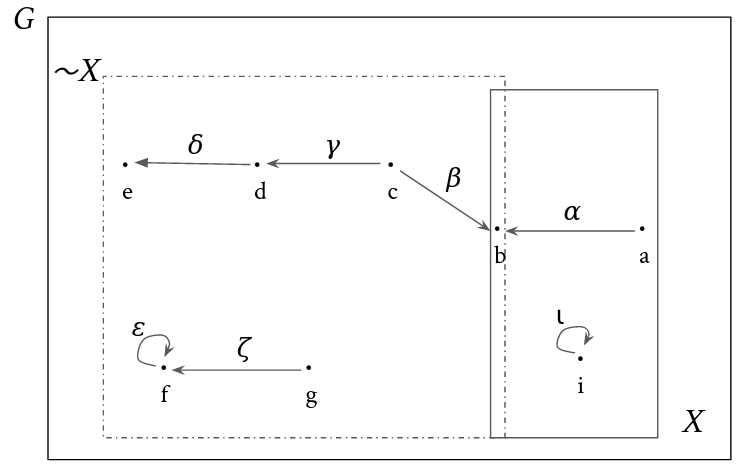}
	\end{center}
	Incidentally, notice that the boundary\index{boundary} of $X$, $\partial X = \sim X \wedge X$, is not the empty subgraph (which functions as $0$),\index{graph} but is the sole vertex $b$. Intuitively, it makes sense to think of the vertex $b$ as the ``boundary" of $X$, since it liaisons between the ``inside" of $X$ and the ``outside," as there is an arrow, namely $\beta$ coming in to $X$ from the outside of $X$ via the target $b$. 
	\par    
	Let us now compute $\Diamond_1 X$, i.e., $\sim \neg X$:
	\begin{center}
		\includegraphics*[scale=0.24]{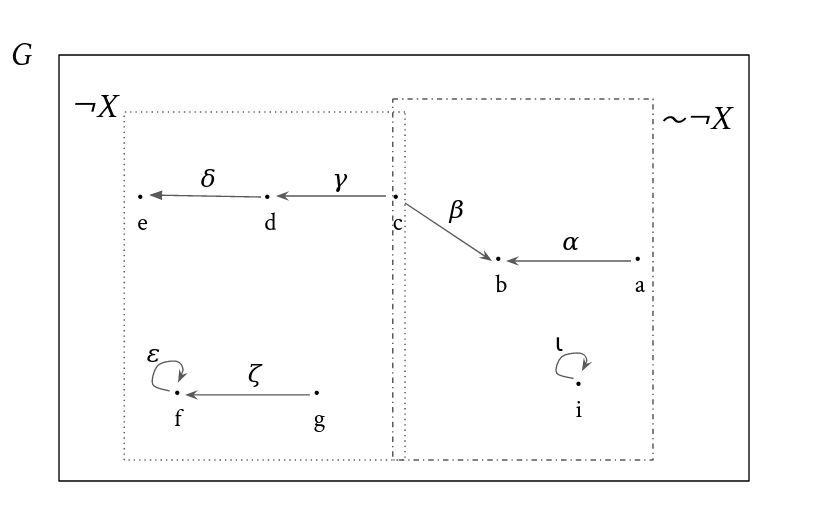}
	\end{center}  
	To compute $\Diamond_2 X$, we must first compute $ \neg \Diamond_1 X$
	\begin{center}
		\includegraphics*[scale=0.24]{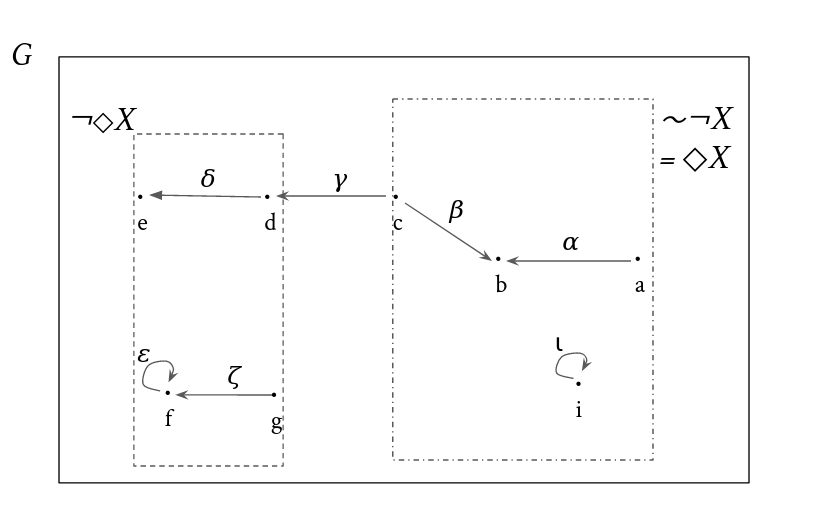}
	\end{center}  
	and then $\sim \neg \Diamond_1 X$, i.e., $\Diamond_2 X$, 
	\begin{center}
		\includegraphics*[scale=0.24]{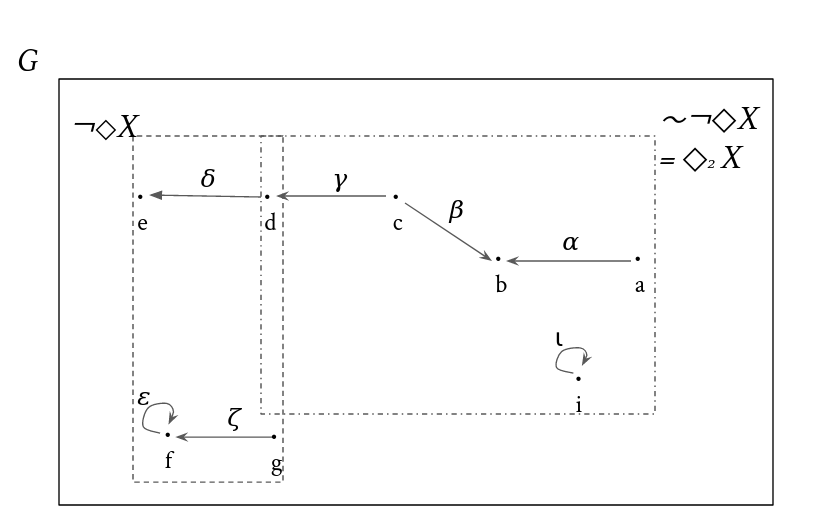}
	\end{center}  
	Running this one more time,\index{graph} we first get $\neg \Diamond_2 X$, 
	\begin{center}
		\includegraphics*[scale=0.24]{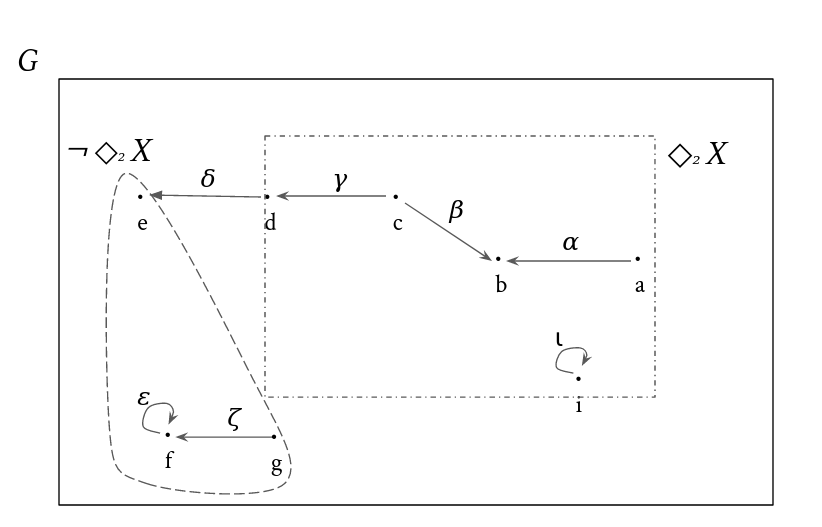}
	\end{center}  
	and then $\sim$ of this, i.e., $\sim \neg \Diamond_2 X = \Diamond_3 X$, 
	\begin{center}
		\includegraphics*[scale=0.24]{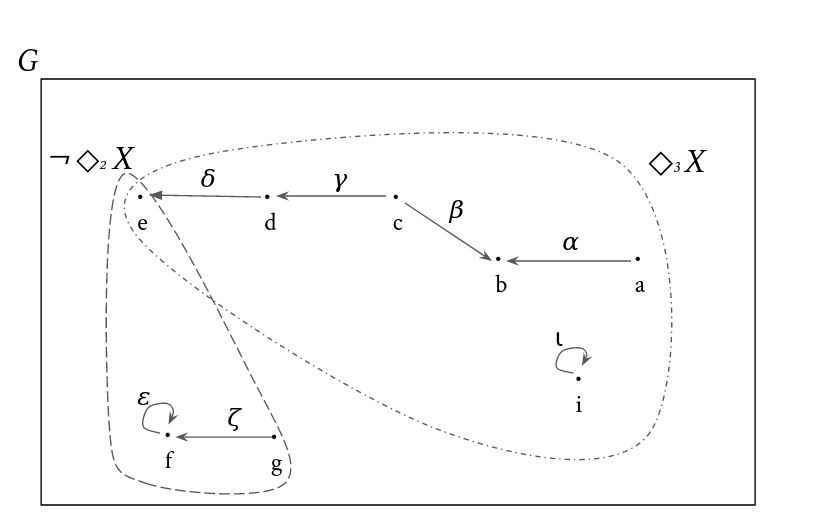}
	\end{center}  
	At this point, something interesting happens. If we run $\neg$ on $\Diamond_3 X$, we get
	\begin{center}
		\includegraphics*[scale=0.24]{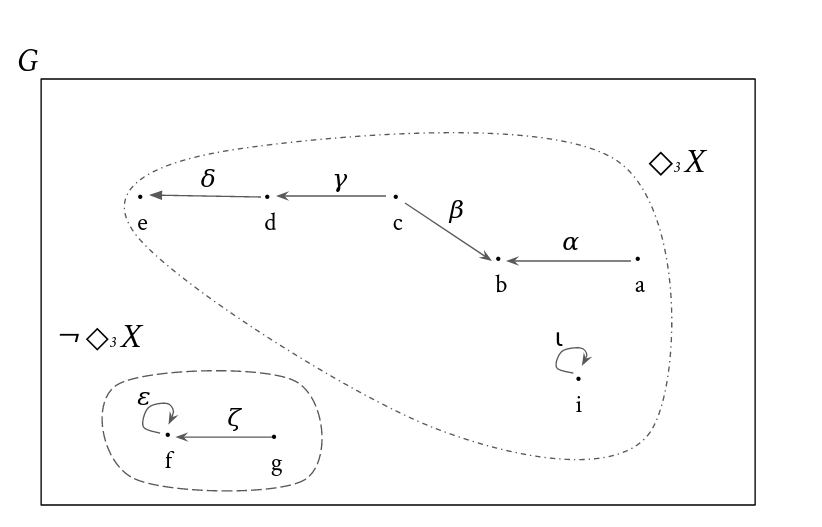}
	\end{center}  
	and then $\sim$ of the above, which gives us $\Diamond_4 X$, is revealed to be the same as $\Diamond_3 X$.
	Thus, this iterative operation stabilizes at $\Diamond_3 = \Diamond_4$, and we have that $\Diamond X = \Diamond_3 X$, having captured, in a subgraph, all those elements of the graph\index{graph} $G$ that can be reached from $X$ through some path, which is clearly something of a ``picture" of the \textit{possibility}\index{possibility} of $X$ (or, perhaps more accurately stated, of what is \textit{possible for} $X$). This illustrates a more general feature as well, namely that every arrow or vertex that is connected to $X$ via some path will end up in $\Diamond X$ after a finite number of steps. In general, applying the operator $\Diamond_n$ to the subgraph should be thought of as capturing those elements connected with $X$ within $n$ paths. Notice also that taking $\neg$ of $\Diamond_3$ is the same as taking $\sim$ of $\Diamond_3$. As such, $\Diamond X$ has no boundary (and no edge going out of it), as $\partial \Diamond X = \Diamond X \wedge \sim \Diamond X$ is the empty subgraph. In general, in the land of graphs, taking the boundary of a subgraph $X$ yields the subgraph whose elements are connected to the ``outside" of $X$.\par 
	We can perform similar computations, in reverse order, to compute $\Box X$, which will supply us with a subgraph whose elements are those that are \textit{not} connected to the outside. First, we compute $\neg$ of $\sim X$, to get $\neg \sim X = \Box_1 X$, 
	\begin{center}
		\includegraphics*[scale=0.24]{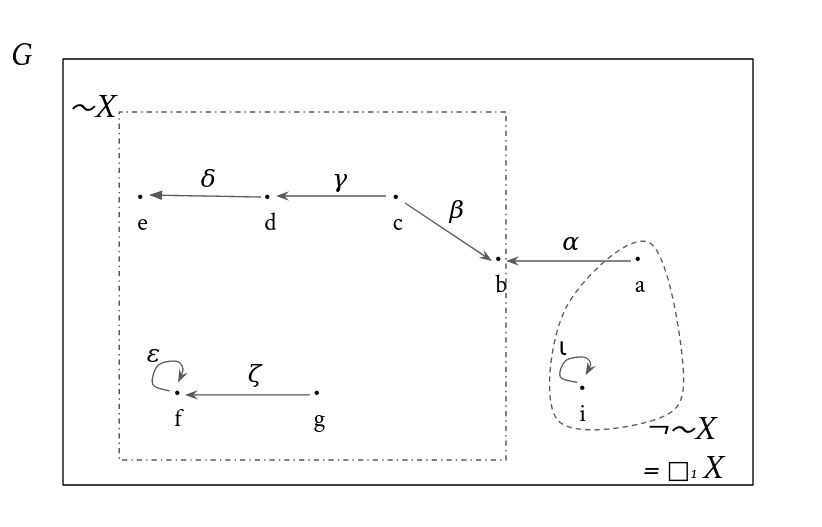}
	\end{center}       
	Then, taking $\sim$ of this, we get $\sim \Box_1 X$, and stripping away one more layer, by taking $\neg$ of the result, we end up with $\neg \sim \Box_1 X = \Box_2 X$:
	\begin{center}
		\includegraphics*[scale=0.24]{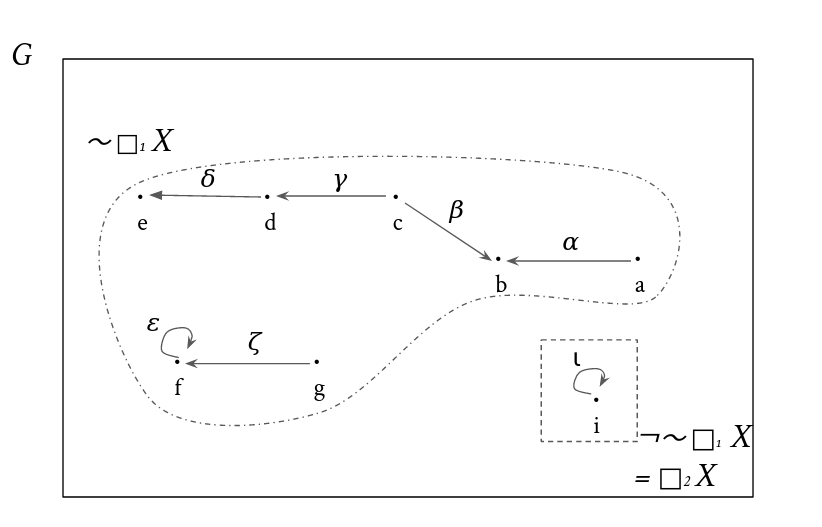}
	\end{center}  
	and we have stability, as any further iterations $\Box_{2 +n} X$ will just reduce back to $\Box_2 X$. Notice that what is left, $\Box X$, just consists of those elements of $X$ that are not connected to the ``outside" (of $X$ in $G$), which seems to align with some intuitions we have about the notion of ``necessity"\index{necessity} for $X$.  \par 
	Altogether then, and for a general graph,\index{graph} $\Diamond X$ supplies the elements of the ambient graph $G$ that can be reached from $X$ via some path, while $\Box X$ has those elements in $X$ not connected, via any path, to the outside. Note, finally, that both the subgraphs $\Diamond X$ and $\Box X$ are complemented sums of connected components.   
	\begin{exercise} Consider the following graph\index{graph} $G$ of routes, with subgraph $X$ corresponding to some region of the northeast (including the nodes Boston, New York, and Long Island, together with the indicated routes between them): 
	\begin{center}
		\includegraphics*[scale=0.2]{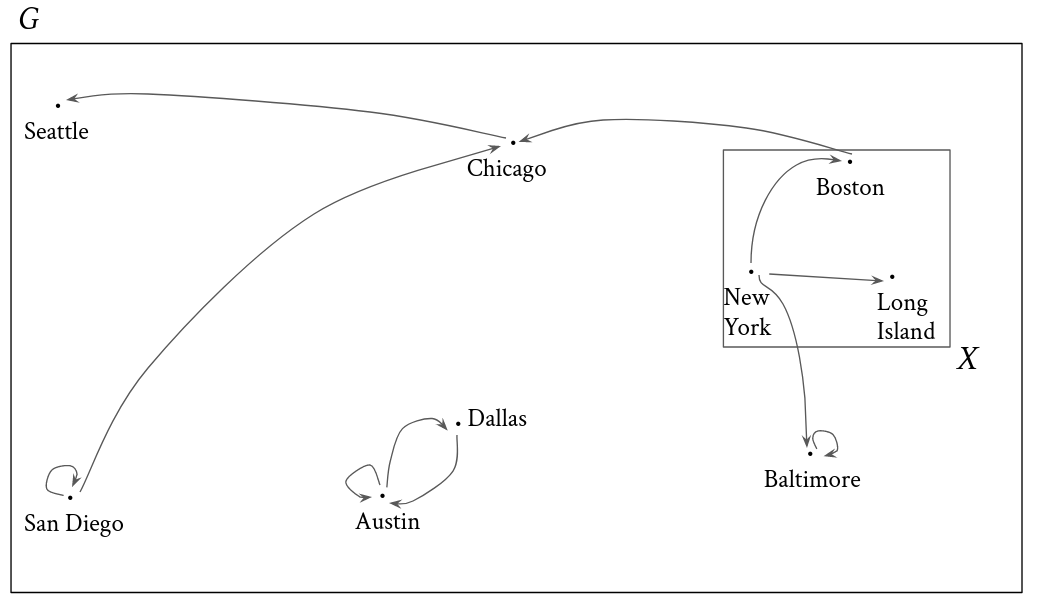}
	\end{center}
Compute $\Diamond X, \Box X,$ and $\partial X$. Then consider, via this example, how the boundary operator $\partial$ interacts with the modal operators. 
\end{exercise} \par \noindent 
\textit{Solution:} First, notice that 
	\begin{center}
		\includegraphics*[scale=0.2]{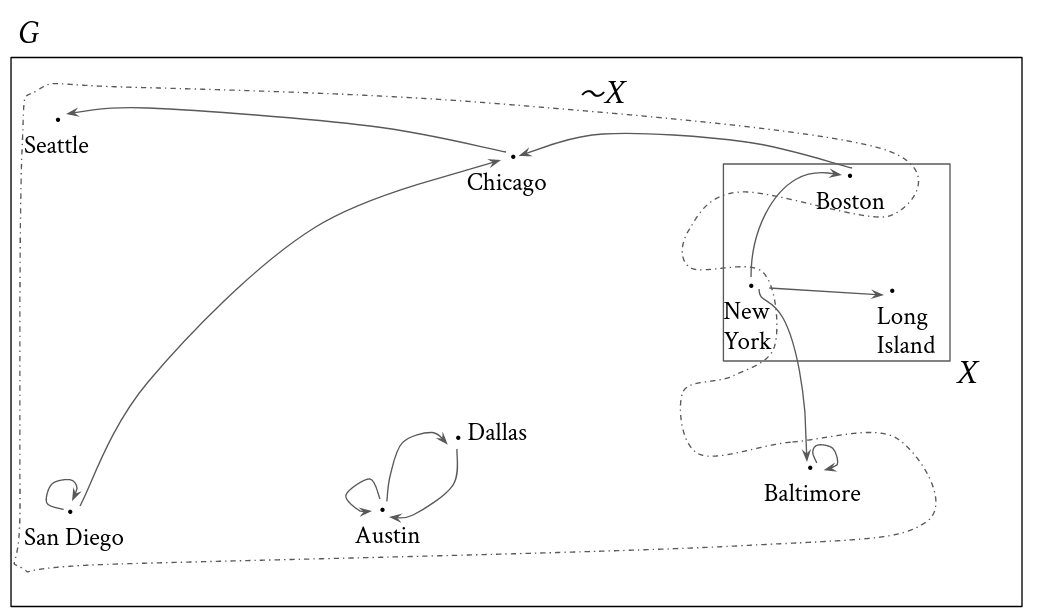}
	\end{center}
and 
	\begin{center}
		\includegraphics*[scale=0.2]{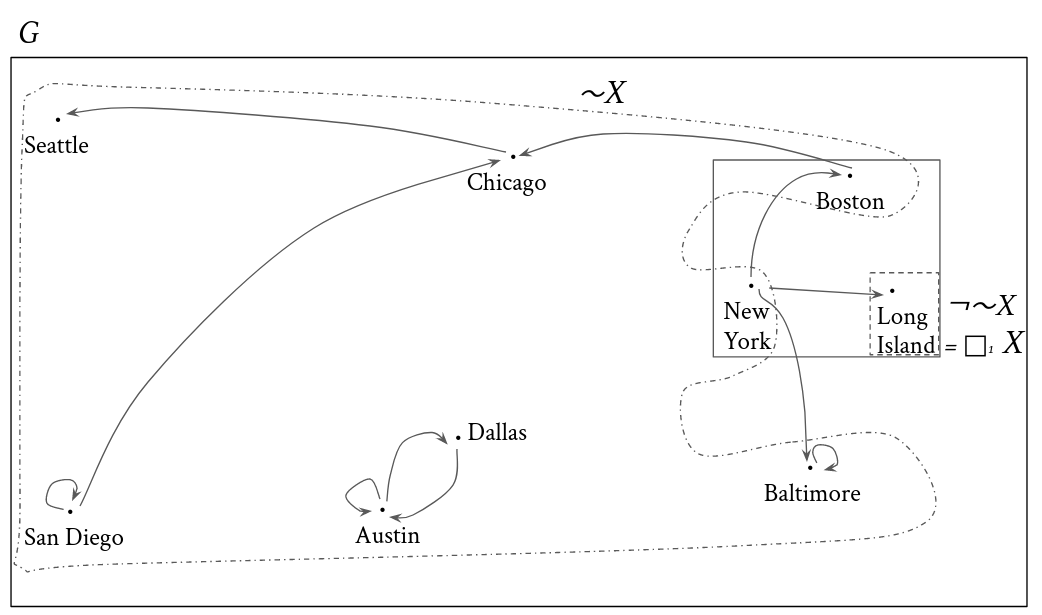}
	\end{center}
	Things stabilize here, with $\Box_1 = \Box$. $\Box X$ includes those parts of $X$ that have no connection to the ``outside" of $X$. The meaning of $\Box X = \{\text{ Long Island }\}$ is this: if one is in $X$ and ends up in Long Island, one will never get out of $X$---having arrived there, one is \textit{necessarily}\index{necessity} in $X$. \par 
	What about $\Diamond X$, or ``possibly" $X$?\index{possibility}   
	\begin{center}
		\includegraphics*[scale=0.2]{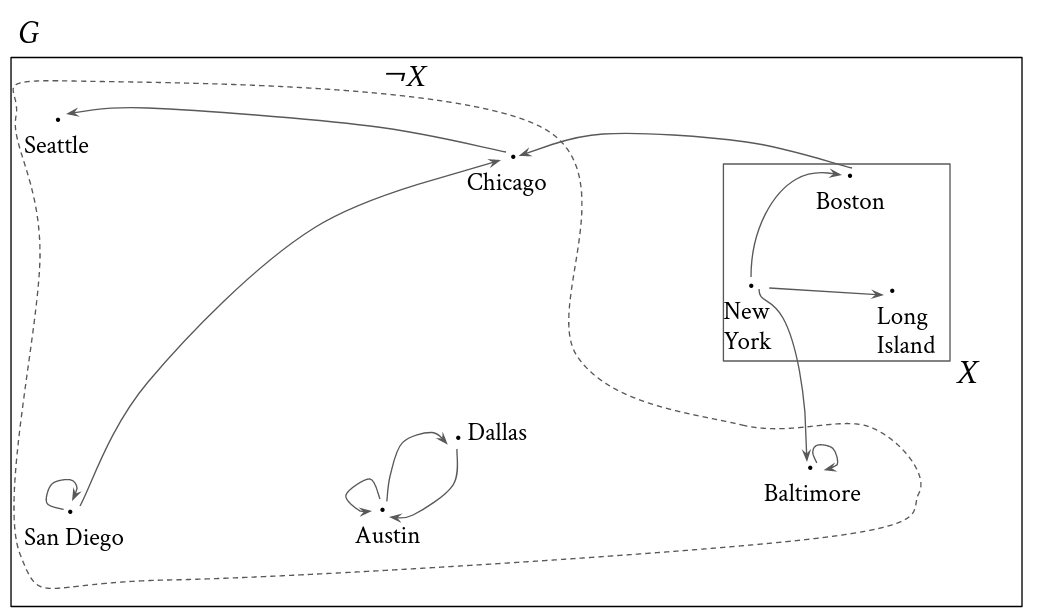}
	\end{center}
	\begin{center}
		\includegraphics*[scale=0.2]{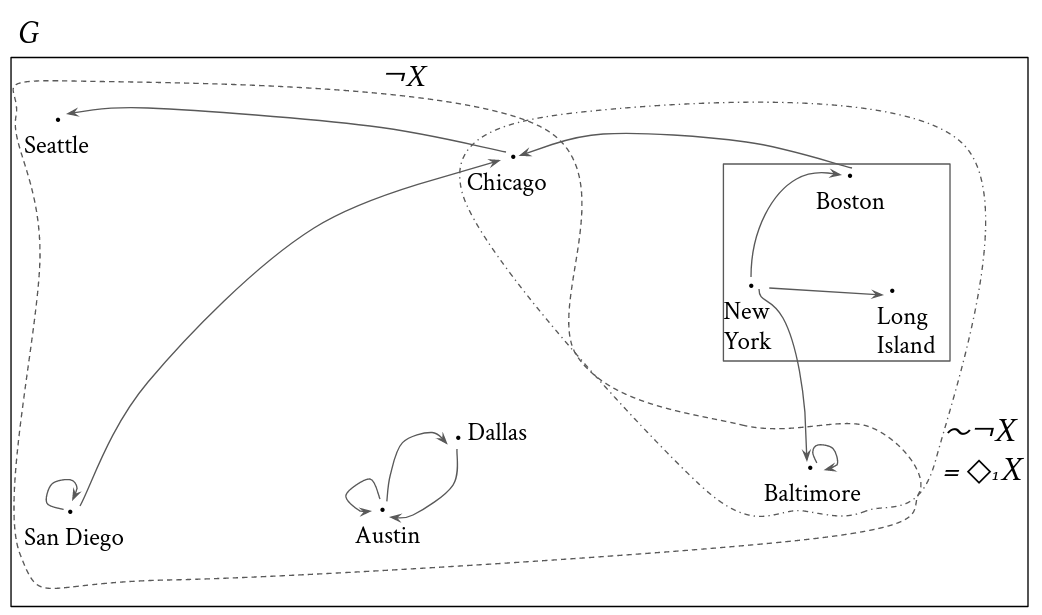}
	\end{center}
	\begin{center}
		\includegraphics*[scale=0.2]{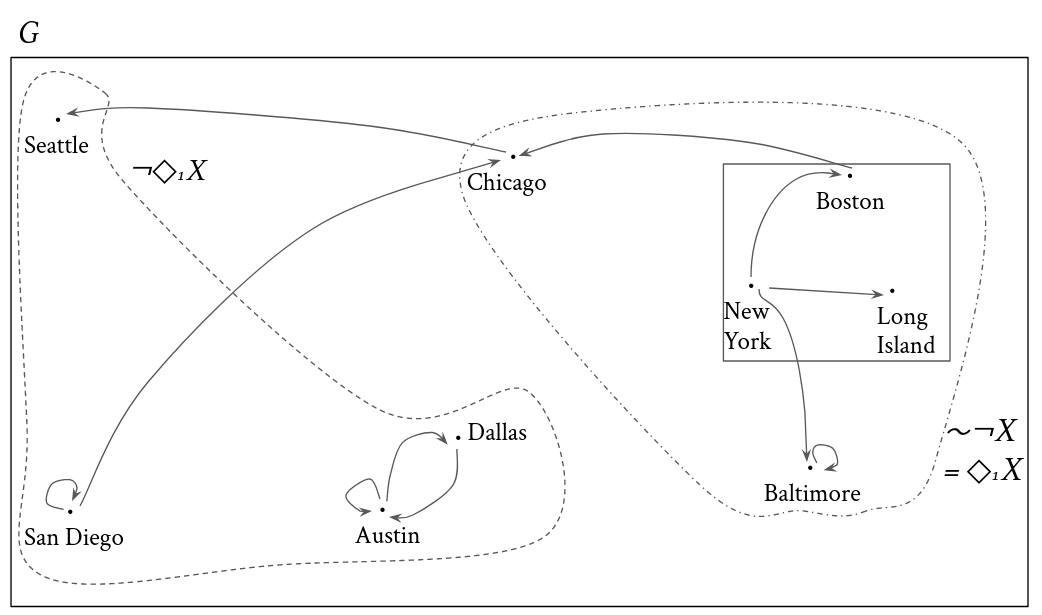}
	\end{center}
	\begin{center}
		\includegraphics*[scale=0.2]{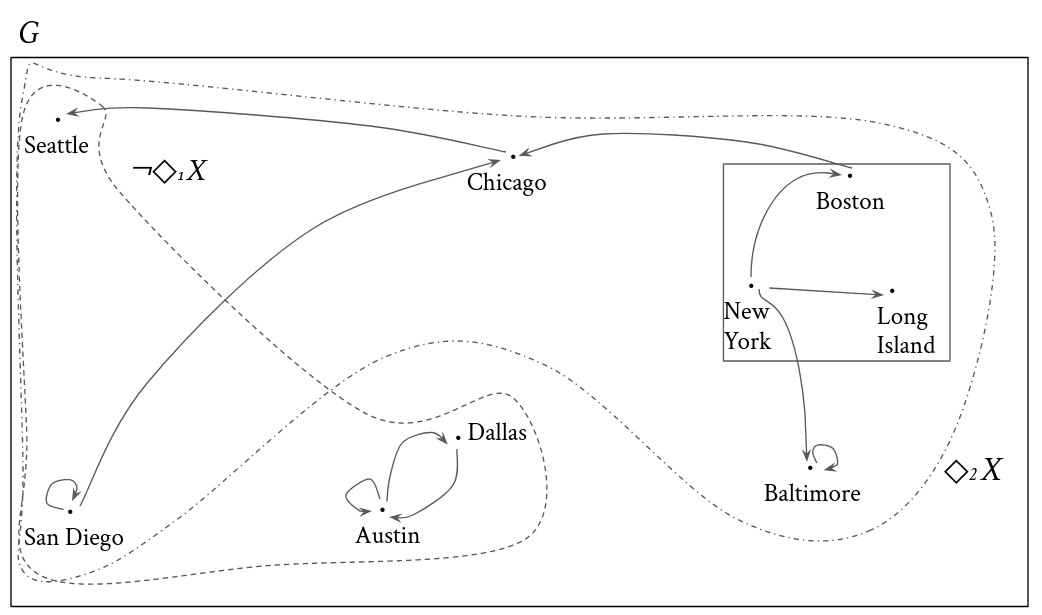}
	\end{center}
	\begin{center}
		\includegraphics*[scale=0.2]{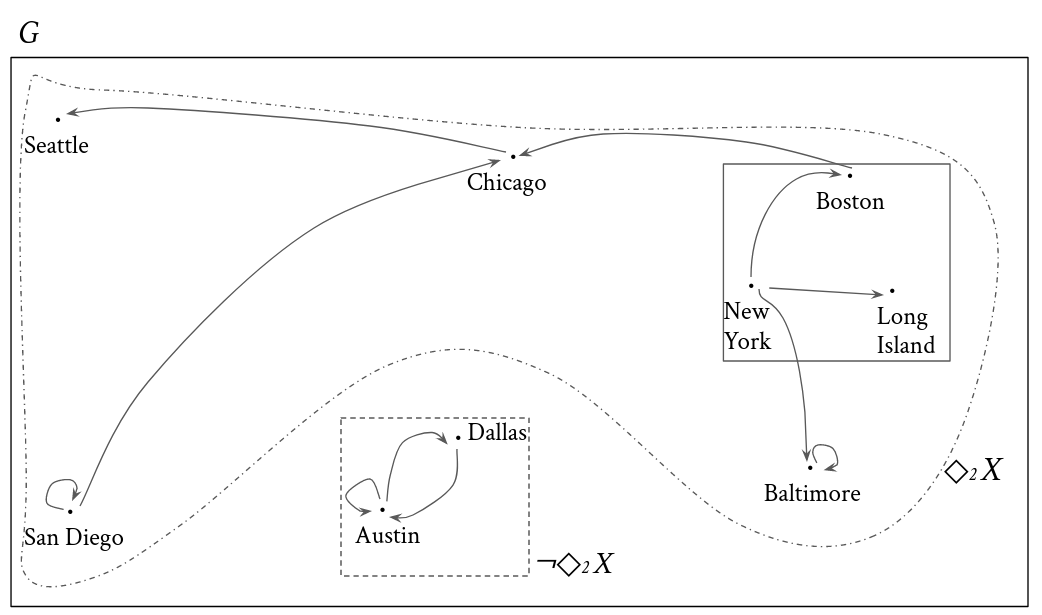}
	\end{center}
	Taking $\sim$ of $\neg \Diamond_2 X$ just returns $\Diamond_2 X$, so we achieve stability at $\Diamond_3 = \Diamond_2$, and so $\Diamond_2 X$ gives us $\Diamond X$, a picture of the ``possibility" of $X$. Intuitively, this makes sense, as $\Diamond X$ supplies those parts of $G$ that can be connected, directly or indirectly, with some part of $X$. For instance, then, even though San Diego is not directly reachable from any city in $X$ or any city reachable by $X$, anyone in San Diego can get to Chicago, and Chicago \textit{is} reachable from $X$. In other words, a person from $X$ might meet someone from San Diego in Chicago or Seattle---and so, for someone from $X$, San Diego may form part of its picture of reality. Dallas and Austin, by contrast, are inaccessible to $X$---given the graph above, someone from anywhere in $X$ could never meet anyone from Dallas or Austin. \par    
	Finally, consider $\partial X = \sim X \wedge X$. As one can see, this is the ``vertices" Boston and New York. Intuitively, this makes sense, as these cities are those parts of $X$ that mediate between the ``inside" of $X$ (as parts of $X$) and the ```outside" in $G$. Seeing how the boundary operator $\partial$ interacts with the modal operators can further solidify the intuitiveness of the reading of $\Diamond$ as ``possibility" and $\Box$ as ``necessity," even in contexts like that of graphs. As one can easily verify, 
	\begin{equation*}
	\partial \Box X = \sim \Box X \wedge \Box X = \Box X
	\end{equation*}   
	which confirms the intuition that the ``boundary"\index{boundary} of what is necessarily $X$---i.e., those parts of $X$ that have no connection to the ``outside" of $X$---is just trivially $\Box X$ itself. Also, 
	\begin{equation*}
	\partial \Diamond X = \sim \Diamond X \wedge \Diamond X = \emptyset,
	\end{equation*}  
	or, more accurately, the empty subgraph. Intuitively, this realizes the idea that the ``world" of what is \textit{not possible} for $X$ has empty overlap with what's \textit{possible} for $X$. 
\end{example} 
Many more examples of important adjunctions will arise organically throughout the book. For now, we conclude the chapter by sketching two examples that are left deliberately somewhat vague, but that are correct ``in spirit" (and, in fact, can be developed to be formally correct as well). They are meant to provide some hopefully engaging examples of adjunctions, while encouraging more fastidious readers to work out the unspecified details on their own. 
\begin{example}
	There is a connection between how the world \textit{appears} to an agent and what that agent \textit{believes} to hold of their world. But `appears' and `believes' are not quite inverses of one another. Instead, we might conjecture that
	\begin{equation*}
	\text{appearance} \dashv \text{belief}, 
	\end{equation*}
	in the sense that there is an adjunction (in a slogan, ``belief as the right adjoint of appearance")
	\begin{equation*}
	\frac{f_{\alpha}(m) \leq m'}{m \leq B_{\alpha}(m')} 
	\end{equation*} 
	realizing, effectively, how all that appears to hold at, or given, state $m$ entails state $m'$ if and only if whenever $m$ holds in the ``real world," this entails that all that agent $\alpha$ \textit{believes} to hold on the assumption that $m'$ holds does in fact hold.  Thus, in general, $B_{\alpha}(m)$ would stand for agent $\alpha$'s belief at $m$ and will consist of those propositions that agent $\alpha$ \textit{believes to hold} whenever $m$ holds. 
\end{example}
\begin{example}
Both small-scale and large-scale projects, such as research or development projects, require resources. Resource allocation (through grants, investment funding, contracts, etc.), requires a detailed plan for how those resources are to be spent, especially as the project increases in scale. If $\textbf{Rsrc}$ is a category\index{category!resource} consisting of relevant resources, so that objects are resources (like, e.g., for simplicity, different-sized checks or bags of money) and morphisms are given by a natural relation between those resources (e.g., $\leq$ in the case of a uniform money-valuation of the different resource objects); and if $\textbf{ProjPlan}$ is a category consisting of project tasks,\index{category!project plan} given some natural ordering (e.g., by order of priority in the carrying out of the plan); then we might consider the functor 
\begin{equation*}
V: \textbf{Rsrc} \rightarrow \textbf{ProjPlan}
\end{equation*} 	
that maps a resource $r$ to the collection of plans $p_i$ that are \textit{viable} given that resource, and the functor 
\begin{equation*}
N: \textbf{ProjPlan} \rightarrow \textbf{Rsrc}
\end{equation*}
taking a project task $p$ to all those resources that are necessary to complete the task (which, depending on how $\textbf{Rsrc}$ is structured, say in a simple case of ``costs," might just amount to returning an interval bounded by the \textit{least} cost for which the task could be carried out, and including all other more ``ample" amounts). \par 
We would probably not expect $V$ and $N$ to construct \textit{inverses} to one another, for we do not expect that, for any given resource $r$, a list of necessary resources for those plans that are deemed viable given $r$ would be \textit{equal} to $r$ (though we might expect that, among the resources, $r \leq NVr$). Similarly, we would not expect that, for a given project task $p$, the result of applying $N$ to $p$ and then $V$ to $Np$, would always \textit{equal} the same task $p$; yet we would expect $VNp \leq p$ in $\textbf{ProjPlan}$. This suggests that we have an adjunction, 
		\begin{center} 
	\tikzset{
		,no line/.style={%
			,draw=none
			,commutative diagrams/every label/.append style={/tikz/auto=false}
		}
	}
	\begin{tikzcd}[column sep =large]
		\textbf{Rsrc} \arrow[shift left = 2]{r}[name=U]{V} & \textbf{ProjPlan} \arrow[shift left]{l}[name=L]{N} \arrow[from=L, to=U, no line, pos=.5]{}{\perp}. 
	\end{tikzcd}
\end{center}  
\end{example}  
\section{Final Thoughts on Fundamentals}
There are many more category-theoretic results and constructions that we could pursue. But with a good working understanding of categories, functors, natural transformations, and adjunctions---built on a number of carefully selected examples and applications---the reader is already well-equipped to tackle more advanced matters, and to approach sheaves in particular. \par   
But before the introduction of sheaves, we promised to consider presheaves more systematically and to dwell on some of the fundamental ideas. Doing so will help solidify a number of important ideas that will resurface in the discussion of sheaves, in addition to presheaves and presheaf categories being of considerable intrinsic interest. A closer consideration of presheaves, in addition to the presentation of new examples of presheaves that will resurface in subsequent chapters, is thus the subject of the next chapter.  
\chapter{Presheaves Revisited} 
\section{Seeing Structures as Presheaves}
When we work with a mathematical structure, it is common to try to approach it in terms of its elements. In general, it is very natural to want to ``break things down" by decomposing more complicated structures into their components---and elements, like dots, are one sort of ``component" we seem especially ready to recognize as such. But in certain settings, one needs to consider \textit{figures of a more general shape} than ``points." Points, after all, might be regarded as just a particularly simple kind of ``shape." For instance, suppose you are presented with the structure $X$: 
\begin{center}
	\includegraphics*[scale=0.25]{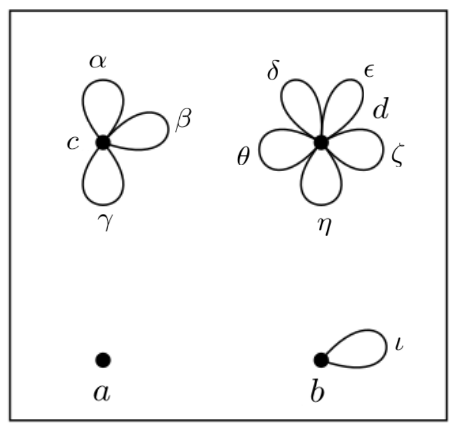}
\end{center}
This $X$ depicts what is called a \textit{bouquet}.\index{bouquet} Figures in the bouquet $X$ with the shape ``point" are maps $\bullet \rightarrow X$, each of which map ``names" a point in $X$: 
\begin{center}
	\includegraphics*[scale=0.25]{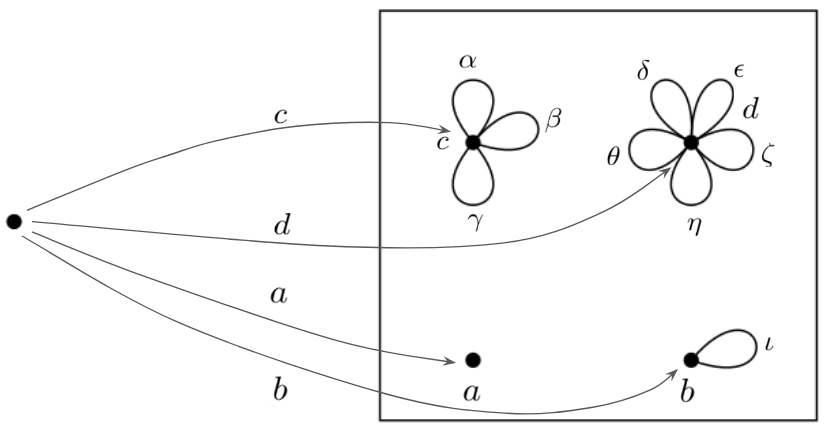}
\end{center}
In other words, we use the `generic shape' $\bullet$ to locate and name all the distinct ``point-like" figures of $X$, i.e., via a particular map $\bullet \xrightarrow{a} X$ we name with $a$ one of the various components of $X$ that ``looks like" $\bullet$. Altogether, the data of our ``points" in $X$ really just amounts to a set 
\begin{equation*}
X(\bullet) = \{a,b,c,d\},
\end{equation*}
which you might read as saying ``$X$ realizes $a, b, c,$ and $d$ as its figures of shape $\bullet$." \par 
But you could not hope to understand all that this structure $X$ is just by considering the point-like figures! After all, points are not the only sort of figural component in $X$. And, as such, there are many distinct bouquets that may even have the same set of points (yet will look rather different!). So we also need a way of picking out those figures in $X$ whose shape is that of a ``loop." Similar to with our points, we can ``pick out" and name our loops via maps from the generic shape \includegraphics*[scale=0.15]{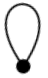} into $X$:  
\begin{center}
	\includegraphics*[scale=0.25]{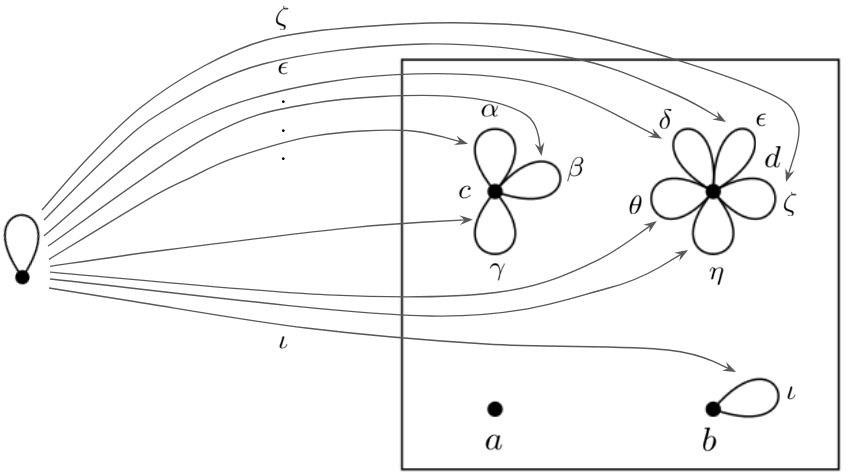}
\end{center}
In other words, the data here is captured by the set
\begin{equation*}
X(\includegraphics*[scale=0.15]{LoopSmall.png}) = \{ \alpha, \beta, \gamma, \delta, \epsilon, \zeta, \eta, \theta, \iota \}.
\end{equation*}
But do we then have enough information to reconstruct $X$? Well, many bouquets\index{bouquet} may have the same set of loops, but the home point at which they ``live" may be different. You would not regard these as the same thing. To fully capture $X$, then, we also need a way of extracting the data of where each loop-shaped figure lives, i.e., how the loop-shaped figures relate to the point-shaped figures. Corresponding to the inclusion of the generic shape ``point" in the generic shape ``loop" 
\begin{equation*}
\bullet \xrightarrow{i} \includegraphics*[scale=0.15]{LoopSmall.png}
\end{equation*}
we should then have a map 
\begin{align*}
X(\includegraphics*[scale=0.15]{LoopSmall.png}) & \xrightarrow{X(i)} X(\bullet) \\
l & \mapsto p
\end{align*}
taking a loop $l$ to the point $p$ at which it is stationed, and so informing us about which loops get stationed at which points. For instance, this will tell us that 
\begin{equation*}
X(i)(\alpha) = c,
\end{equation*}
or ``the loop-shaped component named $\alpha$ lives at the point-shaped component named $c$." The equations telling us which point each of the loops are assigned to just supply us with what are called the \textit{incidence relations}. \par  
With all this information---the set $X(\bullet)$ and $X(\includegraphics*[scale=0.15]{LoopSmall.png})$, together with a map describing how the loop elements in the latter set are sent to the points in the former set---it would seem that we will be able to recover the whole of the information of $X$ itself.\index{bouquet} \par 
But described in this way---and this is the point!---what else have we been saying other than that $X$ is just a presheaf 
\begin{equation*}
X: \mathbb{B}^{op} \rightarrow \textbf{Set},
\end{equation*} 
where the indexing category $\mathbb{B}$ is 
\begin{center} 
	$\mathbb{B} :=$ 
	\begin{tikzcd}[framed]
		\bullet \arrow[r, "i"] & \includegraphics*[scale=0.15]{LoopSmall.png}  
	\end{tikzcd} ? 
\end{center} 
\noindent 
Thus if our domain category is regarded as consisting of some ``shapes," \index{bouquet} related in some particular way (such as points included in pointed loops), then the result of ``realizing" or ``figuring" those shapes, via a presheaf $X$, can be imagined as a ``container" holding onto the various realizations or ``figures" of the different shapes, where the relation between these is respected. \par  
In a similar way, suppose we instead took as our indexing category of ``shapes" the single object $*$ and all the morphisms generated by iterations of $\sigma$, i.e., the free monoid\index{monoid!free} on one generator $(\sigma$),     
\begin{center}  \[\mathbb{E} := 
	\begin{tikzcd}[framed] * \arrow[loop, "{\sigma}", swap]
	\end{tikzcd}
	\] 
\end{center}
Then, $\textbf{Set}^{\mathbb{E}^{op}}$, the category of presheaves on $\mathbb{E}$, is none other than the category of evolutive sets or dynamical systems.\index{category!of dynamical systems} Objects $X$ of $\textbf{Set}^{\mathbb{E}^{op}}$ arise as dynamical systems (evolutive sets) or automata, where $X$ supplies the set of possible states, and the given endomap $\sigma$ gives rise to the evolution of states (think the change in internal state that results after the passage of one unit of time, or as a result of pressing the ``button" $\sigma$ on the outside of a machine). In other words, if $X$ is a presheaf on $\mathbb{E}$, we think of it as a container containing a set of figures (shaped in the form of dots, corresponding to instantiations of the object $*$ of $\mathbb{C}$, and in the form of arrows between certain of those dots, corresponding to the endomap $\sigma$ of $\mathbb{E}$), with a process taking each element to a next stage or next element. In other words, we end up with an $X$ such as 
\begin{center}  
	\begin{tikzpicture}[framed, scale=0.7]
	\tikzset{vertex/.style = {shape=circle,draw, fill=black, minimum size=3pt, inner sep =0pt}}
	\tikzset{edge/.style = {->,> = latex'}}
	\node[vertex] (a) [label=above:{$a$}] at  (0,2) {};
	\node[vertex] (b) [label=below:{$e$}] at  (1.5,0) {};
	\node[vertex] (a4) [label=below:{$f$}] at  (3,0.3) {};
	\node[vertex] (a5) [label=above:{$d$}] at  (4.5,0.5) {};
	\node[vertex] (c) [label=below:{$i$}] at  (6,-1) {};
	\node[vertex] (a1) [label=above:{$b$}] at (1.5,1.75) {};
	\node[vertex] (a2) [label=above:{$c$}] at (3,1.25) {};
	\node[vertex] (a6) [label=above:{$g$}] at (0,-1.75) {};
	\node[vertex] (a7) [label=above:{$h$}] at (1.5,-1.755) {};
	
	\draw[edge] (b) to (a4);
	\draw[edge] (a6) to[bend left=50] (a7);
	\draw[edge] (a7) to[bend left=50] (a6);
	
	\draw[edge] (a)  to (a1);
	
	\draw[edge] (a4) to (a5);
	
	\draw[edge] (a1) to (a2);
	\draw[edge] (a2) to (a5);
	\draw[edge] (a5) to[out =300, in=40, looseness=30] (a5)[above];
	\draw[edge] (c) to[out =300, in=40, looseness=30] (c)[above];
	\end{tikzpicture} \end{center}
Similarly, if we instead took as our ``shape" indexing category the category of $n$-evolving sets, i.e., $\mathbb{E}_n$, freely generated by $n$ non-identity morphisms:
\begin{center}
	\[\mathbb{E}_n :=	
	\begin{tikzcd}[framed, scale=0.7] * \arrow[loop left, "{\sigma_1}"] \arrow[loop above, "{\sigma_2}"] \arrow[loop right, "{\cdots}"] \arrow[loop below, "{\sigma_n}"]
	\end{tikzcd} \] 
\end{center} 
the container of $\mathbb{E}_n$-shaped \index{monoid!free} figures would have figures similar to the above picture, except with (up to) $n$ different processes carrying one $*$-figure to the next, e.g., 
\begin{center}
	\includegraphics*[scale=0.25]{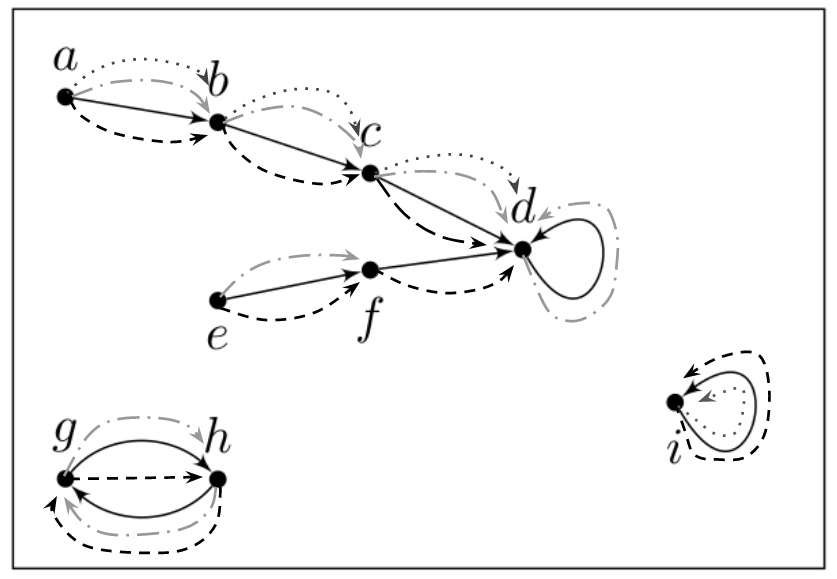}
\end{center} 
A similar approach can be taken to more distinct indexing categories as well. For instance, recall the particular subcategory $\textbf{A}$ of $\textbf{Qua}$,\index{category!qua} 
	\begin{center}
	\begin{tikzcd}
		& \boxed{\text{a scf}}qua\boxed{\text{a scf}} \\
		\boxed{\text{a scf}}qua\boxed{\text{a s}} \arrow[ur] & \boxed{\text{a scf}}qua\boxed{\text{a c}} \arrow[u] & \boxed{\text{a scf}}qua\boxed{\text{a f}} \arrow[ul] \\
		& & \boxed{\text{a scf}}qua\boxed{\text{a p}} \arrow[u] & \boxed{\text{a scf}}qua\boxed{\text{a h}} \arrow[ul] 
	\end{tikzcd}
\end{center}
Applying the \textit{interpretation} functor $X: \textbf{A}^{op} \rightarrow \textbf{Set}$ then results in a ``container" of $\textbf{A}$-shaped figures, e.g., 
\begin{center}
	\includegraphics*[scale=0.38]{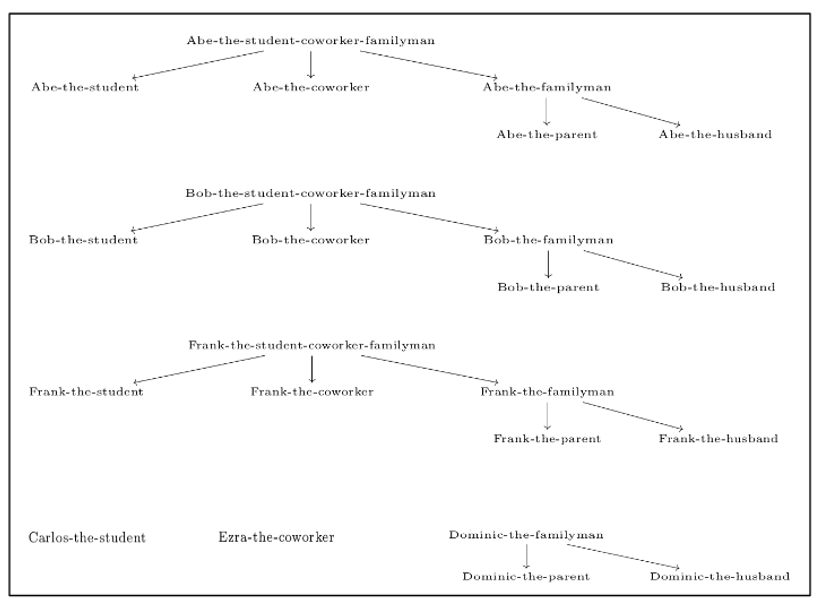}
\end{center}
As before, we could regard $X(S)$, for instance, as picking out the set of people in the ``container" who are ``shaped" like a \textit{student}. If we organize things a bit, grouping together those people that are picked out as conforming to the same ``shape" (role), then what is fundamentally going on here can be displayed more sensibly as   
\begin{center}
	\includegraphics*[scale=0.28]{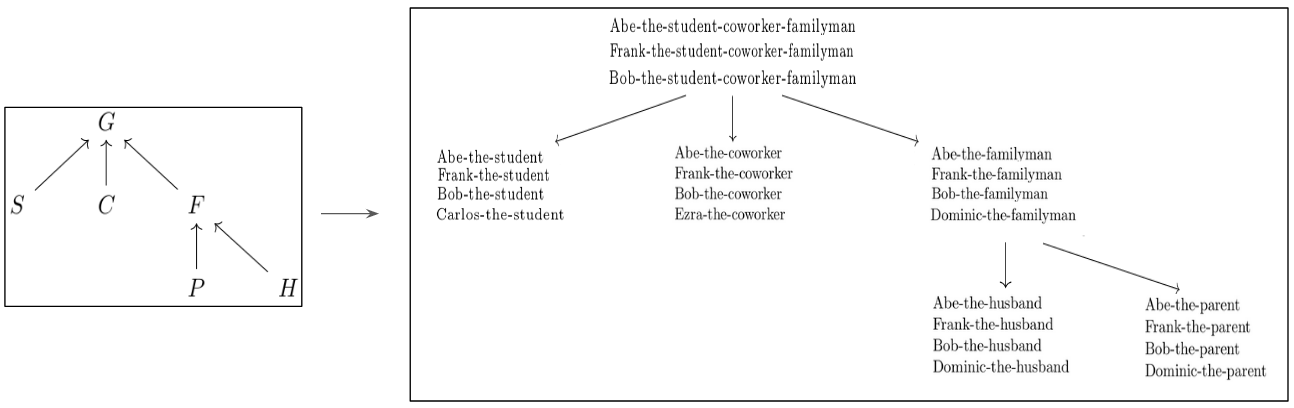}
\end{center}
Similar to the previous examples, the idea here is that the domain category supplies the ``shape" according to which the ``figures" or instantiations are organized, and where the overall realization of such figures as figures of such a shape is accomplished via the functor. \par 
Building on these examples, for a presheaf $P: \textbf{C}^{op} \rightarrow \textbf{Set}$ in general, it is often natural to think of the result of applying $P$ as leaving us with some sort of ``container" of $\textbf{C}$-shaped ``figures," where the various objects $c$ of \textbf{C} are thought of as supplying the ``generic figures" or ``shapes" that are then \textit{instantiated} or \textit{figured}, i.e., $P(c)$ is some particular set of instantiations or figures of $c$-shape. 
A functor fundamentally is just a transformation $P$ that turns objects and maps in one ``world" $\textbf{A}$ into objects and maps in another world $\textbf{B}$, and does so in such a way that certain equations are satisfied (where these code for the transformation preserving structure, or compatibility of the transformation with the composition of maps in $\textbf{A}$ and $\textbf{B}$). Another way of speaking of such a functor $P$ is thus as a \textit{realization} of $\textbf{A}$ in $\textbf{B}$. \par   
This general understanding of presheaves might be thought of in terms of Plato's\index{Plato} notion of the \textit{form} or \textit{shape} (\textit{eidos}) of something, that structural schema according to which the concrete realizations thereof are organized. This form also supports a great variety of realizations or manifestations, and the plurality of particular ``manifestations" (\textit{phantasmata}) of it populate a ``world" that acts as some sort of ``receptable" of those forms. The process by which the manifestations are unfolded according to the structural schema of the form is what Plato would call the ``participation" (\textit{methexis}). The form is held to be invariant, its components sufficiently generic, and altogether it is fundamentally ``simpler" (and so, in the end, more \textit{intelligible}) than its many realizations or manifestations. \par 
Applied to presheaves, the gist here is that the ``generic figures" supplied by the domain category \textbf{C} act as something like the ``form," while the value assignments $P(c)$ for each object of \textbf{C} supply something like the ``appearances" or ``manifestations" of that form, and the presheaf action enforces relationships between the various manifestations based on the invariant relationships between the component objects of the form. The presheaf $P$ itself, on this way of seeing things, would then be nothing other than the process of manifestation or participation of the form in concrete particulars, and understanding how the concrete manifestations and their components respect among themselves the relations that obtain between the components of the form itself is just to understand the general functoriality of the functor $P$.  
\begin{center}
	\includegraphics*[scale=0.23]{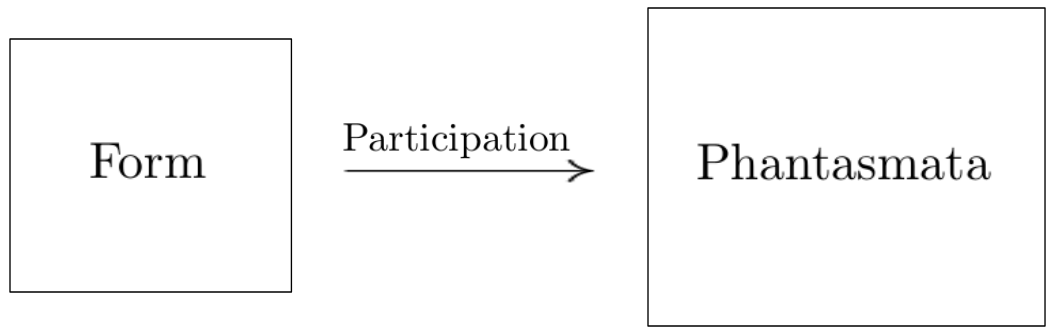}
\end{center}
Presheaves accordingly supply a uniform framework for capturing, in an at once compressed and illuminating way, many structures that appear in math, and that can otherwise appear, in their traditional presentation, in a rather complicated form (thus frequently leading to a complicated description). \par   
So far, we have mostly dwelt on how the presheaf operates \textit{on objects}. Obviously, as a functor, we must also consider the (right) action\index{action} specified by the (contravariant) functor, i.e., how it acts \textit{on morphisms}. The basic idea will be that, following the ``figures" of ``generic shape" $c$ interpretation, a morphism in \textbf{C} from one object to another will give rise to a ``change of figures," where this means, more precisely, that if we have a figure $x$ of shape $c$, i.e., $x \in P(c)$, and a figure $y$ of shape $c'$, i.e., $y \in P(c')$, then asking about the effect of changes of figures amounts, at the level of the presheaf, to asking to what extent the figures are incident or overlap (and what this overlap structure looks like) or otherwise relate. Using a more ``temporal" indexing category, and adopting the interpretation of $P(c)$ as a set existing \textit{at stage} $c$, then the morphisms of \textbf{C}, upon being acted on by the presheaf, would amount to ``varying the stage," so that, overall, the functor can be interpreted as supplying a picture of a set \textit{varying through time}. In the next section, we will look more closely at the ``incidence relation" interpretation, as well as some other approaches, such as the variable-set interpretation---as always, via examples. This will enable us to start to think more systematically about presheaves and their action.
\section{The Presheaf Action}
It is common to see in the literature on presheaves references to \textit{right} \textbf{C}-sets (which is the same as \textit{left} $\textbf{C}^{op}$-sets). Similarly, one can read of a presheaf in terms of its \textit{right action}.\index{action} We will think of there being \textit{four characteristic kinds of cohesivity or variability} presented by presheaf categories in accordance with four main ways the right action of the presheaves in question can be found to operate. But before discussing these interpretations of a presheaf (illustrating them each through select examples), it may be useful to further explain the reference to the presheaf action as a `right action', in case it is not already clear why one sees this referred to as an \textit{action} (and, moreover, why the action is \textit{right}). 
\subsection{Right Action Terminology}
A presheaf is ultimately just a \textit{functor} (one with a particular variance). At least with the usual set-valued presheaves, in applying a given functor $P$ to each of the objects of the domain category \textbf{C}, we just get a bunch of \textit{sets}, $P(c), P(c')$, etc., indexed by the objects of \textbf{C}. The functoriality of the given presheaf $P$, then, just means that for every map $f: c' \rightarrow c$ in \textbf{C}, we will have a \textit{function} (since we are in \textbf{Set}) $P(f): P(c) \rightarrow P(c')$ going the other way. So we just have a function that takes the elements $x$ of the \textit{set} $P(c)$, i.e., the set $P$ ``seen at stage" $c$ (or ``seen in the shape of" $c$), to elements of the \textit{set} $P(c')$, i.e., $P$ ``seen at stage" $c'$ (or ``seen in the shape of" $c'$). In other words, for each element $x$ of $P(c)$ and each map $c' \xrightarrow{f} c$ of \textbf{C}, there is an associated element $xf$ of $P(c')$.\footnote{The reason for writing the element $xf$ this way will be explained in a moment.} Now, the \textit{contravariance} of the functor of course means that the functor applied to a composite $f \circ g$, where $c'' \xrightarrow{g} c' \xrightarrow{f} c$, should be the same as the functor first acting on $f$ then acting on $g$. In other words, in terms of the elements $x \in P(c)$, $f$, $g$, and $xf$ as above, whenever $c'' \xrightarrow{g} c' \xrightarrow{f} c$, we must have $x(fg) = (xf)g$ in $P(c'')$. Moreover, functors must respect identities. But all of this data essentially means that we are dealing with what, in other settings, one would call a \textit{right action} of  \textbf{C} on the underlying set (formed by the presheaf $P$), and where this right action expresses the ``incidence relations" or transitions among the various figures $x, x'$, etc.
\par 
In other contexts, if we have some mapping $X \times A \rightarrow X$, it is common to refer to such a map as a \textit{right action}\index{action} of $A$ on $X$. Usually $A$ is some monoid or group (and $X$ is a set).\index{monoid} The basic idea is that $A$ is thought of as furnishing a set of ``buttons" that control the states of $X$, while the given action $X \times A \xrightarrow{\alpha} X$ is then regarded as supplying us with the data of a state-machine or automaton. Considering a particular ``button" $\textbf{1} \xrightarrow{a} A$ then give rise to an endomap of $X$, specifically $\alpha(\--, a)$, where this means that for each element $x$ of $X$, its image $\alpha(x, a)$ under the action map $\alpha$ is just a new element of $X$. ``Pressing" $a$ once takes a particular state $x$ into the state $\alpha(x,a)$; pressing it twice takes $x$ to $\alpha(\alpha(x,a),a)$; and so on. Of course, we can also press a different button (i.e., take a different object $a'$ of $A$). Combining things, we can press one button and then another. This will mean: supposing we are in state $x$ and the button $a$ is pressed and then the button $a'$, the resulting state will be $\alpha(\alpha(x,a), a')$. As is common, we can choose, notationally, to represent the result of this action $\alpha(x,a) = a \cdot x.$ \par    
In a similar fashion, with presheaves we speak of a \textit{right action} of \textbf{C} on a set $P$ that is partitioned into ``sorts" coming from the objects of \textbf{C} (i.e., parameterized by the objects of \textbf{C}). Being a ``right action" means that whenever we have an arrow $f: c' \rightarrow c$ in \textbf{C} and an element $x \in P(c)$, i.e., an element of the set $P$ of sort $c$, then $xf$ yields an element of $P$ of sort $c'$ subject to the conditions 
\begin{equation*}
\begin{split}
x \text{Id}_c & = x ; \\ 
x(fg) & = (xf)g \text{ whenever } c'' \xrightarrow{g} c' \xrightarrow{f} c \in \textbf{C}. 
\end{split}
\end{equation*}
We write the action in the form of concatenation, i.e., $xf$ is short for $x \cdot f$ where the action $\alpha: \textbf{Set} \times \textbf{C} \rightarrow \textbf{Set}$ is defined as $\alpha(x,f) = x \cdot f$ and the set in question is actually just the disjoint union $\amalg_{c \in Ob(\textbf{C})} P(c)$. \par 
The idea, then, is that given an element $x \in P(c)$ for some $c \in \textbf{C}$, such an $x$ will be acted on by all the morphisms $c' \xrightarrow{f} c$ in \textbf{C}, and will do so in such a way that composite morphisms act as above. In asking what the value of a function $f: c' \rightarrow c$ in \textbf{C} at such an element $x$ will look like, we are asking about $P(f)(x)$. Regarding this in terms of a right action $\alpha(x,f) = x \cdot f$, we have for composite maps, $\alpha(x, fg) := x \cdot (fg) = \alpha(\alpha(x,f), g) = (x \cdot f) \cdot g$. If we agree, notationally, then, to let $x \cdot f = P(f)(x)$,\footnote{One sometimes calls this the \textit{restriction of}\index{restriction} $x$ \textit{along} $f$, and denotes it by either a $|$ or a dot, as in $P(f)(x) = x|f = x \cdot f$. That $f$ gets written on the right of $x$ not only is meant to reveal the underlying right action, but it is a good notational choice since it accords with the induced notation for a composite arrow $f \circ g$ as $x \cdot (f \circ g) = (x \cdot f) \cdot g$.} it is evident that the contravariance of the functor $P$ is equivalent to specifying that \textbf{C} acts\index{action} (and does so \textit{on the right}) on $P$ (regarded as a set). This is evident since 
\begin{equation*}
\begin{split} 
\alpha(x, fg) & = \alpha(\alpha(x,f), g) \\  
x \cdot (fg) & = \alpha (x \cdot f, g)  \\ 
x \cdot (fg) & = (x \cdot f) \cdot g  \\
P(fg)(x) & = P(g)(P(f)(x)) \\ 
P(f \circ g)(x) & = P(g) \circ P(f)(x) .   
\end{split} 
\end{equation*}
The totality of (right) \textbf{C}-sets, i.e., presheaves on \textbf{C}, then induce \textbf{C}-natural morphisms, i.e., natural transformations between presheaves on \textbf{C}, that are covariant.\footnote{One shouldn't worry too much about these remarks, but the above should provide some motivation for why one can often read of presheaves as \textit{right} \textbf{C}-sets (or \textit{left} $\textbf{C}^{op}$-sets), and for the notation decisions commonly found that writes $P(f)(x)$ as $x\cdot f$. Note that the variance is usually \textit{understood} by the context, so it is also common to read just about \textbf{C}-sets.} Having established the reasoning behind that terminological and notational choice, let us now consider more closely the various interpretations this presheaf action takes on in practice.  
\subsection{Four Ways of Acting as a Presheaf} 
We will think of there being \textit{four characteristic kinds of cohesivity or variability} presented by presheaf categories in accordance with four main ways the right action\index{action} can be found to operate:
\begin{enumerate}  
\item as \textit{processual},\index{process} e.g., as passing from sets indexed by one stage to sets indexed by another. This relates to the notion of a \textbf{C}-variable set, modeling sets evolving through time.\footnote{{Objects of \textbf{C} play the role of stages; for every $c$ in \textbf{C}, the set $P(c)$ is the set of elements of $P$ at stage $c$, while the morphisms model transitions between stages. A presheaf on \textbf{C}, then, is just a set varying over the category $\textbf{C}^{op}$.}} 
\item as \textit{extracting boundaries} (or picking out components),\index{boundary!extraction} e.g., graphs with source and target map, simplices picking out lower dimensional boundaries. For something like a topological space that consists of `points' and `edges' and `triangles', etc., in changing figures, we pass from higher-dimensional figures to lower, so that, for instance, the action works by extracting the end-points of an edge or extracting the edges of a triangle, etc.   
\item as \textit{consistency conditions on how different ``probes" of a space relate},\index{probe} where a presheaf $X$ in general is regarded as something like a rule assigning to each object $U$ of $\textbf{C}$ (each ``test space") the set $X(U)$ of admissible maps from $U$ into the space $X$, which space is thus progressively ``probed" by the constituent shapes of the domain category, thereby being ``modeled by" such probes and their interactions.     
\item as \textit{restriction},\index{restriction} e.g., whenever some sort of topology is involved, where the data specified over or about a ``larger" region can be restricted to the data specified over a region included in the former region.\footnote{It is not uncommon to see presheaves and sheaves introduced exclusively via this fourth approach, but the first three perspectives are also important to consider, especially since the first two often involve examples with finitely generated categories (and, as such, provide a good stock of simple and computationally tractable examples), and the third achieves a level of generality that, were it pursued to its end, would ultimately let us speak of sheaves in ``higher dimensions."} 
\end{enumerate} 
 We illustrate these four action\index{action} perspectives, in order, via specific examples. 
	\begin{example}	
		We discussed earlier how presheaves on \textbf{C}, as functors, can also be thought of in general as providing a set of figures with the shape of the indexing category for each object in \textbf{C} and a \textit{process} operator for each morphism in \textbf{C}. For each $a$ in $\textbf{C}^{op}$, the resulting set $F(a)$ is a set of elements of $F$ at stage $a$, while each arrow between objects in $\textbf{C}^{op}$ induces a transition map between the varying set $F$ at stage $a$ and the varying set $F$ at stage $b$ (for an arrow from $b$ to $a$), so that, altogether, we are regarding the objects of $\textbf{C}$ as playing the role of stages of $F: \textbf{C}^{op} \rightarrow \textbf{Set}$ and $F$ itself as a set that ``varies" through the stages.  \par 
		This perspective of the action as exemplifying a kind of \textit{process}\index{process} is nicely illustrated by considering presheaves on a variety of finite indexing categories. For instance, consider the case of finitely free monoids.\index{monoid!free} We saw at the beginning of the chapter that if $\mathbb{E}$ is free monoid on one generator $\sigma$, or the additive monoid of natural numbers, then $\textbf{Set}^{\mathbb{E}^{op}}$, the category of presheaves on $\mathbb{E}$, is none other than the category of evolutive sets or dynamical systems. Objects of $\textbf{Set}^{\mathbb{E}^{op}}$ consist of a set $X$ equipped with a ``process" endomap. For objects $X$ of $\textbf{Set}^{\mathbb{E}^{op}}$, in other words, $X$ supplies the set of possible states, and the given endomap $\sigma$ gives rise to the evolution of states. Referring back to our earlier such $X$, 
		\begin{center}  
			\begin{tikzpicture}[framed, scale=0.7]
			\tikzset{vertex/.style = {shape=circle,draw, fill=black, minimum size=3pt, inner sep =0pt}}
			\tikzset{edge/.style = {->,> = latex'}}
			\node[vertex] (a) [label=above:{$a$}] at  (0,2) {};
			\node[vertex] (b) [label=below:{$e$}] at  (1.5,0) {};
			\node[vertex] (a4) [label=below:{$f$}] at  (3,0.3) {};
			\node[vertex] (a5) [label=above:{$d$}] at  (4.5,0.5) {};
			\node[vertex] (c) [label=below:{$i$}] at  (6,-1) {};
			\node[vertex] (a1) [label=above:{$b$}] at (1.5,1.75) {};
			\node[vertex] (a2) [label=above:{$c$}] at (3,1.25) {};
			\node[vertex] (a6) [label=above:{$g$}] at (0,-1.75) {};
			\node[vertex] (a7) [label=above:{$h$}] at (1.5,-1.755) {};
			
			\draw[edge] (b) to (a4);
			\draw[edge] (a6) to[bend left=50] (a7);
			\draw[edge] (a7) to[bend left=50] (a6);
			
			\draw[edge] (a)  to (a1);
			
			\draw[edge] (a4) to (a5);
			
			\draw[edge] (a1) to (a2);
			\draw[edge] (a2) to (a5);
			\draw[edge] (a5) to[out =300, in=40, looseness=30] (a5)[above];
			\draw[edge] (c) to[out =300, in=40, looseness=30] (c)[above];
			\end{tikzpicture} \end{center}
		the idea is that with this picture we are displaying the presheaf consisting of $X(*) = \{a,b,c,d,e,f,g,h \}$ and where $\sigma$ (i.e., $* \rightarrow *$) acts for instance on the figure $a$ (i.e., on $* \xrightarrow{a} X$) to produce the figure $b: * \rightarrow X$, i.e., $X(\sigma): X(*) \rightarrow X(*)$ takes the particular $*$-figure given the name `a' to the particular $*$-figure given the name `b'.\footnote{The notation $* \xrightarrow{a} X$ is fully justified by the Yoneda results. For any object of $\textbf{Set}^{\textbf{C}^{op}}$, i.e., some presheaf $F$, and any object $c$ of $\textbf{C}$, the set of elements of $F$ of sort or type $c$ can be naturally identified with the set of $\textbf{Set}^{\textbf{C}^{op}}$-morphisms from $\textbf{C}(\--, c)$ to $F$, which is precisely what justifies the abuse of notation that alternately treats the elements of $F$ of sort $c$ as a morphism $c \rightarrow F$ in $\textbf{Set}^{\textbf{C}^{op}}$, letting any $c \xrightarrow{x} F$ be interpreted as a particular figure in $F$ of sort $c$.} \par
		Then a map in this entire presheaf category from a presheaf $X$ (with endomap named $\alpha$) to another object (presheaf) $Y$ (with endomap $\beta$) will be an ``equivariant map" in $\textbf{Set}^{\mathbb{E}^{op}}$, i.e., a natural transformation $(X, \alpha) \xrightarrow{f} (Y, \beta)$ that preserves the structure, where this means it satisfies $f \circ \alpha = \beta \circ f$. \par  
		We also saw how the same story is easily generalized to the category of $n$-evolving sets, i.e., $\mathbb{E}_n$, freely generated by $n$ non-identity morphisms, so that the container of $\mathbb{E}_n$-sets would have figures similar to the above picture, except with (up to) $n$ different processes carrying one $*$-figure to the next. \par 
		We could of course also further consider finitely generated monoids, e.g., $\mathbb{E}_{1,R}$, where certain relations are imposed on the indexing category. For instance, taking $\mathbb{E}_{1,R}$ with one object and non-identity morphism $\sigma$ obeying some relation $R$---say the relation $\sigma^2 = id_*$---then the resulting category of presheaves on $\mathbb{E}_{1,R}$ gives rise to what is usually called, in other contexts, the category of \textit{involution sets}. We could generalize this to any presentation of a monoid $M$, $\mathbb{E}_{n,R}$, where $M = \langle n | R \rangle$, for $n$ generators (i.e., sigmas), and $R$ relations, ultimately leading to the notion that the usual Cayley graph\index{graph!Cayley} for a group is just a presheaf on the category $\mathbb{E}_{n, R}$.\par 
		In this context, we can take the opportunity to highlight that presheaves on a monoid are just equivalent to the usual \textit{right actions}\index{action} on a set by a monoid (which is in part responsible for the `right action' terminology). Recall that a monoid, viewed as a category $\mathcal{M}$ with just one object $*$, will imply that a set-valued functor on $\mathcal{M}$ yields just \textit{one set}, $F(*) \in Ob(\textbf{Set})$. We must also supply, though, a function from $Hom_{\mathcal{M}}(*,*)$ to $Hom_{\textbf{Set}}(F(*), F(*))$, i.e., from $M$ to $Hom_{\textbf{Set}}(F(*), F(*))$. In general, given a set $A$, for any sets $X, Y$, there is a bijection 
		\begin{equation*}
		Hom_{\textbf{Set}}(X \times A, Y) \xrightarrow{\cong} Hom_{\textbf{Set}}(X, Y^A)
		\end{equation*} where $Y^A := Hom_{\textbf{Set}}(A, Y)$ the set of functions from $A$ to $Y$. Moving between these two equivalent formulations via the bijection is sometimes called ``currying." Our function from $M$ to $Hom_{\textbf{Set}}(F(*), F(*))$ just belongs to $Hom_{\textbf{Set}}(M, F(*)^{F(*)})$. Currying, this is the same as a function $\circlearrowright: M \times F(*) \rightarrow F(*)$. Functors preserve identities by definition, so the monoid action law that $e \circlearrowright s = s$ is satisfied, while the composition law for functions provides the other monoid action law. This shows that each monoid action is just nothing other than a set-valued functor. Depending on the variance of the functor from a monoid $M$ to \textbf{Set}, we get the left (covariant) or right (contravariant) $M$-sets.\footnote{Just as for monoids, a group action on a set $S \in Ob(\textbf{Set})$ is just a functor $G \rightarrow \textbf{Set}$ that sends the single object of $G$ to the set $S$. Right $G$-sets are the same as the presheaf category $\textbf{Set}^{G^{op}}$.} \par 
		The variability provided by the right action in each of the above examples is fundamentally \textit{processual}. This perspective is even more clear in an important related example, where we consider sets varying over some ``time-like" linearly ordered category, such as over $\textbf{N}$, the linearly ordered set of natural numbers, regarded as a category. Then, objects in \textbf{Sets$^{\textbf{N}}$} are just sets varying through $n$ successive stages, i.e., $(n-1)$-tuples of maps: \par
		\begin{equation*}
		X: 	X_0 \xrightarrow{f_0} X_1 \xrightarrow{f_1} X_2 \xrightarrow{f_2} \cdots X_{n-2} \xrightarrow{f_{n-2}} X_{n-1}.  
		\end{equation*} 
		The functor $\textbf{N} \rightarrow \textbf{Set}$ picks a sequence $X_0 \rightarrow X_1 \rightarrow \cdots$ of sets $X_n$ and functions $X_n \rightarrow X_{n+1}$. A morphism between two such objects (sequences) is a sequence of functions, e.g., 
		\begin{center} 
		\begin{tikzcd}
			X_0 \rar \dar & X_1 \rar \dar & X_2 \cdots \dar \\
			Y_0 \rar      & Y_1 \rar      & Y_2 \cdots
		\end{tikzcd} 
		\end{center} 
		such that each individual square commutes. More generally, we have an $\textbf{N}$-indexed family of functions $(f_i: X_i \rightarrow Y_i)_{i \in N}$ compatible with the maps, i.e., whenever $i \leq j$, this square commutes: 
		\begin{center}  
		\begin{tikzcd}
			X_i \arrow[d, "f_i", swap] \arrow[r, "\alpha_{ij}"] & X_j \arrow[d, "f_j"] \\
			Y_i \arrow[r, "\beta_{ij}", swap] & Y_j
		\end{tikzcd}
		\end{center}  
		This is equivalently just describing a natural transformation of the functors $X$ and $Y$. So these transition functions $\alpha_{ij}: X_i \rightarrow X_j$ for each $i \leq j$, should satisfy
		\begin{itemize}
			\item $\alpha_{ik} = \alpha_{jk} \circ \alpha_{ij}$ when $i\leq j \leq k$
			\item $\alpha_{ii} = id_{X_i} \forall i$
		\end{itemize}
		Composing with multiple arrows would look something like 
		\begin{center} 
		\begin{tikzcd}
			X \arrow[d, "f"] & X_0 \rar["\alpha_{01}"] \dar & X_1 \rar["\alpha_{12}"] \dar & X_2 \arrow[r, "\alpha_{23}"] \dar & X_3 \dar \cdots \\
			Y \arrow[d, "g"] & Y_0 \rar["\beta_{01}"]  \dar & Y_1 \rar["\beta_{12}"] \dar & Y_2 \arrow[r, "\beta_{23}"] \dar & Y_3 \dar \cdots \\
			Z & Z_0 \arrow[r, "\gamma_{01}"] & Z_1 \arrow[r, "\gamma_{12}"] & Z_2 \arrow[r, "\gamma_{23}"] & Z_3 \cdots \\
		\end{tikzcd} 
		\end{center} 
		where we require that each individual square commutes. The basic idea here is that once an element is in a set, e.g., $ x \in X_t$, it \textit{remains there}, i.e., $\alpha_{tt'}(x) \in X_{t'}$. However, certain elements $a, b \in X_t$ could become identified in the long run (so the $\alpha$'s do not have to be injective); additionally, new elements can appear over time (something that is expressed by the fact that the maps do not have to be surjective). \par 
		But there is of course not really any need to restrict attention to linear orders. Thus, we could consider the functor category \textbf{Set}$^{\mathcal{P}}$ of sets varying over a \textit{poset} $\mathcal{P}$.\footnote{This is closely connected to the notion of Kripke models (and thus intuitionistic logic), a connection we take up briefly later on.} Here, too, it is entirely sensible to regard the resulting functor objects as $\mathcal{P}$-variable sets, since we have sets ``varying" according to the shape of the order supplied by $\mathcal{P}$. More generally, the idea of a set varying over an ordered (poset or preordered) set is really all just a specialization of the general idea of a set ``varying over" some arbitrary small category. 
	\end{example}
		\begin{example}
			We might also consider presheaves on categories with more than one object. For instance, consider the presheaf from the very beginning of the chapter, where the indexing category was 
				\begin{center} 
				$\mathbb{B} :=$ \begin{tikzcd}[framed]
					V \arrow[r, "v"] & L 
				\end{tikzcd}
			\end{center} 
		with identities left implicit. In the indexing category $\mathbb{B}$, the object $V$ of course stands for vertex and $L$ for loop, and the single (non-identity) morphism $v$ amounts to taking vertices to the vertices of the loops. The presheaves on this $\mathbb{B}$ yielded ``bouquets,"\index{bouquet} or those structures with any number of loops stationed at vertices.  A particular presheaf $X$ on this indexing category, then, gives all the data of a particular bouquet $X$, entirely described by a set $X(V)$ of vertices and a set $X(L)$ of loops, together with a function $X(L) \xrightarrow{X(v)} X(V)$ that acts to pick out the vertex of each of the loops. The action $\gamma \cdot v = c$ or $X(v)(\gamma) = c$, where $\gamma \in X(L)$, just functions to ``extract" the appropriate vertex (boundary) of the loop in question. Finally, a natural transformation from one bouquet (presheaf on $\mathbb{B}$) $X$ to another bouquet (presheaf on $\mathbb{B}$) $Y$ will just amount to a rule that sends loop-figures in $X$ to loop-figures in $Y$, point-figures of $X$ to point-figures of $Y$, and does so in such a way that it preserves the incidence relations. In other words, $t: X \Rightarrow Y$ is a map 
			\begin{center}
				\begin{tikzcd}
					\includegraphics*[scale=0.15]{LoopSmall.png} &  X(\includegraphics*[scale=0.15]{LoopSmall.png}) \arrow[r, "t"] \arrow[d, "X(i)", swap] &  Y(\includegraphics*[scale=0.15]{LoopSmall.png}) \arrow[d, "Y(i)"] \\ 
					\bullet \arrow[u, "i"] & X(\bullet) \arrow[r, "t", swap] & Y(\bullet)
				\end{tikzcd}
			\end{center}
			recovering the appropriate notion of a mapping between bouquets. \par 
			For our purposes, the thing to note in the above example is how the presheaf action is one that amounts to an operation of \textit{boundary extraction}. The presheaf action operates by extracting from a loop-figure the vertex-figure to which it is attached, an operation it is very natural to think of as ``taking the boundary," or ``extracting" the simpler parts that form the components of the given (``higher-dimensional") figures. \par 
			For another example of this ``boundary extraction" type, consider a (directed, multi-)graph $X$ 
			\begin{center} 
			\begin{tikzpicture}[framed, scale=0.8]
			\tikzset{vertex/.style = {shape=circle,draw, fill=black, minimum size=3pt, inner sep =0pt}}
			\tikzset{edge/.style = {->,> = latex'}}
			\node[vertex] (a) [label=above:{$a$}] at  (0,0) {};
			\node[vertex] (b) [label=above:{$d$}] at  (4.5,1) {};
			\node[vertex] (a4) [label=above:{$e$}] at  (6,1) {};
			\node[vertex] (a5) [label=above:{$f$}] at  (7.5,1) {};
			\node[vertex] (c) [label=below:{$g$}] at  (6,-1) {};
			\node[vertex] (a1) [label=above:{$b$}] at (2,0) {};
			\node[vertex] (a2) [label=above:{$c$}] at (4,0) {};
			
			\draw[edge] (b) to["${\eta}$"] (a4);

			\draw[edge] (a)  to[bend left, "${\alpha}$"] (a1);
			\draw[edge] (a)  to[bend right=28,"${\beta}$"] (a1);
			\draw[edge] (a1) to[bend left = 70, "${\gamma}$"] (a);
			
			\draw[edge] (a4) to["${\theta}$"] (a5);
			
			\draw[edge] (a1) to[bend left, "${\delta}$"] (a2);
			\draw[edge] (a2) to[bend left, "${\epsilon}$"] (a1);
			\draw[edge] (a2) to[out =300, in=40, looseness=32, "${\zeta}$", style={pos=0.55},swap] (a2);
			\draw[edge] (a5) to[out =300, in=40, looseness=30, "${\iota}$",swap] (a5)[above];
			\end{tikzpicture} 
			\end{center} 
			In a similar way to what we have been doing with the other examples, we can regard this as a functor, i.e., as being generated by a presheaf on a particular indexing category. Moreover, there is then the obvious action representing the ``boundary extraction" of the source and target vertices (boundaries) from a given arrow. More explicitly, take for indexing category the category consisting of two non-identity arrows (the identities again being left understood),
			\begin{center} 
				$\mathcal{G} :=$ \begin{tikzcd}[framed]
					V \arrow[r, shift left = 1ex, "s"] \arrow[r, "t", swap] & A 
				\end{tikzcd}
			\end{center} 
			where the arrows $s,t$ go from an object $V$ (think node or vertex) and another object called $A$ (think arrow). Regarding $X$ as a presheaf on $\mathcal{G}$, then, is straightforward: the presheaf $X: \mathcal{G}^{op} \rightarrow \textbf{Set}$ just assigns a set of $V$(or vertex)-shaped objects, a set of $A$(or arrow)-shaped objects, and functions $X(A) \xrightarrow{X(s)} X(V)$ and $X(A) \xrightarrow{X(t)} X(V)$, where the function $X(s)$ just assigns to each arc its source vertex and the function $X(t)$ picks out each arc's target vertex, thus giving us again an action that can naturally be thought of as performing a sort of ``boundary extraction." \par 
			More explicitly, for our given graph $X$ displayed above, the data (including some of the action data) is 
			\begin{equation*}
			\begin{split}
			X(V) & = \{a,b,c,d,e,f,g\} \\ 
			X(A) & = \{\alpha, \beta, \gamma, \delta, \epsilon, \zeta, \eta, \theta, \iota\} \\ 
			X(s)(\alpha) & = X(s)(\beta) = a , X(s)(\gamma) = X(s)(\delta) = b , X(s)(\epsilon) = X(s)(\zeta) = c, \dots . 
			\end{split}
			\end{equation*}
			By taking presheaves like $X$ for our objects, and natural transformations between such functors for our morphisms (which preserve the incidence relations), we recover the usual graph morphisms, i.e., graph homomorphisms---showing that the presheaf category $\textbf{Set}^{\mathcal{G}^{op}}$ is none other than the category \textbf{dGrph} (again, also sometimes just denoted \textbf{Grph}, or in other contexts, \textbf{Quiv}, for the category of ``quivers"), consisting of irreflexive directed graphs. We can perform a similar analysis for other sorts of graphs, for instance reflexive graphs. Recall that a graph is reflexive provided each vertex $v$ is assigned a designated edge $v \rightarrow v$. (Equivalently, in terms of quivers, a reflexive quiver has a designated identity edge $id_X: X \rightarrow X$ on each object $X$.) For reflexive graphs, we would take for our indexing category  
			\begin{center} 
				$\mathcal{G}' :=$ \begin{boxedtikzcd}[row sep = huge]
					V \arrow[r, bend left = 30, "s"] \arrow[r, bend right = 30, "t", swap] & A \arrow[l, "l"] 
				\end{boxedtikzcd}
			\end{center} 
		consisting of two non-identity arrows (the identities being understood), just as in $\mathcal{G}$, but now with an extra, third arrow going in the other direction from the two already given arrows. This indexing category is subject to the following equations: 
			\begin{equation*}
			\begin{split}
			l \circ t = Id_V = l \circ s . 
			\end{split}
			\end{equation*}
			$\textbf{Set}^{\mathcal{G}'^{op}}$ recovers the category of reflexive (directed, multi-)graphs, $\textbf{rGrph}$ (or $\textbf{rQuiv}$).\footnote{With such a category, there is then an obvious forgetful functor $U: \textbf{rQuiv} \rightarrow \textbf{Quiv}$ from reflexive graphs to irreflexive graphs, where this acts by neglecting the structural component $l$.} Note also that maps of reflexive graphs, i.e., natural transformations between the presheaf objects of $\textbf{Set}^{{G}'^{op}}$, must then not only respect the source and target maps (as was the case with irreflexive graphs, using $\mathcal{G}$), but also the extra map $l$.\par 
			We can generalize all of this to $n$-uniform hypergraphs taking values in multi-sets or still other graph structures. Moreover, as we discussed in an earlier section, each category can be regarded as a directed graph with some structure. We can thus generalize this and consider the $n$-dimensional analogue of a directed graph, i.e., via `globular' shapes. 
			\begin{definition}
			For $n \in \mathbb{N}$, an $n$-\textit{globular set}\index{globular set} $X$ is a diagram 
			\begin{center} 
			\begin{tikzcd} 
				X(n) \arrow[r, shift left = 1ex, "s"] \arrow[r, "t", swap] & X(n-1) \arrow[r, shift left = 1ex, "s"] \arrow[r, "t", swap] & \cdots \arrow[r, shift left = 1ex, "s"] \arrow[r, "t", swap] & X(1)  \arrow[r, shift left = 1ex, "s"] \arrow[r, "t", swap] & X(0) \end{tikzcd}
			\end{center} 
		of sets and functions such that $s(s(x)) = s(t(x))$ and $t(s(x)) = t(t(x))$ for all $m \in \{2,\dots,n\}$ and $x \in X(m)$. \par 
		But an $n$-globular set can also be defined as a presheaf on the category $\mathbb{G}_n$ generated by the objects and arrows 
		\begin{center} 
		\begin{tikzcd} 
			n \arrow[r, shift left = 1ex, leftarrow, "{\sigma_n}"] \arrow[r, "{\tau_n}", swap, leftarrow] & n-1 \arrow[r, shift left = 1ex, "{\sigma_{n-1}}", leftarrow] \arrow[r, "{\tau_{n-1}}", swap, leftarrow] & \cdots \arrow[r, shift left = 1ex, "{\sigma_2}", leftarrow] \arrow[r, "{\tau_2}", swap, leftarrow] & 1  \arrow[r, shift left = 1ex, "{\sigma_1}", leftarrow] \arrow[r, "{\tau_1}", swap, leftarrow] & 0 
		\end{tikzcd}
		\end{center}
	which moreover satisfy the equations $\sigma_m \circ \sigma_{m-1} = \tau_m \circ \sigma_{m-1}$ and $\sigma_m \circ \tau_{m-1} = \tau_m \circ \tau_{m-1}$ for all $m \in \{2,\dots, n\}$.  
			\end{definition}
	In short, the category of $n$-globular sets can also be defined as the presheaf category $\textbf{Set}^{\mathbb{G}_n^{op}}$. For $X$ an $n$-globular set, we call elements of $X(m)$ the $m$-cells of $X$, e.g., $a \in X(0)$ is a dot or vertex labeled by $a$, $f \in X(1)$ is an arrow with a source and target boundary, $\alpha \in X(2)$ looks just like a natural transformation arrow satisfying certain relations, a 3-cell $x \in X(3)$ an arrow between natural transformation type arrows, etc. In this way, various prominent mathematical constructions including simplicial sets, cubical sets, globular sets\index{globular sets} can be construed as examples of presheaf categories.\par 
	The basic idea in all this is that one selects a category $\textbf{C}$ of cell shapes with morphisms `face inclusions' and `degeneracies'; then, similar to above, one can produce a presheaf category \textbf{Set}$^{\textbf{C}^{op}}$, and the ``boundary extraction" action perspective will generally fit such situations. But this also encourages the (third) view that for a category $\textbf{C}$, its objects can be regarded as \textit{spaces} of a certain sort and its morphisms as structure-preserving morphisms between those spaces, leading to the view that presheaves (valued now not necessarily in \textbf{Set}, but just in some category of ``spaces") on such a category give rise to spaces modeled on \textbf{C} in the sense that they are ``testable" or ``probed" by the objects of \textbf{C}.\index{probe} \par 
	The basic idea with this third perspective can be roughly sketched as follows.\footnote{This general perspective is largely due to Lawvere; see, for instance, \cite{lawvere_taking_2005}.} If \textbf{D} is taken to be, e.g., the category of sets, or simplicial sets, or certain topological spaces, or bornological linear spaces,\footnote{One usually starts with thinking about topological spaces, i.e., the category \textbf{Top}, but really we just need that it is a category and that this category supports some notion of how certain objects can be \textit{covered} by other objects. Much more on this in later chapters.} and if we regard the indexing category $\textbf{C}$ as some category supplying the ``shapes" or ``generic (geometrical) figures"---i.e., whose objects are regarded as ``test spaces," and whose morphisms are the right structure-preserving maps between them---then the presheaf category $\textbf{D}^{\textbf{C}^{op}}$ will be a (generally large) category that will include more general objects that are ``probable" or ``testable" with the help of $\textbf{C}$. Altogether, this can be thought of as amounting to a space ``modeled on" (the objects of) \textbf{C}. \par 
	In referring to ``probes" of a hypothetical space $X$ (just think some generic ``space" for now, not necessarily a topological space) with some ``test space," we are really thinking of all the ways of mapping into $X$ using the objects of \textbf{C}. In other words, if you start with a test space $U$ (an object in \textbf{C}) and are returned a set $X(U)$, we are thinking of this as designating the set of ways $U$ can be mapped into $X$; in this way, with such a set, you can imagine that you have received all the probes or ways of ``testing" $X$ with $U$.\index{probe} Really, then, in describing a generalized space modeled on the objects of \textbf{C}, it looks like we are starting to describe a presheaf $X$ on \textbf{C}, where we regard each such presheaf as a rule assigning to each $U \in \textbf{C}$ the set $X(U)$ of admissible maps from $U$ into the space $X$. However, alone, these probes will not usually give you a very thorough or discriminating understanding of the space $X$. To attain a more complete picture, you also need to know about how the different tests or probes of the space relate to one another. If you have a map $f: U \rightarrow V$ from one object to another, then given some (probe) element $V \xrightarrow{p} X$ of $X(V)$, pre-composing with $f$ (acting on the right) will of course induce a map (of sets, simplicial sets, etc.) going in the other direction $X(U) \xleftarrow{X(f)} X(V)$ that will tell you how $V$-probes of the space $X$ change into $U$-probes; and with \textit{that} information, you can get an accurate picture of what $X$ is. Generally speaking, one of the purposes of doing this is that by ``probing" a ``big" space with a number of ``smaller" or simpler test spaces, we can not only model parts of the space into which we are mapping, but we can ultimately look to piecing together the ``small tests" into information about tests with bigger test spaces, ultimately arriving at a ``picture" of the overall space of interest. This is of special importance to us, since we will see that ultimately the information such probes or tests gathers will turn out to be most useful precisely when the presheaves are in fact sheaves, i.e., satisfy some further consistency conditions. In short, this third (admittedly more subtle) perspective invites us to regard presheaves on \textbf{C} as general (probeable) spaces modeled on \textbf{C}.\footnote{The reader intrigued by this perspective may find the extended discussion in \cite{nlab_authors_motivation_2019} particularly illuminating.}              
	\end{example}
	\begin{example}
		To illustrate the last (but arguably most significant) perspective on the presheaf action---namely, action by \textit{restriction}---we\index{restriction} can begin by considering the construction of a presheaf on the lattice or partial order of open sets $\mathscr{O}(X)$, for $X$ a topological space. A presheaf on $X$ is just a functor $F: \mathscr{O}(X)^{op} \rightarrow \textbf{Set}$. For each open $U \subseteq  X$, we then think of the set $F(U)$ as the set that results from assigning set-values or data throughout or ``over" all of $U$. An open subset $V \subseteq  U$ can be seen in terms of an inclusion arrow $V \hookrightarrow U$ when regarded in the poset category $\mathscr{O}(X)$, so applying the (contravariant) functor $F$ will give us a function that passes from the data assigned throughout or specified over $U$ (the generally ``larger" region) to the data assigned throughout the sub-region $V$, in a process aptly called the \textit{restriction}, and typically denoted $\rho_{V,U}: F(U) \rightarrow F(V)$ (or, $\rho^U_V$, or $F(V \hookrightarrow U)$). Especially when the particular application involves looking at all the \textit{functions} of a certain type (e.g., continuous) defined throughout that region $U$, given an element $f \in F(U)$, one sometimes denotes $\rho^{U}_V(f)$ by $f|_V$, and speaks of the \textit{restriction} of $f$ from $U$ to $V$, as this is treated like the usual restriction of a function along a part of its domain. \par 
		As a first illustration of such a functor, we can consider the set of all continuous real-valued functions, i.e., functions from $U \subseteq  X$ to $\mathbb{R}$. Importantly, when there is an inclusion of opens $V \subseteq  U$, we will have a restriction\index{restriction} function $\rho^U_V: \textbf{Top}(U, \mathbb{R}) \rightarrow \textbf{Top}(V, \mathbb{R})$, which just sends $f: U \rightarrow \mathbb{R}$ to $f|_V: V \rightarrow \mathbb{R}$. The presheaf here thus acts to \textit{restrict} the collection of functions given over some region (say, $(0,6)$) down to the open subsets of that region (say, $(2,4)$ in particular), as in the following:
		\begin{center}
			\includegraphics[scale=0.5]{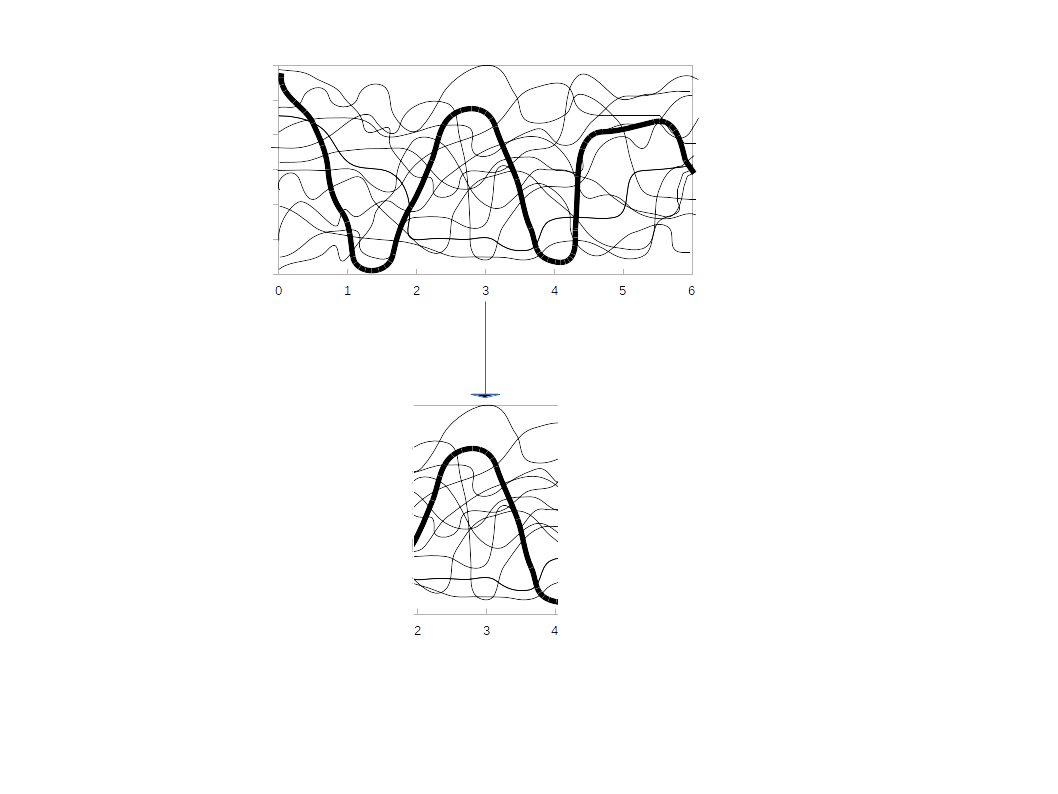}
		\end{center} \vspace*{-1cm} 
		The action of this presheaf is thus given by \textit{restriction}, an action that is clearly functorial.\footnote{Another standard way of producing presheaves on a space arises by taking the ``local sections" of a continuous function $p: E \rightarrow X$, via the local section functor. A local section\index{section!local} of $p$ is a continuous function $s: U \rightarrow E$ from an open subset $U$ of $X$ to $E$, such that $p \circ s (x) = x$ for all $x \in U$. If we let $\Gamma(p)(U) = \{s: U \rightarrow E \hspace*{0.5em} | \hspace*{0.5em} s \text{ is a local section of } p \}$, then by considering that whenever $V \subseteq  U$, we can restrict local sections over $U$ to local sections over $V$, we see that this defines a presheaf on $X$. We will have a lot more to say about this local section approach in a later chapter.}\par 
For another restriction-type example, but of a rather different flavor, take for regions the set $J$ of jurisdictions with their sub-jurisdictions, i.e., $(J, \subseteq )$ a preorder.\footnote{This example comes from \cite{spivak_category_2014}.} We can consider that within the set of \textit{possible laws}---where laws are just propositions, i.e., objects of the preorder $\textbf{Prop}$ regarded as a category whose objects are logical propositions and whose morphisms are \textit{proofs} that one statement implies another---some of these laws are being followed by all people in the region. To each jurisdiction $V$, then, we can assign a set $R(V)$ consisting of whatever laws are being respected by all the people throughout $V$. (In other words: laws are being assigned \textit{locally} to each jurisdiction; after all, a law is dictated to be valid only within a specific region.) If $V$ is a sub-jurisdiction of $U$, i.e., $V \subseteq  U$, then any law respected throughout $U$ is obviously respected throughout $V$, so we can \textit{restrict} from the laws respected throughout $U$ to those respected throughout $V$. Clearly any law respected throughout the state of Illinois will be respected throughout any county included in Illinois, so we can regard such a law given over Illinois from the restricted perspective of a county therein. But observe that the converse is not true! A law respected throughout a part of Illinois need not be respected throughout all of Illinois.\par 
Here we have a local assignment of data to the ``space" of jurisdictions that moreover obeys the property that whenever we have a region $V$ included in $U$, then the action of the presheaf works in the opposite direction: it takes data assignments given throughout $U$ and \textit{restricts}\index{restriction} them to (the same) data assignments now given throughout $V$, the smaller region. The idea to keep in mind here is this: if you have some data (like a ``list" of those laws being respected by everyone) assigned to some region (like Illinois), and you have another list of those laws being respected by everyone in some subregion included in Illinois (like Cook County), then you will expect that the list of laws respected by everyone throughout Cook County will be (equal if not) larger than the list of laws respected by everyone throughout Illinois. In a larger region, there are more chances for the data ``not to fit," e.g., for someone to fail to respect that law, than there is in a smaller region.\par   
The main take-away is that in the previous construction, we have made use of two key ingredients: (1) a local assignment of data to a space (each of the laws in the ``respected laws data" is expected to hold throughout all of the jurisdiction region to which it is assigned); and (2) a natural operation of ``restriction" (induced by the natural ``inclusion" relation governing the overall space of jurisdictions) that allows us to move from the data assigned throughout a region to data assigned throughout subregions. Formally, these two ingredients just specify what we need to have a functor $R$ that is \textit{contravariant}, i.e., we have been describing a presheaf $R: J^{op} \rightarrow \textbf{Prop}$.
	\end{example} 
With an eye towards sheaves, the restriction-style action is in some sense the most decisive of the four perspectives considered, or at least the most immediately relevant to the subsequent initial presentation of sheaves in terms of sheaves on topological spaces. Thus, it pays to understand it very well. The notion of restriction, and its relation to some of the other key basic categorical concepts---for instance, that restrictions are not ``right cancellable" in general, and it is easy to construct many examples of \textit{distinct mappings} that have equal restrictions to the same part, i.e., mappings $f, g$ such that $f|_i = g|_i$ but $f \neq g$---could be dwelt on at greater length. For now, a few general observations concerning restriction may be worth stressing.\par 
Given an inclusion $V \hookrightarrow U$, restriction tells us that we can take some $x \in F(U)$ and ``restrict" that data assignment down to that part of $U$ that makes up $V$, and this will leave us with a viable data assignment on $V$. It should be easy to see, both intuitively and precisely (on the model of function restriction), how this amounts to a ``restriction."\index{restriction} However, it is important to recognize that if one understands (as some newcomers to (pre)sheaves occasionally seem to do) the language of ``restriction" to suggest that, at the level of the presheaf itself, we are passing from an (in general) ``bigger" set of data to a ``smaller" one---``restricting our attention" as it were---strictly speaking, at the level of the maps between the presheaf data $F(U)$ and $F(V)$, this is \textit{not} what is going on. This should already be evident from close consideration of the ``laws" example, where the set of respected laws $R(U)$ specified over a larger region $U$ will actually typically be \textit{smaller} than the set of respected laws $R(V)$ specified over the sub-region $V \subseteq U$. Same with the continuous functions: it is ``easier" to be continuous on a smaller region, i.e., over a bigger region, there are more opportunities to fail to be continuous. It is perhaps an unnecessary warning, but the point is that, at the level of the presheaf maps themselves, e.g., moving from $R(U)$ to $R(V)$, we are not generally dealing with a restriction in the sense of moving from a bigger \textit{data set} (set of value assignments) down to a smaller data set. Confining our attention to $U$ and $V$ as regions or components of a space, it makes sense to think of these objects (regions) as \textit{constraints} of sorts, according to which $V$, a ``smaller" region, amounts to a weaker constraint on any data specified locally over the regions. It should be evident that given any inclusion of a ``smaller" region into a ``larger" region, more data will generally be able to satisfy the weaker constraint (corresponding to the smaller region) than will satisfy the stronger one (corresponding to the larger region). Yet, at the level of a particular data assignment, we can regard the presheaf maps as amounting to a restriction of that data along inclusions of sub-regions. The next and final example of this chapter should help to further clarify this.
\begin{example}
	Consider time intervals as objects and morphisms given by inclusions, yielding a category we will denote $\mathcal{T}$.\footnote{Really, looking ahead to the fact that there is a sheaf lurking in this example, we would want to be constructing a \textit{site} here (using basically the same data), where we also specify that an interval $[t,u]$ is \textit{covered} by a collection of intervals $\{[t_i, u_i] \hspace*{0.2em}|\hspace*{0.2em} i \in I \}$ provided $\cup_i [t_i, u_i] = [t,u]$. There are also some subtleties in this construction that we freely ignore for the moment. We will return to this example, and such matters, in a later chapter.} For concreteness, suppose we consider the period spanning from January 1 (at midnight, or 0:00) of 2018 until June 1 (at 0:00) of 2018. Then, suppose this period is decomposed into various 2 month intervals 
	\begin{equation*} 
	[Jan1, March1], [Feb1, April1], [March1, May1], [April1, June1], 
	\end{equation*}
	that together ``cover" the entire half-year period from $Jan1$ through to the end of $May$.\footnote{For now, you can think of this notion of ``covering" in a naive way. We will be precise about this sort of thing in the future.} There are, moreover, the natural overlapping subintervals, producing the following overall structure on the system of intervals ordered by inclusion (as indicated by the inclusion arrows):  
	\begin{center}
		\includegraphics*[scale= 0.3]{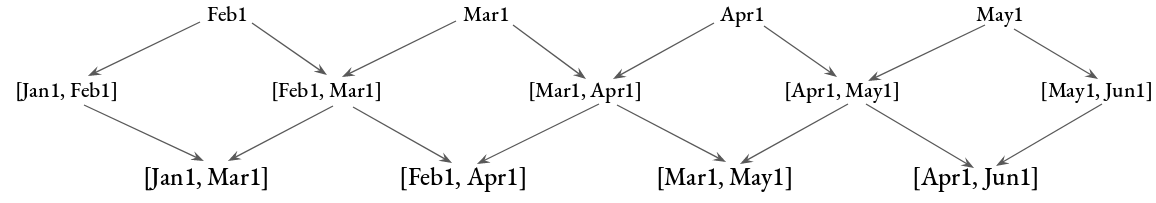}
	\end{center}
and where the single dates like Feb1 are short for the ``degenerate interval" $[Feb1, Feb1]$ representing the ``instant" 0:00 on February 1. \par 
Now, for a particular company with some (generally fluctuating) stockpile of products, we can define a contravariant functor $S: \mathcal{T}^{op} \rightarrow \textbf{Set}$ that assigns to each time period $[t,u]$ the products that are in stock throughout the \textit{entire interval} $[t,u]$ (where, for the moment and for simplicity, we just imagine that a product is simply present or absent, say, as if the only important question was whether they had at least one item of the product or had none, ignoring the question of quantity). $S([t,u])$ is thus just the set of products the company has in stock throughout the entirety of the time-period $[t,u]$. \par 
Then, for any inclusion of time periods $i: [t,u] \hookrightarrow [v,w]$, the functor $S$ acts (contravariantly) \textit{by restriction},\index{restriction} mapping each stocked item onto itself. Clearly, any product present in the company's store throughout the bigger time interval $[v,w]$ must be present as well throughout any sub-interval $[t,u]$. But this tells us that the ``list" of products recorded as present throughout the bigger time interval $[v,w]$ will in general likely be shorter or smaller than the list of products assigned to the smaller time interval $[t,u]$. \par 
A particular presheaf $S$ on such a $\mathcal{T}$ might then be given by something like the following (where each of the $A, B, C,$ etc., sitting over each interval-object, represents one of the products held by the company throughout the entire interval):\footnote{This presheaf will turn out to already be a sheaf, as we will see. The idea is already accessible, however, as long as one appreciates that sheaves employ the notion of \textit{covers}, and that being a sheaf basically involves a `gluing condition' that requires that whenever a property is locally true or valid throughout a cover of an object (time intervals in this case), then it holds over the entire object. Letting the ``subintervals" $\{[t_i, u_i] \hspace*{0.2em}|\hspace*{0.2em} i \in I \}$ cover the interval $[t,u]$, it is practically immediate that if a product is present in the company's stockpile throughout each of the intervals $[t_i, u_i]$, then it will have to be present throughout all of $[t,u]$.} 
\begin{center}
	\includegraphics*[scale=0.3]{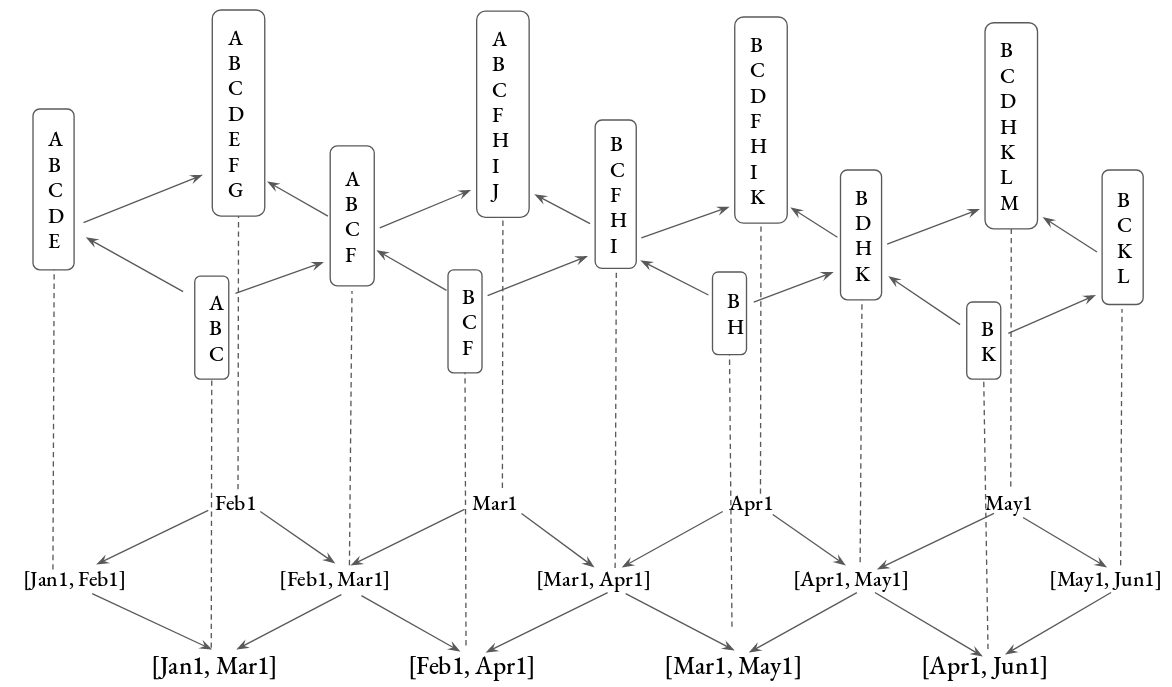}
\end{center}
As one can easily see by inspecting the diagram, the sets of value assignments (of products present) throughout each interval are generally smaller over larger regions (time intervals); thus, strictly attending just to the diagram sitting ``above," the arrows in fact generally go from smaller sets of data to larger ones.\footnote{This relates to the warning discussed just before this example.} \par 
A final thing to realize is that, in general, restriction\index{restriction} along an inclusion is \textit{not necessarily either surjective or injective} (despite what a naive understanding of the language of restriction might seem to imply, for instance, suggesting at least surjectivity). An easy counter-example is provided by the following.\footnote{This counter-example is derived from \cite{lawvere_sets_2003}.} As was seen in an earlier example, restriction of continuous functions is continuous. However, take the (poset) category $\textbf{A}$ that consists of just two objects, $U$ and $C$ with the single non-identity (inclusion) map $U \rightarrow C$, where $U$ is the open interval $(0,1)$ and $C$ is the closed interval $[-1,1]$. Now let the presheaf $F: \textbf{A}^{op} \rightarrow \textbf{Set}$ act on objects as follows: $F(U)$ is the set of all continuous real-valued functions given over the open interval $(0,1)$; and $F(C)$ is the set of all continuous real-valued functions over the closed interval $[-1,1]$. Then the induced single presheaf action $F(C) \rightarrow F(U)$ is clearly by restriction; yet, it should be evident that this particular restriction process can be neither surjective nor injective.\footnote{Non-surjective since there are functions that remain continuous over $(0,1)$ while having discontinuities at either or both ``endpoints"---in particular, at $0$---so that such functions cannot come from any continuous functions specified over all of $[-1,1]$; non-injective since there exist distinct continuous functions given over all of $[-1,1]$, each of whose restrictions down to $(0,1)$ are identical.}    
\end{example}  
\subsection{Philosophical Pass: The Four Action Perspectives}
Many important mathematical structures and categories---including some of those already discussed, e.g., dynamical systems, bouquets, graphs, hypergraphs, and simplicial sets---arise as a presheaf category consisting of contravariant functors on some given indexing category, landing in \textbf{Set}, where the result of applying the functor to the objects of the indexing category yields something we naturally think of as ``containers" (in \textbf{Set}) holding on to various manifestations or figures each of which conform to the ``shape" or ``form" or ``stage" determined by the ``generic figures" populating the indexing category (one for each of its objects), and where the ``changes of figure" indicated by the indexing category (given by its morphisms) are respected by the figures of the container. While this perspective is perhaps most appropriate, or easy to countenance, when the objects and morphisms of the indexing category $\textbf{C}$ have some sort of geometrical interpretation, it is a surprisingly useful perspective even in more general cases. \par
As for the four perspectives on the presheaf action: the fundamental idea shared by each of these is that the domain category \textbf{C} plays the role of specifying the general internal structure or schema---in the form of the figure-types or shapes, the glue, the nature of the internal dynamic, or the locality of data assignments in which all the sets in $\textbf{Set}^{\textbf{C}^{op}}$ must participate. In each case, the resulting entire presheaf category has for its objects all the different instantiations or ``realizations" that each exemplify or mobilize a particular way of realizing the general form supplied by \textbf{C}. The other way of looking at this is that the domain category plays the role of a parameter specifying in an invariant form how (temporal, dynamic, geometric) variation is to take place, while the target category (\textbf{Set}) serves as the container or arena holding on to all the particular values or results of trying to ``participate in" or ``realize" this form of variation. One might accordingly think of a presheaf itself as mediating between the invariance or fixedness ``outside of time" belonging to the domain category, on the one hand, and its multifarious concrete presentations or manifestations ``in time," on the other.  \par   
It is worth emphasizing that while all the presheaf functors above are valued in \textbf{Set}, which can be useful in taming many problems or otherwise complicated structures, philosophically-speaking, presheaves are anything but the static and qualitatively-barren objects the usual naive set-theoretical perspective on sets as ``bags of fixed objects" might encourage us to believe. \par 
Considering presheaves as coming equipped with an action that is \textit{processual}\index{process} recaptures a dynamic perspective in which objects are not regarded as static collections, but instead what is presented to us are sets seen as evolving through stages, either merely temporally or in accordance with an internal dynamic (as in case of motions or evolutions of a certain shape, subject to certain equations, fixed points, etc.). Against the generally ``discrete" and static context of sets, this restores a more ``continuous" and dynamic perspective. The \textit{boundary-extraction}\index{boundary!extraction} perspective, for its part, reveals the incidence relations, which relations altogether describe something like how the overall structure "holds together" or ``coheres." As against the usual set-theoretical perspective of sets of objects grouped together more or less arbitrarily into a set that cannot internally differentiate objects or discern important qualitative (or dimensional) features of those objects or their modes of relation, this perspective restores a kind of ``continuity" in telling us how the various components of a structure can be regarded as being ``glued" or ``stitched" together from other (lower-dimensional) components. The third perspective lets us regard a space in terms of all the ways of ``probing"\index{probe} it from the ``outside" and thinking about the entire space in terms of how these various probes behave with respect to one another. This perspective, similar to that of the ``relationism" of Yoneda, is a more ``continuous" one in the sense of insisting on understanding something in terms of all the relations or perspectives on it. Finally, acting via \textit{restriction},\index{restriction} presheaves open onto a range of relationships between the parts of a whole. In general, such a relationship emerges as ``regular" in the sense that in passing from data specified over some containing region to data over a sub-region, there remains a kind of conformity or identity of the data specified over the parts in relation to the same rule or function describing the containing region. This perspective thus opens onto the notion of a conformity of parts of a whole to a single rule or idea (as opposed to the usual set-theoretical consideration of a ``whole" independently of the specific way, beyond whether or not a part ``belongs," the whole enforces relationships among the parts). \par 
In short, while we are able to benefit from certain ``nice" and ``tame" properties of \textbf{Set}, the presheaf perspective lets us recapture a generally more dynamic, nuanced, and ``continuous" perspective on many structures of interest.
\chapter{First Look at Sheaves}
		A number of the examples of presheaves already introduced---including in particular the examples falling under the heading of ``restriction" (e.g., the continuous functions,  laws respected throughout jurisdictions, and the company stockpile examples)---are \textit{sheaves}. This chapter turns, at long last, to a first presentation of sheaves, in the course of which will be seen a number of initial and intuition-building examples. We will build up gradually towards more and more involved examples. \par 
		While we are very close to sheaves already, in order to provide a definition of a sheaf in the setting of a topological space we need to introduce one last notion, namely that of a \textit{covering}. We will return to (and generalize) this notion of a cover in later chapters, making use of the more general notion of a \textit{site}. For now, we will just work with the easiest examples of sites, via the ``usual" notion of a cover, as it appears in the context of topological spaces and their open sets. 
		\begin{definition}
			Given $X$ a topological space and $U \subseteq  X$ an open set of the space, consider $V_1,\dots, V_n$ open subsets (think ``subregions") of $U$, i.e., for all $1 \leq i \leq n$, we have $V_i \subseteq  U$. Then the $V_i$ are said to collectively \textit{cover},\index{cover} or provide a \textit{covering} of, $U$ if every point that is in $U$ is in a $V_i$ for some $i$. 
		\end{definition} \noindent 
		Another (only slightly more sophisticated) way of saying this is that for $\mathscr{O}(X)$, the poset of open subsets of $X$, ordered by inclusion, an $I$-indexed family of open subsets $V_i \hookrightarrow U$ \textit{covers} $U$ provided the full diagram consisting of the sets $V_i$ together with the inclusions of all their pairwise intersections 
		\begin{center} 
		\begin{tikzcd} 
		V_i \arrow[r, hookleftarrow] & V_i \cap V_j \arrow[r, hookrightarrow] & V_j
		\end{tikzcd} 
		\end{center} 
		has $U$ for its colimit (it is okay, for now, to just think of this in terms of $U$ being the coproduct $\bigcup_{i \in I} V_i$ in $\mathscr{O}(X)$). \par 
		 Roughly, one can think of a covering of a given object $U$ as some sort of \textit{decomposition} of that object into simpler ones, the resulting simpler ``pieces" of which, when taken altogether, can be used to recompose all of $U$. In terms of covers of a \textit{set} $U$, this has a very simple description: it is just a family of subsets $\{V_i\}_{i \in I}$ whose union is $U$ itself, i.e., $\bigcup_{i \in I} V_i = U$. At the outset, it is perfectly fine to just think of this in terms of specifying a collection of subregions that can be ``laid over" a given region in such a way that the entire region is thereby covered, where an entirely obvious but still decisive observation is that such subregions making up the cover can \textit{overlap} one another. The naive image to keep in mind is that we have a region $U$ that we want to cover with some collection of ``pieces" into which it may be regarded as being ``decomposed." Suppose we have some $V_1 \subseteq  U$ and $V_2 \subseteq  U$:
		\begin{center} 
		\includegraphics*[scale=0.23]{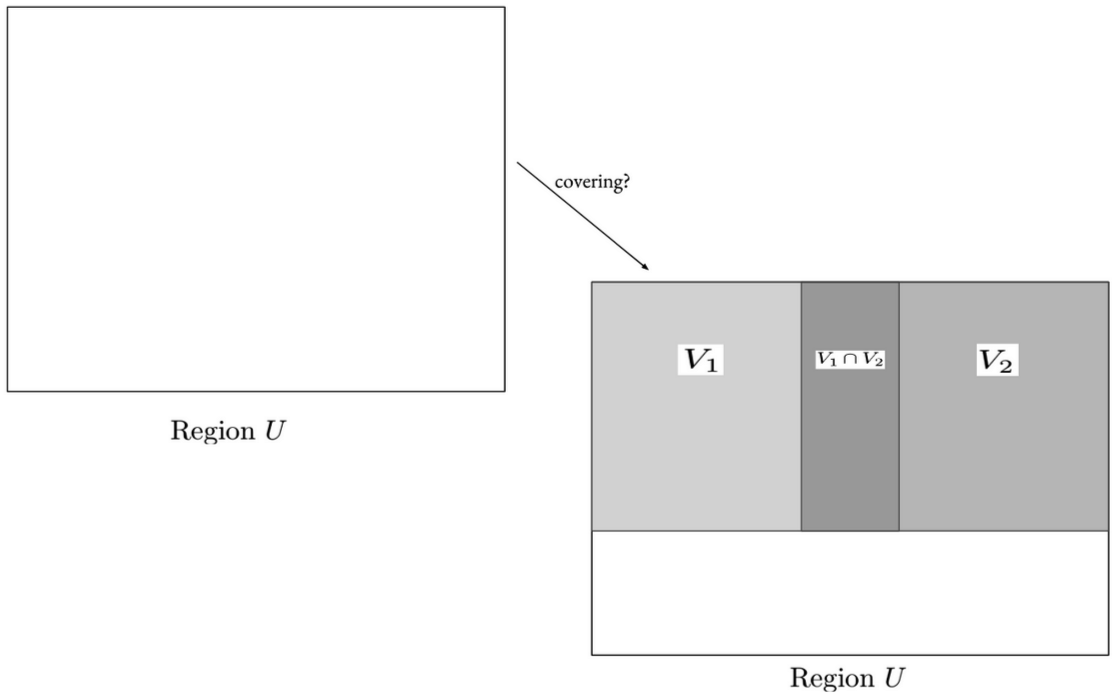} 
		\end{center} \par \noindent
		Clearly, $V_1$ and $V_2$ collectively fail to cover $U$, yet we can observe that there is a subregion where $V_1$ and $V_2$ overlap, which we call $V_1 \cap V_2$ and regard as specifying another ``piece." Since $V_1$ and $V_2$ collectively cover more of $U$ than either does individually, we should also consider the larger region (the entire northern half of $U$) that results from joining $V_1$ and $V_2$. 
		We might continue in this manner, working our way up to a collection of subregions of $U$ that actually cover all of $U$. For instance, we might have another $V_3$, laid on top of the entire southern half of the region (and partly overlapping with each of $V_1$ and $V_2$), such that the entire region $U$ is now covered by the collection $\{V_1, V_2, V_3\}$. Altogether, the data of such a system of open sets, ordered by subset inclusion, will have the structure of a poset (this means, in particular, that we can regard $\mathscr{O}(X)$ as a category). In our particular case, this could be displayed by the diagram: 
			\begin{center}
			\begin{tikzpicture}[yscale=0.65, xscale=0.85]
			\node (max) at (0,4) {$U$};
			\node (a) at (-2.4,2) {$V_1 \cup V_2$};
			\node (b) at (0,2) {$V_1 \cup V_3$};
			\node (c) at (2.4,2) {$V_2 \cup V_3$};
			\node (d) at (-2,0) {$V_1$};
			\node (e) at (0,0) {$V_2$};
			\node (f) at (2,0) {$V_3$};
			\node (g) at (-2.4,-2) {$V_1 \cap V_2$};
			\node (h) at (0,-2) {$V_1 \cap V_3$};
			\node (i) at (2.4,-2) {$V_2 \cap V_3$};
			\node (min) at (0,-4) {$V_1 \cap V_2 \cap V_3$};
			\draw[->] (a) -- (max);
			\draw[->] (b) -- (max);
			\draw[->] (c) -- (max);
			\draw[->] (d) -- (a);
			\draw[->] (d) -- (b);
			\draw[->] (e) -- (a);
			\draw[->] (e) -- (c);
			\draw[->] (f) -- (b);
			\draw[->] (f) -- (c);
			\draw[->] (min) -- (g);
			\draw[->] (min) -- (h);
			\draw[->] (min) -- (i);
			\draw[->] (g) -- (d);
			\draw[->] (g) -- (e);
			\draw[->] (h) -- (d);
			\draw[->] (i) -- (e);
			\draw[->] (h) -- (f);
			\draw[->] (i) -- (f);
			
			\node (dom) at (0,-6) {$\mathscr{O}(X)$};
		\end{tikzpicture}
		\end{center}   
		revealing the components of the space, together with their relevant inclusion relationships as members of a cover of the entire space.\par 
		Sheaves on a topological space---the first sort of sheaf we will consider---will just be particular presheaves on the open subsets $\mathscr{O}(X)$, presheaves that satisfy a further property. What ultimately distinguishes presheaves and sheaves is that sheaves are a special kind of presheaf, one that is ``sensitive to" the information or structure of a covering. The notion of a presheaf together with the notion of a covering of a topological space, supply us with all the ingredients needed to offer a first pass at the sheaf concept. 
		\section{Sheaves: The Topological Definition} 
		The definition of a sheaf that we will presently be concerned with---sheaves on topological spaces---is a definition very much motivated by the action-as-\textit{restriction} presheaf perspective. It may be thought of as involving a choice of data from each of the sets (as assigned by the presheaf to each piece of the cover), that moreover forms a locally compatible family (meaning that it respects the restriction mappings and the chosen data items agree whenever two pieces of the covering overlap) and together induce or extend to a \textit{unique} choice over the entire space being covered. \par 
		Since our sheaf candidate will already be a presheaf, to determine whether or not a given presheaf is a sheaf will just amount to testing for certain properties of a set-valued functor $F: \mathscr{O}(X)^{op} \rightarrow \textbf{Set}$. That we have a presheaf on $\mathscr{O}(X)$ means, first of all, that to every open set $U$ in $\mathscr{O}(X)$, we will have a set $F(U)$. Terminologically (a terminology to be explained further on in this chapter), we will sometimes refer to the set $F(U)$ as the \textit{stalk} at (or over) $U$, or as the \textit{set of sections over} $U$, and to an element $s$ in this $F(U)$ as a \textit{section over} $U$.\footnote{Strictly speaking, elements of $F(U)$, i.e., value assignments specified over $U \subseteq X$, are called \textit{local sections}\index{section!local} of the sheaf $F$ over $U$, to distinguish these from elements of $F(X)$, i.e., value assignments given over the \textit{entire} space, which are called \textit{global sections}\index{section!global} of $F$. More broadly, whenever local information, e.g., elements like functions $f$ and $g$ given over certain domains, restricts to the same element in the intersection of their domains, then such $f$ and $g$ are called \textit{sections}.\index{section} We will explore and better motivate all this terminology in more detail later on; meanwhile, we may make use of the convenience of such language. For now, one thing to realize is that, for an arbitrary presheaf on a space, the set of global sections of the presheaf on the overall space may be different from the set of local sections given over all the open subsets; the ``gluing" axiom (discussed below) is there precisely to enforce that this difference disappears.} Moreover, corresponding to every inclusion of open sets $V \hookrightarrow U$, that $F$ is a presheaf means that we will have a \textit{restriction} $F(U) \rightarrow F(V)$. Supposing $s \in F(U)$ is a section over $U$, it is common to denote its restriction to $V$ by $s|_V$, i.e., $F(V \hookrightarrow U): F(U) \rightarrow F(V)$  takes $s \mapsto s|_V$ for each $s \in FU$, where we treat this like the usual restriction of a function. Observe that whenever we have three nested open sets $W \subseteq V \subseteq U$, restriction will be transitive, i.e., $(s|_V) |_W = s|_W$. \par 
		Altogether, with these two pieces of data, 
		\begin{equation}
		U \mapsto FU, \hspace{2em} \{V \subseteq U\} \mapsto \{FU \rightarrow FV \text{ via } s \mapsto s|_V \}
		\end{equation}
		we are just reiterating that we have a functor from $\mathscr{O}(X)^{op}$ to $\textbf{Set}$. Together with the notion of coverings, we are now ready to define a sheaf. 
			\begin{definition}
			(\textit{Definition of a Sheaf}) Assume given $X$ a topological space, with $\mathscr{O}(X)$ its partial order of open sets, and $F: \mathscr{O}(X)^{op} \rightarrow \textbf{Set}$ a presheaf. Then given an open set $U \subseteq  X$ and a collection $\{U_i\}_{i \in I}$ of open sets covering $U = \cup_{i \in I} U_i$, we can define the following \textit{sheaf condition}:\index{sheaf!defined} 
			\begin{itemize}
				\item Given a family of sections $a_1, \dots, a_n$, where each $a_i \in F(U_i)$ is a value assignment (section) over $U_i$, whenever we have that for all $i, j$,
				\begin{equation*} a_i|_{U_i \cap U_j} = a_j|_{U_i \cap U_j}, 
				\end{equation*}
				then there exists a \textit{unique} value assignment (section) $a \in F(U)$ such that $a|_{U_i} = a_i$ for all $i$. 
			\end{itemize}  
			Whenever there exists such a unique $a \in F(U)$ for every such family, we say that $F$ satisfies the sheaf condition\index{sheaf!condition} for the cover $U = \cup_{i \in I} U_i$.
			The presheaf $F$ is then a \textit{sheaf} (full stop) whenever it satisfies this sheaf condition for \textit{every} cover. 
		\end{definition}
			Let us break this definition down into four, more easily ``digestible," steps. The idea is this: given a presheaf on some space and a covering, the definition of a sheaf begins by making use of what is sometimes called a \textit{matching family}:
			\begin{definition}
				A \textit{matching family}\index{matching family} $\{a_i\}_{i \in I}$ of sections over $\{U_i\}_{i \in I}$ consists of a section $a_i$ in $F(U_i)$ for each $i$---chosen from the entire set $F(U_i)$ of all sections over $U_i$---such that for every $i, j$, we have 
				\begin{equation*}
				a_i|_{U_i \cap U_j} = a_j|_{U_i \cap U_j}.
				\end{equation*} 
				In other words: given a data assignment $a_i$ throughout or over region $U_i$ and a data assignment $a_j$ over region $U_j$, if there is agreement or consistency between the different data assignments when these are restricted to the sub-region where $U_i$ and $U_j$ overlap, then together the data assignments $a_i, a_j$ give a matching family. As the definition requires that it holds for every $i, j$, the idea is that we can build up ``large" matching families of sections via such pairwise checks for agreement. 
			\end{definition}  
			Digesting the definition of matching family is the first step in grasping the definition of a sheaf. \par 
			Next, the definition specifies what is sometimes called a \textit{gluing} (or the \textit{existence} condition).\index{sheaf!gluing condition} Given a matching family for our cover of the space $U$, we call a section over $U$ itself a \textit{gluing} if, whenever this data assignment over all of $U$ is restricted back down to each of the subregions or pieces that make up the cover of the entire object, it is equal to the original local data assigned to each subregion. The definition stipulates that such a gluing $a \in F(U)$ exists.\par 
			Not only does such a gluing exist, but to have a sheaf, we require that there is a \textit{uniqueness} condition.\index{sheaf!uniqueness condition} Specifically, there exists a \textit{unique} section $a \in F(U)$ such that $a|_{U_i} = a_i$ for all $i$. In other words, if $a, s \in F(U)$ are two sections of $F(U)$ such that $a|_{U_i} = s|_{U_i}$ for all $i$, i.e., they are equivalent along all their restrictions, then in fact $a = s$ (i.e., they must be the same, so we have \textit{at most one} $a$ with restrictions $a|_{U_i} = a_i$). \par 
			With the notion of a matching family and that of a unique gluing, we can form the notion of the sheaf condition, and the definition is basically complete. The idea being that if \textit{for every matching family, there exists a unique gluing}, then we say that the presheaf $F$ satisfies the \textit{sheaf condition}. A presheaf will then be a sheaf whenever it satisfies this sheaf condition for every cover. That is all---we have defined what a sheaf is! \par 
	Stepping back, we can accordingly break down the definition of a sheaf on a topological space into a particular presheaf that moreover satisfies two conditions with respect to a cover: (1) existence (or gluing); and (2) uniqueness (or locality). 
	\begin{definition}
		(\textit{Definition of a Sheaf (again)})\index{sheaf!defined}
		Given a presheaf $F: \mathscr{O}(X)^{op} \rightarrow \textbf{Set}$, an open set $U$ with open cover by $\{U_i\}_{i \in I}$, and an $I$-indexed family $s_i \in F(U_i)$, then $F$ is a sheaf provided it satisfies both: 
			\begin{enumerate}
				\item (\textit{Existence/Gluing})\index{sheaf!gluing condition} If, for each $i$, there is a section $s_i \in F(U_i)$ satisfying that for each pair $U_i$ and $U_j$ the restrictions of $s_i$ and $s_j$ to the overlap $U_i \cap U_j$ match (or are ``compatible")---in the sense that   
				\begin{equation*}
				s_i x = s_j x
				\end{equation*}
				for all $x \in U_i \cap U_j$ and all $i, j$---then there \textit{exists} a section $s \in F(U)$ with restrictions $s |_{U_i} = s_i$ for all $i$. (Here, such an $s$ is then called the \textit{gluing}, and the $s_i$ are called \textit{compatible}.) 
				\item (\textit{Uniqueness/Locality})\index{sheaf!uniqueness condition} if $s, t \in F(U)$ are such that 
				\begin{equation*} 
				s|_{U_i} = t|_{U_i}
				\end{equation*}
				for all $i$, 
				then $s = t$. (In other words, there is \textit{at most one} $s$ with restrictions $s|_{U_i} = s_i$.) 
			\end{enumerate} 
	Together, these two axioms assert that compatible sections can be uniquely glued together. 
	\end{definition} 
\noindent If $F$ and $G$ are sheaves on a space $X$, then a morphism $f: F \rightarrow G$ will just be a natural transformation between the underlying presheaves. This lets us define $\textbf{Sh}(X)$ the \textit{category of sheaves} on $X$, which has sheaves for objects and natural transformations for morphisms. There are in general far more presheaves on a space than there are sheaves on the space. This category of sheaves on $X$ will be a (full) subcategory of the category of presheaves on $X$, giving the inclusion functor 
\begin{equation*}
\iota: \textbf{Sh}(X) \rightarrow \textbf{Set}^{\mathscr{O}(X)^{op}}. 
\end{equation*}
\subsection{A Sheaf as Restriction-Collation}
\label{sec: restrict}
Before launching into examples, let us briefly consider how the description of a sheaf (of sets) on a topological space can be motivated by simple observations concerning functions. We know that specifying a topology on a set $X$ lets us define which functions are continuous, such as the continuous functions from the space $X$ (or some open $U \subseteq  X$) to the reals $\mathbb{R}$. Whether or not a function $f: U \rightarrow \mathbb{R}$ is continuous is something that can be \textit{determined locally}. But what exactly does this mean? This fundamentally amounts to saying two things:\footnote{This perspective of restriction-collation is derived from \cite{maclane_sheaves_1994}.}     
\begin{enumerate}
	\item \textbf{Restriction (or Identity)}:\index{restriction} If $f: U \rightarrow \mathbb{R}$ is continuous, and $V \subseteq U$ is open, then restricting $f$ to $V$, i.e., $f|_V : V \rightarrow \mathbb{R}$, yields a continuous function as well.   
	\item \textbf{Uniquely collatable (or Gluability)}: If $U$ is covered by open sets $U_i$, and the functions $f_i: U_i \rightarrow \mathbb{R}$ are continuous for all $i \in I$, then there will be \textit{at most one} continuous $f: U \rightarrow \mathbb{R}$ with restrictions $f |_{U_i} = f_i$ for all $i$. Furthermore, this $f$ will exist iff the given $f_i$ \textit{match} on all the overlaps $U_i \cap U_j$ for all $i, j$, i.e., $f_i x = f_j x$ for all $x \in U_i \cap U_j$.
\end{enumerate}  
\noindent 
One might accordingly think about the ``localness" of a given function's property (such as continuity) as involving two sorts of compatibility conditions or constraints tending in two different directions (the first ``downward" and the second ``upward"): (1) that which requires that information specified over a larger set is compatible whenever \textit{restricted} to information over a smaller open set; (2) that which involves conditions on the \textit{assembly} of matching information on smaller opens into information given over larger open sets. One might also think of the first condition as the ``localizing" part, and the second condition as the ``globalizing" part. \par 
While continuous functions provide a particularly natural example of these sorts of requirements, there is no need to restrict ourselves to \textit{continuous} functions. Various things such as differentiability, real analyticity, and other structures on a space $X$ (including involving things that are not even functions, but are ``function-like") are in fact ``determined locally" in the same sort of way. The underlying idea here is that certain functions (or things that behave like functions), thought of as having some property $P$, are defined on the open sets in such a way that one can check for this property in a neighborhood of every point of the space---this is fundamentally what makes it \textit{local}---and then each inclusion $V \subseteq U$ of open sets in $X$ will determine a function $P(V \hookrightarrow U): PU \rightarrow PV$, for which we just write $t \mapsto t|_V$ for each $t \in PU$, treating it like the usual restriction of a function (which restriction is, moreover, transitive). Altogether, this just says that we have defined a functor $P: \mathscr{O}(X)^{op} \rightarrow \textbf{Set}$; and saying that $P$ is such a presheaf (functor) simply expresses the first (restriction) condition given above. The second condition mentioned above, unique collatability or gluability, can in turn be described category-theoretically in terms of an equalizer diagram for a corresponding covering.\footnote{We will see how this works in a moment.} Accordingly, the two requirements of restriction and unique collatability supply the model for how to define a sheaf more generally (as a functor for which the corresponding equalizer diagram, containing the information of the open sets and the cover, is an equalizer for all coverings). This description ultimately enables the sheaf construction for a wide class of structures. But this motivation in terms of certain properties of classes of functions being checked locally is a particularly useful perspective to keep in mind as one thinks about the construction of sheaves in general.\par 
	Let us now dive right into some examples of sheaves. Over the course of the book, we will provide a multitude of examples, ranging from the more intuitive to the computationally-explicit and involved. With the next examples, we start with a couple of simple sheaves, occasionally omitting some of the details, meant only to develop some initial \textit{intuition} for the sheaf concept and to leave the reader with a number of suggestive pictures and guiding examples. In the sections and chapters that follow, more elaborate and complicated examples are given.
	\section{Examples}
	\begin{example}
		We return to the presheaf of continuous real-valued functions on a topological space $X$, as discussed in the previous chapter. This is a sheaf, specifically a sheaf of real algebras associating to each open $U \subseteq  X$ the algebra $F(U)$ of real-valued continuous functions defined there. Not only can we restrict functions down to any open subset, but we can also glue together local assignments whenever they agree on overlapping regions, producing a global assignment, i.e., a consistent assignment \textit{over the entire region} that will agree with the local assignments when restricted back down to each subregion. Uniqueness in this case is automatic from the fact that we are dealing with \textit{functions}. \par 
		Overall, this process of the pairwise compatibility checks and the subsequent gluing is nicely captured by an image of the following sort (where, for each piece, we just depict the choice that has been made from the overall set of all continuous functions over that region):      
		\begin{center}  
		\includegraphics[scale=0.45]{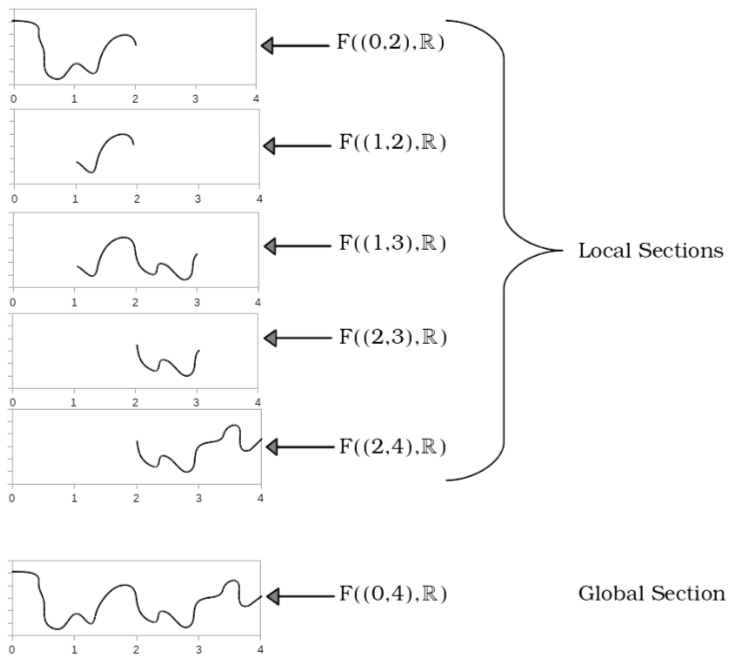}
		\end{center} 
	\end{example}
	\begin{example}
	Revisiting the example of the presheaf of laws being respected throughout a jurisdiction (a geographic area over which some legal authority extends): for $X$ the entire world, to each jurisdiction $U \subseteq  X$ we assigned the set $R(U)$ of laws being respected throughout the region $U$. Is this presheaf $R$ a sheaf? Well, we can check: given some law respected throughout $U$ and another law respected throughout $W$, do they amount to the same law on the sub-region where $U$ and $W$ overlap?\footnote{In fact, as we will see later on, we do not strictly need that they are exactly the same law, just that there is a consistent system of ``translation" between the sets of laws, i.e., a set of isomorphisms translating between each such pairs of sets of laws.} (If there is no overlap, then this is trivially satisfied.) Now repeat this check for each such pair of overlapping regions.  \par 
	For instance, on $U$ there might be a law that stipulates ``no construction near sources of potable water," while on $W$ a law might stipulate ``no construction in public parks." If it turns out that on the overlapping subregion $U \cap W$ all public parks are near sources of potable water (and vice versa), then the laws agree on that overlapping region, and thus can be ``glued" together to form a \textit{single} law about construction that holds throughout the union $U \cup W$ of the two.  \par 
This might seem like a rather harmless or trivial construction, but consider that the global sections of such a sheaf $R$ would tell you exactly those laws that are respected by everyone throughout the planet. This would be a useful piece of information! (For instance, it might reveal the sorts of shared values that are ultimately respected, in one form or another, by every society.) The process of ``checking" for agreement on overlapping regions is straightforward, but the resulting observations or data assignments one can now make concerning the entire space, via the global sections, can be very powerful and far-reaching.   
\end{example}
\begin{example}
	Returning to the presheaf $S: \mathcal{T}^{op} \rightarrow \textbf{Set}$ of a company's stockpile of products, we can describe a sheaf here. If $\{[t_i, u_i] \hspace*{0.25em}| \hspace*{0.25em} i \in I \}$ covers the entire interval $[t,u]$---in our particular case, the interval from January 1 until June 1---it is practically immediate that if a product is present in the company's stockpile throughout each of the pieces $[t_i, u_i]$ of the cover, then it will have to be present throughout all of $[t,u]$, and for any inclusion $i$ of intervals, $S(i)$ is the restriction function mapping each product onto itself (any product found throughout the larger interval must clearly be present throughout a sub-interval). Global sections will then be given by those products that are consistently present throughout the entire time period, such as product $B$, as depicted below: 
	\begin{center}
		\includegraphics*[scale=0.25]{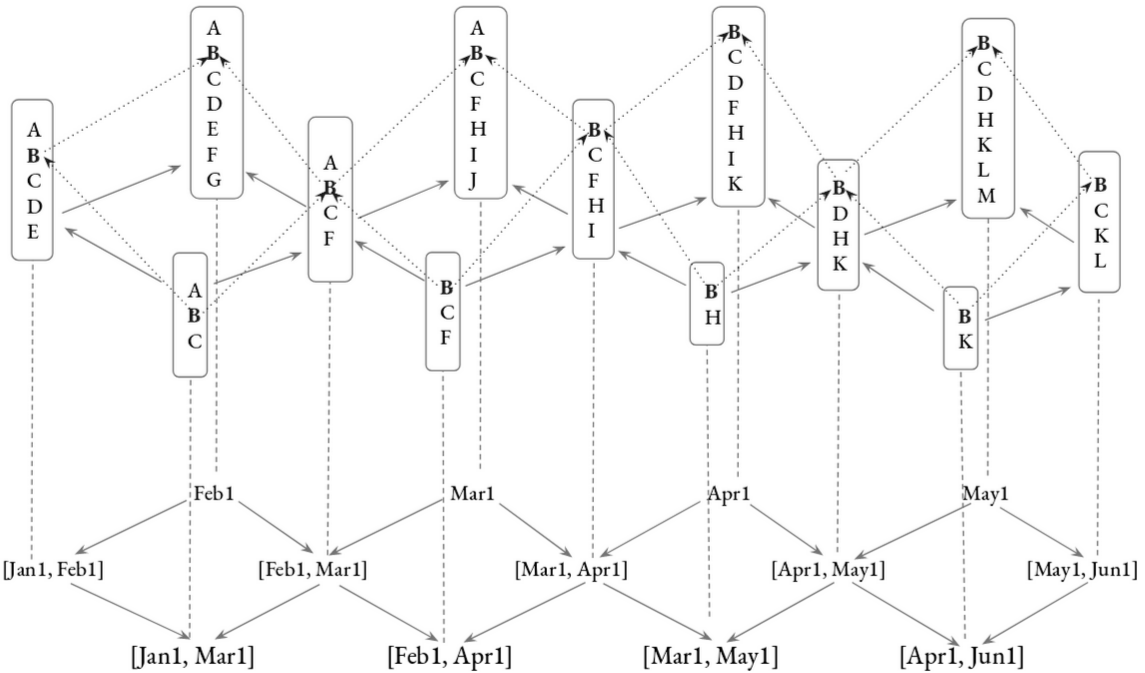}
	\end{center}  
	This may not be a very ``exciting" sheaf, but its simplicity can be useful in helping one achieve an initial working understanding of the difference between local sections that can extend to global sections and local sections that are purely local and satisfy certain local compatibility checks but cannot be glued together into a global section. To see this, suppose, for instance, we had instead selected the product $C$, which indeed appears to be present throughout much of the overall time period. It is certainly present throughout all of $[Jan1, Apr1] = [Jan1, Mar1] \cup [Feb1, Apr1]$, and it is also present $[May1, Jun1]$. However, as can be seen by inspection 
	\begin{center}
		\includegraphics*[scale=0.25]{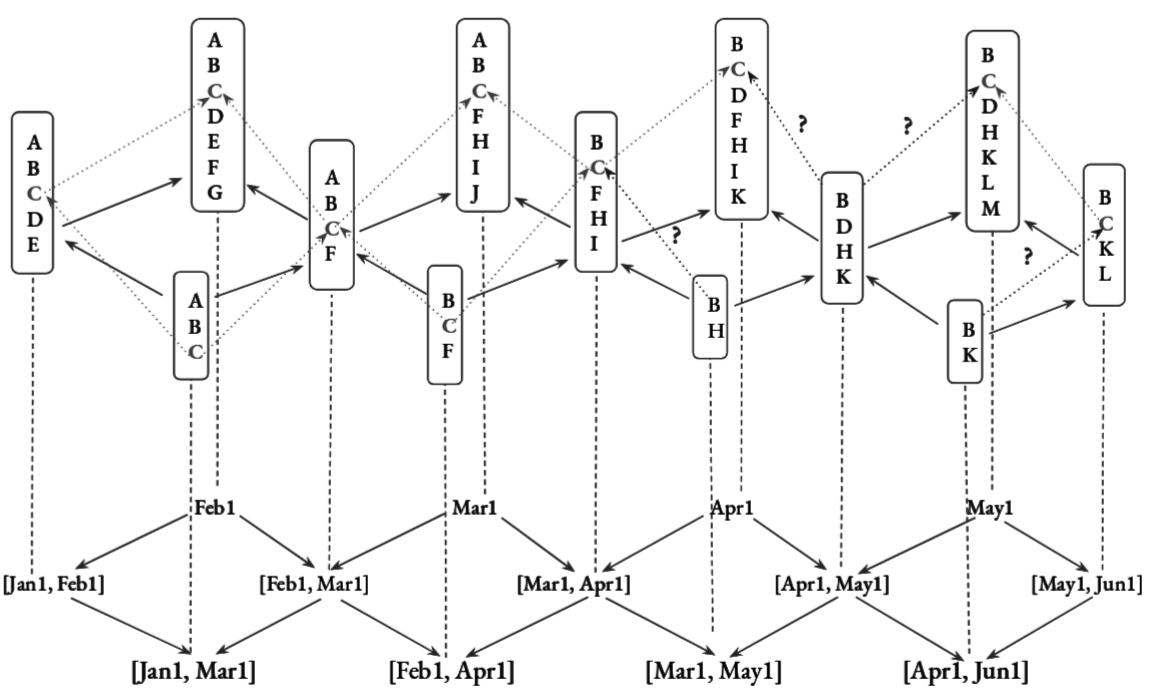}
	\end{center}
	the persistence of the product $C$ (witnessed in particular by the associated arrows to itself) throughout all of $[Jan1, Apr1]$ is ``proven" by the existence of local maps between all the sub-intervals covering this region. However, the (non-)maps indicated with question-marks have that question mark because there is not in fact any selection from the set of products given over, e.g., $S([Mar1, May1] = \{B, H\}$, that could get mapped, under the prescribed action of $S$, to $\{C\}$ in $S([Mar1, Apr1])$. Since there is, however, a map from $S([Mar1, Apr1])$ to $S([Apr1])$ that lands in $C$, the non-existence of the previous map tells us, in particular, that at some point in the period from $(Apr1, May1]$, the product $C$ ceased to be present in the company's stockpile of products. Additionally, as there is a map from $S([May1, Jun1])$ to $S([May1])$ sending $C$ to itself, we know that really the product $C$ could only have been removed from the stockpile in the period strictly between $Apr1$ and $May1$. \par 
	The point is that there is no way of gluing together the combined local sections of $C$ ``on the left" (from $Jan1$ through $Apr1$) to the local section ``on the right" (from $May1$ to $Jun1$). This gap in the presence of $C$ in the stockpile at some point in the time period between $Apr1$ and $May1$ is witnessed by the non-existence of any maps at the presheaf level involving $C$, that would have let us pass from ``one side" to the other. This makes the local section corresponding to the selection of $C$ with its associated maps \textit{strictly local}, since they cannot assemble into a global section, i.e., be glued together into a section specified over \textit{all of} the diagram and covering the time period from January 1 until June 1.\par 
	Notice also how, among the strictly local sections, some sections can be ``more local" than others, in the sense indicated by how, for instance, the product $D$ is only present throughout $[Jan1, Feb1]$ and $[Apr1,May1]$, as witnessed by the following restriction maps:   
	\begin{center}
		\includegraphics*[scale=0.25]{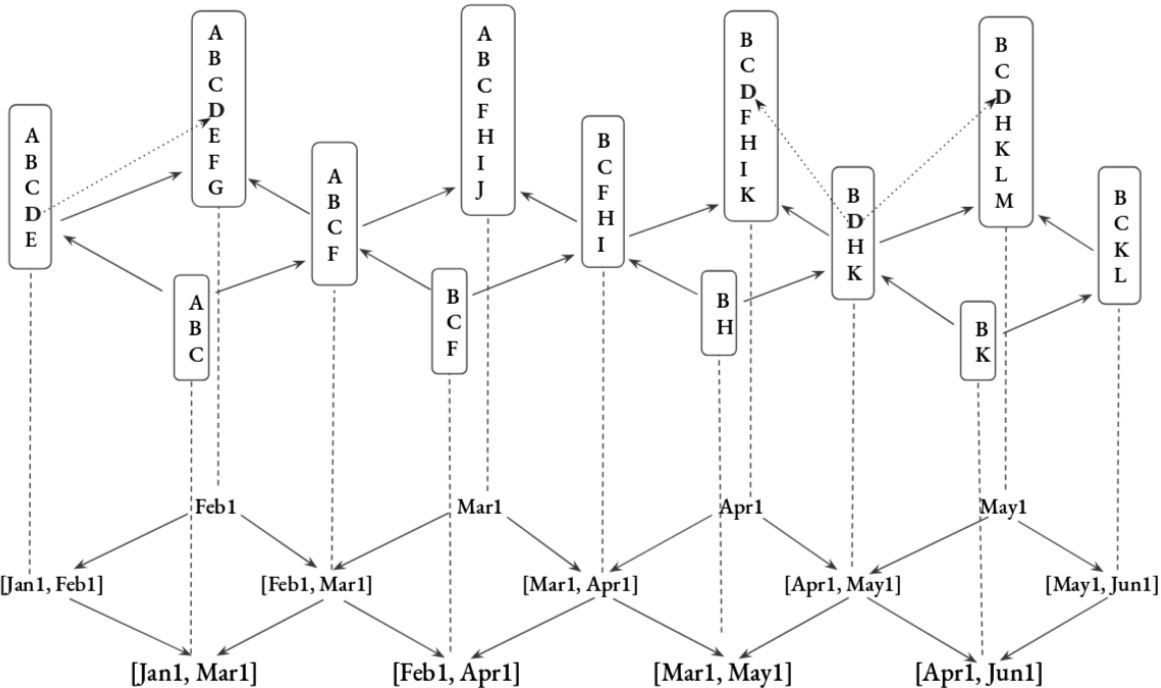}
	\end{center}     
	And we are informed of its moments of ``disappearance" precisely by the non-existence of ($D$-valued) maps into $S([Jan1, Feb1])$, $S([Feb1])$, $S([Apr1, May1])$, $S([Apr1])$, and $S([May1])$. \par 
	One final nuance is worth noting before moving on: observe how, for instance, the product $I$ shows up in the set $S([Mar1]) = \{A,B,C,F,H,I,J\}$, indicating that it is present in the company's stockpile throughout the ``instant" $[Mar1, Mar1]$. One can also see that it is in fact present throughout all of $[Mar1, Apr1]$. In general, viewing such sheaves ``internally," in terms of the contents of the variable sets assigned to each piece of the covers, such sets can be thought of as getting described by their behavior on small intervals. On any small interval, the products present in the stockpile over that period will either include $I$ or not include it. However, if we imagine decomposing the overall region more ``finely" and imagine looking at very small intervals containing $[Mar1]$ (take them however small you like), they may contain a sub-interval on the left over which $I$ is not present and a sub-interval on the right over which it is present. Thus, even over such small intervals containing $[Mar1]$, it is not correct to say that the product $I$ is either present or not present. The point is: over an arbitrary interval containing March 1 at 0:00, it is incorrect to say that the product $I$ is present and it is incorrect to say that the product $I$ is not present.\footnote{This might seem like an artificial feature of the nature of an ``instant," but it in fact suggests a pivotal feature of sheaves (to be discussed in more detail in later sections): that, in general, the logic of sheaves is \textit{intuitionistic} (not classical), where this means that the law of excluded middle $\alpha \vee \neg \alpha$ does not hold.}       
\end{example}
\begin{example}
	Recall the functor $nColor: \textbf{UCGraph}^{op} \rightarrow \textbf{Set}$\index{functor!nColoring}, first introduced in \ref{example: ncoloring}, that takes an undirected connected graph (recall that a graph is connected if there is a path between any two vertices in the graph) to the set of $n$-colorings of its vertices subject to the condition that no adjacent vertices are assigned the same color. In the case of undirected connected graphs, we can define a subgraph $G$ of a graph $H$ as a graph such that $\text{Edges}(G) \subseteq  \text{Edges}(H)$, and the further fact that $\text{Nodes}(G) \subseteq  \text{Nodes}(H)$ follows automatically since $G$ is assumed to be connected. Thus, in this context, to define a cover of a graph $G$ it suffices to specify a family of subgraphs $\{G_i \hookrightarrow G \hspace*{0.3em}| \hspace*{0.3em} i \in I\}$ satisfying the condition that 
	\begin{equation*}
	\bigcup_{i \in I} \text{Edges}(G_i) = \text{Edges}(G).
	\end{equation*} It can be shown that, using subgraph covers of a graph as above, on a given connected graph, we can in fact form a sheaf from the presheaf $nColor$. For concreteness, we exhibit this in the case of a $3$-coloring of the connected graph $K_3$ (with its subgraphs, ordered in the natural way). We first display what this functor assignment looks like over a particular subgraph of the connected graph $K_3$, then we display the full diagram conveying the sheaf over the space of subgraphs. The pictures are very explicit and take care of all the details; by attending to the pictures, the reader should be able to ``see" how there are actually two distinct $3$-coloring sheaves here, each got by selecting one of the two colorings (solutions) on all of $K_3$ and then restricting that down all the way through the inclusions.\footnote{The reader will note, however, that we do not represent all possible colorings, but only those colorings that have already fixed the coloring of the vertex 1 as blue. The rest, however, display colorings that are ultimately isomorphic to these two, since of the six 3-colorings of $K_3$, there are only two non-isomorphic ones.}  
	\par
	\begin{tikzpicture}[scale=.80,every node/.style={minimum size=1cm},on grid]
	\pagestyle{empty}
	\tikzstyle{edge_style} = [draw=black, line width=2, ultra thick]
	\tikzstyle{node_style1} = [circle,draw=blue,fill=blue!20!,font=\sffamily\small\bfseries, scale = 0.8]
	\tikzstyle{node_style2} = [circle,draw=red,fill=red!20!,font=\sffamily\small\bfseries, scale = 0.8]
	\tikzstyle{node_style3} = [circle,draw=green,fill=green!20!,font=\sffamily\small\bfseries, scale =0.8]
	\tikzstyle{node_style} = [circle,draw=black,fill=white!20!,font=\sffamily\small\bfseries, scale = 0.7]
	\begin{scope}[
	yshift=-83,every node/.append style={
		yslant=0.5,xslant=-1},yslant=0.5,xslant=-1
	]
	\fill[white,fill opacity=0.9] (-1,-1) rectangle (4,4);
	
	\draw[black,very thick] (-1,-1) rectangle (3,3);
	
	\end{scope}
	
	\begin{scope}[
	yshift=0,every node/.append style={
		yslant=0.5,xslant=-1},yslant=0.5,xslant=-1
	]
	\fill[white,fill opacity=.9] (-1,-1) rectangle (3,3);
	\draw[black,very thick] (-1,-1) rectangle (3,3);
	\draw[step=5mm, black] (-1,-1) rectangle (3,3);
	\end{scope}
	
	\begin{scope}[
	yshift=90,every node/.append style={
		yslant=0.5,xslant=-1},yslant=0.5,xslant=-1
	]
	\fill[white,fill opacity=.9] (-1,-1) rectangle (3,3);
	
	\draw[black, very thick] (-1,-1) rectangle (3,3);
	\end{scope}
	
	\begin{scope}[
	yshift=170,every node/.append style={
		yslant=0.5,xslant=-1},yslant=0.5,xslant=-1
	]
	\fill[white,fill opacity=0.9] (-1,-1) rectangle (3,3);
	
	\draw[black, very thick] (-1,-1) rectangle (3,3);
	\end{scope}
	
	\begin{scope}[
	yshift=-240,every node/.append style={
		yslant=0.5,xslant=-1},yslant=0.5,xslant=-1
	]      
	
	\draw[black] (0,0) rectangle (4,4); 	
	
	\end{scope}
	
	\draw[-latex,thick] (6.2,2) node[right]{$\mathsf{c_2}$}
	to[out=180,in=90] (4,2);
	
	\draw[-latex,thick](5.8,-.3)node[right]{$\mathsf{c_1}$}
	to[out=180,in=90] (3.9,-1);
	
	\draw[-latex,thick](5.9,5)node[right]{$\mathsf{c_3}$}
	to[out=180,in=90] (3.6,5);
	
	\draw[-latex,thick](5.9,8.4)node[right]{$\mathsf{c_4}$}
	to[out=180,in=90] (3.2,8);
	
	\draw[-latex,thick] (-4,7) node[right]{}
	to[out=180,in=90] (-4,-2);
	\node (t1) at (-6,3) {``Stalk"};
	\node (t2) at (-6,2) {over};
	\node (t3) at (-6,1) {$G$};
	\node (t4) at (4.5,9.5) {\text{nColor}($G$)};
	\draw[-latex,thick](4,-6)node[right]{$G \subseteq K_3$}
	to[out=180,in=90] (2,-5);	
	
	\node[node_style] (w) at (-2,-6.2) {$2$};
	\node[node_style] (d) at (1.2,-5.7) {$1$};
	\node[node_style] (e) at (0.2,-7.5) {$3$};
	\draw (w) -- (d) -- (e);	
	
	\node[node_style2] (w2) at (-2,-1.8) {$2$};
	\node[node_style1] (d2) at (1.2,-1.3) {$1$};
	\node[node_style2] (e2) at (0.2,-3.1) {$3$};
	\draw (w2) -- (d2) -- (e2);	
	
	\node[node_style3] (w3) at (-2,1) {$2$};
	\node[node_style1] (d3) at (1.2,1.7) {$1$};
	\node[node_style3] (e3) at (0.2,0) {$3$};
	\draw (w3) -- (d3) -- (e3);	
	
	\node[node_style2] (w2) at (-2,4.4) {$2$};
	\node[node_style1] (d2) at (1.2,4.9) {$1$};
	\node[node_style3] (e2) at (0.2,3.1) {$3$};
	\draw (w2) -- (d2) -- (e2);	
	
	\node[node_style3] (w2) at (-2,7.2) {$2$};
	\node[node_style1] (d2) at (1.2,7.7) {$1$};
	\node[node_style2] (e2) at (0.2,5.9) {$3$};
	\draw (w2) -- (d2) -- (e2);	
	\end{tikzpicture}
	
	\newcommand*\object[1][]{%
		\draw [blue, ultra thick, -{Stealth[bend]}, #1] (0,0) .. controls (1,-1) and (2,1) .. (3,0);
	}
	\begin{figure} 
		\centering 
		\begin{tikzpicture}[scale=.30,every node/.style={minimum size=1cm},on grid, rotate=32]
		\tikzstyle{edge_style} = [draw=black, line width=2, ultra thick]
		\tikzstyle{node_style1} = [circle,draw=blue,fill=blue!20!,font=\sffamily\small\bfseries, scale = 0.23]
		\tikzstyle{node_style2} = [circle,draw=red,fill=red!20!,font=\sffamily\small\bfseries, scale = 0.23]
		\tikzstyle{node_style3} = [circle,draw=green,fill=green!20!,font=\sffamily\small\bfseries, scale =0.23]
		\tikzstyle{node_style} = [circle,draw=black,fill=white!20!,font=\sffamily\small\bfseries, scale = 0.23]
		\begin{scope} 
		[xshift=-15cm]
		\begin{scope}
		[xshift=-7cm,yshift=-6cm]
		
		\begin{scope} 
		[xshift=-2cm, yshift=-1.2cm]
		\begin{scope}[
		xshift=650, yshift=-450,every node/.append style={
			yslant=0.5,xslant=-1},yslant=0.5,xslant=-1
		]
		\fill[white,fill opacity=0.9] (-1,-1) rectangle (4,4);
		
		\draw[black,very thick] (-1,-1) rectangle (3,3);
		
		\end{scope}
		
		\begin{scope}[
		xshift=650,yshift=-390,every node/.append style={
			yslant=0.5,xslant=-1},yslant=0.5,xslant=-1
		]
		\fill[white,fill opacity=.9] (-1,-1) rectangle (3,3);
		\draw[black,very thick] (-1,-1) rectangle (3,3);
		\draw[step=5mm, black] (-1,-1) rectangle (3,3);
		\end{scope}
		
		\begin{scope}[
		xshift=650,yshift=-330,every node/.append style={
			yslant=0.5,xslant=-1},yslant=0.5,xslant=-1
		]
		\fill[white,fill opacity=.9] (-1,-1) rectangle (3,3);
		
		\draw[black, very thick] (-1,-1) rectangle (3,3);
		\end{scope}
		
		\begin{scope}[
		xshift=650,yshift=-270,every node/.append style={
			yslant=0.5,xslant=-1},yslant=0.5,xslant=-1
		]
		\fill[white,fill opacity=0.9] (-1,-1) rectangle (3,3);
		
		\draw[black, very thick] (-1,-1) rectangle (3,3);
		\end{scope}
		
		\begin{scope}[
		xshift=650,yshift=-600,every node/.append style={
			yslant=0.5,xslant=-1},yslant=0.5,xslant=-1
		]      
		\draw[black,fill=gray!20] (0,0) rectangle (4,4);

		\end{scope} 
		\begin{scope}[
		xshift=650,yshift=-210,every node/.append style={
			yslant=0.5,xslant=-1},yslant=0.5,xslant=-1
		]
		\fill[white,fill opacity=0.9] (-1,-1) rectangle (3,3);
		
		\draw[black, very thick] (-1,-1) rectangle (3,3);
		\end{scope}
		\begin{scope}[
		xshift=650,yshift=-140,every node/.append style={
			yslant=0.5,xslant=-1},yslant=0.5,xslant=-1
		]
		\fill[white,fill opacity=0.9] (-1,-1) rectangle (3,3);
		
		\draw[black, very thick] (-1,-1) rectangle (3,3);
		\end{scope}
		
		\node[node_style] (w) at (21,-19.2) {$2$};
		\node[node_style] (d) at (24.2,-18.2) {$1$};
		\node[node_style] (e) at (22.8,-20.2) {$3$};
		\draw[thick] (w) -- (d);	
		
		\node[node_style2] (w2) at (21,-15) {$2$};
		\node[node_style1] (d2) at (24.5,-14.2) {$1$};
		\node[node_style2] (e2) at (22.8,-16) {$3$};
		\draw (w2) -- (d2);	
		
		\node[node_style3] (w3) at (20.5,-12.8) {$2$};
		\node[node_style1,opacity=0.5] (d3) at (23.8,-12.3) {$1$};
		\node[node_style3] (e3) at (22.8,-14) {$3$};
		\draw (w3) -- (d3);	
		
		\node[node_style2, opacity=0.5] (w2) at (20.6,-10.7) {$2$};
		\node[node_style1, opacity=0.5] (d2) at (24.5,-9.9) {$1$};
		\node[node_style3] (e2) at (22.8,-12) {$3$};
		\draw (w2) -- (d2);	
		
		\node[node_style3] (w2) at (20.5,-8.7) {$2$};
		\node[node_style1, opacity=0.5] (d2) at (24.5,-7.8) {$1$};
		\node[node_style2] (e2) at (22.8,-9.7) {$3$};
		\draw (w2) -- (d2);	
		
		\node[node_style3] (w2) at (20.5,-6.7) {$2$};
		\node[node_style1] (d2) at (24.5,-6) {$1$};
		\node[node_style1] (e2) at (22.8,-7.8) {$3$};
		\draw (w2) -- (d2);	
		
		\node[node_style2] (w2) at (20.5,-4) {$2$};
		\node[node_style1] (d2) at (24.5,-3.3) {$1$};
		\node[node_style1] (e2) at (22.8,-5.1) {$3$};
		\draw (w2) -- (d2);	
		
		\end{scope} 
		\begin{scope} 
		[xshift=0.1cm, yshift=-2.6cm]
		
		\begin{scope}[
		xshift=200, yshift=-210,every node/.append style={
			yslant=0.5,xslant=-1},yslant=0.5,xslant=-1
		]
		\fill[white,fill opacity=0.9] (-1,-1) rectangle (4,4);
		
		\draw[black,very thick] (-1,-1) rectangle (3,3);
		
		\end{scope}
		
		\begin{scope}[
		xshift=200,yshift=-150,every node/.append style={
			yslant=0.5,xslant=-1},yslant=0.5,xslant=-1
		]
		\fill[white,fill opacity=.9] (-1,-1) rectangle (3,3);
		\draw[black,very thick] (-1,-1) rectangle (3,3);
		\draw[step=5mm, black] (-1,-1) rectangle (3,3);
		\end{scope}
		
		\begin{scope}[
		xshift=200,yshift=-90,every node/.append style={
			yslant=0.5,xslant=-1},yslant=0.5,xslant=-1
		]
		\fill[white,fill opacity=.9] (-1,-1) rectangle (3,3);
		
		\draw[black, very thick] (-1,-1) rectangle (3,3);
		\end{scope}
		
		\begin{scope}[
		xshift=200,yshift=-30,every node/.append style={
			yslant=0.5,xslant=-1},yslant=0.5,xslant=-1
		]
		\fill[white,fill opacity=0.9] (-1,-1) rectangle (3,3);
		
		\draw[black, very thick] (-1,-1) rectangle (3,3);
		\end{scope}
		
		\begin{scope}[
		xshift=200,yshift=30,every node/.append style={
			yslant=0.5,xslant=-1},yslant=0.5,xslant=-1
		]
		\fill[white,fill opacity=0.9] (-1,-1) rectangle (3,3);
		
		\draw[black, very thick] (-1,-1) rectangle (3,3);
		\end{scope}
		
		\begin{scope}[
		xshift=200,yshift=90,every node/.append style={
			yslant=0.5,xslant=-1},yslant=0.5,xslant=-1
		]
		\fill[white,fill opacity=0.9] (-1,-1) rectangle (3,3);
		
		\draw[black, very thick] (-1,-1) rectangle (3,3);
		\end{scope}
		
		\begin{scope}[
		xshift=200,yshift=-360,every node/.append style={
			yslant=0.5,xslant=-1},yslant=0.5,xslant=-1
		]      
		\draw[black,fill=gray!20] (0,0) rectangle (4,4);   	
		
		\end{scope}
		\node[node_style] (w) at (4.5,-10.5) {$2$};
		\node[node_style] (d) at (8,-9.7) {$1$};
		\node[node_style] (e) at (6.8,-11.7) {$3$};
		\draw[thick] (d) -- (e);	
		
		\node[node_style2] (w2) at (4.5,-6.5) {$2$};
		\node[node_style1] (d2) at (8,-5.7) {$1$};
		\node[node_style2] (e2) at (6.8,-7.7) {$3$};
		\draw (d2) -- (e2);	
		
		\node[node_style3] (w3) at (4.5,-2.4) {$2$};
		\node[node_style1, opacity=0.5] (d3) at (8.3,-1.8) {$1$};
		\node[node_style3] (e3) at (6.8,-3.6) {$3$};
		\draw (d3) -- (e3);	
		
		\node[node_style2] (w2) at (4.5,0) {$2$};
		\node[node_style1, opacity=0.5] (d2) at (8.4,0.8) {$1$};
		\node[node_style3] (e2) at (6.8,-1.2) {$3$};
		\draw (d2) -- (e2);	
		
		\node[node_style1] (w2) at (4.5,1.9) {$2$};
		\node[node_style1, opacity=0.5] (d2) at (8.4,2.5) {$1$};
		\node[node_style3] (e2) at (6.8,0.7) {$3$};
		\draw (d2) -- (e2);	
		
		\node[node_style1] (w2) at (4.5,4.2) {$2$};
		\node[node_style1] (d2) at (8.4,5) {$1$};
		\node[node_style2] (e2) at (6.8,3) {$3$};
		\draw (d2) -- (e2);	
		
		\node[node_style3] (w2) at (4.5,-4.5) {$2$};
		\node[node_style1, opacity=0.5] (d2) at (8.4,-3.7) {$1$};
		\node[node_style2] (e2) at (6.8,-5.7) {$3$};
		\draw (d2) -- (e2);	
		\end{scope} 
		\begin{scope} 
		[xshift=-.5cm, yshift=-2cm]
		\begin{scope}[
		xshift=410, yshift=-330,every node/.append style={
			yslant=0.5,xslant=-1},yslant=0.5,xslant=-1
		]
		\fill[white,fill opacity=0.3] (-1,-1) rectangle (4,4);
		
		\draw[black,very thick] (-1,-1) rectangle (3,3);
		
		\end{scope}
		
		\begin{scope}[
		xshift=410,yshift=-270,every node/.append style={
			yslant=0.5,xslant=-1},yslant=0.5,xslant=-1
		]
		\fill[white,fill opacity=.9] (-1,-1) rectangle (3,3);
		\draw[black,very thick] (-1,-1) rectangle (3,3);
		\draw[step=5mm, black] (-1,-1) rectangle (3,3);
		\end{scope}
		
		\begin{scope}[
		xshift=410,yshift=-210,every node/.append style={
			yslant=0.5,xslant=-1},yslant=0.5,xslant=-1
		]
		\fill[white,fill opacity=.9] (-1,-1) rectangle (3,3);
		
		\draw[black, very thick] (-1,-1) rectangle (3,3);
		\end{scope}
		
		\begin{scope}[
		xshift=410,yshift=-150,every node/.append style={
			yslant=0.5,xslant=-1},yslant=0.5,xslant=-1
		]
		\fill[white,fill opacity=0.9] (-1,-1) rectangle (3,3);
		
		\draw[black, very thick] (-1,-1) rectangle (3,3);
		\end{scope}
		
		\begin{scope}[
		xshift=410,yshift=-90,every node/.append style={
			yslant=0.5,xslant=-1},yslant=0.5,xslant=-1
		]
		\fill[white,fill opacity=0.9] (-1,-1) rectangle (3,3);
		
		\draw[black, very thick] (-1,-1) rectangle (3,3);
		\end{scope}
		
		\begin{scope}[
		xshift=410,yshift=-20,every node/.append style={
			yslant=0.5,xslant=-1},yslant=0.5,xslant=-1
		]
		\fill[white,fill opacity=0.9] (-1,-1) rectangle (3,3);
		
		\draw[black, very thick] (-1,-1) rectangle (3,3);
		\end{scope}
		
		\begin{scope}[
		xshift=410,yshift=-480,every node/.append style={
			yslant=0.5,xslant=-1},yslant=0.5,xslant=-1
		]      
		\draw[black,fill=gray!20] (0,0) rectangle (4,4);   	
		
		\end{scope}
		\node[node_style] (w) at (12,-14.7) {$2$};
		\node[node_style] (d) at (15.9,-14.2) {$1$};
		\node[node_style] (e) at (14.2,-15.8) {$3$};
		\draw[thick] (e) -- (w);	
		
		\node[node_style2] (w2) at (12,-8.7) {$2$};
		\node[node_style1] (d2) at (15.9,-8.3) {$1$};
		\node[node_style3] (e2) at (14.2,-9.8) {$3$};
		\draw (e2) -- (w2);	
		
		\node[node_style3] (w3) at (12,-6.5) {$2$};
		\node[node_style1] (d3) at (15.9,-6) {$1$};
		\node[node_style1] (e3) at (14.2,-7.7) {$3$};
		\draw (e3) -- (w3);	
		
		\node[node_style2] (w2) at (12,-4.4) {$2$};
		\node[node_style1] (d2) at (15.9,-3.9) {$1$};
		\node[node_style1] (e2) at (14.2,-5.5) {$3$};
		\draw (e2) -- (w2);	
		
		\node[node_style1] (w2) at (12,-2.4) {$2$};
		\node[node_style1] (d2) at (15.9,-1.9) {$1$};
		\node[node_style2] (e2) at (14.2,-3.5) {$3$};
		\draw (e2) -- (w2);	
		
		\node[node_style1] (w2) at (12,0) {$2$};
		\node[node_style1] (d2) at (15.8,1) {$1$};
		\node[node_style3] (e2) at (14.2,-1.1) {$3$};
		\draw (e2) -- (w2);	
		
		\node[node_style3] (w2) at (12,-10.6) {$2$};
		\node[node_style1, opacity=0.5] (d2) at (15.9,-10) {$1$};
		\node[node_style2] (e2) at (14.2,-11.7) {$3$};
		\draw (e2) -- (w2);	
		\end{scope}
		\end{scope}
		
		\begin{scope}
		[xshift=19cm, yshift=-11cm, rotate=0, anchor=south]
		\begin{scope} 
		[xshift=6cm, yshift = -4.8cm]
		\begin{scope}[
		xshift=-190, yshift=320,every node/.append style={
			yslant=0.5,xslant=-1},yslant=0.5,xslant=-1
		]
		\fill[white,fill opacity=0.9] (-1,-1) rectangle (4,4);
		
		\draw[black,very thick] (-1,-1) rectangle (3,3);
		
		\end{scope}
		
		\begin{scope}[
		xshift=-190,yshift=380,every node/.append style={
			yslant=0.5,xslant=-1},yslant=0.5,xslant=-1
		]
		\fill[white,fill opacity=.9] (-1,-1) rectangle (3,3);
		\draw[black,very thick] (-1,-1) rectangle (3,3);
		\draw[step=5mm, black] (-1,-1) rectangle (3,3);
		\end{scope}
		
		\begin{scope}[
		xshift=-190,yshift=450,every node/.append style={
			yslant=0.5,xslant=-1},yslant=0.5,xslant=-1
		]
		\fill[white,fill opacity=.9] (-1,-1) rectangle (3,3);
		
		\draw[black, very thick] (-1,-1) rectangle (3,3);
		\end{scope}
		
		\begin{scope}[
		xshift=-190,yshift=520,every node/.append style={
			yslant=0.5,xslant=-1},yslant=0.5,xslant=-1
		]
		\fill[white,fill opacity=0.9] (-1,-1) rectangle (3,3);
		
		\draw[black, very thick] (-1,-1) rectangle (3,3);
		\end{scope}
		
		\begin{scope}[
		xshift=-190,yshift=160,every node/.append style={
			yslant=0.5,xslant=-1},yslant=0.5,xslant=-1
		]      
		\draw[black,fill=gray!20] (0,0) rectangle (4,4);   	
		
		\end{scope}

		\begin{scope} 
		[xshift=-22cm, yshift=1cm]
		\node[node_style] (w) at (12.3,6.4) {$2$};
		\node[node_style] (d) at (16,7) {$1$};
		\node[node_style] (e) at (15.2,5.1) {$3$};
		\draw[thick] (w) -- (d) -- (e);	
		
		\node[node_style2] (w2) at (12.3,11.1) {$2$};
		\node[node_style1] (d2) at (16.3,11.5) {$1$};
		\node[node_style2] (e2) at (15.2,9.6) {$3$};
		\draw (w2) -- (d2) -- (e2);	
		
		\node[node_style3] (w3) at (12.8,13.2) {$2$};
		\node[node_style1] (d3) at (16.3,13.7) {$1$};
		\node[node_style3] (e3) at (15.2,12) {$3$};
		\draw (w3) -- (d3) -- (e3);	
		
		\node[node_style2] (w2) at (12.3,15.4) {$2$};
		\node[node_style1] (d2) at (16.3,15.9) {$1$};
		\node[node_style3] (e2) at (15.2,14.1) {$3$};
		\draw (w2) -- (d2) -- (e2);	
		
		\node[node_style3] (w2) at (12.3,18.2) {$2$};
		\node[node_style1] (d2) at (16.3,18.7) {$1$};
		\node[node_style2] (e2) at (15.2,16.9) {$3$};
		\draw (w2) -- (d2) -- (e2);	
		\end{scope}
		\end{scope}
		\begin{scope}
		[xshift=5cm, yshift=-12cm]
		\begin{scope}[
		xshift=30, yshift=430,every node/.append style={
			yslant=0.5,xslant=-1},yslant=0.5,xslant=-1
		]
		\fill[white,fill opacity=0.9] (-1,-1) rectangle (4,4);
		
		\draw[black,very thick] (-1,-1) rectangle (3,3);
		
		\end{scope}
		
		\begin{scope}[
		xshift=30,yshift=490,every node/.append style={
			yslant=0.5,xslant=-1},yslant=0.5,xslant=-1
		]
		\fill[white,fill opacity=.9] (-1,-1) rectangle (3,3);
		\draw[black,very thick] (-1,-1) rectangle (3,3);
		\draw[step=5mm, black] (-1,-1) rectangle (3,3);
		\end{scope}
		
		\begin{scope}[
		xshift=30,yshift=550,every node/.append style={
			yslant=0.5,xslant=-1},yslant=0.5,xslant=-1
		]
		\fill[white,fill opacity=.9] (-1,-1) rectangle (3,3);
		
		\draw[black, very thick] (-1,-1) rectangle (3,3);
		\end{scope}
		
		\begin{scope}[
		xshift=30,yshift=620,every node/.append style={
			yslant=0.5,xslant=-1},yslant=0.5,xslant=-1
		]
		\fill[white,fill opacity=0.9] (-1,-1) rectangle (3,3);
		
		\draw[black, very thick] (-1,-1) rectangle (3,3);
		\end{scope}
		
		\begin{scope}[
		xshift=30,yshift=270,every node/.append style={
			yslant=0.5,xslant=-1},yslant=0.5,xslant=-1
		]      
		\draw[black,fill=gray!20] (0,0) rectangle (4,4);

		\end{scope}

		\begin{scope} 
		[xshift=-20cm, yshift=8cm]
		\node[node_style] (w) at (18.3,3.2) {$2$};
		\node[node_style] (d) at (22,4) {$1$};
		\node[node_style] (e) at (21,1.8) {$3$};
		\draw[thick] (d) -- (w) -- (e);	
		
		\node[node_style2] (w2) at (18.3,7.9) {$2$};
		\node[node_style1] (d2) at (22.3,8.2) {$1$};
		\node[node_style1] (e2) at (21.2,6.6) {$3$};
		\draw (d2) -- (w2) -- (e2);	
		
		\node[node_style3] (w3) at (18.3,9.9) {$2$};
		\node[node_style1] (d3) at (22.3,10.2) {$1$};
		\node[node_style1] (e3) at (21,8.6) {$3$};
		\draw (d3) -- (w3) -- (e3);	
		
		\node[node_style2] (w2) at (18.3,11.9) {$2$};
		\node[node_style1] (d2) at (22,12.5) {$1$};
		\node[node_style3] (e2) at (21,10.6) {$3$};
		\draw (d2) -- (w2) -- (e2);	
		
		\node[node_style3] (w2) at (18.3,14.5) {$2$};
		\node[node_style1] (d2) at (22,15.5) {$1$};
		\node[node_style2] (e2) at (21,13.2) {$3$};
		\draw (d2) -- (w2) -- (e2);	
		\end{scope}
		\end{scope}
		\begin{scope}
		[xshift=5.8cm, yshift=-19cm]
		\begin{scope}[
		xshift=210, yshift=520,every node/.append style={
			yslant=0.5,xslant=-1},yslant=0.5,xslant=-1
		]
		\fill[white,fill opacity=0.9] (-1,-1) rectangle (4,4);
		
		\draw[black,very thick] (-1,-1) rectangle (3,3);
		
		\end{scope}
		
		\begin{scope}[
		xshift=210,yshift=580,every node/.append style={
			yslant=0.5,xslant=-1},yslant=0.5,xslant=-1
		]
		\fill[white,fill opacity=.9] (-1,-1) rectangle (3,3);
		\draw[black,very thick] (-1,-1) rectangle (3,3);
		\draw[step=5mm, black] (-1,-1) rectangle (3,3);
		\end{scope}
		
		\begin{scope}[
		xshift=210,yshift=640,every node/.append style={
			yslant=0.5,xslant=-1},yslant=0.5,xslant=-1
		]
		\fill[white,fill opacity=.9] (-1,-1) rectangle (3,3);
		
		\draw[black, very thick] (-1,-1) rectangle (3,3);
		\end{scope}
		
		\begin{scope}[
		xshift=210,yshift=710,every node/.append style={
			yslant=0.5,xslant=-1},yslant=0.5,xslant=-1
		]
		\fill[white,fill opacity=0.9] (-1,-1) rectangle (3,3);
		
		\draw[black, very thick] (-1,-1) rectangle (3,3);
		\end{scope}
		
		\begin{scope}[
		xshift=210,yshift=370,every node/.append style={
			yslant=0.5,xslant=-1},yslant=0.5,xslant=-1
		]      
		\draw[black,fill=gray!20] (0,0) rectangle (4,4);

		\end{scope} 
		\begin{scope} 
		[xshift=-21cm, yshift=16cm]
		\node[node_style] (w) at (26,-1.3) {$2$};
		\node[node_style] (d) at (29.5,-0.5) {$1$};
		\node[node_style] (e) at (28.5,-2.6) {$3$};
		\draw[thick] (d) -- (e) -- (w);	
		
		\node[node_style1] (w2) at (26,3) {$2$};
		\node[node_style1] (d2) at (29.3,3.8) {$1$};
		\node[node_style2] (e2) at (28.5,1.7) {$3$};
		\draw (d2) -- (e2) -- (w2);	
		
		\node[node_style1] (w3) at (26,5) {$2$};
		\node[node_style1] (d3) at (29.3,5.8) {$1$};
		\node[node_style3] (e3) at (28.5,3.7) {$3$};
		\draw (d3) -- (e3) -- (w3);	
		
		\node[node_style2] (w2) at (26,7) {$2$};
		\node[node_style1] (d2) at (29.3,7.8) {$1$};
		\node[node_style3] (e2) at (28.5,5.8) {$3$};
		\draw (d2) -- (e2) -- (w2);	
		
		\node[node_style3] (w2) at (26,9.5) {$2$};
		\node[node_style1] (d2) at (29.3,10.3) {$1$};
		\node[node_style2] (e2) at (28.5,8.1) {$3$};
		\draw (d2) -- (e2) -- (w2);	
		\end{scope}
		\end{scope}
		\end{scope}
	
		\begin{scope}
		[xshift=-7.5cm, yshift=-11cm]
		\begin{scope}[
		xshift=-100, yshift=-470,every node/.append style={
			yslant=0.5,xslant=-1},yslant=0.5,xslant=-1
		]
		\fill[white,fill opacity=0.9] (-1,-1) rectangle (4,4);
		
		\draw[black,very thick] (-1,-1) rectangle (3,3);
		
		\end{scope}
		
		\begin{scope}[
		xshift=-100,yshift=-410,every node/.append style={
			yslant=0.5,xslant=-1},yslant=0.5,xslant=-1
		]
		\fill[white,fill opacity=.9] (-1,-1) rectangle (3,3);
		\draw[black,very thick] (-1,-1) rectangle (3,3);
		\draw[step=5mm, black] (-1,-1) rectangle (3,3);
		\end{scope}
		
		\begin{scope}[
		xshift=-100,yshift=-350,every node/.append style={
			yslant=0.5,xslant=-1},yslant=0.5,xslant=-1
		]
		\fill[white,fill opacity=.9] (-1,-1) rectangle (3,3);
		
		\draw[black, very thick] (-1,-1) rectangle (3,3);
		\end{scope}
		
		\begin{scope}[
		xshift=-100,yshift=-290,every node/.append style={
			yslant=0.5,xslant=-1},yslant=0.5,xslant=-1
		]
		\fill[white,fill opacity=0.9] (-1,-1) rectangle (3,3);
		
		\draw[black, very thick] (-1,-1) rectangle (3,3);
		\end{scope}
		
		\begin{scope}[
		xshift=-100,yshift=-230,every node/.append style={
			yslant=0.5,xslant=-1},yslant=0.5,xslant=-1
		]
		\fill[white,fill opacity=0.9] (-1,-1) rectangle (3,3);
		
		\draw[black, very thick] (-1,-1) rectangle (3,3);
		\end{scope}
		
		\begin{scope}[
		xshift=-100,yshift=-170,every node/.append style={
			yslant=0.5,xslant=-1},yslant=0.5,xslant=-1
		]
		\fill[white,fill opacity=0.9] (-1,-1) rectangle (3,3);
		
		\draw[black, very thick] (-1,-1) rectangle (3,3);
		\end{scope}
		
		\begin{scope}[
		xshift=-100,yshift=-110,every node/.append style={
			yslant=0.5,xslant=-1},yslant=0.5,xslant=-1
		]
		\fill[white,fill opacity=0.9] (-1,-1) rectangle (3,3);
		
		\draw[black, very thick] (-1,-1) rectangle (3,3);
		\end{scope}
		
		\begin{scope}[
		xshift=-100,yshift=-50,every node/.append style={
			yslant=0.5,xslant=-1},yslant=0.5,xslant=-1
		]
		\fill[white,fill opacity=0.9] (-1,-1) rectangle (3,3);
		
		\draw[black, very thick] (-1,-1) rectangle (3,3);
		\end{scope}
		
		\begin{scope}[
		xshift=-100,yshift=10,every node/.append style={
			yslant=0.5,xslant=-1},yslant=0.5,xslant=-1
		]
		\fill[white,fill opacity=0.9] (-1,-1) rectangle (3,3);
		
		\draw[black, very thick] (-1,-1) rectangle (3,3);
		\end{scope}
		\begin{scope}[
		xshift=-100,yshift=-630,every node/.append style={
			yslant=0.5,xslant=-1},yslant=0.5,xslant=-1
		]      
		\draw[black,fill=gray!20] (0,0) rectangle (4,4);

		\end{scope}
		\node[node_style] (w) at (-5.6,-20) {$2$};
		\node[node_style] (d) at (-2,-19.6) {$1$};
		\node[node_style] (e) at (-3.5,-21.3) {$3$};

		\node[node_style2] (w2) at (-5.6,-15.5) {$2$};
		\node[node_style1, opacity=0.6] (d2) at (-2,-15.1) {$1$};
		\node[node_style2] (e2) at (-3.5,-16.7) {$3$};

		\node[node_style3] (w3) at (-5.6,-13.2) {$2$};
		\node[node_style1] (d3) at (-2,-13.3) {$1$};
		\node[node_style3] (e3) at (-3.5,-14.6) {$3$};

		\node[node_style2] (w2) at (-5.6,-11) {$2$};
		\node[node_style1] (d2) at (-1.9,-10.8) {$1$};
		\node[node_style3] (e2) at (-3.5,-12.2) {$3$};

		\node[node_style3] (w2) at (-5.6,-9) {$2$};
		\node[node_style1] (d2) at (-2,-8.8) {$1$};
		\node[node_style2] (e2) at (-3.5,-10.4) {$3$};
		
		\node[node_style3] (w2) at (-5.6,-7) {$2$};
		\node[node_style1] (d2) at (-2,-6.8) {$1$};
		\node[node_style2] (e2) at (-3.5,-8.4) {$3$};
		
		\node[node_style3] (w2) at (-5.6,-5) {$2$};
		\node[node_style1] (d2) at (-2,-4.8) {$1$};
		\node[node_style2] (e2) at (-3.5,-6.2) {$3$};
		
		\node[node_style3] (w2) at (-5.6,-3) {$2$};
		\node[node_style1] (d2) at (-2,-2.8) {$1$};
		\node[node_style2] (e2) at (-3.5,-4.2) {$3$};
		
		\node[node_style3] (w2) at (-5.6,-1) {$2$};
		\node[node_style1] (d2) at (-2,-0.6) {$1$};
		\node[node_style2] (e2) at (-3.5,-2.1) {$3$};
		
		\node[node_style3] (w2) at (-5.5,1.2) {$2$};
		\node[node_style1] (d2) at (-2,1.6) {$1$};
		\node[node_style2] (e2) at (-3.5,0.1) {$3$};
		\end{scope}
		
		\begin{scope} 
		[xshift=52cm, yshift=-20cm]

		\begin{scope}[
		xshift=-460,yshift=750,every node/.append style={
			yslant=0.5,xslant=-1},yslant=0.5,xslant=-1
		]
		\fill[white,fill opacity=.9] (-1,-1) rectangle (3,3);
		
		\draw[black, very thick] (-1,-1) rectangle (3,3);
		\end{scope}
		
		\begin{scope}[
		xshift=-460,yshift=830,every node/.append style={
			yslant=0.5,xslant=-1},yslant=0.5,xslant=-1
		]
		\fill[white,fill opacity=0.9] (-1,-1) rectangle (3,3);
		
		\draw[black, very thick] (-1,-1) rectangle (3,3);
		\end{scope}
		
		\begin{scope}[
		xshift=-460,yshift=530,every node/.append style={
			yslant=0.5,xslant=-1},yslant=0.5,xslant=-1
		]      
		\draw[black,fill=gray!20] (0,0) rectangle (4,4);
		
		\end{scope}
		\node[node_style] (w3) at (-19,20.7) {$2$};
		\node[node_style] (d3) at (-15,21) {$1$};
		\node[node_style] (e3) at (-16.5,19.4) {$3$};
		\draw[thick] (d3) -- (e3) -- (w3) -- (d3);	
		
		\node[node_style2] (w3) at (-19,27.5) {$2$};
		\node[node_style1] (d3) at (-15,27.8) {$1$};
		\node[node_style3] (e3) at (-16.5,26.2) {$3$};
		\draw[thick] (d3) -- (e3) -- (w3) -- (d3);	 
		
		\node[node_style3] (w3) at (-19,30.2) {$2$};
		\node[node_style1] (d3) at (-14.8,30.9) {$1$};
		\node[node_style2] (e3) at (-16.5,28.9) {$3$};
		\draw[thick] (d3) -- (e3) -- (w3) -- (d3);	
		
		\begin{scope}
		[xshift=4cm, yshift=-2cm]
		\draw[->, dashed, gray, thick] (-62,-8) node[right, xshift=1.2cm, yshift=1.3cm] {inclusion} to (-52,-3.5);
		
		\draw[->, dashed, gray, opacity=0.5, thick] (-46,0.2) node[right, xshift=2cm, yshift=2.5cm] {inclusion} to (-31.5,7.5);
		\draw[->, bend right, dashed, gray, thick] (-21.5,8.7) node[right, xshift=-0.2cm, yshift=1.8cm] {inclusion} to (-19.5,19.8);
		
		\draw[->, very thick] (-60.5,15.5) node[left, xshift=-0.4cm, yshift =-0.5cm] {Restriction} to (-64.5,13.5);
		
		\draw[->, very thick] (-42,24) node[left, xshift=-0.7cm, yshift=-1cm] {Restriction} to (-52,19);
		\draw[->, very thick] (-25.5,31)  node[left, xshift=-0.7cm, yshift=-1cm] {Restriction} to (-34.5,26.7);
		\end{scope} 
		\end{scope}  
		\end{scope}
		\end{tikzpicture}
		\caption*{N-Coloring Sheaf for $n=3$} \label{fig:M1}
	\end{figure}
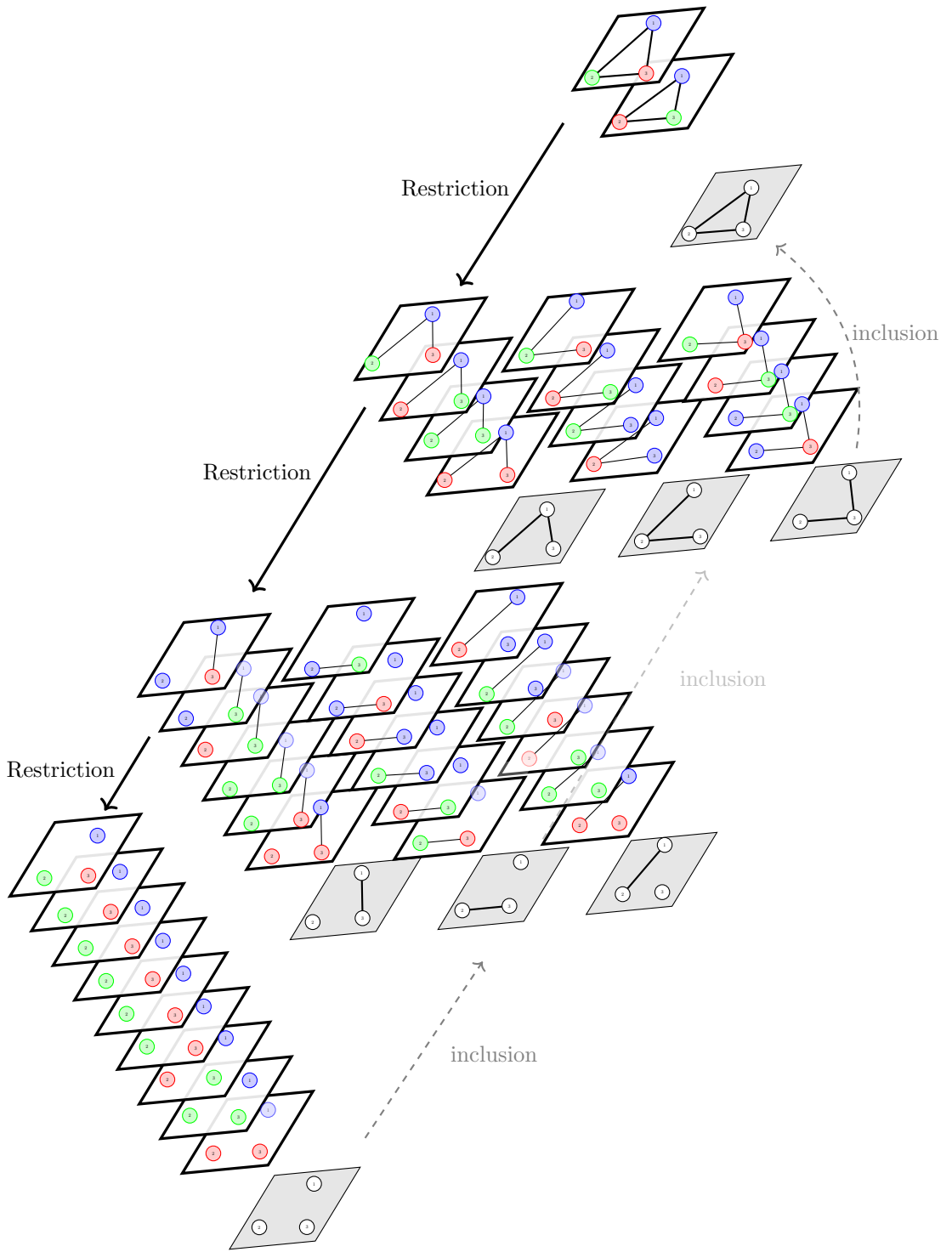  
\end{example}
\begin{example} 	
	For the next example, we consider a satellite, or various satellites, making passes over portions of a region of earth, collecting data as it goes. For concreteness, consider some specific portion of the earth, say Alaska, or that part of Alaska where the Bering Glacier lies, as a topological space $X$. Then given an open subset $U \subseteq  X$, we can let $S(U)$ denote the set of functions from $U$ to $C$, where $C$ might be the set interval of wavelengths in the light spectrum, or some geo-referenced (perhaps timestamped) intensity-valued image data, or some other data corresponding to the data feed of the satellites (or the processing thereof). This presheaf $S$ is in fact a sheaf, since we can indeed fuse together the different data given over the open sets of $X$, forming a larger patched-together image of the glacier. For concreteness, assume we are given the following selection of three satellite images of the Bering Glacier, chosen from among the (possibly very large) sets of images assigned to each region:\footnote{The images come from Landsat 8, here: {https://earthobservatory.nasa.gov/IOTD/view.php?id=4710}.} 
	\begin{center}
		\includegraphics[scale=0.27]{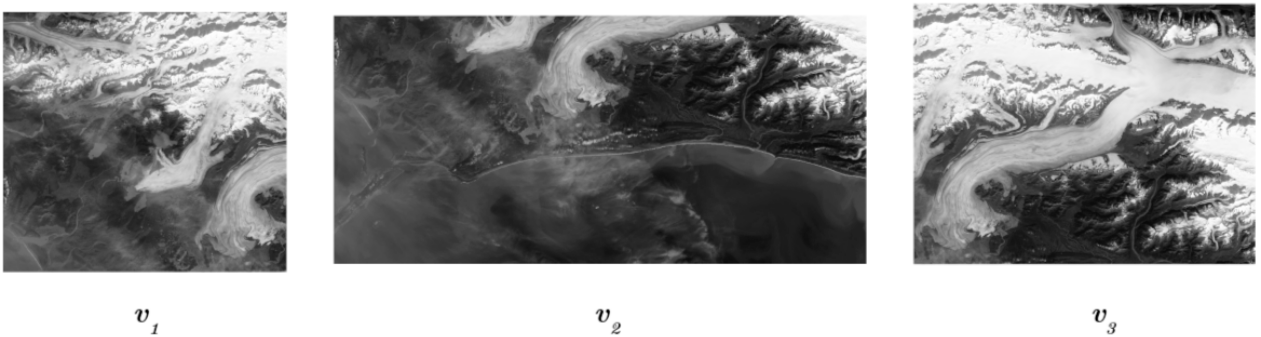}  
	\end{center}
	Each of the $v_i \in S(U_i)$ correspond to value assignments throughout or over certain subsets, $U_1, U_2, U_3$ of $X$, which together cover some subset $U \subseteq  X$---say the region of Alaska corresponding to the glacier. In terms of the data ``sitting over" each of these regions, as provided by each of the satellites in the form of individual images, we can notice that the restriction of $v_1$ to the region $U_1 \cap U_2$ is equal to the restriction of $v_2$ to the same subset $U_1 \cap U_2$, and so on, all the way down to their common restriction to $U_1 \cap U_2 \cap U_3$:
	\begin{center}
		\includegraphics[scale=0.37]{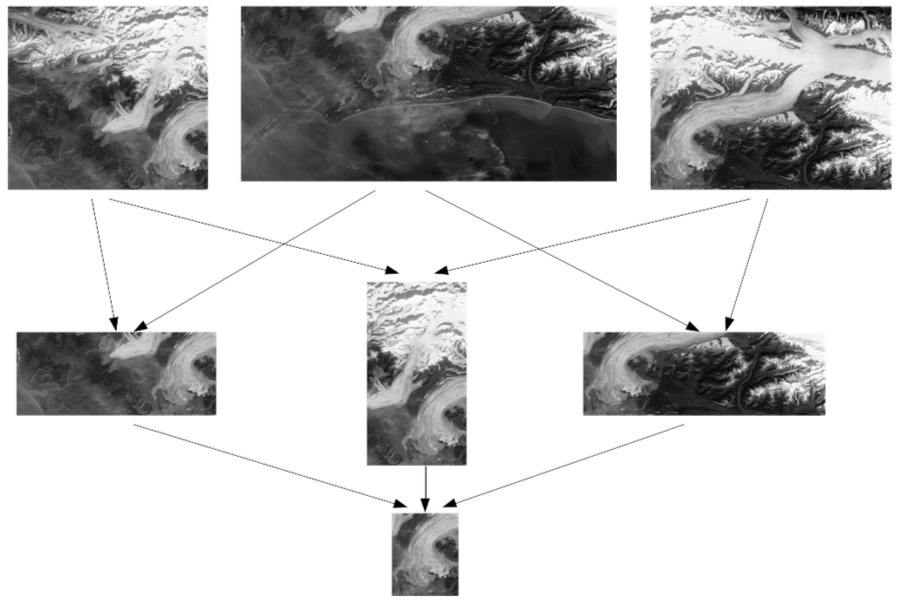}
	\end{center} 
	One can thus immediately \textit{see} that the sheaf condition is met, which means that we can in fact patch together the given local pieces or sections over the members of the open covering of $U$ to obtain a section over all of $U = U_1 \cup U_2 \cup U_3$. In summary, we have the following inclusion diagram (on the left) describing the underlying topology, paired with the sheaf diagram (on the right) with its corresponding restriction maps (notice the change in direction): 
	\begin{center}
		\begin{tikzpicture}[yscale=0.6, xscale=0.8]
		\node (max) at (0,4) {$U$};
		\node (a) at (-2.4,2) {$U_1 \cup U_2$};
		\node (b) at (0,2) {$U_1 \cup U_3$};
		\node (c) at (2.4,2) {$U_2 \cup U_3$};
		\node (d) at (-2,0) {$U_1$};
		\node (e) at (0,0) {$U_2$};
		\node (f) at (2,0) {$U_3$};
		\node (g) at (-2.4,-2) {$U_1 \cap U_2$};
		\node (h) at (0,-2) {$U_1 \cap U_3$};
		\node (i) at (2.4,-2) {$U_2 \cap U_3$};
		\node (min) at (0,-4) {$U_1 \cap U_2 \cap U_3$};
		\draw[->] (a) -- (max);
		\draw[->] (b) -- (max);
		\draw[->] (c) -- (max);
		\draw[->] (d) -- (a);
		\draw[->] (d) -- (b);
		\draw[->] (e) -- (a);
		\draw[->] (e) -- (c);
		\draw[->] (f) -- (b);
		\draw[->] (f) -- (c);
		\draw[->] (min) -- (g);
		\draw[->] (min) -- (h);
		\draw[->] (min) -- (i);
		\draw[->] (g) -- (d);
		\draw[->] (g) -- (e);
		\draw[->] (h) -- (d);
		\draw[->] (i) -- (e);
		\draw[->] (h) -- (f);
		\draw[->] (i) -- (f);
		
		\node (max1) at (8,4) {$S(U)$};
		\node (a1) at (5.5,2) {$S(U_1 \cup U_2)$};
		\node (b1) at (8,2) {$S(U_1 \cup U_3)$};
		\node (c1) at (10.5,2) {$S(U_2 \cup U_3)$};
		\node (d1) at (6,0) {$S(U_1)$};
		\node (e1) at (8,0) {$S(U_2)$};
		\node (f1) at (10,0) {$S(U_3)$};
		\node (g1) at (5.5,-2) {$S(U_1 \cap U_2)$};
		\node (h1) at (8,-2) {$S(U_1 \cap U_3)$};
		\node (i1) at (10.5,-2) {$S(U_2 \cap U_3)$};
		\node (min1) at (8,-4) {$S(U_1 \cap U_2 \cap U_3)$};
		\draw[<-] (a1) -- (max1);
		\draw[<-] (b1) -- (max1);
		\draw[<-] (c1) -- (max1);
		\draw[<-] (d1) -- (a1);
		\draw[<-] (d1) -- (b1);
		\draw[<-] (e1) -- (a1);
		\draw[<-] (e1) -- (c1);
		\draw[<-] (f1) -- (b1);
		\draw[<-] (f1) -- (c1);
		\draw[<-] (min1) -- (g1);
		\draw[<-] (min1) -- (h1);
		\draw[<-] (min1) -- (i1);
		\draw[<-] (g1) -- (d1);
		\draw[<-] (g1) -- (e1);
		\draw[<-] (h1) -- (d1);
		\draw[<-] (i1) -- (e1);
		\draw[<-] (h1) -- (f1);
		\draw[<-] (i1) -- (f1);
		\node (dom) at (0,-6) {$\mathscr{O}(X)^{op}$};
		\node (codom) at (8,-6) {\textbf{Set}};
		\draw[->] (2,-6) to node[pos=0.5, above, font=\large] (6,-6) {$S$} (6,-6);
		\end{tikzpicture}
	\end{center} 
	In terms of the actual images, the sheaf diagram on the right is pictured below, where we can think of the restriction maps as performing a sort of ``cropping" operation, corresponding to a reduction in the size of the domain of the sensor, while the gluing operation corresponds to patching or gluing the images together along their overlaps all the way up to the topmost image (which of course corresponds to the section or assignment over all of $U$).     
	\par 
	\begin{center}
		\includegraphics[scale=0.42]{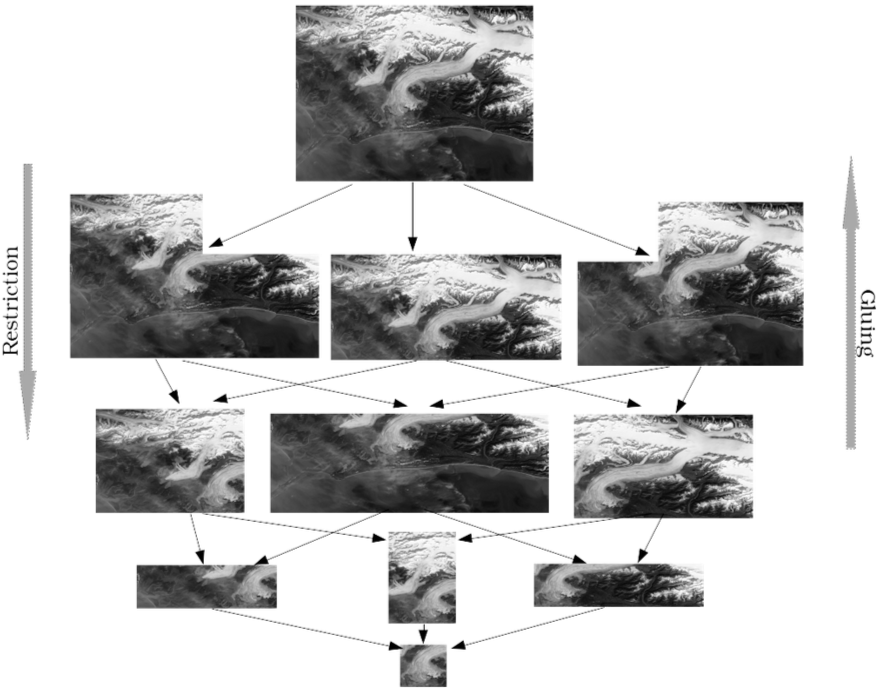}
	\end{center} 
\end{example} \noindent 
This mosaic\index{mosaic} example gives a particularly concrete illustration and motivation for an alternative definition of a sheaf, namely as a presheaf functor $F: \mathscr{O}(X)^{op} \rightarrow \textbf{Set}$ that moreover \textit{preserves limits}---where,\index{limit} because we use the opposite category for domain in defining the presheaf, this means that \textit{colimits} get\index{colimit} sent to limits in \textbf{Set}. We can see that in the lattice of open subsets of $X$, for an $I$-indexed family of open subsets $U_i \subseteq U$ (in the particular case described above, $I = 3$) that \textit{covers}\index{cover} $U$---in the exact sense that the entire diagram comprised of the sets $U_i$ and the inclusions of their pairwise intersections $U_i \cap U_j$ has $U$ for its colimit---the contravariant functor $S$ given above \textit{preserves} this colimit in the sense that it sends it to a \textit{limit} in \textbf{Set}. In terms of the universal characterization of these notions developed earlier, one can basically immediately see that while all arrows ``fall into" $U$ (think of $U$ as the nadir of a cone), any other possible object in this poset will have to pass through $U$, i.e., $U$ is initial; likewise, it is patently visible that the summit $S(U)$ will be terminal among cones. More formally, 
\begin{definition}
	(\textit{Yet Another Definition of a Sheaf on a Topological Space})\index{sheaf!defined} Given a presheaf $F: \mathscr{O}(X)^{op} \rightarrow \textbf{Set}$ from the poset of open sets of $X$ to $\textbf{Set}$, and defining an $I$-indexed family of open subsets $U_i \subseteq U$ as a \textit{cover} for $U$ when the entire diagram consisting of all the $U_i$ together with the inclusions of their pairwise intersections $U_i \cap U_j$---i.e., \begin{tikzcd}
		U_j \arrow[r, hookleftarrow] & U_i \cap U_j  \arrow[r, hookrightarrow] & U_i	
	\end{tikzcd}---has $U$ for its colimit, then such a presheaf $F$ is a \textit{sheaf} (of sets) provided it preserves these colimits, sending them to limits in \textbf{Set}. \par 
This means, in effect, that for any open cover $\{U_i\}_{i \in I}$ of $U$ (colimit), the following is an equalizer diagram \begin{center} 
		\begin{tikzcd}[column sep = 3.8em]
			F(U) \arrow[r, rightarrowtail, "{F(U_i \hookrightarrow U)}"] & \prod_{i \in I} F(U_i) \arrow[r, shift right = 1.2, "p"] \arrow[r, "q", swap] & \prod_{i,j \in I} F(U_i \cap U_j)
	\end{tikzcd} \end{center} 
	in \textbf{Set} (recall that an equalizer diagram is a \textit{limit} diagram).\index{equalizer} Here, for $t \in FU$, and letting the equalizer map $F(U_i \hookrightarrow U)$ be denoted by $e$, it is to be understood that 
	\begin{equation*}
	e(t) = \{t|_{U_i} \hspace*{0.3em}| \hspace*{0.3em} i \in I \}
	\end{equation*}
	and for a family $t_i \in FU_i$, we will have 
	\begin{equation*}
	p({t_i}) = \{t_i |_{(U_i \cap U_j)} \}, \hspace{2em} q({t_i}) = \{t_j |_{(U_i \cap U_j)} \},
	\end{equation*}
	the map $p$ involving $F(U_i \cap U_j \hookrightarrow U_i)$ composed with the appropriate projection map, while $q$ involves $F(U_i \cap U_j \hookrightarrow U_j)$ together with its projection map. 
\end{definition} 
An arrow \textit{into} a product is entirely determined by the components, namely its composition with the projections of the product. Thus, the maps $e$, $p$ and $q$ of the equalizer diagram above are precisely the unique maps making the ``unfolded" diagrams below commute for all $i, j \in I$ (where the vertical maps are the relevant projections of the products): 
\begin{center} 
	\begin{tikzcd}
		& FU_i \arrow[d, leftarrow] \arrow[rr, "{F(U_i \cap U_j \hookrightarrow U_i)}"] & & F(U_i \cap U_j) \arrow[d, leftarrow] \\
		FU \arrow[r, dashed, "e"] \arrow[ur] \arrow[dr] & {\prod_i FU_i} & & {\prod_{i,j} F(U_i \cap U_j)} \arrow[ll, shift right = 1.5, "p", dashed, yshift=-1ex, swap, leftarrow] \arrow[ll, "q", dashed, shift left = 0.3, yshift = -1ex, leftarrow] \\
		& FU_j \arrow[u, leftarrow] & & F(U_i \cap U_j) \arrow[ll, leftarrow, "{F(U_i \cap U_j \hookrightarrow U_j)}"] \arrow[u, leftarrow, swap]
	\end{tikzcd}
\end{center} 
The utility of this alternate description is that it furnishes us with a completely \textit{categorical} description of the equalizer diagram, which means that the above definition of a sheaf will work even when we replace \textbf{Set} with other suitable categories (specifically, those with all small products). In other words, we might just as well have provided a definition of sheaves $F: \mathscr{O}(X)^{op} \rightarrow \textbf{D}$ of \textbf{D}-objects on a space $X$; there are many prominent candidates for \textbf{D} in this more general definition, some of which we will meet in later sections and chapters, e.g., giving rise to sheaves of abelian groups, vector spaces, rings, $R$-modules. \par 
The next two examples are more ``for fun," meant to emphasize or reinforce certain aspects of the \textit{idea} of a sheaf.
\begin{example}
	The 20th-century pianist Glenn Gould\index{Glenn Gould} was one of the first to ardently defend the merits of studio recording and use of the \textit{tape splice}\index{splice} in the creative process, against those who held fast to the supposedly more ``moral" or ``pure" tradition of the live concert performance (and who accordingly thought that the only purpose of the splice would be to rectify performance mishaps or to alleviate the pressure of the ``one-take" approach demanded by the concert form).\footnote{By splices, one means an edit point representing the confluence of distinct takes or inserts (i.e., recorded performance of a portion of the score).} \par 
	Gould challenged the view that the only legitimate continuity of a unified interpretation could come from the one-takeness of traditional performance, proposing instead that the listener's ``splice prerogative" and the performer's newfound editorial control in the recording studio would bestow upon creators an even more demanding ethic concerning matters of architecture and integrity of vision. Gould claimed that new, explicit, and more demanding forms of continuity were to be found in this montage-based approach: ``splicing builds good lines, and it shouldn't much matter if one uses a splice every two seconds or none for an hour so long as the result \textit{appears} to be a coherent whole." Just as one does not demand or expect that the filmmaker shoot a film in one shot, Gould believed that one should not expect that the coherence or continuity of an interpretation of a musical piece can only be secured by the inexorable linearity of time and the single take---the musician has just as much a right to montage as the filmmaker.  \par 
	Gould went as far as to test, with a controlled experiment involving eighteen participants, whether listeners (including laymen and recording experts) could detect the ``in point" of any splice in certain selections of recordings, each of which selection had drastically different splice densities (in some cases, none).\footnote{The results can be found in his essay \cite{page_grass_1990}.} What he found, in short, was that ``the tape does lie and nearly always gets away with it." While originally (with analog magnetic tape splicing) the tape splice involved careful (and literal) cutting of the physical tape with scissors or a blade and (literal) gluing or taping of it to another section of tape (possibly from an entirely different recording session)---whenever qualities, such as tempo, of the two recordings to be joined could be made compatible enough to permit a seamless joining---Gould foresaw the power inherent in the more general \textit{idea} of splicing and montage: this new approach could provide a more analytically acute dissection of the minute connections ultimately defining the coherence of a particular piece of music, displaying the architectural coherence of it less dogmatically than one would have to in relying on the ``in-built continuity" allegedly belonging to the one-take concert ideal, focusing instead on breaking a score down into its smallest parts, recording many ``takes" of such sections of the score, and then gluing together the results of certain of those distinct performances, whenever they could be made compatible on their overlap, all with the aim of producing a single, unified performance of the entire score. Gould's insistence on the virtues of the tape-splice and on the importance of montage in recording practice nicely captures something\index{Glenn Gould} akin to the fundamental ``spirit" of the sheaf construction.\footnote{Incidentally, Gould's closeness to the ``sheaf philosophy" is evidenced in a number of aspects of his life, not just his approach to his art, for instance via his impressive insistence on forever integrating as many disparate planes and partial pieces of information as possible into a single coherent experience, e.g., his alleged habit of simultaneously listening to all of the conversations going on in a caf\'e. This mentality is perhaps most famously illustrated by the various accounts of him purposefully dividing his concentration across multiple channels in order to better understand something, as for instance when he claimed that he discovered he could best understand Schoenberg's Opus 23 if he listened to it while simultaneously playing the news on the radio, or when he mastered a demanding section of a Beethoven sonata only after placing a radio and television next to the piano and turning them up as loud as they would go as he worked through that passage. This embodies something like the ``sheaf philosophy": integration and coherence not through an enforced isolationism, but precisely through complete immersion in the dense texture that arises by decomposition into pieces, careful choices made locally, and the resulting appreciation of the need to make explicit the most minute of links between parts of a whole, as one gradually, piece by piece, assembles a more global or unified perspective. In this connection, we could also mention one of his descriptions of his famous ``contrapuntal radio" programs from the 70's, in which he claimed to try ``to have situations arise cogently from within the framework of the program in which two or three voices could be overlapped, in which they would be heard talking---simultaneously, but from different points of view---about the same subject."} \par 
	Moreover, while the single-take approach to recording and unifying the musical idea is essentially \textit{deductive} (and purports to be neutral in its simple ``transmission of the facts"), reducing the individual (voice, line, note) to its participation in a prefabricated idea of totality and relying on a dubious notion of some ``immediate" continuity, montage/splicing (like the sheaf construction) is fundamentally \textit{inductive}, allowing the individual component materials of a work to create their own formal structure ``from the bottom up" via insistence on the transparent and explicit unfolding of the principles by which the component parts can be patched together locally. According to Gould, it is precisely through the initial discontinuity induced by the decomposition into parts and cutting process in montage/splicing that the task of making explicit the principle of their reorganization/patching into a unified totality is allowed to emerge, and is no longer regarded as something \textit{a priori} or to be taken for granted. Just as in the sheaf construction, this approach essentially involves both cutting (decomposition/discontinuity) of the ``space" (the score) and local patching or gluing (recomposition/continuity) of data (distinct partial recordings), gradually building up to a unique data assignment over the entire space (recording of the entire score).\index{Glenn Gould} \par 
	For a more concrete, if very rough and simplistic, idea of how the splicing or montage approach to recording might be seen as akin to the construction of a sheaf, consider a score consisting of 32 measures. We might then consider that the ``space" of the score has been decomposed into three principal parts or pieces: (A) spanning from measure 1 to the end of measure 16; (B) spanning from the beginning of measure 8 until the end of measure 24; and (C) spanning from the beginning of measure 16 until the final measure. Together, these portions collectively cover the entire 32-measure score, and there are the obvious overlapping measures. We can now imagine that to each section (A)-(C), there corresponds a (possibly very large) set of distinct recordings. If, for some selection of individual recordings from each of the three regions (A)-(C), the select recordings can be made to agree on their overlap---via some system of translation functions, e.g., slowing down one recording to match the tempo of another---then they can be spliced together into a unique recording of the entire work.    
\end{example}
\begin{example}
	Detectives\index{detectives} collect certain information pertaining to a crime that purportedly occurred in some area during a certain time interval. This information will most likely be heterogeneous in nature, i.e., they may have camera footage of some part of the scene, some eyewitness testimony, some roughly time-stamped physical data, etc. These various pieces of data are all considered to be local in the sense that they are assumed to concern (or be valid throughout) a delimited region of space-time, e.g., time-stamped camera footage of one of the parking lot's exits or an eyewitness testimony claiming to have heard a scream coming from the southern end of the parking lot sometime between 8:00 pm and 8:30 pm. In terms of the underlying space-time regions to which these various pieces of information correspond, the various pieces of information may very well agree or be compatible on the overlaps, e.g., an eyewitness's testimony with respect to a particular half-hour interval and location might be checked against the camera feed concerning that same time interval and area. In general, the various pieces of data over the same interval may corroborate one another or contradict one another, either entirely or in some particular respect or with respect to some sub-region of their overlap. It is not always as simple as verifying whether or not they provide the \textit{same} information. It may happen, for instance, that the parking lot is constructed in such a way that certain barriers acoustically account for why the witness heard the scream coming from the southern end of the parking lot, when in fact it could only have come from the western end (which is where the camera shows the victim in conflict during that time). It is the job of the detective to find the appropriate ``translation functions" making sense of these at first (potentially) conflicting local pieces of data and then use these functions to ``glue" together, step by step, the data that can be made to locally cohere into a coherent and self-consistent account of what occurred over the entire spatio-temporal interval in question. In a rough sense, then, given a presheaf assigning information (camera data, propositions, etc.) locally over some collection of space-time regions, the detective is looking to build a sheaf over the entire space-time interval covered by all those regions.\par 
	This example also invites us to consider, in a very preliminary fashion, just one of the many interesting features of a sheaf: namely that in addition to the fact that sheaves let us determine global unknowns or solutions from data given merely locally, whenever we do indeed have a sheaf, say $F$, this may enable the prediction of missing or under-specified data (or at least the specification of what data will be \textit{possible}) with respect to some sub-region, given some selection of data over another region.\footnote{The presentation of this observation in the next paragraph closely follows \cite{spivak_category_2014}, 429-30.} For instance, given some time-interval and some particular area, say, $U = $(8:00, 9:00) $\subseteq  X = $(7:00, 9:40) concerning the parking lot in question, we might consider some subset $A$ of all the data we have over that interval, i.e., $F(U)$. Now consider another region $V = $(7:30, 8:30) $\subseteq  X = $(7:00, 9:40). Diagrammatically, we thus have the following: 
	\begin{center}  
		\begin{tikzcd}
			& F(X) \arrow[d, "{\rho_{U,X}}"] \arrow[r, "{\rho_{V,X}}"] & F(V) \\
			A \arrow[r, "i", swap] & F(U) 
		\end{tikzcd}
	\end{center}  
	We know that because we are dealing with \textit{sets}, we can form the pullback or fiber product, i.e., $A \times_{F(U)} F(X) := \{(a,u,x) \hspace*{0.3em}| \hspace*{0.3em} i(a) = u = \rho_{U,X}(x) \}$ (more on this construction in a later chapter):
	\begin{center}  
		\begin{tikzcd}
			(\rho_{U,X})^{-1}(A) \arrow[d] \arrow[r] \arrow[rr, dashrightarrow, bend left = 40] & F(X) \arrow[d, "{\rho_{U,X}}"] \arrow[r, "{\rho_{V,X}}"] & F(V) \\
			A \arrow[r] & F(U) 
		\end{tikzcd}
	\end{center}     
	Then the image of the top composite (dashed) will yield a subset of $F(V)$ which informs us about the possible value assignments throughout $V$, \textit{given} what we know to be the case (namely $A$) throughout $U$. Moreover, we could further consider maps into $A$ and continue forming pullbacks; since the left-most square will be a pullback iff the composite large rectangle forms a pullback, we can paste together pullbacks and further refine or constrain these predictions in a controlled way.
\end{example}
\subsection{Philosophical Pass: Sheaf as Local-Global Passage}
A sheaf\index{sheaf!as local-global passage} is not to be situated in either the local (restriction) or the global (collation) registers, but rather is to be located in the passage forged between these two, in the translation system or glue that mediates between the two registers. The transit from the local to the global secured via the sheaf gluing (collatability) condition provides a deep but also precisely controllable connection between continuity (via the emerging system of translation functions guaranteeing coherence or compatibility between the local sections) and generality (global sections). By separating something into parts, i.e., by specifying information locally, considering coverings of the relevant region, and enabling the decomposition or refinement of value assignments into assignments over restricted parts of the overall region (restriction condition), we are presented with a problem, a problem that in a sense can only first appear with such a ``downward" movement towards greater refinement. Without having separated something into parts, we may appear to have a sort of trivial or default cohesion of parts, where, without being recognized in their separation, the parts yet remain implicit and so the glue binding them together or the rule allowing one to transit from one part to another in a controlled fashion is simply not visible. However, having decomposed or discretized something into parts, we are at once presented with this separation of parts and the problem of finding and making explicit the glue that will serve to bind them together. A sheaf is a way of taking information that is locally defined or assigned and decomposing those assignments in a controlled fashion into assignments over smaller regions so as to draw out the specific manner of effecting translations or gluings that obtain between those particular assignments with respect to their overlapping regions, and then using this now explicit system of gluings to build up a unique and comprehensive value assignment over the entire network of regions. In this sense, a sheaf equally involves both (i) controlled decomposition (discreteness), and (ii) the recomposition (continuity) of what is partial into an architecture that makes explicit the special form of cooperation and harmony that exists between the decomposed items, items that may have previously been detached, or which may have only \textit{appeared} to ``stick together" because we had not bothered to look closely enough.\par 
Via the restriction/localization step, sheaf theory teaches us that we do not command a more global or integrated vision by renouncing the local nature of information or distinct planes and textures of reality or by glossing over the minute passages between things. Instead, it forces us to first become masters of the smallest link and, precisely through that control of the passages between the local parts, forge a coherent (``collatable") vision of the largest scope. \par 
Phenomenologically speaking, data or observations are frequently presented to us in ``zones," ``fragmented" or isolated in some way. These items can be thought of as various light-beams (perhaps of specific hues or brightness) cast over (and covering only parts of) a vast landscape, some of which may overlap. Even if this data clearly emerges as \textit{evolving} over some region, it remains indexed or determined in some way by a particular ``zone" or context. One interpretation of this initial ``particularity" would be to suggest that the very fact that certain information initially presents itself in this local and bounded fashion is an indication that we are dealing with various discrete approximations, presented piecemeal, to phenomena that may in fact ``really" be continuous. Whether or not that is the case, it is not difficult to accept that in its presentation to us in fragmented form, this step in the process is closely allied with the discrete (in a very general sense of the word). For centuries, the modes of restoring continuity to such partial information have been more or less haphazard. A sheaf removes this aspect of haphazardness. Significant is the at once \textit{progressive} and \textit{necessary} nature of the sheaf concept: how by gradually (progressively) covering fragments of reality, and then systematically gluing them together into unique global solutions (necessary), the construction of sheaves encourages us to shift away from our standard ontologies or descriptions of reality as anchored in some ``absolute" towards a more ``contrapuntal and synthetic" perspective capable of registering ``relative universals."\footnote{\cite{zalamea_synthetic_2013} contains an insightful discussion of precisely this latter perspective on sheaves.}\par 
With the sheaf construction, a global vision is not \textit{imposed} on the local pieces, obliterating the local nature of the presented information via some ``sham" generalization, but emerges progressively, step-by-step, through the unfolding of precise translation systems guaranteeing the compatibility of the various components. A sheaf does not attempt to suppress the richness and polyphony of data in its particularity and relative autonomy, coercing a kind of standardized agreement as so many past models of generality have done. A sheaf is like a master composer who is not content to have her harmony prefabricated for her by habitual associations, or who would achieve harmony only at the expense of suppressing all contrapuntal impulses and polyphony, imposing it ``from above," or restraining the local freedom of each voice to roam with some independence from the constraints that bind it in the name of some prefabricated schema; rather, the sheaf-like composer achieves harmony only progressively, first by letting each component part unfold, in its relative autonomy, its own laws, then by insisting on making explicit even the most minute of links and transits between the laws of movement of each of the parts, securing locally smooth passages for each transition, and from the glue or constraints that emerge out of this process, begins to build up a larger ensemble, step by step. It is not a compromise between the local and the global in the name of some idealogical preference for the more global or universal. Sheaves earn their place as true mediators by virtue of their complete realization of the idea that---to paraphrase Hegel---true mediation comes about only from preserving the extremes as such, and true universality comes about only by sinking as deeply as possible into the particular. \par    
The next section considers, in more detail now, three examples---sheaves in the context of manifolds, analytic continuation, ``cross-sections" of a bundle---that were especially significant in the early development of sheaf theory. 
\subsection{Three Historically Significant Examples}
\begin{example}\label{example: manifold}
	Euclidean space, the space of classical geometry, is quite ``nice." In many areas of geometry and physics, one has to deal with fairly complicated structures, and it is desirable to have a description of these things in terms of simpler properties found in Euclidean space. This is where the concept of a \textit{manifold}\index{manifold} comes in. Manifolds are, roughly, topological spaces that look locally like $\mathbb{R}^n$, where \textit{local} here means that every point of a manifold will have an open neighborhood that admits a one-to-one map onto an open set of $\mathbb{R}^n$. In other words, even if globally a manifold does not ``look like" a Euclidean space, locally it will resemble Euclidean space near each of its points.  \par 
	One of the reasons for moving to manifolds is that there are a number of useful instruments available to us in the context of $\mathbb{R}^n$, such as those of integral and differential calculus, that we would like to import to the study of other (more complicated) spaces. Sometimes we have a topological space on which we would like to employ the instruments of, e.g., calculus, and we find that such spaces are locally like open subsets of the Euclidean space $\mathbb{R}^n$ even while they do not provide coordinates valid everywhere. The basic idea with manifolds is that we can cope with such spaces by transferring the instruments that are available in $\mathbb{R}^n$ to small open sets and then patching those sets together, in a sense ``recovering" important aspects of the original topological space, while taking advantage of the ``usual tools." \par 
	In \ref{sec: restrict}, we discussed how functions of various sorts---e.g., continuous, infinitely differentiable, real analytic, holomorphic functions---are not only all continuous, but the condition for such a continuous function to belong to its given class is in fact, in each case, a \textit{local} condition, meaning that whether it to has a property is in fact equivalent to it having that property in the neighborhood of every point in its domain of definition. In short, manifolds and the functions on them are distinguished by being constructed by the pasting together of pieces that have a particular ``nice" property locally. \par 
	All of this---the issue of locality and the notion of ``patching together" local parts---should suggest that sheaves are lurking somewhere in the background. Indeed, sheaves arise in a particularly natural way in the context of manifolds, and were accordingly pivotal in the early historical development of sheaves.\par 
	Formally, a manifold is defined as follows: 
	\begin{definition} An $n$-dimensional topological \textit{manifold}\index{manifold} $M$ is a second countable (its topology has a countable base) Hausdorff space\footnote{A Hausdorff space is one in which distinct points are contained in disjoint open neighborhoods.} such that each point $q \in M$ has an open neighborhood $V$ homeomorphic to an open set $W \subseteq \mathbb{R}^{n}$, i.e., for each point in $q \in M$, there is an open neighborhood of this point such that there exists a map $\phi$ from this open neighborhood into $\mathbb{R}^n$, which must further satisfy: (1) $\phi$ is invertible, i.e., we have $\phi^{-1}: \phi(V) \rightarrow V$; (2) $\phi$ is continuous; (3) $\phi^{-1}$ is continuous.
	\end{definition}  \noindent 
If the reader has never worked with manifolds, this definition may be difficult to parse at first. If the reader is intimidated by this definition, it is fine to just think of a sphere or surface of the globe, or a finite cylinder. \par 
The homeomorphism $\phi: V \rightarrow W (\subseteq \mathbb{R}^{n})$ given in the definition is usually called the \textit{chart map} for the cover (on which more below). This map, together with the open set $V$, gives the pair $(V, \phi)$, called a \textit{chart}\index{chart} for $M$. Note that $\phi$ is just a map $\phi(q) = (\phi^1 (q), \phi^2 (q),\dots, \phi^n (q))$ (with $n$ entries), where each $\phi^j$ is just a map $\phi^j : V \rightarrow \mathbb{R}$ for $j = 1, 2,..., n$. In other words, the result is $n$ many real numbers, $x_1, \dots,x_n$, the result of $n$ component maps acting on the point. These component maps are called the \textit{coordinate maps}, and provide the \textit{local coordinates}.
	\par  
	All we require is that, for every point of the manifold, there exists a chart that contains the point. The globe cannot be represented by any one single flat map, so in using flat maps or charts to navigate the Earth's surface, one must collect many (occasionally overlapping) charts together into an atlas. Similarly, in general it is not possible to describe a manifold via one chart alone. But the whole manifold may be covered by a collection of charts. A collection of charts $\{(V_i, \phi_i) | i \in I\}$ is called an \textit{atlas}\index{atlas} of the manifold $M$ if $M = \bigcup_{i \in I} V_i$. More formally, 
	\begin{definition}
		An \textit{atlas} for a manifold $M$ is an indexed set $\{\phi_i: V_i \rightarrow W_i \}$ of charts such that together all of the elements of the domains $V_i$ cover $M$. 
	\end{definition} \noindent 
	As an example of this, the torus ($3$-d) is a $2$-d manifold, since it can be covered by charts to $\mathbb{R}^2$. Something like a cross, or branching lines, on the other hand, does not yield a manifold, for at the branch points it is not possible to find an invertible continuous map on an open neighborhood of those points. \par 
	The proper understanding of manifolds\index{manifold} is nicely suggested by the case of the $2$-sphere manifold $S^2$, in which context we can provide an atlas consisting of two charts via \textit{stereographic projection}.\index{stereographic projection} Of course, this sphere can be described as the subset $\{(x,y,z) \hspace*{0.3em}| \hspace*{0.3em} x^2 + y^2 +z^2 = 1\}$ in $\mathbb{R}^3$. But we can also describe it \textit{intrinsically}, meaning without reference to the ambient space $\mathbb{R}^3$,  in terms of the points via some parameterization. It is not possible to cover the sphere with a single chart, since the sphere is a compact space, and the image of a compact space under a continuous map is compact, while $\mathbb{R}^n$ is non-compact---so there cannot be a homeomorphism between $S^n$ and $\mathbb{R}^n$. However, it \textit{is} possible to cover the sphere by two charts, as will be seen.\par 
	In order to develop this further, we will first need to review the relevant notions involved in stereographic projection. The basic idea is that we want to ``picture" the sphere as a plane by ``projecting" it onto the plane. We could do this by first isolating the south pole of the sphere, $S = [0,0,-1]$. To parameterize the sphere we consider a point in the equatorial plane given by $z=0$, namely $P = [r,s,0]$, and then draw a line from the south pole through this point. Such a line will intersect the sphere somewhere, say $Q$, and the resulting line from $S$ to $Q$ will be unique and will intersect the plane in exactly one point. We can describe this line in vector form as $ Q = [0,0,-1] + \lambda (S-P) = [0,0,-1] + \lambda(r,s,1)$ (where round brackets denote a vector and square brackets a point). This is of course equal to a point $[\lambda r, \lambda s, \lambda - 1]$. We are looking for what values of $\lambda$ land us on the sphere, i.e., whenever 
	\begin{align*}
	& (\lambda r)^2 + (\lambda s)^2 + (\lambda - 1)^2 = 1 \\
	& \iff \lambda^2(r^2 + s^2 +1) - 2\lambda = 0 \\
	& \iff \lambda = 0 \text{ or } \lambda = \frac{2}{(1+r^2+s^2)}.
	\end{align*}
	Substituting the non-trivial value of $\lambda$ in to our equation for $Q$, we get the point where the line from $S$ to $P$ meets the sphere
	\begin{equation*}
	Q = \Big[\frac{2r}{1+r^2+s^2}, \frac{2s}{1+r^2+s^2}, \frac{1-r^2-s^2}{1+r^2+s^2} \Big]
	\end{equation*}    
	and where for every value $r,s$, we get a point on the sphere, and every point on the sphere, with one exception, is obtained in this way. The sole exception which prevents this from giving a bijection is obviously the point $S = [0,0,-1]$, corresponding to the line tangent to the south pole. \par 
	In other words, then, for every point on the plane, i.e., every $r,s$ from $[r,s,0]$, there is an associated point on the sphere (with the exception of $S = [0,0,-1]$). We thus have a map $\phi^{-1}_s: \mathbb{R}^2 \rightarrow S^2 -\{S\}$, the inverse projection map, or the rational parameterization of the sphere. We will denote the sphere with the south pole removed $S^2 - \{S\}$ by $U_s$. But now note that $S$, $P$, and $[x,y,z]$ all belong to the same line. Thus, $[r,s,0] = \lambda [x,y,z] + (1 - \lambda)[0,0,-1]$, which gives that $\lambda = \frac{1}{z+1}$, and thus that $r = \frac{x}{1+z}, s = \frac{y}{1+z}$. This is in fact the \textit{projection map}, $\phi_s: U_s \rightarrow \mathbb{R}^2$. On all the points for which they are defined, this map and the inverse map are clearly continuous;\index{manifold} thus we have homeomorphisms. This moreover means that we have produced a coordinate chart for $U_s$. We can do exactly the same sort of thing starting with the north pole, with the expected result that all points are covered except for a sole exception, the \textit{north} pole itself. Similar to before, we can let $U_n$ denote $S^2 - \{N\}$. Then we have another homeomorphism pair $\phi^{-1}_n:\mathbb{R}^2 \rightarrow U_n$ and $\phi_n: U_n \rightarrow \mathbb{R}^2$ . \par 
	Note that together, $U_s$ and $U_n$ provide a \textit{cover}\index{cover} of $S^2$, and that together with their respective charts, they give us an atlas for all of $S^2$. Moreover, on the intersection $U_n \cap U_s$ = $S^2 - \{S,N\}$, i.e., the sphere without its poles, we have that $\phi_i(U_n \cap U_s) = \mathbb{R}^2 - \{(0,0)\}$. The idea is that we can obtain \textit{all of} $S^2$ by taking these two homeomorphic copies of $\mathbb{R}^2$, namely $U_s$ and $U_n$, or rather $\phi_s(U_s) \cong W_s \subseteq \mathbb{R}^2$ and $\phi_n(U_n) \cong W_n \subseteq \mathbb{R}^2$, and pasting the $\phi_s(U_s \cap U_n) \cong W_{sn} \subseteq W_s$ to $\phi_n(U_s \cap U_n) \cong W_{ns} \cong W_n$ together via appropriate transition functions (discussed in a moment). It is straightforward to define the transition function between the two projections which enables such gluing. In this way, we will have thus constructed the atlas $\{(U_s, \phi_s), (U_n, \phi_n)\}$, and an atlas determines $M = S^2$ as a topological space. \par 
	Stepping back, a little more generally now, assume given two charts $(U_i, \phi_i)$ and $(U_j, \phi_j)$. The general idea here is captured in something like the following picture: 
	\begin{center} 
		\includegraphics[scale=0.33]{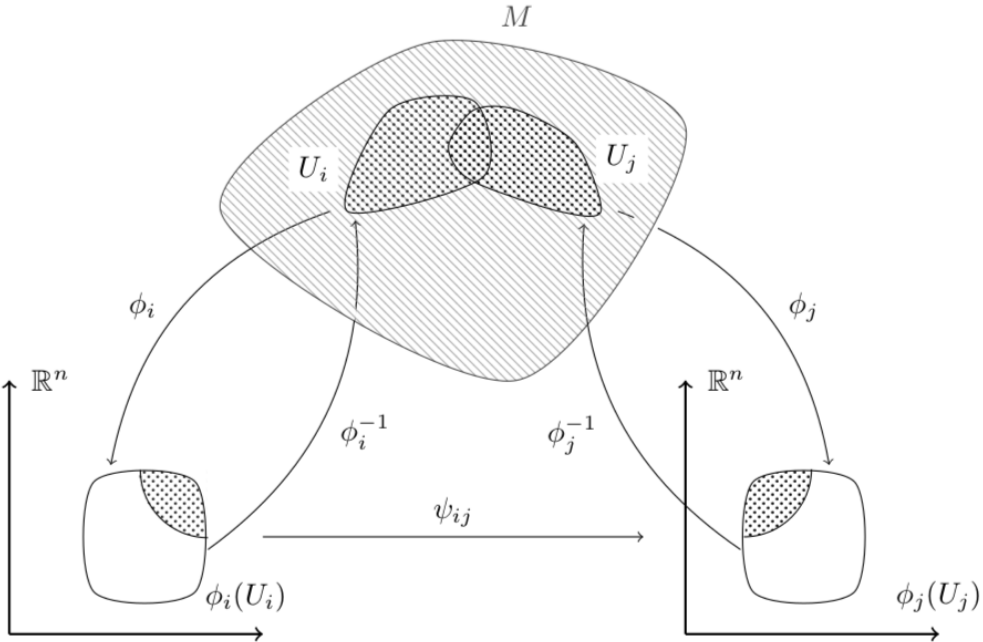}
	\end{center} 
	Two charts $\phi_i$ and $\phi_j$ of an atlas\index{atlas} might overlap\index{manifold} on the set $U_i \cap U_j$. But then by composition of $\phi_i$ with the inclusion $U_i \cap U_j \subseteq U_i \rightarrow W_i$, we get a homeomorphism $\phi_{ij}: U_i \cap U_j \xrightarrow{\cong} W_{ij}$ that goes from the overlap to an open set $W_{ij} \subseteq W_i \subseteq \mathbb{R}^n$. However, while a point in the intersection is mapped to one set of coordinates via $\phi_i$, it is in general mapped to another set of coordinates in $\phi_j$. Accordingly, $\phi_j$ gives an in principle different homeomorphism $\phi_{ji}: U_i \cap U_j \xrightarrow{\cong} W_{ji}$ from the intersection to a different open set $W_{ji} \subseteq W_j \subseteq \mathbb{R}^n$. In the above picture, the dotted regions on the bottom left and right figures in $\mathbb{R}^n$ represent $\phi_{ij} (U_i \cap U_j) \cong W_{ij}$ and $\phi_{ji} (U_i \cap U_j) \cong W_{ji}$, respectively. Notice that, as homeomorphisms, these are not just \textit{maps}, but maps with inverses. \par 
	The question then is: how do we translate or transition between the two homeomorphisms on the overlap? We are looking for a map $\psi_{ij}: W_{ij} \rightarrow W_{ji}$, a rule that enables us to get from a point in one set of coordinates to that same point in terms of the other coordinates. In our example of the manifold $S^2$, such a transition rule is evident and easy to construct explicitly. To find this more generally, notice that the map $\phi_{ij}$ is invertible (it describes a homeomorphism). So we can take a point $q \in \phi_{ij}(U_i \cap U_j) \subseteq \mathbb{R}^n$ and apply the inverse map to it to get back to the manifold. We then apply $\phi_{ji}$, landing us back in $\mathbb{R}^n$, but this time with the coordinates specified by $\phi_j$. In other words, we have constructed the composite map $(\phi_{ji} \circ \phi^{-1}_{ij})(q)$, producing coordinates of $q$ in the other chart, and we can of course do this for any such $q$ in the overlapping region. Thus we are really mapping one open set of $\mathbb{R}^n$ to another, and doing so homeomorphically (the composite of homeomorphisms is a homeomorphism).
	\par 
	With such a \textit{chart transition map}\index{chart!transition map} we thus have the following:
	\begin{center} 
		\begin{tikzcd}
			& U_i \cap U_j \arrow[dl, "{\phi_{ij}}", swap] \arrow[dr, "{\phi_{ji}}"] \\
			\phi_{ij} (U_i \cap U_j) \arrow[rr, "{\phi_{ji} \circ \phi_{ij}^{-1}}", swap] & & \phi_{ji} (U_i \cap U_j)
		\end{tikzcd} 
	\end{center} 
	where, as in the following, the two bottom objects are contained in $\mathbb{R}^n$. 
	\begin{center} 
		\hspace*{-0.2cm}	
		\includegraphics[scale=0.32]{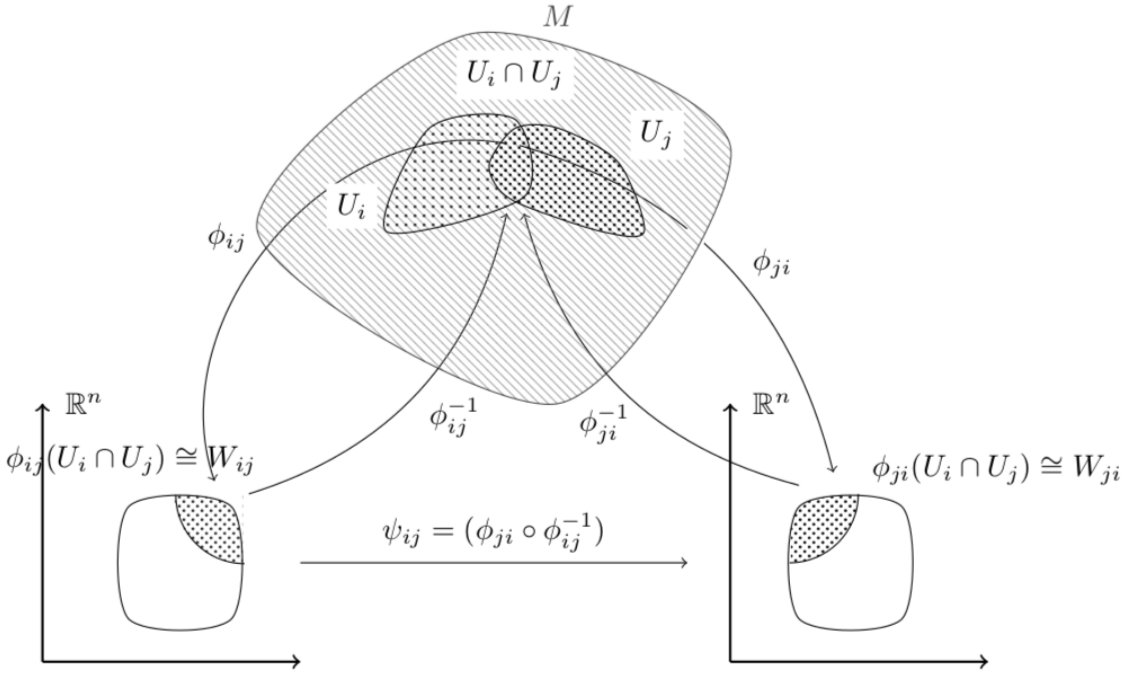}
	\end{center} 
	This transition map $\psi_{ij}$ gives\index{manifold} us our desired \textit{change of coordinates}, a map from an open part of $\mathbb{R}^n$ to another open part of $\mathbb{R}^n$. Since we already had continuous inverses $\phi_i^{-1}$ and $\phi_j^{-1}$, the transition map comprised of the restricted maps is guaranteed to be continuous. Note also that in changing coordinates, from the ``perspective" of $p \in U_i \cap U_j$ (the ``real world"), this $p$ ``does not care" what coordinates it is in, or that you have changed coordinates.\par
	The really important point here is that in general the chart transition maps contain the instructions for \textit{how to glue together the various charts of an atlas}. In other words, with these maps, you have all the \textit{global} information you need. In the case of our example with $S^2$, there were only two charts that needed gluing. But in the general case, the conditions detailed above hold for all $i,j$ in some index set $I$, ultimately allowing an entire manifold $M$ to be covered. This process of obtaining the entire manifold\index{manifold} by taking all the $W_i \subseteq \mathbb{R}^n$ and pasting the $\phi_i(U_i \cap U_j) \cong W_{ij} \subseteq W_i$ to $\phi_j(U_i \cap U_j) \cong W_{ji} \subseteq W_j$ together by the transition functions can be formally described with the following coequalizer diagram: 
	\begin{center}
		\begin{tikzcd}
			\coprod_{i,j} U_i \cap U_j \arrow[r, shift right = 2, "{\beta}", swap] \arrow[r, "{\alpha}"] & \coprod_{i} U_i \arrow[r, "{\gamma}"] & M 
		\end{tikzcd}
	\end{center} 
	Here $\coprod_i$ denotes the coproduct in the category \textbf{Top}; $\gamma$ takes each point $x \in U_i$ and sends it to the same $x \in M$; the maps $\alpha$ takes each point $x_{ij}$ in the intersection $U_i \cap U_j$ to the same $x_{ij}$ in $U_i$, while $\beta$ takes each point $x_{ij}$ in the intersection $U_i \cap U_j$ to the same $x_{ij}$ in $U_j$. $M$ is the coequalizer\index{coequalizer} of these maps $\alpha$ and $\beta$ in the category \textbf{Top}. Given how consideration of locally-defined functions over open sets ordered by inclusion reverses the direction of the arrows, and given the definition of a sheaf in terms of an \textit{equalizer diagram},\index{equalizer} the close connection to sheaves should be apparent.\par 
		Returning to $S^2$ again, a further important feature of all this is that we can now test a function on $S^2$ for certain properties, like continuity or differentiability, by trying it on each of the two parts (each of the two charts) separately. Moreover, a function will be continuous on $V \subseteq M$ when its composite with $\phi^{-1}$ is continuous on $W \subseteq \mathbb{R}^n$. Via the following important \textit{direct image sheaf}\index{sheaf!direct image} construction, we will be able to appreciate how the chart will determine a particular sheaf of continuous functions on $V$. 
		\begin{definition}
			If $f: X \rightarrow Y$ is a continuous map of spaces, then each sheaf $F$ on $X$ yields a sheaf $f_{*}F$ on $Y$ defined, for $V$ open in $Y$, by $(f_{*}F) (V) = F(f^{-1}(V))$. In other words, since $f$ is continuous, it induces an order-preserving map, $f^{-1}$. Composition with $f^{-1}$ then gives $f_{*}F$ defined as the composite functor 
			\begin{equation}
			\mathscr{O}(Y)^{op} \xrightarrow{f^{-1}} \mathscr{O}(X)^{op} \xrightarrow{F} \textbf{Sets}
			\end{equation}
			This sheaf is usually called the \textit{direct image} of $F$ under $f$.\footnote{This direct image sheaf is also sometimes called the \textit{pushforward} of a sheaf.} 
		\end{definition} \noindent 
		Notice that the map $f_{*}$ is in fact a functor 
		\begin{equation}
		f_{*}: \textbf{Sh}(X) \rightarrow \textbf{Sh}(Y)
		\end{equation}
		and that $(fg)_{*} = f_* g_*$, which means that by defining $\textbf{Sh}(f) = f_*$, we have that \textbf{Sh} becomes a functor on the category of all (small) topological spaces. And, in particular, if $f: X \rightarrow Y$ is a homeomorphism, then $f_*$ gives an isomorphism of categories between sheaves on $X$ and sheaves on $Y$.\footnote{As the reader might expect, a continuous map of spaces $f$ will induce another functor going in the other direction on the associated categories of sheaves. Performing the reverse operation yields the \textit{inverse image} sheaf $f^*$, or the \textit{pullback} of a sheaf along a map.\index{sheaf!inverse image (pullback)} Pushforwards and pullbacks are especially useful in providing us with a canonical way to register the effects of changing base spaces. We will look closer at these two functors in a later chapter.} \par 
		Returning to our example, then, our chart will in fact determine a particular sheaf $C_V$ of continuous functions on $V$, specifically as the direct image $C_V = (\phi^{-1})_* C_W$. The coordinate projections $\mathbb{R}^n \rightarrow \mathbb{R}$, first restricted to $W$ then composed with $\phi^{-1}$, give us the local coordinates for the chart $\phi$, i.e., the $n$-coordinate functions $x_1, \dots, x_n: V \rightarrow \mathbb{R}$. Going the other way, these $n$ functions determine the chart as the continuous map $V \rightarrow \mathbb{R}^n$ with components $x_i: V \rightarrow \mathbb{R}$, the map to $\mathbb{R}^n$ being then restricted to its image $W \subseteq \mathbb{R}^n$. \par 
		We now make a key general observation, one that can be applied to our present situation. Consider how if we let $U$ be an open set in a general space $X$, then any sheaf $F$ on $X$, when restricted to open subsets of $U$, clearly gives us a sheaf $F|_{U}$ on $U$. And so, in this way, $U \mapsto \textbf{Sh}(U)$ and $V \subseteq U \mapsto (F|_U \mapsto F|_V)$ defines a contravariant functor on $\mathscr{O}(X)$. Moreover, since the notion of a sheaf is ``local," this suggests that this functor \textit{itself} might be a sheaf.
		\begin{theorem}
			\label{sheafcollation} 
			If $X = \bigcup W_k$ is an open covering of the space $X$, and if, for each $k$, $F_k$ is a sheaf of sets on $W_k$ such that 
			\begin{equation}
			F_k|_{(W_k \cap W_l)} = F_l|_{(W_k \cap W_l)}
			\end{equation} 
			for all indices $k$ and $l$, then there will exist a sheaf $F$ on $X$, unique up to isomorphism, with isomorphisms $F|_{W_k} \cong F_k$ for all indices $k$, which match on the above equation.\footnote{A proof of this can be found in \cite{maclane_sheaves_1994}. Instead of going through this, we note, following Serre, that more generally we could replace the equality in the above theorem with \textit{isomorphisms} of sheaves, showing that sheaves are in fact \textit{collatable up to isomorphism}. Explicitly (since this is actually the version that most concerns us in the present case): If $X$ again has an open cover with the $W_k$, and if, for each index $k$, $F_k$ is a sheaf on $W_k$, so for all $j$ and $k$ 
				\begin{equation}
				\theta_{jk}: F_j|_{(W_j \cap W_k)} \cong F_k|_{(W_j \cap W_k)}
				\end{equation} 
				is an isomorphism of sheaves, and for each $i, j, k$ we have 
				\begin{equation*}
				\theta_{ik} = \theta_{jk} \circ \theta_{ij}
				\end{equation*}
				whenever this is defined; then it can be shown, similar to in the less general case, that there will exist a sheaf $F$ (unique up to isomorphism) on $X$ and isomorphisms $\phi_k: F|{W_k} \rightarrow F_k$ (unique up to isomorphism) such that 
				\begin{equation*}
				\phi_j = \theta_{ij} \phi_i
				\end{equation*}
				 when this is defined.}  
		\end{theorem} \noindent 
		Many sheaves can be built up in this way from the local pieces $F_k$. \par \noindent 
		It is basically immediate from this theorem how we can produce the sheaf of continuous functions on $S^2$, or the sheaf of differentiable functions on $S^2$, etc., and thus ultimately sheaves of smooth structures on $S^2$. \par
	An important point that emerges in this discussion of sheaves on manifolds\index{manifold} is that in the standard definition of sheaf the reader might be misled into thinking that what we need on overlapping regions is \textit{equality} ``on the nose"; but this is not quite right. In fact, all that is actually needed in patching together different information on the same overlapping region is just a consistent \textit{translation system} that can take us from one data-piece to the other and back, i.e., we need a set of isomorphisms allowing us to translate between them, not simple equality. This system of translations typically leads to what is called \textit{descent data}, and can be thought of, intuitively, as providing some sort of generalized matching. \par     
	Before moving on, it is worth noting that the ``meaning" of the manifold only really emerges once we have glued together the charts according to these transition maps. The idea is to study manifolds by ``pushing them down" onto the charts in $\mathbb{R}^n$ with their transitions, thus allowing us to work with far simpler maps from some portion of $\mathbb{R}^n$ to another. Often it is desirable to define properties, like continuity or differentiability, of a ``real-world object" (a map from $\mathbb{R}$ to $M$) by judging suitable conditions on a chart-representative of that real-world object. Say we start with a ``path" or curve $\gamma: \mathbb{R} \rightarrow M$ in a manifold.\index{manifold} We then take an open set $U$ enclosing it, yielding a map $\gamma: \mathbb{R} \rightarrow U$. We can now ``push down" the real world curve ($\gamma$, in bold) to a chart via some chart map $\phi$, to get its image in $\mathbb{R}^n$: 
	\begin{center}  
		\begin{tikzcd}
			\mathbb{R} \arrow[r, "\gamma", thick] & U \arrow[d, "\phi"] \\
			& \phi (U) \subseteq  \mathbb{R}^n 
		\end{tikzcd}
	\end{center}  
	But then, of course, instead of looking at the curve in the real world, we can just focus on the map $\phi \circ \gamma$: 
	\begin{center}  
		\begin{tikzcd}
			\mathbb{R} \arrow[r, "\gamma"] \arrow[dr, "{\phi \circ \gamma}", swap, thick] & U \arrow[d, "\phi"] \\
			& \phi (U) \subseteq  \mathbb{R}^n 
		\end{tikzcd} 
	\end{center}  
	which represents the original curve as some curve in $\mathbb{R}^n$ down in the chart. If we want to know about the continuity of the ``real world" curve, we instead look at continuity in this chart representative of it; but one should not confuse this chart representative for the ``real world" (for instance, one cannot define differentiability on the ``real-world" trajectory given by $\gamma$), or forget that $\phi$ could be ill-defined. In order to overcome this latter issue, and ensure that the chart is not arbitrary, we make sure that a given property (e.g., continuity) does not change under change to another chart. Formally, what we mean is: 
	\begin{center} 
		\begin{tikzcd}
			& \psi(U) \\
			\mathbb{R} \arrow[r, "\gamma"] \arrow[dr, "{\phi \circ \gamma}", swap] \arrow[ur, "{\psi \circ \gamma}"]& U \arrow[d, "\phi"] \arrow[u, "\psi"] \\
			& \phi (U) \subseteq  \mathbb{R}^n 
		\end{tikzcd} 
	\end{center} 
	Then if we know, for instance, that $\phi \circ \gamma$ is continuous, we also want to know if we can say the same for $\psi \circ \gamma$. But this arrow is given by another path around the diagram, in particular,
	\begin{equation*}
	\psi \circ \gamma = \psi \circ (\phi^{-1} \circ \phi) \circ \gamma = (\psi \circ \phi^{-1}) \circ \phi \circ \gamma. 
	\end{equation*}
	But we know that $(\psi \circ \phi^{-1})$ is continuous, and we know that $\phi \circ \gamma$ is continuous, and the composition of continuous functions is continuous. So $\psi \circ \gamma$ must be continuous as well. The map $\phi \circ \gamma$ is thus well-defined and continuous, and since the downward map $\phi$ is continuous, we can show that we can ``lift" continuity up to $\gamma$. But notice that if we had said that $\phi \circ \gamma$ was differentiable, on the other hand, that does not automatically guarantee that $\psi \circ \gamma$ is differentiable. Everything depends, in other words, on the nature of the transition map.\footnote{In terms of differentiability: the transition map might preserve continuity but introduce an edge, preventing differentiability. But if we just ``rip out" all those charts in our atlas that are only continuous, but not differentiable, then we get a \textit{restricted atlas} where all of the transition functions are differentiable.}\par 
	We can be more general. Two charts $(U, \phi)$ and $(V, \psi)$ of a topological manifold are said to be $blank$-compatible if either $U \cap V = \emptyset$ or whenever $U \cap V \neq \emptyset$, we have that the following transition maps 
	\begin{align*}
	\psi \circ \phi^{-1}: \phi(U \cap V) \rightarrow \psi(U \cap V) \\ 
	\phi \circ \psi^{-1}: \psi(U \cap V) \rightarrow \phi(U \cap V) 
	\end{align*} 
	have the $blank$ property. We then take a restriction of the maximal atlas\index{atlas} to get an atlas $\mathscr{A}_{blank}$, and say that this atlas is a $blank$-compatible atlas as long as any two charts in $\mathscr{A}_{blank}$ are $blank$-compatible. A $blank$-manifold\index{manifold} is then a triple $(M, \mathscr{O}, \mathscr{A}_{blank})$.\par 
	Of course, $blank$ could be, for instance, the property of differentiability. But it can be many other things as well. For instance, it could represent: $C^0$ (trivially, because \textit{every} atlas is a $C^0$ atlas); $C^1$ (i.e., differentiable once, the result of which is continuous); $C^k$ (i.e., $k$-times continuously differentiable); $C^{\infty}$ (i.e., continuously differentiable arbitrarily many times); $C^{\omega}$ (i.e., real analytic, meaning there exists a (multi-dimensional) Taylor expansion); $\mathbb{C}^{\infty}$ (i.e., complex differentiable, meaning that each continuous map from $\mathbb{R}^n$ to $\mathbb{R}^n$ satisfies the Cauchy-Riemann equations---this gives a \textit{complex manifold}). In this list, our atlases are becoming more and more restrictive, as we place stronger and stronger conditions on the transition functions. A well-known theorem informs us that, as long as $k \geq 1$, any $C^k$-atlas $\mathscr{A}_{C^{k}}$ of a topological manifold will contain a $C^{\infty}$-atlas. There are some topological manifolds where you simply cannot remove a chart or charts in such a way that all that remains is still an atlas and has continuously differentiable transition functions. However, what the theorem says is that whenever you can achieve that the transition functions are at least once continuously differentiable for an atlas, you can be sure that by removing more and more charts you will eventually have a $C^{\infty}$-atlas. This guarantees that we may always, without loss of generality, consider $C^{\infty}$-manifolds, which are called \textit{smooth manifolds} (``always" meaning as long as we guarantee $C^1$).\par 
	In terms of the general picture given above, we will be able to use the homeomorphism $\phi_i$ to transfer smoothness on $W_i \subseteq \mathbb{R}^n$ to smoothness on $U_i \subseteq M$, giving us the sheaf $C_i^k$ of all smooth functions on the open subsets of $U_i$. Then the smoothness expected to hold of the transition functions guarantees that the sheaves $C^k_i$ and $C^k_j$ will agree when restricted to the overlap $U_i \cap U_j$. Thus, via the result from \ref{sheafcollation}, we know that these sheaves $C^k_i$ themselves can be collated, yielding the sheaf $C^k$ of \textit{all} smooth functions on opens of $M$. Moreover, this sheaf provides an instance of the general notion of a \textit{subsheaf},\index{sheaf!sub-} the description of which in a sense nicely exhibits the local character of a sheaf.   
	\begin{definition}
		If $F$ is a sheaf on $X$, then a subfunctor $S \subseteq F$ is a \textit{subsheaf} iff, for every open set $U$, every element $f \in FU$, and every open covering $U = \bigcup U_i$, we have that $f \in SU$ iff $f|_{U_i} \in SU_i$ for all $i$.
	\end{definition} \noindent 
Our particular sheaf will be a subsheaf of the sheaf of continuous functions; obviously, being a sheaf, its restriction to each of the $U_i$ gives a sheaf $C^k_i$. This is just to say that each smooth manifold supports what is sometimes called a \textit{structure sheaf}, namely the sheaf $C^k$ of smooth functions. This is, importantly, not just a sheaf of \textit{sets}, but a sheaf of $\mathbb{R}$-algebras, a sheaf of rings. \par 
For the Euclidean $n$-space $X = \mathbb{R}^n$, more generally, there are a number of examples of sheaves, which then give rise to a number of subsheaves. For $U$ open in $\mathbb{R}^n$, we let $C^k U$ be the set of all $f: U \rightarrow \mathbb{R}$ having continuous partial derivatives of all orders up to (and including) order $k$. This gives a functor $C^k : \mathscr{O}(X)^{op} \rightarrow \textbf{Set}$ with values in \textbf{Set} (or in $\mathbb{R}$-$\textbf{Mod}$), which in fact yields a sheaf. Notice that \textit{each} $C^k$ will be a sheaf on $\mathbb{R}^n$, leading to a nested sequence of subsheaves on $\mathbb{R}^n$: 
		\begin{equation*}
		C^{\infty} \subseteq \cdots \subseteq C^k \subseteq C^{k-1} \subseteq \cdots \subseteq C^1 \subseteq C^0 = C.
		\end{equation*}
There are other constructions found in the theory of manifolds\index{manifold} that in fact lead to sheaves, for instance, via the notion of a tangent bundle for a smooth ($C^{\infty}$) manifold $M$. The story is pretty similar for various types of manifolds. But perhaps the most important, historically, involves complex analytic manifolds, which brings us to the next example.
\end{example}
\begin{example}
	Historically, the idea of a sheaf largely stems from problems surrounding the \textit{analytic continuation} of functions (together with the development of Riemann surfaces). Analytic continuation\index{analytic continuation} involves the attempts to extend the given domains of definition of functions to larger domains. \par 
	Recall that a function $f: U \rightarrow \mathbb{C}$, where $U$ is open in $\mathbb{C}^n$, is \textit{analytic} (or \textit{holomorphic}) if it is described by a convergent power series in a neighborhood of each point $p \in U$. Another way of saying this is that a function is analytic precisely when its Taylor series about a point, say $x_0$, converges to the function in some neighborhood for every $x_0$ in its domain, i.e., locally it is given by a convergent power series. Now, let $V$ be an open subset of $\mathbb{C}^n$, and for each open subset $U \subseteq V$, let $A_V (U)$ be the set of all analytic functions on $U$. Then $A_V$ will in fact give us a sheaf (specifically, of $\mathbb{C}$-algebras, and of rings).\par 
	Let us explore the connection between analytic continuation and sheaves more closely. We begin with a pair $(D, f)$, where $D \subseteq \mathbb{C}$ is a domain and $f: D \rightarrow \mathbb{C}$ is analytic. We say that another such pair $(D_1, f_1)$ is a \textit{direct analytic continuation} or an analytic extension of $(D, f)$ if $D \cap D_1 \neq \emptyset$ and $f|_{D \cap D_1} = f_1|_{D \cap D_1}$. As one should now expect, this is called the \textit{gluing condition}. In such a case, we then define 
	\begin{equation*}
	g: D \cup D_1 \rightarrow \mathbb{C} 
	\end{equation*}
	by 
	\begin{equation*}
	g|_D = f, \hspace{2em} g|_{D_1} = f_1 \hspace*{0.5em} .
	\end{equation*} 
	Analyticity is locally defined, so we can observe that $g$ must be analytic, and we say that $g$ is obtained by \textit{gluing $f$ and $f_1$}. One of the interesting notions here is that $f$ may be defined very differently than how $f_1$ is defined. In general, an analytic function could be given by a power series, by a formula, or by an integral. The main question that presents itself in this setting is then: what is the \textit{largest} open set to which a function can be extended, and can we describe such an extended function in some way?\par 
	What follows is a simple, concrete illustration of these ideas. Take 
	\begin{equation*} 
	D = \{z : \hspace*{0.3em} \abs*{z} < 1\}
	\end{equation*} 
	and $f(z) = 1 + z + z^2+ z^3 + \cdots$, the geometric series. Then $f(z)$ is a power series, centered at $z$. The radius of convergence is the unit disk, i.e., inside this disk the function always represents an analytic function (if you write out the Taylor expansion at zero, you will get back this power series). Now consider 
	\begin{equation*} 
	D_1 = \mathbb{C}\setminus \{1\} 
	\end{equation*} 
	and $f_1 = \frac{1}{1-z}$. This is a (reciprocal of a) polynomial, and so as long as the denominator does not vanish, it will be analytic everywhere, i.e., an \textit{entire} function. But writing out the Taylor expansion of this function at zero gives exactly the same thing as for $f$, and in fact $(D_1, f_1)$ is an analytic continuation\index{analytic continuation} of $f$. We can now observe that while $f$ lived only strictly inside the unit disk, $f_1$ (to which $f$ is equal!), lives \textit{everywhere} except at points \textit{on} the unit circle. The function thus extends to all of the complex plane (except for the points on the unit circle)! The moral of the story is that just because the power series representation lives only in a certain region does not mean that the analytic function that it represents lives only there. It is worth noting, however, that analytic extension is not always possible.\par 
	Now, the gluing condition stated above already suggests the close connection with sheaves. We also have the required uniqueness. If $(D_1, f_1)$ and $(D_1, f_2)$ are two analytic extensions of $(D, f)$, defined on the same domain $D_1$, then $f_1 = f_2$. For note that it is immediate from the fact that both $f_1$ and $f_2$ are analytic extensions of $f$ that if you take $f_1|_{D \cap D_1}$, this is equal to $f|_{D \cap D_1}$ which is in turn equal to $f_2|_{D \cap D_1}$. To check that the two analytic functions are equal on a domain (specifically the restricted domain in this case), all you have to do is find one convergent sequence in the domain that converges to a point in the domain, and verify that for each point of the sequence the two analytic functions take the same value.\par 
	So far, we only mentioned a direct analytic extension from one domain to another. But we can extend this idea from a pair of domains to any finite number of domains. The idea is that you do not insist that they \textit{all} intersect---only that there is an ordering such that every member intersects with the next. More formally: suppose we have pairs $(D_{\alpha}, f_{\alpha})$, where the $D_{\alpha}$ are open connected domains and the $f_{\alpha}: D_{\alpha} \rightarrow \mathbb{C}$ are analytic, as well as a total ordering of the indexing set consisting of the $\{\alpha\}$. For simplicity, let $A$ be finite, i.e., $A = \{\alpha_1 < \alpha_2 < \cdots < \alpha_m\}$. Now, for all $i$, let $D_{\alpha_i} \cap D_{\alpha_{i+1}} \neq \emptyset$ and 
	\begin{equation*}
	f_{\alpha_i}|_{D_{\alpha_i} \cap D_{\alpha_{i+1}}} = f_{\alpha_{i+1}}|_{D_{\alpha_i} \cap D_{\alpha_{i+1}}}
	\end{equation*}
	In other words, you just have a chain of direct analytic continuations.\index{analytic continuation} In this case, we say that $(D_{\alpha_m}, f_{\alpha_m})$ is the (indirect) analytic continuation of $(D_{\alpha_1}, f_{\alpha_1})$. This process of generating a chain of direct analytic continuations, i.e., interlocking or overlapping regions of extension of various locally-defined functions whereby the domain of that function is extended step by step, might be pictured in the following way (where to each $D_{\alpha_i}$ is of course attached a function $f_{\alpha_i}$): 
	\begin{center}
		\hspace*{-1cm}
		\includegraphics[scale=0.6]{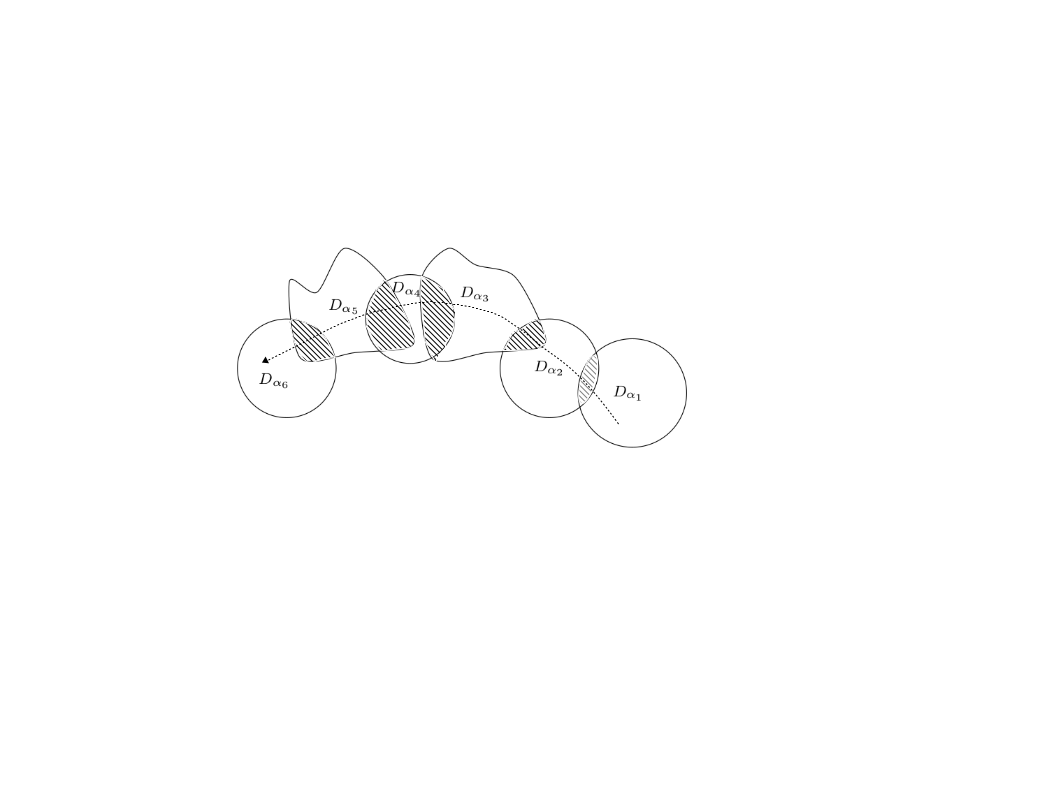}
		\vspace*{-6cm}
	\end{center}
	Now, something interesting can happen in such situations. We can have such a chain of analytic continuations such that it happens that the final domain is the same as the first domain, but the function we get will be completely different. What is in fact happening in such a case is that these functions are two so-called \textit{branches} of an analytic function. (So this process of analytic continuations\index{analytic continuation} helps you find all possible branches of an analytic function.) Note also that this process is only non-trivial and useful if there are singularities. If the function were already entire, there would be nothing to extend. \par 
	Stepping back a bit, and recalling the previous example: smooth and complex analytic manifolds are particular examples of what are called \textit{ringed spaces}, where a ringed space $X$ is a topological space equipped with a fixed sheaf $R$ of rings, in which setting such a sheaf is called the structure sheaf. In addition to leading to the notion of a Riemann surface, the consideration of analytic functions of one complex variable leads to the notion of a \textit{bundle}. The next example provides a closer inspection of this notion, and covers a few notable abstract results that emerge in that context.  
\end{example}
\begin{example}
	In the previous example, we mentioned that the consideration of holomorphic functions of one complex variable can lead to the notion of a \textit{bundle}.\index{bundle} In what follows, we introduce the concept of a bundle, a notion that in a sense can be thought of as capturing the underlying set-theoretic structure of the sheaf idea. We will begin by considering a bundle of sets (over the discrete base space $I$).\par 
	We assume for the moment that we are starting out with a collection $\mathscr{Y}$ of pairwise disjoint sets, sets we assume have some set $I$ of indices or labels. Then a bundle is the \textit{entire} structure depicted in the following picture: \par  
	\begin{center} 
		\begin{tikzpicture}[scale=0.95,x=1cm,y=1cm,circle dotted/.style={dash pattern=on .05mm off 6mm,
			line cap=round}]
		\filldraw[fill=white, draw=black,thick] (0.2,1.1) node{} ellipse (3.9 and 0.7);
		\node(text) at (-5.89,1) {\large Base Space: $I$};
		\path[black,->,>=stealth, thick] (-4.3,4.2) edge node[auto, swap]{\large $p$} (-4.3,1.8);
		\filldraw[fill=white, draw=black,thick] (-2.8,6) node{} ellipse (0.2 and 1.8);
		\node(text) at (-2.75, 1.1) {\footnotesize $\bullet{i}$};
		\path[black,-] (-2.8,4.2) edge node[auto]{} (-2.8,1.1);
		\node(text) at (-2.8, 8) {$Y_i$};
		\draw[line width = 1mm,circle dotted] (-2.8,4.8) -- (-2.8,7.8);
		\filldraw[fill=white, draw=black,thick] (-1,6) node{} ellipse (0.2 and 1);
		\node(text) at (-0.95, 1.2) {\footnotesize $\bullet{j}$};
		\path[black,-] (-1,5) edge node[auto]{} (-1,1.25);
		\node(text) at (-1, 7.2) {$Y_j$};
		\draw[line width = 1mm,circle dotted] (-1,5.2) -- (-1,7);
		\filldraw[fill=white, draw=black, thick] (0.5,6) node{} ellipse (0.25 and 1.4);
		\node(h) at (0.57, 0.82) {\footnotesize $\bullet{k}$};
		\path[black,-] (0.5,4.6) edge node[auto] {} (0.5,0.82);
		\node(text) at (0.5, 7.6) {$Y_k$};
		\draw[line width = 1mm,circle dotted] (0.5,5.2) -- (0.5,7.3);
		\filldraw[fill=white, draw=black,thick] (2,5.5) node{} ellipse (0.15 and 0.7);
		\node(h2) at (2.05, 1.3) {\footnotesize $\bullet{l}$};
		\path[black,-] (2,4.8) edge node[auto]{} (2.02,1.3);
		\node(text) at (2, 6.4) {$Y_l$};
		\draw[line width = 1mm,circle dotted] (2,5.1) -- (2,6.2);
		\filldraw[fill=white, draw=black,thick] (3.5,6) node{} ellipse (0.2 and 1.7);
		\node(text) at (3.65, 1) {\footnotesize $\bullet{m}$};
		\path[black,-] (3.5,4.3) edge node[auto]{} (3.5,1);
		\node(text) at (3.5, 8) {$Y_m$};
		\draw[line width = 1mm,circle dotted] (3.5,4.9) -- (3.5,7.5);
		\draw[decorate,decoration={brace,amplitude=0.3cm},xshift=-1.6cm,yshift=0pt] (-2.6,4.3) -- (-2.6,7.9) node[black,midway,xshift=-1.5cm] {\large $p^{-1}(I) \cong Y$} node[black, midway, xshift=-1.5cm, yshift=0.5cm] {\large Stalk Space:};
		\draw[-latex, dashed] (4.8,6.5) node[right]{}
		to[out=200,in=0] (3.6,6.2);
		\draw[-latex, dashed] (4.8,6.5) node[right]{}
		to[out=190,in=0] (3.6,6.8);
		\draw[-latex, dashed] (4.8,6.5) node[right]{}
		to[out=195,in=10] (3.55,7.4);
		\draw[-latex, dashed] (4.8,6.5) node[right]{}
		to[out=190,in=10] (3.6,5.5);	
		\draw[-latex, dashed] (4.8,6.5) node[right]{}
		to[out=190,in=10] (3.6,4.9);
		\node (g) at (5.3,6.5) {\footnotesize Germs};
		\node (g1) at (5.3,6.2) {\footnotesize at $m$};
		\draw[-latex, dashed] (4.5,9.4) node[right]{}
		to[out=200,in=10] (3.6,8.15);
		\node (y) at (4.3,9.9) {\footnotesize Stalk/Fiber over m:};
		\node (y2) at (4.3,9.6) {\footnotesize $p^{-1}(\{m\}) = \{y : p(y) = m\}$};
		\end{tikzpicture}
	\end{center} 
Let us unpack this picture a bit, and clarify some of the terminology. $Y_i$ is called the \textit{stalk}, or \textit{fiber},\index{stalk} over $i$. The map $p: Y \rightarrow I$ is distinguished by the fact that if $y \in Y$, then there exists precisely one $Y_i$ such that $y \in Y_i$, in which case we set $p(y) = i$. This just means that every member of the $Y_i$ gets sent to $i$. Thus, the stalk or fiber $Y_i$ can be got as the inverse image under the projection map $p$ of $\{i\}$: 
\begin{equation*}
p^{-1}(\{i\}) = \{y \hspace*{0.3em}| \hspace*{0.3em} p(y) = i \} = Y_i .
\end{equation*} 
A \textit{stalk}\index{stalk} over some $i \in I$ can also be seen as consisting of a collection of \textit{germs},\index{germ} where the germs at $i$ are just the members of $Y_i$, represented by dots in the above picture. For simplicity, focusing on just a part of $I$, namely $U = \{j,k,l\}$, suppose we have
\begin{center}
	\includegraphics*[scale=0.27]{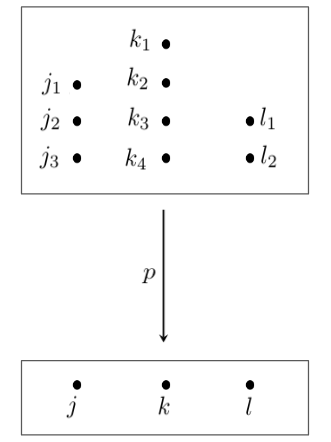}
\end{center}
where, of course, each element of the space ``above" is mapped via $p$ to the element of $I$ directly below it, e.g., $p(j_1) = p(j_2) = p(j_3) = j, p(l_1) = p(l_2) = l$, etc. Then the fiber or stalk over $l$, for instance, is of course just $\{l_1, l_2\}$, while the fiber over $j$ is just $\{j_1, j_2, j_3\}$, and so on. Referring back to the main picture, the set $Y = \{y \hspace*{0.3em}| \hspace*{0.3em} y \in Y_i \text{ for some } i\}$ consisting of all the elements in \textit{any} fiber over $I$, is then called the \textit{stalk space} (or \textit{l'espace \'etal\'e}) of the bundle. And the entire structure is then called a \textit{bundle} (of sets) over the base space $I$.\index{bundle} \par 
This bundle construction is possible whenever there are functions. More specifically, given $p: Y \rightarrow I$ an arbitrary function from a set $Y$ to $I$, we can define the bundle $\mathscr{Y}$ of sets over $I$ with stalk space $Y$ by first defining $Y_i$ as the preimage $p^{-1}(\{i\})$ for each $i \in I$, and then 
\begin{equation*}
	\mathscr{Y} = \{p^{-1}(\{i\}) \hspace*{0.3em}| \hspace*{0.3em} i \in I \} = \{Y_i \hspace*{0.3em}| \hspace*{0.3em}i \in I \} .
\end{equation*}
	But then it should be obvious that a bundle of sets over $I$ is basically just a function with codomain $I$, and so equivalent to set-valued functors defined on $I$, where $I$ is a discrete category. In other words, this is nothing other than the comma (slice)\index{category!slice} category $(\textbf{Set}  \downarrow I)$ of functions with codomain $I$, though it is also common to denote the category of bundles\index{bundle} over $I$ by $\textbf{Bn}(I)$. The reason for discussing such things is that a sheaf can be defined as a bundle with some extra topological structure. We first redefine a bundle in a topological context as follows: 
	\begin{definition} For any topological space $X$, a \textit{bundle over} $X$ is a topological space $Y$ equipped with a continuous map $p: Y \rightarrow X$. 
	\end{definition} \noindent 
As we did a moment ago, you can continue to imagine $Y$ as ``sitting above" $X$, and the map $p$ as \textit{projecting} the points of $Y$ onto their ``shadows" $p(y) \in X$. The bundle is all of this, i.e., a triple $(Y, p, X)$, where $Y$ is the total space, $X$ the base space, and $p$ the projection map. \par 
Notice from the definition that this is just to say that (topological) bundles\index{bundle} are the objects of the slice category $(\textbf{Top} \downarrow X)$,\index{category!slice} where an arrow $f: p \rightarrow p'$ is a continuous (in fact, open) map $f: Y \rightarrow Y'$ such that $p' \circ f = p$, where $p'$ is a map from $Y'$ to $X$. The close connection between the resulting category $(\textbf{Top} \downarrow X)$ (or sometimes just $\textbf{Top}(X)$) and the category of sheaves over $X$ will emerge in what follows. It will turn out that \textit{every} sheaf can be regarded as arising from a bundle; and, conversely, every sheaf gives rise to a bundle. To approach this important connection, we first introduce another important notion: 
	\begin{definition}
		A \textit{cross-section} (or just \textit{section})\index{section} of a bundle $p: Y \rightarrow X$ is a continuous map $s: X \rightarrow Y$ such that $p \circ s = \text{Id}_X$. 
	\end{definition} 
	\noindent  
	In general, a section $s$ for a map $p$ can be regarded as a procedure that at once picks out an element from each of the fibers of $p$. 
	In terms of the earlier simplified case, where we were dealing with a set bundle\index{bundle} over a discrete space, we can depict one such (cross) section $s_1$ as follows: 
	\begin{center}
		\includegraphics*[scale=0.3]{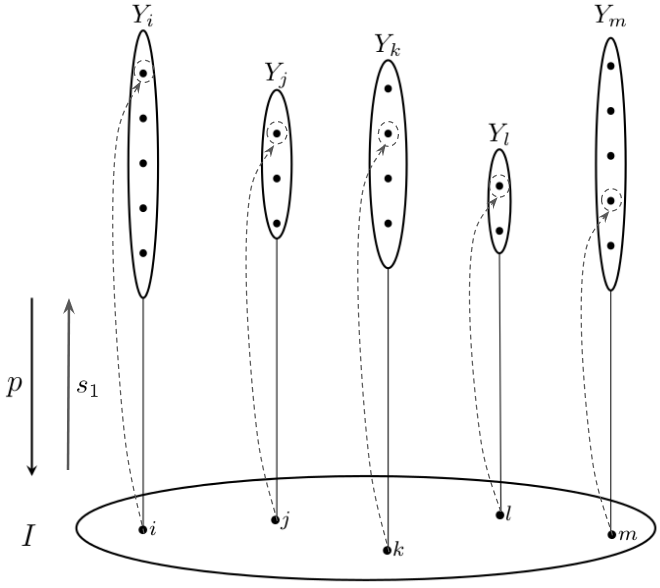}
	\end{center}
Then another section $s_2$ might be given by
\begin{center}
	\includegraphics*[scale=0.3]{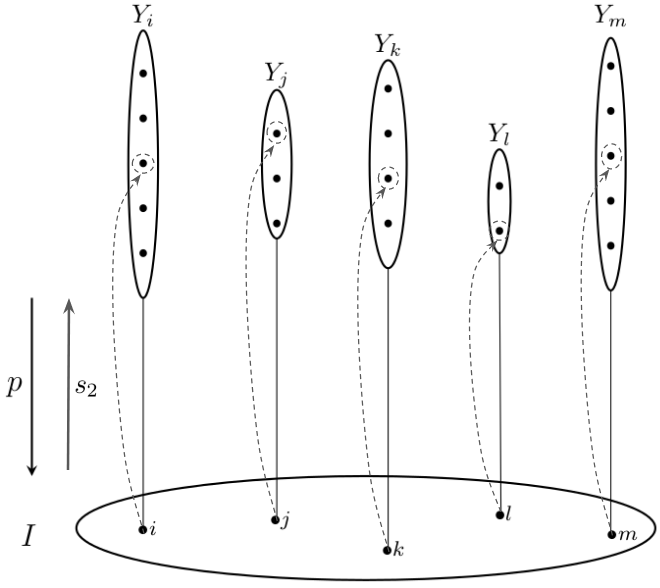}
\end{center}
Still more evocatively, we might picture a particular section with something like 
\begin{center}
	\includegraphics[scale=0.35]{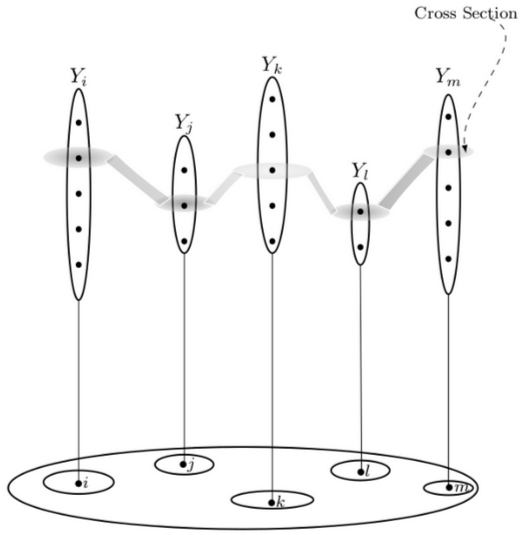}
\end{center}
where the ribbon connecting the individual selections from each fiber is meant to anticipate the gluing process for the same sheaf, i.e., where the various individual selections are glued together (via the topology of $Y$) to create a section of a larger open set (the entire ribbon representing a global section). Altogether, we think of a \textit{bundle} as the indexed family of fibers $p^{-1}(x)$, one per point $x \in X$, glued together by the topology of $Y$. 
\subsection{Bundles to (Pre)Sheaves}
A little more generally now, suppose we have a bundle\index{bundle} $(Y, p, X)$. If we consider some open subset $U$ of the base space $X$ for the bundle $p: Y \rightarrow X$, then $p$ clearly restricts to a map $p_U: p^{-1}(U) \rightarrow U$, and this map will itself yield a bundle (now over $U$). Then, the following diagram, with horizontal arrows as inclusions, will be a pullback in \textbf{Top}, the ``best" way of completing two given morphisms into a commutative square:  
	\begin{center} 
		\begin{tikzcd}
			p^{-1}U \arrow[r] \arrow[d, "{p_U}", swap] & Y \arrow[d, "p"] \\
			U \arrow[r, "i", swap] \arrow[ur, dashed, "s"] & X
		\end{tikzcd}
	\end{center} 
But this lets us define a cross-section (section)\index{section} $s$ of this bundle $p_U$: namely, a section of the bundle $p$ over $U$ is a continuous map $s: U \rightarrow Y$ such that the composite $p \circ s$ is equal to the \textit{inclusion} $i: U \rightarrow X$, i.e., $p \circ s = \text{Id}_U$.\footnote{Note that it may very well happen that a map admits a locally-defined section over a subset $U \subseteq X$, but not a global section, i.e., a section over all of $X$.}
In the case of discrete spaces (where every subset is open, making $p$ automatically continuous), things are especially easy to present. For instance, suppose we have $p$ on $U = \{j,k,l\}$ as earlier 
\begin{center}
	\includegraphics*[scale=0.27]{FiberThingy.png}
\end{center}
Then a section will of course just map each point of $U$ to a point in $Y$ that sits directly above it. For instance, one section $s_1$ for such a $p$ over $U$ will be given by 
\begin{center}
	\includegraphics*[scale=0.25]{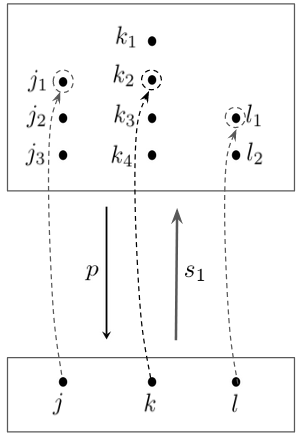}
\end{center}
while another $s_2$ might be given by 
\begin{center}
	\includegraphics*[scale=0.25]{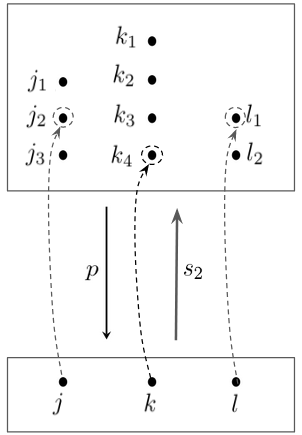}
\end{center}
Observe that there would be a total of 36 distinct sections\index{section} over $U = \{j,k,l\}$, given such a map $p$. Similarly, for more general spaces, here is a picture of sections over an open $U \subseteq  X$
\begin{center}
	\includegraphics*[scale=0.25]{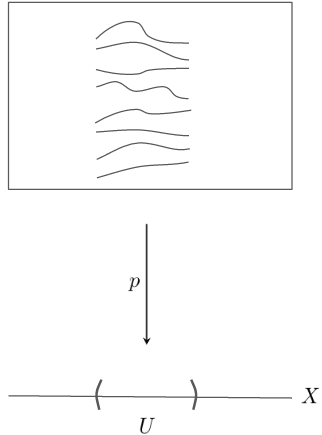}
\end{center} 
Collecting the sections together lets us define 
	\begin{equation}
	\Gamma_p U = \{s \hspace*{0.3em}| \hspace*{0.3em} s: U \rightarrow Y \text{ and } p \circ s = i: U \subseteq X\}
	\end{equation}
	the set of \textit{all} cross-sections over $U$. But now observe that whenever $V \subseteq  U$, we will have the induced (restriction) operation $\Gamma_p U \rightarrow \Gamma_p V$, restricting a function to a subset of its domain. Via this assignment on objects (open sets) $U$, and the induced restriction operation, altogether this just tells us that $\Gamma_p (\--)$ defines a presheaf 
	\begin{equation*} 
	\Gamma_p: \mathscr{O}(X)^{op} \rightarrow \textbf{Set}. 
	\end{equation*}
	On objects $U$, $\Gamma_p(U)$ supplies all the sections over $U$, and a given element $s \in \Gamma_p(U)$ will just be a choice of an element from each fiber over $U$. On arrows, the presheaf $\Gamma_p$ acts by \textit{restriction}\index{restriction} (for every subset inclusion $V \subseteq  U$). Given an inclusion $V \subseteq  U$, and given a section $s$ over $U$, we just restrict this to what $s$ does on $V$, 
	\begin{align*}
	\Gamma_p(U) & \rightarrow \Gamma_p(V) \\
	s & \mapsto s|_V .  
	\end{align*}
	It should be clear that, given a function $s$ on $U$, one can test ``locally" whether or not $s$ amounts to a section. But this locality suggests something more: namely that $\Gamma_p$ is not just a presheaf, but in fact forms a sheaf on $X$, called the \textit{sheaf of cross-sections} (or \textit{sheaf of sections})\index{section} of the bundle\index{bundle} $p$. \par 
	To see that $\Gamma_p$ is in fact a sheaf, not just a presheaf, we need to see that it satisfies the sheaf condition for every cover. Let us explore how this works for the particularly simple case of the discrete sets we have been working with. Suppose given the sets $U_1 = \{j,k\}$ and $U_2 = \{k,l\}$, so that, together, we have a covering of the set $U = U_1 \cup U_2 = \{j,k,l\}$. We have a presheaf $\Gamma_p$ and a cover. We can then define a matching family for the cover: this will, of course, just be a section $t_1$ given over $U_1$, such as 
	\begin{center}
		\includegraphics*[scale=0.3]{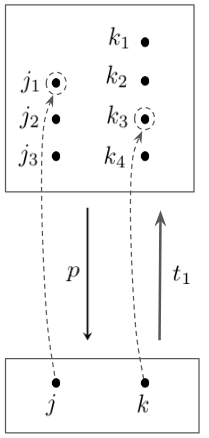}
	\end{center}
	and a section $t_2$ over $U_2$, such as 
	\begin{center}
		\includegraphics*[scale=0.3]{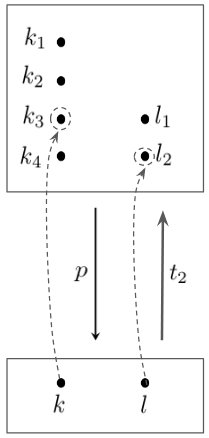}
	\end{center}
	where these agree on the overlapping set $U_1 \cap U_2 = \{k\}$, as these particular sections $t_1$ and $t_2$ do (since they both map $k$ to $k_3$). Pairs of sections that match on the overlap, such as $t_1$ and $t_2$, can then be \textit{glued} together to yield a single section $t \in \Gamma_p(U_1 \cup U_2)$ over the entire set $U = U_1 \cup U_2$: 
	\begin{center}
		\includegraphics*[scale=0.3]{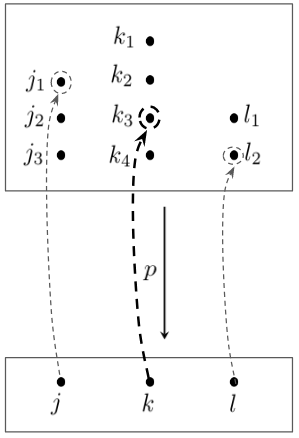}
	\end{center}   
	Moreover, observe how this section $t$ is such that $t|_{U_1} = t_1$ and $t|_{U_2} = t_2$. In this way, we can build up sections over a space by gluing together sections given over local parts.\par 
	This example suggests how we can build a sheaf from the presheaf $\Gamma_p$. This is a general procedure: proceeding in fundamentally the same way as above, every bundle over $X$ will give rise to a sheaf on $X$. Moreover, we know that in the category of bundles $\textbf{Bn}(X)$ on $X$, given (objects) bundles $p: Y \rightarrow X$ and $p': Y' \rightarrow X$, a morphism $p \rightarrow p'$ from the first to the second is just a continuous map $f: Y \rightarrow Y'$ making the triangle $p' \circ f = p$ commute. But each such map $p \rightarrow p'$ of bundles over $X$ (or maps in the slice category) will induce a map $\Gamma_p \rightarrow \Gamma_{p'}$ of (pre)sheaves on $X$, just as one would expect (a morphism of presheaves). This means we actually have a functor $\Gamma$ going from the category of bundles on $X$ to presheaves on $X$: 
	\begin{equation*}
	\Gamma: \textbf{Bn}(X) \rightarrow \textbf{Set}^{\mathscr{O}(X)^{op}}. 
	\end{equation*}
	In fact, this actually gives a functor from bundles to \textit{sheaves}, a fact we will highlight in a moment. Earlier, though, we mentioned that we could also move in the other direction, i.e., that \textit{every} sheaf is in fact a sheaf of sections\index{section} of a suitable bundle.\index{bundle} The next section is devoted to seeing how every (pre)sheaf on $X$ can be regarded as a (pre)sheaf of sections of some bundle. 
\subsection{(Pre)Sheaves to Bundles}	
	Suppose we start with a presheaf 
	\begin{equation*}
	P: \mathscr{O}(X)^{op} \rightarrow \textbf{Set}. 
	\end{equation*}
	We want to use $P$ to construct a bundle over $X$. We will want to construct a collection of sets, indexed by the points $x \in X$, take a union of these sets, and then place a topology on this, leaving us with a space from which we can then define a map down to $X$. This will give us our bundle. \par 
	To see how to construct the relevant sets in our desired collection, first observe that the presheaf $P$ acts on \textit{open sets}, so it does not yet give us sets for the \textit{points} of $X$. This is where the notion of a `germ' comes in. The full demonstration that every sheaf is a sheaf of cross-sections\index{section} of a suitable bundle\index{bundle} relies heavily on this idea of a \textit{germ} of a function. Earlier, with our base space discrete, we thought of \textit{germs}\index{germ} as basically elements of a set. But the more general idea of germs is that functions that agree in a neighborhood of the given ``germ point" are to be treated as equivalent. Two continuous functions $f$ and $g$ are said to have the same germ at a point $x$ provided they agree in an open neighborhood of $x$. Intuitively, this language of ``germs" at a point $x$ can be thought of as naming what data (such functions) look like under a microscrope zeroing in on $x$. Notice that 
	\begin{equation*}
	(\text{germ}_x f = \text{germ}_x g) \Rightarrow (fx = gx),
	\end{equation*} 
	i.e., 
	\begin{quote}
		if two functions have the same germ at a point, they must have the same value at that point,
	\end{quote}
	but the converse does not necessarily hold, i.e., 
	\begin{quote}
		just because two functions have the same value \textit{at} a point does not mean that they will agree \textit{near} or \textit{around} that point (i.e., are \textit{locally equal}).   
	\end{quote}  
	This notion is suggested by the following picture of the functions $y = x^2, y = |x|, y = x$, and $y = -x^2$:  
	\begin{center} 
		\includegraphics[scale=0.17]{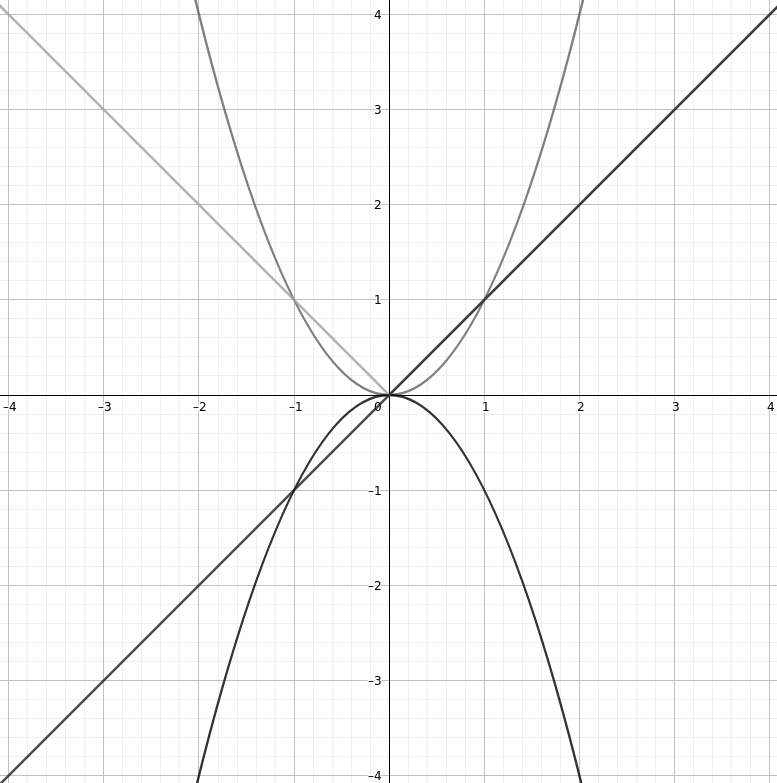}
	\end{center} 
	All of the functions above have the same value at zero, but only $y = x^2$ and $y = -x^2$ \textit{agree in a small neighborhood} around 0, so of the four, only those two will have the \textit{same germ} at zero.\index{germ} Two functions provide the same germ precisely when they become equal when we restrict down to some neighborhood of $x$. The germ of a function can accordingly be thought of as what you get when you focus your microscrope further and further in on the point $x$. Relating this to the example involving analytic continuation, two holomorphic functions $h, k: U \rightarrow \mathbb{C}$ are similarly said to have the same germ at a point $a \in U$ if their power series expansions around $a$ are the same. In other words, $h$ and $k$ agree on some neighborhood of $a$. \par 
	Applied to our present situation: take any presheaf $P: \mathscr{O}(X)^{op} \rightarrow \textbf{Set}$ on a space $X$, a point $x$, two open neighborhoods $U$ and $V$ of $x$, and two elements $s \in PU$ and $t \in PV$. Then we say that 
	\begin{definition} 
	$s$ (in $PU$) and $t$ (in $PV$) \textit{have the same germ at} $x$ when there exists some open set $W \subseteq U \cap V$ in the intersection, where $x \in W$ and $s$ and $t$ agree with respect to this set: 
	\begin{equation*} 
	s|_{W} = t|_{W} \hspace*{0.3em }\in PW.
	\end{equation*} 
	\end{definition} \noindent 
	The relation of ``having the same germ at $x$" is an equivalence relation, and we call the equivalence class of any one such $s$ the \textit{germ}\index{germ} of $s$ at $x$, denoted $\text{germ}_x s$. But this lets us define
	\begin{equation}
	P_x = \{\text{germ}_x s \hspace*{0.3em}| \hspace*{0.3em} s \in PU, x \in U, U \text{ open in } X\},
	\end{equation} 
	the set of \textit{all} germs at $x$. By taking all such functions that have the same germ (at a point) and identifying them, we just get back the stalk $P_x$ of all germs at $x$. In a sense, the very notion of a \textit{stalk}\index{stalk} $P_x$ of a sheaf $P$ is a generalization of the germ of a function, informing us about the properties of a sheaf ``near" a point $x$. Recall that the presheaf $P$ does not yield sets for the points of $X$, but rather for the open sets $U$ of $X$. We will want to narrow in on smaller and smaller open neighborhoods $U' \subseteq  U$ of a point $x$. But in looking at smaller and smaller neighborhoods $U' \subseteq  U$ of a point $x$, it will not suffice to take any single neighborhood, since a smaller one can always be taken. So we need to take some sort of limit. If $U' \subseteq  U$, we know that we have the induced presheaf restriction map $FU \rightarrow FU'$. As we range over all the neighborhoods of $x$, then, we will want to take the \textit{colimit} of all the sets $FU$. Accordingly, the stalk is also defined (over all open $U \subseteq  X$ that contain the given point $x$) as the colimit\index{colimit} (or \textit{direct limit})
	\begin{equation}
	P_x = \varinjlim_{x \in U} P(U).
	\end{equation} 
	An element of the stalk\index{stalk} will then be given by a section over a neighborhood of $x$, where two such sections will be regarded as equivalent provided their restrictions agree on a smaller neighborhood. \par 
	More explicitly, considering the restriction of the functor $P$ to the open neighborhoods of $x$, the functions $\text{germ}_x: PU \rightarrow P_x$ will form a cone as in the following diagram (since $\text{germ}_x s = \text{germ}_x (s|_{U'})$ whenever $x \in U' \subseteq U$ and $s \in PU$): 
	\begin{center} 
		\begin{tikzcd}
			& PU \arrow[ddl, "{\tau_U}", swap] \arrow[d] \arrow[ddr, "{germ_x}"] \\ 
			& PU' \arrow[dl] \arrow[dr] \\
			L && P_x \arrow[ll, "t", dashed]
		\end{tikzcd} 
	\end{center} 
	The morphism $PU \rightarrow P_x$ just takes a section $s \in PU$ defined on an open neighborhood $U$ of $x$ to its germ at $x$, generalizing the usual notion of a germ (of functions). If the functions $\{\tau_U: PU \rightarrow L\}_{x \in U}$ form another cone over $P_x$, then by definition of having the ``same germ," there will be a unique function $t: P_x \rightarrow L$ such that $t \circ \text{germ}_x = \tau$. Altogether, this just tells us that the set $P_x$ of all germs\index{germ} at $x$ is the colimit,\index{colimit} with $\text{germ}_x$ the colimiting cone, of the functor $P$ restricted to open neighborhoods of $x$. \par 
	Now, we also have that any morphism $h: P \rightarrow Q$ of presheaves, i.e., any natural transformation of functors, will induce at each point $x \in X$ a unique function $h_x: P_x \rightarrow Q_x$ such that the following diagram commutes for any open set $U$ with $x \in U$: 
	\begin{center} 
		\begin{tikzcd}
			PU \arrow[r, "{h_U}"] \arrow[d, "{germ_x}", swap] & QU \arrow[d, "{germ_x}"] \\
			P_x \arrow[r, dashed, "h_x", swap] & Q_x
		\end{tikzcd} 
	\end{center}  
	But then notice that the assignments $P \mapsto P_x$, $h \mapsto h_x$ just describe a functor $\textbf{Set}^{\mathscr{O}(X)^{op}} \rightarrow \textbf{Sets}$, a functor you can think of as ``taking the germ at $x$." \par
	Now that we have our sets $P_x$ of germs, we can range over the $x \in X$ and further combine the various sets $P_x$ of germs into the disjoint union (which we will call $\Lambda_P$) over $x \in X$: 
	\begin{equation}
	\Lambda_P = \coprod_x P_x = \{\text{all } \text{germ}_x s \hspace*{0.3em}| \hspace*{0.3em} x \in X, s \in PU\}.
	\end{equation}
	Using this $\Lambda_P$, we can define a unique function 
	\begin{equation*}
	p: \Lambda_P \rightarrow X
	\end{equation*}
	that projects each germ $\text{germ}_x s$ down the point $x$ where it is taken. With such a $p$ and the set $\Lambda_P$, we are making progress towards of description of $P$ as a bundle over $X$. However, it remains to put a topology on $\Lambda_P$, in order to speak about the continuity of $p$, and so finish the construction. \par 
	Notice that each $s \in PU$ determines a function from $U$ to $\Lambda_P$ taking all $x \in U$ to the $\text{germ}_x s$, a function that is in fact a section of $p$. In this manner, each element $s$ of the original presheaf can be replaced by an actual function to the set $\Lambda_P$ of germs. For a basis for the topology, we first take open sets around each point in $\Lambda_P$. Recall that a point in $\Lambda_P$ is just a germ, specifically a point $q \in \Lambda_P$ will lie in some set $P_x$ for some $x$, making it the germ at $x$ of some function $s \in PU$, where $U$ is an open neighborhood of the point $x$. But this $s$ will have germs at other points as well, i.e., $\text{germ}_y s$ for $y \in U$. The idea is that we use the fact that each $t \in PU$ will determine a function 
	\begin{equation*}
	\overline{t}: U \rightarrow \Lambda_P, \hspace*{2em} \overline{t}x = \text{germ}_x t, x \in U,
	\end{equation*}
	where $\overline{t}$ is itself a section of $p$.
	Our base of open sets for $\Lambda_P$ is thus got by taking all the image sets $\overline{t}(U) \subseteq  \Lambda_P$. Then an open set of $\Lambda_P$ will be a union of such sections $\overline{t}$. With such a topology, $p: \Lambda_P \rightarrow X$ (and also every function $\overline{t}$) will be continuous (in fact, a homeomorphism). Suppose $s \in PU$ and $t \in PV$, and let these determine sections $\overline{s}$ and $\overline{t}$, respectively, that agree at some point $x \in U \cap V$. Then, by the definition of a germ, we know that the set of all point $y \in U \cap V$ with $\overline{s}y = \overline{t}y$ will be an open set $W \subseteq  U \cap V$ with $\overline{s}|_W = \overline{t}|_W$. Thus, each $\overline{s}$ is continuous. \par 
	Moreover, if $h: P \rightarrow Q$ is a natural transformation between presheaves, the disjoint union of the functions $h_x: P_x \rightarrow Q_x$ yields a map $\Lambda_P \rightarrow \Lambda_Q$ of \textit{bundles}, that is moreover continuous. Therefore, altogether, we have described a functor 
	\begin{align*}
	\Lambda: \textbf{Set}^{\mathscr{O}(X)^{op}} & \rightarrow \textbf{Bn}(X) \\ 
	P & \mapsto \Lambda_P 
	\end{align*}
	from presheaves to bundles.\par 
	We have thus sketched how we can turn presheaves into bundles. In the prior section, we saw how to turn bundles into (pre)sheaves. In the next section, we come to the important take-away of all this, which involves what happens when we relate and then compose these functors. 
	\subsection{The Bundle-Presheaf Adjunction}
	\begin{theorem}
		For any space $X$, the bundle\index{bundle} functor (assigning to each presheaf $P$ the bundle of germs of $P$)
		\begin{equation*}
		\Lambda: \textbf{Set}^{\mathscr{O}(X)^{op}} \rightarrow \textbf{Bn}(X)
		\end{equation*}
		is left adjoint to the sections functor (assigning to each bundle $p: Y \rightarrow X$ the sheaf of all sections of $Y$) 
		\begin{equation*}
		\Gamma: \textbf{Bn}(X) \rightarrow \textbf{Set}^{\mathscr{O}(X)^{op}}.
		\end{equation*}
		Recall that whenever you have an adjoint pair, with left adjoint $L: \textbf{C} \rightarrow \textbf{D}$ and right adjoint $R: \textbf{D} \rightarrow \textbf{C}$, this comes with a `unit' map 
		\begin{equation*}
		\eta: \text{id} \Rightarrow RL
		\end{equation*} 
		and a `counit' 
		\begin{equation*}
		\epsilon: LR \Rightarrow \text{id}.
		\end{equation*}
		That $\Lambda \dashv \Gamma$ thus means that we will have unit and counit natural transformations  
		\begin{equation*}
		\eta_P: P \rightarrow \Gamma \Lambda P, 
		\end{equation*}
		for $P$ a presheaf, 
		and 
		\begin{equation*}
		\epsilon_Y: \Lambda \Gamma Y \rightarrow Y,
		\end{equation*}
		for $Y$ a bundle. 
	\end{theorem}
We will not give a proper proof of this here,\footnote{A proof can be found in \cite{maclane_sheaves_1994}, II.6.} but instead focus on these unit and counit maps. Consider for a given presheaf $P$ on $X$ the sheaf $\Gamma \Lambda_P$ of sections of the bundle $\Lambda_P \rightarrow X$, formed by first running the ``take the germ" functor $\Lambda$ and then following this with the ``sections" functor $\Gamma$. For each open subset $U$ of $X$, there will be a function $\eta_U: PU \rightarrow \Gamma \Lambda_P (U)$, taking $s \in PU$ to $\overline{s}$, and where restriction of $s$ to the opens of $U$ will agree with $\eta$, thus informing us that we have in fact just described a natural transformation
	\begin{equation}\label{equation: eta}
	\eta: P \rightarrow \Gamma \circ \Lambda_P.
	\end{equation}
	The next result is one of the main ``punch lines" of all this.  
	\begin{theorem}
		If the presheaf $P$ is a sheaf, then the $\eta$ given in \ref{equation: eta} will be an \textit{isomorphism} $P \cong \Gamma \Lambda_P$, i.e., \textit{every sheaf is a sheaf of sections}.
	\end{theorem} 
\begin{proof} 
	Suppose we have a presheaf $P: \mathscr{O}(X)^{op} \rightarrow \textbf{Set}$. We know that for each open subset $U$ of $X$, the function $\eta_U: PU \rightarrow \Gamma \Lambda_P (U)$, takes $s \in PU$ to $\overline{s}$, where, recall, 
	\begin{equation*}
	\overline{s}: U \rightarrow \Lambda_P, \hspace*{2em} \overline{s}x = \text{germ}_x s, x \in U.
	\end{equation*} 
	To show that $\eta_U$ is an isomorphism, let us first show that it is injective, i.e., that 
	\begin{equation*}
	\text{ if } \overline{s} = \overline{t}, \text{ then } s = t, 
	\end{equation*}
	for $s, t \in PU$. But $\overline{s} = \overline{t}$ just means that the germ of $s$ and the germ of $t$ agree on all points of $U$, i.e., $\text{germ}_x s = \text{germ}_x t$ for each point $x \in U$. Thus, for each $x$, there will moreover be an open set $V_x \subseteq  U$ such that $s|_{V_x} = t|_{V_x}$. But notice that the $V_x$ supply a cover of $U$, which means that in 
	\begin{equation*}
	PU \rightarrow \prod_x PV_x
	\end{equation*} 
	the elements $s$ and $t$ will have the same image. But then, supposing $P$ is in fact a sheaf, we know that there can be at most one such element, i.e., $s = t$, which shows $\eta_U$ injective. \par 
	To finish the proof, one must then show that $r: U \rightarrow \Lambda_P$ an arbitrary section of the bundle of germs over an open $U \subseteq  X$ is in the image of $\eta$, altogether showing that $\eta$ is an isomorphism. Details are left to the reader.\footnote{A proof can be found in \cite{maclane_sheaves_1994}, II.5.} 
\end{proof}
We can say even more about this special map $\eta$. Given a presheaf $P$, a sheaf $F$, and $\theta: P \rightarrow F$ any morphism of presheaves, then there will be a unique map $\sigma: \Gamma \Lambda_P \rightarrow F$ of sheaves making  
\begin{center}
	\begin{tikzcd}
		P \arrow[r, "\eta"] \arrow[dr, "\theta", swap] & \Gamma \Lambda_P \arrow[d, dashed, "{\sigma}"] \\ 
		& F
	\end{tikzcd}
\end{center} 
commute, i.e., $\sigma \circ \eta = \theta$. Another way of saying this is that the morphism $\eta$ is \textit{universal} from $P$ to sheaves.\footnote{See \cite{maclane_sheaves_1994}, II.5 for further details. All of this is in Grothendieck\index{Alexander Grothendieck} for the first time.} \par 
This important functor 
	\begin{equation}
	\Gamma \Lambda: \textbf{Set}^{\mathscr{O}(X)^{op}} \rightarrow \textbf{Sh}(X)
	\end{equation} 
	is known as the associated sheaf functor, or the \textit{sheafification functor}.\index{sheafification} For now, you can think of it as taking each presheaf $P$ to the ``best approximation" $\Gamma \Lambda_P$ of $P$ by a sheaf. The previous remark about the universality of $\eta$ further entails that this sheafification functor $\Gamma \Lambda$ is left adjoint to the (full subcategory) inclusion functor 
	\begin{equation*}
	\textbf{Sh}(X) \rightarrow \textbf{Set}^{\mathscr{O}(X)^{op}}. 
	\end{equation*}
	Constructing sheaves in this way as sheaves of cross-sections of a bundle suggests the further idea that a sheaf $F$ on $X$ can be replaced by the corresponding bundle $p: \Lambda F \rightarrow X$. But in constructing the topology on $\Lambda_P$ a moment ago, we mentioned that every function $\overline{t}: U \rightarrow \Lambda_P$ will not only be continuous, but will actually be a local homeomorphism. This basically means that each point of $\Lambda_P$ will have an open neighborhood that is mapped via $p$ homeomorphically onto an open subset of $X$. This leads to the following definition: 
	\begin{definition}
		A bundle $p: E \rightarrow X$ is said to be \textit{\'etal\'e (or \'etal\'e over $X$)}\index{bundle!etale} when $p$ is a \textit{local homeomorphism}, where this means that to each $e \in E$, there is an open set $V$, with $e \in V \subseteq E$, such that $p(V)$ is open in $X$ and $p|_V$ is a homeomorphism (bi-continuous isomorphism) $V \rightarrow p(V)$.
	\end{definition} \noindent 
When we have $\Lambda_P = E$, this can be imagined as a space `sitting over' $X$, where its open sets ``look like" the opens down in $X$. If we now let $\textbf{Etale}(X)$, the category of \'etal\'e bundles over $X$ denote the (full) subcategory of the category $\textbf{Bn}(X)$ of bundles over $X$, then we actually have the following powerful result:
\begin{proposition}
	The adjoint functors $\Lambda$ and $\Gamma$ from before 
			\begin{center} 
		\tikzset{
			,no line/.style={%
				,draw=none
				,commutative diagrams/every label/.append style={/tikz/auto=false}
			}
		}
		\begin{tikzcd}[column sep =large]
			\textbf{Set}^{\mathscr{O}(X)^{op}} \arrow[shift left = 2]{r}[name=U]{\Lambda} & \textbf{Bn}(X) \arrow[shift left]{l}[name=L]{\Gamma} \arrow[from=L, to=U, no line, pos=.5]{}{\perp}. 
		\end{tikzcd}
	\end{center}  
restrict, by restricting these functors to the subcategories $\textbf{Sh}(X)$ and $\textbf{Etale}(X)$, to an equivalence $\textbf{Sh}(X) \simeq \textbf{Etale}(X)$, i.e., 
	\begin{center} 
	\tikzset{
		,no line/.style={%
			,draw=none
			,commutative diagrams/every label/.append style={/tikz/auto=false}
		}
	}
	\begin{tikzcd}[column sep =large]
		\textbf{Set}^{\mathscr{O}(X)^{op}} \arrow[shift left = 2]{r}[name=U]{\Lambda} & \textbf{Bn}(X) \arrow[shift left]{l}[name=L]{\Gamma} \arrow[from=L, to=U, no line, pos=.5]{}{\perp} \\
		\textbf{Sh}(X) \arrow[u, rightarrowtail, "i"] \arrow[shift left = 2]{r}[name=S]{\Lambda_0} & \textbf{Etale}(X) \arrow[u, rightarrowtail, "i"] \arrow[shift left]{l}[name=T]{\Gamma_0} \arrow[from=T, to=S, no line, pos=.5]{}{\perp} 
	\end{tikzcd}
\end{center}  
And while the sheafification\index{sheafification} functor $\Gamma \Lambda$ is left adjoint to the inclusion of sheaves into presheaves, $\Lambda_0 \Gamma_0$ is right adjoint to the inclusion of \'etal\'e bundles into bundles. 
\end{proposition}
The idea is that just as we saw how $P$ is a sheaf precisely when $\eta_P$ is an isomorphism, it can be shown that a bundle $(Y, p, X)$ is \'etal\'e precisely when the counit morphism $\epsilon_Y$ is an isomorphism. The proof of the proposition above basically follows from general categorical reasoning regarding how the adjunction at the top restricts to consideration of subcategories. One could also use the equivalence of categories between $\textbf{Sh}(X)$ and $\textbf{Etale}(X)$ to imply the fact we discussed earlier, namely that every sheaf can be viewed as a sheaf of cross-sections.\footnote{We refer the reader to \cite{maclane_sheaves_1994}, II.5-6 for more details about the equivalence of these categories. The reader should note, however, that throughout \cite{maclane_sheaves_1994}, the authors write \'etale, when they mean \'etal\'e (the former being something else entirely).}
\end{example}
\subsection{Take-Aways}
In the last section, we described the basic adjunction involving presheaves and bundles. Just as sheaves are a special sort of ``nice" presheaf, \'etal\'e bundles are a special sort of ``nice" bundle. That there is an \textit{equivalence} of the subcategories of sheaves (of sets) on a space $X$ and the category of \'etal\'e bundles is like saying that ``to be a nice presheaf is the same thing as being a nice bundle." One of the advantages of this perspective, allowing us to regard sheaves as \'etal\'e spaces (and conversely), is that certain constructions may be simpler to define in one of the two settings, e.g., pullbacks are very easily defined in the context of local homeomorphisms, so the \textit{pullback sheaf} is more easily defined in this setting (conversely, the \textit{direct image sheaf} is simpler to define in the context of sheaves seen as a set-valued functor). \par 
The relationships explored in the previous section also allows us to take three equivalent ways of viewing morphisms between sheaves. Specifically, a morphism $h: F \rightarrow G$ of sheaves $F, G$ can be described (equivalently) in terms of: (1) just a natural transformation $h: F \rightarrow G$ of functors; (2) a continuous map $h: \Lambda F \rightarrow \Lambda G$ of bundles over $X$; and (3) as a family $h_x: F_x \rightarrow G_x$ of functions on the fibers over each $x \in X$ such that, for each open set $U$ and each $s \in FU$, the function $x \mapsto h_x(s(x))$ is a continuous $U \rightarrow \Lambda G$. The point of view given by the equivalence of (1) and (3), in particular, namely of sheaf maps $h: F \rightarrow G$ in terms of stalk maps $h_x: F_x \rightarrow G_x$, can be rather useful, for it allows many facts about sheaves to be checked ``at the level of stalks"---something that cannot be done in general for presheaves (which also reflects, conceptually, the local nature of sheaves). But the more general definition of sheaves in terms of functors satisfying certain properties is arguably superior in that it allows us to consider many non-topological cases where there is no notion of an \'etal\'e space. Thus, while it is valuable to see the important close conceptual (and historical) connection between sheaves and \'etal\'e bundles, we will usually just insist on the more general approach, where a morphism $F \rightarrow G$ of sheaves will just be a natural transformation of functors, so that with $\textbf{Sh}(X)$ we will just think of the category that has for objects all sheaves $F$ (in this case, of sets) on $X$, and for morphisms the natural transformations between them. \par
Thus far, we have confined our attention to sheaves $F$ on a topological space $X$. For such sheaves on spaces, we have been exploring basically two important candidate descriptions: 
\begin{enumerate} 
	\item the \textit{restriction-collation} description, and 
	\item the \textit{section} description.
\end{enumerate} 
	The first description was motivated by structures, such as classes of functions with certain ``nice" properties (like continuity), that are defined ``locally" on a space. The previous example introduced and developed the importance of the second of these two perspectives. The idea here was that we took a sheaf to be some sort of principled way of assigning to each point $x$ of the underlying space a set $F_x$ consisting of all the ``germs" at $x$ of the functions being considered (where these ``germs" are equivalence classes identifying sets that look ``locally the same" in neighborhoods of $x$), after which these sets $F_x$ then get patched together by a topology to form a space (or bundle) projected onto $X$ (a suitable function for this sheaf then being a ``cross section" of the projection of this bundle). According to this perspective, the sheaf $F$ can ultimately be thought of as a set $F_x$ that varies with the points $x$ of the space. \par 
	After learning a new concept and seeing some descriptions and preliminary examples, it is important to think a bit about when and how this can go wrong (and acquire a store of ``non-examples"). Before ending this chapter and introducing more involved examples of sheaves (as well as some more advanced aspects of sheaves), we pause to consider what is \textit{not} a sheaf, i.e., when and why a construction fails to satisfy the sheaf conditions.  After that, this chapter will end with a brief but important discussion of a general result allowing us to blur the distinction between presheaves and sheaves in the special case of posets. 
\subsection{What is \textit{Not} a Sheaf}
Even when structures are determined locally, sometimes local properties alone do not suffice to determine global properties. In such cases, we will not have a sheaf. A common example given to illustrate this is the set $BU$ of all \textit{bounded} functions on $U$ to $\mathbb{R}$---this will give a functor on $U$, but not a sheaf. The reason for this is that while the collation of functions that are bounded does indeed define a unique function on $U$, such a function may be unbounded. \par
\begin{exercise}
	Make sense of the previous sentence by giving an example of a collection of bounded functions on subintervals whose collation is not bounded.
\end{exercise}
If a structure is not even determined locally, specifically in the sense that it does not even obey the first (restriction) condition, then it certainly cannot be a sheaf. An intuitive example of this might be given by the game of Scrabble$^{TM}$,\index{Scrabble} where one thinks of this as follows: the $15 \times 15$ board with its squares labeled in some sensible way (with $x$ and $y$ ``coordinates" of a tessellation, so that $(1,1)$ would indicate the leftmost top corner square), may be regarded as a topological space, with a notion of covers. Then one might attempt to regard the assignment of sets of legal (English) word-forming letter combinations to subsets of the grid of squares (satisfying a further constraint capturing how words are to be ``read" down and to the right) as a functor. However, while each inclusion of ``opens" in the underlying grid of squares would have to induce a function restricting the word-forming letters assigned over a bigger region of the board to the word-forming letter assignments over a sub-region, in general not every sub-word of a word is a word, so it is not clear how to make this work. Even if we agreed to treat individual letters as (legal) words, it is evident that the inclusions of open subsets will sometimes determine a ``restriction" to a particular part of a word that does indeed form a word (now a \textit{different} word), but on other occasions such a process will not result in a word at all. For instance, confining our attention for simplicity to a small $3 \times 3$ region of the board, and displaying a portion of the ``opens" ordered by inclusion
\begin{center}
	\includegraphics*[scale=0.23]{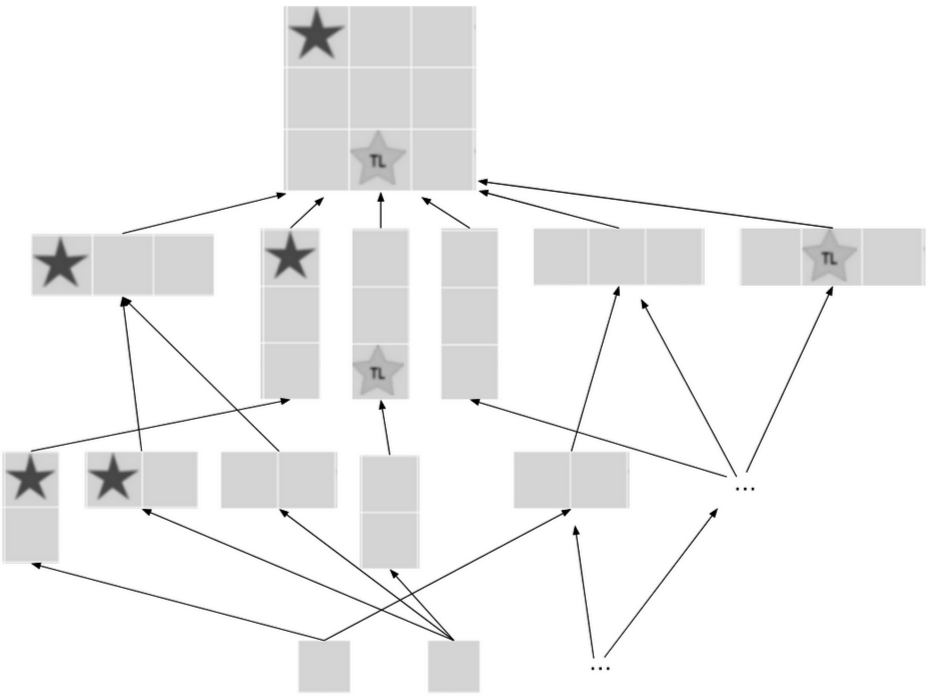}
\end{center}
we might then regard a particular selection of possible letter assignments as follows 
 \begin{center}
 	\includegraphics*[scale=0.23]{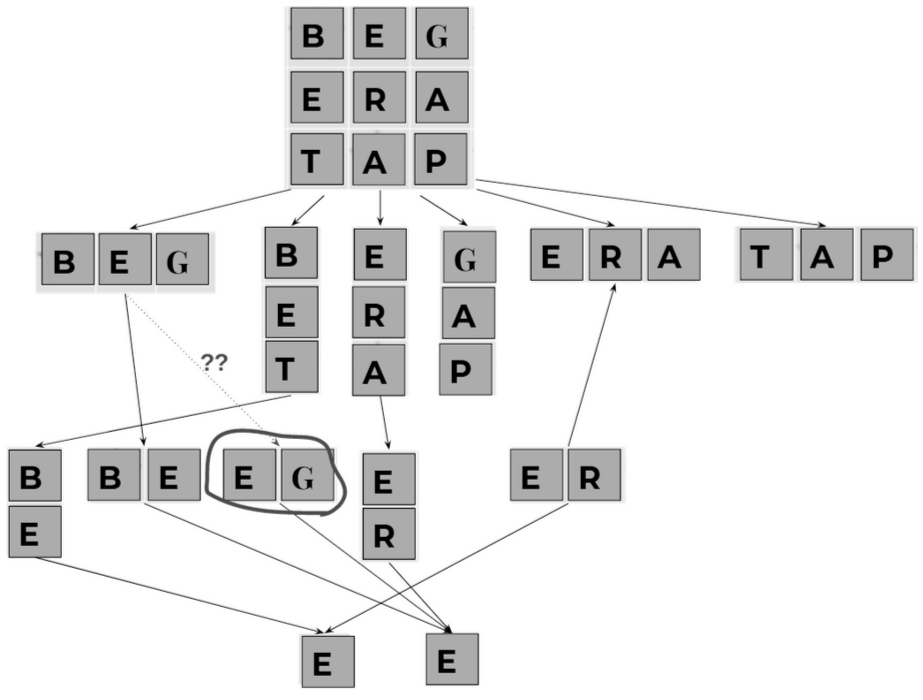}
 \end{center}
Here, with such an assignment, even though the letter assignments over the entire $3 \times 3$ portion result in valid words (in every possible 3-letter combination, e.g., `beg', `bet', `era', `gap', `tap'), some of whose parts even themselves form words (like `be' in `beg'), it is not clear what to do with the (failed) ``restriction" from the word `beg' down to `eg', which is not a word. 
\par    
Even when we do have a presheaf, in general a presheaf can itself fail to be a sheaf in two (fundamentally independent) ways:  
\begin{itemize}
	\item \textbf{Non-locality}: If a presheaf has a section $s \in F(U)$ that cannot be constructed from sections over smaller open sets in $U$---via a cover, for instance---then $F$ fails to be a sheaf. 
	\item \textbf{Inconsistency}: If a presheaf has a pair of sections $s \neq t \in F(U)$ such that when restricted to every smaller open set they define the same section, then $F$ fails to be a sheaf. In other words, informally, the presheaf has local sections that ``ought to" patch together to give a unique global section, but do not.  
\end{itemize}
It is often thought that the second sort of failure is somehow easier to understand. But that does not mean that there are not examples of the first sort of failure. A standard illustration of non-locality is given by the following example. Consider $X$ the topological space consisting of two points $p, q$, endowed with the discrete topology (i.e., every set is open). Then $X$ consists of the open sets $\{p,q\}, \{p\}, \{q\}, \emptyset$, ordered by inclusion. We can form the \textit{constant presheaf} $P$ on $X$ which assigns a set (or abelian group) to each of the four open sets and the identity map to each of the nine restriction maps (five plus the four trivial self-maps). For concreteness, let this presheaf assign $\mathbb{Z}$ to each of the sets.     
\begin{center} 
\begin{tikzpicture}[yscale=0.7, xscale=0.9]
\node (max) at (0,4) {$\{p,	q\}$};
\node (a) at (-2.4,2) {$\{p\}$};
\node (c) at (2.4,2) {$\{q\}$};
\node (e) at (0,0) {$\emptyset$};

\draw[->] (a) -- (max);
\draw[->] (c) -- (max);
\draw[->] (e) -- (a);
\draw[->] (e) -- (c);
\draw[->] (e) -- (max);

\node (max1) at (8,4) {$P(\{p,q\}) = \mathbb{Z}$};
\node (a1) at (5.5,2) {$P(\{p\}) = \mathbb{Z}$};
\node (c1) at (10.5,2) {$P(\{q\}) = \mathbb{Z}$};
\node (e1) at (8,0) {$P(\emptyset) = \mathbb{Z}$};

\path[<-, draw] (a1) edge node[above, font=\footnotesize] {$id$} (max1);
\path[<-, draw] (c1) edge node[above, font=\footnotesize] {$id$} (max1);
\path[<-, draw] (e1) edge node[left, font=\footnotesize] {$id$} (a1);
\path[<-, draw] (e1) edge node[right, font=\footnotesize] {$id$} (c1);
\path[<-, draw] (e1) edge node[right, font=\footnotesize] {$id$} (max1);
\node (dom) at (0,-2) {$\mathscr{O}(X)^{op}$};
\node (codom) at (8,-2) {\textbf{Set}};
\draw[->] (2,-2) to node[pos=0.5, above, font=\large] (6,-2) {$P$} (6,-2);
\end{tikzpicture} 
\end{center} 
This presheaf $P$ does indeed satisfy the gluing axiom. However, it fails to satisfy the locality/identity axiom, specifically with respect to the assignment on the empty set. The empty set is covered by the empty family of sets; but clearly any two local sections of $P$ are equal when restricted to their common intersection in the empty family. If the locality axiom were satisfied, then any two sections of $P$ over the empty set will be equal---however, this need not be true.\par  
Another example of non-locality is given by the following: take an open set $X \subseteq \mathbb{C}^n$, and for open $U \subseteq X$, let
\begin{equation}
S(U) := \{f: U \rightarrow \mathbb{C} \hspace*{0.3em}| \hspace*{0.3em} f \text{ is holomorphic} \}.
\end{equation}      
We then define the restriction maps by stipulating that for $V \subsetneqq U$, we set the restriction $\rho_{UV} = 0$, and set $\rho_{UU} = \text{id}$. This $S$ is a presheaf but it is not a sheaf precisely because it obviously has a section in $S(U)$ that cannot be built from sections over smaller open sets in $U$ (for which $f$ must be the $0$ map).\par 
To illustrate the second type of failure, inconsistency, consider again the constant presheaf $P$ from before. We construct a new presheaf $G$ over the same $X$ with the same discrete topology, which is just like $P$ except that we now let $G(\emptyset) = \{*\}$, where $\{*\}$ is a one-element set (the terminal object). We retain $\mathbb{Z}$ as our value assignment for the remaining non-empty sets. Now, however, for each inclusion of opens that has the empty set for domain, $G$ will assign the unique map $0$; otherwise, it just assigns the identity map as before.
\begin{center} 
\begin{tikzpicture}[yscale=0.7, xscale=0.9]
\node (max) at (0,4) {$\{p,	q\}$};
\node (a) at (-2.4,2) {$\{p\}$};
\node (c) at (2.4,2) {$\{q\}$};
\node (e) at (0,0) {$\emptyset$};

\draw[->] (a) -- (max);
\draw[->] (c) -- (max);
\draw[->] (e) -- (a);
\draw[->] (e) -- (c);
\draw[->] (e) -- (max);

\node (max1) at (8,4) {$G(\{p,q\}) = \mathbb{Z}$};
\node (a1) at (5.5,2) {$G(\{p\}) = \mathbb{Z}$};
\node (c1) at (10.5,2) {$G(\{q\}) = \mathbb{Z}$};
\node (e1) at (8,0) {$G(\emptyset) = \{*\}$};

\path[<-, draw] (a1) edge node[above, font=\footnotesize] {$id$} (max1);
\path[<-, draw] (c1) edge node[above, font=\footnotesize] {$id$} (max1);
\path[<-, draw] (e1) edge node[left, font=\footnotesize] {$0$} (a1);
\path[<-, draw] (e1) edge node[right, font=\footnotesize] {$0$} (c1);
\path[<-, draw] (e1) edge node[right, font=\footnotesize] {$0$} (max1);
\node (dom) at (0,-2) {$\mathscr{O}(X)^{op}$};
\node (codom) at (8,-2) {\textbf{Set}};
\draw[->] (2,-2) to node[pos=0.5, above, font=\large] (6,-2) {$G$} (6,-2);
\end{tikzpicture} 
\end{center} 
This presheaf $G$ satisfies the locality/identity axiom, but now it fails to be a sheaf on account of not satisfying the gluing axiom. The entire set $X = \{p,q\}$ is covered by $\{p\}$ and $\{q\}$, which individual sets obviously have empty intersection. By definition, the sections on $\{p\}$ and $\{q\}$ will just be an element of $\mathbb{Z}$, i.e., an integer. By selecting a section $m$ (an integer) for our section over $\{p\}$ and $n$ (another integer) for our section over $\{q\}$, such that $m \neq n$, we can easily see this violation. $m$ and $n$ must restrict to the same element over $\emptyset$ on account of the action of the trivial $0$ restriction map; but then, if the gluing axiom were satisfied, because $m$ and $n$ restrict to the same element over their (trivial) intersection, we would need the existence of a unique section $s$ over the union of the two sets, i.e., in $G(\{p,q\})$, which moreover restricts back to $m$ on $\{p\}$ and to $n$ on $\{q\}$. But the restriction maps from $G(\{p,q\})$ along $\{p\}$ and $\{q\}$, being the identity map in both cases, forces that $s = m$ and $s = n$, from which $m = n$, contradicting the assumption that the sections (integers) $m$ and $n$ were different. \par      
Before moving on, we take the opportunity to briefly mention a common issue that may arise in the construction of sheaves in practice (construed as sheaves of sections), but one that is somewhat distinct from the failures of the previous two sorts. It concerns a situation where we may in fact be dealing with a sheaf, but \textit{given} certain selections from the stalks or local sections, we may find that we simply cannot extend those assignments to produce a global section, since no matter what we assign to the remaining open set(s), we will run into inconsistency with respect to the other sections. This may simply be a problem with failing to select the ``right" elements from the sets of possible values assigned to each underlying open. In the next chapter, we will see a number of explicit instances of this.\par 
Finally, let us briefly look at a particularly interesting example (due to \cite{goguen_sheaf_1992}) of a presheaf that is not a sheaf, for reasons distinct from the two discussed above.  
\begin{example}
Taking our indexing (domain) category to be some \textit{base} for some data assignment or observations, where this base is a poset, then the particular base consisting of intervals of natural numbers beginning with an ``initial time" 0---where the various intervals beginning from 0 may represent periods of continuous observation of the system---can be described as
\begin{equation*}
\mathscr{I}_0 (\omega) = \{\emptyset, \{0\}, \{0,1\}, \{0,1,2\}, \dots \} \cup \{\omega \},
\end{equation*} 
with $\omega$ representing the domain for observations or data assignment over an infinite time. Using such a base, we might consider a ``fair scheduler" $F$ for events $a,b$, i.e., if $a$ occurs, then at some point $b$ must occur, and vice versa. More formally, this means 
\begin{equation*}
F(\omega) = (a^{+}b^{+} + b^{+}a^{+})^{\omega},
\end{equation*}
 i.e., we have a concatenation of strings where the components consist of some number of $a$'s followed by some number of $b$'s or some number of $b$'s followed by some number of $a$'s, on to infinity. For each $n$ in the set of natural numbers $\omega$, we will have that $F(\{0,1,\dots,n-1 \}) = \{a, b\}^n$, which is just to say that we might have any combination of $a$'s and $b$'s (including all $a$'s or all $b$'s) for some finite interval. The point, however, is that while $F$ is indeed a presheaf, if we had that $F$ was a sheaf, then the sheaf condition will imply that $F(\omega) = \{a,b \}^{\omega}$, i.e., that in the limit any combination is also possible.  However, this contradicts the definition of $F$ as a fair scheduler. Thus, it is not a sheaf. Though, if we require the indexing set to be finite, then this $F$ does satisfy the finite sheaf and gluing conditions. This (non)example is interesting because it suggests that interesting phenomena can appear ``at infinity" that do not show up in the finite approximations. 
\end{example} 
Returning to more general considerations, as the section on bundles already suggested, there is a standard procedure for completing a presheaf to make it a sheaf. Since there are two fundamental ways a presheaf can fail to be a sheaf, this process, usually dubbed ``sheafification," can be roughly thought of as doing one of two things: (1) it \textit{discards} those extra sections that make the presheaf fail to satisfy the locality condition; (2) it \textit{adds} those missing sections which, had they been present, would allow the local sections to glue together into a unique global section, satisfying the gluability condition.  In other words, with respect to the second of these two, we are adding functions to the global set that restrict to compatible functions on each of the opens, and then, recursively, we continue adding the restrictions of the newly-generated global functions. We have already seen what this sheafification abstractly looks like in the setting of bundles. In Chapter 5, we will look more closely at this process of sheafification whereby an arbitrary presheaf can be turned into a sheaf of the same type. \par 
For now, let us just sketch how the sheaf $G$ from a moment ago would be ``sheafified". Basically, this involves expanding $G(\{p,q\})$ to $\mathbb{Z} \oplus \mathbb{Z}$, thus defining a new sheaf $H$, then letting the restriction maps be the appropriate projection maps $\pi_i: \mathbb{Z} \oplus \mathbb{Z} \rightarrow \mathbb{Z}$, thereby defining $H(\{p\}) = image(\pi_1) = \mathbb{Z}$ and $H(\{q\}) = image(\pi_2) = \mathbb{Z}$. Everything else can be defined just as it was for $G$. The resulting functor $H$ now satisfies the gluability condition and so is a sheaf, usually called the \textit{constant sheaf} on $X$ valued in $\mathbb{Z}$.  
\subsection{Presheaves and Sheaves in Order Theory}
As is so often the case when working with categories, things are greatly simplified when dealing with orders. As we saw in the first chapter, highly abstract results and concepts in category theory can have a particularly friendly showing when specialized to orders. Sometimes, posets are especially ``nice" to us, in that important general distinctions (such as that between presheaves and sheaves) can be ``collapsed" when dealing with posets. This in turn can sometimes make the more general distinction easier to grasp. \par 
Before concluding this chapter, we will briefly cover a highly useful result that relates presheaves and sheaves on a poset (regarded as a category):  
\begin{quote}
	\textit{Presheaves on a poset are} equivalent \textit{to sheaves over that poset, once the latter has been equipped with a suitable topology (the ``Alexandrov topology").}\index{topology!Alexandrov} 
\end{quote}
We will sketch how this works. First recall from \ref{definition: downset} the notion of a \textit{downset}\index{downset} (and its dual, an \textit{upper set}).\index{upper set} Recall also how we defined \textit{principal downsets}, denoted $\mathcal{D}_p$ (or just $\downarrow p$): these were sets of the form 
\begin{equation*}
\downarrow p := \{q \in \mathcal{P}: q \leq p \},  
\end{equation*}
for $p \in \mathcal{P}$ (and dually for the principal upper sets). It is fairly straightforward to show how the principal downsets can generate a topology, called the Alexandrov topology (or, sometimes, \textit{lower Alexandrov topology}, to distinguish it from the topology generated by the principal upper sets).\footnote{This can be described as a functor from $\textbf{PreOrd}$ to $\textbf{Top}$.} While the poset just supplies us with points $p \in \mathcal{P}$, the topology generated by the principal down (upper) sets supplies us with a way of looking at points now in terms of opens (the language a topological sheaf will understand).  Recall that any downset can be written as the union of principal downsets\index{downset!principal} (taking unions of downsets is the same as taking colimits of representables in the poset category of downsets), so using the principal down (upper) sets as a basis, we can form the collection of all downsets of a poset $\mathcal{P}$, denoted by $\mathcal{D}(\mathcal{P})$, and this will define a topology on $\mathcal{P}$, where we take $\mathcal{D}(\mathcal{P})$ as our open sets $\mathscr{O}(\mathcal{P})$. As $\mathcal{D}(\mathcal{P})$ is closed under arbitrary intersections, the ``closed sets" (i.e., the upper sets of $\mathcal{P}$) will be closed under arbitrary unions. But this means that the upper sets also forms a topology, usually called the \textit{upper Alexandrov topology}. Thus, dually, denoting by $\mathcal{U}(\mathcal{P})$ the collection of all upper sets of $\mathcal{P}$, this yields another Alexander topology, the upper Alexandrov topology on $\mathcal{P}$. \par     
With these notions in hand, it can be shown that for $\mathcal{P}$ a poset, presheaves on $\mathcal{P}$ are the same as sheaves on $\mathcal{D}(\mathcal{P})$, while (dually) copresheaves (variable sets) on $\mathcal{P}$ are the same as sheaves on $\mathcal{U}(\mathcal{P})$. 
\begin{theorem}
	\begin{equation*}
	\textbf{Set}^{\mathcal{P}^{op}} \simeq \textbf{Sh}(\mathcal{D}(\mathcal{P})),
	\end{equation*}
	and 
	\begin{equation*}
	\textbf{Set}^{\mathcal{P}} \simeq \textbf{Sh}(\mathcal{U}(\mathcal{P})).
	\end{equation*}
\end{theorem} \noindent 
In other words, this tells us (in the first case) that, given a presheaf defined on a poset $\mathcal{P}$, this can be regarded as a sheaf when $\mathcal{P}$ is equipped with the natural Alexandrov topology (induced by the downset completion construction). Instead of going through a full proof of this, we will sketch one way of looking at why this is true. Because of the example to follow, we will also focus on $\textbf{Set}^{\mathcal{P}} \simeq \textbf{Sh}(\mathcal{U}(\mathcal{P}))$.\par 
Recall, from the first chapter, the very close relationship between $\downarrow p$ and $p$, via the downset embedding (which specialized the Yoneda embedding); and we just saw that the principal downsets of $\mathcal{P}$ can be shown to form a topology. There is a similar result for upper sets, where $p \leq p'$ iff $\uparrow p' \subseteq  \uparrow p$ (note the reversal of order). Likewise, we just mentioned that via such a construction, there is the natural upper Alexandrov topology on $\mathcal{P}$, the basis of which is given by the principal upper sets $\uparrow p$ (note that this is the \textit{smallest} open set that will contain $p$).\par  
Suppose you have a functor (copresheaf, or variable set) $F^* \in \textbf{Set}^{\mathcal{P}}$. (You can of course also think of this as a presheaf on $\mathcal{P}^{op}$.) Let $p$ be an element of $\mathcal{P}$. We can take this functor to a sheaf $F: \mathscr{O}(\mathcal{P})^{op} \rightarrow \textbf{Set}$, where $\mathscr{O}(\mathcal{P}) = \mathcal{U}(\mathcal{P})$, by defining $F$ on a basis of the Alexandrov topology by taking 
\begin{equation*}
F(\uparrow p) = F^*_p,
\end{equation*}
where $F^*_p$ is just the image of $F^*(p)$, and where of course $\uparrow p \in \mathcal{U}(\mathcal{P})$. Another way of thinking of this is that, given a functor $F^*: \mathcal{P} \rightarrow \textbf{Set}$, via the inclusion functor 
\begin{align*}
\iota: & \mathcal{P} \rightarrow \mathcal{U}(\mathcal{P})^{op} \\ 
&  p \mapsto \uparrow p
\end{align*} 
we want to know whether we can produce a sheaf 
\begin{center} 
	\begin{tikzcd}
		\mathcal{P} \arrow[dr, "\iota", swap] \arrow[rr, "F^*"] & &  \textbf{Set} \\ 
		& \mathscr{U}(\mathcal{P})^{op} \arrow[ur, dashed, "{?}", swap]
	\end{tikzcd}
\end{center} 
The most concise categorical way of accomplishing all this would be to use what is called Kan extensions, specifically the right Kan extension of $F^*$ along $\iota$, denoted $\text{Ran}_{\iota} F^*$, to assign data to the opens in our poset ``in a nice way" such that it is a sheaf, and then declare that when this happens, then $F^*$ itself can be seen as a sheaf. But instead, we will just describe things in more elementary terms. The idea is that, having defined $F(\uparrow p) = F^*_p$, we extend this to the general opens of $\mathcal{U}(\mathcal{P})$ by recognizing that for each $V \in \mathcal{U}(\mathcal{P})$, the set $\{\uparrow p: p \in V\}$ will cover $V$, and we have 
\begin{equation*}
F(V) = \varprojlim_{p \in V} F(\downarrow p) = \varprojlim_{p \in V} F^*_p. 
\end{equation*}
Passing from a sheaf $F$ to the corresponding presheaf $F^*$ on $\mathcal{P}$ is straightforward, using the same identification $F^*_p = F(\uparrow p)$. \par 
Note that whenever $p \leq p'$, we know that $\uparrow p' \subseteq  \uparrow p$, so we will have the maps 
\begin{equation*}
F^*_{p, p'}: F(\uparrow p) \rightarrow F(\uparrow p'), 
\end{equation*} where this is the image of the inclusion $\uparrow p' \hookrightarrow \uparrow p$ under the (contravariant) sheaf functor $F$, i.e., maps $F^*_{p'} \rightarrow F^*_{p}$. Note how the inclusion $\iota$ reverses the order, and then the underlying presheaf action of the sheaf $F$ reverses order once more. \par 
Altogether, this result informs us that we can pass freely from a sheaf over the natural topology $\mathcal{U}(\mathcal{P})$ induced on $\mathcal{P}$ and a ``plain" copresheaf (variable set) on that $\mathcal{P}$. By taking $\mathcal{P}^{op}$ instead, and recognizing that the upper sets of $\mathcal{P}^{op}$ are the same as the down sets of $\mathcal{P}$, we also get that we can move between a sheaf over the natural lower Alexandrov topology given by $\mathcal{D}(\mathcal{P})$ and a ``plain" presheaf on $\mathcal{P}$. 
Of course, in general cases, with categories that are not posets, presheaves are not automatically going to yield a sheaf. In making a presheaf a sheaf in the general case, we are in a sense demanding that it ``awaken to" the features of the underlying topology; specifically, in trying to find the sheaves, we want to restrict attention to those presheaves that are ``sensitive" to the structure supplied by the cover. The above result letting us blur the distinction between presheaves and sheaves, when $\mathcal{P}$ is a poset now equipped with its natural Alexandrov topology, in a sense tells us that presheaves in this setting are automatically ``sensitive" to the structure of covers. A presheaf on a poset can already be regarded as a sheaf (with respect to its Alexandrov topology).\par 
The next chapter will make good use of this result. This chapter ends with a brief look at a different sort of poset where this result is of some utility.  
\begin{example}
	We start by defining a \textit{(time)frame} as a structure $\mathcal{T} = (T, \leq)$ consisting of a non-empty set $T$ of ``times" (instants, events) on which the relation $\leq$ forms a reflexive and transitive order, i.e., $\mathcal{T}$ is a pre-order. We have not insisted that a (time)frame be a partial order, i.e., that $\leq$ be antisymmetric, so the equivalence relation defined on $T$ by $t \approx s$ iff $t \leq s$ and $s \leq t$ will be non-trivial in general. We call the $\approx$-equivalence classes the \textit{clusters} of $\mathcal{T}$, ordered by setting $\hat{t} \leq \hat{s}$ iff $t \leq s$, where $\hat{t}$ is the cluster containing $t$. Such an order is an antisymmetric order, enabling us to regard the (time)frame as a poset of clusters. Finally, a frame is then said to be \textit{directed} provided any two elements have an upper bound, i.e., for all $t, s \in T$, $\exists v \in T$ such that $t \leq v$ and $s \leq v$. \par 
	Starting in \ref{section: modalities}, we discussed some features of modalities.\index{modalities} Propositional modal logic consists of sentences constructed from sentence letters $p, q, r, \dots$, Boolean connectives, and the modal operator $\Box$ (which is typically interpreted, in the context of tense logics, as `it will always be'). The Greek philosopher Diodorus of Megara\index{Diodorus} is often thought to have held that the modalities `necessity' and `possibility' were definable in terms of time. For instance, Diodorus held that the necessary should be understood as that which \textit{is} (now) and \textit{will always be} the case. One might then define $\Diamond$ (interpreted as `it will (at some time) be') as $\neg \Box \neg$. In terms of contemporary modal logics, while it is common for tense logics to regard time as an irreflexive ordering (so that ``at all future times" does not include the present moment), Diodorus's approach suggests a temporal interpretation of $\Box$ that uses reflexive orderings instead, and we can apply this to non-linear time structures. The reflexivity of $\leq$ in our time frame will give $\Box$ the Diodorean intepretation of `\textit{is} (now) and \textit{always will be}'.  \par 
	As \cite{goldblatt_diodorean_1980} first showed, one can develop a Diodorean logic of $n$-dimensional (for $n \geq 2$) Minkowskian special-relativistic spacetime, and show that this is exactly the modal logic \textbf{S4.2}. In more detail, if $x = \langle x_1, \dots, x_n \rangle$ is an $n$-tuple of real numbers, then let 
	\begin{equation*}
	\mu(x) = x_1^2 + \cdots x_{n-1}^2 - x_n^2.
	\end{equation*}
	Then, for $n \geq 2$, $n$-dimensional spacetime is the frame 
	\begin{equation*}
	\mathbb{T}^n = (\mathbb{R}^n, \leq),
	\end{equation*} 
	with $\mathbb{R}^n$ the set of all real $n$-tuples, and the order $\leq$ defined, for all $x$ and $y$ in $\mathbb{R}^n$, as follows: 
	\begin{equation*}
	\begin{split} 
	x \leq y & \text{ iff } \mu(y-x) \leq 0 \text{ and } x_n \leq y_n \\
	& \text{ iff} \sum_{i=1}^{n-1} (y_i - x_i)^2 \leq (y_n - x_n)^2 \text{ and } x_n \leq y_n. 
	\end{split}  
	\end{equation*}
	The frame $\mathbb{T}^n$ is partially-ordered and directed. 
	Then the usual Minkowski spacetime of special relativity is given by $\mathbb{T}^4$, for which a point will represent a spatial location $\langle x_1, x_2, x_3 \rangle$ at time $x_4$. In this case, the interpretation of $x \leq y$ is that a signal may be sent from event $x$ to event $y$ at a speed at most that of the speed of light, entailing that $y$ ``comes after," or is in the ``causal future" of, $x$. In other words, the (reflexive) relation is given by `can reach with a lightspeed-or-slower signal'. For simplicity, we can look at the frame $\mathbb{T}^2$, and easily visualize the future cone $\{z: x \leq z \}$ for a point $x = \langle x_1, x_2 \rangle$, where the future cone contains all points on or above the directed rays of slopes $\pm 1$ beginning from $x$ (where we assume a coordinate system for which speed of light is one unit of distance per unit of time). 
	\begin{center}
		\includegraphics*[scale =0.4]{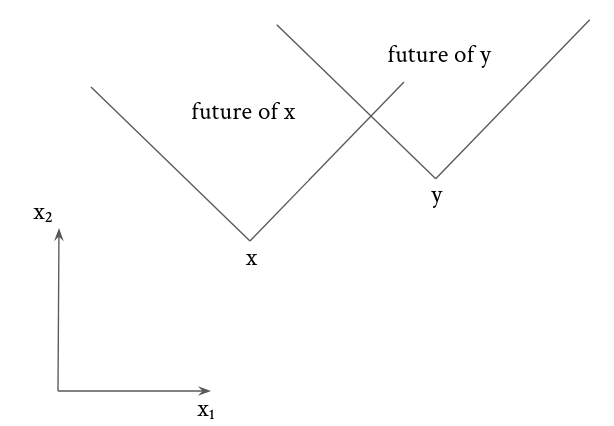}
	\end{center} \par \noindent 
	One can observe that the future cones of any two points will eventually intersect, which entails that the underlying order is \textit{directed}, in the sense that for any two locations $x, y$, there is a third that is in the future of both $x$ and $y$. This directedness forces the Diodorean interpretation of $\Box$ to validate the $\textbf{S4.2}$ axiom schema $\Diamond \Box A \rightarrow \Box \Diamond A$, where $\textbf{S4.2}$ just arises by adding that axiom schema to the usual axioms of $\textbf{S4}$. \cite{goldblatt_diodorean_1980} showed that, in fact, \textit{each} of the frames $\mathbb{T}^n$ has the logic $\textbf{S4.2}$ for its Diodorean modal logic, regardless of spatial dimensions.\par 
	Now, following Goldblatt, we may also call $T' \subseteq  T$ \textit{future-closed} under $\leq$ provided whenever $t \in T'$ and $t \leq s$, then also $s \in T'$.\footnote{Note that this is just to say that $T'$ is an upper set!\index{upper set}} In this case, $\mathcal{T}' = (T', \leq)$ will be a \textit{subframe} of $\mathcal{T}$, and by the transitivity of $\leq$, for each $t$ the set $\{s: t \leq s\}$ will be the \textit{base} of a subframe. We can denote the collection of future-closed subsets of $\mathcal{T}$ by $\mathcal{T}^+$.\footnote{We could also define, as usual, a $\mathcal{T}$-\textit{valuation} as a function $V: \Phi \rightarrow \mathcal{T}^+$, where $\Phi$ is the set of atomic formulae, with typical member $p$. A valuation sends each sentence letter $p$ to a future-closed subset $V(p) \subseteq  T$, interpreted as the set of times at which $p$ is ``true." We also have that $t \in V(\Box A)$ iff $t \leq s$ implies $s \in V(A)$, allowing the valuation to be extended to all sentences (using the Boolean connectives as well). A \textit{model} based on $\mathcal{T}$ is then defined as a pair $\mathcal{M} = (\mathcal{T}, V)$, where $V$ is a $\mathcal{T}$-valuation. If $\mathcal{T}'$ is a subframe of $\mathcal{T}$, then for any sentence $A$, $A$ is valid on the frame $\mathcal{T}$ (i.e., $A$ is true in every model based on $\mathcal{T}$) only if $A$ is valid on the frame $\mathcal{T}'$.} 
	But since the ordered collection $\mathcal{T}^+$ of all future-closed subsets of $\mathcal{T}$ constitutes an (Alexandrov) topology on $T$, we can use the exact correspondence 
	\begin{equation*}
	\textbf{Set}^{\mathcal{P}} \simeq \textbf{Sh}(\mathcal{U}(\mathcal{P}))
	\end{equation*}
	 to move freely between sheaves on $\mathcal{T}^+$ and the usual ``variable set" perspective of $\textbf{Set}^{\mathcal{T}}$. 
\end{example}
\chapter{Sheaf Cohomology through Examples} 
In this chapter, we start to look at some more involved and computationally-explicit examples, working up towards an extended introduction to sheaf cohomology, presented via a particularly computational example. Roughly, if sheaves represent local data---or, more precisely, represent how to properly ensure that what is locally the case everywhere is in fact globally the case---sheaf cohomology\index{sheaf!cohomology} can be thought of as a tool for systematically exploring, representing, and relating \textit{obstructions} to such passages from the local to the global. Moreover, in sheaf theory more generally, one could argue (as does Grothendieck,\index{Alexander Grothendieck} for instance) that individual sheaves are only of secondary importance---the real power of sheaf theory emerges from the use of constructions involving various sheaves, linked together via sheaf morphisms. Sheaf cohomology will allow us to begin to appreciate such a perspective. 
\section{Simplices and their Sheaves}
Of the many ways to represent a topological space, a particularly computationally-friendly way is to perform a triangulation with entities called \textit{simplices},\index{simplices} decomposing the space into simple pieces (thought of as being `glued together') whose common intersections or boundaries are lower-dimensional pieces of the same kind. With simplices come certain simplicial maps that, moreover, approximate continuous maps. In this way, simplices play a role in bridging the gap between continuous figures and their discrete representation and approximation via decompositions of spaces into discrete parts. More than that, as we will see, they allow us to develop profound connections between algebra and geometry. Simplices are a powerful and easy-to-use device for understanding qualitative features of data collections, and in general they can be thought to represent $n$-ary relations between $n$ vertices. \par 
	Basically, we use collections of simplices---for now just think of points, line segments, generalized triangles or tetrahedra generalized to arbitrary dimensions---to build up what are called \textit{simplicial complexes}.\index{complex!simplicial} A (geometrical) simplicial complex $K$ is a collection of simplices such that (i) every face of a simplex of $K$ is in $K$, and (ii) the intersection of any two simplices of $K$ is a face of each of them. Roughly, then, one can think of a simplicial complex $K$ as comprised of generalized triangles of various dimensions, glued together along common faces. This is really part of a more general story involving \textit{cell complexes}\index{complex!cell} (including cubical complexes, multigraphs, etc.), where one can roughly think of a cell complex as a collection of \textit{closed disks} of various dimensions which are moreover glued together along their boundaries. But we will instead focus on the more computationally-tractable combinatorial counterpart to the already simplified notion of a simplex: that of \textit{abstract simplicial complexes}. Here, we re-encode the information of a simplicial complex via the more computationally-friendly notion of an abstract simplicial complex (or ASC), where this is basically just a collection of subsets of ``vertices" (elements), closed under the operation of taking subsets. This captures in a purely combinatorial way the geometrical notion of simplicial complex.\footnote{What has been lost is how the simplex is embedded in, say, Euclidean space; however, this specification retains all the data needed to reconstruct the complex up to homeomorphism.}   
	\begin{definition}
		An \textit{abstract simplicial complex} (ASC)\index{complex!abstract simplicial} $K$ on a set $A$ is a collection of ordered finite non-empty subsets $K \subseteq  \mathbb{P}(A)$ that is closed under taking subsets (sublists), i.e., every subset of a set in $K$ is also in $K$. In other words, we must have 
		\begin{itemize}
			\item for each $x \in A$, the singleton $\{x\} \in K$; and 
			\item if $\sigma \in K$ and $\tau \subseteq  \sigma$, then $\tau \in K$. 
		\end{itemize}\par \noindent 
		Terminologically, each member of $K$ is called a \textit{simplex} or a \textit{face} (or sometimes a \textit{cell}). A face with $n+1$ elements is called an $n$-dimensional face (or a $n$-face, or $n$-simplex) of $K$.\footnote{But, as one would expect, a $0$-face is usually just called a \textit{vertex}, and a $1$-face an \textit{edge}.} If all of the faces of an abstract simplicial complex $K$ are of dimensional $n$ or less, $K$ is said to be an $n$-dimensional simplicial complex, i.e., the dimension of an ASC is the maximal dimension of its constituent simplices.\par
		A \textit{simplicial map} $f: K \rightarrow K'$ from an abstract simplicial complex on $A$ to an abstract simplicial complex on $B$ is a function induced on simplices by a usual function $A \rightarrow B$, so that the image of any element of $K$ is an element of $K'$.\footnote{Simplicial maps between simplicial complexes are the natural equivalent of continuous maps between topological spaces.} \par 
		Altogether, this data in fact lets us define the category $\textbf{SCpx}$ that has (abstract) simplicial complexes as objects and simplicial maps as morphisms.
	\end{definition} 
ASCs are particularly easy to describe, but one might worry that certain topologically-valuable information gets lost in encoding things in this simplified, set-theoretical fashion. We will ultimately be interested in certain topological information, so it makes sense to want to perform, for a given ASC $K$, what is called the geometrical \textit{realization} $|K|$ of $K$, allowing $K$ to be realized (basically, ``pictured") as some (generalized) triangles glued together in suitable ways, living in a subspace of $\mathbb{R}^n$. For every simplicial complex $K$, there in fact exists a unique (up to simplicial isomorphism) geometric realization $|K|$, so we will mostly not bother to distinguish between ASCs and their geometric realizations.\footnote{Ensuring the uniqueness of this is precisely the reason for considering \textit{ordered} lists in the definition of an ASC, instead of just unordered sets.} Such a realization better explains why we think of objects (sets) in an ASC $K$ as ``faces" or ``simplexes," since via the realization, three-element sets correspond to filled-in triangles, two-element sets to edges, singletons to vertices, etc. For instance, given a set $A = \{a,b,c\}$ for which we have the ASC $K = \{\{a,b\}, \{b,c\}, \{a,c\}, \{a\}, \{b\}, \{c\}\}$, then its realization $|K|$ will be the (hollow) triangle (with a natural orientation). Drawing such a picture yields the usual geometrical notion of a simplicial complex, obtained by ``gluing" together the \textit{standard simplices} along the boundaries; the usual approach then employs these $n$-simplices to probe a topological space via continuous maps into the space.\footnote{More details on these matters can be found in \cite{ghrist_elementary_2014} or \cite{hatcher_algebraic_2001}.} \par 
In more detail, observe that an oriented $0$-simplex thus corresponds to a point $P$, while an oriented $1$-simplex is a directed line segment $P_1P_2$ connecting the points $P_1$ and $P_2$, where we assume that we are traveling in the direction from $P_1$ to $P_2$, i.e., $P_1P_2 \neq P_2P_1$ (however, $P_1P_2 = -P_2P_1$). An oriented 2-simplex will be a triangular region $P_1P_2P_3$, with a prescribed order of movement around the triangle. An oriented $3$-simplex is given by an ordered sequence $P_1P_2P_3P_4$ of four vertices of a solid tetrahedron. Similar definitions hold for $n > 3$. By gluing together various simplices along their boundaries, we get simplicial complexes such as the following: 
\begin{center} 
		\begin{tikzpicture}[scale=0.8]
		
		\coordinate (a) at (4,3.5);
		\coordinate (b) at (3,.8);
		\coordinate (z) at (7.3,0);
		\coordinate (c) at (5,0);
		\coordinate (d) at (6,1.2);
		\coordinate (e) at (1.2,0);
		\coordinate (g) at (-0.5,0);
		
		\draw[thick, fill=black!30,path fading=south] (a) -- (c) -- (d) -- cycle;
		\draw[very thick] (c) -- (d);
		\draw[thick, fill=black!20] (b) -- (c) -- (e) -- cycle;
		\draw[thick] (c) -- (d) -- (z) -- cycle;
		\draw[very thick] (a) -- (c);
		\draw[thick] (a) -- (b);
		\draw[thick, fill=black!30,path fading=south] (a) -- (c) -- (b) -- cycle;
		\draw[thick, fill=black!30] (b) -- (e);
		\draw[thick, dash dot dot] (b) -- (d);
		\draw[thick] (e) -- (g);
		
		\fill[black!20, draw=black, thick] (a) circle (3pt) node[black, above right] {$a$};
		\fill[black!20, draw=black, thick] (b) circle (3pt) node[black, above left] {$d$};
		\fill[black!20, draw=black, thick] (c) circle (3pt) node[black, below right] {$c$};
		\fill[black!20, draw=black, thick] (d) circle (3pt) node[black, above right] {$b$};
		\fill[black!20, draw=black, thick] (e) circle (3pt) node[black, below left] {$e$};
		\fill[black!20, draw=black, thick] (g) circle (3pt) node[black, below left] {$f$};
		\fill[black!20, draw=black, thick] (z) circle (3pt) node[black, below left] {$z$};
		\end{tikzpicture} 
		\end{center}
While it is perhaps useful to visualize things in this way, and while this perspective can be important for connections to other concepts, the alert reader might have observed that given the way ASCs were defined, they should already come with a natural topology, letting us ``by-pass" the geometric realization, and associate to each abstract simplicial complex a particular topological space. This will be useful to us in the construction of sheaves on such spaces towards which we are building.\par  
To see this, first observe that ASCs come with a canonical partial order on faces, given by face subset inclusion (or attachment) relation between vertices, edges, and higher dimensional faces. We can use this fact to define the \textit{face (or cell) category}\index{category!face (cell)} which has for objects the elements of $K$, a cell complex, and (setwise) inclusions of one element/cell of $X$ into another for its morphisms; if $a$ and $b$ are two faces in a complex $X$ with $a \subseteq b$ and $|a| \leq |b|$, we will write $a \rightsquigarrow b$ and say that $a$ is \textit{attached} to $b$.\footnote{Note that technically to make the following construction work we need to use the more general notion of \textit{cell complexes}, the full definition of which can be found in Chapter 4 of \cite{curry_sheaves_2013}; but since ultimately the realization $|K|$ of an ASC $K$ on a finite set, which is the sort of thing we will be dealing with in our example, can be shown to be a cell complex, and simplicial maps $f: K \rightarrow K'$ induce a cellular map $|f|: |K| \rightarrow |K'|$, we will not worry too much about the distinction.} We will then identify a complex with its face poset, writing the incidence relation $a \rightsquigarrow b$. Building on the graphical construction of a complex, the attachments between the faces or cells of a complex can then be displayed in attachment diagrams, where the links represent attachments going from lower dimensional cells to higher (and where any additional attachments that arise as compositions of attachments are left implicit). Observe that the attachment diagram of a graphical complex is just a set partially ordered via the attachment relations, i.e., it is a poset. This face relation poset can be displayed with a Hasse diagram. For instance, assume we have the following simplicial complex $K$, which we imagine has been realized thus:\footnote{Note that the rightmost simplex ($abcd$) is \textit{not} meant to depict a hollow \textit{tetrahedron}; each of the four component triangles are to be thought of as lying in the plane. We have simply spaced it this way to make the sheaf diagrams we build on top of this in a moment a little more readable.} 
\begin{center} 
	\begin{tikzpicture}
	
	\coordinate (a) at (3.2,3.3);
	\coordinate (b) at (3,.8);
	\coordinate (c) at (5,0);
	\coordinate (d) at (6,1.2);
	\coordinate (e) at (1.2,0);
	\coordinate (g) at (-0.5,0);
	
	\draw[thick] (a) -- (c) -- (d) -- cycle;
	\draw[thick, fill=black!20] (b) -- (c) -- (e) -- cycle;
	\draw[very thick] (a) -- (c);
	\draw[thick] (a) -- (b);
	\draw[thick, fill=black!30] (b) -- (e);
	\draw[thick] (b) -- (d);
	\draw[thick] (e) -- (g);
	
	\fill[black!20, draw=black, thick] (a) circle (3pt) node[black, above right] {$a$};
	\fill[black!20, draw=black, thick] (b) circle (3pt) node[black, above left] {$d$};
	\fill[black!20, draw=black, thick] (c) circle (3pt) node[black, below right] {$c$};
	\fill[black!20, draw=black, thick] (d) circle (3pt) node[black, above right] {$b$};
	\fill[black!20, draw=black, thick] (e) circle (3pt) node[black, below left] {$e$};
	\fill[black!20, draw=black, thick] (g) circle (3pt) node[black, below left] {$f$};
	\end{tikzpicture} 
\end{center} 
We can form the diagram of the face-subset relations, where the edges and higher faces are given the natural names (and are assumed to be ordered lexicographically):   
\begin{center} 	
	\begin{tikzpicture}[yscale=0.75, xscale=0.5]
	\node (b1) at (5.5,2) {$cde$};
	\node (d1) at (-3.5,0) {$ab$};
	\node (d) at (-2,0) {$ac$};
	\node (e) at (-0.5,0) {$ad$};
	\node (f) at (1,0) {$bc$};
	\node (d2) at (2.5,0) {$bd$};
	\node (d3) at (4,0) {$cd$};
	\node (e1) at (5.5,0) {$ce$};
	\node (f1) at (7,0) {$de$};
	\node (d4) at (8.5,0) {$ef$};
	\node (g) at (-2.4,-2) {$a$};
	\node (h) at (0,-2) {$b$};
	\node (i) at (2.4,-2) {$c$};
	\node (j) at (4.4,-2) {$d$};
	\node (k) at (6.4,-2) {$e$};
	\node (l) at (8.4,-2) {$f$};
	\node (min) at (3,-4) {$\emptyset$};
	
	\draw[->] (min) -- (g);
	\draw[->] (min) -- (h);
	\draw[->] (min) -- (i);
	\draw[->] (min) -- (j);
	\draw[->] (min) -- (k);
	\draw[->] (min) -- (l);
	\draw[->] (g) -- (d);
	\draw[->] (g) -- (d1);
	\draw[->] (g) -- (e);
	\draw[->] (h) -- (d1);
	\draw[->] (i) -- (d);
	\draw[->] (h) -- (f);
	\draw[->] (h) -- (d2);
	\draw[->] (i) -- (f);
	\draw[->] (i) -- (d3);
	\draw[->] (j) -- (e);
	\draw[->] (j) -- (d2);
	\draw[->] (j) -- (d3);
	\draw[->] (k) -- (e1);
	\draw[->] (k) -- (f1);
	\draw[->] (k) -- (d4);
	\draw[->] (l) -- (d4);
	
	\draw[->] (d3) -- (b1);
	\draw[->] (e1) -- (b1);
	\draw[->] (f1) -- (b1);
	\end{tikzpicture} \end{center} 
But it is more revealing to display this in the form of an attachment diagram, as follows:
\begin{center} 
	\begin{tikzpicture}[yscale=1.2, xscale=1.2,every label/.append style={above=0cm}]
	\node (a) at (2.7,4.2) {$a$};
	\node (d) at (2.2,1.5) {$d$};
	\node (c) at (4.7,0) {$c$};
	\node (b) at (5.9,1.8) {$b$};
	\node (e) at (0.1,0) {$e$};
	\node (f) at (-2,0) {$f$};
	\node (ab) at (4.9,2.6) {$ab$};
	\node (ad) at (2.4,2.6) {$ad$};
	\node (cd) at (3.55,0.8) {$cd$};
	\node (bc) at (5.3,0.9) {$bc$};
	\node (ed) at (1.1,0.8) {$de$};
	\node (cde) at (2.3,0.7) {$cde$};
	\node (ef) at (-0.9,0) {$ef$};
	\node (ce) at (2.3,0) {$ce$};
	\node (bd) at (4.2,1.65) {$bd$};
	\node (ac) at (3.7,2.2) {$ac$};
	
	\draw[->] (a) -- (ab);
	\draw[->] (b) -- (ab);
	\draw[->] (a) -- (ad);
	\draw[->] (d) -- (ad);
	\draw[->] (b) -- (bd);
	\draw[->] (b) -- (bc);
	\draw[->] (c) -- (bc);
	\draw[->] (c) -- (cd);
	\draw[->] (d) -- (ad);
	\draw[->] (d) -- (cd);
	\draw[->] (d) -- (bd);
	\draw[->] (a) -- (ac);
	\draw[->] (c) -- (ac);
	\draw[->] (d) -- (ed);
	\draw[->] (d) -- (ed);
	\draw[->] (e) -- (ed);
	\draw[->] (c) -- (ce);
	\draw[->] (e) -- (ce);
	\draw[->] (ed) -- (cde);
	\draw[->] (cd) -- (cde);
	\draw[->] (ce) -- (cde);
	\draw[->] (e) -- (ef);
	\draw[->] (f) -- (ef);					
	\end{tikzpicture}
\end{center}
It is this face/cell incidence poset that we are regarding as a category, $\mathcal{F}_K$.\index{category!face (cell)} In other words: to a simplicial complex, we can associate a category, which is just the face incidence poset viewed as a category. We can then put the Alexandrov topology\index{topology!Alexandrov} on this poset of face-relations. Then, given this Alexandrov topology on a poset, the usual topological notions (such as of interiors, boundaries, and closures) can be easily understood in terms of the poset itself.\par 
The basic idea is this: pick a simplex; then look at all the other simplices that include that one as a face, i.e., higher dimensional simplices adjacent to it; then regard such ``upper sets" as the open sets. The sets of the form $\uparrow x = \{y \in \mathcal{P}: x \leq y \}$, i.e., the principal upper sets, form a basis for the topology. We can also define the closure of $x$ by $\overline{x} = \{y \in \mathcal{P}: y \leq x \}$; and, provided the poset $\mathcal{P}$ is finite, a basis of closed sets is given by these $\overline{x}$.\footnote{A dual topology thus arises by considering, for a given simplex, all the other simplices that are attached to (included in) it---these are the downsets, which also can serve as the opens for this topology. In the Alexandrov topology, arbitrary intersections of opens are open and arbitrary unions of closed sets are closed; therefore, by exchanging opens with closed sets, we can pass from any Alexandrov space to its dual topology. \par The Alexandrov topology construction is really appropriate when the elements of the underlying poset $\mathcal{P}$ represent finite pieces of information, i.e., are compact, something that is typical for many combinatorial and computer science applications; however, if $\mathcal{P}$ includes infinite elements, then the ``Scott topology" is called for (see Chapter 7 of \cite{vickers_topology_1996} for details on this).} \par 
In the cellular context, for $\sigma$ a cell of a cell complex, the analogue of the principal upper set construction is called the open star of $\sigma$, where this is denoted st$(\sigma)$, and consists of the set of cells $\tau$ such that $\sigma \rightsquigarrow \tau$, i.e., captures all the higher-dimensional cells containing that cell. Then, in terms of the topology, st$(\sigma)$ will be the smallest open set of cells that contain $\sigma$. Taking all the stars and the union of all the stars will give a topology for the complex/simplex---this is the Alexandrov topology. Every intersection of opens in the Alexandrov topology\index{topology!Alexandrov} on a poset $\mathcal{P}$ is open. Thus, a star over $A \subseteq  X$ is then defined to be the intersection of the collection of all open sets containing $A$. The resulting collection of stars will be a basis for the Alexandrov topology.  While stars need not exist in general, in the Alexandrov topology on a poset, there will exist a star of every subset. There is accordingly a dictionary between cellular complexes and Alexandrov spaces, which can be seen by considering that for a cell complex, every cell $\Delta_{\sigma}$ has a star, where this is a set consisting of all those cells $\Delta_{\tau}$ such that $\Delta_{\sigma} \leq \Delta_{\tau}$. \par 
	Note that with respect to the inclusions in the face relation poset $\mathcal{F}_K$, the containment relation for the open sets (stars) in the Alexandrov topology is order-reversing. In other words, the Alexandrov construction will yield (just as we saw last chapter) an order-reversing inclusion functor $\mathcal{P} \rightarrow \mathscr{O}(\mathcal{P})^{op}$.\footnote{More generally, we have described a contravariant functor $Alxd: \textbf{PreOrd} \rightarrow \textbf{Top}$, one that, upon applying $Alxd$ to $X^{op}$, yields a space that has for open sets unions of simplices.} \par 
	It turns out that a sheaf (of sets) over an abstract simplicial complex $K$ can be defined as a covariant functor from the face category\index{category!face (cell)} $\mathcal{F}_K$ of its associated face-relation poset to \textbf{Set}. More generally, given a functor from a preorder or poset $\mathcal{P}$ to a category \textbf{D}, this functor can be used to produce a sheaf on the Alexandrov topology (via a right Kan extension). With this construction, the sheaf gluing axiom for any cover is automatically satisfied! In other words, given a poset endowed with the Alexandrov topology, as anticipated in the last chapter, we do not even need to distinguish between presheaves and sheaves. \par 
The following definition follows \cite{shepard_cellular_1985}, who defined sheaves for more general \textit{cell complexes},\index{complex!cell} which are just a collection of closed disks of certain dimensionality that are glued together along the boundaries. Via its realization $|K|$, an ASC\index{complex!abstract simplicial} $K$ is just a cell complex, so while the definition is more general, it can be applied to a simplicial complex\index{complex!simplicial} (which is what we will work with in the coming example). Just as with ASCs, since cell complexes are built up from simple pieces (cells), the associated attachment diagram exhibiting the relations between cells contains the information of the cell complex itself. Attending to the face poset in particular, then, we will define a \textit{cellular sheaf}, following Shepard, as a covariant functor from the face category of a complex $K$ to some other category \textbf{D}. For concreteness, for the remainder we consider \textbf{D} = \textbf{Vect}.\index{category!of vector spaces}
	\begin{definition}
		A \textit{cellular sheaf} (of vector spaces)\index{sheaf!cellular} on a cell complex $X$ is 
		\begin{itemize}
			\item an assignment of a vector space $F(\sigma)$ to each cell $\sigma$ of $X$,\footnote{This vector space $F(\sigma)$ is then the \textit{stalk}\index{stalk} of $F$ at (or over) $\sigma$.} together with 
			\item a linear transformation 
			\begin{equation*} 
			F_{\sigma \rightsquigarrow \tau}: F(\sigma) \rightarrow F(\tau)
			\end{equation*} for each incident cell pair $\sigma \rightsquigarrow \tau$. 
		\end{itemize}
		These linear maps have to further satisfy the identity relation $F_{\sigma \rightsquigarrow \sigma} = $ id and the usual composition condition, namely
		\begin{equation*}
		\text{ if } \rho \rightsquigarrow \sigma \rightsquigarrow \tau , \text{ then } F_{\rho \rightsquigarrow \tau} = F_{\sigma \rightsquigarrow \tau} \circ F_{\rho \rightsquigarrow \sigma}. 
		\end{equation*}
	\end{definition} 
The reader may be wondering if there is a typo in the direction of the maps described in this definition of a cellular \textit{sheaf}. A sheaf, after all, is a particular presheaf, so one would have expected (order-reversing) \textit{restriction} maps.\index{restriction} But this is not a typo, and in fact goes to the heart of the underlying construction. The reason for the seemingly ``wrong" direction of the arrows---typically, sheaf restriction maps reverse the direction of the arrows, and cosheaves preserve them---is explained (as correct) by a result we have already encountered. Recall that, in general, when dealing with a poset $\mathcal{P}$, we can regard a sheaf on $\mathcal{P}$, once this has been equipped with the upper Alexandrov topology, as a plain old copresheaf (covariant functor) on $\mathcal{P}$ (which could, in turn, be regarded as a presheaf on $\mathcal{P}^{op}$).\footnote{Recall also, that taking $\mathcal{P}^{op}$ instead, sheaves on $\mathcal{P}$, where this is equipped with the \textit{lower} Alexandrov topology (using that $\mathcal{D}(\mathcal{P}) \simeq \mathcal{U}(\mathcal{P}^{op})$), are just presheaves on $\mathcal{P}$.} But the face (or cell) category $\mathcal{F}_{X}$ with which we identify a complex $X$ is a poset. So, using our general result\footnote{Earlier, we just described it in terms of set-valued functors; but in fact, the equivalences hold for (pre)sheaves valued in a category $\textbf{D}$, provided that category is complete and co-complete.} 
\begin{equation*}
\textbf{Sh}(\mathcal{U}(\mathcal{P})) \simeq \textbf{D}^{\mathcal{P}}
\end{equation*} 
letting us move freely between sheaves (valued in $\textbf{D}$) on the upper Alexandrov topology placed on $\mathcal{P}$ and covariant $\textbf{D}$-valued functors (copresheaves) on $\mathcal{P}$, we know that we can freely regard a plain old \textit{covariant functor} $\mathcal{P} \rightarrow \textbf{D}$ as a \textit{sheaf} on $\mathcal{U}(\mathcal{P})$, i.e., on the upper Alexandrov topology on $\mathcal{P}$. The inclusion taking a poset into its upper sets is order-reversing, and the underlying functor of a sheaf (now on the upper sets) is itself order-reversing---their composition, equivalent to the original covariant functor, is accordingly covariant. This accounts for why the definition of a cellular \textit{sheaf} seems to just contain the data of a definition of a covariant functor---that is indeed all it says! The theorem discussed last chapter lets us conflate these two perspectives. In short, we have the slogan: cellular sheaves are covariant functors from the face category into some other category \textbf{D}, e.g., a covariant functor from $\mathcal{F}_X$ to $\textbf{Vec}$ is already just a sheaf of vector spaces. \par    
Dually, we could define: 
	\begin{definition}
		A \textit{cellular cosheaf} (of vector spaces)\index{cosheaf} on a cell complex $X$ is 
		\begin{itemize}
			\item an assignment of a vector space $F(\sigma)$ to each of the cells $\sigma$ of $X$---this vector space $F(\sigma)$ is called the \textit{costalk} of $F$ at (or over) $\sigma$---together with 
			\item a linear transformation $F_{\sigma \rightsquigarrow \tau}: F(\tau) \rightarrow F(\sigma)$ for each incident cell pair $\sigma \rightsquigarrow \tau$. 
		\end{itemize}
		These maps---called the \textit{corestriction maps}---have\index{corestriction} to further satisfy the identity relation $F_{\sigma \rightsquigarrow \sigma} = $ id and the usual composition condition, namely
		\begin{equation*}
		\text{ if } \rho \rightsquigarrow \sigma \rightsquigarrow \tau , \text{ then } F_{\rho \rightsquigarrow \tau} = F_{\rho \rightsquigarrow \sigma} \circ F_{\sigma \rightsquigarrow \tau}. 
		\end{equation*}
	\end{definition}
In the example to follow, we will focus on cellular sheaves; later in the chapter, an example with cellular cosheaves is presented. \par 
Let us also highlight a few more explicit definitions of the corresponding cellular notions that are more or less as you would expect for any sheaf.   
	\begin{definition}
		For $F$ a cellular sheaf on $X$, we define a \textit{global section} $x$ of\index{section!global} $F$ to be a choice $x_{\sigma} \in F(\sigma)$ for each cell $\sigma$ of $X$, where this satisfies 
		\begin{equation*}
		x_{\tau} = F_{\sigma \rightsquigarrow \tau} x_{\sigma} 
		\end{equation*}
		for all $\sigma \rightsquigarrow \tau$. 
	\end{definition} \noindent 
Fundamentally, the data of a cellular sheaf on a complex $X$ amounts to a specification of spaces of local sections on a cover of $X$ (namely, the one given by open stars of cells). 
Ultimately, we will be able to form the category of \textit{all sheaves} $\textbf{Sh}(X)$ over a fixed complex $X$, adopting the only notion of morphisms that there could be, namely as the natural transformations between the corresponding functors.\footnote{In the next definition, we focus on the corresponding notion of morphism for sheaves of vector spaces, but of course, if we were to work with sheaves valued in some other category $\textbf{D}$, then we would require that the maps defined below be the appropriate structure-preserving map, e.g., for sheaves of groups, just homomorphisms.} Explicitly, 
\begin{definition}
	A morphism $f: F \rightarrow G$ of sheaves on a cell complex $X$ is an assignment of a linear map 
	\begin{equation*}
	f_{\sigma}: F(\sigma) \rightarrow G(\sigma)
	\end{equation*}
	to each cell $\sigma$ of $X$, where for each attachment $\sigma \rightsquigarrow \tau$, the usual (natural transformation) compatibility condition holds, making the diagram 
	\begin{center}
	\begin{tikzcd}
			F(\sigma) \arrow[d, "{F(\sigma \rightsquigarrow \tau)}", swap] \arrow[r, "{f_{\sigma}}"] & G(\sigma) \arrow[d, "{G(\sigma \rightsquigarrow \tau)}"] \\ 
			F(\tau) \arrow[r, "{f_{\tau}}", swap] & G(\tau)
	\end{tikzcd}
	\end{center} \par \noindent 
	commute. \par 
	A \textit{sheaf isomorphism}, then, is also defined in the inevitable way: as such a morphism where each of the $f_{\sigma}$ are isomorphisms.\footnote{It is also easy to show that a morphism between sheaves of vector spaces, such as that given above, induces a linear map between the spaces of global sections of the sheaves; moreover, isomorphic sheaves will have isomorphic spaces of global sections.} 
\end{definition} 
	In brief: if $X$ is a cell (or simplicial) complex with the associated face poset category $\mathcal{F}_X$, we identify the complex with its face category, and then a cellular sheaf is just a vector-valued functor $F: \mathcal{F}_X \rightarrow \textbf{Vect}$, while a cellular cosheaf is a vector-valued functor $F: \mathcal{F}_X^{op} \rightarrow \textbf{Vect}$. Thus, to avoid confusion, realize that in the above definitions of cellular (co)sheaves, $X$ was really short for $\mathcal{F}_X$ to which we associate the cell complex, so that a cellular sheaf on $X$ (per the definition) is really a sheaf on $\mathcal{F}_X$ where this has been equipped with the Alexandrov topology, which in turn is uniquely determined by a functor $\mathcal{F}_X \rightarrow \textbf{Vect}$ (and \textit{this} last functor is what the definition is describing). \par  
 	In the example that follows, we illustrate these ideas in a concrete fashion by constructing a cellular sheaf on a particular simplicial complex. 
\begin{example}  	
 Recall our simplicial complex $X$ from before, which we imagine has been realized thus:
 	\begin{center} 
		\begin{tikzpicture}
		
		\coordinate (a) at (3.2,3.3);
		\coordinate (b) at (3,.8);
		\coordinate (c) at (5,0);
		\coordinate (d) at (6,1.2);
		\coordinate (e) at (1.2,0);
		\coordinate (g) at (-0.5,0);
		
		\draw[thick] (a) -- (c) -- (d) -- cycle;
		\draw[thick, fill=black!20] (b) -- (c) -- (e) -- cycle;
		\draw[very thick] (a) -- (c);
		\draw[thick] (a) -- (b);
		\draw[thick, fill=black!30] (b) -- (e);
		\draw[thick] (b) -- (d);
		\draw[thick] (e) -- (g);
		
		\fill[black!20, draw=black, thick] (a) circle (3pt) node[black, above right] {$a$};
		\fill[black!20, draw=black, thick] (b) circle (3pt) node[black, above left] {$d$};
		\fill[black!20, draw=black, thick] (c) circle (3pt) node[black, below right] {$c$};
		\fill[black!20, draw=black, thick] (d) circle (3pt) node[black, above right] {$b$};
		\fill[black!20, draw=black, thick] (e) circle (3pt) node[black, below left] {$e$};
		\fill[black!20, draw=black, thick] (g) circle (3pt) node[black, below left] {$f$};
		\end{tikzpicture} 
		\end{center} 
together with its associated attachment diagram displaying the poset of face relations (i.e., its face category $\mathcal{F}_X$): 
	\begin{center} 
	\begin{tikzpicture}[yscale=1.2, xscale=1.2,every label/.append style={above=0cm}]
		\node (a) at (2.7,4.2) {$a$};
		\node (d) at (2.2,1.5) {$d$};
		\node (c) at (4.7,0) {$c$};
		\node (b) at (5.9,1.8) {$b$};
		\node (e) at (0.1,0) {$e$};
		\node (f) at (-2,0) {$f$};
			\node (ab) at (4.9,2.6) {$ab$};
			\node (ad) at (2.4,2.6) {$ad$};
			\node (cd) at (3.55,0.8) {$cd$};
			\node (bc) at (5.3,0.9) {$bc$};
			\node (ed) at (1.1,0.8) {$de$};
			\node (cde) at (2.3,0.7) {$cde$};
			\node (ef) at (-0.9,0) {$ef$};
			\node (ce) at (2.3,0) {$ce$};
			\node (bd) at (4.2,1.65) {$bd$};
			\node (ac) at (3.7,2.2) {$ac$};
			
				\draw[->] (a) -- (ab);
				\draw[->] (b) -- (ab);
				\draw[->] (a) -- (ad);
				\draw[->] (d) -- (ad);
				\draw[->] (b) -- (bd);
				\draw[->] (b) -- (bc);
				\draw[->] (c) -- (bc);
				\draw[->] (c) -- (cd);
				\draw[->] (d) -- (ad);
				\draw[->] (d) -- (cd);
					\draw[->] (d) -- (bd);
					\draw[->] (a) -- (ac);
					\draw[->] (c) -- (ac);
					\draw[->] (d) -- (ed);
					\draw[->] (d) -- (ed);
					\draw[->] (e) -- (ed);
					\draw[->] (c) -- (ce);
					\draw[->] (e) -- (ce);
						\draw[->] (ed) -- (cde);
						\draw[->] (cd) -- (cde);
						\draw[->] (ce) -- (cde);
							\draw[->] (e) -- (ef);
							\draw[->] (f) -- (ef);					
	\end{tikzpicture}
	\end{center}
	Let us see what a sheaf (of vector spaces) on $X$ (or rather, on the corresponding face attachment diagram $\mathcal{F}_X$) looks like. Well, we need to spell out all the data of the topology, covers, and the sheaf conditions, right? No! Using the main theorem, it will suffice to just describe a vector-valued covariant functor on this diagram---and this will already contain all the data of a sheaf! This is one of the many instances where very abstract category-theoretic results, which may be difficult to understand at first, can make our lives a lot easier in practice.\par 
	Thus, following the definition of a cellular sheaf, for the values of the sheaf over cells, this will just amount to the specification or assignment of values (spaces) to each of the cells of the simplices, data that comes in the form of vectors, e.g.,  
		\begin{center}  
		\begin{tikzpicture}[yscale=1.25, xscale=1.15]
		\node (a) at (3.75,4.5) {$\mathbb{R}^2$};
		\node (d) at (2.2,1.65) {$\mathbb{R}$};
		\node (c) at (4.7,0) {$\mathbb{R}^2$};
		\node (b) at (6,2) {$\mathbb{R}^3$};
		\node (e) at (-0.1,0) {$\mathbb{R}^3$};
		\node (f) at (-2.5,0) {$\mathbb{R}^3$};
		\node (ab) at (5.1,3) {$\mathbb{R}^2$};
		\node (ad) at (2.8,2.7) {$\mathbb{R}$};
		\node (cd) at (3.5,0.85) {$\mathbb{R}^2$};
		\node (bc) at (5.3,0.9) {$\mathbb{R}$};
		\node (ed) at (1.1,0.8) {$\mathbb{R}^2$};
		\node (cde) at (2.3,0.75) {$\mathbb{R}$};
		\node (ef) at (-1.4,0) {$\mathbb{R}^2$};
		\node (ce) at (2.3,0) {$\mathbb{R}^2$};
		\node (bd) at (3.9,1.8) {$\mathbb{R}$};
		\node (ac) at (4.3,2.2) {$\mathbb{R}^2$};
		
		\draw[->] (a) -- node {}(ab);
		\draw[->] (b) -- node {}(ab);
		\draw[->] (a) -- node {}(ad);
		\draw[->] (d) -- node {}(ad);
		\draw[->] (b) -- node {}(bd);
		\draw[->] (b) -- node {}(bc);
		\draw[->] (c) -- node {}(bc);
		\draw[->] (c) -- node {}(cd);
		\draw[->] (d) -- node {}(cd);
		\draw[->] (d) -- node {}(bd);
		\draw[->] (a) -- node {}(ac);
		\draw[->] (c) -- node {}(ac);
		\draw[->] (d) -- node {}(ed);
		\draw[->] (e) -- node {}(ed);
		\draw[->] (c) -- node {}(ce);
		\draw[->] (e) -- node {}(ce);
		\draw[->] (ed) -- node {}(cde);
		\draw[->] (cd) -- node {}(cde);
		\draw[->] (ce) -- node {}(cde);
		\draw[->] (e) -- node {}(ef);
		\draw[->] (f) -- node {}(ef);
		\end{tikzpicture} 
	\end{center}  
	But we need maps as well. The maps, for their part, may be thought of as representing some sort of local constraints or as enforcing certain relations between the data. In general, when the stalks $F(\sigma)$ have structure---for instance, here they are vector spaces---then a sheaf of that type is obtained when the restriction maps preserve this structure. In other words, we should have a function (in our particular case of vector space assignments, these will be given by a linear map) assigned to each inclusion of faces in such a way that the diagram commutes, i.e., the composition of functions throughout the diagram is path independent. The following sheaf diagram nicely displays all these ideas: 
	\begin{center}  
		\begin{tikzpicture}[yscale=1.44, xscale=1.32]
		\node (a) at (3.75,4.5) {$\mathbb{R}^2$};
		\node (d) at (2.2,1.65) {$\mathbb{R}$};
		\node (c) at (4.7,0) {$\mathbb{R}^2$};
		\node (b) at (6,2) {$\mathbb{R}^3$};
		\node (e) at (-0.1,0) {$\mathbb{R}^3$};
		\node (f) at (-2.5,0) {$\mathbb{R}^3$};
		\node (ab) at (5.1,3) {$\mathbb{R}^2$};
		\node (ad) at (2.8,2.7) {$\mathbb{R}$};
		\node (cd) at (3.5,0.85) {$\mathbb{R}^2$};
		\node (bc) at (5.3,0.9) {$\mathbb{R}$};
		\node (ed) at (1.1,0.8) {$\mathbb{R}^2$};
		\node (cde) at (2.3,0.75) {$\mathbb{R}$};
		\node (ef) at (-1.4,0) {$\mathbb{R}^2$};
		\node (ce) at (2.3,0) {$\mathbb{R}^2$};
		\node (bd) at (3.9,1.8) {$\mathbb{R}$};
		\node (ac) at (4.3,2.2) {$\mathbb{R}^2$};
		
		\tiny 
		\draw[->] (a) -- node[label={[xshift=0.5cm]90:$\begin{psmallmatrix} 1 & 0 \\ -1 & 2 \end{psmallmatrix}$}] {}(ab);
		\draw[->] (b) -- node[label={[xshift=0.5cm]90:$\begin{psmallmatrix} 1 & 0 & 1\\ 0 & -1 & -1 \end{psmallmatrix}$}] {}(ab);
		\draw[->] (a) -- node[label={[xshift=-0.5cm]90:$\begin{psmallmatrix} 0 & -2 \end{psmallmatrix}$}] {}(ad);
		\draw[->] (d) -- node[label={[xshift=-0.1cm]90:$\begin{psmallmatrix} 1 \end{psmallmatrix}$}] {}(ad);
		\draw[->] (b) -- node[label={[xshift=0.2cm]90:$\begin{psmallmatrix} 2 & 0 & 2 \end{psmallmatrix}$}] {}(bd);
		\draw[->] (b) -- node[label={[xshift=0.4cm, yshift=-0.7cm]90:$\begin{psmallmatrix} 1 & 2 & 1 \end{psmallmatrix}$}] {}(bc);
		\draw[->] (c) -- node[label={[xshift=0.4cm, yshift=-0.6cm]90:$\begin{psmallmatrix} 1 & 1 \end{psmallmatrix}$}] {}(bc);
		\draw[->] (c) -- node[label={[xshift=-0.2cm, yshift=-0.7cm]90:$\begin{psmallmatrix} -1 & -1 \\ 3 & 1 \end{psmallmatrix}$}] {}(cd);
		\draw[->] (d) -- node[label={[xshift=0.25cm, yshift=-0.2cm]90:$\begin{psmallmatrix} 0.5 \\ 1 \end{psmallmatrix}$}] {}(cd);
		\draw[->] (d) -- node[label={[yshift=0cm]90:$\begin{psmallmatrix} -3 \end{psmallmatrix}$}] {}(bd);
		\draw[->] (a) -- node[label={[xshift=-0.3cm, yshift=-0.4cm]90:$\begin{psmallmatrix} 1 & 0 \\ 0 & 1 \end{psmallmatrix}$}] {}(ac);
		\draw[->] (c) -- node[label={[xshift=0.4cm, yshift=-0.4cm]90:$\begin{psmallmatrix} 3 & 3 \\ 1 & 1 \end{psmallmatrix}$}] {}(ac);
		\draw[->] (d) -- node[label={[xshift=-0.2cm]90:$\begin{psmallmatrix} 3 \\ 1 \end{psmallmatrix}$}] {}(ed);
		\draw[->] (e) -- node[label={[xshift=-0.2cm, yshift=-0.1cm]90:$\begin{psmallmatrix} 2 & 0 & 1\\ 0 & 3 & -1 \end{psmallmatrix}$}] {}(ed);
		\draw[->] (c) -- node[label={[yshift=-0.8cm]90:$\begin{psmallmatrix} 1 & -1 \\ -1 & 2 \end{psmallmatrix}$}] {}(ce);
		\draw[->] (e) -- node[label={[yshift=-0.8cm]90:$\begin{psmallmatrix} 2 & -3 & 2\\ 1 & 0 & 7.5 \end{psmallmatrix}$}] {}(ce);
		\draw[->] (ed) -- node[label={[yshift=-0.2cm]90:$\begin{psmallmatrix} 1 & -1 \end{psmallmatrix}$}] {}(cde);
		\draw[->] (cd) -- node[label={[yshift=-0.2cm]90:$\begin{psmallmatrix} 2 & 1 \end{psmallmatrix}$}] {}(cde);
		\draw[->] (ce) -- node[label={[xshift=0.5cm, yshift=-0.5cm]90:$\begin{psmallmatrix} 1 & 0 \end{psmallmatrix}$}] {}(cde);
		\draw[->] (e) -- node[label={[yshift=-0.8cm]90:$\begin{psmallmatrix} 2 & 0 & 2\\ 1 & -1 & 1 \end{psmallmatrix}$}] {}(ef);
		\draw[->] (f) -- node[label={[yshift=-0.8cm]90:$\begin{psmallmatrix} 0 & 1 & 1\\ 1 & -1 & 0 \end{psmallmatrix}$}] {}(ef);
		\end{tikzpicture} 
		\end{center}  
		A sheaf is generated by its values specified on individual simplices of $X$, i.e., by local sections\index{section!local} specified on the vertices. But a sheaf is not just this data. The restriction\index{restriction} maps of a sheaf are an essential part of the construction, as they encode how any local sections can be extended into sections over a larger part of the diagram (ultimately throughout the entire diagram), and so it is precisely via the restriction maps that it is made explicit how the local sections can be glued together. The sheaf assignment given over \textit{all} of $X$ will be specified by a collection of local sections that can be extended along all the restriction maps to higher-dimensional faces. There may be some flexibility or freedom in the actual data assignments over a vertex, but they are not entirely arbitrary, for the restriction maps encode how local assignments---values specified on certain parts of the diagram---can be extended to other parts of the diagram, and so the maps will constrain the assignments in various ways. If the reader would like to get a good ``working" understanding of the important distinction between a local section and a global section, it would be useful to closely consider what happens in the concrete case when we assign, for instance, the value $\begin{psmallmatrix}
		1 \\ 0 \\ -1
		\end{psmallmatrix}$ to the stalk over the vertex $e$, versus what happens when we assign, for instance, the value $\begin{psmallmatrix}
		-1 \\ 2 \\ 2
		\end{psmallmatrix}$ to the same vertex. In the first case, one finds that we can ``extend" or propagate this particular selection along \textit{some} of the arrows to those stalks highlighted in gray; but then there is a problem, an obstruction to our continuing this process any further upward along the edges of the diagram:  
		\begin{center} 
			\begin{tikzpicture}[yscale=1.44, xscale=1.32]
			\footnotesize
			\node (a) at (3.75,4.5) {$\mathbb{R}^2$};
			\node[gray] (d) at (2.2,1.65) {$\textbf{?}$};
			\node[gray] (c) at (4.7,0) {$\begin{psmallmatrix} 1 \\ 1  \end{psmallmatrix}$};
			\node (b) at (6,2) {$\mathbb{R}^3$};
			\node[gray] (e) at (-0.1,0) {$\begin{psmallmatrix} 1 \\ 0 \\ -1  \end{psmallmatrix}$};
			\node[gray] (f) at (-2.5,0) {$\begin{psmallmatrix} 1 \\ 1 \\ -1  \end{psmallmatrix}$};
			\node (ab) at (5.1,3) {$\mathbb{R}^2$};
			\node (ad) at (2.8,2.7) {$\mathbb{R}$};
			\node[gray] (cd) at (3.5,0.85) {$\begin{psmallmatrix} -2 \\ -4  \end{psmallmatrix}$};
			\node (bc) at (5.3,0.9) {$\mathbb{R}$};
			\node[gray] (ed) at (1.1,0.8) {$\begin{psmallmatrix} 1 \\ 1  \end{psmallmatrix}$};
			\node[gray] (cde) at (2.3,0.75) {$\begin{psmallmatrix} 0  \end{psmallmatrix}$};
			\node[gray] (ef) at (-1.4,0) {$\begin{psmallmatrix} 0 \\ 0  \end{psmallmatrix}$};
			\node[gray] (ce) at (2.3,0) {$\begin{psmallmatrix} 0 \\ -6.5  \end{psmallmatrix}$};
			\node (bd) at (3.9,1.8) {$\mathbb{R}$};
			\node (ac) at (4.3,2.2) {$\mathbb{R}^2$};
			\tiny 
			\draw[->] (a) -- node[label={[xshift=0.5cm]90:$\begin{psmallmatrix} 1 & 0 \\ -1 & 2 \end{psmallmatrix}$}] {}(ab);
			\draw[->] (b) -- node[label={[xshift=0.5cm]90:$\begin{psmallmatrix} 1 & 0 & 1\\ 0 & -1 & -1 \end{psmallmatrix}$}] {}(ab);
			\draw[->] (a) -- node[label={[xshift=-0.5cm]90:$\begin{psmallmatrix} 0 & -2 \end{psmallmatrix}$}] {}(ad);
			\draw[->] (d) -- node[label={[xshift=-0.1cm]90:$\begin{psmallmatrix} 1 \end{psmallmatrix}$}] {}(ad);
			\draw[->] (b) -- node[label={[xshift=0.2cm]90:$\begin{psmallmatrix} 2 & 0 & 2 \end{psmallmatrix}$}] {}(bd);
			\draw[->] (b) -- node[label={[xshift=0.4cm, yshift=-0.7cm]90:$\begin{psmallmatrix} 1 & 2 & 1 \end{psmallmatrix}$}] {}(bc);
			\draw[->] (c) -- node[label={[xshift=0.4cm, yshift=-0.6cm]90:$\begin{psmallmatrix} 1 & 1 \end{psmallmatrix}$}] {}(bc);
			\draw[->] (c) -- node[label={[xshift=-0.2cm, yshift=-0.7cm]90:$\begin{psmallmatrix} -1 & -1 \\ 3 & 1 \end{psmallmatrix}$}] {}(cd);
			\draw[->] (d) -- node[label={[xshift=0.25cm, yshift=-0.2cm]90:$\begin{psmallmatrix} 0.5 \\ 1 \end{psmallmatrix}$}] {}(cd);
			\draw[->] (d) -- node[label={[yshift=0cm]90:$\begin{psmallmatrix} -3 \end{psmallmatrix}$}] {}(bd);
			\draw[->] (a) -- node[label={[xshift=-0.3cm, yshift=-0.4cm]90:$\begin{psmallmatrix} 1 & 0 \\ 0 & 1 \end{psmallmatrix}$}] {}(ac);
			\draw[->] (c) -- node[label={[xshift=0.4cm, yshift=-0.4cm]90:$\begin{psmallmatrix} 3 & 3 \\ 1 & 1 \end{psmallmatrix}$}] {}(ac);
			\draw[->] (d) -- node[label={[xshift=-0.2cm]90:$\begin{psmallmatrix} 3 \\ 1 \end{psmallmatrix}$}] {}(ed);
			\draw[->] (e) -- node[label={[xshift=-0.2cm, yshift=-0.1cm]90:$\begin{psmallmatrix} 2 & 0 & 1\\ 0 & 3 & -1 \end{psmallmatrix}$}] {}(ed);
			\draw[->] (c) -- node[label={[yshift=-0.8cm]90:$\begin{psmallmatrix} 1 & -1 \\ -1 & 2 \end{psmallmatrix}$}] {}(ce);
			\draw[->] (e) -- node[label={[yshift=-0.8cm]90:$\begin{psmallmatrix} 2 & -3 & 2\\ 1 & 0 & 7.5 \end{psmallmatrix}$}] {}(ce);
			\draw[->] (ed) -- node[label={[yshift=-0.2cm]90:$\begin{psmallmatrix} 1 & -1 \end{psmallmatrix}$}] {}(cde);
			\draw[->] (cd) -- node[label={[yshift=-0.2cm]90:$\begin{psmallmatrix} 2 & 1 \end{psmallmatrix}$}] {}(cde);
			\draw[->] (ce) -- node[label={[xshift=0.5cm, yshift=-0.5cm]90:$\begin{psmallmatrix} 1 & 0 \end{psmallmatrix}$}] {}(cde);
			\draw[->] (e) -- node[label={[yshift=-0.8cm]90:$\begin{psmallmatrix} 2 & 0 & 2\\ 1 & -1 & 1 \end{psmallmatrix}$}] {}(ef);
			\draw[->] (f) -- node[label={[yshift=-0.8cm]90:$\begin{psmallmatrix} 0 & 1 & 1\\ 1 & -1 & 0 \end{psmallmatrix}$}] {}(ef);
			\end{tikzpicture} 
		\end{center}
		The reader will note that there is simply \textit{no value} that might be placed at the stalk over vertex $d$ (thus the `?') that would allow us to continue with this extension process. If we assigned $(-4)$ to the stalk over $d$, this would indeed be consistent with the map $\begin{psmallmatrix}
		0.5 \\ 1 \end{psmallmatrix}$ proceeding down and to the right and landing in the stalk over $cd$, which in turn lands us, perfectly consistently, in the stalk over $cde$ with the value $(0)$, as required by the other restriction maps. However, by following what happens to that same assignment $-4$ under the action of the map down and to the left via $\begin{psmallmatrix} 3 \\ 1 \end{psmallmatrix}$, we see that we would get $\begin{psmallmatrix}
		-12 \\ -4 \end{psmallmatrix}$ which, when further mapped under $\begin{psmallmatrix}
		1 & -1 \end{psmallmatrix}$ yields $(-8)$. We thus cannot assign $-4$---or \textit{anything} for that matter---to the stalk over $d$, given our original assignment over $e$. \par 
		The original assignment at the stalk over $e$, then, is said to describe a strictly \textit{local section}, one that importantly cannot be extended globally, i.e., over the entire complex. By beginning with other values at the same (or, if one desires, another) vertex, the reader can explore various other solutions that are merely local versus those that manage to be global. In this way, the reader will not only discover that some local solutions or sections are ``more local" than others, but will discover other ``types" of obstructions to the extension of local sections.\index{section!local} For instance, whereas with our test assignment above it turned out that there was simply \textit{no value at all} that could be assigned to the stalk over vertex $d$, while maintaining consistency with the other stalk assignments required by the restriction maps, another (less serious) issue one frequently encounters is that a specific assignment on one stalk ends up requiring two \textit{different} assignments at some stalk.\footnote{We might here take the opportunity to mention that in a variety of applications of sheaves, including beyond sheaves of vector spaces on complexes, it is possible that \textit{exact equality} of assignments will be unattainable or the least valuable thing to consider. There are a few ways to develop machinery to accommodate this, but it is the beyond the scope of this book to touch on them in detail. Instead, for one particularly friendly approach, the reader is referred to Robinson's work, who has proposed to deal with this situation via formalizing a ``consistency structure" with a corresponding notion of \textit{distance} between assignments (say, with the structure of a pseudo-metric space). Moreover, the notion of \textit{pseudosections} can be developed, and Robinson has shown that in fact pseudosections are already sections, just with respect to a different sheaf; for instance, pseudosections of a sheaf over an ASC $X$ are veritable sections of another sheaf over the barycentric subdivision of $X$. See \cite{robinson_pseudosections_2015}, \cite{robinson_sheaves_2016}, and \cite{praggastis_maximal_2016} for more details. Pushing this a little further, we could analyze data using the \textit{consistency radius} of the sheaf, i.e., the maximum distance between the value in a stalk and the values propagated along the restriction maps. By imposing such a consistency structure on the sheaf, this could tell you ``how far" a particular data instance was from conforming to the consistency requirements stipulated by the structure encoded by the entire sheaf. In other words, given a particular data assignment, it could be used to inform how to find the \textit{closest global section} (where ``closest," of course, would be given by, say, the pseudo-metric placed on the assignments). See \cite{robinson_assignments_2018}.} \par 
		In contrast to the above failures, one observes that by seeding the stalk over vertex $e$ with the value $\begin{psmallmatrix}
		-1 \\ 2 \\ 2
		\end{psmallmatrix}$, we encounter no such obstruction to the extension of this assignment to a consistent assignment of values over the \textit{entire} diagram, thus yielding what is appropriately called a \textit{global section}.\index{section!global} A global section is just a selection of value assignments from each of the stalks over all the cells that is consistent with \textit{all} of the restriction maps of the diagram:
		\begin{center} 
		\begin{tikzpicture}[yscale=1.6, xscale=1.38]
		\footnotesize
		\node (a) at (3.75,4.5) {$\begin{psmallmatrix} 3 \\ 1  \end{psmallmatrix}$};
		\node (d) at (2.2,1.65) {$(-2)$};
		\node (c) at (4.7,0) {$\begin{psmallmatrix} -1.5 \\ 2.5  \end{psmallmatrix}$};
		\node (b) at (6,2) {$\begin{psmallmatrix} 1 \\ -1 \\ 2 \end{psmallmatrix}$};
		\node (e) at (-0.1,0) {$\begin{psmallmatrix} -1 \\ 2 \\ 2  \end{psmallmatrix}$};
		\node (f) at (-2.5,0) {$\begin{psmallmatrix} 2 \\ 3 \\ -1  \end{psmallmatrix}$};
		\node (ab) at (5.1,3) {$\begin{psmallmatrix} 3 \\ -1  \end{psmallmatrix}$};
		\node (ad) at (2.8,2.7) {$\begin{psmallmatrix} 1  \end{psmallmatrix}$};
		\node (cd) at (3.5,0.85) {$\begin{psmallmatrix} -1 \\ -2  \end{psmallmatrix}$};
		\node (bc) at (5.3,0.9) {$\begin{psmallmatrix} 1  \end{psmallmatrix}$};
		\node (ed) at (1.1,0.8) {$\begin{psmallmatrix} -6 \\ -2  \end{psmallmatrix}$};
		\node (cde) at (2.3,0.75) {$\begin{psmallmatrix} -4  \end{psmallmatrix}$};
		\node (ef) at (-1.4,0) {$\begin{psmallmatrix} 2 \\ -1  \end{psmallmatrix}$};
		\node (ce) at (2.3,0) {$\begin{psmallmatrix} -4 \\ 6.5  \end{psmallmatrix}$};
		\node (bd) at (3.9,1.8) {$\begin{psmallmatrix} 6  \end{psmallmatrix}$};
		\node (ac) at (4.3,2.2) {$\begin{psmallmatrix} 3 \\ 1  \end{psmallmatrix}$};
		\tiny 
		\draw[->] (a) -- node[label={[xshift=0.5cm]90:$\begin{psmallmatrix} 1 & 0 \\ -1 & 2 \end{psmallmatrix}$}] {}(ab);
		\draw[->] (b) -- node[label={[xshift=0.5cm]90:$\begin{psmallmatrix} 1 & 0 & 1\\ 0 & -1 & -1 \end{psmallmatrix}$}] {}(ab);
		\draw[->] (a) -- node[label={[xshift=-0.5cm]90:$\begin{psmallmatrix} 0 & -2 \end{psmallmatrix}$}] {}(ad);
		\draw[->] (d) -- node[label={[xshift=-0.1cm]90:$\begin{psmallmatrix} 1 \end{psmallmatrix}$}] {}(ad);
		\draw[->] (b) -- node[label={[xshift=0.2cm]90:$\begin{psmallmatrix} 2 & 0 & 2 \end{psmallmatrix}$}] {}(bd);
		\draw[->] (b) -- node[label={[xshift=0.4cm, yshift=-0.7cm]90:$\begin{psmallmatrix} 1 & 2 & 1 \end{psmallmatrix}$}] {}(bc);
		\draw[->] (c) -- node[label={[xshift=0.4cm, yshift=-0.6cm]90:$\begin{psmallmatrix} 1 & 1 \end{psmallmatrix}$}] {}(bc);
		\draw[->] (c) -- node[label={[xshift=-0.2cm, yshift=-0.7cm]90:$\begin{psmallmatrix} -1 & -1 \\ 3 & 1 \end{psmallmatrix}$}] {}(cd);
		\draw[->] (d) -- node[label={[xshift=0.25cm, yshift=-0.2cm]90:$\begin{psmallmatrix} 0.5 \\ 1 \end{psmallmatrix}$}] {}(cd);
		\draw[->] (d) -- node[label={[yshift=0cm]90:$\begin{psmallmatrix} -3 \end{psmallmatrix}$}] {}(bd);
		\draw[->] (a) -- node[label={[xshift=-0.3cm, yshift=-0.4cm]90:$\begin{psmallmatrix} 1 & 0 \\ 0 & 1 \end{psmallmatrix}$}] {}(ac);
		\draw[->] (c) -- node[label={[xshift=0.4cm, yshift=-0.4cm]90:$\begin{psmallmatrix} 3 & 3 \\ 1 & 1 \end{psmallmatrix}$}] {}(ac);
		\draw[->] (d) -- node[label={[xshift=-0.2cm]90:$\begin{psmallmatrix} 3 \\ 1 \end{psmallmatrix}$}] {}(ed);
		\draw[->] (e) -- node[label={[xshift=-0.2cm, yshift=-0.1cm]90:$\begin{psmallmatrix} 2 & 0 & 1\\ 0 & 3 & -1 \end{psmallmatrix}$}] {}(ed);
		\draw[->] (c) -- node[label={[yshift=-0.8cm]90:$\begin{psmallmatrix} 1 & -1 \\ -1 & 2 \end{psmallmatrix}$}] {}(ce);
		\draw[->] (e) -- node[label={[yshift=-0.8cm]90:$\begin{psmallmatrix} 2 & -3 & 2\\ 1 & 0 & 7.5 \end{psmallmatrix}$}] {}(ce);
		\draw[->] (ed) -- node[label={[yshift=-0.2cm]90:$\begin{psmallmatrix} 1 & -1 \end{psmallmatrix}$}] {}(cde);
		\draw[->] (cd) -- node[label={[yshift=-0.2cm]90:$\begin{psmallmatrix} 2 & 1 \end{psmallmatrix}$}] {}(cde);
		\draw[->] (ce) -- node[label={[xshift=0.5cm, yshift=-0.5cm]90:$\begin{psmallmatrix} 1 & 0 \end{psmallmatrix}$}] {}(cde);
		\draw[->] (e) -- node[label={[yshift=-0.8cm]90:$\begin{psmallmatrix} 2 & 0 & 2\\ 1 & -1 & 1 \end{psmallmatrix}$}] {}(ef);
		\draw[->] (f) -- node[label={[yshift=-0.8cm]90:$\begin{psmallmatrix} 0 & 1 & 1\\ 1 & -1 & 0 \end{psmallmatrix}$}] {}(ef);
		
		\end{tikzpicture} 
		\end{center} \par \noindent 
		In last chapter's discussion of sheaves of sections, from the various examples it was clear that some local sections of a sheaf (of \textit{sets}) will not extend to global sections. As we just saw, in the cellular sheaf of vector spaces, the same sort of thing occurs, i.e., sections can remain strictly local, when they cannot be defined across all the faces of the simplex or when they conflict with the constraints of the restriction maps. There might be interesting local solutions among the variable sets of solutions to a local problem, but only those solutions that can be consistently propagated along the entire diagram, respecting the sheaf restriction maps, will provide us with a global section or solution.
	\end{example}
\subsection{Sheaf Morphisms and Some Operations on Sheaves}
		Earlier, in discussing the definition of cellular sheaves, we mentioned the category $\textbf{Sh}(X)$ of all cellular sheaves on a fixed cellular space $X$, where we observed that the morphisms of this category are, inevitably, just natural transformations between the functors defining the cellular sheaves. But suppose we no longer fix the cellular space, so that we are considering cellular sheaves on (possibly) different spaces. We would like to extend the notion of a morphism of cellular sheaves to provide maps between sheaves on different spaces. \par 
		Suppose that we have two attachment diagrams, each corresponding to two different cell complexes. For concreteness, we consider the following two: the one on the left consists of three vertices with two connecting edges; the one on the right of two vertices and a single connecting edge. A simplicial (or, more generally, cellular) map is a map going from the one on the right to the one on the left:
		\par \begin{center} 	
			\begin{tikzcd}[row sep= small]
				v_0 \arrow[d] & & & & w_0 \arrow[dd] \\
				v_0 v_1 \\
				v_1 \arrow[u] \arrow[d] & {} & & \arrow[ll, dashed] & w_0 w_1 \\
				v_1 v_2 \\
				v_2 \arrow[u] & & & & w_1 \arrow[uu]
			\end{tikzcd}
		\end{center} 
		 Then if we represent the sheaves as ``sitting over" the attachment diagrams representing the simplicial complexes, a \textit{sheaf morphism}---displayed\index{sheaf!morphism} in bold, and going in the opposite direction as the simplicial (cellular) map---will look something like: 
		\begin{center} 
			\begin{tikzcd}[row sep=1.9em,column sep =1.4em]
				\mathbb{Z}^3 \arrow[bend left, thick]{rrrrr}{\begin{psmallmatrix} 1 \amsamp 1 \amsamp 1 \\ 1 \amsamp 0 \amsamp 3 \end{psmallmatrix}} \arrow{dr}{\begin{psmallmatrix} 1 \amsamp 0 \amsamp 2 \\ 1 \amsamp 3 \amsamp 1 \end{psmallmatrix}} \arrow[dd, leftarrow,gray] && & \checkmark & &
				\mathbb{Z}^2 \arrow[dd, near start,leftarrow,gray] \arrow{ddrrrr}{\begin{psmallmatrix} 3 \amsamp 2 \end{psmallmatrix}} \\
				& \mathbb{Z}^2 \arrow[crossing over, bend left =10, near start, thick]{drrrrrrrr}{\begin{psmallmatrix} 4 \amsamp 1 \end{psmallmatrix}} \\ 
				v_0 \arrow[rrrrr, near end,dashed, leftarrow] \arrow[dr] & & \mathbb{Z}^2 \arrow[dd, leftarrow,gray] \arrow[ul] \arrow[dr] \arrow[rrrrrrr, bend left =10, thick] & & & w_0 \arrow[ddrrrr] & & & & \mathbb{Z} \\
				& v_0 v_1 \arrow[uu,crossing over, near end,gray]&& \mathbb{Z} \arrow[urrrrrr, bend left = 10, thick] \arrow[dd,leftarrow,gray]  \\
				& & v_1 \arrow[dr] \arrow[ul] & & \mathbb{Z} \arrow[dd,leftarrow, gray] \arrow[ul] \arrow[rrrrrrrrr, bend left =20, thick] & & & & & w_0 w_1 \arrow[uu,gray] \arrow[lllllll,dashed] & & & & \mathbb{Z} \arrow[dd,leftarrow,gray] \arrow[uullll] \\
				& & & v_1 v_2 \\
				& & & & v_2 \arrow[ul] & & & & & & & & & w_1 \arrow[uullll] \arrow[lllllllll,dashed]
			\end{tikzcd}  
		\end{center}
		In other words, a sheaf morphism takes data in the stalks over two sheaves and relates them through linear maps in such a way that the resulting diagram commutes.\footnote{A lot is left unspecified in the above diagram, displaying how this might look only for the top commuting square.} Of course, we could consider a morphism between sheaves on a fixed space $X$ as a special case of this approach, but we have really been heading towards the more general case of a sheaf morphism involving different spaces. 
		\begin{definition}
			A \textit{sheaf morphism}\index{sheaf!morphism} $s: F \rightarrow G$ from a sheaf $F$ over a space $Y$ to a sheaf $G$ over the space $X$ consists of the following data: 
			\begin{itemize}
				\item a cellular map $f: X \rightarrow Y$, 
				\item a collection of (linear) maps $l_{\sigma}: F(f(\sigma)) \rightarrow G(\sigma)$ such that for each attachment map $\sigma \rightsquigarrow \tau$ in $X$, the following diagram commutes: 
					\begin{center}  
					\begin{tikzcd}
						F(f(\sigma)) \arrow[d, "{F(f(\sigma) \rightsquigarrow f(\tau))}", swap] \arrow[r, "{l_{\sigma}}"] & G(x) \arrow[d, "{G(\sigma \rightsquigarrow \tau)}"] \\
						F(f(\tau)) \arrow[r, "{l_{\tau}}", swap] & G(\tau) 
					\end{tikzcd}
				\end{center}  
			\end{itemize}
		\end{definition} \noindent 
	The reader may note how this in fact makes use of the \textit{pullback sheaf}\index{sheaf!inverse image (pullback)} notion, first mentioned in \ref{example: manifold}, where for a map $f: X \rightarrow Y$, and a sheaf $F$ on $Y$, the pullback $f^*F$ will be a sheaf on $X$ defined by $(f^*F)(\gamma) = F(f(\gamma))$ and $(f^*F)(\sigma \rightsquigarrow \tau) = F(f(\sigma) \rightsquigarrow f(\tau))$. \par 
	Such notions (and others, such as the pushforward sheaf) are clearly useful for switching base spaces. Sheaf morphisms can also be composed, under certain conditions, leading to \textit{sequences of sheaves}, linked together by sheaf morphisms. Certain sequences will even exhibit special algebraic properties, like \textit{exactness}, which will have significance in a number of applications. We take up these matters in a moment, after a brief discussion of a few further operations on sheaves. \par 
	We have seen a few constructions such as that of subsheaves, pullback sheaves, and  pushforward sheaves, where ``old sheaves" are used to generate new ones. Here is a very brief look at just a few other important things one can do \textit{to} or \textit{with} sheaves, to generate new sheaves; specifically, we focus on indicating a few of the \textit{algebraic} operations one can perform on sheaves. These notions can be defined in greater generality, but we will stick to the cellular context. 
	\begin{definition}
	For $F$ and $G$, two sheaves of vector spaces on a cell complex $X$, we can define $F \oplus G$, their \textit{direct sum},\index{sheaf!direct sum} in the natural way: 
	\begin{equation*}
	(F \oplus G)(\gamma) = F(\gamma) \oplus G(\gamma) 
	\end{equation*} 
	and 
	\begin{equation*}
	(F \oplus G)(\sigma \rightsquigarrow \tau)(v, w) = (F(\sigma \rightsquigarrow \tau)v, G(\sigma \rightsquigarrow \tau)w) 
	\end{equation*}
	for $v \in F(\gamma)$ and $w \in G(\gamma)$.  
	\end{definition} \noindent 
In a similar fashion, we could define $F \otimes G$, the \textit{tensor product}\index{sheaf!tensor product} of two sheaves $F, G$, in the expected way; but there is a subtlety here when we try to think of this in terms of sheaves, and we will not be making use of this, so instead we will just indicate an example of a direct sum of sheaves. 
\begin{example} 
Consider a network, i.e., a 1-d cell complexes with oriented edges. This is just a directed acyclic graph in which each vertex has a finite `indegree' and `outdegree', and where we assume there is a special designated vertex representing the connection of the network to the "outside world." \cite{ghrist_elementary_2014} provides some nice applications of cellular sheaves over such networks, called \textit{flow sheaves}, where these represent the flow of a commodity (as in supply chains or various information or transportation of goods moving through networks). The underlying graph supports certain viable flow values, and one of the purposes of the sheaf is to encode these feasibility conditions or constraints. Basically a sheaf on such a network will supply some algebraic structure encoding a particular collection of logical, numerical, stochastic, or other constraints on the ``flows" or transport of commodities through a network. One of the advantages of using sheaves, in this setting, is that we can easily generalize beyond numerical constraints on a network to other (perhaps noisy or logical) constraints. \par 
Here is a very rough sketch of how this works.\footnote{The reader who desires a less rough sketch of these notions is invited to look at \cite{ghrist_elementary_2014}.} A \textit{flow} (or \textit{flow sheaf}) on a network $X$ is an assignment of coefficients (e.g., in $\mathbb{Z}$ or $\mathbb{N}$) to each edge of $X$, in such a way that a particular ``conservation" condition is met (namely the sum of the incoming edge flow values equals the sum of the outgoing edge values, except at the special `external' vertex, where they are not conserved). Each such value can be imagined as representing an amount of a commodity or resource in transit at a location of the underlying graph. Restrictions $F(v \rightsquigarrow e)$ are then projections onto components. The \textit{direct sum} of two flow sheaves $F \oplus G$ could then be used to represent the transportation of two \textit{different} resources being transported along the same given network. In other words, the sum $(F \oplus G)$ would represent the number of both sorts of resources, so $(F \oplus G)(e) = \mathbb{N}^2$ would represent the number of $F$-items and $G$-items being carried along the edge $e$. 
\end{example}
	\section{Sheaf Cohomology}
		\begin{example} 
		The impatient reader might well be wondering at this point: ``Okay, I understand what a sheaf is already! But what good is all this?" One glib answer, following Hubbard, might be that ``without cohomology theory, they aren't good for much"!\footnote{See \cite{hubbard_teichmuller_2006}, 383.} While this seems too pessimistic---after all, even if all sheaves on their own did was organize a wealth of particular and highly disparate constructions involving local data into a powerfully general framework, this would be immensely valuable---it is a perspective that gets at something important. If sheaves represent local data---or, more precisely, represent how to properly ensure that what is locally the case everywhere is in fact, more than that, globally the case---sheaf cohomology is a device that lets us extract global information from local data and systematically explore, represent, and relate obstructions to the extension of the local to the global. In this way, sheaf cohomology can cope with situations where the local-to-global passage breaks down, and this is of immense value, since we would like to be able to handle and talk about structures that somehow ``fall short" of assembling into sheaves. Moreover, as we mentioned at the outset of the chapter, we would like to appreciate a fact that Grothendieck\index{Alexander Grothendieck} insisted on, namely that individual sheaves are only of secondary importance---the real power of sheaf theory emerging from the use of constructions involving various sheaves, linked together via sheaf morphisms. Sheaf cohomology\index{sheaf!cohomology} will allow us to glimpse this. \par 
		In line with our categorical approach thus far, in our presentation sheaf cohomology will emerge as a \textit{functor}, specifically as one with domain the category of sheaves (together with their sheaf morphisms) and codomain the category of vector spaces. The cellular sheaves that we will continue to work with, together with their cohomology, have the nice property that it is just as easy to compute with them as to compute the ``usual" cellular cohomology of a cell complex, something one can learn about in more elementary contexts. Sheaf cohomology is particularly important to understand, though, so we do not assume that the reader already knows or recalls all they need to know about the basic notions of (co)homology. Over the next few pages, we accordingly build up to sheaf cohomology by first reviewing the basic notions of (co)homology with respect to ordinary simplices.\footnote{Though, again, technically we ought to be working with the more general \textit{cell complexes} and their cellular maps.}
	\subsection{Primer on (Co)Homology}
	Given oriented simplices (or cell complexes), as described above, a very natural thing to look at is the \textit{boundary}\index{boundary} of a given $n$-simplex. As one might expect, the boundary of a $1$-simplex $P_1P_2$ is simply the vertices of the edge, however we must now carefully attend to the issue of orientation. Taking the boundary, an operation denoted by $\partial$, is more precisely defined as taking the formal ``difference" between the endpoint and the initial point, i.e., $\partial_1(P_1P_2) = P_2 - P_1$. Similarly, the boundary of a $2$-simplex $P_1P_2P_3$ is given by 
	\begin{equation*}
	\partial_2(P_1P_2P_3) = P_2P_3 - P_1P_3 + P_1P_2.
	\end{equation*} 
	This in fact corresponds to traveling around what we intuitively think of as the boundary of a triangle in the direction indicated by the orientation arrow. The boundary of a $3$-simplex is then defined as 
	\begin{equation*}
	\partial_3(P_1P_2P_3P_4) = P_2P_3P_4 - P_1P_2P_3 + P_1P_2P_4 - P_1P_2P_3.
	\end{equation*}
	The pattern should be clear, allowing us to define the boundary operator $\partial_n$ more generally for $n > 3$: 
	\begin{equation}
	\partial_k(\sigma) = \Sigma_{i = 0}^k (-1)^i (v_0,\dots,\widehat{v_i},\dots,v_k), 
	\end{equation}
	where the oriented simplex $(v_0,\dots, \widehat{v_i}, \dots, v_k)$ is the $i$-th face of $\sigma$ obtained by deleting its $i$-th vertex. Notice that each individual summand (i.e., the positive terms) of the boundary of a simplex is just a \textit{face} of the simplex. \par  
	We can associate some \textit{groups} to a given complex $X$. The group $C_n(X)$ of oriented $n$-chains of $X$ is defined to be the free abelian group generated by the oriented $n$-simplices of $X$. Every element of $C_n(X)$ is a finite sum $\Sigma_i m_i \sigma_i$, where the $\sigma_i$ are $n$-simplices of $X$ and $m_i \in \mathbb{Z}$. Then the addition of chains is carried out by algebraically combining the coefficients of each occurrence in the chains of a given simplex. For instance, considering the surface of a tetrahedron $S$ (oriented in an obvious way), the elements of $C_2(S)$ will look like $m_1P_2P_3P_4 + m_2P_1P_3P_4 + m_3P_1P_2P_4 + m_4P_1P_2P_3$, while an element of $C_1(S)$ will look like $m_1P_1P_2 + m_2P_1P_3 + m_3P_1P_4 + m_4P_2P_3 + m_5P_2P_4 + m_6P_3P_4$, etc.\par 
	Now observe that if $\sigma$ is an $n$-simplex, then applying the boundary operator to $\sigma$ will land us inside the group of $(n-1)$-chains, i.e., $\partial_n(\sigma) \in C_{n-1}(X)$.\footnote{For the record, if we define $C_{-1}(X) = \{0\}$, the trivial group of one element, then $\partial_0(\sigma) \in C_{-1}(X)$.} But now observe that since $C_n(X)$ is a free abelian group---thus enabling us to describe a \textit{homomorphism} of such a group by specifying its values on generators---it is clear that $\partial_n$ describes a boundary homomorphism mapping $C_n(X)$ into $C_{n-1}(X)$. In other words, 
	\begin{equation}
	\partial_n \Big(\sum_i m_i \sigma_i \Big) = \sum_i m_i \partial_n (\sigma_i).
	\end{equation}
	For instance, 	 
	\begin{align*}
	\partial_1(3P_1P_2 - 4P_1P_3 + 5P_2P_4) & =  3 \partial_1 (P_1P_2) - 4\partial_1(P_1P_3) + 5\partial_1(P_2P_4) \\
	& = 3(P_2-P_1) -4(P_3-P_1) + 5(P_4-P_2) \\
	& = P_1 - 2P_2 - 4P_3 + 5P_4. 
	\end{align*} 
	But since we have a homomorphism, we are naturally drawn to look at two things: the \textit{kernel} and the \textit{image}. The kernel of $\partial_n$ will consist of those $n$-chains with boundary zero, and so the elements of the kernel are just $n$-cycles. We sometimes denote the kernel of $\partial_n$, the group of $n$-cycles, by $Z_n(X)$. So for instance, if $q = P_1P_2 + P_2P_3 +P_3P_1$, then $\partial_1(q) = (P_2 - P_1) + (P_3 - P_2) + (P_1 - P_3) = 0$. Note that this $q$ corresponds to a cycle around the triangle with vertices $P_1, P_2,P_3$ (oriented in the obvious way). Furthermore, we can consider the image under $\partial_n$, namely the group of $(n-1)$-boundaries, which consists of those $(n-1)$-chains that are boundaries of $n$-chains. We sometimes denote this group $B_{n-1}(X)$.\par 
	Homomorphisms compose and it is a well-known fact that the composite homomorphism $\partial_{n-1}\partial_{n}$ taking $C_n(X)$ into $C_{n-2}(X)$ in fact takes everything into zero, i.e., for each $c \in C_n(X)$, we have that $\partial_{n-1}(\partial_n(c)) = 0$, or $\partial^2 = 0$. A corollary of this is that $B_n(X)$ is a \textit{subgroup} of $Z_n(X)$, allowing us to form the quotient or factor group $Z_n(X)/B_n(X)$ which we denote $H_n(X)$ and call the $n$-dimensional \textit{homology group} of $X$. While perhaps obvious, it is important to realize that this quotient simply puts an equivalence relation on $Z_n$ with respect to $B_n$, i.e., $\omega \sim \sigma \iff \omega - \sigma \in B_n$, and so is technically represented as some coset. \par 
	The important thing to realize now is that we can form the \textit{sequence of chain groups} linked together by such boundary homomorphisms, a sequence we call the \textit{chain complex}: \par 
	\begin{center} 
	\begin{tikzcd}
	\cdots \arrow[r, "{\partial_{n+2}}"] & C_{n+1} \arrow[r, "{\partial_{n+1}}"] & C_n \arrow[r, "{\partial_{n}}"] & C_{n-1} \arrow[r,"{\partial_{n-1}}"] & \cdots
	\end{tikzcd} 
	\end{center} 
	\par 
	Moreover, if $C = \langle C, \partial \rangle$ is the previous (in principle doubly-infinite) sequence of abelian groups together with the collection of homomorphisms satisfying the condition that each map descends by one dimension and that $\partial^2 = 0$, then we can extend all of the above reasoning to the sequences themselves and immediately see that under these conditions the image under $\partial_k$ will be a subgroup of the kernel of $\partial_{k-1}$. In brief, then, we can define the kernel $Z_k(C)$ of $\partial_k$ as the group of $k$-cycles, the image $B_k(C) = \partial_{k+1}[C_{k+1}]$ as the group of $k$-boundaries, and then the factor group $H_k(C) = Z_k(C)/B_k(C) = \text{ker } \partial_k / \text{image } \partial_{k+1}$ as the $k$-th homology group of $C$. In other words, $H_k$ gives all the vectors that are annihilated in stage $k$ that were not already present in stage $k+1$. If for all $k$ in a sequence we have that the image under $\partial_k$ is equal to the kernel of $\partial_{k-1}$, then we have what is called an \textit{exact sequence}. While exact sequences are chain complexes, the converse is not true, since a chain complex need only satisfy that the image (of the prior map) is contained in the kernel (of the subsequent map). The important thing to realize here is that \textit{homology just measures the difference} between the image and the kernel maps. \par 
	For simplicial complexes $X$ and $Y$, a map $f$ from $X$ to $Y$ induces a mapping (i.e., homomorphism) of homology groups $H_k(X)$ into $H_k(Y)$. This arises by considering that for certain triangulations of $X$ and $Y$, the map $f$ will give rise to a homomorphism $f_k$ of $C_k(X)$ into $C_k(Y)$ which moreover commutes with $\partial_k$, i.e., $\partial_k f_k = f_{k-1}\partial_k$. \par  
	We can dualize this entire account to get an account of \textit{co}homology, and we do so briefly now since it will be important for what follows. Consider a simplicial complex $X$. For an oriented $n$-simplex $\sigma$ of $X$, we can define the \textit{coboundary}\index{coboundary} $\delta^{(n)}(\sigma)$ of $\sigma$ as the $(n+1)$-chain summing up all of the $(n+1)$-simplices $\tau$ that have $\sigma$ as a face. In other words, we are summing those $\tau$ that have $\sigma$ as a summand of $\partial_{n+1}(\tau)$. For instance, if we let $X$ be the simplicial complex consisting of the solid tetrahedron, then $\delta^{(2)}(P_3P_2P_4) = P_1P_3P_2P_4$, while $\delta^{(1)}(P_3P_2) = P_1P_3P_2 + P_4P_3P_2$.\par 
	We can also define the group $C^{(n)}(X)$ of $n$-cochains as the same as the group $C_n(X)$. However, the coboundary maps $\delta^{(n)}$ go the other way from the boundary maps, i.e., we have $\delta^{(n)}: C^{(n)} \rightarrow C^{(n+1)}$, defined by 
	\begin{equation}
	\delta^{(n)} \Big(\sum_i m_i \sigma_i \Big) = \sum_i m_i \delta^{(n)}(\sigma_i).
	\end{equation}   
	 We can then build up sequences of cochain groups into cochain complexes, just as one would expect. Cochain complexes in general can be thought of as looking at how objects are related to larger superstructures instead of to smaller substructures (as was the case for chain complexes). Just as before, we could show that $\delta^2 = 0$, i.e., that $\delta^{(n+1)}(\delta^{(n)}(c))= 0$ for each $c \in C^{(n)}(X)$. Moreover, we can define the group $Z^{(n)}(X)$ of $n$-cocycles of $X$ as the kernel of the coboundary homomorphism $\delta^{(n)}$, the group $B^{(n)}$ of $n$-coboundaries of $C^{(n)}(X)$ as the image of $\delta^{(n-1)}$, i.e., $\delta^{(n-1)}[C^{(n-1)}(X)]$, and since we have that $\delta^2 = 0$, again $B^{(n)}(X)$ will be a subgroup of $Z^{(n)}(X)$. This last fact allows us to define the $n$-dimensional \textit{cohomology group} $H^{(n)}(X)$ of $X$ as $Z^{(n)}(X)/B^{(n)}(X)$, i.e., the kernel of the map ``going out" mod the image of the map ``coming in." \par 
	 (Co)Chain complexes can be thought of as representing a cell complex within the context of linear algebra, expressing the action of taking the boundary of a cell in terms of a linear transformation. If we put this latter approach together with the development of ASCs from before, we get what is usually called simplicial (co)homology. We can thus start with some arbitrary ASC $X$ or its realization and turn it into a chain complex (for instance). This then allows us to compute its homology. Note that what we are doing here is moving \textit{via functors} from \textbf{Top} (after having placed the appropriate topology on $X$) to the category of chain complexes \textbf{Chn}, finally landing in the category \textbf{Vect}. The idea here is that we can use the algebraic properties exhibited by the composite of functors $H_{\bullet}$ to view the topological properties of the original ASC. Initially, the $C_k$ and their maps may be vector spaces over some field like $\mathbb{F}_2$ (the field of two elements, $0$ and $1$), and so the simplicial homology of $X$, denoted $H_{\bullet}(X; \mathbb{F}_2)$ will take coefficients in $\mathbb{F}_2$ (``on" or ``off"). But we could also let (co)homology take coefficients elsewhere, for instance, in $\mathbb{R}$ (thereby describing simplices' intensities, say, as opposed to the simple ``on-off" of $\mathbb{F}_2$). Eventually, the idea here is to abstract further and let (co)homology take \textit{sheaves as coefficients}.\par  
	When we start computing simplicial homology, we must choose and fix an ordering on the list of vertices in each simplex. We use coefficients other than $\mathbb{F}_2$, like $\mathbb{R}$, and the boundary maps will be used to track orientation.\footnote{Really, though, we would like to generalize beyond \textit{field} coefficients, say to $\mathbb{Z}$ coefficients, making each $C_k$ a $\mathbb{Z}$-module. In this case, the chains will record finite collections of simplices with some orientation and multiplicity. We then move to a chain complex over an $R$-module, where $R$ is a ring, the boundary maps being module homomorphisms. Then we have a sufficiently general definition: a chain complex $\mathscr{C} = (C_{\bullet}, \partial)$ is any sequence of $R$-modules $C_k$ with homomorphisms $\partial_k: C_k \rightarrow C_{k-1}$ satisfying $\partial_k \circ \partial_{k+1} = 0$.} In brief, for an arbitrary simplicial complex $X$, $C_k (X)$ will be a vector space whose dimension is the number of $k$-simplices of X, i.e., the vector space whose basis (i.e., each row as one of the $k$-simplices with some coefficients as well) is the list (given some fixed ordering) of $k$-simplices of $X$. The boundary maps $\partial_k: C_k(X) \rightarrow C_{k-1}(X)$ go down in dimension, going from the chain space built of the $k$-dimensional simplices to the chain space built of the $k-1$ simplices. In other words, the map acts on the ordered set $[v_0,\dots,v_k]$ and should be a linear map.  Indeed it takes some linear combination of $k$-simplices and returns some linear combination of $(k-1)$-simplices. We use the formula given above (with the hat) to compute with this operation.\par 
	As a very simple example illustrating these ideas, consider the following abstract simplicial complex $X = \{[v_0],[v_1],[v_0v_1]\}$, realized geometrically as 
		\begin{center} 
			\begin{tikzpicture}
			\coordinate (v0) at (0,0); 
			\coordinate (v1) at (2,0);
			\draw[thick] (v0) -- (v1) node[pos=0.5,above]{$[v_0v_1]$} node[pos=0.5, below, gray] {$\rightarrow$};
			\fill[black!20, draw=black, thick] (v0) circle (3pt) node[black, below left] {$v_0$};
			\fill[black!20, draw=black, thick] (v1) circle (3pt) node[black, below right] {$v_1$};
			\end{tikzpicture} \end{center} \par \noindent 
	Then the chain complex associated to this is given by 
	\par 
	\begin{center} 
	\begin{tikzcd}
		C_2 \arrow[r, "{\partial_{2}}"] & C_1 \arrow[r, "{\partial_{1}}"] & C_0 \arrow[r, "{\partial_{0}}"] & C_{-1} \\
		0 \arrow[r, "{\partial_{2}}"] & \mathbb{R} \arrow[r, "{\partial_{1}=\begin{psmallmatrix} -1 \\ 1 \end{psmallmatrix}}"] & \mathbb{R}^2 \arrow[r, "{\partial_{0}}"] & 0
	\end{tikzcd} \end{center} \par \noindent 
	By inspection, one can see that $H_2(X)$ (and all higher homology groups) must be zero (the trivial group). Moreover, the kernel of the $\partial_1$ map is trivial, and the image of $\partial_2$ (a $1\times 0$ matrix map) is zero, so $H_1(X)$ must be zero. As for $H_0$, on the other hand, by inspection we observe that the kernel of $\partial_0$ is everything, i.e., $\mathbb{R}^2$. The image of $\partial_1$ is spanned by the $\partial_1 = \begin{psmallmatrix} -1 \\ 1
	\end{psmallmatrix}$ matrix. Taking the quotient, then, we see that the dimension of $H_0$ must be $1$. Strictly speaking, $H_0$ is a coset, and has a vector space that is isomorphic to a single copy of $\mathbb{R}$. Thinking of it in terms of its coset representation, parameterized by one free parameter, one can think of this as \textit{identifying} the two points on account of the fact that they happen to be connected via an edge. As $H_0$ effectively represents the connected components, this should make good sense. Since the dimension of $H_1 (X)$ can be thought of as picking out the number of \textit{cycles} in the graph (i.e, among the vertices and edges) that are not filled in by 2-dimensional simplices, it should also be intuitively clear that $H_1$ ought to be trivial in this example. If we had found an $H_1$ not equal to zero, say for another simplex, this might be indicating that the simplices all fit together in some fashion, but that they cannot be glued together into one big construct, on account of some kind of obstruction or ``hole." The dimension of $H_k(X)$ is usually referred to as $k$-cycles, and for non-trivial values this can be thought of as picking out $(k+1)$-dimensional ``voids" or holes.\footnote{We could further develop this story in a number of directions, for instance going on to define the Betti numbers, which indicate various levels of obstructions, stringing them together as $k$ varies; but we leave the curious reader to pursue these matters on their own.} \par 
	We now return to simplicial or cellular maps.\footnote{We ignore more sophisticated issues here, like maps between different dimensions.} Just as $X$ can be unfolded into a chain complex $C_{\bullet}(X)$, cellular maps $f: X \rightarrow Y$ between cell complexes $X$ and $Y$ can be unfolded to yield a sequence $f_{\bullet}$ of homomorphisms $C_k(X) \rightarrow C_k(Y)$. Since $f$ is continuous, it induces a \textit{chain map} $f_{\bullet}$ that ``plays nicely" with the boundary maps of $C_{\bullet}(X)$ and $C_{\bullet}(Y)$, meaning that the following diagram is made to commute:
	\begin{center}  
	\begin{tikzcd}
		\cdots \arrow[r] & C_{n+1}(X) \arrow[r, "{\partial}"] \arrow[d, "{f_{\bullet}}"] & C_n (X) \arrow[r, "{\partial}"] \arrow[d, "{f_{\bullet}}"] & C_{n-1}(X) \arrow[r, "{\partial}"] \arrow[d, "{f_{\bullet}}"] & \cdots \\ 
		\cdots \arrow[r] & C_{n+1}(Y) \arrow[r, "{\partial'}"] & C_n (Y) \arrow[r, "{\partial'}"] & C_{n-1}(Y) \arrow[r, "{\partial'}"] & \cdots
	\end{tikzcd}
	\end{center} 
	Since the squares in this diagram commute, meaning the chain map respects the boundary operator, we have that $f$ will act not just on chains but on cycles and boundaries as well, which entails that it induces the homomorphism $H(f): H_{\bullet}(X) \rightarrow H_{\bullet}(Y)$ on homology. We are dealing with functors! In particular, the functoriality of homology means that the induced homomorphisms will reflect, algebraically, the underlying properties of continuous maps between spaces. Chain maps, in respecting the boundary operators, send neighbors to neighbors, and thus capture an essential feature of the underlying continuous maps. Via such functors, we can build up big ``quiver" diagrams with composable maps between chain complexes, i.e., between one complex representation to the next, etc. \par 
	 	Overall, the idea is then that homology should algebraically capture changes that happen to the underlying complex structure. Homology proceeds by first replacing topological spaces with complexes of algebraic objects. It then takes a hierarchically ordered sequence of these parts ``chained together," i.e., a chain complex, as input and returns the global features. Then other topological concepts---like continuous functions and homeomorphisms---have analogues at the level of chain complexes. \textit{Cohomology} is just the homology of the cochain complex.\par 
	 	In terms of the ``big picture," then, homology can be seen as a way of translating topological problems into algebraic ones. Specifically, it will (in the general approach) translate a topological problem into a problem about modules over commutative rings; but when we can take coefficients in a \textit{field}, this actually amounts to a translation of the topological problem into one of linear algebra. And this is one of the principal motivations! Linear algebra is generally much simpler than topology, in part for reasons having to do with how dimension can classify finite-dimensional vector spaces (thus the centrality of ``dimension formula" that relates the kernel of a linear transformation to its image). The ``long exact sequences" we started to look at are basically fancy versions of the dimension formula. In the context of sheaf theory, we can develop powerful way of building long exact sequences of cohomology spaces from short exact sequences of sheaves, where such sequences can already tell us a lot on their own.  
\subsection{Cohomology with Sheaves}
	We now proceed, at last, to the construction of the cochain complex for a sheaf $F$, where each $C^k$ will be comprised of the stalks, stacked together, over the $k$-simplices, and the coboundary maps (denoted with $d^k$) are built by gluing together a bunch of restriction maps. Once we have a sheaf, the tactic for computing all sections at once is to build a chain complex (where we go up in dimension, i.e., really we have a \textit{cochain} complex) and examine its zero-th homology. The difference between the (co)homology of spaces and sheaf (co)homology\index{sheaf!cohomology} mostly just has to do with the fact that with sheaves we are again looking at functions on a space, but the range of these functions is allowed to vary, i.e., we will have a collection of possible outputs, where the output space of the functions can change as we move around the domain space. \par 
	The basic idea here can be nicely explained as follows.\footnote{The next two paragraphs closely follow \cite{robinson_topological_2014}.} Suppose we have a simplicial complex, say, for simplicity 
	\begin{center} 
	\begin{tikzpicture}
	\coordinate (v1) at (0,0); 
	\coordinate (v2) at (2,0);
	\draw[thick] (v1) -- (v2);
	\fill[black!20, draw=black, thick] (v1) circle (3pt) node[black, below left] {$v_1$};
	\fill[black!20, draw=black, thick] (v2) circle (3pt) node[black, below right] {$v_2$};
	\end{tikzpicture} \end{center} \par \noindent 
with the attachment diagram above and corresponding data of the sheaf given below: 
\par 
\begin{center} 
 \begin{tikzpicture}
 \node (v1) at (-1.5,0) {$v_1$};
 \node (v2) at (3.5,0) {$v_2$};
 \node(e) at (1,0) {$e$};
 \draw[->] (v1) -- (e);
 \draw[->] (v2) -- (e);
 
  \node (F1) at (-1.5,-1) {$F(v_1)$};
  \node (F2) at (3.5,-1) {$F(v_2)$};
  \node(Fe) at (1,-1) {$F(e)$};
  \draw[->] (F1) -- (Fe) node[pos=0.5,above]{$F(v_1 \rightsquigarrow e)$};
  \draw[->] (F2) -- (Fe) node[pos=0.5,above]{$F(v_2 \rightsquigarrow e)$};
 \end{tikzpicture}.
 \end{center} 
 Now suppose that $s$ is a global section of this sheaf. Then obviously we must have 
 \begin{equation}
 F(v_1 \rightsquigarrow e) s(v_1) = s(e) = F(v_2 \rightsquigarrow e) s(v_2),
 \end{equation}
 where $F(v_i \rightsquigarrow e)$ is the restriction map, $s(v_1)$ is a section belonging to $F(v_1)$, and $s(v_2)$ is a section belonging to $F(v_2)$. But now we need only observe that this equation in fact holds in a vector space, which means that we can rewrite it $F(v_1 \rightsquigarrow e) s(v_1) - F(v_2 \rightsquigarrow e) s(v_2) = 0$, or in matrix form: 
 \begin{equation}
 \begin{psmallmatrix}
 +F(v_1 \rightsquigarrow e) \Big| - F(v_2 \rightsquigarrow e) \end{psmallmatrix} \begin{psmallmatrix}
 s(v_1) \\ s(v_2)
 \end{psmallmatrix} = 0.
 \end{equation} 
 Note that these $F(v_i \rightsquigarrow e)$ entries are just the restriction maps, so in the context of sheaves of vector spaces, they will in general be (potentially large) matrices the entries of which will be given by the linear (restriction) maps.\footnote{Note that, in the context of the (co)homology groups defined earlier in terms of equivalence relations and cosets, saying that the difference of the two restriction maps is equal to their value along the edge (i.e., $s(e)$), is effectively the same as saying that they their difference can be regarded as \textit{zero}.}\par 
 Extending this reasoning to arbitrary simplicial complexes (which we assume comes with a listing of vertices in a particular order, say lexicographic for concreteness), we can observe that computing the space of global sections of a sheaf is equivalent to computing the kernel of a particular matrix. Moreover, the matrix 
  \begin{equation*}
  \begin{psmallmatrix}
  +F(v_1 \rightsquigarrow e) \Big| - F(v_2 \rightsquigarrow e) \end{psmallmatrix}
  \end{equation*}
  generalizes into the coboundary map 
  \begin{equation*}
  \delta^{(k)}: C^{(k)}(X; F) \rightarrow C^{(k+1)}(X;F)
  \end{equation*}
  which takes an assignment $s$ on the $k$-faces to another assignment $\delta^{(k)}(s)$ whose value at a $(k+1)$-face $b$ is 
  \begin{equation}
  (\delta^{(k)}(s))(b) = \sum_{\text{all k-faces } a \text{ of } X} [b:a] F(a \rightsquigarrow b) s(a),
  \end{equation} 
  where $[b:a]$ is defined to be $0$ if $a$ (a $k$-simplex) is not a face of $b$ (a $(k+1)$-simplex) and $(-1)^n$ if you have to delete the $n$-th vertex of $b$ to get back $a$.\footnote{We start counting at 0.} This makes sense since $a$ and $b$ must differ by one dimension, so either $a$ is not a face of $b$, or $a$ is a face of $b$ in which case they will differ by exactly one vertex, i.e., one need only delete one of the vertices of $b$ to get $a$. The sign mechanism tied to the deleted vertex allows us to track and respect the chosen orientation given to the complex.\par  
 The matrix kernel construction enables us to extend this approach to higher dimensions and over simplices that are arbitrarily larger. Then we can show that we have a cellular cochain complex for a sheaf $F$ on some simplicial complex $X$. To form the cochain spaces $C^{(k)}(X;F)$, we just collect the stalks over vertices (the domain) and edges (the codomain) together via direct sum, so that an element of $C^{(k)}(X;F)$ comes from the stalk at each $k$-simplex: 
	\begin{equation}
	C^{(k)}(X;F) = \bigoplus F(a),
	\end{equation}
	where $a$ is a $k$-simplex. Moreover, the coboundary map $\delta^{(k)}: C^{(k)}(X;F) \rightarrow C^{(k+1)}(X;F)$ is defined as the block matrix where for row $i$ and column $j$, the $(i,j)$-th entry is given by $[b_i: a_j]F(a_j \rightsquigarrow b_i)$, where the $[b_i:a_j]$ term is either $0, +1$, or $-1$, depending on the relative orientation of $a_j$ and $b_i$, assuming one is a face of the other. Carrying on in this way, we have defined the cellular cochain complex: 
	\par 
	\begin{tikzcd}
		\cdots \arrow[r, "{\delta^{k-2}}"] & C^{k-1}(X;F) \arrow[r, "{\delta^{k-1}}"] & C^{k}(X;F) \arrow[r, "{\delta^{k}}"] & C^{k+1}(X;F) \arrow[r,"{\delta^{k+1}}"] & \cdots
	\end{tikzcd} \par \noindent 
	Using the cellular sheaf cohomology group definition, 
	\begin{equation*} 
	H^k(X;F) = \text{ker }\delta^k / \text{image }\delta^{k-1}, 
	\end{equation*} 
	we will recover all the cochains that are \textit{consistent} in dimension $k$ (kernel part) but that did not yet show up in dimension $k-1$ (image part). \par 
	All of this is perhaps better illustrated with our running example. 
	\begin{example} 
	We reproduce our running example below for convenience: 
	\begin{center}
	\begin{tikzpicture}[yscale=1.4, xscale=1.25]
	\node (a) at (3.75,4.5) {$\mathbb{R}^2$};
	\node (d) at (2.2,1.65) {$\mathbb{R}$};
	\node (c) at (4.7,0) {$\mathbb{R}^2$};
	\node (b) at (6,2) {$\mathbb{R}^3$};
	\node (e) at (-0.1,0) {$\mathbb{R}^3$};
	\node (f) at (-2.5,0) {$\mathbb{R}^3$};
	\node (ab) at (5.1,3) {$\mathbb{R}^2$};
	\node (ad) at (2.8,2.7) {$\mathbb{R}$};
	\node (cd) at (3.5,0.85) {$\mathbb{R}^2$};
	\node (bc) at (5.3,0.9) {$\mathbb{R}$};
	\node (ed) at (1.1,0.8) {$\mathbb{R}^2$};
	\node (cde) at (2.3,0.75) {$\mathbb{R}$};
	\node (ef) at (-1.4,0) {$\mathbb{R}^2$};
	\node (ce) at (2.3,0) {$\mathbb{R}^2$};
	\node (bd) at (3.9,1.8) {$\mathbb{R}$};
	\node (ac) at (4.3,2.2) {$\mathbb{R}^2$};
	\tiny 
	\draw[->] (a) -- node[label={[xshift=0.5cm]90:$\begin{psmallmatrix} 1 & 0 \\ -1 & 2 \end{psmallmatrix}$}] {}(ab);
	\draw[->] (b) -- node[label={[xshift=0.5cm]90:$\begin{psmallmatrix} 1 & 0 & 1\\ 0 & -1 & -1 \end{psmallmatrix}$}] {}(ab);
	\draw[->] (a) -- node[label={[xshift=-0.5cm]90:$\begin{psmallmatrix} 0 & -2 \end{psmallmatrix}$}] {}(ad);
	\draw[->] (d) -- node[label={[xshift=-0.1cm]90:$\begin{psmallmatrix} 1 \end{psmallmatrix}$}] {}(ad);
	\draw[->] (b) -- node[label={[xshift=0.2cm]90:$\begin{psmallmatrix} 2 & 0 & 2 \end{psmallmatrix}$}] {}(bd);
	\draw[->] (b) -- node[label={[xshift=0.4cm, yshift=-0.7cm]90:$\begin{psmallmatrix} 1 & 2 & 1 \end{psmallmatrix}$}] {}(bc);
	\draw[->] (c) -- node[label={[xshift=0.4cm, yshift=-0.6cm]90:$\begin{psmallmatrix} 1 & 1 \end{psmallmatrix}$}] {}(bc);
	\draw[->] (c) -- node[label={[xshift=-0.2cm, yshift=-0.7cm]90:$\begin{psmallmatrix} -1 & -1 \\ 3 & 1 \end{psmallmatrix}$}] {}(cd);
	\draw[->] (d) -- node[label={[xshift=0.25cm, yshift=-0.2cm]90:$\begin{psmallmatrix} 0.5 \\ 1 \end{psmallmatrix}$}] {}(cd);
	\draw[->] (d) -- node[label={[yshift=0cm]90:$\begin{psmallmatrix} -3 \end{psmallmatrix}$}] {}(bd);
	\draw[->] (a) -- node[label={[xshift=-0.3cm, yshift=-0.4cm]90:$\begin{psmallmatrix} 1 & 0 \\ 0 & 1 \end{psmallmatrix}$}] {}(ac);
	\draw[->] (c) -- node[label={[xshift=0.4cm, yshift=-0.4cm]90:$\begin{psmallmatrix} 3 & 3 \\ 1 & 1 \end{psmallmatrix}$}] {}(ac);
	\draw[->] (d) -- node[label={[xshift=-0.2cm]90:$\begin{psmallmatrix} 3 \\ 1 \end{psmallmatrix}$}] {}(ed);
	\draw[->] (e) -- node[label={[xshift=-0.2cm, yshift=-0.1cm]90:$\begin{psmallmatrix} 2 & 0 & 1\\ 0 & 3 & -1 \end{psmallmatrix}$}] {}(ed);
	\draw[->] (c) -- node[label={[yshift=-0.8cm]90:$\begin{psmallmatrix} 1 & -1 \\ -1 & 2 \end{psmallmatrix}$}] {}(ce);
	\draw[->] (e) -- node[label={[yshift=-0.8cm]90:$\begin{psmallmatrix} 2 & -3 & 2\\ 1 & 0 & 7.5 \end{psmallmatrix}$}] {}(ce);
	\draw[->] (ed) -- node[label={[yshift=-0.2cm]90:$\begin{psmallmatrix} 1 & -1 \end{psmallmatrix}$}] {}(cde);
	\draw[->] (cd) -- node[label={[yshift=-0.2cm]90:$\begin{psmallmatrix} 2 & 1 \end{psmallmatrix}$}] {}(cde);
	\draw[->] (ce) -- node[label={[xshift=0.5cm, yshift=-0.5cm]90:$\begin{psmallmatrix} 1 & 0 \end{psmallmatrix}$}] {}(cde);
	\draw[->] (e) -- node[label={[yshift=-0.8cm]90:$\begin{psmallmatrix} 2 & 0 & 2\\ 1 & -1 & 1 \end{psmallmatrix}$}] {}(ef);
	\draw[->] (f) -- node[label={[yshift=-0.8cm]90:$\begin{psmallmatrix} 0 & 1 & 1\\ 1 & -1 & 0 \end{psmallmatrix}$}] {}(ef);
	
	\end{tikzpicture}
	\end{center}  
	For such a sheaf $F$ over our given complex $X$, we have the following:  
		\begin{center} 
		\begin{tikzcd}[row sep=tiny]
			H^0(X;F) \arrow[r, "{\delta^{0}}"] & H^1(X;F) \arrow[r, "{\delta^{1}}"] & H^{2}(X;F) \arrow[r, "{\delta^{2}}"] & H^{3}(X;F) \\ \mathbb{R}^2 \oplus \mathbb{R}^3 \oplus \mathbb{R}^2 \oplus & \mathbb{R}^2 \oplus \mathbb{R}^2 \oplus \mathbb{R} \oplus \mathbb{R} \oplus \\ 
			\mathbb{R} \oplus \mathbb{R}^3 \oplus \mathbb{R}^3  \arrow[r, "{\delta^{0}}"]  & \mathbb{R} \oplus \mathbb{R}^2 \oplus \mathbb{R}^2 \oplus \mathbb{R}^2 \oplus \mathbb{R}^2 \arrow[r, "{\delta^{1}}"] & \mathbb{R} \arrow[r, "{\delta^{2}}"] & 0 \\
			\begin{psmallmatrix}
			s(a) \\ s(b) \\ s(c) \\ s(d) \\ s(e) \\ s(f)
			\end{psmallmatrix} \arrow[r, "{\delta^{0}}"] & \begin{psmallmatrix}
			s(ab) \\ s(ac) \\ s(ad) \\ s(bc) \\ s(bd) \\ s(cd) \\ s(ce) \\ s(de) \\ s(ef)
		\end{psmallmatrix} \arrow[r, "{\delta^{1}}"] & \begin{psmallmatrix} s(cde) \end{psmallmatrix} \arrow[r, "{\delta^{2}}"] & 0
		\end{tikzcd} 
		\end{center}
		and where $\delta^0$ is given by 
\vspace{-0.5cm}
\begin{center} 
\tiny 
\[
\renewcommand\arraystretch{1.4}
\begin{blockarray}{rccccccccc}
&&& a & b & c & d & e & f & \\
&&& \downarrow & \downarrow & \downarrow & \downarrow & \downarrow & \downarrow & \\
\begin{block}{rc[c@{}c|c|c|c|c|c@{}c]}
\text{[ab]}  & \to && -F(a \rightsquigarrow ab) & F(b \rightsquigarrow ab) & 0 & 0 & 0 & 0 &\vphantom{\smash[b]{\bigg|}} \\
\BAhhline{~~~------~}
\text{[ac]}  & \to && -F(a \rightsquigarrow ac) & 0 & F(c \rightsquigarrow ac) & 0 & 0 & 0 & \\
\BAhhline{~~~------~}
\text{[ad]}  & \to && -F(a \rightsquigarrow ad) & 0 & 0 & F(d \rightsquigarrow ad) & 0 & 0 & \\
\BAhhline{~~~------~}
\text{[bc]}  & \to && 0 & -F(b \rightsquigarrow bc) & F(c \rightsquigarrow bc) & 0 & 0 &  0 &\\
\BAhhline{~~~------~}
\text{[bd]} & \to && 0 & -F(b \rightsquigarrow bd) & 0 & F(d \rightsquigarrow bc) & 0 & 0 & \vphantom{\smash[t]{\bigg|}}\\
\BAhhline{~~~------~}
\text{[cd]}  & \to && 0 & 0 & -F(c \rightsquigarrow cd) & F(d \rightsquigarrow cd) & 0 & 0 & \\
\BAhhline{~~~------~}
\text{[ce]}  & \to && 0 & 0 & -F(c \rightsquigarrow ce) & 0 & F(e \rightsquigarrow ce) & 0 & \\
\BAhhline{~~~------~}
\text{[de]}  & \to && 0 & 0 & 0 & -F(d \rightsquigarrow de) & F(e \rightsquigarrow de) & 0 & \\
\BAhhline{~~~------~}
\text{[ef]}  & \to && 0 & 0 & 0 & 0 & -F(e \rightsquigarrow ef) & F(f \rightsquigarrow ef) & \\
\end{block}
\end{blockarray}
\]
\end{center}
and $\delta^1$ by 
\begin{center} 
	\tiny
	\[
	\renewcommand\arraystretch{0.9}
	\begin{blockarray}{rccccccccccc}
	&& [ab] & [ac] & [ad] & [bc] & [bd] & [cd] & [ce] & [de] & [ef] \\
	&& \downarrow & \downarrow & \downarrow & \downarrow & \downarrow & \downarrow & \downarrow & \downarrow & \downarrow \\
	\begin{block}{rc[c|c|c|c|c|c|c|c|c@{}c]}
	\text{[cde]} \to && 0 & 0 & 0 & 0 & 0 & F(cd \rightsquigarrow cde) & -F(ce \rightsquigarrow cde) & F(de \rightsquigarrow cde) & 0 &\vphantom{\smash[b]{\bigg|}} \\
	\end{block}
	\end{blockarray}
	\]
\end{center} 
and $\delta^2$ is trivial. Computing the cohomologies here will require doing some linear algebra and row-reduction, which is left to the reader. Observe that in general by adding in more higher-dimensional consistency checks, i.e., ``filling in" data corresponding to the ``missing" higher-dimensional simplices of the underlying simplex, we would be able to reduce the kernel of $\delta^1$, and thus get closer to reducing the non-trivial $H^1$. 
\end{example}  
It is worth emphasizing a more general truth, namely that the space of global sections of a sheaf $F$ on a cell complex $X$ will be isomorphic to $H^0(X; F)$. Moreover, since this $H^k$, i.e., cohomology with sheaves as coefficients, is a functor, this means that when we have sheaf morphisms between sheaves, these will induce linear maps between their cohomology spaces, allowing us to extend this still further. Ultimately, it is also worth noting that while we have been focused on cellular sheaves, \textit{general sheaf cohomology} can be shown to be isomorphic to cellular sheaf cohomology, so this sort of example is really part of a much more general story. And the vector space version of homological algebra we have used in order to provide a concrete example, while useful for many applications and for building intuition, does not display the full power: ultimately one would like to (and could) retell the story using rings, modules, and other categories. \par 
Via cohomology, global compatibilities between pieces of local data can display global qualitative features of the data structure. Via sheaf (co)homology, we might, for instance, be better able to isolate potential ``holes" in data collections. At a very high level, the general idea of all this is suggested by the following picture, illustrating how we get a comparatively simple algebraic representation of features of spaces: 
	\begin{center}
	\hspace*{-3em}
		\includegraphics*[scale =0.37]{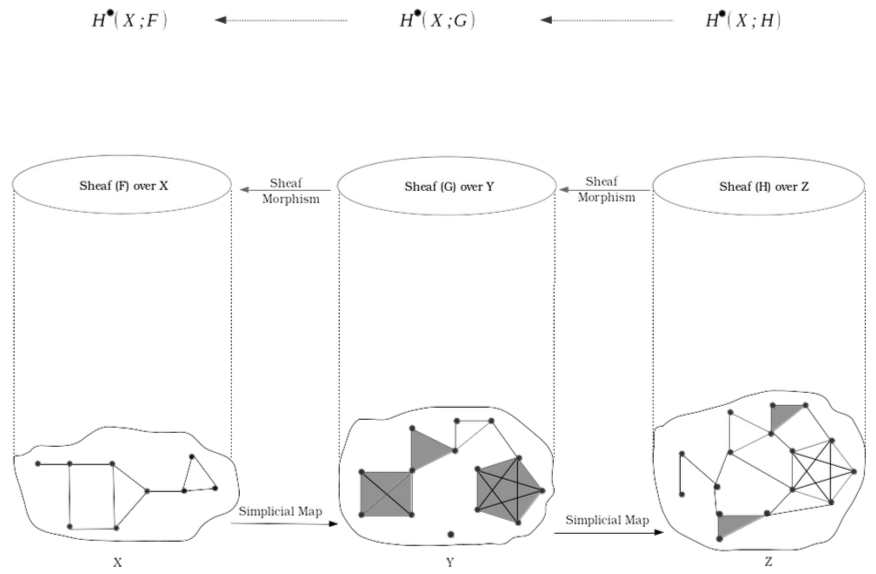}
	\end{center}  
\end{example}
\subsection{Philosophical Pass: Sheaf Cohomology}
If the sheaf compatibility conditions require controlled transitions from one local description to another, enabling progressive patching of information over overlapping regions until a unique value assignment emerges over the entire region---something that is captured by the vanishing of the group $H^0$ (yielding our global sections)---higher (non-vanishing) cohomology groups basically detect and summarize (in an algebraic fashion) obstructions to such local patching and consistency relations among various dimensional subsystems. In other words, it can be thought of as measuring (for some cover) how many incompatible (purely local) systems we would have to ``discard" in order to be left with only the compatible systems. In this way, sheaf cohomology\index{sheaf!cohomology} moreover allows us to examine the relationship between information valid globally and the underlying topology of the space. \par 
A proper discussion of the possible invariants that emerge in cohomology, especially as we ascend in dimension, would require a much longer and more detailed discussion. Instead, here are some very general reflections on the idea of sheaf cohomology. Referencing the high-level picture on the previous page, in forming sheaves (middle level) over the discrete approximations of spaces via their triangulations (bottom level), a more ``continuous" perspective is recaptured. However, the non-vanishing cohomology groups (top level) give an algebraic (more ``discrete") representation of something like the \textit{resistance} of certain information (assigned to a part of a space) to integration into a more global system. In short, if the collation condition in the sheaf construction aligns them with ``continuity" in the sense that it ensures smooth passage from the local to the global, cohomology with sheaves is something like its discrete counterpart providing us an algebraic summarization of when and how such local-global passages might be blocked. In this respect, sheaf cohomology could intuitively be thought of as capturing, in a dialectic between the continuous and discrete, the non-globalizability or non-extendibility of a given information structure in relation to other overlapping data structures. Both in its algebraic representation and in this general interpretation, then, the non-vanishing cohomology groups might be thought of as giving us a picture of just how ``non-integrated" a system of information over a space may be. On the other hand, vanishing cohomology groups indicate the mutual compatibility or ``globalizability" of local information systems (since they tell us about the global sections). In this way, sheaf cohomology emerges as a tool for representing (algebraically) what might be thought of as the degree of generality (or lack thereof) of a given system of measurement or interlocking ways of assigning information to a space. \par 
Some years before the invention of sheaf theory, Charles Peirce\index{Charles Peirce} argued that ``continuity is shown by the logic of relations to be nothing but a higher type of that which we know as generality. It is relational generality."\footnote{See \cite{peirce_collected_1997}, CP, 6.190.} Such a suggestive, if somewhat cryptic, remark provokes us to take a closer look at the connections between generality and continuity that emerge in the context of sheaves. We know that a sheaf enables a collection of local sections to be patched together uniquely given that they agree (or that there exists a translation system for making them agree) on the intersections. Consider the satellite image ``mosaic"\index{mosaic} sheaf introduced in Chapter 3. Recall the way in which the sheaf (collation/gluing) condition ensures a systematic passage from local sections (images of parts of the glacier) to a unique global section (the image of the entire glacier). Where the localizing step of the sheaf construction might be thought of as analytic, decomposing an object into a multitude of individual parts (local), the gluing steps are synthetic in restoring systematic relations between those parts and thereby securing a unique assignment over the entire space (global). The global sections of such a sheaf should not be thought of as a single (topmost) image, but rather as the entire network of component parts welded together via certain compatibility relations or constraints. In this connection, `generality' can be understood in terms of the systematic passages from the local to the \textit{global}, a passage that is strictly \textit{relational}, in that the action of the component restriction maps is precisely an enforcing of certain relations or mutual constraints between the local sections (that are then built up, along the lines of these relations, into a global section), and these are an ineliminable part of the construction. \par 
One might further think of the indexing (domain) category in the (pre)sheaf construction as providing a particular `context' specifying the possible scope of the generality of a given sheaf diagram---just \textit{how} global the global section is, in relation to other possible spaces. In this way, the degree of generality achieved by a particular sheaf construction can be thought to depend upon the form of the indexing category. In so far as such distinct systems for the production of generalities can be themselves compared via natural transformations, one might also think of this as introducing yet another (higher-order) layer of relationality into the notion of generality. In these ways, via the sheaf concept, Peirce's suggestive idea that 
\begin{equation*}
\text{Continuity} = \text{Relational Generality}
\end{equation*}  
is given a particularly powerful interpretation. 
\subsection{A Glimpse into Cosheaves}
In the cellular case, we can easily talk about \textit{sheaf homology} by just reversing the direction of the arrows, technically producing a \textit{cosheaf}\index{cosheaf} (with its ``corestriction" maps). Simply by observing that the vector space dual of every extension/corestriction map in a cosheaf produces a sheaf over the poset, we arrive at homology for a cosheaf. More explicitly, this reversal gives us a \textit{cosheaf} $\hat{F}$ of vector spaces on a complex, assigning to each simplex $\sigma$ a vector space $\hat{F}(\sigma)$ but to each face attachment $\sigma \rightsquigarrow \tau$ of $\tau$ a corestriction (or extension)\index{corestriction} map, i.e., a linear transformation $\hat{F}(\sigma \rightsquigarrow \tau): \hat{F}(\tau) \rightarrow \hat{F}(\sigma)$, which reverses the direction of the sheaf maps. As before, this cosheaf must respect the composition and identity rules. Similar to the sheaf case, a cosheaf can come to serve as a system of coefficients for homology that varies as the space varies. However, the globality of the cosheaf sections will be found in the top dimension (unlike how the global sections were built up from the \textit{vertices} in the case of cellular sheaves).\footnote{More generally, following \cite{curry_sheaves_2013}, we define a \textit{pre-cosheaf} as expected: namely as a functor $\hat{F}: \mathscr{O}(X) \rightarrow \textbf{D}$ (when we come to the cosheaf conditions, technically we ought to insist that $\textbf{D}$ be not just any category, but rather a category with enough (co)limits), and specify that whenever $V \subseteq U$, we can define the extension or co-restriction map for the cosheaf as $r_{U,V}: \hat{F}(V) \rightarrow \hat{F}(U)$. We again need the notion of open cover, a notion we can now think of in terms of a function from the \textit{nerve} construction, which basically acts on covers to produce the ASC consisting only of those finite subsets of the cover whose intersection is non-empty. More formally, if we suppose $\mathscr{U} = \{U_i\}$ is an open cover of $U$, then we can take the \textit{nerve} of the cover to yield an ASC $N(\mathscr{U})$, which will have for elements the subsets $I = \{i_0, \dots, i_n \}$ for which it holds that $U_I = U_{i_0} \cap \cdots \cap U_{i_n} \neq \emptyset$. $N(\mathscr{U})$ is then the category with objects the finite subsets $I$ where $U_I \neq \emptyset$, and unique arrows from $I$ to $J$ whenever $J \subseteq  I$. Finite intersections of opens are open, so we thus get the functors $\iota_{\mathscr{U}}: N(\mathscr{U}) \rightarrow \mathscr{O}(X)$ and $\iota_{\mathscr{U}}^{op}: N(\mathscr{U})^{op} \rightarrow \mathscr{O}(X)^{op}$. In general, as the colimit of a cover $N(\mathscr{U}) \rightarrow \mathscr{O}(X)$ is the union $U = \bigcup_i U_i$, the data we associate to $U$ ought to be expressible as the colimit of data assigned to the nerve. With this expectation in mind, we can say that $\hat{F}$ is a \textit{cosheaf on} $\mathscr{U}$ if the unique map from the colimit of $\hat{F} \circ \iota_{\mathscr{U}}$ to $\hat{F}(U)$, given by 
	\begin{equation}
	\hat{F}[\mathscr{U}] := \varinjlim_{I \in N(\mathscr{U})} \hat{F}(U_I) \rightarrow \hat{F}(U)
	\end{equation} 
	is in fact an isomorphism. Then $\hat{F}$ is simply a \textit{cosheaf} if for every open set $U$ and every open cover $\mathscr{U}$ of $U$, the map $\hat{F}(\mathscr{U}) \rightarrow \hat{F}(U)$ is an isomorphism. For simplicity, suppose $\textbf{D} = \textbf{Set}$ and take a cover $\mathscr{U} = \{U_1, U_2\}$ of $U$ by just two open sets. The sheaf condition would stipulate that for two functions or sections $s_1$, $s_2$ (from $F(U_1)$, $F(U_2)$ respectively) to give an element in $U = U_1 \cup U_2$, the sections must agree on the overlap $U_1 \cap U_2$. This constraint serves to pick out the consistent choices of elements over the local sections that can then be glued together into a section over the larger set. With a similar setup, i.e., $\textbf{D} = \textbf{Set}$ and $U = U_1 \cup U_2$, the cosheaf condition requires not that we find consistent choices, but rather that we use \textit{quotient objects}. We do indeed form the union of the two sections, but in the process we identify those elements that would be double-counted on account of coming from the intersection. Formally, 
	\begin{equation}
	\hat{F}(U) \cong (\coprod_{i = 1,2} \hat{F}(U_i)) / \sim 
	\end{equation}
	where $s_1 \sim s_2$ iff there exists an $s_{12}$ (a section over the intersection) such that $s_1 = r_{U_1, U_{12}}(s_{12})$ and $s_2 = r_{U_2, U_{12}}(s_{12})$. This makes sense, since in accordance with duality, we would expect that the equalizer definition of the sheaf condition would be converted, in passing to cosheaves, into an underlying coequalizer diagram. For more details on this, see \cite{curry_sheaves_2013}.} 
	\par 
A sheaf was constructed in such a way that the values of its sections on larger sets in the Alexandrov topology will determine values on smaller sets. A cosheaf basically reverses this dependence. While so far we have seen many examples of sheaves in contexts where it makes sense to perform \textit{restrictions} of assignments of data from larger spaces to data over smaller spaces, building up global assignments from the local pieces, i.e., a ``bottom-up" approach, roughly cosheaves can be thought of as proceeding ``top-down," involving \textit{extensions} of data given over smaller spaces to larger spaces.\footnote{In the context of simplices, because of the reversal of direction involved in the topology given on the face relation construction, such extensions will go from higher dimensional simplices to lower.} While in some sense the paradigmatic example of a sheaf was given by the restriction of continuous functions, the paradigmatic example of a cosheaf might be given by the cosheaf of compactly supported continuous functions where, instead of restricting along inclusions, we extend by zero (in the other direction).\par 
Importantly, while in the cellular context the difference between sheaf and cosheaf is somewhat immaterial, simply a matter of which direction makes the most sense for the framing of the problem, things can be far more subtle in the context of sheaves and cosheaves over opens sets for a continuous domain. In insisting upon the more general functorial perspective, allowing us to make use of duality, one might suspect that the differences between sheaves and cosheaves are merely formal and not worth discussing. For instance, a sheaf is a particular functor that commutes with limits in open covers. As one might expect, a cosheaf is a functor that preserves colimits in open covers. However, in more general contexts than the cellular one, especially with open sets coming from a continuous domain, the differences can reflect much more than a preference for direction of arrows.\footnote{For instance, when working with the Alexandrov topology on a poset (or when working with locally finite topological spaces), we can ignore the distinction between pre(co)sheaves and (co)sheaves; however, while for general topological spaces, there is a sheafification\index{sheafification} functor that allows us to pass from a presheaf to the unique smallest sheaf consistent with the given presheaf, there is no analogous cosheafification functor for general topological spaces. For more on cosheaves, and a number of interesting connections and differences with sheaves, we refer the reader to Justin Curry's thesis, \cite{curry_sheaves_2013}. This reference, as well as \cite{robinson_sheaf_2016}, contains more details on some of the dualities in the sheaf-cosheaf perspective, as well as instances of asymmetry (when certain constructions are natural for sheaves but not for cosheaves).} \par   
The next example explores a particularly fascinating connection between sheaves and cosheaves in the context of probabilities and Bayes nets.   
\begin{example}  
Imagine we are given a set of random variables $X_0, X_1, \dots, X_n$.\footnote{The idea for this example was inspired by \cite{robinson_sheaf_2016}.} We can consider the set $P(X_0, X_1,\dots, X_n)$ of all joint probability distributions over these random variables, i.e., the non-negative measures or generalized functions with unit integral. Now, there is a very natural map, one that will be familiar to anyone with some exposure to probability theory: 
\begin{equation}
P(X_0, X_1, \dots, X_n) \rightarrow P(X_0, X_1, \dots, X_{n-1}).
\end{equation}
This map is accomplished via \textit{marginalization}, i.e., we have 
\begin{equation}
f(X_0, X_1,\dots, X_{n-1}) = \int f(X_0, X_1, \dots, X_n) d X_n,
\end{equation}
and where there are similar maps for marginalizing out the other random variables. The important point is that this in fact yields a cosheaf on the complete $n$-simplex. We elaborate on this with an example. \par 
Assume given the random variables $X_0 = W$, $X_1 = S$, $X_2 = R$ (the reason for renaming these random variables thus will be clear in a moment). Then the space of probability measures $P(W,S,R)$ is a function from the cartesian product of the random variables to the (non-negative) reals such that the integral is zero. Now, we can perform the marginalization operation, for instance, marginalizing out the $R$ by integrating along $R$, yielding a map $P(W,S,R) \rightarrow P(W,S)$. We can do this for each of the random variables, and continue down in dimension until we reach the measurable functions over a single random variable. In other words, we have: 
	\begin{center} 
	\small
	\begin{tikzcd} 
		&  P(W,S,R) \arrow[bend left]{dr} \arrow[bend right]{dl}
		\arrow{d}
		\\
		
		P(S,R) \arrow{d} \arrow[dr] 
		& P(W,R) \arrow{dl} \arrow[dr] &  
		P(W,S) \arrow[dl] \arrow{d} \\ 
		
		P(R)
		&  
		P(S) 
		& P(W)
	\end{tikzcd}
\end{center} 
which in fact represents the attachment diagram of a complete $2$-simplex, 
\begin{center} 
\begin{tikzpicture}

\coordinate (b) at (3,1);
\coordinate (c) at (5,0);
\coordinate (e) at (1,0);

\draw[thick, fill=black!20] (b) -- (c) -- (e) -- cycle;

\fill[black!20, draw=black, thick] (b) circle (3pt) node[black, above] {$W$};
\fill[black!20, draw=black, thick] (c) circle (3pt) node[black, below right] {$S$};
\fill[black!20, draw=black, thick] (e) circle (3pt) node[black, below left] {$R$};

\end{tikzpicture} 
\end{center}
The commutative diagram, together with the appropriate marginalization maps, gives a cosheaf on this abstract simplicial complex.\par 
Now, the reader may have already wondered about maps going the other way, namely: 
\begin{equation}
P(X_0,X_1,\dots, X_{n-1}) \rightarrow P(X_0,X_1,\dots,X_n),
\end{equation}
an operation that is parameterized by functions $C$
\begin{equation}
P(X_0, X_1,\dots, X_n) = C(X_0,X_1,\dots,X_n) f(X_0,X_1,\dots,X_{n-1}),
\end{equation}
where usually one writes the arguments to $C$ as follows 
\begin{equation}
C(X_n \vert X_0, X_1, \dots, X_{n-1}).
\end{equation}
The reader may recognize that we are just describing \textit{conditional probabilities} and Bayes's rule. The key observation is that such conditional probability maps yield a sheaf over a portion of the $n$-simplex.\par 
The reader may be familiar with the construction of a \textit{Bayes net},\index{Bayes net} given by a directed acyclic graph (with an induced topology) together with local conditional probabilities. A Bayes net encodes joint distributions and does so as a product of local conditional distributions, i.e., 
\begin{equation}
P(X_1, X_2, \dots, X_n) = \prod_{i=1}^{n} P(X_i \lvert \text{ parents} (X_i) ).
\end{equation} 
For the sake of concreteness, we consider the following very simple example of a Bayes net, the standard one given in most introductions to the device, involving probabilities of ``grass being wet" given that it rained or that the sprinkler was running, etc., illustrated with some sample probability assignments: 
\begin{center} 
\includegraphics*[scale=0.3]{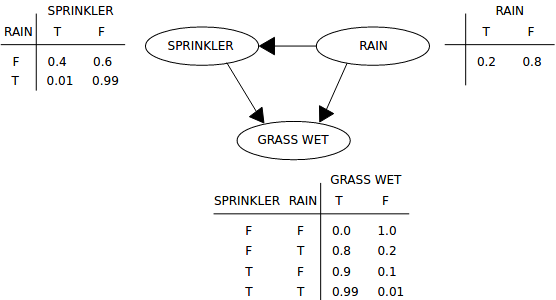}
\end{center}  
The perhaps surprising result is that the content of this Bayes net is in fact entirely captured by the paired sheaf-cosheaf construction given below, where the marginalization cosheaf\index{cosheaf} is given by the entire diagram (arrows going down the page), while the paired conditional probability sheaf is given in bold (up the page) over a part of the underlying attachment diagram.    
	\begin{center} 
		\tiny
		\begin{tikzcd}[cramped, column sep = -12em, row sep=huge, xscale = 0.6] 
			&  \begin{psmallmatrix} P(w,s,r) \\ P(w,s,\neg r) \\ P(w,\neg s,r) \\ P(w,\neg s,\neg r) \\ P(\neg w,s,r) \\ P(\neg w,s,\neg r) \\ P(\neg w,\neg s,r) \\ P(\neg w,\neg s,\neg r) \end{psmallmatrix} \arrow[bend left]{ddddr}{\begin{psmallmatrix} 1 \amsamp 1 \amsamp 0 \amsamp 0 \amsamp 0 \amsamp 0 \amsamp 0 \amsamp 0 \\ 0 \amsamp 0 \amsamp 1 \amsamp 1 \amsamp 0 \amsamp 0 \amsamp 0 \amsamp 0 \\ 0 \amsamp 0 \amsamp 0 \amsamp 0 \amsamp 1 \amsamp 1 \amsamp 0 \amsamp 0 \\ 0 \amsamp 0 \amsamp 0 \amsamp 0 \amsamp 0 \amsamp 0 \amsamp 1 \amsamp 1 \end{psmallmatrix}}
			
			\arrow[bend right, swap]{ddddl}{\begin{psmallmatrix} 1 \amsamp 0 \amsamp 0 \amsamp 0 \amsamp 1 \amsamp 0 \amsamp 0 \amsamp 0 \\ 0 \amsamp 1 \amsamp 0 \amsamp 0 \amsamp 0 \amsamp 1 \amsamp 0 \amsamp 0 \\ 0 \amsamp 0 \amsamp 1 \amsamp 0 \amsamp 0 \amsamp 0 \amsamp 1 \amsamp 0 \\ 0 \amsamp 0 \amsamp 0 \amsamp 1 \amsamp 0 \amsamp 0 \amsamp 0 \amsamp 1 \end{psmallmatrix}} 
			\arrow{dddd}{\begin{psmallmatrix} 1 \amsamp 0 \amsamp 1 \amsamp 0 \amsamp 0 \amsamp 0 \amsamp 0 \amsamp 0 \\ 0 \amsamp 1 \amsamp 0 \amsamp 1 \amsamp 0 \amsamp 0 \amsamp 0 \amsamp 0 \\ 0 \amsamp 0 \amsamp 0 \amsamp 0 \amsamp 1 \amsamp 0 \amsamp 1 \amsamp 0 \\ 0 \amsamp 0 \amsamp 0 \amsamp 0 \amsamp 0 \amsamp 1 \amsamp 0 \amsamp 1 \end{psmallmatrix}}
			\\ \\ \\ \\ 
			
			\begin{psmallmatrix} P(w,s,r) + P(\neg w,s, r) \\ P(w,s,\neg r) + P(\neg w, s,\neg r) \\ P(w,\neg s,r) + P(\neg w,\neg s,r) \\ P(w,\neg s,\neg r) + P(\neg w,\neg s,\neg r)  \end{psmallmatrix} \arrow[swap]{ddddd}{\begin{psmallmatrix} 1 \amsamp 0 \amsamp 1 \amsamp 0 \\ 0 \amsamp 1 \amsamp 0 \amsamp 1 \end{psmallmatrix}} \arrow[ddddddr] \arrow[ultra thick]{uuuur}{P(W|S,R)}
			& \begin{psmallmatrix} P(w,s,r) + P(w,\neg s, r) \\ P(w,s,\neg r) + P(w, \neg s,\neg r) \\ P(\neg w,s,r) + P(\neg w,\neg s,r) \\ P(\neg w, s,\neg r) + P(\neg w,\neg s,\neg r) \end{psmallmatrix} \arrow[swap]{dddddl}{\begin{psmallmatrix} 1 \amsamp 0 \amsamp 1 \amsamp 0 \\ 0 \amsamp 1 \amsamp 0 \amsamp 1 \end{psmallmatrix}} \arrow{dddddr}{\begin{psmallmatrix} 1 \amsamp 1 \amsamp 0 \amsamp 0 \\ 0 \amsamp 0 \amsamp 1 \amsamp 1 \end{psmallmatrix}}
			\arrow[ultra thick, bend left]{uuuu}{P(S|W,R)} &  
			\begin{psmallmatrix} P(w,s,r) + P(w, s, \neg r) \\ P(w,\neg s,r) + P(w, \neg s,\neg r) \\ P(\neg w,s,r) + P(\neg w, s,\neg r) \\ P(\neg w, \neg s, r) + P(\neg w,\neg s,\neg r)  \end{psmallmatrix} \arrow[ddddddl] \arrow{ddddd}{\begin{psmallmatrix} 1 \amsamp 1 \amsamp 0 \amsamp 0 \\ 0 \amsamp 0 \amsamp 1 \amsamp 1 \end{psmallmatrix}} \\ \\ \\ \\ \\
			
			\begin{psmallmatrix} P(w,s,r) + P(w,\neg s, r) + P(\neg w,s,r) + P(\neg w, \neg s,r) \\ P(w,s,\neg r) + P(w,\neg s,\neg r) + P(\neg w, s,\neg r) + P(\neg w,\neg s,\neg r)  \end{psmallmatrix} \arrow[bend left = 37ex, ultra thick]{uuuuu}{P(S|R)} \arrow[bend right = 30ex, ultra thick, swap]{uuuuur}{P(W|R)}
			&  & 
			\begin{psmallmatrix} P(w,s,r) + P(w,s, \neg r) + P(w,\neg s,r) + P(w, \neg s,\neg r) \\ P(\neg w,s, r) + P(\neg w, s,\neg r) + P(\neg w, \neg s, r) + P(\neg w,\neg s,\neg r)  \end{psmallmatrix} 
			\\
			& \begin{psmallmatrix} P(w,s,r) + P(\neg w,s, r) + P(w,s,\neg r) + P(\neg w, s,\neg r) \\ P(w,\neg s, r) + P(\neg w, \neg s, r) + P(w, \neg s, \neg r) + P(\neg w,\neg s,\neg r)  \end{psmallmatrix}
		\end{tikzcd}
	\end{center} 
 The paired sheaf-cosheaf construction contains all the data of a Bayes net; and, in fact, a solution to the Bayes net is just a global section that is a section for \textit{both the sheaf and the cosheaf}! \par  
 Before ending this chapter and moving into discussion of toposes, we take the opportunity to make an important but frequently overlooked observation. Consider the (co)sheaf construction above. Now, consider that relatively small Bayes nets, say, one with only 8 or 9 nodes, are already rather simple compared to those that will be of use in practice. One might thus be suspicious of just how complicated the corresponding (co)sheaf might look for even only slightly more involved examples, not to mention the issue of storing the relevant sections for such sheaves. This indeed seems to be a real issue. One might also suspect that computing the sheaf cohomology (and global sections) on ``monster" (extremely large) sheaves would be extremely difficult. As the discussion of this Bayes net sheaf-cosheaf construction suggests, most ``real-life" (co)sheaves may very well turn out to be monsters in the sense of being so large as to cause difficulties in storage, representation, or computation---difficulties we have simply avoided by confining our attention to more modest constructions. This is an issue that deserves to be recognized and pondered.\footnote{For some ways to reduce the difficulty, in one setting, see \cite{smith_introduction_2014} (527-533); see \cite{curry_sheaves_2013} (64ff.) for some ways to think about pre-processing the input data so as to deal with the ``too many sections" problem. The reader may also find \cite{curry_discrete_2015} highly relevant, as this shows how you can ``collapse" the data structure if your restriction maps are nice enough. Thanks are due to Michael Robinson for pointing me in the direction of this paper and observing the connection.}       
\end{example}
\printbibliography 
\printindex
\end{document}